\documentclass[11pt,reqno]{amsart}
\usepackage[margin=1in]{geometry}
\usepackage{amsmath,amssymb}
\usepackage{algpseudocode}
\usepackage{algorithm}
\usepackage{mathrsfs}
\usepackage{algorithmicx}
\usepackage{caption}
\usepackage{subcaption}
\usepackage{amsthm}
\usepackage{esint}
\numberwithin{equation}{section}
\usepackage{amsfonts}
\usepackage{graphicx}
\usepackage{hyperref}
\usepackage[capitalise]{cleveref}
\usepackage{makecell}
\usepackage{longtable}
\usepackage{comment}
\usepackage{enumerate}

\usepackage{setspace}
\usepackage{xcolor}
\usepackage{fancyhdr}
\usepackage{listings}
\usepackage{cite}
\usepackage{tabularx}
\usepackage{epstopdf}    
\usepackage{epsfig}
\usepackage{float}
\usepackage{listings}
\usepackage{appendix}
\theoremstyle{plain}

\newtheorem{theorem}{Theorem}[section]
\newtheorem{proposition}[theorem]{Proposition}

\newtheorem{lemma}[theorem]{Lemma}
\newtheorem{cor}[theorem]{Corollary}
\theoremstyle{definition}
\newtheorem{definition}[theorem]{Definition}

\newtheorem{remark}[theorem]{Remark}

\let\pa=\partial
\let\al=\alpha
\let\b=\beta
\let\d=\delta
\let\g=\gamma
\let\e=\varepsilon
\let\ze=\zeta
\let\kp =\kappa
\let\lam=\lambda

\let\s=\sigma
\let\f=\frac
\let \les = \lesssim
\let \gtr = \gtrsim
\let\om=\omega

\let\th =\theta
\let\Th =\Theta

\let\pr = \prime
\let \vp =\varphi

\let\B = \Big
\let\D=\Delta
\let\Lam=\Lambda
\let\S=\Sigma

\let\Om=\Omega
\let\td = \tilde

\let\wh=\widehat
\let \mr = \mathring

\let\pa=\partial

\def\cA{{\mathcal A}}
\def\cB{{\mathcal B}}

\def\cD{{\mathcal D}}
\def\cE{{\mathcal E}}
\def\cF{{\mathcal F}}

\def\cI{{\mathcal I}}
\def\cJ{{\mathcal J}}
\def\cK{{\mathcal K}}
\def\cL{{\mathcal L}}
\def\cM{{\mathcal M}}
\def\cN{{\mathcal N}}
\def\cO{{\mathcal O}}
\def\cP{{\mathcal P}}

\def\cR{{\mathcal R}}
\def\cS{{\mathcal S}}
\def\cT{{\mathcal T}}

\def\cX{{\mathcal X}}
\def\cY{{\mathcal Y}}
\def\cZ{{\mathcal Z}}

\def\tw{\widetilde{w}}
\def\tB{\tilde{B}}

\def\tH{\widetilde{H}}

\def\tb{\widetilde{b}}

\def\cM{{\mathcal M}}
\def\cN{{\mathcal N}}

\def\na{\nabla}
\def\la{\langle}
\def\ra{\rangle}

\def\one{\mathbf{1}}
\def\udb{\underbrace}

\def\uds{\underset}

\def\Re{\mathrm{Re}}

\newcommand{\bseq}{\begin{subequations}}
\newcommand{\eseq}{\end{subequations}}

\newcommand{\beq}{\begin{equation}}
\newcommand{\eeq}{\end{equation}}
\newcommand{\bal}{\begin{aligned} }
\newcommand{\eal}{\end{aligned}}
\newcommand{\bit}{\begin{itemize} }
\newcommand{\eit}{\end{itemize}}
\newcommand{\bga}{ \begin{gathered} }
\newcommand{\ega}{ \end{gathered} }
\newcommand{\ben}{\begin{eqnarray}}
\newcommand{\een}{\end{eqnarray}}
\newcommand{\beno}{\begin{eqnarray*}}
\newcommand{\eeno}{\end{eqnarray*}}

\renewcommand{\div}{\mathrm{div}}

\newcommand{\BZ}{\mathbb{Z}}
\newcommand{\aaa}{\mathbf{a}}
\newcommand{\ee}{\mathbf{e}}

\newcommand{\FF}{\mathbf{F}}
\newcommand{\GG}{\mathbf{G}}

\newcommand{\JJ}{\mathbf{J}}

\newcommand{\R}{\mathbb{R}}
\newcommand{\C}{\mathbb{C}}

\newcommand{\UU}{\mathbf{U}} 
 
\newcommand{\WW}{\mathbf{W}}

\newcommand{\Id}{\mathrm{Id}}

\newcommand{\tf}{\tilde{f}}

\newcommand{\sgn}{\mathrm{sgn}}
\newcommand{\supp}{\mathrm{supp}}

\makeatletter
\def\th@plain{%
  \thm@notefont{}%
  \itshape %
}
\def\th@definition{%
  \thm@notefont{}%
  \normalfont %
}
\def\paragraph{\vspace{0.1in}\@startsection{paragraph}{4}%
  \z@\z@{-\fontdimen2\font}%
  {\normalfont\bfseries}}
\makeatother

\newcommand{\ang}[1]{\la #1 \ra}

\definecolor{labelkey}{rgb}{0,0,1}

\let \mr = \mathring
\let \mw = \mathrm
\let \msf = \mathsf

\let \tf = \tfrac

\def\iin{\mathrm{in}}
\def \EEs {\mathscr{E}}
\def \TTs {\mathscr{T}}
\def \bbJ {\boldsymbol J}
\def \bbY {\boldsymbol{\cY}_{\etab}}

\def \kkl {\mathsf{k_L}}
\def \kk {\mathsf{k}}
\def \cgam {\bar C_{\g}}
\def \cnon {\bar C_{\cN}}

\def \aaa { a }
\def \bbb  { b}

\def \mass {  \varrho}
\def \mome {\mathsf{ m } }
\def \eee {\mathsf{e}}

\def \Lmic {\cL_{\msf{mic},s} }

\def \FFp {F_{m, \msf{pos}}}

\def \tFM{\td F_M}
\def \tFm{\td F_m}

\def \tFMa { \td F_{M, 1} }
\def \tFMb {\td F_{M, 2}}

\def \sc{\mathsf{C}}
\def \bth{{\bar \Theta}}
\def \ths{{\bth _s}}
\def \bu{{\bar \UU}}
\def \bc{{\bar{\mathsf C}}}
\def \cs{{\bc _s}}

\def \bpi{{\boldsymbol \Pi}}

\def \brho{ \bar \rho_s}
\def \bp{\bar P_s}

\def \rhos {\bar \rho_s}
\def \Lams{\Lam}

\def \cMM{\cM _1} 
\def \bcx{\bar c_x}
\def \bcv{\bar c_v}

\def \vc{{\mathring{V}}}
\def \wc{{\mathring{W}}}

\def \dcm{d _{\cM}}

\def \tu { \td{\UU} }
\def \trho {\td \rho}
\def \tp {\td P}
\def \tb {\td B}
\def \tw {\widetilde \WW}
\def \hw {\wh \WW}

\def \eM {\cE _\cM}
\def \erho {\cE_{\rho}}
\def \ec{\cE_{\sc}}
\def \eu {\cE_{\UU}}
\def \ep{\cE_P}

\def \bF {\mathbf F}

\def \sss {\mathsf{g}}

\def \es {{\e _s}}
\def \rs {{R _s}}

\def \lame {\lam_{\xeta}}
\def \lamu {\lam_{\mathsf{u}}}
\def \lams {\lam_{\mathsf{s}}}

\def \etab {{\bar \eta}}

\def \xeta {{\underline{\eta}}}
\def \cXX {\cX_{\xeta}}
\def \cXb {\cX_{\etab}}

\def \cXC{ \cX_{\C, \xeta} }

\def \cYe {\cY_{\eta}}
\def \cYY {\cY_{\xeta}}
\def \cYb {\cY_{\etab}}
\newcommand{\cYL}[1]{\cY _{\Lam, #1}}
\def \cYE {\cYL{\eta}}

\def \cZZ { Z }
\def \cZZQ {Z_{\Lam}}
\def \cZZQR {Z_{\Lam, R}}
\def \cZZQinf { Z_{\Lam, \infty}}
\def \lvv{ {V} }

\def \cYQ{ \cYL{\etab} }
\def \tFR { \td F_{\th, R} }
\let  \rsa =\rightsquigarrow

\newcommand{\bS}{\operatorname{\boldsymbol \Sigma}}

\def \xis {\xi_*}

\def \rE { \om }

\def \kmu { \kp_2 }

\def \id { {\rm Id} }

\def \init {\mathsf{in}}
\def \Lmic {\cL _{\mathsf{mic}}}

\def \el{ \ell }
\def \ups {\upsilon}

\def \vpi {\varpi}

\def \Znonneg {\mathbb Z_{\geq 0}}

\renewcommand{\div}{\operatorname{div}}
\renewcommand{\cO}{O}

\author[J.~Bedrossian]{Jacob Bedrossian}
\address{Department of Mathematics, University of California, Los Angeles, CA 90095.}
\email{\href{jacob@math.ucla.edu}{jacob@math.ucla.edu}}

\author[J.~Chen]{Jiajie Chen}
\address{Department of Mathematics, University of Chicago, Chicago, IL 60637.}
\email{\href{jiajiechen@uchicago.edu}{jiajiechen@uchicago.edu}}

\author[M.~Gualdani]{Maria Pia Gualdani}
\address{Department of Mathematics, The University of Texas at Austin,  Austin, TX 78712.}
\email{\href{gualdani@math.utexas.edu}{gualdani@math.utexas.edu}}

\author[S.~Ji]{Sehyun Ji}
\address{Department of Mathematics, University of Chicago, Chicago, IL 60637.}
\email{\href{jise0624@uchicago.edu}{jise0624@uchicago.edu}}

\author[V.~Vicol]{Vlad Vicol}
\address{Courant Institute of Mathematical Sciences, New York University, New York, NY 10012.}
\email{\href{vicol@cims.nyu.edu}{vicol@cims.nyu.edu}}

\author[J.~Yang]{Jincheng Yang}
\address{Department of Applied Mathematics and Statistics, Johns Hopkins University, Baltimore, MD 21218.}
\email{\href{jincheng@jhu.edu}{jincheng@jhu.edu}}

\date{\today}

\title[]{Finite time singularities in the Landau equation with very hard potentials}

\thanks{\textit{Acknowledgment}. The authors would like to thank American Institute of Mathematics (AIM) for hosting the workshop ``Integro-differential equations in many-particle interacting systems'', where this project \cite{AIM2025} was initiated. AIM receives major funding from the US National Science Foundation and the Fry Foundation. JB gratefully acknowledges support from NSF DMS-2108633 and NSF DMS-2510949. JC is partially supported by NSF Grant DMS-2408098. MPG is supported by NSF DMS-2511625 and thanks the College of Natural Science at UT Austin for its support through the SPARK Grant. VV is partially supported by the Collaborative NSF grant DMS-2307681 and a Simons Investigator Award. JY is partially supported by NSF Grant DMS-1926686. The authors would like to thank Luis Silvestre for insightful discussions.
}

\begin{document}

\begin{abstract} 

We consider the inhomogeneous Landau equation with $\gamma \in (\sqrt{3},2]$ and construct smooth, strictly positive initial data that develop a finite time singularity. The $C^{\alpha}$‑norm of the distribution function blows up for every $\alpha>0$, whereas its $L^{\infty}$‑norm remains uniformly bounded. In self‑similar variables, the solution becomes asymptotically hydrodynamic—the distribution function converges to a local Maxwellian, while the hydrodynamic fields develop an asymptotically self‑similar implosion whose profile coincides with a smooth imploding profile of the compressible Euler equations.  To our knowledge, this provides the first example of a collisional kinetic model which is globally well-posed in the homogeneous setting, but admits finite time singularities for inhomogeneous data.
\end{abstract}

\maketitle
\tableofcontents

\newpage

\section{Introduction}

We consider the three-dimensional inhomogeneous Landau equation
\beq\label{eq:Landau}
\bal
& \pa_t f + v \cdot \na_x f = \f{1}{\e_0} Q(f, f) ,
\eal
\eeq
where  $f(t,x,v)\ge 0 $ is the distribution function, $x,v\in \mathbb{R}^3 $ denote the spatial and velocity variables, and  $\glossary Q$ is the collision operator defined as\footnote{For simplicity of notation, we omit the $(t, x)$-dependence of $f$ in~\eqref{eq:Q}.}
\bseq\label{eq:Q}
\beq
\bal
Q(f, g)(v) &= \div _v \int_{\R ^3} \boldsymbol \Phi (v - w) \B(f (w) \na _{v} g (v) - g (v) \na _{w} f(w) \B) d w , \\
& = \div _{v} (A [f] \na _v g - \div _v A [f] g) (v),
\eal
\eeq
where $\boldsymbol \Phi$ and $A$ are given by
\beq 
\bal
\label{eqn: defn-bPhi-A}
\boldsymbol \Phi (v) &:= \f1{8 \pi} |v| ^{\gamma + 2} (\Id - \bpi _v), \quad \quad A [f] (v) := \boldsymbol \Phi * f (v)  . %
\eal
\eeq
\eseq
Here $\bpi _v = \f v{|v|} \otimes \f v{|v|}$ is the projection along the $v$ direction whenever $v \ne 0$ and $\bpi _v = O$ if $v = 0$.
The constant $\varepsilon_0$ denotes the Knudsen number. 

The physically relevant case is $\gamma =-3$, typically called the Landau--Coulomb model. Nonetheless, a range of other values of $\gamma$ has been explored in the mathematical literature. 
In this work, we consider the inhomogeneous Landau equation for $\gamma \in (\sqrt{3}, 2]$, 
and construct smooth, strictly positive initial data that develop a finite time singularity. 
To our knowledge, this provides the first example of a collisional kinetic model that is globally well-posed for homogeneous data \cite{desvillettes2000hardpotential} yet admits finite time singularities in the inhomogeneous setting. See Section \ref{sec:intro:Landau:previous} below for more discussion on the well-posedness problem for these equations.

Our construction of the singular solutions is inspired by the hydrodynamic limit of kinetic equations to the compressible Euler equations and by the known smooth imploding solutions in the compressible Euler equations.
In \cite{merle2022implosion1,buckmaster2022smooth,shao2025blow}, the authors constructed families of smooth, self-similar imploding 
profiles for the 3D compressible (isentropic) Euler equations with 
various adiabatic exponents and blowup rates, denoted by $r>1$. 
Specifically, they construct blowup solutions on $t\in[0,1)$ %
with the following leading order structures
\begin{align*} 
    \rho (t, x) &\sim (1 - t)^{-\f{3 (r - 1)}{r }} \bar \rho \left(\f x{(1 - t) ^{1/r}} 
    \right), \\
    (\rho \mathbf U) (t, x) &\sim (1 - t) ^{-\f{4 (r - 1)}{r}} (\bar \rho \bu) \left(\f x{(1 - t) ^{1/r}} 
    \right),
\end{align*}
where $(\rho, \mathbf U)$ is the density and velocity of the compressible Euler equations and $(\bar\rho, \bu)$ denotes a smooth self-similar imploding profile for the density and velocity.

The blowup we construct for the Landau equation is \emph{asymptotically hydrodynamic} in self-similar variables. That is, 
in self-similar variables, the hydrodynamic fields converge to an imploding profile for Euler and the distribution function converges to the corresponding local Maxwellian. 
Hence, \emph{informally} the distribution looks like the following to leading order near the blowup
\begin{align}
\label{eqn: leading-order-at-blowup}
f(t,x,v) \sim \mathcal{M}_{\bar{\rho},\bu,\bar{\Theta}}\left( \frac{x}{(1-t)^{1/r}}, (1-t)^{1 - 1/r} v \right), 
\end{align}
where $\mathcal{M}_{\bar{\rho}, \bu,\bar{\Theta}}$ defined in \eqref{eq:local_max} denotes the local Maxwellian with the hydrodynamic fields given by the self-similar imploding profile, and $\bar{\Theta} = \f35 \bar \rho^{2/3}$.

\subsection{Main result} We define the mass $\varrho$, momentum $\mome$, and energy density $\eee$ associated to $f$ by
\beq\label{eq:density}
 (\mass, \mome, \eee )(t,x)  :=	\int f(t, x ,v)( 1 , v , |v|^2 ) d v .
\eeq
Our main result is the following. 
\begin{theorem}
\label{thm:blowup}
Fix  $\g \in (\sqrt{3}, 2]$  in  \eqref{eq:Landau}. Let $( \bar \rho, \bar \UU, \bar \Th, r )$ %
denote the imploding profile for the 3 dimensional  compressible Euler equations constructed in \cite{shao2025blow}, with a blowup speed $r  > (\gamma+3)/(\gamma+2)$.

 There exists a small $\e_*>0$, such that for any $0< \e_0 \leq \e_*$, the following statement holds: 
 there exists a positive initial data $f_{\iin} \in C^{\infty}$ with Gaussian decay $f_{\iin}(x, v)  
 \leq c \exp(- C |v|^2)$ and  uniformly bounded momentum $\mome_{\iin}(x)$ and energy $ \eee_{\iin}(x)$ 
 that decay to $0$ as $|x| \to \infty$, and with a  mass density $ 0< c_1 \leq \mass_{\iin}(x) \leq c_2$ for some constants $c_1, c_2>0$
\footnote{
The initial mass density $\varrho_{\iin}$ may be chosen to equal a constant outside a compact set in $\R^3$. See Section \ref{sec:non_stab}. 
}
, such that the corresponding positive solution $f$ to \eqref{eq:Landau} develops a finite time singularity at a time $T=1$ in the following sense: 
\vspace{0.2cm}

\begin{enumerate}[\upshape (a)]
\item 
Regularity: as $t \to 1^-$, the $C ^{\al}$-norm of $f(t, \cdot, v = 0)$ blows up for any $\al > 0$, 
while  %
$\| f(t) \|_{L^{\infty}_{x,v}}$ remains uniformly bounded. Moreover, for any $v$, the spatial gradient $\na_x f $ blows up at $(0, v)$ in the following sense: 
$ \sup_{ |x| \leq (1-t)^{1/r}} |\na_x f(t, x, v)| \to \infty$ as  $ t \to 1^-$. 
Furthermore, the solution remains smooth away from $x=0$: for any $x \neq 0$, $v \in \R^3$, and multi-indices $\al, \b$ with $|\al| + |\b| \leq 16$,\footnote{
Given any $\msf{n} \geq 0$, we can construct a blowup solution $f$ that satisfies all the properties in Theorem \ref{thm:blowup}
and estimate \eqref{eq:smooth_away_0:1} with $|\al| + |\b| \leq \msf{n}$. 
One only needs to modify $\kk$ in \eqref{def:kk} with $\kk \geq 2 d + \msf{n}$ in the proof. In \eqref{eq:smooth_away_0:1}, we fix $\msf{n}=16$. 
} we have
\beq\label{eq:smooth_away_0:1}
\sup_{t \in [0, 1)}|\pa_x^\al \pa_v^\b f(t, x, v )| \leq C(\e_0, x, v) < \infty .
\eeq

\item Hydrodynamic limit in self-similar variables: in self-similar variables, %
$f$ converges to the local Maxwellian associated with the Euler profile :
\beq\label{eq:blow_asym:micro}
		\lim_{t \to 1^-}  f(t, (1-t)^{ \f{1}{r} } X, (1-t)^{ -\f{r-1}{r} } V)
	= \cM_{\bar \rho, \bar \UU, \bar \Th}(X,V), 
\eeq
for any fixed $X, V \in \R^3$, where $\cM_{\bar \rho, \bar \UU, \bar \Th}(X, V)
=\bar \rho( X) \f{1}{(2 \pi \bar \Th( X))^{3/2}} \exp\left(-\f{|V - \bar \UU( X)|^2}{2 \bar \Th( X)} \right)$.

\item Implosion in hydrodynamic fields: as $t \to 1^-$, the hydrodynamic fields $(\varrho, \mome, \eee)$ \eqref{eq:density} all blow up 
 at $x =0$. The blowup is asymptotically self-similar in the sense that
 \begin{subequations}\label{eq:blow_asym:macro}
\begin{align} 
 \lim_{t \to 1^-}   (1-t)^{ \f{3(r-1)}{r} } \mass ( t, (1-t)^{ \f{1}{r} } X) & =  \bar \rho (X),   \\
 \lim_{t \to 1^-}   (1-t)^{  \f{4(r-1)}{r} } \mome( t, (1-t)^{ \f{1}{r} } X)   & = (\bar \rho \bu) (X),  \\
 \lim_{t \to 1^-}   (1-t)^{  \f{5(r-1)}{r} } \eee ( t, (1-t)^{ \f{1}{r}  } X)  & = (  3 \bar \rho  \bar \Th    +  \bar \rho  |\bu|^2 )(X) ,
\end{align}
 \end{subequations}
for any fixed $X \in \R^3$, where $r > 1$.
\end{enumerate}

\end{theorem}

We establish Theorem \ref{thm:blowup} by exploiting the connection between the Landau equation and the compressible Euler equations. 
Instead of performing a Hilbert expansion similar to 
\cite{caflisch1980fluid,guo2010acoustic}, we develop a framework to establish nonlinear (finite codimension) stability of the local Maxwellian in self-similar variables and justify the 
connection between the two equations up to the blowup time $T=1$ for a small, \textit{fixed} $\e_0$. See Section \ref{sec:idea} for  more discussions.
To simplify notation, we may use the scaling symmetries of the Landau equation to fix the initial time at $t=0$ and the blowup time at $ t = T = 1$.
Below, we list a few remarks on the results in Theorem \ref{thm:blowup}.

\begin{remark}[Range of $\g$]
\label{rem:gamma:lower:bound}

In this work, we provide a \textit{proof of concept}\footnote{Our contribution may be phrased as follows: if you only use certain properties of the collision kernel $Q$ appearing in~\eqref{eq:Landau}, then finite time blowup cannot be ruled out. In some sense, this is akin to Tao's result for ``averaged Navier--Stokes''~\cite{tao2016finite}; while there is no hydrodynamic meaning to the model in~\cite{tao2016finite}, Tao's paper shows that if you only use certain properties of the bilinear nonlinear term, then you cannot rule out finite time blowup.} in kinetic equations, showing how one can ``lift'' compressible Euler singularities to the Landau equation via the hydrodynamic limit in self-similar variables. This confirms a scenario vaguely alluded to in~\cite[Section 8.1]{Silvestre2023}, albeit only for $\gamma > \sqrt{3}$. We expect that the admissible range of $\gamma$ may be extended if a broader class of implosion profiles for compressible Euler equations is shown to exist. See further discussion at the end of Section~\ref{sec:Euler_review}.

The condition $\gamma > \sqrt{3}$ arises solely from the \textit{existing}  class of 
smooth imploding profiles for the 3D compressible Euler equations  with monatomic gas used to lift singularities to the kinetic level. 
Since smooth imploding profiles in this class are known to exist for $r < 3 - \sqrt{3}$, this restriction combined with the condition $r > (\gamma+3)/(\gamma+2)$ yields $\gamma > \sqrt{3}$.

\end{remark}

\begin{remark}[Set of initial data]
\label{rem:initial:data:set}
To simplify the analysis, we consider solution with radial symmetry: $f(t,Q x , Q v) = f(t,x, v)$ for any orthogonal matrix $Q \in SO(3)$, which is preserved by equation \eqref{eq:Landau}. %
The associated hydrodynamic fields $(\mass, \mome, \eee)$ \eqref{eq:density} are radially symmetric in $x$.

Within the radially symmetric class, as to be explained in Remark~\ref{rem:data}, the initial data can be decomposed into $F_{\iin} = \cM + \cMM^{1/2} ( \cF_M(\tw) + \tFm)$, where $\cM$ is the modified local Maxwellian defined in \eqref{eq:localmax2}, $\cF_M(\tw)$ is a perturbation associated with the hydrodynamic fields $\tw$ (similar to $\mass, \mome, \eee$ \eqref{eq:density}), and $\tFm$ is the micro-perturbation. There exists an open set $X_2$ (a ball in a weighted Sobolev space) and a finite codimension set $X_1$ such that the positive initial data in Theorem~\ref{thm:blowup} may be taken such that $ \tw  \in X_1$ and $ \tFm \in X_2$. We can localize the unstable directions of the blowup profile to ensure that the initial mass density $\mass_{\iin}$ admits a uniform positive lower bound. In addition, we can construct a blowup solution whose initial hydrodynamic fields $(\mass, \mome, \eee)(x)$
decay to $0$ as $|x| \to \infty$.

Extending the blowup construction to non-radial perturbations of $\cM$ would follow the framework developed here, combined with the global weighted $H^k$ stability estimates of implosion with non-radial perturbations made in \cite{chen2024vorticity}. 
\end{remark}

\begin{remark}[Tail fattening at $x = 0$]\label{rem:tail_flat}

Let $ \bc $ be the sound speed profile  defined in \eqref{eq:Euler_profi}
and $\mu(\cdot)$ be the Gaussian defined in \eqref{eq:gauss}. 
Recall that the initial data satisfies $f(0,x,v) \les \exp(-C|v|^2)$. 
In the proof of Theorem \ref{thm:blowup} in Section \ref{sec:non_stab}, %
at $x = 0$, for any $ t \in [0, 1)$ and $v \in\R^3$, we establish
\[
  \left|f(t, 0, v) - \mu \left( \f{v}{ \bc(0) \cdot (1-t)^{-\f{r - 1}{r}}} \right) \right| \les \e_0^{\ell} \mu^{ \f12 } \left( \f{v}{ \bc(0) \cdot (1-t)^{-\f{r - 1}{r}}}\right), \qquad r > 1, \    \ell = 10^{-4}.
\]
In particular, by choosing $\e_0$ small, $f(t, 0, \cdot)$ is close to a constant as $t \to 1^-$:
\beq\label{eq:solu_lim:0}
  \limsup_{ t \to 1^-} \bigl|f(t, 0, v) -\mu(0) \bigr| \les \e_0 ^{\ell} \mu (0) ^\f12, \qquad \forall v \in \R ^3. 
\eeq
Here $\e _0 ^\ell \mu (0) ^\f12 \ll \mu (0)$ when $\e _0$ is small.
\end{remark}

\begin{remark}[Smoothness away from $x=0$]\label{rem:smooth_away_0}

We establish a more quantitative version of \eqref{eq:smooth_away_0:1} in the proof of Theorem \ref{thm:blowup} in Section \ref{sec:non_stab}.
There exists a large parameter $R_0= R_0(\e_0)$, a function $ \msf{c}_{R_0}(x) \asymp \min \{|x|, R_0 \} ^{-(r - 1)}$ (defined in \eqref{def:cRx}) and an absolute constant $ C_{\bu}$ such that, for 
\beq\label{eq:def_small_vring}
\mr{v} := \f{ v - C_{\bu} |x|^{-r} x}{  \msf{c}_{R_0} (x) }
\eeq
(see also \eqref{def:small_vring}), and for any multi-indices $\al, \b$ with $|\al| + |\b| \leq 16 $, 
and any $x\neq 0, v \in \R ^3$, we have %
\[
\sup_{t \in [0,1 )} |\pa_x^{\al} \pa_v^{\b} f(t, x, v)| \les_{ \e_0} |x|^{-|\al| } \exp(- C |\mr{v}|^2), 
\]
where $C$  is some absolute constant independent of $\al, \b, \e_0$. 

\end{remark}

\begin{remark}[Limiting solution]\label{rem:limit_solution}
Let $\mu( \cdot)$ be the Gaussian defined in \eqref{eq:gauss}. For any fixed $x \neq 0, v \in \R^3$, in the proof of Theorem \ref{thm:blowup} in Section \ref{sec:non_stab}, we establish that the blowup solution $f$ is close to $\mu(\mr{v})$ (see Figures~\ref{fig:1} and~\ref{fig:2} below) in the following sense:
\beq\label{eq:limit_solu}
\bal 
\limsup_{t \to 1^-} |   f(t, x, v) - \mu( \mr{v} )| & \les  \e_0^{\ell} 
  \mu( \mr{v} )^{\f12} .
\eal 
\eeq
For $\mr{v}$ defined in \eqref{eq:def_small_vring}, we can obtain $|\mr{v}| \les  |v| \cdot |x|^{r-1} + 1$ with $r >1$. If $  |v| \cdot |x|^{r-1} \leq c (\log \f{1}{\e_0} )^{1/2}$ for some small $c$ and $\e_0$ is small, the error term is smaller: $\e_0^{\ell} \mu( \mr{v})^{1/2} \les \e_0^{\ell/2} \mu( \mr{v}) \ll  \mu( \mr{v})$. 
Thus, along the surface $\{  (x, v): |\mr{v}|^2 = C \}$ with $C \ll \log \e_0^{-1}$, $f$ is close to a constant.  Using the formula of $\msf{c}_{R_0}(x)$ \eqref{def:cRx}, when $|x| \leq R_0$, we compute $\mr{v}, |\mr{v}|^2$ 
\[
 \mr{v} = c_1 v |x|^{r-1} -c_2 \f{x}{|x|},
 \quad |\mr{v}|^2 = c_1^2 |v|^2 |x|^{2r-2} - 2 c_1 c_2 ( v \cdot x) |x|^{r-2} + c_2^2, 
\]
where $c_1 = C_{\bc}^{-1}, c_2 = \f{C_{\bu}}{C_{\bc}}$ are the constants associated with profile. See \eqref{def:cRx} and 
\eqref{eq:UC_refine_asym}.
Note that when $x = 0$, as $t\to 1^-$, $f(t, 0, v)$ is close to $\mu(0)$ \eqref{eq:solu_lim:0}. 
For a fixed $v$, since $|\mr{v}|^2 \to c_2^2$ as $x\to 0 $ and $c_2$ may not be $0$, the limiting function of $f$ may not be continuous at $x = 0$.

\begin{figure}[htbp]
  \centering
  \begin{minipage}[b]{0.48\textwidth}
    \centering
\includegraphics[width=\textwidth]{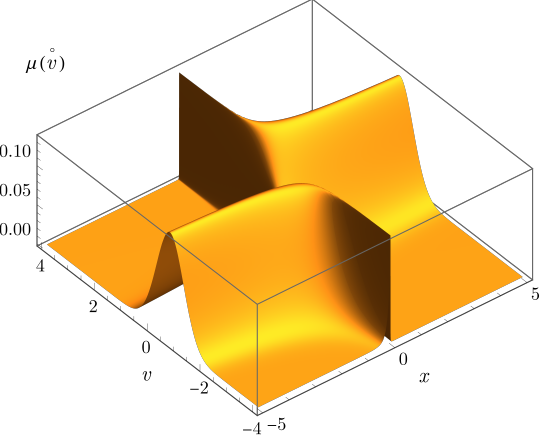}
  \end{minipage}
  \hfill
  \begin{minipage}[b]{0.48\textwidth}
    \centering   \includegraphics[width=\textwidth]{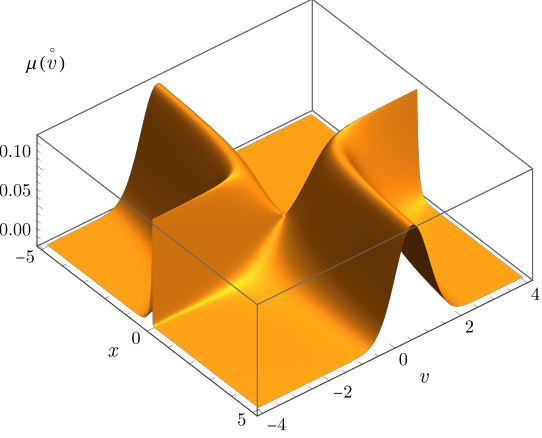}
  \end{minipage}
  \caption{A cartoon of the limiting density $\mu(\mr{v})$ with $x \neq 0 $ when both $x$ and $v$ are one-dimensional, with $\mr v$ defined in \eqref{eq:def_small_vring}. We have taken $r=1.26<3-\sqrt{3}$ and set $C _\bc = C _\bu = 1$, $R _0 \to \infty$ in~\eqref{eq:def_small_vring} and \eqref{def:cRx}. The two images represent the same function $\mu(\mr{v})$ from two different perspectives (rotated by $90^\circ$ in the horizontal plane).
  }
  \label{fig:1}
\end{figure}

\begin{figure}[htbp]
  \centering \includegraphics[width=0.75\textwidth]{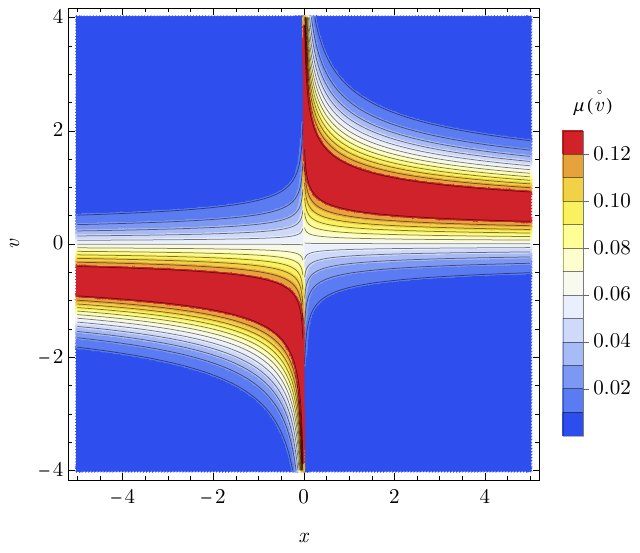} 
  \caption{A contour plot (top-down view) of the same function $\mu(\mr{v})$ from Figure~\ref{fig:1}. The red regions represent higher values of $\mu(\mr{v})$, while the blue regions represent values close to zero.}
  \label{fig:2}
\end{figure}

\end{remark}

\begin{remark}[$L^{\infty}$ bound $\not \Rightarrow$ smoothness]
The singularity formation established here concerns the macroscopic quantities---namely mass, momentum, and energy. At the blowup time, the distribution function $f$ remains bounded in 
$L^\infty$, whereas the $C^\alpha$ norm blows up. A De Giorgi second lemma (boundedness $\Rightarrow$ H\"older regularity) fails, partly because boundedness of $f$ does not guarantee boundedness of the coefficients $A[f]$ and $\div A[f]$ in \eqref{eq:Q}. Unlike the homogeneous case, in the inhomogeneous setting, the bounds for $A[f]$ and its derivative require control on the hydrodynamic fields, which is lost at the blowup time. A recent work by Golding and Henderson~\cite{GoldingHenderson25} suggests that for the Coulomb potential $\gamma = -3$ the condition $f\in L^1_t L^\infty_{x,v}$ is, however, enough to guarantee global regularity, which would rule out a finite time singularity of this type in that case.

\end{remark}

\begin{remark}[Local well-posedness]
A general local well-posedness theory for the Landau equation \eqref{eq:Landau} with $\gamma \in ( \sqrt{3}, 2]$ is not covered in this manuscript. Our work provides, however, local well-posedness for \eqref{eq:Landau} with any $\e _0 > 0$ and $\g \in (\sqrt 3, 2]$ for initial data near the local Maxwellian associated with our profile. See Proposition \ref{prop:LWP_landau}. The argument used to prove Proposition \ref{prop:LWP_landau} is inspired by the 
method used in \cite{henderson2019local}.
We also prove that a Gaussian lower bound for $f$ propagates in time via a barrier argument.
\end{remark}

\subsection{Related work}\label{sec:review}

We review related work on well-posedness results for kinetic equations, singularities in the compressible Euler equations, and hydrodynamic limits.

\subsubsection{Well-posedness results for the Landau and Boltzmann equations} 
\label{sec:intro:Landau:previous} 

    The Landau equation is one of the most important mathematical models in collisional plasmas. It was derived from the Boltzmann equation by Landau \cite{Landau1936Coulomb} in 1936 to model the grazing collisions that dominate in the charged particle collisions present in plasma physics (in the case $\gamma = -3$). In this regime, the classical Boltzmann collisional operator becomes singular and reduces to \eqref{eq:Q}.

    The mathematical investigation of these equations began with the work of DiPerna and Lions in \cite{diperna1989cauchy} for the inhomogeneous Boltzmann equation. The case of the inhomogeneous Landau equation was  discussed in \cite{lions1994CauchyLandau} by Lions later. In these works, the authors provided the first complete Cauchy theory for very weak solutions, the so-called renormalized solutions. Further progress was made in the late 1990s by Villani \cite{villani1998hsolution}, who introduced the notion of (very) weak solutions for the spatially homogeneous case which satisfy the Boltzmann $H$-theorem and relative bounds on the entropy $\int f \log f \;dv$. These are hence referred to as $H$-solutions. 
    
    Since then, significant effort has been devoted to the understanding of the {\em{spatially homogeneous}} Landau equation. In 1998, Villani \cite{villani1998MM} proved the global well-posedness of the space-homogeneous Landau equation for Maxwellian molecules $(\gamma=0)$, and shortly after, with Desvillettes, he extended the global well-posedness to the case of hard potentials ($\gamma>0$) in \cite{desvillettes2000hardpotential} and \cite{desvillettes2000hardpotential2}.
    The global existence of smooth solutions of the spatially homogeneous Landau equation for the case of moderately soft potentials ($\gamma \in [-2,0)$) was established later by Wu \cite{wu2014softpotential} and Silvestre \cite{silvestre2017upperboundlandau}. While our understanding for the cases of moderately soft potentials and hard potentials was satisfactory, the cases of the Coulomb potential ($\gamma=-3)$ and very soft potentials $\gamma \in (-3,-2)$ remained elusive. The main question was to understand whether for $\gamma <-2$ the drift term would be controlled by the diffusion, implying existence of smooth solutions for arbitrarily large times. 

 Partial progress was made on the regularity of the solutions, analogous to the results for the Navier--Stokes equation of the partial regularity theory of Caffarelli--Kohn--Nirenberg \cite{caffarelli1982partial}. Golse, Gualdani, Imbert, and Vasseur \cite{golse2022timepartial} proved that the Hausdorff dimension of the singular set in time is at most $1/2$ for the space-homogeneous Landau equation with the Coulomb potential. This result was later extended by Golse, Imbert, Ji, and Vasseur \cite{golse2024velocitytimepartialregularity} to partial regularity in velocity-time space. 
    In parallel, the conditional regularity results for the space-homogeneous case were obtained by Silvestre \cite{silvestre2017upperboundlandau} and Alonso, Bagland, Desvillettes, and Lods \cite{alonso2024prodiserrin}, in analogy with the Prodi--Serrin criteria for the Navier--Stokes equation. Other types of conditional regularity results were also established by Gualdani and Guillen \cite{gualdani2019Apweight}. Recently, Guillen and Silvestre \cite{GuillenSilvestreLandau} made a breakthrough, proving the global well-posedness for initial data with Maxwellian tails, based on the discovery of a new monotone quantity, the Fisher information. Subsequently, the global existence of smooth solutions was extended to broader classes of initial data by several authors, and currently it is known for weighted $L^1$ initial data (see \cite{golding2025smoothsolution,ji2024dissipationestimatesfisherinformation,desvillettes2024productionfisherinformationlandaucoulomb,he2025existenceuniquenesssmoothingestimates}).

For the space-inhomogeneous Landau equation, the well-posedness theory was initiated by Guo \cite{guo2002landau}, who constructed a unique global solution for $\gamma\ge -3$ for initial data close to the global Maxwellian in a high-order Sobolev space with a fast decaying velocity tail, namely $H^8(\mathbb{T}_x^3 \times \R^3_v ; \mu(v)^{-1/2} d v )$. 
For $\gamma \in[-3,-2)$, this result was improved in terms of the decay rate of the velocity tail in \cite{Strain2006almost}, \cite{Strain2008exponential}, which also include the study of various types of Boltzmann and Landau equations in the perturbative regime. The key ingredient underlying these stability estimates near the global Maxwellian is coercivity estimates for the linearized Landau collision operator \cite{guo2002landau}. 
\footnote{
The proof of Guo in \cite{guo2002landau} relies on a compactness argument and therefore yields
non-explicit constants in the coercivity estimates. Constructive proofs of coercivity
estimates for the Boltzmann and Landau collision operators, with explicit and computable
constants, are established by Baranger--Mouhot in \cite{baranger2005spectralgap},
by Mouhot \cite{Mouhot2006coercivity}, and by Mouhot--Strain \cite{mouhot2007spectral}.
}
Carrapatoso and Mischler \cite{carrapatoso2017nearMaxwellian} substantially improved these results for initial data close to the global Maxwellian in $H^2_x L^2_v(\mathbb{T}_x^3 \times \R^3_v ; (1+|v|^2)^{k/2} dv)$ for $\gamma \in [-3,-2).$ It is worth noting that these works consider the space-inhomogeneous case, but the main difficulty arises from analyzing the (weak) dissipation generated by the collision operator, which only acts on the velocity variable $v$.
We also remark that there is a development of the global well-posedness theory for solutions near vacuum, which can be understood as another type of singularity, by Luk \cite{luk2019vacuumsoft} for $\gamma \in (-2, 0)$ and by Chaturvedi \cite{chaturvedi2022vacuumhard} for $\gamma \in [0,1)$. See also \cite{guo2001vlasovVacuum,guo2002vlasovMaxwellian,guo20212vlasovLandau,Chaturvedi2023damping,Yan2003VlasovMaxwellBoltzmann} and the references therein for related works regarding the Vlasov--Poisson--Landau, Vlasov--Poisson--Boltzmann systems, and Vlasov--Maxwell--Boltzmann, where the kinetic equation is coupled with a self-consistent field.

For large initial data, the conditional regularity results were first established for the inhomogeneous Boltzmann equation with moderately soft potentials by Imbert and Silvestre in \cite{imbert2020conditionalregualarityBoltzmann} (and later improved in \cite{imbert2022Boltzmannconditional}). Moving to the inhomogeneous Landau equation, the first conditional regularity was obtained by Cameron, Snelson and Silvestre \cite{cameron2018landauconditional} for the case of moderately soft potentials (See also \cite{henderson2020smoothening}). A similar conditional regularity result for hard potentials was later established by Snelson in \cite{snelson2020hardpotential}.
 These results are conditioned on the $L^\infty$ boundedness of the hydrodynamic densities (mass, energy, and entropy). From this, one can expect that one possible type of singularity formation for the space-inhomogeneous Landau equation is the implosion of hydrodynamic quantities. 
 In a series of papers \cite{henderson2019local}, \cite{henderson2020local}, Henderson, Tarfulea, and Snelson relaxed the assumptions on the hydrodynamic densities in conditional regularity results for moderately soft potentials. Moreover, for the Coulomb and very soft potentials, they obtained a continuation criterion under an additional assumption on the $L^{\infty}_{t,x} L_v^p$ norm of $f$ for $p>\frac{3}{3+\gamma}$ (with $p=\infty$ when $\gamma=-3$). In more recent work, Snelson and Solomon 
\cite{snelson2023continuation} improved this criterion, by lowering the required exponent to $p>\frac{3}{5+\gamma}$.

Very interestingly, Golding and Henderson \cite{GoldingHenderson25} recently established a new continuation criterion for the Landau equation with a Coulomb potential 
($\gamma = -3$) based on a fundamentally different quantity, the $L^1_t L^\infty_{x,v}$ norm of $f$, which does not rely on the hydrodynamic quantities. %
For $\gamma = -3$, this criterion rules out the same type of singularities constructed in our work, where $\| f \|_{L^{\infty}_{x,v}}$ remains uniformly bounded.
However, \cite{GoldingHenderson25} does not rule out a Type I blowup with rate $c_f=-1$ (see \eqref{eq:hydro}) in 
the case $\gamma=-3$, and it would be interesting to explore a Type I blowup with a profile $F = \cM + g$ with slow decay in $v$ \cite{bedrossian2022non}, where $\cM$ is a local Maxwellian and $g$ is a small perturbation that \textit{does not} converge to zero asymptotically, 
\footnote{Such a scenario is exploited in \cite{chen2023nearly} to construct nearly self-similar
blowup for a model problem of the homogeneous Landau equation with Coulomb potential, in
which an additional nonlinear term $a f^2$ with $0<a\ll1$ is included.
}
corresponding to the degenerate case of \eqref{eq:hydro}.

In the literature on the Boltzmann and Landau equations, the parameter $\gamma$ is often restricted to the range $[-3,1]$, corresponding to the power-law interaction potentials from the Coulomb potential for the Landau equation ($\gamma=-3$) to the hard sphere model for the Boltzmann equation ($\gamma=1$).
However, this upper bound $\gamma=1$ is not a fundamental obstruction and can be extended up to $2$ in many settings. For example, Desvillettes and Villani note in the introduction of \cite[Section 1]{desvillettes2000hardpotential} that their restriction to $\gamma \in (0,1]$ is made for simplicity, and it can be easily relaxed to $\gamma \in (0,2)$.
A similar situation arises in \cite[Theorem 1.2]{GuillenSilvestreLandau}, where it is stated only for the range $\gamma \in[-3,1]$, even though the monotonicity of the Fisher information holds for a substantially larger range that covers $\gamma \in [-3,2]$. More precisely, the restriction $\gamma \le 1$ was only used in \cite[Theorem 2.4]{GuillenSilvestreLandau}, which concerns the propagation of moments. 
For $\gamma>0$, however, 
as noted in the paragraph following  \cite[Theorem 2.4]{GuillenSilvestreLandau}, one in fact expects the generation of moments, which is even stronger than the propagation of moments. It was explicitly mentioned in remarks following \cite[Lemma 2]{desvillettes2000hardpotential} that their argument works for $\gamma \in (0,2]$.
We also remark that, in the recent lecture notes by Villani \cite[Section 4]{villani2025fisher}, he argues that the natural physical restriction for $\gamma$ is $\gamma\leq 2$, and provides a figure covering this range.

There have been some efforts to study potential singularity formation in collisional kinetic equations via model problems and self-similar methods. In \cite{andreasson2004blowup}, the authors established finite time blowup for a model of the Boltzmann equation without the loss term in the collision operator. In \cite{chen2023nearly}, %
Chen studied the homogeneous Landau equation with $ \f{1}{\e_0} Q(f,f)$ replaced by a modified collision operator $\frac{1}{\e_0} Q(f,f) + \div_v( \div_v A[f]) f$ \eqref{eq:Landau}. For \emph{any} sufficiently small $\e_0>0$, %
by perturbing the global Maxwellian and adapting a perturbation-of-equilibrium idea from \cite{chen2020slightly}, Chen established nearly self-similar blowup for the model equation with very soft potentials. 
In \cite{bedrossian2022non}, the authors ruled out Type I self-similar blowups for the Landau equation \eqref{eq:Landau} for any $\gamma \in [-3,-2]$ over a range of possible blowup speeds, assuming at least integrability in velocity on the inner profile (along with certain other structural assumptions).

\subsubsection{Finite time singularities for the compressible Euler equations}\label{sec:Euler_review}

The compressible Euler equations are the fundamental macroscopic model governing the motion of inviscid fluids and gases. A central problem in the analysis of these equations is the formation of finite time singularities from smooth initial data. While the question of global regularity versus finite time blowup has been extensively studied, the precise \emph{nature} of the singularity and the precise mechanism by which smooth solutions break down, has only been understood in certain regimes, which broadly speaking fall into two fundamentally different categories.

The classical mechanism is \textit{shock formation}; in this setting, as the smooth solution evolves, the density, pressure, and velocity all remain bounded, but their gradients blow up. In  {one space dimension}, this phenomenon is pretty well understood, ranging from the pioneering work of Riemann~\cite{Riemann1860}, up to the modern theory of hyperbolic systems of conservation laws, which provides a framework for shock formation, propagation, and interaction (see e.g.~the book of Dafermos~\cite{dafermos2005hyberbolic}). The  {multi-dimensional problem}, however, presents formidable difficulties due to the geometry of the steepening wavefronts. The theory eventually culminated in the recent work of Shkoller and Vicol~\cite{ShVi2023}, who have addressed the problem of maximal globally hyperbolic development of smooth and non-vacuous initial data, providing a complete description of how smooth solutions form  their first gradient singularity, and then continue beyond, through a succession of gradient catastrophes. For a detailed account of the literature on multi-D shocks for the Euler equations, we refer the interested reader to the summary in~\cite{ShVi2023}. Despite this progress, \emph{shock development in multiple space dimensions}, remains a fundamental open problem.
 
A qualitatively different blowup mechanism is provided by \textit{implosion singularities}. Therein, (some of) the primary flow variables  themselves become unbounded at the singularity. Such solutions describe the self-focusing collapse of a fluid/gas onto a point, and represent strong amplitude singularities, rather than  gradient catastrophes. The remainder of this section  focuses on this latter class of singularities, which has seen remarkable progress in recent years.  

The study of imploding solutions to the compressible Euler equations dates to the seminal 1942 work of Guderley~\cite{Gu1942}, who constructed self-similar solutions describing a converging spherical shock wave into a quiescent medium, which collapses at the origin. 
At the moment of collapse $t=0$, the pressure and velocity diverge at the origin, but the density does not.

A  rigorous mathematical construction of Guderley's converging shock was obtained recently by Jang, Liu, and Schrecker~\cite{JangLiuSchrecker2025}. In a different direction, Cialdea, Shkoller, and Vicol~\cite{cialdea2025Guderley} proved that Guderley's imploding shock solution arises dynamically from classical, shock-free initial data.  

We wish to emphasize that the Guderley solution suffers from a major drawback, which prevents us from using it in the construction of this paper: in the quiescent core, namely at radii $< (1 - t)^{1/r}$ for a specific similarity exponent $r$, the \textit{pressure} and the \textit{speed of sound} are assumed to \textit{vanish identically}, which is inconsistent with the regularizing effects/positivity gaining effects of the Landau collision operator (see e.g.~\cite{henderson2019local}).
 
A fundamentally new type of implosion singularity was discovered by Merle, Rapha\"el, Rodnianski, and Szeftel in~\cite{merle2022implosion1}. Unlike the Guderley solution, which contains a shock discontinuity, the implosions in~\cite{merle2022implosion1} remain $C^\infty$ \emph{smooth}
\footnote{Jenssen and collaborators~\cite{Jenssen2025, JenssenTsikkou2018} have studied a different class of amplitude blowup solutions without shock discontinuities, which are \emph{continuous but not smooth}. %
} 
until the singular time, at which point both density and velocity blow up at a single point; these solutions are also  radially symmetric, but as opposed to Guderley, they are isentropic. 
The analysis in~\cite{merle2022implosion1} requires a delicate understanding of the ODE phase portrait governing self-similar profiles, showing that for special values of $r$, 
the solution curve in the phase portrait 
passes smoothly through the so-called ``sonic point''. The price for this smooth transition through the sonic point is that the resulting profiles are only stable up to a \textit{finite-dimension instability}. It is essential to recognize how \emph{unstable} these solutions are; while the unstable manifold of these self-similar profiles is finite-dimensional, the precise number of unstable directions has not been rigorously determined, and the detailed structure of the unstable manifold remains unknown.\footnote{Numerical investigations by Biasi~\cite{biasi2021self} suggest that certain unstable directions lead to shock formation before the implosion can occur, demonstrating that generic perturbations destroy the smooth implosion mechanism entirely.} This extreme instability means that it is essentially impossible to observe these smooth implosions in physical experiments or direct numerical simulations of compressible Euler; more pertinent to the present work, these instabilities account for severe technical difficulties in implementing a finite-codimension stability argument for the macroscopic part of the Landau distribution function.

The construction in~\cite{merle2022implosion1} holds for a generic set of adiabatic exponents (countable complement), excluding the physically important monatomic gas exponent $5/3$, which corresponds to %
a degenerate case in the structure of the phase portrait of~\cite{merle2022implosion1}; the exponent $5/3$ is the one relevant for the analysis in this paper. Recently, Buckmaster, Cao-Labora, and G\'omez-Serrano~\cite{buckmaster2022smooth} have extended the construction of smooth isentropic implosion profiles to cover \emph{all} adiabatic exponents $> 1$, and established the 
profile properties for stability analysis in the diatomic case with an adiabatic exponent $7/5$. Very recently, Shao, Wang, Wei, and Zhang~\cite{shao2025blow} have specifically addressed the degenerate monatomic case with an adiabatic exponent $5/3$, establishing the existence of smooth isentropic self-similar implosion profiles for %
a sequence of blowup speeds $r_n$. These values approach (from the left) the limiting value $r_* = 3- \sqrt{3}$
and correspond to the largest possible self-similar exponent $r$ for smooth isentropic implosion with radially symmetric profiles.
See Section~\ref{sec:isen_euler} below. The limitation $r< r_* = 3-\sqrt{3}$ is the only cause for the limitation~$\gamma> \sqrt{3}$ in Theorem~\ref{thm:blowup}; see also Remark~\ref{rem:gamma:lower:bound}. 
Moreover, the authors~\cite{shao2025blow} established the profile properties for stability analysis (see Lemma \ref{lem:profile0}), which we will use in our stability analysis.

The constructions in~\cite{merle2022implosion1,buckmaster2022smooth,shao2025blow} are inherently radially symmetric, raising the question of whether smooth implosions exist without radial symmetry. 
In ~\cite{cao2024non}, Cao-Labora, G\'omez-Serrano, Shi, and Staffilani proved that the existing radial implosion profiles are stable for non-radial perturbation in a finite codimension set. 
In a different direction Chen, Cialdea, Shkoller, and Vicol~\cite{chen2024Euler} proved that a radially symmetric (hence irrotational) implosion may be ``lifted'' as an axisymmetric imploding solution to the 2D compressible Euler equations, which exhibits \emph{vorticity blowup} in finite time. 
A distinctive feature is that the swirl velocity enjoys full stability, rather than the finite-codimension stability. This result was subsequently generalized to dimensions $d \geq 3$ by Chen \cite{chen2024vorticity}. At the technical level, the analytical framework put forth in~\cite{chen2024Euler} and~\cite{chen2024vorticity} 
establishes the global weighted $H^k$ stability estimates of implosion based on the primary flow variables and plays a key role for the analysis in the current paper. See Section~\ref{sec:idea} below.

We emphasize that the smooth implosion mechanism discovered for the compressible Euler system has been shown to have profound implications for singularity formation in related equations. The works~\cite{merle2022implosion1,buckmaster2022smooth,shao2025blow} establish %
implosion for  compressible Navier--Stokes by showing that the Euler self-similar profile dominates the dynamics in appropriate parameter regimes, with viscosity treated as a perturbation. Via the Madelung transform, the defocusing NLS maps to a system resembling compressible Euler with quantum pressure. Merle, Rapha\"el, Rodnianski, and Szeftel~\cite{merle2022implosion3} used their Euler implosion profiles to prove finite time blowup for the energy-supercritical defocusing NLS, resolving a longstanding open problem. Non-radial extensions were obtained in \cite{CGSS_NLS2024}.  Shao, Wei, and Zhang~\cite{shao2024self} constructed self-similar imploding solutions for the relativistic Euler equations and used them to prove blowup for the supercritical defocusing nonlinear wave equation~\cite{shao2025blowup} with 
complex-valued solution. Subsequent work by Buckmaster and Chen~\cite{buckmaster2024blowup} 
established the blowup result 
in dimension $d=4$ and for the nonlinearity $p=7$, which corresponds to the end-point case of the blowup mechanism for the wave equation with a radially symmetric, complex-valued solution.

We note that all known \emph{smooth implosions}~\cite{merle2022implosion1,buckmaster2022smooth,shao2025blow} share three fundamental limitations: 
the blowup profiles are \emph{radially symmetric}, they apply only to the \emph{isentropic} Euler equations, and they are \emph{unstable} in the sense that only a non-quantitative, finite codimensional manifold of initial data leads to implosion. It remains a formidable open challenge to discover (1) smooth imploding solutions with non-radial profiles, and (2) smooth imploding solutions for the \emph{full} (non-isentropic) compressible Euler equations that are also \emph{stable}; that is, attracting an open set of initial data (with respect to a suitable topology). Such stable, non-isentropic implosions, if they exist, would need to be fundamentally different. For the physically important case of adiabatic exponent equaling $5/3$ (monatomic gas), in three dimensions, the discovery of these solutions would have immediate consequences, potentially extending the range of possible blowup speeds $r$, and hence the admissible values of $\gamma$ (cf.~\eqref{eq:hydro}) for which the results established in this paper may be established.

\subsubsection{Hydrodynamic limits }
The classical problem of hydrodynamic limits is generally focused on studying the limit $\epsilon \to 0$ in equations such as   
\begin{align*}
\partial_t f + v \cdot \nabla_x f = \frac{1}{\epsilon}Q(f,f)    
\end{align*}
over a fixed time window $t\in [0,T]$, which can be extended to a variety of kinetic models arising in gas dynamics,  plasmas, and a variety of corresponding macroscopic models; see e.g. \cite{saint2009hydrodynamic,bardos1991fluid,golse2005hydrodynamic,duan2023compressible,duan2024compressible,guo2010global,caflisch1980fluid,guo2010acoustic} and the references therein. Traditionally, this limit is performed using a Hilbert or Chapman--Enskog expansion (the latter providing the next order viscous corrections) and often work with smooth solutions, though sometimes the limits are only weak solutions. 
While our result is inspired by the idea of hydrodynamic limits, 
our approach differs from previous work. 
In particular, we fix a small $\e_0$ in~\eqref{eq:Landau} rather than taking $\e \to 0$ \cite{bardos1991fluid},
and obtain estimates for the Landau equation and its associated compressible Euler equations \emph{up to} the blowup time $T$, rather than only \emph{before} $T$ \cite{caflisch1980fluid,guo2010acoustic,duan2024compressible}.

Another closely related variation that is more related to this work, are the kinetic-level construction of rarefaction waves \cite{duan2021small}, contact waves \cite{duan2022asymptotics}, and especially inner shock-layers \cite{caflisch1982shock,metivier2009existence,albritton2024kinetic,wynter2024shock}, where a weak shock behaves similar to a sort of hydrodynamic limit at small, but \emph{fixed} $\epsilon$ (essentially the jump-size becomes the analogue of $\epsilon$). 
The weak shock proofs proceed by constructing an exact traveling wave solution of the Landau or Boltzmann equations by perturbing around a traveling wave solution of the compressible Navier--Stokes equations. These proofs are not dynamical and do not prove stability up to translations, only construct traveling waves. 
While these shock proofs are closer to our work than a traditional hydrodynamic limit, our result and proof proceed quite differently from any of these existing works. In our case, the limit only occurs in re-scaled time $s \to \infty$ and in a region shrinking to the origin $(x, v) = (0,0)$, involves the formation of a finite time imploding singularity in the macroscopic equations, and necessitates a fully dynamical argument which proves finite codimensional stability of the blowup.

\paragraph{Organization}

The rest of the paper is organized as follows. In Section \ref{sec:proof_idea}, 
we introduce the self-similar ansatz and outline the proof of Theorem \ref{thm:blowup}.
Section~\ref{sec:Euler_setup} discusses the properties of the Euler imploding profile 
and introduces the equations for the macro-perturbation. 
In Section \ref{sec:lin_macro}, we establish the linear stability estimates for the macro-perturbation. 
Sections \ref{sec:trilinear_Q} and \ref{sec:lin_micro} are devoted to estimates of the collision operator and to linear stability estimates for the micro-perturbation.
In Section \ref{sec:EE_top}, we estimate the top-order interaction terms between the macro and micro perturbations.
In Section \ref{sec:non_Q}, we estimate the nonlinear terms in the stability analysis.
Building on these estimates, we construct the blowup solution and prove Theorem \ref{thm:blowup} in Section \ref{sec:solu}.
Finally, in Section \ref{app:LWP_fixed}, we establish local well-posedness results for the fixed-point equations introduced in Section 
\ref{sec:solu} and for the Landau equation. Additional technical estimates and derivations are deferred to the Appendix.

\section{Self-similar ansatz and outline of the proof}\label{sec:proof_idea}

In this section, we develop the framework that reduces constructing finite time singularities in the Landau equation \eqref{eq:Landau} to establishing nonlinear stability of the local Maxwellian in self-similar variables. 
We first derive the self-similar ansatz for the Landau equation \eqref{eq:Landau} under a Type II scaling, 
and then consider its formal hydrodynamic limit to the compressible Euler equations. This converts the original equation into equation \eqref{eqn: linearized-td-F}. Then, we introduce the functional spaces and analytic framework for the stability analysis. In Section \ref{sec:idea}, we outline the steps for proving nonlinear stability and Theorem \ref{thm:blowup}.

\subsection{Self-similar blowup ansatz}

For any function $f$ and $ l > 0$, the collision operator $Q$ \eqref{eq:Q} satisfies the following scaling property in $v$
\[
	Q(f_{l}, f_{l})(v) = l^{-(\g+3)} Q(f, f)( l v),
	\quad f_{ l }(v) = f( l v).
\]

We consider the self-similar ansatz
\begin{align}\label{eq:SS}
	\bga
	f(t,x, v) = (1-t)^{c_f} F\B(s , \f{x}{(1-t)^{c_x}}, \f{v}{(1-t)^{c_v}} \B), \\
	s = - \log(1-t), \quad X = \f{x}{(1 -t)^{c_x}}, \quad V = \f{v}{( 1 -t)^{c_v}} ,
	\ega
\end{align}
where $c_f, c_x, c_v$ are time-independent blowup exponents.
\footnote{
In self-similar analysis, 
$c_f, c_x, c_v$ are commonly referred to as modulation parameters. They can be chosen to eliminate unstable or neutrally stable 
directions of the blowup profile that arise from scaling symmetries of the equation. 
Since in this work we establish only finite codimension stability—rather than full stability, we choose time-independent parameters $c_f, c_x, c_v$ to simplify the analysis.
}
The physical time $t=0$ corresponds to $s=0$ in the self-similar variables, 
and the blowup time $t=1$ corresponds to $s=\infty$. We have
\[
	\bal
	\pa_t f &= (1 -t)^{c_f - 1} (- c_f F + c_x X \cdot \na_X + c_v V \cdot \na_V + \pa_s F), \\
	v \cdot \na_x f &= (1 -t)^{c_f - c_x + c_v} V \cdot \na_x F, \\
	Q(f, f)(v) &= (1 -t)^{2 c_f + (\g + 3) c_v} Q(F, F) .
	\eal
\]

Choosing
\bseq\label{eq:LC_ss}
\beq
c_f - 1 = c_f - c_x + c_v \iff c_x = c_v + 1,
\eeq
and using the above identities, we get the self-similar equation of $F$
\beq\label{eq:self-similar-landau}
\pa_s F + (c_x X \cdot \na_X + V \cdot \na_X + c_v V \cdot \na_V) F = c_f F + \f1\es Q(F, F),
\eeq
where the Knudsen number in the self-similar equation is given by
\beq\label{eq:knud}
\es =\e _0 (1-t)^{- (c_f + (\g+3) c_v + 1)}
= \e_0 \operatorname e^{s (c_f + (\g+3) c_v + 1)} .
\eeq
\eseq

If $c_f + (\g+3) c_v + 1 = 0$, then the transport terms  and the collision term in \eqref{eq:self-similar-landau} have the same scaling, and we obtain the rate for Type I blowup. If, on the other hand
\beq\label{eq:hydro}
c_f + (\g+3) c_v + 1 < 0, 
\eeq
then we obtain $\es \to 0$ as $s \to \infty$; thus, we have a Type II scaling, and we are formally in a kind of hydrodynamic limit as $s \to \infty$, i.e., the self-similar equation becomes asymptotically collision dominated. 
Throughout this paper, we consider exponents satisfying the Type II scaling~\eqref{eq:hydro}.

\subsection{Local Maxwellian and the compressible Euler equations}

Since $\es \to 0$, to leading order as $s\to \infty$ one would formally expect the solution to be a local Maxwellian, with hydrodynamic fields governed by the Euler equations.
Hence, we construct the profile $\bar F$ as a local Maxwellian, with fields that can a priori depend on time $s$
\beq\label{eq:local_max}
\bar F = \cM_{\rho, \UU, \Th} = \rho(s, X) \f{1}{(2 \pi \Th(s, X))^{3/2}} \exp\B(-\f{|V - \UU(s, X)|^2}{2 \Th(s, X)} \B) .
\eeq
A direct computation yields
\beq\label{eq:moments}
\bal
\int \bar F d V &= \rho , &
\quad \int \bar F V d V &= \rho \UU, \\
\int \bar F V \otimes V d V 
&= \rho(\Th \Id + \UU \otimes \UU) , \quad &
\int \bar F |V|^2 V d V &= \rho \UU (5 \Th + |\UU|^2).
\eal
\eeq
We determine $\rho, \Th, \UU, c_x, c_v, c_f$ by integrating \eqref{eq:self-similar-landau} against $1, V_i, |V|^2$, which yields the self-similar equations for the full compressible Euler equations
\bseq\label{eq:Euler1}
\beq
\bal
(\pa_s + c_x X \cdot \na) \rho + \na \cdot (\rho \UU) &= (c_f + 3 c_v) \rho, \\
(\pa_s + c_x X \cdot \na) (\rho \UU) + \na \cdot (\rho (\Th \Id + \UU \otimes \UU)) &= (c_f + 4 c_v) \rho \UU, \\
(\pa_s + c_x X \cdot \na) (\rho (3 \Th + |\UU|^2))
+ \na \cdot (\rho \UU (5 \Th + |\UU|^2)) &= (c_f + 5 c_v) \rho (3 \Th + |\UU|^2).
\eal
\eeq
We introduce the pressure $P = \rho \Th$, 
\footnote{
	We use the law $P = \rho \msf{R} \Th$ with $\msf{R}=1$.
}
so we can rewrite the above system equivalently for $(\rho, \UU, P)$
\beq
\bal
\label{eqn: rhoUP}
& [ \pa_s + (c_x X + \UU)\cdot \na] \rho + \rho (\na \cdot \UU) = (c_f + 3 c_v) \rho, \\
& [ \pa_s + (c_x X + \UU)\cdot \na] \UU + \f1\rho \na P = c_v \UU, \\
& [ \pa_s + (c_x X + \UU)\cdot \na] P + \f53 P (\na \cdot \UU) = (c_f + 5 c_v) P.
\eal
\eeq
\eseq

Let $\kp = \f{5}{3}$ be the adiabatic exponent for monatomic gases.\footnote{
	We use $\kp = \f{5}{3}$ for the adiabatic exponent rather than 
    the canonical notation $\g$, since $\gamma$ is used to denote the exponent for the collision kernel \eqref{eq:Q}.    
} The pressure can be expressed by the ideal gas law in terms of the density $\rho$ and the specific entropy $\mathsf{s}$, or the \textit{pseudo entropy} $B$, via
\bseq
\label{eq:Euler_var}
\beq
P = \rho \Th = \f{1}{\kp} \rho^{\kp} \operatorname e^{\msf{s}}, \quad  B = \operatorname e^{\msf{s}}, \quad 
\kp = \f{5}{3}, 
\eeq
so that $P = \f{1}{\kp} \rho^{\kp} B$. Here $\Theta$ denotes temperature. 
We further introduce the sound speed $\sc$
\footnote{
	The sound speed $\sc$ differs from the rescaled sound speed $\s$ in \cite{chen2024Euler,chen2024vorticity}, and we do not use 
    the rescaled sound speed in this paper.
}
\beq
\sc = \sqrt{\f{ d P}{d \rho}} = \rho^{\f{\kp-1}2} B^\f12.%
\eeq
\eseq
Then $\rho$ and $P$ can be expressed in term of $\sc$ and $B$ as $\rho = \sc ^\f2{\kp - 1} B ^{-\f1{\kp - 1}}$, $P = \f{1}{\kp} \sc^{\f{2\kp}{\kp -1}} B^{- \f{1}{\kp-1}}$. Note that $\f P\rho = \f1\kp \sc ^2$, and
\[
	\f1\rho \na P = \f1\kp \sc ^2 \na \log P = \f{2}{\kp - 1} \sc \na \sc - \f{1}{(\kp-1) \kp} \sc^2 \f{\na B}{B}.
\]
Then we can rewrite \eqref{eq:Euler1} as the compressible Euler in terms of the unknowns $(\UU, \sc, B)$ as 
\beq\label{eq:Euler2}
\bal
&[ \pa_s + (c_x X + \UU)\cdot \na] \sc + \f{\kp - 1}2 \sc (\na \cdot \UU) = c_v \sc, \\
& [ \pa_s + (c_x X + \UU) \cdot \na] \UU + \f{2}{\kp - 1} \sc \na \sc = c_v \UU + \f{1}{(\kp-1) \kp} \sc^2 \f{\na B}{B}, \\
& [ \pa_s + (c_x X + \UU) \cdot \na] B = (1 - \kappa) c_f B .
\eal
\eeq

\subsubsection{Smooth implosion for the Isentropic Euler equations}\label{sec:isen_euler}

Taking $B \equiv 1$ and $c_f = 0$, we get the isentropic Euler system
\beq\label{eq:Euler3}
\bal
&[ \pa_s + (c_x X + \UU)\cdot \na] \sc + \f13 \sc (\na \cdot \UU) = c_v \sc, \\
& [ \pa_s + (c_x X + \UU) \cdot \na] \UU + 3 \sc \na \sc = c_v \UU.
\eal
\eeq
The recent work \cite{shao2025blow} extends the construction in \cite{merle2022implosion1,buckmaster2022smooth} on smooth imploding blowup solutions for the 3D isentropic Euler equations, in the case of monatomic gases ($\kappa = 5/3$); moreover,~\cite{shao2025blow} constructed a sequence of
smooth radially symmetric profiles $(\UU _n, \sc_n, B_n \equiv 1, c_{x, n}, c_{v, n})$ with
\beq\label{eq:r_limit}
c_{x, n} = \f{1}{r_n} ,\quad c_{v, n} = \f{1}{r_n} - 1 ,\quad
r_n \to (r_{*})^-, \quad %
r_{*} = \f{6}{3 + \sqrt{3}}=3-\sqrt{3},
\eeq

To achieve the Type II blowup condition that suggests a hydrodynamic limit to the 3D isentropic Euler equations, we need $\es \to 0$ (defined in \eqref{eq:knud}) as $s \to +\infty$. That is, $\es = \e_0 \operatorname e^{-\om s}$ with $\om$ defined as
\bseq\label{eq:hydro2}
\beq
\label{eqn: defn-omega}
\omega (\g, r) := -c _f - (\g+3) c_v - 1 = - (\g + 3) \left(\f{1}{r} - 1\right) - 1 = \g + 2 - \f{\g + 3}{r} > 0.
\eeq
In the limit $c_{v, n} = \f{1}{r_n}-1 \to \f{1}{r_*} - 1$, the constraint \eqref{eqn: defn-omega} reduces to 
\[
	\g > \sqrt{3} .
\]

\eseq

For each $\g > \sqrt{3}$, we choose $n$ to be large enough so that $3 - \sqrt{3} - r_n$ is small enough to ensure $\omega(\g,r_n)>0$ in~\eqref{eqn: defn-omega}. For such a chosen $n$, for the rest of the paper we fix the smooth radial profile $(\UU_n, \sc_n)$ which solves \eqref{eq:Euler3} with exponents determined by~\eqref{eq:r_limit} in terms of $r_n$. For simplicity, we denote %
\bseq\label{eq:Euler_profi}
\beq 
 \bu = \UU_n, \quad   \bar\sc = \sc_n,  \quad  r = r_n  , 
\eeq
and simplify $(c_{x, n}, c_{v,n})$ in \eqref{eq:r_limit} as $(\bcx, \bcv)$. We denote
\beq
    \label{eq:Euler_profi_cfcvcx}
\bar c_f = 0,  \quad 
	\bcv = \f{1}{r} - 1, \quad \bcx = \f{1}{r},
\eeq
and denote $\bu = \bar U \ee _R, \bar B \equiv 1$. The relation \eqref{eq:Euler_var} reduces to
\beq
    \label{eq:Euler_profi_rhoPth}
    \bar \rho = \bar{\sc}^3, \quad 
    \bar P = \f{1}{\kp} \bar \rho^{\kp} 
    = \f{1}{\kp} \bar{\sc}^5, \quad 
    \bar \Th = \f{1}{\kp} \bar \rho^{2/3}
    = \f{1}{\kp} \bar{\sc}^2 .
\eeq
\eseq

\subsubsection{Modified Euler profile}\label{sec: modified-profile}

\let\oldchi\chi
\renewcommand{\chi}{\raisebox{\depth}{$\oldchi$}}

To construct a blowup solution with non-vacuous density for large $|X|$, we modify the tail of the Euler profile. For $R_0 \gg 1$ which will be chosen to be sufficiently large, we define the time-dependent cut-off radius $\rs$ by  
\beq
\label{eq:S_radial}
\rs := R _0 \operatorname e ^{\bcx s} = R _0 \operatorname e ^{s/r}.
\eeq
Let $\chi \in C _c^{\infty} (\R ^3)$ be a radial cut-off function with $\one _{B _1} \le \chi \le \one _{B _2}$, so that $\chi _R (X) := \chi (X / R)$ becomes a cut-off function between $B _R$ and $B _{2 R}$,
where $B_a$ denotes the ball $\{ |X|: |X| < a\}$. We modify the profile in \eqref{eq:Euler_profi} based on the cutoff profile $( \bu, \cs)$ where
\begin{subequations}
\label{eq:Euler_profi_modi}
    \begin{equation}
    \cs := \bar \sc \chi _{\rs} + \rs ^{-r + 1} (1 - \chi _{\rs}),
    \end{equation}
and for consistency with \eqref{eq:Euler_var} we define
    \begin{equation}
    \bar \rho_s := \cs^3, 
    \quad 
    \bar \Th_s  := \f{1}{\kp} \cs^2, 
    \quad  
    \bar P_s := \f{1}{\kp} \cs^5
    .
    \end{equation}
\end{subequations}
The purpose of~\eqref{eq:Euler_profi_modi} is to replace the sound speed $\bar \sc$ by a (time-dependent) constant in the far-field, so that the profiles of 
$\brho, \bp$ are positive constants for large $|X|$. From the relation \eqref{eq:SS}, $|X| \approx \rs$ corresponds to $|x| \approx R_0$ in the physical variable.

\begin{remark}[Far-field profiles]\label{rem:far_profi}
From  \eqref{eq:dec_U}, we have $\bc \asymp \rs^{-(r-1)}$ for $|X| \in [\rs, 2 \rs]$. 
The term 
$\rs^{-(r-1)}$ in the second part in \eqref{eq:Euler_profi_modi} captures the correct scale of $\bc(X)$ 
 for $|X| \in [\rs, 2 \rs]$. 
In addition, instead of $\rs^{-(r-1)}$, we can choose another far-field profile in $\rs ^{-r + 1} (1 - \chi _{\rs})$ \eqref{eq:Euler_profi_modi} to obtain different far-field asymptotics of the macroscopic part of the blowup solution. 
For example, we can use $\bc$ as the profile and do not need the modification in \eqref{eq:Euler_profi_modi}. 
In that case, the associated density profile $\bar \rho(X)$ would vanish to $0$ as $|X| \to \infty$.  
\end{remark}

\begin{remark}[Growth rate $\rs$]\label{rem:exp_rate}
Note that the growth rate in $\rs$ defined in \eqref{eq:S_radial} is the same as the self-similar spatial rate $(1 - t)^{- \bcx} = \operatorname e^{\bcx s}$ in \eqref{eq:SS}.
From \eqref{eq:SS}, the far-field profile $\cs(X) =\rs ^{-(r - 1)}$ for $|X| > 2 \rs$ in the self-similar variables corresponds to a constant profile %
for $|x| > 2 R_0$ in the original physical variables. This choice of growth rate
is crucial for us to show that both the perturbation and residual error of the profile are relatively small. See Lemma \ref{lem: cutoff_error}.
\end{remark}

Next, we introduce a \textit{normalized relative velocity} $\vc$, which plays a fundamental role in our analysis, and is defined by
\begin{align}
\label{eqn: defn-vring}
\vc := \f{V - \bu}{\cs}.
\end{align}
To fix notation, we also let $\mu$ denote a specific  Gaussian:
	\beq\label{eq:gauss}
	\mu (x) = \left(\frac\kappa{2 \pi}\right) ^\frac32 \exp \left(
	- \kmu |x| ^2
	\right), \quad \kp = \f{5}{3},
	\quad \kmu = \f{\kp}{2} = \f{5}{6}.
	\eeq
The parameter $\kp_2$ appears naturally in the local Maxwellian \eqref{eq:local_max}.

In the rest of the paper, we denote by $ \cM, \cMM$ the time-dependent local Maxwellians 
\eqref{eq:local_max} 
\bseq\label{eq:localmax2}
\beq
    \cMM = \cM_{1, \bar \UU, \bar \Th_s}, \quad 
    \cM = \cM_{\brho, \bar \UU, \bar \Th_s} = \brho \cMM.
\eeq
Using the notation $\mu(\cdot)$ from \eqref{eq:gauss} and the cutoff profiles from \eqref{eq:Euler_profi_modi}, we can rewrite the local Maxwellian in terms of $\vc$
\beq
  \cMM = \cs^{-3} \mu(\vc), 
  \quad \cM =  \cs^3 \cMM = \mu(\vc).
\eeq
\eseq
Thus, the variable $\vc$ can be viewed as the normalized $V$ adapted to the local Maxwellian profile.

\paragraph{Error of the profile}

Since the modified profile does not solve the isentropic Euler equations exactly, we introduce the following  micro-error $\eM$ and macro-error $(\erho, \eu, \ep )$ associated with the above-defined profiles 
\bseq\label{eq:error0}
\beq
\bga
    \eM := (\pa_s + \bcx X \cdot \na_X + \bcx V \cdot \na _V + V \cdot \na_X) \cM, \\
    \erho := \cs^{-3} \la \eM, 1 \ra_\lvv, \quad 
    \eu   := \cs^{-4} \la \eM, V - \bu \ra_\lvv, \quad  
    \ep   := \cs^{-5} \left\la \eM, \f{1}{3} |V - \bu|^2 \right\ra _\lvv,
\ega
\eeq
where $\la \cdot, \cdot \ra_V  $ is defined as 
\beq\label{eq:inner_V}
 \la f , \, g \ra_V := \int f(V) g(V) d V .
\eeq
We introduce the relative error $\ec$ in solving the $\mathsf{C}-$equation \eqref{eq:Euler3} using the modified profile $(\cs, \bu)$: %
\begin{equation}\begin{aligned}
    \ec &=  \cs^{-1} \B( [\pa_s + ( \bcx X + \bu) \cdot \na] \cs + \f13 \cs (\na \cdot \bu) -  \bcv \cs \B) \\
    &= [\pa_s + ( \bcx X + \bu) \cdot \na] \log \cs + \f13 (\na \cdot \bu) - \bcv.
\end{aligned}\end{equation}
\eseq
The errors $(\erho, \eu, \ep , \ec )$  are supported in the far-field $|X| \geq \rs$ and have appropriate decay as $|X|\to \infty$. We estimate  these errors in Lemma~\ref{lem: cutoff_error} of Appendix~\ref{app:error_est}.  We defer the computation of $\eM$ to Lemma~\ref{lem:ds_t}.

\subsection{Decomposition of the perturbation}\label{sec:lin_decomp}

Let  $ \cM, \cMM$ be the local Maxwellians defined in \eqref{eq:localmax2}. Our goal is to construct a global solution to the self-similar equation \eqref{eq:LC_ss} near the local Maxwellian $\cM$. To this end, we decompose the full solution to \eqref{eq:LC_ss} as
\footnote{
	We renormalize the perturbation by $\cMM^{1/2}$ rather than $\cM^{1/2}$ since the density in $\cMM$ is $\rho \equiv 1$ and it is more convenient to define orthogonality using $\cMM^{1/2}$.
}
\beq\label{eq:pertb_dec}
F = \cM + \cMM^{1/2} \td F .
\eeq
Denote
\beq\label{eq:func_Phi}
\begin{aligned}
    \Phi_0 &= \cMM ^{1/2}, \\
    \Phi_i &= \f{V_i - \bu_i}{\ths^{1/2}} \cMM^{1/2} = \kp ^{1/2} \vc _i \cMM ^{1/2}, &i = 1, 2, 3,\\
    \Phi_4 &= \f{1}{\sqrt 6} \left(
        \f{|V - \bu| ^2}{\ths} - 3 
    \right) \cMM ^{1/2} = \f{\kp}{\sqrt 6} \left(
        |\vc| ^2 - \f95 
    \right) \cMM ^{1/2}.
\end{aligned}
\eeq
By elementary computation, we have
\[
	\la \Phi_i , \Phi_j \ra_\lvv = \d_{ij} .
\]

For a function $g \in L ^2 (V)$, we use $\cP _M$ and $\cP _m$
\footnote{
We use the calligraphic font to denote the operator $\cP$ and $P$ for the pressure to avoid confusion.
}
to denote the projection onto the macro and micro parts
\bseq\label{def:proj}
\beq
\cP_M g := \sum _{0 \le i \le 4} \la g, \Phi_i \ra _{\lvv} \Phi_i, 
\qquad 
\cP_m g := (I -\cP_M) g.
\eeq
We denote the macro part and the micro part of $\td F$ by:
\beq
\tFM := \cP _M \td F, \qquad \tFm := \cP _m \td F .
\eeq
\eseq

\begin{remark}[\bf Macro- and micro-perturbations]
Throughout the paper, we refer to $\tFM$ as the macro-perturbation and $\tFm$ as the micro-perturbation.
\end{remark}

Denote the transport operator $\cT$ and the symmetric linear collision operator $\cL _\cM$ by 
\bseq\label{def:lin_op}
\beq
\label{def:T-LM}
\bal
\cT g & := %
(V \cdot \na_X + \bcx X \cdot \na_X + \bcv V \cdot \na_V) g, \\
\cL_{\cM} g &:= \cMM^{-1/2} \B[ 
    Q(\cM, \cMM^{1/2} g) + Q(\cMM^{1/2} g, \cM)
\B],
\eal
\eeq
and denote the nonlinear collision operator $\cN$ by
\beq\label{eq:non_nota}
 \cN(f, g) = \cM_1^{-1/2} Q( \cM_1^{1/2} f, \cMM^{1/2} g  ) .
\eeq
We denote the following moments from the micro part
\beq\label{eq:moments_Fm}
\bal
\cI_1(\tFm) &:= \big\la V \cdot \na_X (\cMM^{1/2} \tFm), \vc \big\ra_\lvv, \\
\cI_2(\tFm) & := \B\la V \cdot \na_X (\cMM^{1/2} \tFm)\, , \f13 |\vc| ^2 \B\ra_\lvv .\\
\eal
\eeq
\eseq

Recall the exponents $\bcx, \bar c_f, \bcv$ from\eqref{eq:Euler_profi}. We choose time-independent blowup exponents  $c_x , c_v , c_f$:
\[
	c_f = \bar c_f  = 0 ,\quad c_x = \bcx = r^{-1} , 
	\quad c_v = \bcv = r^{-1} - 1.
\]

Linearizing \eqref{eq:LC_ss} around the local Maxwellian \eqref{eq:localmax2}, we obtain the linearized equation for the perturbation
\bseq
\label{eq:lin}
\beq
\label{eqn: lin-with-M11/2}
\bal
& \pa_s (\cMM^{\f12} \td F) + \cT (\cMM^{ \f12} \td F) 
 = 
\f1\es \B[ 
    Q(\cM, \cMM^{ \f12 } \td F)
	+ Q(\cMM^{ \f12} \td F, \cM)
	+ Q(\cMM^{ \f12} \td F , \cMM^{ \f12 } \td F) \B]
- \eM ,
\eal
\eeq
where $\eM = (\pa _s + \cT) \cM$ is the error defined in \eqref{eq:error0}. We derive the linearized equation of $\td F$ by dividing $\cMM^{1/2}$:
\begin{align*}
    (\pa_s + \cT) \td F + \f12 (\pa _s + \cT) \log \cMM  \cdot \td F = \f1\es \cL _\cM \td F + \f1\es \cN (\td F, \td F) - \cMM ^{-1/2} \eM,
\end{align*}
where the $\log \cMM$ term comes from
\begin{align*}
	\cMM^{-1/2} (\pa_s + \cT) \cMM ^{1/2} = \frac12 (\pa _s + \cT) \log \cMM.
\end{align*}
The leading order term of this term is $-\f32 \bcv$ (see \eqref{eqn: Dab-dcm-bound}), so we introduce $\dcm$ and write 
\begin{align*}
    \f12 (\pa _s + \cT) \log \cMM = \dcm - \f32 \bcv.
\end{align*}
The linearized equation for the perturbation $\td F$ is thus 
\beq 
    \label{eqn: linearized-td-F}
    \left(\pa_s + \cT + \dcm - \f32 \bcv \right) \td F = \f1\es \cL _\cM \td F + \f1\es \cN (\td F, \td F) - \cMM ^{-1/2} \eM .
\eeq
\eseq

We aim to construct a non-trivial global-in-time solution to \eqref{eq:lin} by establishing the nonlinear stability estimates of $\td F$, upon modulating finitely many unstable directions. Using the self-similar transform \eqref{eq:SS} with exponents \eqref{eq:Euler_profi}, we construct a finite time singularity in the Landau equation \eqref{eq:Landau}. We introduce the functional spaces in the next subsection and outline the proofs in Section \ref{sec:idea}. 

\subsection{Functional spaces and weighted derivatives}\label{sec:function_space}

\subsubsection{Weighted derivatives}

In view of the local Maxwellian \eqref{eq:localmax2}, we introduce weighted $X, V$-derivatives to capture the scaling of the profile $\cM$ and the perturbation. Let $\vp_1 = \vp _1 (X)$ be a positive smooth weight to be designed in Lemma \ref{lem:wg}. For any multi-indices $\al, \b \in 
 \Znonneg^3$,  we define (weighted) derivatives 
\footnote{
	Note that we do not have $ D_X^{2 \al} = (D_X^{\al})^2$ or similar identities since we first take derivatives and then multiply the weights. They do agree up to lower order terms: see Corollary \ref{cor: concatenate-derivative}.
}
\beq\label{eq:deri_wg}
	D_X^{\al} := \vp_1^{|\al|} \pa_X^{\al}, \quad
	D_V^{\b} := \cs^{|\b|} \pa_V^{\b}, \quad
	D^{\al, \b} := \vp_1^{|\al|} \cs^{|\b|} \pa_X^{\al} \pa_V^{\b},
	\quad \pa_{X, V}^{(\al, \b)} := \pa_X^{\al} \pa_V^{\b} ,
\eeq
where we denote $\pa_Z^{\th}= \pa_{Z_1}^{\th_1} \pa_{Z_2}^{\th_2} \pa_{Z_3}^{\th_3}$
for $Z = X, V \in \R^3$ and $\th = \al, \b$.  By definition, the weighted derivatives satisfy the usual Leibniz rule. 
    
For $k \geq 0$, we define the tensor $ D^{\leq k} f$ spanned by the mixed derivatives of $f$
\beq\label{eq:mix_deri}
\bal
D^{\leq k} f & :=
\{D^{\al, \b} f\}_{|\al|+ |\b| \leq k},
\quad
| D^{\leq k} f |
= \B(\sum_{|\al| + |\b| \leq k} |D^{\al, \b} f|^2 \B)^{1/2}, \\
\eal
\eeq
and $D ^{< k} f = \one_{k>0} D ^{\le {k - 1}} f$. %
Similarly, we define the tensor $D ^{\preceq (\al, \b)}$ spanned by the mixed derivatives of $f$:
\beq
\bal 
    D ^{\preceq (\al, \b)} &= \{D ^{\al', \b'} \} _{\al' \preceq \al, \b' \preceq \b}, &
    |D ^{\preceq (\al, \b)}| = \B( 
        \sum \nolimits_{\substack{\al' \preceq \al \\\b' \preceq \b}} |D ^{\al', \b'} f| ^2
    \B) ^{1/2}.
\eal
\eeq
$D ^{\prec (\al, \b)}$ are defined analogously, where $(\al', \b') \prec (\al, \b)$ if $\al' \preceq \al$, $\b' \preceq \b$, and at least one of the inequality $\preceq$ is strict $\prec$.

\paragraph{Motivation of $D^{\al,\b}$}
Following  \cite{chen2024vorticity,chen2024Euler}, we use the weight $\vp_1 \asymp \ang X$ to capture the flow structure in the compressible Euler equations for stability analysis.  See details in Section \ref{sec:lin_macro}. We weight $\pa_V$ by the standard deviation in $\cM$ \eqref{eq:localmax2} so that the weighted operator $D_V$ is similar to $\pa_{\vc}$. For example, we have $D _V \vc _i = \ee _i$ is the $i$-th basis vector, and $D _V \mu(\vc) =  (\na \mu) (\vc)$.
     
A crucial property of $D^{\al,\b}$ is that it commutes with the self-similar flow  $\pa_s + \cT$ in equations \eqref{eq:lin}, up to lower order terms that decay faster in $X$. See Lemma \ref{lem: commu-derivative-projection} (2). We use this property crucially to perform sharp decay estimates for $\td F$ 
in $X$ and its higher order derivatives. See Lemma \ref{lem:prod}. Moreover, in many cases, $D^{\al, \b}$ behaves similarly to a constant multiplier and simplifies many estimates. 

\subsubsection{Weighted Sobolev norms: $\s$, $\cX$, $\cY$-norms}
\label{sec: norms} Now we introduce function spaces in $X$ and $V$.

\paragraph{$\cX$-norm} 


For hydrodynamic fields $\WW = (\UU, P, B)$ which are functions of $X$, we equip the following norm.
For any $k \geq 0$, $\eta \in \R$, and some parameter $\vpi_{k, \eta} $ determined in Theorem \ref{thm:coer_est},  we introduce the $\cX _\eta ^{2k}$-norm to analyze $\WW$ and the Euler equations 
\begin{align}
\la (\UU_a, P_a, B_a) , ( \UU_b,  P_b,  B_b) \ra_{\cX_{\eta}^{2k}}
& :=  \int  \sum_{g = \UU, P, B } w_{g} (   \D^k g_a \cdot \D^k  g_b  \; \vp_{2 k}^ 2  + 
\vpi_{k, \eta}
g_a \cdot  g_b ) \; \la X \ra^{\eta}   d X,  \ k \geq 1,  \notag \\
 (w_\UU, w_P, w_B ) & :=  (1, 1, \tf{3}{2}) , 
\label{norm:Xk0}
\end{align}
The $\cX^{2 k+1} _\eta$-norm is defined similarly; see \eqref{norm:Xk}.

\paragraph{$\s$-norm}

For the micro-perturbation, we first introduce the $\s$-norm by generalizing the corresponding norm in \cite{guo2002landau}.
For a function $g$, we define the collision norm in $V$ as (recall $\cM (V) = \mu (\vc)$)
\beq\label{norm:sig}
\bal 
    \| g \| _\s ^2 &:= \int _{\R ^3} A [\cM] \na _V g \cdot \na _V g + \kp _2 ^2 \cs ^{-2} A [\cM \vc \otimes \vc] g ^2 d V \\
    &= \cs ^{\gamma + 5} \int _{\R ^3} A [\mu]  (\vc) \na _V g \cdot \na _V g d V + \kp _2 ^2 \cs ^{\gamma + 3} \int _{\R ^3} A [\mu \vc \otimes \vc] (\vc) g ^2 d V,
\eal 
\eeq 
where we recall the definition of $A [f]$ from \eqref{eqn: defn-bPhi-A}:
\begin{align*}
    A [f] (V) = \int \boldsymbol \Phi (V - W) f (W) d W = \int |V - W| ^{\gamma + 2} (\Id - \bpi _{V - W}) f (W) d W .
\end{align*}
Here $\bpi _V = \f V{|V|} \otimes \f V{|V|}$ is the projection along the $V$ direction. We define $A [f]$ the same way if $f$ is a vector-valued function, and for matrix-valued function $\boldsymbol f$ we define
\begin{align*}
    A [\boldsymbol f] (V) = \int \boldsymbol \Phi (V - W) : \boldsymbol f (W) d W = \f1{8\pi} \int |V - W| ^{\gamma + 2} (\Id - \bpi _{V - W}) : \boldsymbol f (W) d W,
\end{align*}
with matrix product $\boldsymbol f : \boldsymbol g = \sum _{i = 1} ^3 \sum _{j = 1} ^3 f _{i j} g _{i j}$.

\paragraph{$\cY$-norm}
Now we introduce new norms that also take into account the $X$ variable. 
For every $\eta \in \R$, we define 
\bseq\label{norm:Y}
\beq
\bal 
    \| g \| _{\cY _\eta} ^2 &:= \int \ang X ^\eta \| g \| _{L ^2 (V)} ^2 d X, &
    \| g \| _{\cYE} ^2 &:= \int \ang X ^\eta \| g \| _\sigma ^2 d X.
\eal 
\eeq
Then we introduce the $H ^k$ counterparts of these norms:
\beq\label{norm:Y:b}
\bal
    \| g \| _{\cYe ^k} ^2 &:= \sum _{|\al| + |\b| \le k} \nu ^{|\al| + |\b| - k} 
    \f{ |\al|! }{ \al !} \int \ang X ^\eta \| D ^{\al, \b} g \| _{L ^2 (V)} ^2 d X, \\
    \| g \| _{\cYE ^k} ^2 &:= \sum _{|\al| + |\b| \le k} \nu ^{|\al| + |\b| - k}
    \f{ |\al|! }{ \al !}  \int \ang X ^\eta \| D ^{\al, \b} g \| _\sigma ^2 d X,
\eal
\eeq
\eseq
with constant coefficients $\nu ^{|\al| + |\b| - k}$ depending on $\nu \ll 1$ to be determined in Theorem \ref{thm: micro-Hk-main},
where $\al, \b \in  \Znonneg^3$ are multi-indices, $n! = \prod_{1\leq i\leq n} i$ denotes the factorial,  and  $\al ! = \al_1 ! \al_2! \al_3!$ for multi-index $\al = (\al_1, \al_2, \al_3)$.  Note that $\cYe$ and $\cYE$ coincide with $\cYe ^0$ and $\cYE ^0$.
We define $\la  \cdot , \cdot \ra_{\cYe ^k}$, $\la  \cdot , \cdot \ra_{\cYE ^k}$ to be the inner product associated with the norm $\cYe ^k$ and $\cYE ^k$. 

We choose the weight $\nu^{k-k}=1$ for the top-order derivative terms with $|\alpha|+|\beta|=k$ in the $\cY_{\eta}^k$-norm to facilitate the top-order estimates in Section~\ref{sec:EE_top}. The multiplicity constant $   \f{ |\al|! }{ \al !} $ in \eqref{norm:Y:b} arises from the difference between two sums:
\beq\label{eq:two_summation}
\sum_{ |\al| = n} \f{ |\al|!}{ \al !} H( \pa^{\al} _X g_1,  \pa^{\al} _X g_2  ) = \sum_{i_1, .., i_n \in \{1,2,3\} } 
H( \pa_{X_{i_1}} .. \pa_{X_{i_n}} g_1, \pa_{X_{i_1}} .. \pa_{X_{i_n}} g_2 ),
\eeq
for any function $g_i$ and functional $H$. The left hand side sums over different multi-indices $\al = (\al_1, \al_2, \al_3)$ with $|\al| = n$. The proof follows from a simple combinatorial calculation and is omitted.

\paragraph{Critical exponent}
By choosing different exponents $\eta$ in \eqref{norm:Xk0} and \eqref{norm:Y}, we obtain different coercivity estimates of the linear operators in \eqref{eq:lin}. See estimates \eqref{eq:idea_macro_EE},\eqref{eq:idea_micro_EE}. 
       We define %
\beq\label{wg:X_power}
  \etab =  - 3 + 6 (r-1) .
\eeq

Under the self-similar scaling for the density $\rho$: $\rho _\lambda = \lambda ^{ 3 \bcv } \rho \big( \f{X}{ \lambda^{ \bcx}} \big)
= \lambda^{  -3(1-1/r) } \rho \big( \f{X}{ \lambda^{1/r}} \big)$ with $\lambda >0$,
\footnote{
Using \eqref{eq:SS} with $c_f = 0$ and $\mass$ defined in \eqref{eq:density}, one can show that the self-similar ansatz for density is :
$\mass = (1 - t)^{3 c_v } \rho( \f{X}{(1-t)^{c_x}}) $.
}
the exponent $\etab$ can be viewed as weighted $L^2$-critical since the following $\etab$-weighted norm is invariant %
\beq\label{eq:inv_etab}
   \int \rho^2 |X|^{\etab} d X = \int \rho_\lambda ^2 |X|^{\etab} d X .
\eeq

In Section \ref{sec:lin_euler}, we normalize the hydrodynamic variables $\tw = (\tu, \tp, \tb)$ for the perturbation so that they have the same scaling as the density $\rho$.

\subsection{Steps and ideas of the proof}\label{sec:idea}

To establish the nonlinear stability of the perturbation and construct global solutions to \eqref{eq:lin}, our argument proceeds in the following steps.

\subsubsection*{\underline{Step 1.} Decomposing the perturbation.}
\label{sec: step1}

We decompose the perturbation $\tilde F$ into a macroscopic part $\tFM$, a microscopic part $\tFm$,  
and derive the equations for $\tFM$ and $\tFm$ using \eqref{eqn: linearized-td-F}.
For $\es$ sufficiently small, the evolutions of $\tFM$ and $\tFm$ are weakly coupled via the kinetic transport term $V \cdot \na_X \td F$ and the nonlinear terms. This structure allows us to essentially decompose the whole stability analysis into proving the stability of the macro-perturbation and micro-perturbation separately.

\paragraph{Size of perturbation}

We design the $\cX^k_{\eta}$ norm (see~\eqref{norm:Xk0}, or \eqref{norm:Xk}) to analyze the hydrodynamic fields of the macro-perturbation
$\tw = (\tu, \tp, \tb)$ (see \hyperref[sec: step2]{\itshape \underline{Step 2}}), and the $\cYe ^k$ norm (see~\eqref{norm:Y}) to analyze the micro-perturbation; here $k$ indicates the regularity index and $\eta$ indicates the power of $|X|$ in the weight.

Let $\etab$ be the exponent in \eqref{wg:X_power}. 
We choose weights with two exponents $\etab$ and $\xeta$ with $\xeta < \etab$ in order to capture different temporal and spatial decays; we also  choose two regularity indices $k, k+1$ with $k$ sufficiently large.

\paragraph{Exponential decay estimates}
In the norm with faster-decaying weights indicated by $\xeta$,  we aim to establish 
\bseq\label{eq:idea_size}
  \begin{align}
 \| \tw(s) \|_{\cXX ^{2 k + 2} }  , \ \| \tFm(s) \|_{\cYY ^{2 k+ 2}} & < \es^{1/2- \el},
 \quad  \el = 10^{-4},  \label{eq:idea_xeta_H}  \\
\| \tw(s) \|_{\cXX^{2k}}  & < \es^{2/3} ,  \label{eq:idea_xeta_L} 
  \end{align}
for any $s \geq 0$. We fix $\el >0$ to be a sufficiently small absolute parameter and use $\es^{-\el}$ for small $\e_0$ to absorb any large absolute constants. Recall $\es = \e _0 \operatorname e ^{- \om s}$ from \eqref{eqn: defn-omega}.

\paragraph{Relative smallness estimates}

The estimates in \eqref{eq:idea_size} yield non-sharp \textit{spatial} decay estimates for the perturbation at large $|X|$ and are insufficient to close the nonlinear estimates. To overcome this, we also work in the norms with critical decaying weights indicated by $\etab$ and aim to establish
\beq\label{eq:idea_etab}
 \| \tw(s) \|_{\cXb ^{2 k + 2} } , \quad  \| \tFm(s) \|_{\cYb ^{2 k + 2}} < \d^{ \el }, \quad \forall s \geq 0. 
\eeq 
We choose $\e _0 = \d$ to be sufficiently small after we fix the parameters $k, \xeta, \etab, \el$.
\eseq

For $|X|$ sufficiently large, due to the decay of the mass $\rhos (X)$ and variance $\ths (X)$, the 
coercivity estimates obtained from the linear collision operator $\cL_{\cM}$ become much weaker. 
Moreover, since the far-field of the $\etab$-based norm $\cXb^{n}, \cYb^{n}$ ($n = 2 k$ or $2 k + 2$) is almost invariant under the self-similar scaling (based on identities similar to \eqref{eq:inv_etab}), 
the self-similar scaling fields $\bcx X\cdot\nabla_X + \bcv V\cdot\nabla_V$ in \eqref{def:lin_op}-\eqref{eq:lin}
do not generate a damping effect for large $|X|$ in the $\cXb^n, \cYb^n$-estimates. As a result, we can establish only smallness, rather than decay, estimates for the perturbation in the norms in \eqref{eq:idea_etab}. Estimate \eqref{eq:idea_etab} implies that $\td\rho$ is small \textit{relative} to its profile $\bar\rho_s$ and implies \textit{relative smallness} of similar variables 
and their weighted derivatives via the embedding inequalities in Lemma~\ref{lem:prod}.

In the $\xeta$-based norm, with $\xeta < \etab$, we use the stability mechanism from the scaling fields
$ \bcx X \cdot \na_X + \bcv V \cdot \na_V$ in \eqref{def:lin_op}-\eqref{eq:lin}, leading to the exponential decay estimates in \eqref{eq:idea_xeta_L}, \eqref{eq:idea_xeta_H}. 
By contrast, in the $\cX^n_{\eta}$ and $\cY^n_{\eta}$ energy estimates with $\eta>\etab$, 
the scaling fields $\bcx X\cdot\nabla_X + \bcv V\cdot\nabla_V$ induce an \emph{anti-damping effect}, so that
the perturbation in these norms is expected to grow \emph{exponentially fast}.
 Therefore, the choice of the functional spaces in \eqref{eq:idea_size} is crucial for establishing
nonlinear stability.

\subsubsection*{\underline{Step 2}. Finite codimension stability of macro-perturbation}
\label{sec: step2}

The macro-perturbation $\tFM$ is governed by the linearized Euler equations around the isentropic imploding profile $( \brho, \bu, \bp )$. 
Due to the weak coupling with the micro-perturbation, we need to analyze perturbations to density, velocity and pressure $(\trho, \tu, \tp)$, 
which are not covered by  previous stability analyses of isentropic implosions~\cite{merle2022implosion2,chen2024Euler,buckmaster2022smooth,cao2024non}. 
Instead of estimating the system evolving $(\trho, \tu, \tp)$, we introduce a variable $\tB$ related to the entropy, and perform estimates on the system for $\tw = (\tu, \tp, \tB)$, which is symmetric and hyperbolic.  The finite codimension stability of the profile relies on the interior repulsive property 
\eqref{eq:rep1}, and the outgoing property of the profile \eqref{eq:rep2}. 
We generalize the finite codimension stability argument  developed in \cite{chen2024Euler,chen2024vorticity}, and perform weighted linear stability estimates in Section \ref{sec:lin_macro}. 

Using these stability estimates and applying the splitting method \cite{ChenHou2023a} to the perturbation $\tw = \tw_1 + \tw_2$,  up to lower order terms, we obtain 
\bseq\label{eq:idea_macro_EE}
\begin{align}
  \f{1}{2} \f{d}{d s} \| \tw_1 \|_{\cXX^{ 2 k_1} }^2 
  & \leq - \lam_1 \| \tw_1 \|_{\cXX^{ 2 k_1} }^2  
   +  \B\la  \tw_1,  (  - \cI_1 , - \cI_2, \cI_2)(\tFm) \B\ra_{\cXX^{ 2 k_1}} 
 + l.o.t.,
 \label{eq:idea_macro_EE:a}  \\
   \f{1}{2} \f{d}{d s} \| \tw \|_{\cXb^{ 2 k + 2 } }^2 
   & \leq C_k \| \tw \|_{\cXX^{2 k + 2 } }^2  
     +   \B\la  \tw,  (  - \cI_1 , - \cI_2, \cI_2 )(\tFm)  \B\ra_{\cX_{\etab}^{2 k + 2 }} 
   + l.o.t. \label{eq:idea_macro_EE:b}
\end{align}
\eseq
for $k_1 = k, k+1$, where  $\kp = \f{5}{3}$, $\lam_1 > 0$ is independent of $k$ and is defined in  \eqref{eq:eta_constraint}. See the discussion of decay exponents at the end of this Section, e.g. \eqref{eq:eta_constraint}. 
Note that the $\cXb^{2k+2}$-energy estimates \eqref{eq:idea_macro_EE:b} with the critical exponent $\eta=\etab$ do not contain any damping terms.  It is therefore important that the upper bound in \eqref{eq:idea_macro_EE:b} involves the $\xeta$-norm $\|\tw \|_{\cXX^{2k+2}}^2$ rather than $\|\tw \|_{\cXb^{2k+2}}^2$, since the former will be shown to decay exponentially fast in time \eqref{eq:idea_size}.

The perturbation $\tw_2$ captures the potential unstable modes and is small relative to $\tw_1$.
We treat it as a sufficiently smooth forcing, 
and discuss its estimate in \hyperref[sec: step6]{\itshape \underline{Step 6}}.

\subsubsection*{\underline{Step 3}. Full stability of micro-perturbation}
\label{sec: step3}

To control the micro-perturbation $\tFm$, we use the coercivity estimates of 
the linear collision operator $\cL_{\cM}$ defined in \eqref{def:T-LM}, inspired by  \cite{guo2002landau}. Thanks to the small parameter $\es$ in \eqref{eq:lin}, we can treat the transport term $ V \cdot \na_X $ in \eqref{eq:lin} as a perturbation of the coercive linear part.
However, since the density and temperature decay spatially, these coercivity estimates weaken for large $|X|$. In this region, we instead combine the coercivity of $\cL_{\cM}$ and the stability effects generated by the scaling fields $ \bcx X \cdot \na_X + \bcv V \cdot \na_V$ in \eqref{def:lin_op}. 

We develop these estimates in Sections \ref{sec:trilinear_Q} and \ref{sec:lin_micro} and establish
\footnote{
In \eqref{eq:idea_micro_EE:b}, we may replace the term 
$ C_k \es \| \tFm \|_{\cYb^{2\kk+2}}^2$ by the upper bound $ C_k \| \tFm \|_{\cYY^{2\kk+2}}^2$. Both estimates are sufficient to prove nonlinear stability estimates. 
}
\bseq\label{eq:idea_micro_EE}
\begin{align}
  \f{1}{2} \f{d}{d s} \| \tFm \|_{\cYY^{ 2 k + 2} }^2 
   & \leq - \lame  \| \tFm \|_{\cYY^{ 2 k  + 2} }^2 
   - \f{ \cgam }{8 \es} \| \tFm \|_{\cYL{\xeta}^{ 2 k + 2}}^2
 +   \la  \tFm ,  V \cdot \na_X \tFM  \ra_{\cYY^{2 k + 2}}  \notag \\ 
 & \qquad + 
 \es^{-1} \la \cN(\td F, \td F) , \tFm \ra_{\cYY^{ 2 k + 2}} +    C\es  +  l.o.t.,  \label{eq:idea_micro_EE:a}  \\
  \f{1}{2} \f{d}{d s} \| \tFm \|_{\cYb^{ 2 k + 2 } }^2 
    & \leq C_k \es \| \tFm \|_{\cYb^{2 k + 2 } }^2  
    - \f{ \cgam  }{8 \es} \| \tFm \|_{\cYQ^{2 k +2}}^2 
 +   \la  \tFm ,  V \cdot \na_X \tFM  \ra_{\cYb^{ 2 k + 2}}  \notag  \\ 
 & \qquad  + \es^{-1} \la \cN(\td F, \td F) , \tFm \ra_{\cYb^{ 2 k +2 }}  + C \es  + l.o.t.  \label{eq:idea_micro_EE:b}
\end{align}
\eseq
where $\cgam > 0$ and $\lame > \lam_1 >0$ is independent of $k$ and is defined in \eqref{eq:eta_constraint:c}. The term $\f{1}{\es} \| \tFm \|_{\cYL{\xeta}^{2 k + 2 }}^2$ is from the coercivity estimates of $\cL_{\cM}$, and \textit{l.o.t.} denotes lower order terms that can be treated perturbatively.
We refer to \eqref{norm:Y} for the definition of the norms $\cY_{\eta}^k, \cYL{\eta}^k$. Similar to \eqref{eq:idea_macro_EE:b}, the $\cYb^{ 2 k + 2  }$-energy estimate \eqref{eq:idea_micro_EE:b} does not contain any damping terms of the form $- C \| \tFm \|_{\cYb^{ 2 k + 2  } }^2 $.

\subsubsection*{\underline{Step 4}. Coupled estimates of macro and micro-perturbation}
\label{sec: step4}

At the linear level, the equations for the macro $\tFM$ and micro $\tFm$ perturbations are coupled via the transport term $V \cdot \na_X$, which could potentially lead to a loss of derivatives. 
See estimates \eqref{eq:idea_macro_EE} and~\eqref{eq:idea_micro_EE}. 

In Section \ref{sec:EE_top}, we show that the cross term in \eqref{eq:idea_macro_EE} can be rewritten as  
 \beq\label{eq:idea_top_cross}
        \B\la  \tw_1,  (  - \cI_1 , - \cI_2, \cI_2 )(\tFm)  \B\ra_{\cX_{\eta}^{2 k_1}} 
=  \f{1}{\kp}  \la  \td \cF _M (\tw _1) ,  V \cdot \na_X \tFm  \ra_{\cYe ^{2 k_1}}  + l.o.t.  , 
\quad k_1 =k , k+1, 
 \eeq
 up to lower order terms that can be bounded perturbatively; here $\cF _M (\tw _1)$ is the macro-perturbation with hydrodynamic fields $\tw_1$. 
To avoid loss of derivatives, at the top order weighted $H^{2k+2}$ estimates, we choose a specific energy norm
$\kp \| \tw_1 \|_{\cX_{\xeta}^{2 k+ 2 } }^2
+  \| \tFm \|_{\cYY ^{2 k+ 2}}^2$
which satisfies 
\[
\kp \| \tw_1 \|_{\cX_{\xeta}^{2 k + 2 } }^2
+  \| \tFm \|_{\cYY ^{2 k+ 2 }}^2
 = \sum_{|\al| + |\b| = 2 k+ 2} \int   \vp_{k, \xeta}(X) ( | D^{\al, \b} \cF _M (\tw _1) |^2 + | D^{\al, \b} \tFm|^2 )  d X d V+ l.o.t.
\]
for some weight $\vp_{k ,\xeta}$, where \textit{l.o.t.} contains terms involving $X$-derivatives of order at most $2k+1$.
Similarly, we estimate $\kp \| \tw \|_{\cX_{\etab}^{2 k+ 2 } }^2 +  \| \tFm \|_{\cYb^{2 k + 2 }}^2$ to close the energy estimates \eqref{eq:idea_macro_EE:b} and \eqref{eq:idea_micro_EE:b}.  The above structure allows us to combine 
the estimates of the cross terms in \eqref{eq:idea_macro_EE} and \eqref{eq:idea_micro_EE} and to perform integration by parts,
thereby transferring the $X$-derivatives onto the weight.

\subsubsection*{\underline{Step 5}. Estimate of nonlinear terms}
\label{sec: step5}

We aim to treat the nonlinear terms $\cN$ defined in \eqref{eq:non_nota} as a small perturbation of the linear coercive part. 
There are two difficulties. Firstly, the local Maxwellian is spatially dependent with states decaying in $X$. Secondly, the coefficient  $\es^{-1}$ of the nonlinear terms \eqref{eq:lin} grows exponentially.
To overcome the first difficulty, we design careful weighted estimates in Sections \ref{sec:lin_micro} and \ref{sec:non_Q}.
To overcome the second difficulty, in Section \ref{sec:non_Q}, we establish two types of nonlinear estimates for $\cN(f, g)$ \eqref{eq:non_nota}, based on the decomposition 
\footnote{
Note that $\cN$ is a bilinear operator and $\td F = \tFm + \tFM$.
}
\beq\label{eq:idea_decom_N}
  \cN(\td F, \td F) = \cN(\td F, \tFm) + \cN( \tFm, \tFM )
  + \cN(\tFM, \tFM)  
:= \cN_{m} + \cN_{m M} + \cN_{MM}.
\eeq
The first two terms contain the micro-perturbation $\tFm$. 
For the first term, we establish 
\beq\label{eq:idea_decomp_N2}
\bal
   \B|\f{1}{\es} \la \cN( \td F , \tFm ), \tFm \ra_{\cYe^{2 k + 2}} \B|
  & \les 
\f{1}{\es} \| \td F \|_{\cYb^{2 k+2}}  \|  \tFm \|_{\cYE^{2 k+2} }^2 
\les  \f{1}{\es} ( \| \tw \|_{\cXb^{2 k + 2} } + \| \tFm\|_{\cYb^{2 k + 2} } )  
  \|  \tFm \|_{\cYE^{2 k + 2} }^2 , \\
  \eal
\eeq
for $\eta = \xeta, \etab$. %
Using the crucial relative smallness estimates \eqref{eq:idea_etab}, we treat it as a perturbation of the coercive terms in \eqref{eq:idea_micro_EE}. The second term $\cN_{m M}$ is estimated similarly.

However, using similar estimates and the top-order bounds \eqref{eq:idea_xeta_H}, \eqref{eq:idea_etab} is not sufficient 
to bound $\cN_{MM}$, since this leads to 
\[
\bal 
  \es^{-1}  \| \tw \|_{\cX_{\eta_1}^{2 k + 2} } 
  \| \tw \|_{\cX_{\eta_2}^{2 k + 2} } 
  \|  \tFm \|_{\cYE^{2 k + 2} } 
& \les \es^{-1}   \|  \tFm \|_{\cYE^{2 k + 2} }^2 
+ \es^{-1}  \| \tw \|_{\cX_{\eta_1}^{2 k + 2} }^2 
  \| \tw \|_{\cX_{\eta_2}^{2 k + 2} }^2  \\
& \les \es^{-1}   \|  \tFm \|_{\cYE^{2 k + 2} }^2 
+ \es^{-2\el}  \| \tw \|_{\cX_{\eta_2}^{2 k + 2} }^2  ,
\eal 
\]
for suitable $\eta_1, \eta_2$. Owing to the large factor $\es^{-2 \ell}$, the term $\es^{-2\el}  \| \tw \|_{\cX_{\eta_2}^{2 k + 2}  }^2$ cannot be absorbed by the $O(1)$ damping term in \eqref{eq:idea_macro_EE:a}. Thus, 
 extra smallness and faster decay estimates are required to treat $\cN_{MM}$ as a perturbation of the damping terms in \eqref{eq:idea_macro_EE}. 
 \footnote{
 This difficulty does not appear in the stability analysis of the global Maxwellian with a \textit{fixed} $\es$, e.g. \cite{guo2002landau,carrapatoso2017nearMaxwellian}. 
 }

To overcome this difficulty, we estimate the perturbation in a lower order 
weighted Sobolev norm \eqref{eq:idea_xeta_L},
which provides an extra smallness $\es^{2/3}$ compared to $\es^{1/2-}$. 
We further establish 
\footnote{
Let $\tFM$ be the macro-perturbation with hydrodynamic fields $\tw$. %
From  Lemma \ref{lem:macro_UPB}, we have the equivalence $\| \tFM \|_{\cYe^k} \asymp_{k, \eta} \| \tw \|_{\cX_{\eta}^k}$.
}
\beq
\bal
\left|
    \f{1}{\es} \la \cN(\tFM, \tFM), \tFm \ra _{\cYe^{{2 k + 2} }}  
\right| \les & \f{1}{\es} 
\| \tFM \|_{\cYY ^{2 k}}
\| \tFM \|_{\cYY ^{2 k + 2}} 
\| \tFm \|_{\cYE^{2 k + 2}} 
\les \f{1}{\es} 
\| \tw \|_{\cXX ^{2 k}}
\| \tFM \|_{\cYY ^{2 k + 2}} 
\| \tFm \|_{\cYE^{2 k + 2} } 
\\
\overset{ \eqref{eq:idea_xeta_L}}{\les} &
\es^{-1/3}  
\| \tFM \|_{\cYY ^{2 k + 2 }} 
\| \tFm \|_{\cYE^{2 k + 2} }
\les \es^{1/6} \left(
    \es^{-1} \| \tFm   \|_{\cYE^{2 k + 2} }^2 
    + \| \tFM \|_{\cYY ^{2 k + 2}}^2 
\right),
\eal 
\label{eq:idea_non}
\eeq
for $\eta = \xeta$ and $\etab$. This allows us to treat the nonlinear term perturbatively. Note that in \eqref{eq:idea_non},
we gain crucial spatial decay so that we can bound $\tFM$ in the spaces $\cYY ^{2k + 2}$ and 
$\cYY ^{2 k}$ with a weaker spatial weight $\la X \ra^{\xeta}$, rather than the stronger weight $\la X \ra^{\etab}$ 
appearing in $\cYb^{2 k + 2}$. Moreover, by interpolation, one of the two $\tFM$ terms can be placed in the lower order norm $\cYY ^{2 k }$,
which satisfies the sharper estimate \eqref{eq:idea_xeta_L}.

\subsubsection*{\underline{Step 6}. Construction of global solutions to~\eqref{eq:lin}}
\label{sec: step6}

We construct a global solution to \eqref{eq:lin}, which satisfies  estimates \eqref{eq:idea_size}, by combining the a priori estimates established in previous steps and using a fixed point argument. To avoid the potential unstable directions and estimate $\tw_2$ in \hyperref[sec: step2]{\itshape \underline{Step 2}}, we generalize the argument in \cite{chen2024Euler,ChenHou2023a} by splitting the equations, applying Duhamel's formula, and backward-in-time semigroup estimates for $\tw_2$; see \eqref{eq:W2_form}. We consider initial data small enough (relative to $ \e_0 = \d$) in the norms in \eqref{eq:idea_size} and  prove \eqref{eq:idea_size} using a bootstrap argument. Define 
\[
  E_{k+1, \xeta} = \kp \| \tw_1 \|_{\cXX^{2 k+2 }}^2 + \| \tFm \|_{\cYY^{2 k+ 2}}^2, 
  \qquad  E_{k+1, \etab} =\kp \| \tw \|_{\cXb^{2 k+ 2}}^2 + \| \tFm \|_{\cYb^{2 k+2 }}^2 .
\]
We combine estimates \eqref{eq:idea_macro_EE:a} and \eqref{eq:idea_micro_EE:a} and estimates in \hyperref[sec: step4]{\itshape \underline{Step 4}}, \hyperref[sec: step5]{\itshape \underline{Step 5}} to obtain 
\beq\label{eq:idea_xeta_EE}
\bal
  \f{1}{2} \f{d}{d s } E_{k+1, \xeta} 
  \leq -\lam_1 E_{k+1, \xeta} - \f{ \cgam }{4 \es} \| \tFm \|_{\cYL{\xeta}^{2 k + 2 }}^2 + C_k \es .
  \eal
\eeq
We combine \eqref{eq:idea_macro_EE:b} and \eqref{eq:idea_micro_EE:b} and estimates in \hyperref[sec: step4]{\itshape \underline{Step 4}}, \hyperref[sec: step5]{\itshape \underline{Step 5}} to obtain 
\[
   \f{1}{2} \f{d}{ds} E_{k+1, \etab} 
 \leq C_k E_{k+1, \xeta} 
 -  \f{ \cgam  }{4 \es} \| \tFm \|_{\cYQ^{ 2 k + 2 }}^2
  +C_k \es  . 
\]
These estimates imply $E_{k+1, \xeta} \les_k \es, E_{k+1, \etab} \les_k   \d^{1 - 2 \el}$ and improve estimate \eqref{eq:idea_xeta_H}. 
\footnote{
We have $\es = \d \operatorname e^{-\rE s }$ \eqref{eq:EE_para1} and choose $\lam_1 > \rE$ \eqref{eq:eta_constraint}.
}

\paragraph{Extra smallness of $\| \tw_1\|_{\cXX^k}$}

With the above top order estimates, we have an improved estimate at lower order. By interpolating 
the two damping terms in \eqref{eq:idea_xeta_EE}, we exploit the large damping term $\es^{-1} \| \tFm \|_{\cYL{\xeta}^{2k+2}} $  and bound the cross term in \eqref{eq:idea_macro_EE:a} with $k_1 = k$: 
\beq\label{eq:idea_xeta_EE_L}
  \f{1}{2} \f{d}{d s} \| \tw_1 \|_{\cXX^{ 2 k} }^2 
   \leq - \lam_1 \| \tw_1 \|_{\cXX^{ 2 k} }^2  
   + C_k \es^{1/2} \| \tw_1 \|_{\cXX^{ 2 k}} \cdot   \B(  \f{1}{\es}  \| \tFm \|_{\cYL{\xeta}^{ 2 k + 2} }^2 \B)^{1/2} 
   + C_k \es^2 ,
\eeq
where the last term $C_k\es^2$ arises from estimating terms introduced by the profile modification in \eqref{eq:Euler_profi_modi} and is negligible compared to other terms. The small factor $\es^{1/2}$ in the above estimates shows that at lower order, the estimates of the micro and macro perturbation are weakly coupled. Combining estimates \eqref{eq:idea_xeta_EE_L}, \eqref{eq:idea_xeta_EE}, and using the extra small factor $\es^{1/2}$, we improve \eqref{eq:idea_xeta_L}.

With these global estimates, we choose the initial perturbation carefully to ensure non-negativity of initial data and prove Theorem \ref{thm:blowup}.

\begin{remark}[\bf Effect of dissipation]

To illustrate the mechanism behind the improved estimates for
$\|\tw_1\|_{\cXX^{2k}}$ in \eqref{eq:idea_xeta_EE_L}, we consider a simplified
model. Approximating $(\cT+\dcm-\tfrac32\bcv)\td F$ by $\td F$, the error term
$-\cMM^{-1/2}\eM$ by $1$, the dissipative operator $\cL_{\cM}\td F$ by 
a damping term $-\td F$,
and neglecting the nonlinear terms $\cN$, equation
\eqref{eqn: linearized-td-F} may be heuristically reduced to
\[
\pa_s \td F + \td F  = - \es^{-1}\,\td F + 1.
\]
For small $\e_0$ and initial data, the dissipation $-\es^{-1}\td F$ leads to 
$|\td F(s)|\lesssim \es$, which is much smaller than the scale
$\es^{1/2-\ell}$ in \eqref{eq:idea_size}. This suggests that exploiting the
dissipation yields sharper estimates. At the top-order level, however, this
mechanism cannot be fully exploited, since closing the estimates
\eqref{eq:idea_macro_EE},\eqref{eq:idea_micro_EE},\eqref{eq:idea_xeta_EE} requires the specific coupled structure between $\tFm$
and $\tFM$ in \hyperref[sec: step4]{\itshape \underline{Step 4}}. Instead, we exploit the dissipative effect at
the lower-order level in \eqref{eq:idea_xeta_EE_L} to establish
\eqref{eq:idea_xeta_L}.

\end{remark}

\subsubsection{Choice of the parameters}
Below, we discuss several parameters in the energy estimates. %
\paragraph{Parameters of decay and weights}

We discuss the constraint on the exponent $\xeta$ and choose $\lame, \lam_1$, which appeared in the above steps.  Recall the definition of $\rE$ from \eqref{eq:hydro2}, $\bcx$ from \eqref{eq:Euler_profi},
and $\etab$ from \eqref{wg:X_power}. We choose $\xeta$ 
with the following properties 
\bseq\label{eq:eta_constraint}
\begin{align}
\f{\bcx}{4} (\etab - \xeta)  & > \rE > 0,  \label{eq:eta_constraint:a} \\ 
  \etab - \xeta & \leq \f{(1 + \om) r}{2} . 
  \label{eq:eta_constraint:b} 
\end{align}
From Remark \ref{lem:para_bd}, it is not difficult to see that the constraints for $\xeta$ form a non-empty interval.

We define its related decay exponents $\lame$ and choose  $\lam_1$ close to $\lame$ such that 
\beq 
    \lame := \f{\bcx}{4} (\etab - \xeta), \quad   \rE  <  \lam_1   < \lame. 
     \label{eq:eta_constraint:c}
\eeq 
\eseq 

The factors $ \lam_1$ and $\lame$ are related to the spectral gap in the linear stability estimate, for the $\cYY$ and $\cXX$ norms, respectively;~see \eqref{eq:idea_macro_EE}, \eqref{eq:idea_micro_EE}, \eqref{eq:idea_xeta_EE}, and Theorem~\ref{thm:coer_est}.  We impose the lower bound on $\etab - \xeta$ in \eqref{eq:eta_constraint:a} so that 
the linear damping terms, e.g. $-\lam_1 E_{k+1,\xeta}$ lead to decay faster than the error $\es$ \eqref{eq:idea_xeta_EE}. 
\footnote{
If \eqref{eq:eta_constraint:a} does not hold and $\lam_1 < (1/2 - \el) \rE$, estimate \eqref{eq:idea_xeta_EE} implies decay estimate $E_{k+1,\xeta}^{1/2}(s) \les E_{k+1,\xeta}^{1/2}(0) \cdot \operatorname e^{-  \lam_1 s }$, which is slower than $\es^{1/2 -\ell} = C(\e_0) \operatorname e^{-(1/2-\ell) \rE s}$ \eqref{eq:idea_xeta_H}. $\ell$ can be  essentially treated as a small parameter close to $0$.  
}

We impose the upper bounds on $\etab - \xeta$ in \eqref{eq:eta_constraint:b} so that the $\cYY$-norm is not too weak compared to the $\cYb$-norm. This constraint comes from estimate \eqref{eq:idea_non}, where we bound the $\cYb$-estimate of the nonlinear term using the weaker $\cYY$-norm. 
See more details in \eqref{eqn: etaD-constraint2} and Theorem \ref{theo: N nonlin_est}.

\paragraph{Parameters related to $\es, \rs$}

Recall $\rE, \bcv, \bcx$ from \eqref{eq:Euler_profi} and \eqref{eq:hydro2}. 
For some small $\d \in (0,1)$ to be chosen in Theorem \ref{thm:non}, we choose $R_0$ in \eqref{eq:S_radial} 
and $\e_0$ in \eqref{eq:knud} as
\bseq\label{eq:EE_para1}
\beq
  \e_0 = \d , \quad R_0 = \e_0^{ - \ell_r }
  = \d^{-  \ell_r  } , \quad \ell_r = \f{\bcx}{ \rE }.
\eeq
Since $\rs = R_0 \operatorname e^{\bcx s} 
= \e_0^{-\ell_r} \operatorname e ^{\ell_r \rE s}$ and $\es = \e_0 \operatorname e^{-\rE s}$, we get
 \beq
 \es = \d \operatorname e^{ -\rE s} \leq \d \leq 1,  \quad \rs = \es^{ - \ell_r } .
 \eeq
\eseq
for any $s \geq 0$. 

\begin{remark}[Range of parameters]\label{lem:para_bd}
Let $\rE$ be defined in \eqref{eq:hydro2}. For $\g \in (\sqrt{3}, 2]$ and $ r \in  ( r_*- 0.01 , r_*)$ with $r_*  = 3 - \sqrt{3}$, we have 
the following inequality regarding these parameters
\[
 1.25< r < 1.3,  \quad 0 < \rE \leq 5 \cdot \f{2 - \sqrt{3}}{3 - \sqrt{3}} -1 < 0.06 .
\]
As a result
\beq\label{eq:para_bd}
  r \ell_r = \f{1}{\rE} > 2 , \quad \ell_r = \f{1}{r \rE} > 2,
  \quad \rs^{-r} = \es^{ r \ell_r}  \les \es^2,
  \quad \rs \gtr \es^{-2}.
\eeq
\end{remark}

\paragraph{Other parameters}

We have fixed $\g$ for the Landau equation \eqref{eq:Q}, fixed the exponent $r$ for the profile in Section \ref{sec:isen_euler}, 
and determined parameters $\etab, \xeta,\lame,\lam_1$ in \eqref{eq:eta_constraint}. 
We fix the parameter $\ell$ in \eqref{eq:idea_size} to be a small absolute constant.

Our stability estimates involve a few more parameters: $\d$ for $\e_0$ and the size of perturbation \eqref{eq:idea_size}, $\nu$ in the $\cY$-norm \eqref{norm:Y}, and the order of energy estimates $k$ (see \hyperref[sec: step1]{\itshape \underline{Steps 1-6}}). We determine these parameters sequentially:
\[
  k = \kk  \rsa  \nu \rsa \d ,
\]
where each later parameter may depend on the previous ones. We determine $\kk$ in \eqref{def:kk},
 $\nu$ in Theorem \ref{thm: micro-Hk-main}, and $\d$ in Theorem \ref{thm:non}.

\subsubsection{Comparison with Guo's stability estimates in \cite{guo2002landau}}\label{sec:compare_guo}

Part of our stability estimates for the micro-perturbation build on those in \cite{guo2002landau}. 
For instance, we adopt the coercivity estimates for the linearized Landau collision operator and the associated functional framework from
\cite{guo2002landau} in the stability estimates in \hyperref[sec: step3]{\itshape \underline{Step 3}} in Section \ref{sec:idea} for the micro-perturbation in
$V$ at each fixed point $X$; we need to generalize the $\s$-norm introduced in~\cite{guo2002landau} 
to define a $\s$-norm associated with a local Maxwellian  in \eqref{norm:sig}. 
In addition, some of our nonlinear estimates for the collision operator, such as \eqref{eq:idea_decomp_N2}
in \hyperref[sec: step5]{\itshape \underline{Step 5}}, are inspired by those in \cite{guo2002landau}.

Despite these similarities, there are several essential differences between \cite{guo2002landau} and the present work. 
First, our profile \eqref{eqn: leading-order-at-blowup} is a \emph{local} Maxwellian
rather than a global one \cite{guo2002landau}, which necessitates the development of genuinely inhomogeneous estimates.
In particular, we design weighted operators and functional spaces with $X$-weights that depend on the stability estimates for the macro-perturbation and are adapted to the self-similar scaling fields. See Section \ref{sec:function_space}.
Second, we perform stability analysis for $X$ in the whole space, rather than on the torus $ X \in\mathbb{T}^3$ \cite{guo2002landau}. Since coefficients in the coercivity estimates decay in $X$, rather than remaining uniformly bounded away from zero \cite{guo2002landau}, we need to carefully control the spatial decay of the perturbation for large $X$.
We emphasize that the \emph{temporal} decay estimates for the perturbation are \emph{sensitive} to the choice of weights used for the \emph{spatial} decay estimates, as reflected in \eqref{eq:idea_size}, making the control of spatial decay one of the major difficulties.
See the paragraph \emph{Relative smallness estimates} in \hyperref[sec: step1]{\itshape \underline{Step 1}} for further discussion of this difficulty.

Our analysis involves further challenges, including stability estimates for the macro-perturbation, the limit $\es \to 0$ in the self-similar equation, and the construction of a blowup solution from a finite codimension set.

\subsection{Notation}

We collect here the main notation used throughout the paper.  For each variable or operator, we list below it the equation or result in which it is first defined or determined.

We use lowercase letters $f,x,v$ to denote variables in the physical equations, whereas uppercase letters $F,X,V$ denote variables in the self-similar equations. The time variables are denoted by $t$ in the physical equations and by $s$ in the self-similar equations. Lowercase letter $m$ indicate microscopic variables or operators, such as $\tFm, \cP_m$, 
while uppercase letter $M$ indicate macroscopic variables or operators, such as $\tFM, \cP_M$. 

\paragraph{Operators}
We use calligraphic font to denote operators. 
Calligraphic $\cL$ is reserved for linearized operators around the profile, such as 
\[
    \uds{ \eqref{def:T-LM}}{ \cL_\cM}, \quad 
    \udb{ \cL_{E} ,\cL_U, \cL_P, \cL_B }_{ \eqref{eq:lin_euler_limit} }  , \quad
    \udb{ \cL_{E,s} , \cL_{U, s}, \cL_{P, s}, \cL_{B, s} }_{  \eqref{eq:lin_euler} } , \quad  
    \uds{ \eqref{eqn: defn-Lmic}  }{ \Lmic }  .
\]
The following operators are introduced in the linearization in Section \ref{sec:lin_decomp}
and in Section \ref{sec:non_collision_op}
\[
    \uds{\eqref{def:T-LM}}{\cT}, \quad
    \uds{\eqref{eq:non_nota}}{\cN}, \quad
    \uds{\eqref{eq:moments_Fm}}{\cI},
    \quad \uds{ \eqref{N(f,g)} }{ \cN_i}.
\]

We use calligraphic $\cF$ to denote the maps between the macro-perturbation and the variables in the Euler equations :
$
    \underset{\eqref{eq:pertb_macro}}{\cF _E}, \quad
    \underset{\eqref{eq:macro_UPB}}{\cF _M} $. 
    
We use $\cK$ to denote a compact operator defined in Proposition \ref{prop:compact} :
$
    \uds{ \mathrm{Proposition \ } \ref{prop:compact} }{ \cK_{k, \eta} },
    \quad  \uds{ \eqref{def:kk} }{ \cK_{\kk}} . 
$

We use calligraphic $\cP$ to denote projection : $\uds{\eqref{def:proj}}{\cP_m},  \quad \uds{\eqref{def:proj}}{\cP_M} .$

We use $\Pi_{\cdot}$ to denote various projections 
\[
  \uds{  \eqref{eqn: defn-bPhi-A}   }{\bpi_v }, 
  \quad \uds{  \eqref{eqn: defn-matrix} }{ \bpi_{\vc } }, 
  \quad   \uds{   \eqref{eq:W2_form}  }{ \Pi_{ \mathsf{s} } },  \quad  \uds{ \eqref{eq:W2_form} }{\Pi_{\mathsf{u}} } .
\]

\paragraph{Functions and parameters}

We use $F$-functions  $\td F, \tFm, \tFM$ to denote functions related to the solution of the self-similar Landau equation \eqref{eq:LC_ss}.

We introduce the following radial variable, unit vector, and velocity field:
\begin{align}
		\xi &= |X|, &
        \ee _R &= \f{X}{|X|}, &
       \uds{ \eqref{eq:Euler_profi}}{ \bu (X) = \bar U (\xi) \ee _R . }
               \label{eqn: radial}
\end{align}

We use the following parameters related to the profiles,
\[
  \uds{  \eqref{eqn: defn-bPhi-A} }{ \gamma },  \quad 
  \uds{   \eqref{eq:Euler_profi}   }{r}, 
  \quad \uds{  \eqref{eqn: defn-omega}   }{ \rE(r, \gamma)},
  \quad 
  \uds{ \eqref{eq:EE_para1} }{ \e_0 , \ R_0 , \ \ell_r },
\]
 and parameters related to weights and estimates 
\[
  \uds{ \eqref{wg:X_power} }{ \xeta, \ \etab},
  \quad \uds{  \eqref{eq:idea_xeta_H}  }{\ell = 10^{-4}  }, 
  \quad  \uds{ \eqref{def:kk} }{\kkl , \ \kk } ,
  \quad  \uds{  \mathrm{Theorem  \ } \ref{thm: micro-Hk-main}  }{ \nu},
  \quad \uds{ \mathrm{Theorem \ } \ref{thm:non},  \ \ref{thm:blowup}  }{\d} .
\]

We use a ``bar" notation $\bar \cdot$ to denote constants and functions associated with the profile: 
\[
\udb{  
    \bar c _f,\quad
    \bcx,\quad
    \bcv
}_{\eqref{eq:Euler_profi_cfcvcx}},\quad 
\udb{
    \bar \rho,\quad  
    \bu,\quad   
    \bth,\quad  
    \bc,\quad 
    \bar B,\quad 
    \bar P    
}_{\eqref{eq:Euler_profi_rhoPth}} ,
\]
and a ``tilde" notation $\td \cdot$ to denote perturbation variables: 
\[
    \uds{\eqref{eq:pertb_dec}}{\td F},\quad 
    \udb{\tFm, \ \tFM}_{\eqref{def:proj}},\quad 
    \udb{ \td {\rho}, \tu, \tb, \tp }_{\eqref{eq:pertb_macro}},\quad
    \tw = (\tu, \tp, \tb) .
\]

We use subscript  $\bar \cdot _s$ to denote variables with the cutoff-modification introduced in \eqref{eq:Euler_profi_modi} 
\[
    \udb{ \rhos , \quad  \ths , \quad  \bar P_s,  \quad  \cs }_{ \eqref{eq:Euler_profi_modi} }
\]
and to denote variables and operators depending on the self-similar time $s$
\[
  \uds{ \eqref{eq:knud} }{ \es } , \quad 
\uds{  \eqref{eq:S_radial}   }{ \rs} , \quad 
    \udb{ \cL_{E,s} , \  \cL_{U, s},  \ \cL_{P, s},  \  \cL_{B, s} }_{  \eqref{eq:lin_euler} }
    \quad  
\]

We use  calligraphic  $\cM$ to denote a local Maxwellian and $\mu$ to denote the Gaussian function:
\[
\uds{ \eqref{eq:local_max} }{ \cM_{\rho, \UU, \Th} } ,  \quad   \uds{ \eqref{eq:localmax2} }{  \cM,  \  \cMM }  ,
\quad \uds{ \eqref{eq:gauss}  }{ \mu(\cdot)} ,
\]
The calligraphic  $\cE$ is reserved for variables related to errors %
\[
    \udb{ \eM, \   \erho, \  \eu, \  \ep,  \  \ec  }_{ \eqref{eq:error0} } .
\]

We use $\lam$-parameters and $\Lam$ 
\[
    \uds{  \eqref{eq:eta_constraint:c} }{ \lam_1, \ \lame }  , 
    \quad  \uds{ \eqref{eq:lam_eta} }{\lam_{\eta}} , 
    \quad  \uds{ \eqref{eq:decay_para} }{ \lam_s, \ \lam_u } ,
    \quad  \uds{ \eqref{eq:coer_coe} }{\Lam(s, X, V) } 
\]
to denote parameters or functions related to the decay rates and coercivity estimates. 

We use $\varpi$-parameters to denote parameters in the norms, e.g. $\cX$-norm \eqref{norm:Xk} and $Z^j_R$-norm \eqref{norm:cZZ}
\[
  \uds{ \eqref{norm:Xk0} }{ \varpi_{k, \eta } } ,   \quad 
  \uds{ \eqref{norm:Z}  }{\varpi_k^{\prime} }  , \quad  
  \uds{  \eqref{norm:cZZ} }{\varpi_{Z, i}} .
\]

The $\s$, $\cX$, $\cY$ norms are defined in Section \ref{sec:function_space}. 
We define the $\cZ$-norm in \eqref{norm:Z}, the $Y$-space for the fixed point argument 
in \eqref{norm:fix}, and the $Z^j_R$-norm in \eqref{norm:cZZ}.

\paragraph{Symbols}

Angled brackets represent the Japanese bracket $\ang \cdot$ or an inner product or duality pairing $\ang{\cdot, \cdot}$ depending on the context. In particular $\ang{\cdot, \cdot} _V$ is a duality pairing in the $V$ variable defined in $\eqref{eq:inner_V}$.

We write  $p \les q$ to mean that there exists some absolute constant $C > 0$ such that
	$p \leq C q$, and  $p \asymp q$ to mean that $p \les q$ and $q \les p$. We use the notation $A=  B + O_{h}(B^\prime)$ 
    to indicate that there exists $C_{h} > 0$ such that $|A - B| \leq {C_h} B^\prime$. 
 In particular, $A = O_{h}(B^\prime)$ means that $|A| \leq C_h B^{\prime}$  for some constant $C_h>0$ depending on $h$. Throughout the paper, $c,C$ and $C_i$ denote absolute constants that may vary from
line to line, while $\bar C_h$ (with a bar) denotes a fixed constant depending on
$h$. We use the following fixed constants in this paper 
\[
\uds{ \eqref{eq:Euler_var}  }{ \kappa  = \f53 },
\quad \uds{ \eqref{eq:gauss} }{  \kp_2 = \f56 },  
\quad 
  \uds{ \mathrm{Lemma \ } \ref{lem: spectral-gap}  }{ \cgam }  ,
  \quad \uds{ \mathrm{ Theorem \ } \ref{thm:coer_est} } { \bar C_{k, \eta} }, 
  \quad \uds{ \mathrm{Theorem \ } \ref{theo: N nonlin_est}   }{  \bar C_{ \cN }  } .
\]
For any multi-index $\al, \b \in \Znonneg^3$, we write $ \al \preceq \b$ if and only if $\al_i \leq \b_i$. We write $\al \prec \b$ if $\al \preceq \b$ and $\al \neq \b$.

\section{Properties of Euler profile and equations of macro-perturbation}\label{sec:Euler_setup}

In this section, we present the properties of the Euler profile and its modification which are used to 
established finite codimension stability of the macro-perturbation, and derive the equations of macro-perturbation, which is the linearized Euler equations.

\subsection{Properties of the Euler profile}

For the Euler profile $(\bar \UU, \bar{\sc})$ without modification \eqref{eq:Euler_profi},  we first recall the following properties from \cite[Theorem 1.1]{shao2025blow}. 
\begin{lemma}\label{lem:profile0}
 The profile $ \bar \UU = \bar U \ee _R, \bar{\sc} $ are radially symmetric and satisfy  $\bar U(0) = 0$ and
\bseq\label{eq:euler_prop0}
\beq\label{eq:dec_U}
  \bal
  |\na^k \bar \UU | \les_k \la X \ra^{-r +1 - k}, \quad 
   \bar{\sc}  \asymp \la X \ra^{-r+1} , 
  \quad  
    |\na^k \bar{\sc} | \les_k \la X \ra^{-r + 1- k},
  \eal
\eeq
for any $k \geq 0$. There exists $c_0>0$ such that for any $\xi \geq 0$, we have
\begin{align}
\bcx + \partial_\xi \bar{U}(\xi) - |\pa_{\xi} \bar{\sc}(\xi)| &> c_0 .   \label{eq:rep10} 
\end{align}
\eseq
\end{lemma}

Let $\xi_*$ be the unique root of
\beq\label{eq:sonic_pt}
\bcx \xi_* + \bar U(\xi_*) - \bar \sc(\xi_*) = 0.
\eeq
It corresponds to the degenerate point in the phase portrait \cite{merle2022implosion1,buckmaster2022smooth,shao2025blow} and is called the sonic point.

We have the following estimates for the modified profile $ \cs$ \eqref{eq:Euler_profi_modi}. 
\begin{lemma}\label{lem:profile}
The modified profile $\cs$ satisfies 
\bseq
\beq\label{eq:dec_S}
\bal
   \cs  \asymp \la X \ra^{-r+1} + \rs^{-r+1},  \quad 
    |\na^k \cs | \les_k \la X \ra^{-r + 1- k} \les \cs \la X \ra^{- k},
  \eal 
\eeq
for any $k \geq 1$. There exists $c_1 > 0$, $\xi_1 > \xis$, and $R_{0, 1} \gg 1$, such that for any $\xi \geq 0$ and $R_0 \geq R_{0, 1}$, we have 
\footnote{
Note that $\cs$ depend on the parameter $R_0$ via $\rs$. See \eqref{eq:S_radial} and \eqref{eq:Euler_profi_modi}.
}
\begin{align}
\bcx + \partial_\xi \bar{U}(\xi) - |\pa_{\xi} \cs(\xi)| &> c_1,  \quad \xi \in [0, \xi_1] \label{eq:rep1}, \\
\bcx \xi + \bar U(\xi) - \cs(\xi)&> \min \left\{ \bcx  \xi + \bar U - \bc ,  \, \f{1}{2}  \bcx \xi \right\} > 0, \quad  \forall \, \xi > \xis, \label{eq:rep2} \\
\bcx + \xi^{-1} \bar U(\xi)  &> c_1, \quad \xi \geq 0 ,  \label{eq:rep22} 
\end{align}
with implicit constants independent of $s, R_0$ in the definition of $\cs$ \eqref{eq:Euler_profi_modi}, \eqref{eq:S_radial}.
\eseq

 Moreover, there exists constants $ C_{\bu}$ and $ C_{\bc} >  0$ such that for any $|X| \geq 1$, we have the refine asymptotics 
\beq\label{eq:UC_refine_asym}
   | \bc(X) - C_{\bc} |X|^{-(r-1)} | \les |X|^{ - 2 r + 1},
   \quad |\bar U(\xi) - C_{\bu}  |X|^{- (r-1)}   | \les |X|^{-2 r + 1}. 
\eeq

\end{lemma}

\begin{remark}[Repulsive conditions]

As in \cite{chen2024Euler,chen2024vorticity}, to establish finite codimension stability estimates, we only need the interior repulsive condition \eqref{eq:rep1}, which follows from \eqref{eq:rep10} for $\xi \in [0, \xi_*]$.  Establishing the exterior repulsive condition (\eqref{eq:rep10} with $\xi > \xi_*$) can be highly nontrivial, see e.g. \cite{merle2022implosion1}. While we use the full repulsive condition \eqref{eq:rep10} to prove 
the outgoing conditions \eqref{eq:rep2} and \eqref{eq:rep22} below, 
these two conditions follow from the natural barrier functions 
and the sign of the denominator in the ODEs for profile $(\bu, \bc)$, %
and they are much simpler to establish than the exterior repulsive condition. See further discussion in \cite[Remark 2.3]{chen2024Euler}.
\end{remark}

\begin{proof}[Proof of Lemma \ref{lem:profile}]
On the one hand, from the definition of $\cs$ in \eqref{eq:Euler_profi_modi} and the upper bound of $\bar\sc$ in \eqref{eq:dec_U}, we see directly that $\cs \les \bar\sc + \rs ^{-r + 1} \les \ang X ^{-r + 1} + \rs ^{-r + 1}$. On the other hand, $\bar\sc \les \ang X ^{-r + 1} \le \rs ^{-r + 1}$ for $|X| \ge \rs$, together with $\chi _\rs \equiv 1$ for $|X| \le \rs$ we see 
$$\bar\sc \les \bar\sc \chi _{\rs} + \rs ^{-r + 1} (1 - \chi _\rs) = \cs.$$
Similarly, from $\rs ^{-r + 1} \les \ang X ^{-r + 1} \asymp \bar\sc$ when $|X| \le 2 \rs$ and $1 - \chi _{\rs} \equiv 1$ for $|X| \ge 2 \rs$, we know $\rs ^{-r + 1} \les \cs$. Combined, %
we prove the first comparison in \eqref{eq:dec_S}.

By definition of $\cs$ and $\chi_{\rs}$ in \eqref{eq:Euler_profi_modi}, we obtain 
\begin{align*}
    |\na^i \chi_{\rs}| \les_i \one _{\{ \rs \le |X| \le 2 \rs \}} \rs ^{-i}, & \qquad \forall i \geq 1.
\end{align*}
We apply this to the second estimate in \eqref{eq:dec_S}, using \eqref{eq:dec_U} and Leibniz rule for $k \ge 1$: 
\[
\bal
	|\na^k \cs| & \les_k \sum_{ i = 0} ^k |\na^i \bc|  \cdot |\na^{k-i} \chi_{\rs}|
	+ \rs^{-r+1} |\na^k (1 - \chi_{\rs})| \\
	& \les_k |\na ^k \bc| \cdot \chi _{\rs} + \sum_{i = 1} ^k \ang X^{-r + 1 - i} \rs^{-(k-i)} \one _{\{\rs \le |X| \le 2 \rs\}}
	+ \rs^{-r + 1 - k} \one _{\{\rs \le |X| \le 2 \rs\}} \\
    & \les_k \ang X^{-r+1 - k} 
    \les \cs \ang X^{-k}.
\eal
\]

\paragraph{Proof of \eqref{eq:rep1}-\eqref{eq:rep22}}

Firstly, we take $0< c_1 < c_0$,  $\xi_1 > \xis$, and $R_{0, 1} > \xi_1 $. Then we get $\cs = \bc$ for $|X| \leq \xi_1 < R_0 < \rs$. Thus, \eqref{eq:rep1} follows from \eqref{eq:rep10}. 

Since $\bcx \xi + \bar U |_{\xi =0} = 0$ 
and $\bcx \xi + \bar U - \bc |_{\xi = \xis} = 0$, using \eqref{eq:rep10} and integration, we obtain 
\bseq\label{eq:outgo_pf}
\begin{align}
	\bcx \xi + \bar U(\xi) & \geq c_0 \xi  > c_1 \xi, \quad \forall \ \xi \geq 0, 
\label{eq:outgo_pf:a} \\
		\bcx \xi + \bar U(\xi) -\bc(\xi) & \geq c_0 (\xi-\xi_*)  > 0, \quad \forall \ \xi > \xi_* .
\label{eq:outgo_pf:b}		 
\end{align}
\eseq

Estimate \eqref{eq:outgo_pf:a} implies \eqref{eq:rep22}. 

From \eqref{eq:dec_S} and \eqref{eq:dec_U}, 
for $R_0$ sufficiently large and $\xi =|X| \geq \rs \geq R_0$, we obtain 
\beq\label{eq:outgo_pf2}
	\bcx \xi + \bar U - \cs 
	\geq \bcx \xi - C \ang \xi^{-r + 1} - C \rs^{-r+1} \geq 
 \f{1}{2} \bcx \xi  + \f{1}{4} \bcx \rs -C \rs^{-r+1} \geq  \f{1}{2} \bcx \xi.
\eeq

Since $	\bcx \xi + \bar U - \cs  = 	\bcx \xi + \bar U - \bc$ for $\xi \leq \rs$, 
combining \eqref{eq:outgo_pf:b} and \eqref{eq:outgo_pf2}, we prove \eqref{eq:rep2}.

\paragraph{Proof of \eqref{eq:UC_refine_asym}}

Recall that the profile $(\bu, \bc , \bcx, \bcv)$ solves the steady state of \eqref{eq:Euler3}
\beq\label{eq:profile_eqn}
 [ (\bcx X + \bu) \cdot \na] \bc + \f13 \bc (\na \cdot \bu)   = \bcv \bc, \quad 
 [ ( \bcx X + \bu)  \cdot \na] \bu + 3 \bc \na \bc   = \bcv \bu.
 \eeq

Since $\bu(X) =\bar U(\xi) \ee _R$ is radially symmetric, we get
\[
\bcx X \cdot \na \bar U(\xi) -   \bcv \bar U   = \bcx \xi \pa_{\xi} \bar U(\xi) -   \bcv \bar U  = - 3 \bc \pa_{\xi} \bc 
- \bar U(\xi) \pa_{\xi} \bar U(\xi).
\]
Since $ \f{\bcv}{\bcx} = 1 -r$, using the integrating factor $\xi^{r-1} = |X|^{r-1}$
and then dividing $\xi$, we get
\[
 \bcx  \pa_{\xi} ( \xi^{r-1} \bar U(\xi) ) = - \xi^{r-2} ( 3 \bc \pa_{\xi} \bc 
+ \bar U(\xi) \pa_{\xi} \bar U(\xi) ) . 
\]

Using the decay estimates \eqref{eq:dec_U}, we obtain $| \xi^{r-2} ( 3 \bc \pa_{\xi} \bc 
+ \bar U(\xi) \pa_{\xi} \bar U(\xi) )| \les \xi^{-r-1} $, which is $L^1$-integrable in $\xi$. Thus, there exists
$C_{\bu}$ such that 
\[
 |\xi^{r-1} \bar U(\xi) - C_{\bu}| \les \int_{\xi}^{\infty} |\xi^{r-2} ( 3 \bc \pa_{\xi} \bc 
+ \bar U(\xi) \pa_{\xi} \bar U(\xi) )   | \les \xi^{-r}. 
\]
Dividing $\xi^{r-1}$ on both sides, we prove the asymptotics of $\bar U(\xi)$ in \eqref{eq:UC_refine_asym}. 
The asymptotics of $\bc(X)$ in \eqref{eq:UC_refine_asym} is proved similarly.

Since $\bc(X) \gtr \ang X^{-r+1}$ \eqref{eq:dec_S}, we obtain $ C_{\bc} > 0$, where $C_{\bc}$ is the coefficient in \eqref{eq:UC_refine_asym}. 
\end{proof}

\subsection{Linearized Euler equations}\label{sec:lin_euler}

To control the macroscopic terms, we introduce the weighted hydrodynamic fields 
\bseq\label{eq:pertb_macro}
\beq
\bal
  ( \td  \rho, \td \UU, \td P) & := \int  \cMM^{1/2} \td F \cdot \B(  1,  \f{ V - \bu}{\cs}, \f{|V -\bu|^2}{3\cs^2}  \B) d V, \quad  \tB  := \trho - \tp,
\eal
\eeq
and encode the above linear map from $\td F$ to $\tw := (\tu, \tp, \tb)$ as
\beq
  (\tu, \tp, \tb) := \cF_E( \td F ) 
  = \int \cMM^{1/2} \td F \cdot \B( \f{ V - \bu}{\cs}, \f{|V -\bu|^2}{3\cs^2} , 1 - \f{|V -\bu|^2}{3\cs^2}  \B) d V , 
\eeq
\eseq
where $E$ is short for \textit{Euler}. Variables $\tu, \tB$ are similar to the perturbation of the velocity and entropy up to some weights in $(X, s)$. 

Integrating \eqref{eq:lin} against $1, V - \bar \UU, |V - \bar \UU|^2$, we obtain the equations of  
$ \td \rho$
\bseq\label{eq:lin_euler}
\beq
  \pa_s \trho + (\bcx X + \bu ) \cdot \na \trho 
+ \na \cdot ( \cs \tu) = (3 \bcv - \na \cdot \bu) \trho - \cs ^3 \erho,
\eeq
and of $ (\td \UU, \td P, \tb)$
\beq
\bal
  \pa_s \tu & = \cL_{U,s}(\tu, \tp, \tb) - \cI_1(\tFm) - \cs^3  \eu ,  \\
  \pa_s \tp & = \cL_{P,s}(\tu, \tp, \tb)   - \cI_2(\tFm) -\cs^3 \ep ,\\
    \pa_s \tb & =  \cL_{B,s}(\tu, \tp, \tb)     +      \cI_2(\tFm) , \\
\eal
\eeq
where the linearized operators are defined as
\beq
\bal
\cL_{U,s} \tw &:= 
- ( \bcx X + \bu) \cdot \na \tu - \cs \na \tp   
+ \B( 3 \bcv - \f{2}{3} \na \cdot \bu - (\na \bu) - \ec \B) \tu \\
& \qquad - 2 \na \cs \cdot \tp + 3 \cs^{-1} \bar \sc \na \bar \sc (\tp + \tb) ,  \\
\cL_{P, s} \tw & := - (\bcx X + \bu ) \cdot \na \tp - \cs \na \cdot \tu   
+  \B( 3 \bcv -  \na \cdot \bu
- 2 \ec
\B) \tp -  \B( \na \cs +  \f{2}{3} \eu \B) \cdot \tu, \\
\cL_{B, s} \tw &:= - (\bcx X + \bu) \cdot \na \tB
+ (3 \bcv - \na \cdot \bu) \tB
+ 2\ec \tp + \f{2}{3} \eu \cdot \tu  ,
\eal
\eeq
\eseq
the matrix $ \na \tu $ is given by $(\na \tu)_{ij} = \pa_j \tu_i $, and $\cI_i$ depends on the micro part and is defined in \eqref{eq:moments_Fm}. We refer the derivation to Appendix \ref{sec: derivation-euler}. Note that the projection of \eqref{eq:lin} onto the hydrodynamic fields give the \textit{full} linearized Euler equations around the isentropic profile and we \textit{do not} have nonlinear terms. 
Denote $\cL_{E,s} = (\cL_{U,s}, \cL_{P,s}, \cL_{B,s})$.
Here, the subindex $s$ indicates that the operator $\cL_{E,s}$ is time-dependent. 

As $s \to \infty$, the error $ \ec,\eu$ defined in \eqref{eq:error0} becomes $0$
and $(\cs, \brho, \bp) \to (\bar \sc, \bar \rho ,\bar P)$.  Denote by 
\bseq\label{eq:lin_euler_limit}
\beq
\cL_E = \cL_{E,\infty}, \quad  (\cL_U, \cL_P, \cL_B) = (\cL_{U,\infty}, \cL_{P,\infty}, \cL_{B,\infty} )
\eeq
the limiting operator as $s \to \infty$. We have 
\beq
\bal
\cL_{U} \tw &:= -  ( \bcx X + \bu) \cdot \na \tu - \bar \sc \na \tp   
+ \B( 3 \bcv - \f{2}{3} \na \cdot \bu - (\na \bu)  \B) \tu   - 2 \na \bar \sc \cdot \tp   +  3  \na \bar \sc (\tp + \tb) ,  \\
\cL_{P} \tw & := - (\bcx X + \bu ) \cdot \na \tp - \bar \sc \na \cdot \tu   
+ \B( 3 \bcv -  \na \cdot \bu \B) \tp -   \na \bar \sc  \cdot \tu, \\
\cL_{B} \tw & := - (\bcx X + \bu) \cdot \na \tB
+ (3 \bcv - \na \cdot \bu) \tB .
\eal
\eeq
\eseq

In the rest of the work, we estimate the system of $(\tu, \tp, \tb)$ instead of $(\trho, \tu, \tp)$,
as the former is a symmetric hyperbolic system. 
We can rewrite \eqref{eq:lin_euler} schematically as 
\beq\label{eq:lin_euler2}
\bal
  \pa_s ( \tu, \tp, \tb) &= \cL_{E, s} (  \tu , \tp, \tb)  -  (\cI_1, \cI_2, - \cI_2)(\tFm)
   - ( \cs^3 \eu, \cs^3  \ep, 0), \\
   \cL_{E,s} & = (\cL_{U, s}, \cL_{P, s}, \cL_{B, s}) , 
\eal
\eeq

\subsection{Relations between \texorpdfstring{$\tFM$ and $(\tu, \tp, \tb)$}{FM and U, P, B}}
\label{sec:macro_UPB}

In this section, we derive the relations among $\tFM, \tFm$ \eqref{def:proj}, $\cI_i$ defined in \eqref{eq:moments_Fm}, and $(\tu, \tp, \tb, \trho)$ defined in \eqref{eq:pertb_macro}.

Recall $\ths = \f{1}{\kp} \cs^2 = \f{3}{5} \cs^2$ from \eqref{eq:Euler_profi_modi}. For any functions $G$, we have 
\beq\label{eq:mom_G}
\bal
  \la G , \Phi_0 \ra_V & = \la \cMM^{1/2} G, 1 \ra_\lvv, \\
    \la G, \Phi_i \ra_V
    & =     \B\la \cMM^{1/2} G, \f{ V - \bu}{ \ths^{1/2}} \B\ra_\lvv
    =  \f{\cs}{\ths^{1/2}} \cdot \left\la \cMM^{1/2} G, \f{ V - \bu}{ \cs } \right\ra_\lvv
    = \kp^{1/2} \left\la \cMM^{1/2} G, \f{ V - \bu}{ \cs } \right\ra_\lvv , \ i = 1,2,3 , \\
    \la  G, \Phi_4 \ra_\lvv
      & = \f{1}{\sqrt 6} \B\la \cMM^{1/2} G, \f{ |V - \bu|^2 }{ \ths }- 3 \B\ra_\lvv
      =  \f{1}{\sqrt 6} \B\la \cMM^{1/2} G, 3 \kp \cdot \f{ |V - \bu|^2 }{ 3 \cs^2 }- 3 \B\ra_\lvv.
  \eal
\eeq

Applying the above identities with $ G = \tFM$, using the definition \eqref{eq:pertb_macro} and $\kp = \f{5}{3}$, we get 
\beq\label{eq:mom_comp_FM}
\bal
  \la  \tFM, \Phi_0 \ra_\lvv &= \la \cMM^{1/2} \tFM, 1 \ra_\lvv =  \trho = \tp + \tb,  \\
    \la  \tFM, \Phi_i \ra_\lvv & = \kp^{1/2} \tu_i , \ i = 1,2,3 , \\
    \la    \tFM, \Phi_4 \ra_\lvv 
      &        =  \f{1}{\sqrt 6} \B\la \cMM^{1/2} \tFM, 5  \cdot \f{ |V - \bu|^2 }{ 3 \cs^2 }- 3 \B\ra_\lvv
  = \f{1}{\sqrt 6}(5 \tp - 3 \trho)
  =  \f{1}{\sqrt 6}(2 \tp - 3 \tb) .
\eal
\eeq

By definition of the projection \eqref{def:proj}, we have $ \cMM^{1/2} \tFm \perp 1, V_i, |V|^2$. Thus, applying the above computation with $G = \cMM^{-1/2} V \cdot \na_X( \cMM^{1/2} \tFm )$ and using 
$\cI_i$ defined in \eqref{eq:moments_Fm}, we obtain 
\beq\label{eq:mom_comp_tranFm}
\bal
  \la \cMM^{-1/2} V \cdot \na_X ( \cMM^{1/2} \tFm) , \Phi_0 \ra_\lvv & = \div _X \int V  \cMM^{1/2} \tFm d V = 0,  \\
    \la \cMM^{-1/2} V \cdot \na_X ( \cMM^{1/2} \tFm) , \Phi_i \ra_\lvv &= \kp^{1/2} \cI_{1, i}, \\
    \la \cMM^{-1/2} V \cdot \na_X ( \cMM^{1/2} \tFm) , \Phi_4 \ra_\lvv & = \f{1}{\sqrt 6} \B\la V \cdot \na_X (  \cMM^{1/2} \tFm),  
    3 \kp \cdot \f{ | V - \bu|^2 }{ 3 \cs^2 }  \B\ra_\lvv
    = \f{3 \kp}{\sqrt{6}} \cI_2 ,  \\
  \eal
  \eeq
where we have used $ \la G, 1 \ra_V = 0$ in the last identity.

We define the linear operator $\cF_M$ mapping the hydrodynamic fields $(\tu, \tp, \tb)$ to macro-perturbation
\bseq\label{eq:macro_UPB}
\beq
\cF_M( \tu, \tp, \tb ) := ( \tp + \tB ) \Phi_0 + \kp^{1/2} \tu_i \Phi_i 
  + \sqrt{\f{1}{6}} (2 \td P - 3 \tb) \Phi_4.
\eeq
Then we can rewrite $\tFM$ defined in \eqref{def:proj} as 
\beq
  \tFM = \cF_M(\tu, \tp, \tb). 
\eeq
\eseq

Recall the operator $\cF_E$ in \eqref{eq:pertb_macro}.  
For any $\td F$, we denote $(\tu, \tp, \tb):= \cF_E(\td F), \tFM = \cP_M \td F$. By definition of $\cP_M$ in \eqref{def:proj}, we obtain 
\beq\label{eq:macro_UPB_inv}
\bga
   (\tu, \tp, \tb) := \cF_E( \td F ) = \cF_E(\tFM) = \cF_E \circ \cF_M(\tu, \tp, \tb), \\
   \tFM = \cF_M(\tu, \tp, \tb) =  \cF_M \circ \cF_E(\td F) =  \cF_M \circ  \cF_E(\tFM) ,  \\
  \quad \cF_M \circ \cF_E = \cP _M, \quad \cF_E \circ \cF_M = \Id,
\ega
\eeq
Thus, $\cF_M$ and $\cF_E | _{\operatorname{Span} \{\Phi _i \}}$ are inverse operators. We estimate the operator $\cF_M, \cF_E$ in Lemma \ref{lem:macro_UPB}.

\section{Linear stability estimates: macroscopic part }\label{sec:lin_macro}

In this section, we perform linear stability estimates on the hydrodynamic fields 
$(\tu, \tp, \tb)$ in \eqref{eq:lin_euler}. Throughout this section, we simplify the perturbation $(\tu, \tp, \tb)$ as $(\UU, P, B)$.

Firstly, we design the weight $\vp_{2 k}$ for weighted $H^{2 k}$ estimates. 
Recall the sonic point $\xi_*$ defined in \eqref{eq:sonic_pt}. We have the following results similar to \cite[Lemma 3.1]{chen2024Euler}.

\begin{lemma}[Lemma 3.2 \cite{chen2024vorticity}]\label{lem:wg}
There exists a radially symmetric weight $\vp_1$ in the form of 
\footnote{
The parameter $c_2, c_3$ in \eqref{eq:repul_form_wg} corresponds to $(\kp_2, \nu )$ in \cite[Lemma 3.1]{chen2024Euler}. 
The forms \eqref{eq:repul_form_wg:a}, \eqref{eq:repul_form_wg:b} are given in 
\cite[Eqn (3.5), Eqn (3.6)]{chen2024Euler}, respectively. 
}
\bseq\label{eq:repul_form_wg}
\beq\label{eq:repul_form_wg:a}
 \vp_1(y) := \vp_b(y)^{c_2} \vp_f(y),  \qquad \vp_f(y) := 1 + c_3  \la y \ra, 
\eeq
where $ c_2, c_3 > 0$ and $\vp_b \in C^{\infty}$ satisfies 
\beq\label{eq:repul_form_wg:b}
\bal 
\vp_b(y)  &= 1, \quad |y| \leq \xi_*, & \quad    \vp_b(y) =  \tfrac 12 , & \quad |y| \geq R_2 + 1, \\
\pa_{\xi} \vp_b & \leq 0, \quad \forall y \in {\mathbb R}^3,  & \quad 
 \pa_{\xi} \vp_b \leq - c_1  < 0,  & \quad |y| \in [R_1, R_2] , 
\eal 
\eeq
\eseq
for some  $ R_1, R_2$ with $\xi_* < R_1 < R_2$,  and there exists a constant $\mu_1 > 0$, such that 
\begin{subequations}
\beq 
 \f{  (\xi + \bar U)\pa_{\xi} \vp_1 }{\vp_1} 
+ i  \cs  \B| \f{ \pa_{\xi} \vp_1 }{\vp_1} \B|  - 
(1 + \pa_{\xi} \bar U - i   |\pa_{\xi} \cs |  )  \leq - \frac{\mu_1}{\la  \xi \ra}  , \label{eq:repul_wg}
\eeq 
\end{subequations}
for all $\xi \in (0,\infty)$ and any $i=0, 1$.
For $k \geq 0$ and $\vp_1$ satisfying the above properties, we define 
\begin{equation}
\varphi_k(y) = \varphi_1(y)^k.
\label{eq:wg_vp}
\end{equation}
\end{lemma}

The term \eqref{eq:repul_wg} relates to the coefficients of the top order term in later weighted $H^{2k}$ estimate. In \cite{chen2024Euler}, the proof of the above theorem relies on the 
repulsive properties \eqref{eq:rep1}  for $\xi \in [0, \xi_1]$, and the outgoing property \eqref{eq:rep2} of the imploding profile for 2D compressible Euler. In \cite[Lemma 3.1]{chen2024Euler}, the weight $\vp_b$ is not stated as a $C^{\infty}$ function. 
From the proof of  \cite[Lemma 3.1]{chen2024Euler}, $\vp_b$ can be chosen to be any function satisfying \eqref{eq:repul_form_wg:b}. In particular, we can choose $\vp_b \in C^{\infty}$. 
Since the modified profile $(\bu, \cs)$ for 3D compressible Euler satisfies these properties uniformly in $s$, the proof of Lemma \ref{lem:wg} in the current setting is the same and is omitted.

From \eqref{eq:repul_form_wg}, we obtain the following estimates of $\vp_1$
\bseq \label{eq:wg_asym}   
\beq 
 \vp_1(X) \asymp \la X \ra, \quad |\na \vp_1| \les 1 , \quad 
  | D_X^{\al} \vp_1 | \les_\al  \vp_1  , \\ 
  \eeq 
for any multi-index $\al$, 
where $D_X^{\al}$ is defined in \eqref{eq:deri_wg}. 
Since $\ang X = X + O(\ang X^{-1})$, \eqref{eq:repul_form_wg} implies 
  \beq 
  X \cdot \na_X \log \vp_1 
  = O(\ang X^{-2}) + X \cdot \na_X \log (1 + c_3 \ang X)
  = 1 + O(\ang X^{-1}) .
  \eeq 
 \eseq 

\subsection{Weighted \texorpdfstring{$H^k$}{H k} coercivity estimates}
\label{sec:est_L}

\begin{theorem}
\label{thm:coer_est}

Let $\cL_{U, s}, \cL_{P, s}, \cL_{B, s}$ be defined in \eqref{eq:lin_euler} 
and $\etab$ in \eqref{wg:X_power}.  Denote $\cL_{E, s} = (\cL_{U, s}, \cL_{P, s}, \cL_{B, s})$. There exists $k_0 \geq 6$ large enough, %
\footnote{The parameters $k_0, R_{\eta}$ and $\lambda_{\eta}$ depend only on the weight  $\vp_1$ from Lemma~\ref{lem:wg}, on $r>0$, and on the profiles $(\bu, \bar  \sc)$.
} 
such that the following statements hold true.  For any $\eta \in [-6,  \etab)$ 
\footnote{
By choosing $\eta > -3$, we can obtain decay estimates of $f$ from $\| f \|_{\cX^k_{\eta}}$ using 
the embedding in Lemma \ref{lem:prod}. The lower bound $-6$ of $\eta$ is not important as we will only apply Theorem \ref{thm:coer_est} with $\eta$ close to $\etab$. We impose it to avoid tracking some constants depending on $\eta$. 
}
and $k \geq k_0$, there exists 
$R_{\eta} > 0, \bar C_{k, \eta} $ large enough and $\vpi_{k, \eta} = \vpi_k(k_0, R_{\eta}, \bar C_{\eta},  \eta)>0$ such that
\bseq\label{eq:coer_est}
\beq
\bal\label{eq:coer_est:a}
  \la  \cL_{E, s} ( \UU, P, B),  (\UU, P, B) \ra_{\cX_{\eta}^{2 k}} 
  & \leq 
  -  \lam_{\eta} \|  (\UU, P, B) \|_{\cX_{\eta}^{2k}} ^2
  + \bar C_{k,\eta} \int_{|X| \leq R_{\eta} } |( \UU, P, B)|^2 d X ,
  \\
  \eal 
\eeq
where $\lam_{\eta}$ is defined as follows and is independent of $k$ 
\footnote{
  From \eqref{eq:r_limit}, \eqref{eq:Euler_profi}, and \eqref{wg:X_power}, we have $\f{\bcv}{ \bcx } = \f{1/r - 1}{1/r} 
  =-(r-1)$ and $\etab = -3 + 6 (r-1) = -3 - 6 \f{\bcv}{ \bcx }$.
}
\beq\label{eq:lam_eta}
  \lam_{\eta} = - \f{1}{2} \left(  3 \bcv +    \f{\bcx}{2} (\eta + 3) \right)
= \f{\bcx}{4} ( \etab - \eta )
   > 0, \quad 
 \mathrm{ \ for \ } \eta < \etab.
\eeq

For $\eta = \etab$ and any $ k \geq k_0$, there exists $\varpi_{k, \eta} = \varpi_k(k_0, R_1, \bar C,  \eta)>0$
and a constant $\bar C_{k, \etab} > 0$, such that 
\beq\label{eq:coer_est:b}
    \la  \cL_{E, s} ( \UU, P, B),  (\UU, P, B) \ra_{\cX_{\etab}^k} 
   \leq  \bar C_{k, \etab}  \int \la X \ra^{\etab - r} | ( \UU, P, B)|^2 d X .
\eeq
\eseq
Here, the Hilbert spaces $\cX^{n}_{\eta}$ are defined as the completion of the space of $C^\infty_c(\R^3)$ {\em radially-symmetric}
\footnote{
By radially symmetric functions we mean $f(y) = f (|y|)$ and by radially symmetric vectors we mean $\boldsymbol{f}(y) = f(|y|) \f{y}{|y|}$.
} 
scalar/vector functions, with respect to the norm induced by the inner products 
\footnote{It is also convenient to denote $\cX^{\infty} = \cap_{k \geq 0} \cX^{k}$.} 
\begin{align}
\la (\UU_a, P_a, B_a) , ( \UU_b,  P_b,  B_b) \ra_{\cX_{\eta}^{2k}}
& :=  \int  \sum_{g = \UU, P, B } w_{g} (   \D^k g_a \cdot \D^k  g_b  \; \vp_{2 k}^ 2  + 
\vpi_{k, \eta}
g_a \cdot  g_b ) \; \la X \ra^{\eta}   d X,  \ k \geq 1,  \notag \\
\la (\UU_a, P_a, B_a) , ( \UU_b,  P_b,  B_b) \ra_{\cX_{\eta}^{2k + 1}}
& :=  \int  \sum_{g = \UU, P, B } w_{g} (  \na  \D^k g_a \cdot  \na \D^k  g_b  \; \vp_{2 k + 1}^ 2  + 
\vpi_{k, \eta}
g_a \cdot  g_b ) \; \la X \ra^{\eta}   d X,  \ k \geq 0,  \notag \\
\la (\UU_a, P_a, B_a) , (  \UU_b,  P_b,  B_b) \ra_{\cX_{\eta}^0}
& :=  \int  \sum_{g = \UU, P, B } w_{g}  g_a \cdot  g_b  \; \la X \ra^{\eta}   d X,  \notag \\ 
(w_\UU, w_P, w_B ) & :=  (1, 1, \tf{3}{2}) , \label{norm:Xk}
\end{align}
The inner products $\la\cdot,\cdot\ra_{\cX_{\eta}^{ n }}$, and the associated norms, are defined in terms of the constants $ \vpi_{n,\eta}$ (defined in the last paragraph of Section~\ref{sec:est_L}), the weight $\vp_{n} = \vp_1^{n}$ defined in Lemma \ref{lem:wg}. 
In particular, these constants $k_0, R_{\eta}, \bar C_{k,\eta}, \vpi_{k,\eta},  \eta $ are independent of 
$s,\e_0$ and $R_0$ used in \eqref{eq:S_radial}.
\end{theorem}

The special weights $(w_\UU, w_P, w_B)$ are determined by 
the relationship between $\int |\cM^{1/2} F_M|^2 d V$ and $(\tu, \tp, \tb)$. See Lemma \ref{lem:macro_UPB}.

\begin{remark}[The odd order norm]
The norm $\cX^{2k+ 1}$ with odd index is auxiliary. We only perform energy estimates, construct compact operator, and develop semigroup estimates in the norm $\cX^{2k}_{\eta}$. Here, $\cX^{2k}_{\eta}$ corresponds to the norm $\cX^{k}_{\eta}$ in \cite{chen2024Euler}.
\end{remark}

\begin{remark}[Full stability of $\tb$]
The linear evolution of $\tb$ in \eqref{eq:lin_euler} is almost decoupled from $\td \UU, \td P$, 
and one can establish full stability of $\cL_{B, \infty}$ with radial symmetry using weighted estimates. %
\end{remark}

Before proving Theorem~\ref{thm:coer_est}, we note the following simple {\em nestedness property} of the spaces $\cX_\eta^n$, which follows from Lemmas \ref{lem:interp_wg}, \ref{lem:norm_equiv} with $\d_1 = 1, \d_2 = a$ or $\d_2 = b$. 
\begin{lemma}\label{lem:Xm_chain}
For any $n \geq m$ and $a \geq b$, we have $\|f \|_{\cX_{b}^m} \les_{n, a, b} \|f \|_{\cX_{a}^n}$ and $\cX_{a}^n \subset \cX_{b}^m$.
\end{lemma}

Next, we prove  Theorem~\ref{thm:coer_est}. We will drop $\td \cdot$ in the variables $\tu, \tp, \tb$ to simplify the notations, and write $\WW = (\UU, P, B)$.

\begin{proof}[Proof of Theorem \ref{thm:coer_est}]
Applying the operator $\D^k$ to the linearized operators $\cL_{U, s}, \cL_{P, s}, \cL_{B,s}$ 
defined in \eqref{eq:lin_euler}, and using Lemma~\ref{lem:leib} to extract the leading order parts 
with $\geq 2k$-derivatives 
from the terms containing $\na \UU,  \na P, \na B$,  we get 
\begin{subequations}
\label{eq:lin_Hk}
\begin{align}
\D^k \cL_{U,s} \WW & = 
\underbrace{- \, (\bcx  X + \bar \UU) \cdot \na \D^k \UU - \cs \na \D^k P }_{\cT_{U}}
\underbrace{+ \, 3 \bcv   \D^k \UU- 2 k \pa_{\xi}( \bcx \xi + \bar U) \D^k \UU- 2 k  \na \cs \D^k P  }_{ \cD_{U} } \notag \\
& \quad 
\underbrace{- \, \left( \f{2}{3} \na \cdot \bar \UU  + (\na \bu) +\ec \right) \D^k  \UU
- 2 \na \cs \cdot  \D^k P + 3 \cs^{-1} \bar \sc \na \bar \sc  (\D^k P + \D^k B)  }_{ \cS_U}
+ \cR_{U, k}, \\
\D^k \cL_{P,s} \WW & = 
\underbrace{- \,(\bcx X + \bar \UU) \cdot \na \D^k P - \cs \na \cdot (\D^k \UU)  }_{\cT_{P}}
\underbrace{ +\, 3 \bcv \D^k P
- 2 k  \pa_{\xi}( \bcx \xi + \bar U) \D^k P - 2 k  \na \cs \cdot \D^k \UU   }_{ \cD_{P}} \notag \\
& \quad 
\underbrace{ -\, ( \na \cdot \bu +  2\ec ) \cdot  \D^k P 
- \left(\na \cs + \f{2}{3} \eu \right) \cdot \D^k \UU  }_{\cS_{P}} + \cR_{P, k}, \\
\D^k \cL_{B, s} \WW & = 
\underbrace{ -\, (\bcx X + \bar \UU) \cdot \na \D^k B }_{\cT_{B}}
\underbrace{ +\, 3 \bcv   \D^k B
- 2 k \pa_{\xi}( \bcx \xi + \bar U) \D^k B
}_{ \cD_{B}} \notag \\
& \quad  
\underbrace{ -\, ( \na \cdot \bu) \cdot  \D^k B + 2 \ec \D^k P + \f{2}{3} \eu \cdot \D^k \UU }_{\cS_{B}} + \cR_{B, k}, 
\end{align}
\end{subequations}
In~\eqref{eq:lin_Hk} we have denoted by $\cR_{U, k}, \cR_{P,k}, \cR_{B,k}$ {\em remainder} terms which are of lower order (in terms of highest derivative count on an individual term); moreover, we have used the notation $\cT, \cD,\cS$ to single out {\em transport}, {\em dissipative}, and {\em stretching} terms.

Using Lemma \ref{lem:leib} and the decay estimates \eqref{eq:dec_U}, \eqref{eq:dec_S} on $\ths$, $\cs$,
and \eqref{eq:error_ec}, \eqref{eq:error_eu} on $\ec$, $\eu$, we obtain that the remainder terms are bounded as  
\begin{subequations}
\label{eq:lin_lower}
\begin{align}
|\cR_{U, k}| 
&\les_k \sum_{ 0 \leq i \leq 2 k-1}  ( | \na^{ 2 k + 1- i}  \bar \UU | +|\na^{2k - i} \ec| ) \; | \na^i \UU |  
+ |\na^{2 k + 1 - i} \cs | \; |\na^i P|  \notag \\
& \qquad + |\na^{2k-i} (\cs^{-1} \bar \sc \na \bar \sc)| \; |\na^i (P + B)| 
\notag\\
&\les_k \sum_{0 \leq i \leq 2 k - 1} 
\la X \ra^{-2k +i - r}| (  |\na^i \UU | +  | \na^i B |   +  |\na^i P|  ) , \\ 
|\cR_{P, k}| 
&\les_k \sum_{ 0 \leq i \leq 2 k-1} (  | \na^{ 2 k + 1- i}  \bar \UU |
+ |\na^{2k-i} \ec| )\; | \na^i P |  
+ ( |\na^{2 k + 1 - i} \cs | + |\na^{2k-i} \eu| ) \; |\na^i \UU|
\notag\\
&\les_k \sum_{0 \leq i \leq 2 k - 1} 
\la X \ra^{-2k + i - r}| (  |\na^i \UU | +  | \na^i B |   +  |\na^i P|  ) , \\ 
 |\cR_{B, k}|  &\les_k \sum_{ 0 \leq i \leq 2 k-1}   | \na^{ 2 k + 1- i}  \bar \UU | 
 \cdot |\na^i B|
 + |\na^{2k-i} \ec| \; |\na^i P| + |\na^{2k-i} \eu| \; |\na^i \UU| \notag \\
&\les_k \sum_{0 \leq i \leq 2 k - 1} 
\la X \ra^{-2k +i - r}| (  |\na^i \UU | +  | \na^i B |   +  |\na^i P|  ) .
\end{align}
\end{subequations}

Next, in order to bound the left side of~\eqref{eq:coer_est}, we perform weighted $H^{2k}$ estimates with weight given by $\vp_{2k}^2 \la X\ra^{\eta}$, as dictated by the definitions in~\eqref{eq:coer_est}. 
 To this end, we estimate the term
\begin{equation}
 \int ( \D^k \cL_{U, s} \WW \cdot \D^k \UU 
 + \D^k \cL_{P,s} \WW \, \D^k P 
+ w_B \D^k \cL_{B, s} \WW \cdot \D^k B   ) 
 \vp_{2 k}^2 \la X \ra^{\eta} d X,
 \label{eq:main:term:0}
\end{equation}
by appealing to the decomposition in~\eqref{eq:lin_Hk}.

\paragraph{Estimate for $\cT_U, \cT_P, \cT_B$}

We first combine the estimates of $\cT_U, \cT_P$, and then estimate $\cT_B$. Using the identity
\beq\label{eq:lin_cross1}
 \na \D^k P \cdot \D^k \UU + \na \cdot (\D^k \UU) \D^k P = \na \cdot ( \D^k \UU \cdot \D^k P )
\eeq
and integration by parts, we obtain that the contribution of the transport terms 
$\cT_U, \cT_P$ in \eqref{eq:lin_Hk} to the expression~\eqref{eq:main:term:0} is given by 
\[
\bal
I_{\cT_U + \cT_P} & = -\int \B( ( \bcx X + \bar \UU) \cdot \na \D^k \UU \cdot \D^k \UU  
+ ( \bcx X + \bar \UU) \cdot \na \D^k P \cdot \D^k P  \\
& \qquad \qquad \qquad \qquad +   \cs \na \D^k P \cdot \D^k \UU
+\cs \na \cdot ( \D^k \UU \cdot \D^k P )
\B) \vp_{2 k}^2 \la X \ra^{\eta} d X \\ 
& = \int \B(  \f{1}{2}  \f{ \na \cdot ( (\bcx X + \bar \UU) \vp_{2 k}^2 \la X \ra^{\eta} ) }{ \vp_{2 k}^2 \la X \ra^{\eta} }
( |\D^k \UU|^2 + |\D^k P|^2  ) %
+ \f{ \na ( \cs \vp_{2 k}^2 \la X \ra^{\eta} )}{ \vp_{2 k}^2 \la X \ra^{\eta} }\cdot \D^k \UU  \D^k P
\B)  \vp_{2k}^2 \la X \ra^{\eta} d X .
\eal
\]

Recall $\vp_{2 k} = \vp_1^{2 k}$ 
from \eqref{eq:wg_vp} and $r < 2$ form \eqref{eq:r_limit}.  Using the decay estimates in \eqref{eq:dec_U}, the outgoing property $\xi + \bar U>0$ \eqref{eq:rep22}, we obtain
\[
 (\bcx \xi + \bar U)  \f{\pa_{\xi} \la X \ra^{\eta}}{ \la X \ra^{\eta}}
 = (\bcx  + \f{\bar U}{\xi}) \f{\eta \xi^2}{1 + \xi^2}
 \leq \bcx \eta + C ( \la \xi \ra^{-r} + \la \xi \ra^{-2} )
  \leq \bcx \eta + C \la \xi \ra^{-r}. 
\]
Using the above inequality and  ~\eqref{eq:dec_U} (with $k=1$),
we estimate
\[
\bal
\f{1}{2}  \f{ \na \cdot ( (\bcx X + \bar \UU) \vp_{2k}^2 \la X \ra^{\eta}) }{ \vp_{2k}^2 \la X \ra^{\eta} }
& = \f{1}{2} \B( 3 \bcx  + \na \cdot \bar \UU
+ 4k ( \bcx \xi + \bar U) \f{\pa_{\xi} \vp_1 }{\vp_1}
 + ( \bcx \xi + \bar U) \f{\pa_{\xi} \la X \ra^{\eta} }{ \la X \ra^{\eta} } \B) \\
 & \leq 
 \f{\bcx}{2} (\eta + 3) + 
 C \la \xi \ra^{-r}+ 2k ( \bcx \xi + \bar U) \f{ \pa_{\xi} \vp_1}{\vp_1} , \\
\B| \f{ \na( \cs \vp_{2k}^2 \la X \ra^{\eta} )}{ \vp_{2k}^2 \la X \ra^{\eta} } \B|
 &  \leq  |\na \cs | + \cs \Bigl( 4k \f{ |\na \vp_1|}{\vp_1}
 + \f{ |\na \la X \ra^{\eta}| }{ \la X \ra^{\eta} }\Bigr) 
 \leq 4 k \cs  \f{ |\pa_{\xi}\vp_1|}{\vp_1}
 + C \la \xi \ra^{- r} ,
  \eal
\]
where we have used $|\na f| = |\pa_{\xi} f| $ for any radially symmetric function $f$. 
 Combining the above estimates and using $|a b| \leq \f{1}{2}(a^2 + b^2) $
on $\D^k \UU \D^k P$, we get
\bseq\label{eq:lin_est1}
\begin{align}
  I_{\cT_U +\cT_P} 
 & \leq \int \f{1}{2} \B(  \f{ \na \cdot ( ( \bcx X + \bar \UU) \vp_{2 k}^2 \la X \ra^{\eta} ) }{ \vp_{2 k}^2 \la X \ra^{\eta} }
 +  \f{ | \na ( \cs \vp_{2 k}^2 \la X \ra^{\eta} ) | }{ \vp_{2 k}^2 \la X \ra^{\eta} }
\B) ( |\D^k \UU|^2 + |\D^k P|^2  ) \vp_{2 k}^2 \la X \ra^{\eta} d X \notag \\
   \leq \int & \B( %
 \f{\bcx}{2} (\eta + 3)
    +  2k \B( ( \bcx \xi + \bar U) \f{\pa_{\xi} \vp_1}{\vp_1} 
+ \cs \Bigl| \f{\pa_{\xi} \vp_1}{\vp_1} \Bigr| \B) + C \la \xi \ra^{-r} \B) 
( |\D^k \UU|^2 + |\D^k P|^2  ) \vp_{2 k}^2 \la X \ra^{\eta} d X,
\end{align}
with $C>0$ independent of $k$.

The estimate of contribution of $\cT_B$  in \eqref{eq:lin_Hk} to the expression~\eqref{eq:main:term:0} is easier and similar. Using integration by parts, we obtain 
\begin{align}
  I_{\cT_B}
  & =  -\int \B( ( \bcx X + \bar \UU) \cdot \na \D^k B \cdot \D^k B \vp_{2k}^2 \la X \ra^{\eta}
  = \int \f{1}{2}  \f{ \na \cdot ( (\bcx X + \bar \UU) \vp_{2 k}^2 \la X \ra^{\eta} ) }{ \vp_{2 k}^2 \la X \ra^{\eta} }  | \D^k B|^2 \vp_{2k}^2 \la X \ra^{\eta} d X \notag \\
  & \leq  \int \B( %
   \f{\bcx}{2} (\eta + 3)
  +  2k \B( ( \bcx \xi + \bar U) \f{\pa_{\xi} \vp_1}{\vp_1} 
 \B) + C \la \xi \ra^{-r} \B) 
|\D^k B|^2 \vp_{2 k}^2 \la X \ra^{\eta} d X,
\end{align}
\eseq

\paragraph{Estimate of $\cD_U, \cD_{P}, \cD_B,  \cS_U, \cS_{P}, \cS_B$}
Recall  the definitions of the terms $\cD_U, \cD_{P}, \cD_B,  \cS_U, \cS_{P}, \cS_B$ from \eqref{eq:lin_Hk}. Using the estimates of the error term $\eu, \ec$ in Lemma \ref{lem: cutoff_error} and the decay estimates \eqref{eq:dec_U},  we get 
\[
\bal
|\cS_U| & \leq 
   C  \la \xi \ra^{-r} (|\D^k \UU| + |\D^k P | + |\D^k B| ),  \\
|\cS_P| & \leq %
C  \la \xi \ra^{-r} (|\D^k \UU| + |\D^k P | ) , \\
|\cS_B| & \leq %
 C  \la \xi \ra^{-r} (|\D^k \UU| + |\D^k B| ) .
\eal 
\]

For $\cD_U, \cD_P$, using Cauchy--Schwarz inequality for the cross term
\beq\label{eq:lin_cross2}
|\na \cs \D^k P \cdot \D^k \UU | 
+ |\na \cs \cdot   \D^k \UU \cdot  \D^k P | 
\leq |\na \cs| ( |\D^k P|^2 + |\D^k \UU |^2  ),  \quad |\na \cs| = |\pa_{\xi} \cs|,
\eeq
we obtain 
\begin{align}
I_{\cD+\cS} &= \int \B(   ( \cD_U + \cS_U)    \D^k \UU +  
 ( \cD_{P} + \cS_{P})   \D^k P 
+ w_B ( \cD_{B} + \cS_{B})   \D^k B
 \B) \vp_{2 k}^2 \la X \ra^{\eta } d X \notag \\
& \leq  \int \B( 3 \bcv - 2 k ( \bcx + \pa_{\xi} \bar U)  
+ 2 k  |\pa_{\xi} \cs | + C \la \xi \ra^{-r} \B) ( |\D^k \UU|^2 + |\D^k P|^2 ) \vp_{2k}^2 \la X \ra^{\eta}\notag \\
& \quad +  w_B  \B( 3 \bcv - 2 k ( \bcx + \pa_{\xi} \bar U)  
+  C \la \xi \ra^{-r} \B)  |\D^k B|^2 \vp_{2k}^2 \la X \ra^{\eta} d X, %
\label{eq:lin_est2}
\end{align}
with $C>0$ independent of $k$.

\paragraph{Estimates of $\cR_U, \cR_{P}$}
Recall that the remainder terms $ \cR_U, \cR_{P}, \cR_B$ from \eqref{eq:lin_Hk} satisfy~\eqref{eq:lin_lower}. Moreover, using that $\vp_{2 k} \asymp_k \la X \ra^{2k}$ from~\eqref{eq:wg_asym},  we obtain 
\[
\vp_{2 k}^2 \la X \ra^{\eta} \la X \ra^{-2 k + i} \asymp_k \la \xi \ra^{2 k + i + \eta} 
= \la \xi \ra^{2 k + i + 2 \cdot \eta/2}, 
\]
At this stage we apply  Lemma~\ref{lem:interp_wg} and Lemma~\ref{lem:norm_equiv}, with $\d_1 = 1$ and $\d_2= \eta - r $ for $1 \leq i \leq 2k-1$, and an arbitrary $\nu > 0$, to obtain 
\begin{align*}
& \int \la X \ra^{ - 2 k + i} |\na^i F| |\D^k G | \vp_{2k}^2 \la X \ra^{\eta - r} d X  \\
&\qquad \leq   \nu  \| \la X \ra^{2k + (\eta - r) /2} \D^k G \|^2_{L^2} 
+ C_{k ,\eta, \nu } \|\la X \ra^{i + (\eta - r) /2 } \na^i F  \|^2_{L^2} \\
&\qquad\leq  
\nu \| \la X \ra^{2k + (\eta - r) /2} \D^k G  \|_{L^2}^2
+ \nu \| \la  X\ra^{2k + ( \eta - r) /2} \na^{2 k} F   \|_{L^2}^2
+ C_{k, \eta, \nu} \| \la X \ra^{ ( \eta - r )/2} F \|_{L^2}^2 \\
&\qquad\leq   2 \nu \| \la X \ra^{2 k + (\eta - r) /2} \D^k G  \|^2_{L^2} 
+ 2\nu \| \la X \ra^{2k + (\eta - r) /2} \D^k F    \|^2_{L^2} 
+ C_{k, \eta, \nu} \| \la X \ra^{ (\eta-r)/2}  F \|_{L^2}^2 .
\end{align*}
We may apply the above estimates to each term in \eqref{eq:lin_lower}. Using the bound $ \la X \ra^{2( 2k + (\eta - r) /2)} \les_k  \vp_{2k}^2 \la X \ra^{ \eta - r } $, and choosing $\nu = \f{1}{2}$ in the above estimates, we get 
\begin{align}
I_{\cR} &= \int \bigl| \cR_{U, k}  \cdot \D^k \UU + \cR_{P, k} \cdot \D^k P
+ w_B \cR_{B, k} \cdot \D^k B
 | \vp_{2k}^2 \la X \ra^{\eta} d X \notag \\
& \leq \int %
( |\D^k \UU|^2 + |\D^k P|^2 + w_B |\D^k B|^2 ) \vp_{2k}^2 \la X \ra^{\eta - r}
+ C_{k, \eta}  |(\UU, P, B)|^2  \la X \ra^{\eta - r} d X.
\label{eq:lin_est3}
\end{align}

Combining the bounds \eqref{eq:lin_est1}, \eqref{eq:lin_est2}, \eqref{eq:lin_est3} and using the estimates in Lemma~\ref{lem:wg}, 
we arrive at
\begin{align*}
&  \int \left( 
    \D^k \cL_{U,s} \WW \cdot \D^k \UU + 
    \D^k \cL_{P, s} \WW \cdot \D^k P + 
    w_B \D^k \cL_{B,s} \WW \cdot \D^k B 
\right)
\vp_{2k}^2 \la X \ra^{\eta} d X
= I_{\cT} + I_{\cD+\cS}  + I_{\cR} \notag  \\
 & \leq 
\int \B\{ \left(
    3 \bcv + \f{\bcx}{2} (\eta + 3) +  C \la \xi \ra^{-r} 
\right) ( |\D^k \UU|^2 + |\D^k P |^2 + w_B |\D^k B|^2) \vp_{2k}^2 \la X \ra^{\eta}  \\ 
& \quad \qquad + 2k \B( (\xi + \bar U) \f{\pa_{\xi} \vp_1}{\vp_1} 
 + \cs | \f{\pa_{\xi} \vp_1}{\vp_1} |  - (1 + \pa_{\xi} \bar U -  | \pa_{\xi} \cs|)\B)  \;  ( |\D^k \UU|^2 + |\D^k P |^2) \vp_{2k}^2 \la X \ra^{\eta} \\
& \quad \qquad + 2k \B( (\xi + \bar U) \f{\pa_{\xi} \vp_1}{\vp_1}  - (1 + \pa_{\xi} \bar U  ) \B)  \;  
\cdot w_B |\D^k B|^2 \vp_{2k}^2 \la X \ra^{\eta} \B\} d X %
\\
& \qquad + C_{k, \eta} \int  |(\UU, P, B)|^2 \la X \ra^{\eta- r} d X \notag \\
&  \leq \int \left( %
- 2 \lam_{\eta}
- 2k \mu_1 \la \xi \ra^{-1} + a_1 \la \xi \ra^{- r} \right)  ( |\D^k \UU|^2 + |\D^k P|^2 + w_B |\D^k B|^2 ) \vp_{2k}^2 \la X \ra^{\eta} d X \\
& \qquad + C_{k, \eta} \int
|(\UU, P, B)|^2
 \la X \ra^{\eta - r} d X ,
\end{align*}
for some constant $a_1 > 0$ independent of $\eta$, %
where we have used the notation $\lam_{\eta} = - \f{1}{2} (  3 \bcv +    \f{\bcx}{2} (\eta + 3)  ) $ defined in 
\eqref{eq:lam_eta}.  
Since $ r > 1$ (see \eqref{eq:r_limit}), there exists $k_0$ sufficiently large, e.g.~$k_0 = \lceil \f{a_1}{2\mu_1} \rceil + 1$, such that for any $ k \geq k_0$, we get
\[
 - 2k \mu_1 \la \xi \ra^{-1} + a_1 \la \xi \ra^{- \kp_3} 
\leq  - 2 k_0 \mu_1 \la \xi \ra^{-1} + a_1 \la \xi \ra^{- \kp_3} \leq 0.
\]
resulting in 
\begin{align}
& \int (\D^k \cL_{U,s} \WW \cdot \D^k \UU + \D^k \cL_{P,s} \WW \D^k P
+ w_B \D^k \cL_{B, s} \WW \D^k B )   \vp_{2k}^2 \la X \ra^{\eta} d X
\notag\\
& 
\qquad \leq %
  - 2 \lam_{\eta} \int ( |\D^k \UU|^2 + |\D^k P|^2 + w_B |\D^k B|^2) \vp_{2k}^2 \la X \ra^{\eta} 
+ C_{k, \eta}  |(\UU, P, B)|^2  \la X \ra^{\eta - r} d X.
\label{eq:lin_coerHk}
\end{align}

\paragraph{Weighted $L^2$ estimates}
For $k = 0$, we do not have the lower order terms $\cR_{\cdot, 0}$ in \eqref{eq:lin_lower} and we do not need to estimate $I_{\cR}$ as in \eqref{eq:lin_est3}. Combining \eqref{eq:lin_est1} and \eqref{eq:lin_est2} (with $k=0$), we obtain 
\begin{equation}\label{eq:lin_coerL2a}
\bal
& \int \left( 
    \cL_{U,s} \WW \cdot \UU + 
    \cL_{P, s} \WW \cdot P + 
    w_B \cL_{B,s} \WW \cdot B 
\right) \ang X ^\eta d X \\
& \qquad \leq \int \left(
    3 \bcv + \f{\bcx}{2} (\eta + 3) +  C \la \xi \ra^{-r} 
\right) ( | \UU|^2 + |P|^2 + w_B |B|^2 ) \la X \ra^{\eta} d X.
\eal
\end{equation}
for some constant $ C>0$, independent of $k$ and $\eta$.

If $\eta  \in [ -6, \etab)$, we have 
\[
 2 \lam_{\eta} = - \bigl( 3 \bcv +    \f{\bcx}{2} (\eta + 3)   \bigr)   > 0. 
\]

Thus, there exists a sufficiently large $R_{\eta}$ such that for all $\xi= |X| \geq R_{\eta}$ we have
\[
3 \bcv +    \f{\bcx}{2} (\eta + 3)  
+  C \la \xi \ra^{-r} 
= 
- 2 \lam_{\eta} +  C \la \xi \ra^{- r} \leq - 2 \lam_{\eta}  +  C \la R_{\eta} \ra^{- r} \leq 
- \f{3}{2} \lam_{\eta} .
\]
Therefore, combining the two estimates above and taking  $\bar C_{\eta} =  C  w_B \max \{ \la R_{\eta} \ra^{\eta - r}, 1\}$, we arrive at 
\begin{equation}
\bal
& \int \left( 
    \cL_{U,s} \WW \cdot \UU + 
    \cL_{P, s} \WW \cdot P + 
    w_B \cL_{B,s} \WW \cdot B 
\right) \ang X ^\eta d X \\
& \qquad \leq \int %
- \f{3}{2} \lam_{\eta}  ( | \UU|^2 + | P|^2 + w_B |B|^2 )  \la X \ra^{\eta}  + \bar C_{\eta} \one_{|X| \leq R_{\eta}}  | (\UU, P, B)|^2 d X . \label{eq:lin_coerL2}
\eal
\end{equation}

\paragraph{Choosing $\vpi_{k, \eta}$} 

In order to conclude the proof of~\eqref{eq:coer_est},  we combine \eqref{eq:lin_coerHk}, \eqref{eq:lin_coerL2a}, and \eqref{eq:lin_coerL2}. 

If $\eta \in [-6, \etab)$, choosing $\vpi_k $ sufficiently small, e.g.~$\vpi_{k, \eta} = \f{ 4 C_{k, \eta}}
{ |\lam_{\eta}| }$, where $C_{k,\eta}$ is as in \eqref{eq:lin_coerHk}, 
multiplying \eqref{eq:lin_coerL2} with $\vpi_{k, \eta}$ and then adding to~\eqref{eq:lin_coerHk}  we deduce \eqref{eq:coer_est:a} with $\lam_{\eta}$ independent of $k$.

If $\eta = \etab$, we have $\lam_{\eta} = 0$. We choose $\vpi_{k, \eta} = 1$. Multiplying 
\eqref{eq:lin_coerL2a} with $\vpi_{k, \eta}$ and then adding to~\eqref{eq:lin_coerHk}, we deduce \eqref{eq:coer_est:b}.
\end{proof}

\subsection{Compact perturbation}

We follow \cite{chen2024Euler} to construct a compact operator $\cK_{k, \eta}$ such that $\cL_{E, s} - \cK_{k, \eta}$ is dissipative in $\cX_{\eta}^{2k}$. We fix $k_0\geq 6$ and restrict $\eta \in [-6, \etab)$. 
\footnote{
The result similar to Proposition \ref{prop:compact} was first proved in \cite[Proposition 3.4]{chen2024Euler} for the stability analysis for the imploding profile of the 2D isentropic Euler equations. 
}

\begin{proposition}\label{prop:compact}
For any $k \geq k_0 \geq 6$ and $\eta \in [-6, \etab) $, there exists a bounded linear operator $\cK_{k, \eta} \colon \cX_{\eta}^0 \to \cX_{\eta}^{2 k}$ independent of $s , R_0$ in the definition of $\cs$ \eqref{eq:Euler_profi_modi}, \eqref{eq:S_radial} with:
\begin{enumerate}[\upshape(a)]
\item
for any $f \in \cX_{\eta}^{0}$ we have
\[
{\rm supp}(\cK_{k, \eta} f) \subset B(0, 4 R_{\eta}),
\]
where $R_\eta$ is chosen in Theorem~\ref{thm:coer_est} (in particular, it is independent of $k$);

\item the operator $\cK_{k, \eta}$ is compact from $ \cX_{\eta}^{2 k} \to \cX_{\eta}^{2k}$;  
\item the enhanced smoothing property  $\cK_{k, \eta}: \cX_{\eta}^{0} \to \cX_{\eta}^{2 k+ 6}$ holds;
\item the operator $\cL _{E, s} - \cK_{k, \eta}$ is dissipative on $\cX^{2k}_{\eta}$ and we have the estimate
\begin{equation}
    \la (\cL _{E, s} - \cK_{k, \eta}) f ,f \ra_{\cX_{\eta}^{2 k} } \leq - \lam_{\eta} \| f \|_{\cX_{\eta}^{2k} }^2
    \label{eq:dissip}
\end{equation}
for all $f \in \{ (\UU, P, B) \in \cX_{\eta}^{2k} \colon \cL _{E,s} ( \UU, P, B) \in \cX_{\eta}^k\}$, $\cL _{E,s} = (\cL_{U,s}, \cL_{P,s}, \cL_{B,s} )$, 
and any $s \geq0$, where $\lambda_{\eta} >0$ is the  parameter from~\eqref{eq:coer_est} (in particular, it is independent of $k$).
\end{enumerate}
\end{proposition}

In \cite{chen2024Euler}, the compact operator $\cK_k = \bar C_{\eta} \cK_{k, 0} $ is constructed by applying the Riesz representation theorem in the Hilbert space $\cX_{\eta}^{2k}$ to the bilinear form 
\[
\la \cK_{k, 0} f , g \ra_{\cX_{\eta}^{2k}} :=
\int \chi  f \cdot  g d X
\]
for some smooth cutoff function $\chi$ supported in $B(0, 4 R_1)$ and vector value functions $f, g$, 
with $R_1, \bar C, \cX_{\eta}^{2k}$ chosen in the coercivity estimates similar to those in Theorem \ref{thm:coer_est}. 
The proof of the properties of $\cK_{k, \eta}$ in
 \cite{chen2024Euler} is based the Riesz representation theory and the Rellich--Kondrachov compact embedding theorem. When we apply the argument in \cite{chen2024Euler}, since the space $\cX^{2 k}_{\eta}$, and the parameters $R_1, \bar C_{k, \eta}$, are independent of 
$R_0, \d_R,s$ in  \eqref{eq:S_radial} and \eqref{eq:Euler_profi_modi}, the
operator $\cK_{k, \eta}$ constructed by the same argument associated to $\cX^{2k}_{\eta}$ is 
independent of  $R_0, \d_R,s$. Since the proof is the same, we omit it.

\subsection{Semigroup estimates of the limiting operator \texorpdfstring{$\cL_{E}$}{L E}}
\label{sec:semi}

To estimate the unstable part, we  will apply semigroup estimates to 
the limiting operator $\cL_E =\cL_{E,\infty}$ as $s \to \infty$ (see \eqref{eq:lin_euler_limit}), which is time-independent, and then estimate the error $\cL_{E, s} - \cL_{E,\infty}$ perturbatively.
We estimate $\cL_E, \cL_{E,s}$ in $\cX_{ \xeta}^{2 k}$ with $\xeta$ chosen in \eqref{eq:eta_constraint}.

Applying Theorem \ref{thm:coer_est} and Proposition \ref{prop:compact}
with $\eta \rsa \xeta$, we obtain stability estimates of $\cL_E, \cL_{E, s}$ in $\cX_{\xeta}^{2 k}$, and construct the compact operator $\cK_{k, \xeta}$ with the properties in Proposition \ref{prop:compact}. We recall the notations from \eqref{eq:eta_constraint}
\beq\label{norm:X_lim}
\lame = \f{\bcx}{4} (\etab - \xeta) > 0 .
\eeq

\subsubsection{Complex Banach space $\cX_{\C, \eta}^{k}$}

To apply various functional analysis argument, following \cite{chen2024Euler}, we introduce the complex Banach space $\cX_{\C, \eta}^{k}$ associated with $\cX_{\eta}^k$. 
Recall the inner product $\la \cdot ,\cdot \ra_{\cX_{\eta}^k}$ from \eqref{norm:Xk}. For any vector value functions $f, g$, we define the inner product 
\beq\label{norm:Xmc}
 \la f , g \ra_{\cX_{\C, \eta}^k} := \la f , \bar g \ra_{\cX_{\eta}^k}.
\eeq
The Hilbert spaces $\cX_{\C, \eta}^k$ is defined as the completion of 
complex-valued $C_c^{\infty}$ radially-symmetric scalar/vector functions with respect to the norm induced by the above inner products. It is not difficult to show that 
\beq\label{eq:comp_norm}
\| f \|_{\cX_{\C, \eta}^m}^2 = 	\la f, f \ra_{\cX_{\C, \eta}^m} 
	= \| \Re f \|_{\cX_{\eta}^m}^2
	+ \| \mathrm{Im} f \|_{\cX_{\eta}^m}^2,
\eeq
where $\mathrm{Im} f $ denotes the imaginary part of $f$. %

Using linearity and following \cite[Section 3.4-3.5]{chen2024Euler}, for any real bounded linear operator $\cB : \cX_{\eta}^a \to \cX_{\eta}^b$ with some $a, b\geq 0$, we define its extension to
$\cX_{\C, \eta}^a \to \cX_{\C, \eta}^b $ using linearity : 
\beq\label{eq:op_lin}
	\cB (f + i g) = \cB f + i \cB g , \quad \forall f, \ g \in \cX^m .
\eeq
From \eqref{eq:comp_norm}, it is not difficult to show that $\cB$ defined as above is a bounded complex linear operator with 
\beq\label{eq:op_norm}
	 \| \cB \|_{\cX^k_{\C, \eta}} = \| \cB \|_{\cX^k_{\eta}}. 
\eeq

\subsubsection{Construction of the semigroup}\label{sec:semi_construct}

We follow \cite{chen2024Euler} to construct strongly continuous semigroups generated by $\cL_E, 
\cD_k = \cL_E - \cK_{k, \xeta}$ for any $k \geq k_0$.
We define the domains of $\cL_E, \cD_{k} $ as 
\beq\label{eq:domain}
D_k(\cD_k) = D_k(\cL) := \{  (\UU, P, B) \in \cX^{2 k}_{\C, \xeta}, \; \cL_E (\UU, P, B) \in \cX_{\C, \xeta}^{2 k} \}. 
\eeq

 We have the following results for $\cL_E, \cD_k$.

\begin{proposition}\label{prop:semi}
Suppose that $k \geq k_0$ and $\cK_{k, \xeta}$ is the compact operator constructed in Proposition \ref{prop:compact}. The operators $\cL_E, \cD_k = \cL_E - \cK_{k,\xeta} : D_k \subset \cX_{\C, \xeta}^{2k} \to \cX_{\C, \xeta}^{2k}$ generate strongly continuous semigroups 
\[
\operatorname e^{s \cL_E} : \cX_{\C, \xeta}^{2 k} \to  \cX_{\C, \xeta}^{2 k} ,
\quad 
\operatorname e^{s \cD_k} : \cX_{\C, \xeta}^{2 k } \to  \cX_{\C, \xeta}^{2 k} . 
\]
We have the following estimates of the semigroup 
\beq\label{eq:dissp_semi}
 \| \operatorname e^{s \cD_k} \|_{ \cX_{\C, \xeta}^{2 k} \to  \cX_{\C, \xeta}^{2 k} } \leq \operatorname e^{-\lame s} , \quad \| \operatorname e^{s \cL _E} \|_{ \cX_{\C, \xeta }^{ 2 k  } \to  \cX_{\C, \xeta }^{ 2 k  } } \leq \operatorname e^{C_k s}
\eeq
for some $C_k>0$, and the following spectral property of $\cD_k$
\beq\label{eq:spec}
 \{ z: \Re(z) > -\lame \} \subset \rho_{\mw{res}}(\cD_k) ,
\eeq
where $\rho_{\mw{res}}(\cA)$ denotes the resolvent set of an operator $\cA$. 

Moreover, the semigroups map the real Banach space into the real Banach space:
\[
    \operatorname e^{s \cL _E} : \cX_{\xeta}^{2 k} \to  \cX_{ \xeta}^{2 k} ,
    \quad 
    \operatorname e^{s \cD_k} : \cX_{  \xeta}^{2 k } \to  \cX_{  \xeta}^{2 k} .
\]
\end{proposition}

The construction and estimates of the semigroups $\operatorname e^{s \cL}, \operatorname e^{s \cD_k}$ in \cite{chen2024Euler} are based on the following steps.

\begin{enumerate}[(1)]
\item Solve a linear PDE of $(\UU, \S)$ similar to \eqref{eq:lin_euler_limit}, which is a symmetric hyperbolic system, and prove its well-posedness. 
\footnote{
In the linear stability analysis of 
the isentropic Euler equations in \cite{chen2024vorticity}, 
$\S$ is the perturbation of the \emph{rescaled} sound speed. 
We do not use such a variable in this work.
}

\item Apply the coercivity estimates \eqref{eq:coer_est} to obtain uniqueness of the solution and continuous dependence on the initial data. 

\item To further construct the semigroup $\operatorname e^{ s (\cL -\cK_{k, \xeta})}$ based on $\operatorname e^{s \cL}$, one applies the Bounded Perturbation Theorem \cite[Theorem 1.3, Chapter III]{ElNa2000}.

\item Generalize the estimates on the real Banach space $\cX_{\xeta}^{2 k}$ to the complex Banach space 
$ \cX_{\C, \xeta}^{2k}$ using linearity.
\end{enumerate}

In the current setting, the linear PDE with linear operators $\cL_E = (\cL_U, \cL_{P}, \cL_B)$ defined in \eqref{eq:lin_euler_limit} is a symmetric hyperbolic system. Moreover, we develop the same type of coercivity estimate \eqref{eq:coer_est} as that in \cite{chen2024Euler}. Thus, the proof is the same as that in \cite{chen2024Euler}. After we construct the semigroup, the decay estimate of $\operatorname e^{s \cD_k} = \operatorname e^{s ( \cL_E - \cK_{k, \xeta})}$ follows from the dissipative estimates \eqref{eq:dissip} 
for $\cL_E - \cK_{k, \xeta }$, and implies the estimates of the resolvent set \eqref{eq:spec}. Note that the estimates \eqref{eq:dissp_semi}, \eqref{eq:spec} apply to all  $(\cD_k, \cXC^{2k} )$ with $k \geq k_0$ with $\lame$ independent of $k$. Since the argument is the same, we omit the proof here and refer the reader to \cite[Sections 3.4, 3.5]{chen2024Euler}.

\subsubsection{ Hyperbolic decomposition}\label{sec:semi_grow}

We consider $k \geq k_0$. Denote by ${\sigma}(\cA)$ the spectrum of $\cA$ :
\[
	{\sigma}(\cA) =  \{ z \in {\mathbb C}: z - \cA \mbox{ is not bijective in }  \cXC^{2k} \}. 
\]

We follow \cite{chen2024Euler,chen2024vorticity} to obtain decay estimates of $\operatorname e^{s \cL_E}$ and decompose the space $\cXC^{2k}$. 
Based on Proposition \ref{prop:semi}, using operator theories from \cite[Corollary 2.11, Chapter IV]{ElNa2000}  for the growth bound, and  \cite[Theorem 2.1, Chapter XV, Part IV (Page 326)]{GoGoKa2013} for the spectral projection, we can perform the following decomposition of $\cXC^{2 k}$ and $\s(\cL_E)$. 
The arguments are the same as those in \cite[Section 3.5]{chen2024Euler}; we refer the reader to further discussions therein. Below, we summarize the results.

 Recall the parameter $\lame$ from  \eqref{eq:eta_constraint:c} and \eqref{norm:X_lim} in Theorem \ref{thm:coer_est}, $\el$ from \eqref{eq:idea_size}, and $\rE$ from \eqref{eq:hydro2}. We have $\lame > \rE$. Next, we fix parameters $\lams, \lamu$ with 
\beq\label{eq:decay_para0}
\left(\f{2}{3} \rE - \el\right)< \lams <  \lamu  < \f{2}{3} \rE  < \lame.
\eeq

Due to \eqref{eq:spec},  the set 
\beq\label{eq:sig_eta}
\s_{\mathsf{u}} := \s( \cL_E) \cap  \{ z : \Re(z) > - \lamu \}, 
\eeq
only consists of finite many eigenvalues of $\cL_E$ with finite multiplicity. 
Applying the spectral projection, we can decompose $ \cXC^{2 k }$ into the stable part $\cX_{  \mw{st}}^{2 k}$ associated with the spectrum $\s(\cL_E) \backslash \s_{\mathsf{u}}$ and unstable part $\cX_{ \mw{un}}^{2 k}$ associated with $\s_{\mathsf{u}}$:
\beq\label{eq:dec_X}
\bal
  \cX_{\xeta}^{2k } &= \cX^{2k}_{\mw{un}} \oplus \cX^{2k}_{\mw{st}},  \quad 
  \s( \cL_E |_{\cX^k_{\mw{st} } } ) = \s(\cL_E) \backslash \s_{\mathsf{u}}
\subset \{ z: \Re(z) \leq - \lamu \}, 
  \quad \s( \cL_E|_{\cX^k_{\mw{un} }} ) = \s_{\mathsf{u}}.
\eal 
\eeq
We omit the subindex $ \xeta$ in $\cX_{\mw{st}}^{2 k}, \cX^{2 k }_{\mw{un}}$ since we only apply the decomposition to spaces with this parameter.

The space $\cX^{2 k }_{\mw{un}}$ has finite dimension and can be decomposed as follows 
\beq\label{eq:dec_Xu}
\cX^{2 k}_{\mw{un}} = \bigoplus_{ z \in \s_{\mathsf{u}}} \ker( (z- \cL_E)^{\mu_z}), \quad \mu_z < \infty, \quad |\s_{\mathsf{u}}| < +\infty,
\quad  \dim ( \cX^{2k}_{\mw{un}} ) < \infty .
\eeq

We have the following estimates of the semigroup in these two spaces %
\bseq%
\begin{align}
 \| \operatorname e^{s \cL_E} f \|_{\cXC^{2 k } }  & \leq C_k \operatorname e^{ - \lams s } \| f \|_{\cXC^{2 k} }, \quad \forall f \in \cX^{2 k}_{\mw{st}}, 
 \label{eq:decay_stab} \\
  \| \operatorname e^{-s \cL_E} f \|_{\cXC^{2 k } }  & \leq  C_k \operatorname e^{\lamu s} \| f \|_{\cXC^{2 k }},
   \quad \forall f \in \cX^{2 k}_{\mw{un}} ,
   \label{eq:decay_unstab}
\end{align}
for any $ s > 0$, where $\lams, \lamu$ are defined in \eqref{eq:decay_para0}. 
 \eseq

\subsubsection{Smoothness of unstable directions}\label{sec:smooth_unstab}

We use the following Lemma proved in \cite[Lemma 3.9]{chen2024Euler} to show that
the unstable part $\cX_{\mw{un}}^k$ \eqref{eq:dec_Xu} is spanned by smooth functions.%

\begin{lemma}\label{lem:smooth}
Let $\{ X^i\}_{i \geq 0}$ be a sequence of Banach spaces with $X^{i+1} \subset X^{i}$ for all $i \geq 0$.
Assume that for any $i \geq i_0$ we can decompose the linear operator 
$\cA \colon D(\cA) \subset X^i \to X^i$ as 
 $\cA = \cD_i + \cK_i$, where the linear operators $\cD_i$ and $\cK_i$ satisfy
\beq\label{eq:smooth_ass}
\cD_i : D(\cA) \subset X^i \to X^i , \quad  \cK_i :  X^{i-1}  \to  X^i,
\quad \{ z \in {\mathbb C}\colon \Re(z) > - \lam \} \subset \rho_{\mw{res}}(\cD_i),
\eeq
for some $\lam \in \R$. Here, $\rho_{\mw{res}}(\cdot)$ denotes the resolvent set of an operator and $\lam>0$ is independent of $i\geq i_0$. Fix $n\geq 0$ and $z \in {\mathbb C}$ with $\Re(z) > -\lam$. Assume that the functions $ f_0, \ldots , f_n  \in  X^{i_0}$ satisfy 
\[
(z - \cA) f_0 = 0,
\quad (z - \cA) f_{i} = f_{i-1}, \quad \mbox{ for } \quad 1 \leq i \leq n.
 \]
Then, we have $ f_0, \ldots , f_n  \in X^{\infty} := \cap_{i \geq 0} X^i$. 
\end{lemma}

Consider  Lemma \ref{lem:smooth} with $(\cA,  \{ X^i \}_{i\geq 0} , i_0 , \lam) \rightsquigarrow $
$(\cL_E, \{  \cXC^{2i}\}_{i\geq 0}, k_0, \lame)$ and the decomposition $\cL_E = \cD_k + \cK_{k,\xeta}$ for any $k \geq k_0 \geq 6$, where $\cK_{k,\xeta}$ is constructed in Proposition \ref{prop:compact}. Using Lemma \ref{lem:Xm_chain}, Proposition \ref{prop:compact}, and \eqref{eq:spec}, 
we verify the assumption on $X^i$ and \eqref{eq:smooth_ass}  in Lemma \ref{lem:smooth}. 
Applying Lemma \ref{lem:smooth} to 
$\ker( (z- \cL)^{\mu_z} ) \subset \cX_{\mw{un}}^{ 2 k}$ in \eqref{eq:dec_Xu}
with  $z \in \s(\cL_E) \cap \{ z : \Re(z) > -\lamu \} \subset  \{z : \Re(z) > - \lame \}$ (see \eqref{eq:sig_eta}), we obtain that $\cX_{\mw{un}}^{2 k}$ are spanned by smooth functions in $\cXC^{\infty}
= \cap_{i\geq 0} \cXC^i \subset C^{\infty}$. 
From the definition of $\cXC^i$ \eqref{eq:comp_norm},
$f \in \cXC^i$ if and only if $f, \mathrm{Im} f \in \cXX^i$.  It follows 
$\Re ( \cX^{2k}_{\mw{un}} )  \subset \cap_{i \geq 0} \cXX^{ i} $.

For any $k \ge k_0$ and $s\in \R$,  \eqref{eq:dec_Xu} implies that $\cX^{2k}_{\mw{un}}$ is finite-dimensional and invariant under the operator $\operatorname e^{s \cL_E}$. For any $n \geq 0$, $ k\geq k_0$, and $s\in\R$,  
since all norms in the finite-dimensional space $\Re (\cX^{2k}_{\mw{un}} )
\subset \cap_{i \geq 0} \cXX^{ i}$ are equivalent 
and $\cXX^i$ is a norm for $\Re (\cX^{2k}_{\mw{un}} )$, we obtain 
 \beq\label{eq:X_un_smooth}
 \bal 
   \operatorname e^{s \cL_E} f  \subset \cX^{2 k}_{\mw{un}}, \quad 
   \| \Re \, f \|_{ \cXX^{  n}}
 \les_{n, k} \| \Re \, f \|_{\cXX^{2k}} , \quad \forall f \in \cX^{2 k}_{\mw{un}} .
\eal 
 \eeq

\subsubsection{Additional decay estimates of $\cL_E$ }

In order to localize the initial data (see~\eqref{eq:W2_form:a} below), we will need the following decay estimates for %
$\operatorname e^{s \cL_E} g $ with initial data $g$ supported in the far-field. The following Proposition is similar to \cite[Proposition 3.8]{chen2024Euler}.  

\begin{proposition}\label{prop:far_decay}
Let $R_{\xeta}$ be as defined in Theorem~\ref{thm:coer_est} with $\eta= \xeta$ 
and $\cL_E$ be as defined in \eqref{eq:lin_euler_limit}. Consider the linear equations 
\beq\label{eq:lin_semi2}
  \pa_s \WW = \cL_E \WW ,
  \quad \WW = (\UU, P, B),
\eeq
with initial data $\WW_0 = (\UU_0, P_0, B_0)$ with $\supp(\WW_0) \subset B( 0, R)^\complement$ for some $R > 4 R_{\xeta}> \xi_* $. 
For any $k \geq k_0$, the solution $\WW(s)= \operatorname e^{s \cL_E} \WW_0$ satisfies
\begin{equation*}
\supp(\WW(s)) \subset B( 0, 4 R_{\xeta})^\complement, \qquad  \| \operatorname e^{s \cL_E} \WW_0\|_{\cX_{\xeta}^{2k}} \leq \operatorname e^{-\lame s} \| \WW_0 \|_{\cX_{\xeta}^{2k}}.
\end{equation*}
\end{proposition}

The proof is similar to that of \cite[Proposition 3.8]{chen2024Euler}. %

\begin{proof}

Our key observation is that the support of the solution $\operatorname e^{s \cL_E} \WW_0$ is moving away from $X=0$, remaining outside of $B(0, 4 R_{\xeta})$ for all time. Since $\supp(\cK_{k, \xeta} f) \subset B(0,4 R_{\xeta})$ (item~(a) in Proposition~\ref{prop:compact}), we 
get $\cK_{k, \xeta} \WW(s) = 0$ for all time $s$, and the desired decay estimate follows from \eqref{eq:dissip} or \eqref{eq:dissp_semi}.

Based on the above discussion, we only need to show that the solution satisfies $\WW(s, X)= 0$ for all $X \in B(0, 4 R_{\xeta} )$ and $s>0$. Let $\chi$ be a radially symmetric cutoff function with $\chi(X) = 1$ for $|X| \leq 4 R_{\xeta} $, $\chi(X)= 0 $ for $|X| \geq R>4 R_{\xeta}$, and with $\chi(|X|)$ decreasing in $|X|$. Our goal is to show that the weighted $L^2$ norm $\int (|\UU(s)|^2 + |P(s)|^2 + |B(s)|^2 ) \chi$ of the solution $\WW = (\UU, P, B)$ of~\eqref{eq:lin_semi2}, vanishes identically for $s\geq 0$. By assumption, we have that $\int (|\UU(0)|^2 + |P(0)|^2+ |B(0)|^2 ) \chi = 0$, so it remains to compute $\frac{d}{ds}$ of this weighted $L^2$ norm using~\eqref{eq:lin_semi2}.

Recall the decomposition of $\cL_E$ from \eqref{eq:lin_Hk} with $k=0$ and $\cL_{E, s} \WW = \cL_{E, \infty} \WW = \cL_E \WW$. Then, we have  
\[
  \cR_{U, k} =\cR_{P, k} = \cR_{B, k}=0, \quad \eu =\ep= \cE_B = 0, \quad \cs = \bc.
  \]
Performing weighted $L^2$ estimates analogous to the ones in the proof of Theorem~\ref{thm:coer_est}, for the transport terms we obtain 
\begin{align}
&\int  ( \cT_U \cdot \UU   + \cT_{P} \, P + \cT_B B ) \chi d X \notag\\
&\quad 
= - \int \B( ( X + \bar \UU) \cdot \f{1}{2} \na \bigl(  |\UU|^2  + P^2 + B^2 \bigr)
   + \bc  \,  \na P \cdot \UU
+ \bc  (\na \cdot  \UU)  \cdot P
\B) \chi d X \notag \\
&\quad = \int \f{1}{2} \na \cdot \B( (X + \bar \UU ) \chi \B) 
(|\UU|^2 + P^2 + B^2)
+ \na ( \bc \chi) \cdot \UU  P d X \notag \\
&\quad \leq \int \f{1}{2} \B(  ( X + \bar \UU) \cdot \na \chi
 + \chi \na \cdot ( X + \bar \UU)  \B) ( |\UU|^2  + P^2 + B^2) 
 +  \bigl(\bc |\na \chi| + \chi |\na \bc| \bigr) |\UU P| d X.
\label{eq:junk:label:1}
\end{align}
We focus on the terms in~\eqref{eq:junk:label:1} that involve $|\na \chi|$. Since $\chi$ is radially symmetric, we get 
\[
(X + \bar \UU) \cdot \na \chi = (\xi + \bar U) \pa_{\xi} \chi.
\]
Using Cauchy--Schwarz, the fact that $ \pa_{\xi} \chi(X) = 0, |X| \leq 4 R_{\xeta}$ and  $\pa_{\xi} \chi \leq 0$ globally, and using that \eqref{eq:rep2}, \eqref{eq:rep22} yield $ \xi + \bar U(\xi) - \bc(\xi) > 0$ for $\xi = |X| > \xi_*$ (hence for $\xi>4 R_{\xeta}$) 
and $\xi + \bu(\xi) \geq 0$ for any $\xi$, we obtain 
\begin{align*}
\tfrac{1}{2} \bigl(  ( X+ \bar \UU) \cdot \na \chi  \bigr) ( |\UU|^2  + P^2 + B^2) 
 + \bc |\na \chi| |\UU P| 
& \leq \tfrac{1}{2}\bigl( (\xi + \bar U ) \pa_{\xi}\chi + \bc |\pa_{\xi} \chi| \bigr)  ( |\UU|^2  + P^2)  \\
&\leq \tfrac{1}{2}(\xi + \bar U  - \bc ) \pa_{\xi}\chi  ( |\UU|^2  + P^2)  \leq 0.
\end{align*}

For the remaining contributions, resulting from the $\chi$-terms in~\eqref{eq:junk:label:1}, and from the $\cD_U, \cS_U,  \cD_P$, $\cS_P, \cD_B, \cS_B$-terms in \eqref{eq:lin_Hk}, in light of~\eqref{eq:dec_U} and  we have that 
\begin{align*}
&\f12 \chi \na \cdot ( X + \bar \UU)  ( |\UU|^2  + P^2 + B^2) +   \chi |\na \bc | |\UU P |
+ 3 \bcv \chi ( |\UU|^2 + P^2 + B^2) 
+ O( \chi |(\UU, P, B)|^2 )
\notag\\
&\leq C \chi ( |\UU|^2  + P^2 + B^2)
\end{align*}
for some sufficiently large $C>0$ (depending on $ r, \bar \UU, \bar P, \bc $). Thus, we obtain 
\[
\f{1}{2} \f{d}{ds} \int (|\UU|^2 + P^2 + B^2 )\chi d X
= 
\int  ( \cL_U \cdot \UU  + \cL_{P} \cdot  P + \cL_B  \cdot B ) \chi d X
\leq C \int (|\UU|^2 + P^2 + B^2 )\chi d X ,
\]
which implies via Gr\"onwall that %
$\int |(\UU, P, B)(s)|^2 \chi d X = 0$ 
for all $s\geq0$. The claim follows.
\end{proof}

\subsubsection{Regularity parameter}\label{sec:para_k}

Let $k_0$ be the parameter from Theorem~\ref{thm:coer_est}. To construct blowup solution in Section \ref{sec:solu}, we fix regularity parameters
$\kkl$ and $\kk$ 
with the special font
\bseq\label{def:kk}
\beq
\kkl = 2 d + 16, \quad 
\kk \geq \max\{k_0, \kkl\}.
\eeq
We simplify the compact operator $\cK_{\kk, \xeta}$ constructed in Proposition \ref{prop:compact}  with $k = \kk$:
\beq 
 \cK_{\kk} : = \cK_{\kk, \xeta}.
\eeq 
\eseq

\subsection{Estimate of \texorpdfstring{$ \cL_{E, s} - \cL_E$}{L E,s - L E}}

We have the following estimates of the error terms $ (\cL_{E, s} - \cL_E) \WW$.
\begin{proposition}\label{prop:error_diff}
Let $ \xeta$ be chosen in \eqref{wg:X_power}. 
For any $\eta \in [\xeta , \etab]$, $ k \geq 0$, and $(\UU, P, B) \in \cX_{\xeta}^{2k}$, we have
\[
  \| (\cL_{E, s} - \cL_E) ( \UU, P, B) \|_{\cX_{\eta}^{2k} }
\les_{k } \rs^{-r + \f{1}{2}( \eta - \xeta)}
\| (\UU, P, B) \|_{\cX_{\xeta}^{ 2 k + 1}} .
\]
\end{proposition}

The above estimates show that the error terms $(\cL_{E, s} - \cL_E) ( \UU, P, B) $ decay faster  than $(\UU, P, B)$. We will gain  regularity from the compact operator $\cK_{k, \eta}$ in Proposition \ref{prop:compact} in the later fixed point argument. 
See \eqref{eq:non_W} and \eqref{eq:non_fix} in Section \ref{sec:solu}.

\begin{proof}
We drop the dependence of $\cL_{E, s}, \cL_E$ on $(\UU, P, B)$ for simplicity. Denote $\WW = (\UU, P, B)$. From \eqref{eq:Euler_profi_modi} and Lemma \ref{lem: cutoff_error}, we have $\bar \sc - \cs =0, \erho = 0,\eu = 0 $ for $|X| < \rs $. For $|X| \geq \rs$, using \eqref{eq:dec_S} and $\bar \sc \les \cs$, we obtain 
\beq\label{eq:L_err_pf1}
|\na^i (\cs -\bar \sc)| \les_i
\la X \ra^{  - i } \rs^{-(r-1)}, 
\quad 
| \na^i ( (\cs^{-1} \bar \sc - 1) \na \bar \sc  )| 
\les_i \la X \ra^{-r- i}
\eeq

Recall the definition of $\cL_E$ and $\cL_{E, s}$ from \eqref{eq:lin_euler2}, \eqref{eq:lin_euler_limit}.  Using \eqref{eq:dec_S} and \eqref{eq:dec_U},  for any $ l \geq 0$ and $|X| \geq \rs $, we obtain 
\[
\bal
   J_{2l}  :=| \na^{2l} ( \cL_{E,s} - \cL_E) | 
  &\les_l  \sum_{ 0 \leq i \leq 2 l} 
  \B( |\na^{2 l-i}(\cs - \bar \sc )| \;  |\na^{i+1} \WW | 
  +
  |\na^{2 l-i + 1 }(\cs - \bar \sc )| \;  |\na^{i} \WW | 
   \\
   & \quad \quad \quad + ( |\na^{2 l-i}(\cs^{-1} \bar \sc \na \bar \sc - \na \bar \sc)|  
 + |\na^{2 l-i} (\ec ,\eu) | 
 )   |\na^i \WW | \B) . \\
   \eal 
   \]
 Using  Lemma \ref{lem: cutoff_error} and \eqref{eq:L_err_pf1}, for any $|X| \geq \rs $, 
 we obtain 
\[
\bal 
   J_{2l}   & \les_l 
 \sum_{ 0 \leq i \leq 2 l} \la X \ra^{-(2 l-i) } 
 |\na^{i+1} \WW| \rs^{-(r-1)}
 + 
 \la X \ra^{-(2 l-i) - 1 }  |\na^{i} \WW| \rs^{-(r-1)}
 +  \la X \ra^{- ( 2 l-i) - r} 
 |\na^i \WW | \\
 & \les_l %
 \sum_{ 0\leq j \leq 2 l+1}
   ( \la X \ra^{- 2 l - 1 + j} \rs^{-(r-1)} + \la X \ra^{- 2 l + j - r}   ) |\na^j \WW| \\
 & \les_l %
   \rs^{- r + 1} \la X \ra^{-1} \sum_{ 0\leq j \leq 2 l+1}
    \la X \ra^{-  2 l + j }   |\na^j \WW|. 
\eal 
\]

Note that for $|X| \geq \rs $ and $\eta \in [\xeta, \etab]$ (hence $\eta - 2 - \etab < 0$), we have
\[
  \vp_{2l}^2 \la X \ra^{\eta - 2} \les_l \la X \ra^{4 l + \eta- 2}
  \les \la X \ra^{4 l + \xeta} \la X \ra^{\eta - \xeta -2}
  \les \la X \ra^{4 l + \xeta} \rs^{\eta - \xeta -2}.
  \]

  Since $\cL_{E, s} - \cL_E = 0$ for $|X| < \rs $, using the above estimates and applying the interpolation in Lemma \ref{lem:interp_wg} and Lemma \ref{lem:norm_equiv} with $\d_1 =1 ,\d_2 =  \xeta $, we obtain 
\[
\bal
   \int    | J_{2l} |^2 \vp_{2l}^2 \la X \ra^{ \eta }  d X 
 & \les_l 
 \rs^{-2 r  + \eta - \xeta}
 \sum_{0\leq j \leq 2 l +1} \int_{|X| \geq  \rs  } \la X \ra^{4 l - 4 l + 2 j +  \xeta } |\na^j \WW|^2  d X 
 \les_l \rs^{-2 r + \eta - \xeta} \| \WW\|_{\cX^{2 l+1}_{ \xeta}}^2.\\
\eal 
\]
Applying the above estimates with $l=0, k$, we complete the proof. 
\end{proof}

\section{Trilinear estimates of the collision operator in \texorpdfstring{$V$}{V}}
\label{sec:trilinear_Q}

In this section, we estimate the $V$-integral of the nonlinear operator $\cN(f, g)$ defined in \eqref{eq:non_nota}. 
In Section \ref{sec:diffu_matrix}, we estimate the diffusion coefficient matrix. 
In Section \ref{sec:non_collision_op}, we decompose and estimate the collision operator.

\subsection{Estimate of diffusion coefficient matrix}\label{sec:diffu_matrix}

Recall $\bpi _v = \f v{|v|} \otimes \f v{|v|}$. Define a matrix-valued function $\bS$ as 
\bseq \label{eqn: defn-matrix}
\begin{align}
    \bS & := \cs ^{\g + 5} \left(\ang \vc ^\g \bpi _{\vc} + \ang \vc ^{\g + 2} (\Id - \bpi _{\vc}) \right).
\end{align}
Since $\bS$ is a positive definite matrix, by definition, we have
\beq 
         \bS^{1/2} =  \cs ^{ \f{\g + 5}{2}} \left(\ang \vc ^{\f{\g}{2} } \bpi _{\vc} + \ang \vc ^{ \f{\g + 2}{2}} (\Id - \bpi _{\vc}) \right), 
\eeq 
and
\beq
    \cs ^{\g + 5} \ang \vc ^\g \Id \preceq \bS \preceq \cs ^{\g + 5} \ang \vc ^{\g + 2} \Id.
\eeq
\eseq 
Here for matrices $\boldsymbol M _1, \boldsymbol M _2$, $\boldsymbol M _1 \preceq \boldsymbol M _2$ means $\boldsymbol M _2 - \boldsymbol M _1$ is nonnegative definite. $\vc$ is an eigenvector of $\bS$ with eigenvalue $\cs ^{\g + 5} \ang \vc ^\g$, and $\vc ^\top$ is a two-dimensional eigenspace of $\bS$ with eigenvalue $\cs ^{\g + 5} \ang \vc ^{\g + 2}$.

\begin{lemma}
    \label{lem: A-pointwise-bound}
    Let $f$ be a scalar-valued function of $V$. %
    Then for any $N \ge 0$ and $\g \geq 0$, there exists a constant $C _N$ such that for every $V \in \R ^3$:
    \bseq\label{eqn: A-pointwise-bound}
    \begin{align}
         - C _N \cs ^{-3} \| \ang \vc ^{-N} f \| _{L ^2 (V)} \bS 
         & \preceq  A [\cMM ^{1/2} f] (V)   
         \preceq C _N \cs ^{-3} \| \ang \vc ^{-N} f \| _{L ^2 (V)} \bS,
         \label{eqn: A-pointwise-bound:a} \\ 
        \f1{C _0} \bS &\preceq A [\cM] (V) \preceq C _0 \bS.    
        \label{eqn: A-pointwise-bound:b}
    \end{align}
    \eseq
    Moreover, for $i=0, 1,2$, and any $j \geq0$ and $N\geq 0$, we have 
    \begin{align}\label{eqn: divA-pointwise}    
                   \left| D^{\leq j } \na_V^i A [\cMM ^{1/2} f] (V)\right| & \les_{j, N} \| \ang \vc ^{-N} D^{\leq j} f \| _{L ^2 (V)} \cs ^{\g + 2-i} \ang \vc ^{\g + 2-i} .  
    \end{align}
\end{lemma}

Using the embedding \eqref{eq:linf_Y}, if $f$ is a function of $X$ and $V$, then we further obtain that $ \| \ang \vc ^{-N} f (X, \cdot) \| _{L ^2 (V)} \leq  \| f (X, \cdot) \| _{L ^2 (V)} \les  \cs^3 \| f \|_{\cYb^{\kkl}}$ for all $N \ge 0$ and every $X \in \R ^3$.

\begin{proof}
\def\bxi{\boldsymbol{\xi}}
    Without loss of generality, assume $\| \ang \vc ^{-N} f \| _{L ^2 (V)} = 1$. 
    For \eqref{eqn: A-pointwise-bound}, it suffices to show that for any vector $\bxi \in \R ^3$, 
    \begin{align}\label{eq: divA-pointwise-pf1}
        \bxi ^\top A [\cMM ^{1/2} f] \bxi &\leq  C_N \cs ^{-3} \bxi ^\top \bS \bxi, &
        \f1{C _0} \bxi ^\top \bS \bxi & \le \bxi ^\top A [\cM] \bxi \le C _0 \bxi ^\top \bS \bxi.
    \end{align}
    By changing $f$ to $-f$, this estimate implies the lower bound in  \eqref{eqn: A-pointwise-bound:a}. Furthermore, it is sufficient to show this for $\bxi = \vc$ and $\bxi \perp \vc$ since $A [\cMM ^{1/2} f]$ and $A [\cM]$ are symmetric. 

\paragraph{Proof of the first bound in \eqref{eq: divA-pointwise-pf1}}
    First, we show the upper bound of $A [\cMM ^{1/2} f]$. Note that 
    \begin{align*}
        A [\cMM ^{1/2} f] & = \int |V - W| ^{\g + 2} (\Id - \bpi _{V - W}) \cMM ^{1/2} (W) f (W) d W \\
        & = \cs ^{\g + 2} \int |\vc - \wc| ^{\g + 2} (\Id - \bpi _{\vc - \wc}) \cMM ^{1/2} (W) f (W) d W \\
        & \preceq C _N \cs ^{\g + 2} \ang \vc ^{\g + 2} \int 
        \ang \wc ^{\g + 2} (\Id - \bpi _{\vc - \wc}) \cMM ^{1/2} (W) |f (W) | d W.
    \end{align*}
    As a consequence, 
    \beq\bal
    \label{eqn: A-pointwise-bound-abs}
        |A [\cMM ^{1/2} f]| &\les \cs ^{\g + 2} \ang \vc ^{\g + 2} \int 
        \ang \wc ^{\g + 2} \cMM ^{1/2} (W) |f (W)| d W \\
        &\les \cs ^{\g + 2} \ang \vc ^{\g + 2} \| \ang \wc ^{\f{\g + 2}2 + N} \cMM ^{1/2} (W) \| _{L ^2 (W)} \| \ang \wc ^{-N} f (W) \| _{L ^2 (W)} \\
        & \les _N \cs ^{\g + 2} \ang \vc ^{\g + 2}.
    \eal\eeq
    Here we used that $\| \ang \vc ^N \cMM ^{1/2} \| _{L ^2 (V)} \le C _N$ for any $N \ge 0$.

    If $\bxi = \vc$ and $|\vc| \geq 1$, since $|\vc|\gtr \ang \vc$, we have 
    \begin{align*}
        \vc ^\top A [\cMM ^{1/2} f] \vc & \les \cs ^{\g + 2} \ang \vc ^{\g + 2} \int 
        \ang \wc ^{\g + 2} \vc ^\top (\Id - \bpi _{\vc - \wc}) \vc \cMM ^{1/2} (W) |f (W) | d W \\
        & = \cs ^{\g + 2} \ang \vc ^{\g + 2} \int 
        \ang \wc ^{\g + 2} \wc ^\top (\Id - \bpi _{\vc - \wc}) \wc \cMM ^{1/2} (W) | f (W) | d W \\
        & \les \cs ^{\g + 2} |\vc|^2  \ang \vc ^{\g } \| \ang \wc ^{ \g + 2 + 2  + N} \cMM ^{1/2} (W) \| _{L ^2 (W)} \| \ang \wc ^{-N} f (W) \| _{L ^2 (W)} \\
        &\les \cs ^{\g + 2} |\vc|^2 \ang \vc ^{\g } \\
        & = \cs ^{-3} \vc ^\top \bS \vc .
    \end{align*}

    If $|\vc| \le 1$, then $\ang \vc \asymp 1$, so $\ang \vc ^{\g + 2} \les \ang \vc ^\g$, and we use the upper bound \eqref{eqn: A-pointwise-bound-abs} to deduce
    \begin{align*}
        \vc ^\top A [\cMM ^{1/2} f] \vc &\le |\vc| ^2 |A [\cMM ^{1/2} f] \les |\vc| ^2 \cs ^{\g + 2} \ang \vc ^{\g + 2} \les \cs ^{\g + 2} |\vc| ^2 \ang \vc ^\g,
    \end{align*}
    which is the same upper bound as above.

    Now we suppose $\bxi \perp \vc$, we also use \eqref{eqn: A-pointwise-bound-abs} to obtain: 
    \begin{align*}
        \bxi ^\top A [\cMM ^{1/2} f] \bxi &\le |\bxi| ^2 |A [\cMM ^{1/2} f]| \les |\bxi| ^2 \cs ^{\g + 2} \ang \vc ^{\g + 2} = \cs ^{-3} \bxi ^\top \bS \bxi.
    \end{align*}

\paragraph{Proof of the second bound in \eqref{eq: divA-pointwise-pf1} }
    The upper bound of $A [\cM]$ follows directly from the first bound because $\cM = \cMM ^{1/2} \cdot \rhos \cMM ^{1/2}$, and $\| \cMM ^{1/2} \| _{L ^2 (V)} = 1$.
    Next, we show the lower bound of $A [\cM]$. By direct computation,  
    \begin{align*}
        \bxi ^\top A [\cM] \bxi = \cs ^{\g + 5} \int |\vc - \wc| ^\g \left(
            |\bxi| ^2 |\vc - \wc| ^2 - |\bxi \cdot (\vc - \wc)| ^2
        \right) \mu (\wc) d \wc.
    \end{align*}
    Again, we only need to show lower bound for the cases $\bxi = \vc$ and $\bxi \perp \vc$.

    If $\bxi \perp \vc$, we have 
    \begin{align*}
        \bxi ^\top A [\cM] \bxi = \cs ^{\g + 5} \int |\vc - \wc| ^\g \left(
            |\bxi| ^2 |\vc - \wc| ^2 - |\bxi \cdot \wc| ^2
        \right) \mu (\wc) d \wc.
    \end{align*}
    When $|\vc| \ge 1$, we restrict the integral in $\{|\wc| \le \f13 |\vc| \}$, and we can bound from below by 
    \begin{align*}
        \bxi ^\top A [\cM] \bxi &\ge \cs ^{\g + 5} \int _{|\wc| \le \f13 |\vc|} |\vc - \wc| ^\g \left(
            |\bxi| ^2 |\vc - \wc| ^2 - |\bxi \cdot \wc| ^2
        \right) \mu (\wc) d \wc \\
        &\ge \cs ^{\g + 5} \int _{|\wc| \le \f13 |\vc|} \left(\f23 |\vc|\right) ^\g \left(
            |\bxi| ^2 \left(\f23 |\vc|\right) ^2 - |\bxi| ^2 \left( \f13 |\vc| \right) ^2
        \right) \mu (\wc) d \wc \\
        & \gtrsim \cs ^{\g + 5} |\vc| ^{\g + 2} |\bxi| ^2. 
    \end{align*}
    When $|\vc| \le 1$, define a cone in $\vc$'s direction
    \begin{align*}
        \mathcal C _{\vc} = \B\{ \wc \in \R ^3: |\wc \cdot \vc| \ge \f{24}{25} |\wc| |\vc| \B\}.
    \end{align*}
    Since $\vc \perp \bxi$, in this cone, we use  Pythagoras' rule to obtain 
    \[ 
    |\wc|^2 = |\bxi|^{-2} | \bxi \cdot \wc|^2 + |\vc|^{-2} |\wc \cdot \vc|^2 
    \geq  |\bxi|^{-2} | \bxi \cdot \wc|^2 + \tf{24^2}{25^2} |\wc|^2
    \implies |\bxi \cdot \wc| \leq \tf{7}{25} |\bxi | \cdot |\wc|.
    \]
    We restrict the integral in the cone intersecting an annulus:
    \begin{align*}
        \bxi ^\top A [\cM] \bxi &\ge \cs ^{\g + 5} \int _{\wc \in \mathcal C _{\vc}, 2 \le |\wc| \le 3} |\vc - \wc| ^\g \left(
            |\bxi| ^2 |\vc - \wc| ^2 - |\bxi \cdot \wc| ^2
        \right) \mu (\wc) d \wc \\
        &\ge \cs ^{\g + 5} \int _{\wc \in \mathcal C _{\vc}, 2 \le |\wc| \le 3} |\vc - \wc| ^\g \left(
            |\bxi| ^2 |\vc - \wc| ^2 - \left(\f7{25} |\bxi| |\wc|\right) ^2 
        \right) \mu (\wc) d \wc \\
        &\ge \cs ^{\g + 5} \int _{\wc \in \mathcal C _{\vc}, 2 \le |\wc| \le 3} \left(
            |\bxi| ^2 - \f{21 ^2}{25 ^2} |\bxi| ^2 
        \right) \mu (\wc) d \wc \\
        &\gtrsim \cs ^{\g + 5} |\bxi| ^2.
    \end{align*}

    Combined, we have shown 
    \begin{align}
        \label{eqn: A[M]-lower-bound-perp}
        \bxi ^\top A [\cM] \bxi \gtrsim \cs ^{\g + 5} \ang \vc ^{\g + 2} |\bxi| ^2 \gtrsim \bxi ^\top \bS \bxi, \qquad \forall \, \bxi \perp \vc.
    \end{align}
    
    Next, when $\bxi = \vc$, we have
    \begin{align*}
        \bxi ^\top A [\cM] \bxi &= \cs ^{\g + 5} \int |\vc - \wc| ^\g \left(
            |\vc| ^2 |\vc - \wc| ^2 - |\vc \cdot (\vc - \wc)| ^2
        \right) \mu (\wc) d \wc \\
        &= \cs ^{\g + 5} \int |\vc - \wc| ^\g \left(
            |\vc| ^2 |\wc| ^2 - |\vc \cdot \wc| ^2
        \right) \mu (\wc) d \wc.
    \end{align*}
    When $|\vc| \ge 1$, we restrict the integral in an annulus but outside the cone:
    \begin{align*}
        \bxi ^\top A [\cM] \bxi &\ge \cs ^{\g + 5} \int _{\wc \notin \mathcal C _{\vc}, \f13 \le |\wc| \le \f12 |\vc|} |\vc - \wc| ^\g \left(
            |\vc| ^2 |\wc| ^2 - |\vc \cdot \wc| ^2
        \right) \mu (\wc) d \wc \\
        &\ge \cs ^{\g + 5} \int _{\wc \notin \mathcal C _{\vc}, \f13 \le |\wc| \le \f12 |\vc|} \left(\f12 |\vc|\right) ^\g \left(\f7{25}\right) ^2 |\vc| ^2 |\wc| ^2 \mu (\wc) d \wc \\
        & \gtrsim \cs ^{\g + 5} |\vc| ^{\g + 2}.
    \end{align*}
    When $|\vc| \le 1$, we integrate in another annulus but still outside the cone:
    \begin{align*}
        \bxi ^\top A [\cM] \bxi &\ge \cs ^{\g + 5} \int _{\wc \notin \mathcal C _{\vc}, 2 \le |\wc| \le 3} |\vc - \wc| ^\g \left(
            |\vc| ^2 |\wc| ^2 - |\vc \cdot \wc| ^2
        \right) \mu (\wc) d \wc \\
        &\ge \cs ^{\g + 5} \int _{\wc \notin \mathcal C _{\vc}, 2 \le |\wc| \le 3} \left(\f7{25}\right) ^2 |\vc| ^2 |\wc| ^2 \mu (\wc) d \wc  \\
        & \gtrsim \cs ^{\g + 5} |\vc| ^{2}.
    \end{align*}
    Combined we have shown 
    \begin{align*}
        \bxi ^\top A [\cM] \bxi \gtrsim \cs ^{\g + 5} \ang \vc ^{\g} |\vc| ^2 = \bxi ^\top \bS \bxi, \qquad \bxi = \vc.
    \end{align*}
    Together with \eqref{eqn: A[M]-lower-bound-perp}, we have completed the proof for the lower bound of $A [\cM]$.

 Finally, we prove \eqref{eqn: divA-pointwise}. For $i \leq 2$ and 
 any multi-indices $\al, \b$ with  $ |\al| +| \b | = j$, since $\g \geq 0$
 and $D^{\al, \b}$ commutes with $\na_V$ by \eqref{eq:deri_wg}, 
 using \eqref{eq:M11/2} and Leibniz rule, we get
     \begin{align*}
        | D^{\al, \b} \na_V^i A [\cMM ^{1/2} f]| &= \B|\int \na_V^i \boldsymbol \Phi (V - W) \cdot 
        D^{\al, \b} (\cMM ^{1/2} f) (W) d W \B| \\
        &\les_j \int |V - W| ^{\g + 2 - i} \cdot \cMM ^{1/2} (W) \cdot |D^{\leq j} f (X, W)| d W
        \\
        &\les_j \cs ^{\g + 2-i} \int |\vc - \wc| ^{\g + 2-i} \cdot \cMM ^{1/2} (W) \cdot |D^{\leq j} f (X, W)| d W
        \\
        &\les_j \cs ^{\g + 2-i} \int \ang \vc ^{\g + 2-i} \ang \wc ^{\g + 2- i} \cdot \cMM ^{1/2} (W) \cdot |D^{\leq j} f (X, W)| d W 
        \\
        &\les_j \cs ^{\g + 2 -i} \ang \vc ^{\g + 2-i} \| \ang \wc ^{\g + 2 - i} \ang \wc ^N \cMM ^{1/2} \| _{L ^2 (W)} \| \ang \wc ^{-N} D^{\leq j} f \| _{L ^2 (W)}.
    \end{align*}
   Using \eqref{eq:cross_pf2:a} and $\cMM = \cs^{-3} \mu(\vc)$, we know $\| \ang \wc ^{\g + N+2} \cMM ^{1/2} \| _{L ^2 (W)} \le C(N, \gamma)$, so \eqref{eqn: divA-pointwise} is proven.
\end{proof}

\begin{cor}
    \label{lem: A-L2-bound}
    Let $f$ be a scalar-valued function of $V$, and $\vec g, \vec h$ be vector-valued functions of $V$. Then for any $N \ge 0$, it holds that  
    \begin{align*}
        \left|\la A [\cMM ^{1/2} f] \vec g, \vec h \ra _V \right| \les _N \cs ^{-3} \| \ang \vc ^{-N} f \| _{L ^2 (V)} \| \bS ^\f12 \vec g \| _{L ^2 (V)} \| \bS ^\f12 \vec h \| _{L ^2 (V)}.
    \end{align*}
\end{cor}

\begin{proof}
Define  $ M =  \bS^{-1/2} A [\cMM ^{1/2} f] \bS^{-1/2}$. Since $\bS$ is a positive definite matrix by 
\eqref{eqn: defn-matrix}, using \eqref{eqn: A-pointwise-bound}, we obtain 
\[
    - C _N \cs ^{-3} \| \ang \vc ^{-N} f \| _{L ^2 (V)} \Id
     \preceq  M 
    \preceq C _N \cs ^{-3} \| \ang \vc ^{-N} f \| _{L ^2 (V)} \Id
    \implies \  | M(V) | \les C _N \cs ^{-3} \| \ang \vc ^{-N} f \| _{L ^2 (V)}.
\]
Thus, using Cauchy--Schwarz inequality, we prove the desired result 
\[
\bal 
        |\la A [\cMM ^{1/2} f] \vec g, \vec h \ra _V | 
        & = |\la M \bS^{1/2} \vec g, \bS^{1/2} \vec h \ra_V| \les \| M \|_{L^{\infty}(V)}
         \| \bS ^\f12 \vec g \| _{L ^2 (V)} \| \bS ^\f12 \vec h \| _{L ^2 (V)} \\
         & \les   C _N \cs ^{-3} \| \ang \vc ^{-N} f \| _{L ^2 (V)}   \| \bS ^\f12 \vec g \| _{L ^2 (V)} \| \bS ^\f12 \vec h \| _{L ^2 (V)}.
\eal 
\]
\end{proof}

Another direct consequence is that $\sigma$ norm is equivalent to the following weighted $H ^1$ norm.

\begin{cor}\label{cor:sig_norm}
Recall the $\s$-norm from \eqref{norm:sig}. Define 
    \beq\label{eq:coer_coe}
        \Lams (s, X, V) = \cs^{\g+3} \la \vc\ra^{\g+2}.
    \eeq
    Then 
    \begin{align}
        \label{eqn: sigma-norm-approx}
	   \| g (s, X, \cdot) \| _{\sigma} ^2 \asymp \int \Lams |g| ^2 + \la \na _V g, \bS \na _V g \ra d V = \| \Lams ^{1/2} g \| _{L ^2 (V)} ^2 + \| \bS ^{1/2} \na _V g \| _{L ^2 (V)} ^2.
    \end{align}
    In particular, we have 
    \bseq\label{eq:sig_norm_L2}
    \begin{align}
        \| \bS ^\f12 \na g \| _{L ^2 (V)} ^2 & \les \| g \| _\sigma ^2, \label{eq:sig_norm_L2:a} \\
        \| \bS ^\f12 (g \cs ^{-1} \vc) \| _{L ^2 (V)} ^2 
        & \les \| \Lams^{1/2} g \| _{L ^2 (V)} ^2 , \label{eq:sig_norm_L2:b} \\
        \| \bS ^\f12 (g \cs ^{-1} \vc) \| _{L ^2 (V)} ^2 
        & \les \| g \| _\sigma ^2. \label{eq:sig_norm_L2:c}
    \end{align}
\eseq
\end{cor}

\begin{proof}
    The proof of \eqref{eqn: sigma-norm-approx} can be found in \cite[ Corollary 1]{guo2002landau}. \eqref{eq:sig_norm_L2:a} is a direct consequences of \eqref{eqn: sigma-norm-approx}. For \eqref{eq:sig_norm_L2:b}, we can compute it as
    \begin{align*}
        \| \bS ^\f12 (g  \cs ^{-1} \vc) \| _{L ^2 (V)} ^2 
        =  \cs ^{\g + 5 - 2} \| \ang \vc ^{\f\g2} g \vc \| _{L ^2 (V)} ^2 
        & \les \| \Lams^{1/2} g \| _{L ^2 (V)} ^2.
    \end{align*}
    We used $\Lams = \cs^{\g+3} \ang \vc^{\g+2}$ defined in \eqref{eq:coer_coe}. Then \eqref{eq:sig_norm_L2:c} is due to \eqref{eqn: sigma-norm-approx}.
\end{proof}

\subsection{Decomposition and estimates of the collision operator}\label{sec:non_collision_op}

Recall %
$\cN(\cdot)$ form \eqref{eq:non_nota}:
\[
    \cN(f, g) = \cMM ^{-1/2} Q(\cMM ^{1/2} f, \cMM ^{1/2} g).
\]
In the next lemmas we will derive the following equivalent formulation for $\cN(f, g)$: 
\beq
\begin{aligned}\label{N(f,g)} 
\cN (f, g) &= \sum _{1\leq  i\leq 6}  \cN _i (f, g), \\
\text{where} \qquad \cN _1 (f, g) &:= \div \left( A [\cMM ^{1/2} f] \na g \right) \\
\cN _2 (f, g) &:= - \div \left(\div A [\cMM ^{1/2} f] g\right) \\
\cN _3 (f, g) &:= - \kp _2 \cs ^{-1} \div \left( A [\cMM ^{1/2} f] \vc g \right) \\
\cN _4 (f, g) &:=  - \kp _2 \cs ^{-1} \vc ^\top A [\cMM ^{1/2} f] \nabla g \\ 
\cN _5 (f, g) &:= \kp _2 ^2 \cs ^{-2} g \vc ^\top A [\cMM ^{1/2} f] \vc \\
\cN _6 (f, g) &:= \kp _2 \cs ^{-1} g \div A [\cMM ^{1/2} f] \cdot \vc.
\end{aligned}
\eeq
In the rest of this section, the divergence operator $\div$ acts on the $V$ variable.

\begin{lemma} \label{lem: nonlinear}
    Let $f, g, h$ be functions of $V$, then for each $\cN _i$ defined in \eqref{N(f,g)} the following holds:
    \begin{align*}
        |\la \cN _i (f, g), h \ra _V| \les _N \cs ^{-3} \| \ang \vc ^{-N} f \| _{L ^2 (V)} \|g \| _{\sigma} \| h \| _{\sigma}.
    \end{align*}
    Recall $\Lams = \cs^{\g+3} \ang \vc^{\g+2}$ from \eqref{eq:coer_coe}. For $\cN_i, 2\leq i\leq 6$, we have the following improved estimates 
    \begin{align*}
        |\la \cN _i (f, g), h \ra _V| \les \cs ^{-3} \| \ang \vc ^{-N} f \| _{L ^2 (V)} \| \Lams^{1/2} g \| _{L ^2 (V)} \| h \| _{\sigma} ,
        \quad  2\leq i \leq 6.
    \end{align*}
\end{lemma}

\begin{proof}
    By the definition of $\cN$ and $Q$, we split the inner product into two parts:
    \begin{align*}
        \la \cN (f, g), h \ra &= \left\la \cMM ^{-1/2} Q (\cMM ^{1/2} f, \cMM ^{1/2} g), h \right\ra \\
        &= -\left\la 
            A [\cMM ^{1/2} f] \na (\cMM ^{1/2} g) - \div A [\cMM ^{1/2} f] \cMM ^{1/2} g, \na (\cMM ^{-1/2} h)
        \right\ra \\
        &= -\left\la 
            A [\cMM ^{1/2} f] \na (\cMM ^{1/2} g), \na (\cMM ^{-1/2} h)
        \right\ra + \left\la \div A [\cMM ^{1/2} f] \cMM ^{1/2} g, \na (\cMM ^{-1/2} h)
        \right\ra. 
    \end{align*}
    We simplify the first inner product as 
    \begin{align*}
        & \left \la A [\cMM ^{1/2} f] \na (\cMM ^{1/2} g), \na (\cMM ^{-1/2} h) \right \ra
        \\
        & \qquad = \left \la   A [\cMM ^{1/2} f] \cMM^{1/2} (\na g+ g \na \log \cMM^{1/2}) , \cMM^{-1/2} (\na h- h\na \log \cMM^{1/2} ) \right \ra\\
        & \qquad = \left \la   A [\cMM ^{1/2} f] (\na g+ g \na \log \cMM^{1/2}) ,  \na h- h\na \log \cMM^{1/2}  \right \ra.
    \end{align*}
    We carry a similar computation for the second inner product:
    \begin{align*}
        \left\la \div A [\cMM ^{1/2} f] \cMM ^{1/2} g, \na (\cMM ^{-1/2} h)
        \right\ra = \left\la \div A [\cMM ^{1/2} f] g, \na h - h \na \log \cMM ^{1/2}
        \right\ra. 
    \end{align*}
    By direct computation, 
	\begin{align}
        \label{eqn: dv_cMM}
		\na _V \cM &= -\cM \cdot \kappa \cs ^{-1} \vc, & \na _V \cMM ^{1/2} = -\cMM ^{1/2} \cdot \kappa _2 \cs ^{-1} \vc,
	\end{align}
    so $\na \log \cMM ^{1/2} = -\kp _2 \cs ^{-1} \vc$. Therefore, the inner product can be expanded as 
    \beq
    \label{eqn: expanding-N(f,g)}
    \begin{aligned}
        \la \cN (f, g), h \ra &= - \left \la A [\cMM ^{1/2} f] \na g, \na h\right \ra - \left \la A [\cMM ^{1/2} f] \na g, h \kp _2 \cs ^{-1} \vc \right \ra \\
        & \quad + \left \la A [\cMM ^{1/2} f] g \kp _2 \cs ^{-1} \vc, \na h \right \ra + \left \la A [\cMM ^{1/2} f] g \kp _2 \cs ^{-1} \vc, h \kp _2 \cs ^{-1} \vc \right \ra \\
        & \quad + \left \la \div A [\cMM ^{1/2} f] g, \na h \right \ra + \left \la \div A [\cMM ^{1/2} f] g, h \kp _2 \cs ^{-1} \vc \right \ra \\
        &= \la \cN _1 (f, g), h \ra + \la \cN _4 (f, g), h \ra + \la \cN _3 (f, g), h \ra + \la \cN _5 (f, g), h \ra \\
        & \quad + \la \cN _2 (f, g), h \ra + \la \cN _6 (f, g), h \ra.
    \end{aligned}
    \eeq
    This equality holds for any $h \in L ^2 (V)$, so we have proven the decomposition \eqref{N(f,g)}.     

    Thanks to Corollary \ref{lem: A-L2-bound}, we only need to bound a $\bS ^\f12$ weighted norm for $\na g$, $\na h$, $g \vc$, and $h \vc$. To this end, we use \eqref{eq:sig_norm_L2} to bound them by the $\s$ norm or the weighted $L ^2$ norm.
    Applying Corollary \ref{lem: A-L2-bound} to \eqref{eqn: expanding-N(f,g)}, we get
    \begin{align}\label{eq:nonlinear_N1_pf1}
        \left|\la \cN _1 (f, g), h \ra\right| = \left|\left \la A [\cMM ^{1/2} f] \na g, \na h \right \ra \right|
        &\lesssim  \cs ^{-3} \|\ang \vc^{-N}f\|_{L ^2 (V)} \| g \| _{\sigma} \| h \| _{\sigma}. 
    \end{align}
    The estimate of $\cN _3$, $\cN _5$ are the same except that we apply the estimate \eqref{eq:sig_norm_L2:b} 
    instead of \eqref{eq:sig_norm_L2:c} since $\cN_{3},\cN_5$ do not involve $\na_V g$: 
    \[
     \left|\la \cN _{3,5} (f, g), h \ra\right| \lesssim  \cs ^{-3} \|\ang \vc^{-N}f\|_{L ^2 (V)} \| \Lams^{1/2} g \| _{L ^2 (V)} \| h \| _{\sigma}. 
    \]

    Next, we derive bounds for $\cN _2$, $\cN _6$. We invoke pointwise bound of $\div A$ in \eqref{eqn: divA-pointwise} with $i=1, j=0$: 
    \begin{align*}
        \left|\la \cN _2 (f, g), h \ra\right| & = \left| \left\la \div A [\cMM ^{1/2} f] g, \na h
        \right\ra \right| \\
        & \le \cs ^{\g + 1} \| \ang \vc ^{-N} f \| _{L ^2 (V)} \int \ang \vc ^{\gamma + 1} |g| \cdot |\na h| d V \\
        & \le \cs ^{\g + 1} \| \ang \vc ^{-N} f \| _{L ^2 (V)} \| \ang \vc ^{\f{\g + 2}2} g \| _{L ^2 (V)} \| \ang \vc ^\f\g2 \na h \| _{L ^2 (V)} \\
        & = \cs ^{-3} \| \ang \vc ^{-N} f \| _{L ^2 (V)} \| \cs ^{\f{\g + 3}2} \ang \vc ^{\f{\g + 2}2} g \| _{L ^2 (V)} \| \cs ^{\f{\g + 5}2} \ang \vc ^\f\g2 \na h \| _{L ^2 (V)} \\
        & \les \cs ^{-3} \| \ang \vc ^{-N} f \| _{L ^2 (V)} \| \Lams^{1/2} g \| _{L ^2 (V)} \| h \| _\sigma.
    \end{align*}
    Similarly, 
    \begin{align*}
        \left|\la \cN _6 (f, g), h \ra\right| & = \left| \left\la \div A [\cMM ^{1/2} f] g, h \kp _2 \cs ^{-1} \vc
        \right\ra \right| \\
        & \le \cs ^{-3} \| \ang \vc ^{-N} f \| _{L ^2 (V)} \| \cs ^{\f{\g + 3}2} \ang \vc ^{\f{\g + 2}2} g \| _{L ^2 (V)} \| \cs ^{\f{\g + 5}2} \ang \vc ^\f\g2 h \kp _2 \cs ^{-1} \vc \| _{L ^2 (V)} \\
        & \le \cs ^{-3} \| \ang \vc ^{-N} f \| _{L ^2 (V)} \| \cs ^{\f{\g + 3}2} \ang \vc ^{\f{\g + 2}2} g \| _{L ^2 (V)} \| \cs ^{\f{\g + 3}2} \ang \vc ^{\f{\g + 2}2} h \| _{L ^2 (V)} \\
        & \les \cs ^{-3} \| \ang \vc ^{-N} f \| _{L ^2 (V)} \| \Lams^{1/2} g \| _{L ^2 (V)} \| h \| _\sigma.
    \end{align*}

Finally for $\cN_4$, using integration by parts in 
$\na_V$, we obtain
\[
\bal 
& \la \cN_4(f, g), h \ra
 =  - \left \la A [\cMM ^{1/2} f] \na g, h \kp _2 \cs ^{-1} \vc \right \ra   \\ 
& =  \B\la \div A [\cMM ^{1/2} f] g, h \kp _2 \cs ^{-1} \vc \B\ra
+ \B \la   A [\cMM ^{1/2} f] g,  \na_V h \otimes  \kp_2 \cs^{-1}  \vc  \B\ra
+  \B \la   A [\cMM ^{1/2} f] g,   h \kp_2 \cs^{-1} \na_V \vc  \B\ra
\eal 
 \]
where in the last term, we have $\na_V \vc = \cs^{-1} \Id$. The first two terms are the same as the above $\cN_6$-term, 
$\cN_3$-term. For the last term, using $\Lams = \cs^{\g+3} \ang \vc^{\g+2}$ and  Lemma \ref{lem: A-pointwise-bound}, we bound 
\[
\bal 
\B|  \B \la   A [\cMM ^{1/2} f] g,   h \kp_2 \cs^{-1} \na_V \vc  \B\ra \B| 
& \les \cs^{\g+2}  \| \ang \vc^{-N} f \|_{L ^2 (V)} \int \ang \vc^{\g+2} \cs^{-2}  |g| \cdot |h| d V  \\ 
& \les \cs^{\g+2 - 2 - \g - 3} \| \ang \vc^{-N} f \|_{L ^2 (V)} \| \Lams^{1/2} g\|_{L ^2 (V)}   \| h\|_{\s}  .
\eal 
\]
Combining the above estimate and the estimates of $\cN_3, \cN_6$, we prove 
\[
| \la \cN_4(f, g), h \ra| \les \cs ^{-3} \| \ang \vc^{-N} f \|_{L ^2 (V)} \| \Lams^{1/2} g\|_{L ^2 (V)}   \| h\|_{\s} .
\]
Since $ \| \Lams^{1/2} g\|_{L ^2 (V)} \leq \| g \|_{\s}$, we complete the proof.    
\end{proof}

Next, we estimate commutator with derivatives.

\begin{lemma}
    \label{lem: N-derivative-bound}
    Let $f, g, h$ be functions of $V$. For each $1 \le i \le 6$ and any multi-indices $\al, \b$, we have
    \begin{align*}
        |\la D^{\al, \b} \cN _i (f, g) - \cN _i (f, D ^{\al, \b} g), h \ra _V| \les _{\al, \b} \sum _{\substack{\al _1 + \al _2 \preceq \al \\ \b _1 + \b _2 \preceq \b \\ (\al _2, \b _2) \prec (\al, \b)}} \cs ^{-3} \|  D^{\alpha_1, \beta_1} f \| _{L ^2 (V)} \| D^{\alpha_2, \beta_2}g \| _{\sigma} \| h \| _{\sigma}, \\
        |\la D^{\al, \b} \cN _i (f, g) - \cN _i (D ^{\al, \b} f, g), h \ra _V| \les _{\al, \b} \sum _{\substack{\al _1 + \al _2 \preceq \al \\ \b _1 + \b _2 \preceq \b \\ (\al _1, \b _1) \prec (\al, \b)}} \cs ^{-3} \|  D^{\alpha_1, \beta_1} f \| _{L ^2 (V)} \| D^{\alpha_2, \beta_2}g \| _{\sigma} \| h \| _{\sigma}.
    \end{align*}
\end{lemma}

\begin{proof}
    By Leibniz rule and \eqref{eq:M11/2} we have
    \begin{align*}
        \B| D ^{\al, \b} (\cMM ^{1/2} f) \B| &= \B| \sum _{\substack{\al' \preceq \al\\ \b' \preceq \b}} C ^\al _{\al'} \cdot C ^\b _{\b'} \cdot D ^{\al - \al', \b - \b'} \cMM ^{1/2} \cdot D ^{\al', \b'} f \B|  \les  \cMM ^{1/2} \sum _{\substack{\al' \preceq \al\\ \b' \preceq \b}} \ang \vc ^{|\b| + 2 |\al|} |D ^{\al', \b'} f|.
    \end{align*}
    We can write 
\[
\bal
        D ^{\al, \b} \cN _1 (f, g) &= \sum _{\substack{\al _2 \preceq \al , \,  \b _2 \preceq \b}} C ^\al _{\al _2} \cdot C ^\b _{\b _2} \div \left( A [D ^{\al - \al _2, \b - \b _2} (\cMM ^{1/2} f)] \na D ^{\al _2, \b _2} g \right) \\
        &= \sum _{\substack{\al _2 \preceq \al , \, \b _2 \preceq \b}} C ^\al _{\al _2} \cdot C ^\b _{\b _2} \cN _1 \B( 
            \f{D ^{\al - \al _2, \b - \b _2} (\cMM ^{1/2} f)}{\cMM ^{1/2}},
            D ^{\al _2, \b _2} g
        \B),
\eal
\]
    where by Lemma \ref{lem: nonlinear} and \eqref{eq:M11/2},
    \begin{align*}
       & \B| \B\la \cN _1 \B( 
            \f{D ^{\al - \al _2, \b - \b _2} (\cMM ^{1/2} f)}{\cMM ^{1/2}},
            D ^{\al _2, \b _2} g
        \B), h \B\ra _V \B| \\
        & \qquad \les \cs ^{-3}  \B\| \ang \vc ^{-N} \f{D ^{\al - \al _2, \b - \b _2} (\cMM ^{1/2} f)}{\cMM ^{1/2}} \B\| _{L ^2 (V)} \| D ^{\al _2, \b _2} g \| _\sigma \| h \| _\sigma \\
        & \qquad \les \cs ^{-3} \| D ^{\preceq (\al - \al _2, \b - \b _2)} f \| _{L ^2 (V)} \| D ^{\al _2, \b _2} g \| _\sigma \| h \| _\sigma. 
    \end{align*}
    The bounds for $\cN _1$ are proven after taking summations. Similarly, 
    \begin{align*}
        D ^{\al, \b} \cN _2 (f, g) &= \sum _{\substack{\al _2 \preceq \al ,\, \b _2 \preceq \b}} C ^\al _{\al _2} \cdot C ^\b _{\b _2} \cN _2 \B( 
            \f{D ^{\al - \al _2, \b - \b _2} (\cMM ^{1/2} f)}{\cMM ^{1/2}},
            D ^{\al _2, \b _2} g
        \B) , 
    \end{align*}
    and the bound for $\cN _2$ follows.

    For $\cN _3$, since $A[\cMM^{1/2} f] \vc = A[ \cMM^{1/2} f \vc]$, we rewrite it as 
    \[
  \cN_3(f, g) := - \kp _2 \cs ^{-1} \div \big( A [\cMM ^{1/2} f] \vc g \big)  
   = -  \sum_{1\leq i\leq 3} \kp _2 \cs ^{-1} \div \big( A [\cMM ^{1/2} f \vc_i] \ee _i g \big) 
  =   \sum_{1\leq i\leq 3} \cN_{3, i}(f \vc_i , g) ,
     \]
    where $\ee _i$ is the standard basis in $\R^3$ and  we define
   \begin{align}\label{eq:non_decomp_N3}
        \cN _{3, i} ( F, G) := - \kp _2 \cs ^{-1} \div (A [\cMM ^{1/2} F] \ee _i \cdot G),
    \end{align}
which is $\cN _3$ with the $\vc$ factor replaced by the basis vector $\ee _i$. 
For each $ i=1,2,3$, we also have derivatives hitting $\kp _2 \cs ^{-1}$ when applying Leibniz rule: 
\begin{align*}
    D ^{\al, \b} \cN _{3,i} (f \vc_i, g) &=  - \sum _{\substack{\al _2  + \al _4 \preceq \al \\ \b _2  \preceq \b}} C ^{\al} _{\al _2,  \al _4} \cdot C ^\b _{\b _2  } 
        \cdot   \kp _2  D ^{\al _4, 0} (\cs ^{-1})  \cdot 
         \div \B( A [D ^{\preceq (\al, \b)} (\cMM ^{1/2} f \vc_i)] \ee _i \cdot  D ^{\al _2, \b _2} g \B)  \\ 
         & =  \sum _{\substack{\al _2  + \al _4 \preceq \al \\ \b _2  \preceq \b}} C ^{\al} _{\al _2,  \al _4} \cdot C ^\b _{\b _2  } 
         \f{ D ^{\al _4, 0} (\cs ^{-1}) }{\cs^{-1}} \cN_{3, i}\B( \f{  D ^{\preceq (\al, \b)}( \cMM^{1/2} f \vc_i)  }{\cMM^{1/2}}, D^{\al_2, \b_2} g \B) .
\end{align*}

    Note that similar to \eqref{eq:sig_norm_L2:c}, we also have 
    \begin{align*}
        \| \bS ^\f12 (g \kp _2 \cs ^{-1} \ee _i) \| _{L ^2 (V)} ^2 & \le \kp _2 ^2 \cs ^{\g + 3} \| \ang \vc ^{\f{\g + 2}2} g \ee _i \| _{L ^2 (V)} ^2 \les \| g \| _\sigma ^2.
    \end{align*}
    Following the proof of Lemma \ref{lem: nonlinear}, one can show $\cN _{3, i}$ defined in \eqref{eq:non_decomp_N3} enjoys the same bound of $\cN _i$ in Lemma \ref{lem: nonlinear}. %
    Since $|D ^{\al _4, 0} (\cs ^{-1}) | \les \cs ^{-1}$ from Lemma \ref{lem:example_fl} (4), we conclude the proof for $\cN _3$ in the same way as $\cN _1$. The proof for $\cN _4$, $\cN _5$, and $\cN _6$ are similar.
    \end{proof}

\section{Linear stability estimates: microscopic part}
\label{sec:lin_micro}

In this section, we derive the linear stability estimates for the equation that governs the evolution of the microscopic part of the perturbation. 

Recall the $\s$-norm in $V$ from \eqref{norm:sig} and its equivalent formulation in \eqref{eqn: sigma-norm-approx}, 
 and recall the $\cY$-norm in $(X, V)$ from \eqref{norm:Y}.

The following lemma contains several interpolation inequalities. These inequalities will be used later to control the free transport term by the collision term and the self-similar scaling field.

\begin{lemma}
    \label{lem:V/X}
    Recall $\Lams = \cs^{\g+3} \ang \vc^{\g+2}$ defined in \eqref{eq:coer_coe}. %
    Assume $0 < r < 3 - \sqrt 3$. 
     For any $1 < \g \leq 2$, we have
    \bseq\label{eqn: interpolation}
    \begin{align}
        \label{eqn: interpolation-v3-holder}
        \ang X ^{-r} \ang \vc ^3  
        \les \ang X ^{-1} \cs \ang \vc ^3 &\les 
        \Lams ^\f3{\g + 2} .
\end{align}
For any $0\leq \g \leq 2$, we have 
\begin{align}
        \label{eqn: interpolation-v-holder}
        \ang X ^{-1} + \ang X ^{-1} |V| + \ang X ^{-1} \cs \ang \vc & \les 
        \Lams ^\f12.
    \end{align}
    \eseq
\end{lemma}

\begin{proof}
    From the bounds $\cs \gtrsim \ang X ^{-r + 1}$ in \eqref{eq:dec_S}  
    we know $
    \ang X ^{-r} \ang \vc ^3 \les \ang X ^{-1} \cs \ang \vc ^3$.
    By direct computation, 
    \begin{align*}
        \Lams ^\f3{\g + 2} = \cs ^\f{3 (\g + 3)}{\g + 2} \ang \vc ^3.
    \end{align*}
    To prove \eqref{eqn: interpolation-v3-holder}, it suffices to show $\ang X ^{-1} \les \cs ^{\f{3 (\g + 3)}{\g + 2} - 1} = \cs ^{2 + \f3{\g + 2}}$, which holds provided
    \begin{align*}
        \left(2 + \f3{\g + 2}\right) (r - 1) \le 1.
    \end{align*}
    Clearly, it holds for all $\g > 1$ and $r < \f43$. 

    Using the relation $\vc = \f{V - \bu}{\cs}$, the estimates \eqref{eq:dec_S}, and \eqref{eq:dec_U}, we obtain 
    \beq\label{eq:V_est}
        |V| = | \bu + \cs \vc | \les \cs (1 + |\vc|) \les \cs \la \vc \ra.
    \eeq
    Recall $\cs$ is bounded, so $\ang X ^{-1} + \ang X ^{-1} |V| + \ang X ^{-1} \cs \ang \vc \les \ang X ^{-1} \ang \vc$.    
    Using 
    $$\Lams ^\f12 = \cs ^{\f{\g + 3}2} \ang \vc ^{\f{\g + 2}2} \gtrsim \ang X ^{-(r - 1) \f{\g + 3}2} \ang \vc,$$ 
    for any $\g \geq 0$, to prove \eqref{eqn: interpolation-v-holder}
    we only need 
    \begin{align*}
        \f{\g + 3}2 (r - 1) \le 1,
    \end{align*}
    which holds for any $\g \leq 2$ and $r < \f75$.   %
\end{proof}

\subsection{The micro equation and the linear estimate}

In this section we derive the equation for $\tFm$ and $D ^{\al, \b} \tFm$. We also present the $H ^k$ estimate in Theorem \ref{thm: micro-Hk-main}. The proof of this theorem occupies the following subsections. 

Define 
\begin{align}
    \label{eqn: dcm}
    \dcm &= \frac12 (\pa _s + \cT) \log \cMM + \f32 \bcv, & 
    \td d _\cM &= -\frac12 V \cdot \na _X \log \cMM.
\end{align}
where $\cT g = (V \cdot \na_X + \bcx X \cdot \na_X + \bcv V \cdot \na_V) g$ was defined in \eqref{def:T-LM}. The estimates on $\dcm$ and $\td d _\cM$ can be found in Lemma \ref{lem: commu-derivative-projection} (4).

\begin{lemma}\label{lem:Deri-Fm}
    $\tFm = \cP _m \td F$ satisfies
    \begin{align}
    	\label{eqn: micro-part}
    	\left(\pa _s + \cT + \dcm - \f32 \bcv\right) \tFm &+ \cP _m [(V \cdot \na _X + 2 \dcm + \td d _\cM) \tFM] %
        - \cP _M [(V \cdot \na _X - 2 \dcm - \td d _\cM) \tFm]  \notag \\ 
        & = \frac1{\es} \cL _{\cM} (\tFm) 
        - \cP _m [\cMM ^{-1/2} \eM] + \f1\es \cN (\td F, \td F).
    \end{align}
    Here $\eM = (\pa _s + \cT) \cM$ is defined in \eqref{eq:error0}, $
        \cN (f, g) := \cMM ^{-1/2} Q (\cMM ^{1/2} f, \cMM ^{1/2} g)
    $ is defined in \eqref{eq:non_nota}, and $\cL _\cM = \cN (\rhos \cMM ^{1/2}, \cdot) + \cN (\cdot, \rhos \cMM ^{1/2})$ is defined in \eqref{def:T-LM}:
    \begin{align*}
		\cL _{\cM} g & = \cMM ^{-1/2} \left[
		  Q (\cM, \cMM ^{1/2} g) + Q(\cMM ^{1/2} g, \cM)
		\right].
    \end{align*}
\end{lemma}

\begin{proof}
We derive the linearized equation of $\td F$ from \eqref{eq:lin} by dividing $\cMM^{1/2}$:
\beq\label{eq:lin_micro}
\pa_s \td F + \cT \td F + \f12 (\pa _s + \cT) \log \cMM  \cdot \td F = \f{1}{\es} \cL_{\cM} \td F - \cMM ^{-1/2} \eM + \f1\es \cMM ^{-1/2} Q (\cMM ^{1/2} \td F, \cMM ^{1/2} \td F).
\eeq
where $\cL_{\cM}$ is defined in \eqref{def:lin_op}, and the $\log \cMM$ term comes from
\begin{align*}
	\cMM^{-1/2} (\pa_s + \cT) \cMM ^{1/2} = \frac12 (\pa _s + \cT) \log \cMM.
\end{align*}
We write $\f12 (\pa _s + \cT) \log \cMM = \dcm - \f32 \bcv$.
Using the fact that $\cP _m$ commutes with $\cL _{\cM}$, we now project \eqref{eq:lin_micro} to the microscopic part: 
\begin{align}
    \label{eq:lin_micro_projection}
    \cP _m \B[\B(\pa _s + \cT + \dcm - \f32 \bcv\B) \td F\B] = \frac1{\es} \cL _{\cM} (\tFm) - \cP _m [\cMM ^{-1/2} \eM] + \f1\es \cN (\td F, \td F). 
\end{align}
Here we used $\cN (f, g) \perp \Phi _i$ in $L ^2$, so $\cP _m \cN (f, g) = \cN (f, g)$. Use Lemma \ref{lem: commu-derivative-projection} (1), 
\begin{align*}
    \cP _m [(\pa _s + \cT) \td F] = (\pa _s + \cT) \tFm &+ \cP _m [(V \cdot \na _X + \dcm + \td d _\cM) \tFM] \\
    & - \cP _M [(V \cdot \na _X - \dcm - \td d _\cM) \tFm].
\end{align*}
Combine this with \eqref{eq:lin_micro_projection} and 
\begin{align*}
    \cP _m (d _\cM \td F) = d _\cM \tFm + d _\cM \tFM - \cP _M (d _\cM \td F) = d _\cM \tFm + \cP _m [d _\cM \tFM] - \cP _M [d _\cM \tFm],     
\end{align*}
we conclude \eqref{eqn: micro-part}.
\end{proof}

We introduce a linear micro operator $\Lmic$:
\bseq \label{eqn: defn-Lmic}
\begin{align}
        \Lmic g =  \f1\es \cL _\cM g - \left(\cT + \dcm - \f32 \bcv \right) g+ \cP _M [(V \cdot \na _X - 2 \dcm - \td d _\cM) g] .
\end{align}
Note that the operator $\Lmic$ depends on $s$. For simplicity of notation, we suppress this dependence.
Then \eqref{eqn: micro-part} can be written in the following form 
\begin{align}
        (\pa _s - \Lmic) \tFm = - \cP _m [(V \cdot \na _X + 2 \dcm + \td d _\cM) \tFM] 
         + \f1\es \cN (\td F, \td F) - \cP _m [\cMM ^{-1/2} \eM]. 
\end{align}
\eseq

One main goal of this section is to prove estimates for $\pa _s - \Lmic$, $\cP _m [(V \cdot \na _X + 2 \dcm + \td d _\cM) \tFM]$ and $\cP _m [\cMM ^{-1/2} \eM]$. The nonlinear term $\cN$ will be estimated in Section \ref{sec:non_Q}.

\begin{theorem}
    \label{thm: micro-Hk-main}
    Suppose that $1 < \g \le 2$.  For every $k \ge 0$, $\eta \leq \etab$, if $\nu \leq \nu _k$ for the $\nu_k$ determined by Lemma \ref{lem:L_m_Dalphabeta}, then for every $g = \cP _m g$, we have the following estimates of the operator $\Lmic$:
    \bseq\label{eq:micro_lin_est}
    \begin{align}
        \left\la (\pa _s - \Lmic) g, g \right\ra _{\cYe ^k} &\ge \f12  \f{d}{d s} \| g \| _{\cYe ^k} ^2 + \f{\bcx}2 (\etab - \eta) \| g \| _{\cYe ^k} ^2 + \frac{\cgam}{4 \es} \| g \| _{\cYE ^k} ^2  \nonumber
        \\
        & \qquad - C _{k, \eta} \| g \| _{\cYe ^k} ^\f{2 (\g - 1)}{\g + 2} \| g \| _{\cYE ^k} ^\f6{\g + 2} - C _{k, \eta} \| g \| _{\cYe ^k} \| g \| _{\cYE ^k}
        \label{eq:micro_lin_est:lin} \\ 
        & \ge \f12  \f{d}{d s} \| g \| _{\cYe ^k} ^2 + \left( 
            2 \lam_{\eta} - C _{k, \eta} \es  
        \right)
        \| g \| _{\cYe ^k} ^2 + \frac{\cgam}{ 6 \es} \| g \| _{\cYE ^k} ^2 \label{eq:micro_lin_est:lin2} ,
    \end{align}
    where $\lam_{\eta} = \f{\bcx}{4} (\etab - \eta)$ is defined in \eqref{eq:lam_eta}.
    The remaining operators in \eqref{eqn: defn-Lmic} satisfy the following estimates:
    \begin{align}
        \left| \la \cP _m [(2 \dcm + \td d _\cM) \tFM], \tFm \ra_{ \cYe ^k} \right|  
        & \le C _{k} \| \tFm \|_{\cYE ^k}  \| \tFM \|_{\cYe ^k} , 
        \label{eq:micro_lin_est:cross_d}
        \\
        - \la \cP _m [(V \cdot \na _X) \tFM], \tFm \ra _{\cYe ^k}  
        & =  \cO _{k} \left( \nu ^{-\f12}
        \| \tFm \|_{\cYE ^k}  \| \tFM \|_{\cYe ^k}
        \right) 
        \label{eq:micro_lin_est:cross}   \\
        & \qquad - \underbrace{\sum _{|\al| = k}     \f{ |\al|! }{ \al !} \iint (V \cdot \na _X) D ^{\al, 0} \tFM \cdot D ^{\al, 0} \tFm \ang X ^\eta d V d X} _{\text{\upshape cross terms}}, \notag \\
        |\la \cP _m [\cMM ^{-1/2} \eM] , \tFm   \ra_{\cYe ^k} | & 
        \les _k \| \tFm \|_{\cYE ^k}  .   \label{eq:micro_lin_est:err} 
    \end{align}
    \eseq
    In particular, let $\kk$ be the regularity parameter chosen in \eqref{def:kk}. We choose $\nu = \nu_{\kk}$ in the $\cY$-norm \eqref{norm:Y}. 
\end{theorem}

\vspace{0.5cm}
Theorem \ref{thm: micro-Hk-main} is proven by estimating $\cL _\cM$, $\cT$, $\dcm$, and $\cP _M$ part of $\Lmic$ separately in the following subsections. Proof of \eqref{eq:micro_lin_est:lin} and \eqref{eq:micro_lin_est:lin2} is provided at the end of Section \ref{sec: projection-m-free-transport}. The proof of \eqref{eq:micro_lin_est:cross_d} and \eqref{eq:micro_lin_est:cross} can be found in Section \ref{sec: projection-free-transport}. Finally, the proof of \eqref{eq:micro_lin_est:err} is in Section \ref{sec: micro_error}. The remaining cross terms in \eqref{eq:micro_lin_est:cross} will be estimated in Section \ref{sec:EE_top},
and the nonlinear term will be estimated in Section \ref{sec:non_Q}. 

 \vspace{0.5cm}

{\renewcommand{\tFm}g
Before we start with the calculation, let's introduce the term $D^{\al, \b} (\pa_s - \Lmic) \tFm$ and $- D^{\al, \b} 
\cP _m [(V \cdot \na _X + 2 \dcm + \td d _\cM)  \tFM]$. We start with $D^{\al, \b} (\pa_s - \Lmic) \tFm$ that can be rewritten as  
    \beq
    \bal \label{eqn: lin_micro_dab}
     D^{\al, \b} ( -\pa_s + \Lmic) \tFm   &= ( - \pa_s + \Lmic) D^{\al,  \b} \tFm  +   \EEs _{\al, \b },
    \eal 
    \eeq 
where $\EEs_{\al, \b}$ is computed as 
\[
\EEs_{\al, \b} := \f1\es h _0 -h _1  + h _4 + h _5, 
\]
and we denote
\begin{align*}
    h _0 &: = D ^{\al, \b} \cL _\cM (\tFm) - \cL _\cM (D ^{\al, \b} \tFm), \\
    h _1 &: = D ^{\al, \b} (\pa _s + \cT + \dcm) \tFm 
    - (\pa _s + \cT + \dcm) D ^{\al, \b} \tFm, \\
    h _4 &: = D ^{\al, \b} \cP _M [V \cdot \na _X \tFm] - \cP _M [(V \cdot \na _X) D ^{\al, \b} \tFm], \\
    h _5 &:= -D ^{\al, \b} \cP _M [(2 \dcm + \td d _\cM) \tFm] + \cP _M[ (2 \dcm + \td d _\cM) D ^{\al, \b} \tFm].
\end{align*}
The other term writes
\begin{align}\label{eq:Pm_tran_FM}
     - D^{\al, \b}  \cP _m [(V \cdot \na _X + 2 \dcm + \td d _\cM)  \tFM]
     & = -  \cP _m [(V \cdot \na _X + 2 \dcm + \td d _\cM) D^{\al, \b} \tFM]  - h_2 - h_3, 
\end{align}
where 
\begin{align*}
    h _2 &:= D ^{\al, \b} \cP _m [V \cdot \na _X \tFM] - \cP _m [(V \cdot \na _X) D ^{\al, \b} \tFM], \\
    h _3 &:= D ^{\al, \b} \cP _m [(2 \dcm + \td d _\cM) \tFM]
    - \cP _m [(2 \dcm + \td d _\cM) D ^{\al, \b} \tFM]. 
\end{align*}

}

\subsection{Linear collision operator \texorpdfstring{$\cL _\cM$}{L M} estimate}
\label{sec: linear-collision}

In this subsection, we provide coercivity estimate for the first part of $\Lmic$, i.e. the linearized collision operator $\cL _\cM$.
We first recall the spectral gap of the linearized collision operator.

\begin{lemma}
    \label{lem: spectral-gap}
	Let $g$ be a function of $V$. There exists a constant $\cgam > 0$ such that
	\begin{align*}
		\int \cL _\cM g \cdot g d V & \le - \cgam \| \cP _m g \| _{\sigma} ^2.
	\end{align*}
\end{lemma}

\begin{proof}
Recall $\vc = \f{V- \bu}{\cs}$. We make a change of variable
	\begin{align*}
		\mathring g (\vc) = g (\cs \vc + \bu) = g (V), \qquad
		\mathring \cMM (\vc) = \cMM (\cs \vc + \bu) = \cMM (V).
	\end{align*}
Denote $\wc = \frac{W - \bu}{\cs}$ and $\boldsymbol \Phi (v) = \frac1{8 \pi} (\Id - \frac{v \otimes v}{|v|^2}) |v|^{\gamma + 2}$. 
Recall the Gaussian $\mu(\cdot)$ from \eqref{eq:gauss}.  Then
	\begin{align*}
		Q (\cM, \cMM ^{1/2} g) (V) & = \div _{V} \int _{\R ^3} \boldsymbol \Phi (V - W) \na _{V - W} \left[
			\cM (W) \cMM ^{1/2} (V) g (V)
		\right] d W \\
		& = \div _{V} \int _{\R ^3} \cs ^{\gamma + 2} \boldsymbol \Phi (\vc - \wc) \na _{V - W} \left[
			\mu (\wc) \mathring \cMM ^{1/2} (\vc) \mathring g (\vc)
		\right] d W \\
		& = \cs ^{\gamma + 3} \div _{\vc} \int _{\R ^3} \boldsymbol \Phi (\vc - \wc) \na _{\vc - \wc} \left[
			\mu (\wc) \mathring \cMM ^{1/2} (\vc) \mathring g (\vc)
		\right] d \wc \\
		& = \cs ^{\gamma + 3} Q (\mu, \mathring \cMM ^{1/2} \mathring g) (\vc) .
	\end{align*}
 By symmetry, we obtain 
	\begin{align*}
		\cL _{\cM} (g) (V) & = \cMM ^{-1/2} (V) \left(
		Q (\cM, \cMM ^{1/2} g) (V) + Q(\cMM ^{1/2} g, \cM) (V)
		\right) \\
		& = \cs ^{\gamma + 3} \mathring \cMM ^{-1/2} (\vc) \left(
		Q (\mu, \mathring \cMM ^{1/2} \mathring g) (\vc) + Q (\mathring \cMM ^{1/2} \mathring g, \mu) (\vc)
		\right) \\
		& = \cs ^{\gamma + 3} \mu ^{-1/2} (\vc) \left(
		Q (\mu, \mu ^{1/2} \mathring g) (\vc) + Q (\mu ^{1/2} \mathring g, \mu) (\vc)
		\right) \\
		& = \cs ^{\gamma + 3} \cL _{\mu} \mathring g (\vc),
	\end{align*}
where we introduce the linear operator $\cL_{\mu}$ similar to $\cL_{\cM}$ from \eqref{def:T-LM}.

    Multiplying $g$, integrating in $V$, and then performing a change of variable $V \to \vc$, we yield     
	\begin{align*}
		\int \cL _{\cM} g \cdot g d V & = \cs ^{\gamma + 3} \int \cL _\mu \mathring g (\vc) \cdot \mathring g (\vc) d V  = \cs ^{\gamma + 6} \int \cL _\mu \mathring g (\vc) \cdot \mathring g (\vc) d \vc  .
    \end{align*}
Applying the coercivity estimates of $\cL_{\mu}$ \cite[Lemma 5]{guo2002landau} and then changing $\vc \to V$, we obtain 
\begin{align*}
			\int \cL _{\cM} g \cdot g d V & 	 \le - \cs ^{\gamma + 6} \cdot \cgam \int A [\mu (\vc) ] \na _{\vc} \mathring \cP _m \mathring g \cdot \na _{\vc} \mathring \cP _m \mathring g + A [\mu \vc \otimes \vc] (\mathring \cP _m \mathring g (\vc)) ^2 d \vc \\
        &\leq - \cs ^{\gamma + 5}  \cgam \int _{\R ^3} A [\mu (\vc)] \na _V \cP _m g \cdot \na _V \cP _m g d V - \cs ^{\gamma + 3} \cgam \int _{\R ^3} A [\mu \vc \otimes \vc] (\cP _m g) ^2 d V \\
        &= - \cgam \| \cP _m g \| _\sigma ^2.
	\end{align*}
	The constant $\cgam> 0 $ depends on $\gamma \ge -3$ only. Here, $\mathring \cP _m = \operatorname{Id} - \mathring \cP _M$ is the micro projection in the $\vc$ variable.
\end{proof}

\vspace{0.5cm}

We now show the $H ^k$ estimate for the linearized collision operator. %

\begin{lemma} \label{lem:L_m_Dalphabeta}
    For $k \ge 0$, $\eta \in \R$, it holds that 
    \begin{align*}        
        \left\langle \cL _\cM g, g \right \rangle _{\cYe ^k} &= 
        \sum _{|\al| + |\b| \le k} \nu ^{|\al| + |\b| - k} 
            \f{ |\al|! }{ \al !} \iint \ang X ^\eta D ^{\al, \b} \cL _\cM g \cdot D ^{\al, \b} g d V d X
        \\
        &\le -\frac{\cgam}3 \| g \| _{\cYE ^k} ^2 + \cgam \| \cP _M g \| _{\cYE ^k} ^2 + \one _{k > 0} C _k \| g \| _{\cYE ^{k - 1}} ^2.
    \end{align*}
    In particular, if $\nu \leq \nu _k \le 1$ for some $\nu _k > 0$ then 
    \begin{align*}        
        \left\langle \cL _\cM g, g \right \rangle _{\cYe ^k} \le -\frac{\cgam}4 \| g \| _{\cYE ^k} ^2 + \cgam \| \cP _M g \| _{\cYE ^k} ^2.
    \end{align*}
\end{lemma}

Note that this coercivity estimate is only used for stability analysis in Section \ref{sec:solu} and is \textit{not} used 
to prove local well-posedness in Section \ref{app:LWP_fixed}.

\begin{proof}
    The case $k = 0$ is a direct consequence of Lemma \ref{lem: spectral-gap} and the elementary inequality 
    \begin{align*}
        -\| \cP _m g \| _{\sigma} ^2 \le -\frac12 \| g \| _{\sigma} ^2 + \| \cP _M g \| _\sigma ^2. 
    \end{align*}
    For $k > 0$, we need to prove there exists a constant $C _k > 0$ such that for any multi-index $\al$, $\b$ with $|\al| + |\b| \leq k$, the following holds:
	\begin{align}
        \label{eqn: Dab-Lm-g}
		\la D ^{\al, \b} \cL _\cM g, D ^{\al, \b} g \ra _{\cYe} \le - \cgam \| \cP _m D ^{\al, \b} g \| _{\cYE} ^2 &+ \f\cgam6 \| D ^{\al, \b} g \| _{\cYE} ^2 
		+ C _k \| D ^{<|\al| + |\b|} g \| _{\cYE} ^2.
	\end{align}
    Provided this is true, we are left with removing the projection $\cP _m$ from the $\cgam$-term; for that we use the bound
    $$
        - \| \cP _m D ^{\al, \b} g \| _{\cYE} ^2 \le  -\frac{1}{2} \|  D ^{\al, \b} g \| _{\cYE} ^2 +  \| \cP _M D ^{\al, \b} g \| _{\cYE} ^2.
    $$
    In the second term we will commute the operators $\cP _M$ with $D ^{\al, \b}$. Using \eqref{eqn:commutator-projection-sigma} we have
    \begin{align*}        
        \| \cP _M D ^{\al, \b} g - D ^{\al, \b} \cP _M g \| _{\cYE} ^2 &\le C _k \| D ^{< |\al| + |\b|} g \| _{\cYE} ^2.
    \end{align*}
    Combined with \eqref{eqn: Dab-Lm-g}, we deduce 
    \begin{align*}        
		\la D ^{\al, \b} \cL _\cM g, D ^{\al, \b} g \ra _{\cYe} \le - \frac{ \cgam }2 \| D ^{\al, \b} g \| _{\cYE} ^2 + \cgam \| D ^{\al, \b} \cP _M g \| _{\cYE} ^2 
        & + \f\cgam6 \| D ^{\al, \b} g \| _{\cYE} ^2 \\
		& + C _k \| D ^{< |\al| + |\b|} g \| _{\cYE} ^2.
    \end{align*}
    We conclude the proof of the first claim by the definition of $\cYe ^k$. The second claim follows from absorbing the last lower order term using $\| g \| _{\cYE ^{k - 1}} ^2 \le \nu \| g \| _{\cYE ^k} ^2$ and setting $\nu _k \le \cgam / (12 C_k) $.
    
    The rest of the proof is devoted to proving \eqref{eqn: Dab-Lm-g}. Recall that
    \begin{align*}
        h _0 = D ^{\al, \b} \cL _\cM (\tFm) - \cL _\cM (D ^{\al, \b} \tFm)
        &= D ^{\al, \b} \cN (\rhos \cMM ^{1/2}, g) - \cN (\rhos \cMM ^{1/2}, D ^{\al, \b} g) \\ 
        & \qquad + D ^{\al, \b} \cN (g, \rhos \cMM ^{1/2}) - \cN (D ^{\al, \b} g, \rhos \cMM ^{1/2}).
    \end{align*}
    We apply Lemma \ref{lem: N-derivative-bound} to the two commutator terms:
    \begin{align*}
        &\B| \left\la D ^{\al, \b} \cN (\rhos \cMM ^{1/2}, g) - \cN (\rhos \cMM ^{1/2}, D ^{\al, \b} g), D ^{\al, \b} g \right\ra _V \B|\\
        & \qquad \les _{\al, \b} \sum _{\substack{\al _1 + \al _2 \preceq \al , \, \b _1 + \b _2 \preceq \b \\ (\al _2, \b _2) \prec (\al, \b)}} \cs ^{-3} \|  D^{\alpha_1, \beta_1} (\rhos \cMM ^{1/2}) \| _{L ^2 (V)} \| D^{\alpha_2, \beta_2} g \| _{\sigma} \| D ^{\al, \b} g \| _{\sigma} \\
        & \qquad \les _{\al, \b} \sum _{\substack{\al _1 + \al _2 \preceq \al , \,  \b _1 + \b _2 \preceq \b \\ (\al _2, \b _2) \prec (\al, \b)}} \| D^{\alpha_2, \beta_2} g \| _{\sigma} \| D ^{\al, \b} g \| _{\sigma} \\ 
        & \qquad \le \f\cgam{12} \| D ^{\al, \b} g \| _\sigma ^2 + C _k \| D ^{< |\al| + |\b|} g \| _\sigma ^2 , 
    \end{align*}
    \begin{align*}
        & \B| \left\la D ^{\al, \b} \cN (g, \rhos \cMM ^{1/2}) - \cN (D ^{\al, \b} g, \rhos \cMM ^{1/2}), D ^{\al, \b} g \right\ra _V \B| \\
        & \qquad \les _{\al, \b} \sum _{\substack{\al _1 + \al _2 \preceq \al , \, \b _1 + \b _2 \preceq \b \\ (\al _1, \b _1) \prec (\al, \b)}} \cs ^{-3} \| D^{\alpha_1, \beta_1} g \| _{L ^2 (V)} \| D^{\alpha_2, \beta_2} (\rhos \cMM ^{1/2}) \| _{\sigma} \| D ^{\al, \b} g \| _{\sigma} \\
        & \qquad \les _{\al, \b} \sum _{\substack{\al _1 + \al _2 \preceq \al , \, \b _1 + \b _2 \preceq \b \\ (\al _1, \b _1) \prec (\al, \b)}} \cs ^{\f{\g + 3}2} \| D^{\alpha_1, \beta_1} g \| _{L ^2 (V)} \| D ^{\al, \b} g \| _{\sigma} \\ 
        & \qquad \les _{\al, \b} \sum _{\substack{\al _1 + \al _2 \preceq \al , \, \b _1 + \b _2 \preceq \b \\ (\al _1, \b _1) \prec (\al, \b)}} \| D^{\alpha_1, \beta_1} g \| _{\sigma} \| D ^{\al, \b} g \| _{\sigma} \\ 
        & \qquad \le \f\cgam{12} \| D ^{\al, \b} g \| _\sigma ^2 + C _k \| D ^{< |\al| + |\b|} g \| _\sigma ^2.
    \end{align*}
    Here we used $|D ^{\al, \b} \rhos| \les _{\al, \b} \rhos$, $\| D ^{\al, \b} \cMM ^{1/2} \| _{L ^2 (V)} \les _{\al, \b} 1$, and $\| D ^{\al, \b} \cMM ^{1/2} \| _\s \les _{\al, \b} \cs ^{\f{\g + 3}2}$, which follow from \eqref{eq:dec_S} and Lemma \ref{lem:basis}.  Note that the above estimates \textit{do not} depend on $\nu$. Thus 
    \begin{align*}
        \la h _0, D ^{\al, \b} g \ra _V \le \f\cgam6 \| D ^{\al, \b} g \| _\s ^2 + C _k \|D ^{<|\al| + |\b|} g \| _\s ^2.
    \end{align*}
    For the main term, we apply Lemma \ref{lem: spectral-gap}:
    \begin{align*}
        \la \cL _\cM D ^{\al, \b} g, D ^{\al, \b} g \ra _V \le - \cgam \| \cP _m D ^{\al, \b} g \| _\s ^2.
    \end{align*}
    In summary, we have 
    \begin{align*}
        \la D ^{\al, \b} \cL _\cM g, D ^{\al, \b} g \ra _V \le -\cgam \| \cP _m D ^{\al, \b} g \| _\s ^2 + \f\cgam6 \| D ^{\al, \b} g \| _\s ^2 + C _k \| D ^{< |\al| + |\b|} g \| _\s ^2.
    \end{align*}
    Integrating in $X$ with $\ang X ^\eta$ weight, we conclude \eqref{eqn: Dab-Lm-g}. The proof of the lemma is completed.
\end{proof}

\subsection{Transport operator \texorpdfstring{$\cT$}{T} estimate}
\label{sec: scaling-field}

In this subsection, we estimate the second part of $\Lmic$, which includes the transport operator $\cT$ and also the reaction terms $\f32\bcv$, $\dcm$.
Recall 
\begin{align*}
    \cT g &= (\bcx X \cdot \na _X + \bcv V \cdot \na _V + V \cdot \na _X) g, &
    \dcm &= \frac12 (\pa _s + \cT) \log \cMM + \f32 \bcv.
\end{align*}

\begin{lemma}
    \label{lem: T-dM-Yl}
    Let $\eta \in \R$. 
    There exists $C > 0$ such that
    \begin{align}
        \label{eqn: T-Yl-estimate}
        -\left\langle \left(\cT - \f32 \bcv \right) g, g \right\rangle _{\cYe} &\le \frac{\bcx}2 (\eta - \etab) \| g \| _{\cYe} ^2 + C |\eta| \cdot 
       \| g \| _{\cYe} \| g \| _{\cYE}, %
        \\
        \label{eqn: dM-Yl-estimate}
        |\langle \dcm g, g \rangle _{\cYe}| &\le C %
        \| g \| _{\cYe} ^\f{2 (\g - 1)}{\g + 2} \| g \| _{\cYE}^\f6{\g + 2}  . 
    \end{align}
\end{lemma}

\begin{proof}
    We use integration by parts:
    \begin{align*}
		-2 \ang{\cT g, g} _{\cYe} 
        & = -\iint \ang X ^\eta \cT |g| ^2 d V d X \\
		& = \iint [\div _X (\bcx X \ang X ^\eta + V \ang X ^\eta) + \div _V (\bcv V \ang X ^\eta)] |g| ^2 d V d X \\
		& = \iint \left(
		  3 \bcv + 3 \bcx + \eta \bcx \frac{|X| ^2}{\ang X ^2} + \eta \frac{V \cdot X}{\ang X ^2}
		\right) \ang X ^\eta |g| ^2 d V d X \\
        & \leq \iint \left(
		  3 \bcv + 3 \bcx + \eta \bcx 
		\right) \ang X ^\eta |g| ^2 d V d X
        \\
        & \qquad - \eta \bcx \iint \ang X ^{-2} \ang X ^\eta |g| ^2 d V d X + |\eta| \iint |V| \ang X ^{-1} \ang X ^\eta |g| ^2 d V d X. 
	\end{align*}
    The first integral is exactly $(3 \bcv + 3 \bcx + \eta \bcx) \| g \| _{\cYe} ^2 = (\bcx (\eta - \etab) - 3 \bcv) \| g \| _{\cYe} ^2$, because $\bcx \etab = -3 \bcx - 6 \bcv$ from the definition \eqref{wg:X_power}. 

    For the second and third integral, we use \eqref{eqn: interpolation-v-holder}: $\ang X ^{-2} \le \ang X ^{-1} \les \Lams ^\f12$, $|V| \ang X ^{-1} \les \Lams ^\f12$, so
    \begin{align*}
        \iint ( \bcx \ang X^{-2} + |V| \ang X ^{-1} ) \ang X ^\eta |g| ^2 d V d X & \les \iint \Lams ^\f12 \ang X ^\eta |g| ^2 d V d X \\
        & \le \left(\iint |g| ^2 d V d X\right) ^\f12 \left(\iint \Lams \ang X ^\eta |g| ^2 d V d X\right) ^\f12 \\
        & \les \| g \| _{\cYe} \| g \| _{\cYE}.
    \end{align*}
    This proves \eqref{eqn: T-Yl-estimate}. 
    
    As for \eqref{eqn: dM-Yl-estimate}, recall that $\dcm = \cO (\ang X ^{-r} \ang \vc ^3)$ in \eqref{eqn: Dab-dcm-bound}. Using \eqref{eqn: interpolation-v3-holder}, we can control
    \begin{align*}
        \left| \iint \dcm |g| ^2  \ang X ^\eta d V d X \right| & \les \iint |g| ^2 \Lams ^\f3{\g + 2} \ang X ^\eta d V d X \les \| g \| _{\cYe} ^\f{2(\g - 1)}{\g + 2} \| g \| _{\cYE} ^\f6{\g + 2}.
    \end{align*}
    The lemma is thus proven.
\end{proof}

\begin{cor}
    \label{cor: scaling-field-coersive}
    For every $k \ge 0$, $\eta \le 0$, 
    there exists $C _{k, \eta}$ such that 
    \begin{align*}
        \left\langle \left(\pa _s + \cT + \dcm - \f32 \bcv\right) g, g \right\rangle _{\cYe ^k} & \ge \f12 \f{d}{ds} \| g \| _{\cYe ^k} ^2 + \f{\bcx}2 (\etab - \eta) \| g \| _{\cYe ^k} ^2 \\  
        & \qquad - C_{k, \eta} \| g \| _{\cYe ^k} ^\f{2(\g - 1)}{\g + 2} \| g \| _{\cYE ^k} ^\f6{\g + 2} 
        - C _{k, \eta} \| g \| _{\cYe ^k} \| g \| _{\cYE ^k} .
    \end{align*}
\end{cor}

\begin{proof}
    Recall the commutator $h _1$ is defined by
    \begin{align*}
        D ^{\al, \b} \left(\pa _s + \cT + \dcm - \f32 \bcv\right) g = \left(\pa _s + \cT + \dcm - \f32 \bcv\right) D ^{\al, \b} g + h _1. 
    \end{align*}
    Therefore, 
    \begin{align*}
        \left\la D ^{\al, \b} \left(\pa _s + \cT + \dcm - \f32 \bcv\right) g, D^{\al, \b} g \right\ra _{\cYe} & = \left\la \pa _s D ^{\al, \b} g, D ^{\al, \b} g \right\ra _{\cYe} \\ 
        & \qquad + \left\la \left(\cT + \dcm - \f32 \bcv\right)  D ^{\al, \b} g, D ^{\al, \b} g \right\ra _{\cYe} \\
        & \qquad + \la h _1, D ^{\al, \b} g \ra _{\cYe}. 
    \end{align*}
    For the first inner product, it equals 
    \begin{align*}
        \left\la \pa _s D ^{\al, \b} g, D ^{\al, \b} g \right\ra _{\cYe} = \f12 \f{d}{ds} \left\la D ^{\al, \b} g, D ^{\al, \b} g \right\ra _{\cYe} = \f12 \f{d}{ds} \| D ^{\al, \b} g \|^2 _{\cYe}.
    \end{align*}
    For the second inner product, 
    applying Lemma \ref{lem: T-dM-Yl} to $D ^{\al, \b} g$ yields
    \begin{align*}    
        & -\left\langle \left(\cT + \dcm - \f32 \bcv \right) D ^{\al, \b} g, D ^{\al, \b} g \right\rangle _{\cYe} \\
        & \qquad \le \frac{\bcx}2 (\eta - \etab) \| D ^{\al, \b} g \| _{\cYe} ^2 
        + C _{\eta} \| D ^{\al, \b} g \| _{\cYE} ^\f{2 (\g - 1)}{\g + 2} \| D ^{\al, \b} g \| _{\cYe} ^\f6{\g + 2} 
        + C \| D ^{\al, \b} g \| _{\cYE} \| D ^{\al, \b} g \| _{\cYe}.
    \end{align*}
    For the third inner product, we use Lemma \ref{lem: commu-derivative-projection} (2) and (4) :
    \renewcommand{\tFm}{g}
    \begin{align*}
        h_1 &= D ^{\al, \b} (\pa _s + \cT) \tFm - (\pa _s + \cT) D ^{\al, \b} \tFm + D ^{\al, \b} (\dcm \tFm) - \dcm D ^{\al, \b} \tFm \\
        &= \cO (\cs \ang X ^{-1} \ang \vc + \ang X^{-1}) |D ^{\le |\al| + |\b|} \tFm| + \cO ( \ang X ^{-r} \ang \vc ^3) |D^{<|\al| + |\b|}\tFm| 
        \\
        & \les (\Lams ^\f12 + \Lams ^\f{3}{\g + 2})  |D^{\le |\al| + |\b|}\tFm|.
    \end{align*}
    We used the Leibniz rule for the term $\dcm$, and applied \eqref{eqn: interpolation-v3-holder} and \eqref{eqn: interpolation-v-holder}. We conclude 
    \begin{align*}
        \la h _1, D ^{\al, \b} g \ra _{\cYe} \les C _\eta \| D ^{\le |\al| + |\b|} g \| _{\cYE} ^\f{2 (\g - 1)}{\g + 2} \| D ^{\le |\al| + |\b|} g \| _{\cYe} ^\f6{\g + 2} 
        + C \| D ^{\le |\al| + |\b|} g \| _{\cYE} \| D ^{\le |\al| + |\b|} g \| _{\cYe}.
    \end{align*}
    
    Combine the three inner products and summing $\al, \b$ by the definition of $\cYe ^k$ norm \eqref{norm:Y}, we conclude 
    \begin{align*}
        \f12 \f{d}{ds} \| g \| _{\cYe ^k} ^2 + \f{\bcx}2 (\etab - \eta) \| g \| _{\cYe ^k} ^2 &\le \left\langle \left(\pa _s + \cT + \dcm - \f32 \bcv\right) g, g \right\rangle _{\cYe ^k} \\
        & \qquad + C _{k, \eta} \sum _{|\al| + |\b| \le k} \nu ^{|\al| + |\b| - k} \| D ^{\le |\al| + |\b|} g \| _{\cYE} ^\f{2 (\g - 1)}{\g + 2} \| D ^{\le |\al| + |\b|} g \| _{\cYe} ^\f6{\g + 2} \\
        & \qquad + C _{k, \eta} \sum _{|\al| + |\b| \le k} \nu ^{|\al| + |\b| - k} \| D ^{\le |\al| + |\b|} g \| _{\cYE} \| D ^{\le |\al| + |\b|} g \| _{\cYe}.
    \end{align*}
    The proof is completed by H\"older inequality.
\end{proof}

\subsection{Projection \texorpdfstring{$\cP _M$}{P M} estimates}
\label{sec: projection-m-free-transport}

Now we estimate the last part of $\Lmic$, which is the term involving macro projection $\cP _M$. First, we estimate the main terms. 

\begin{lemma}
    \label{lem: PMVnaX}
Suppose $\g \in (1, 2]$. Let $g = \cP _m g$. For any multi-indices $\al, \b$ and $\eta \in \R$ it holds that 
    \begin{align}
        \label{eqn: PMVnaX}
        &\left| \la \cP _M [V \cdot \na _X D^{\al,\b} g], D^{\al,\b} g \ra_{\cYe} \right| 
        \les _{\eta, \al, \b} \| D ^{\le |\al| + |\b|} g \| _{\cYe} \|D ^{\le |\al| + |\b|} g\| _{\cYE} .
    \end{align}
    Moreover, 
    \begin{align}
        \label{eqn: PMdcm-1}
        &\left| \la \cP _M [\dcm D^{\al,\b} g], D^{\al,\b} g \ra_{\cYe} \right|
        \les _{\eta, \al, \b} \| D ^{\le |\al| + |\b|} g \| _{\cYe} ^\f{2 (\g - 1)}{\g + 2} \|D ^{\le |\al| + |\b|} g\| _{\cYE} ^\f6{\g + 2}, \\
        \label{eqn: PMdcm-2}
        &\left| \la \cP _M [\td d _\cM D^{\al,\b} g], D^{\al,\b} g \ra_{\cYe} \right| 
        \les _{\eta, \al, \b} \| D ^{\le |\al| + |\b|} g \| _{\cYe} ^\f{2 (\g - 1)}{\g + 2} \|D ^{\le |\al| + |\b|} g\| _{\cYE} ^\f6{\g + 2} .
    \end{align}
\end{lemma}

\begin{proof}
We use the fact that $\cP_M$ is a projection and integrate by part to obtain
\begin{align*}
    &\left| \la \cP _M [V \cdot \na _X D^{\al,\b} g], D^{\al,\b} g \ra_{\cYe} \right| \\
    & \qquad = \left| \iint \ang X ^\eta \cP _M [V \cdot \na _X (D ^{\al, \b} g)] \cdot D ^{\al, \b} g d V d X \right| \\
    & \qquad = \left| \iint D ^{\al, \b} g \cdot V \cdot \na _X \left(\ang X ^\eta \cP _M [D ^{\al, \b} g] \right) d V d X \right| \\
    & \qquad \le \iint |D ^{\al, \b} g| \cdot |V| \left(|\eta| \ang X ^{\eta - 1} \left| \cP _M [D ^{\al, \b} g] \right| + \ang X ^\eta \left|\na _X \cP _M [D ^{\al, \b} g]\right| \right) d V d X \\
    & \qquad \les _\eta \iint |D ^{\al, \b} g| \cdot |V| \ang X ^{-1} \left(\left| \cP _M [D ^{\al, \b} g] \right| + \left|\ang X \na _X \cP _M [D ^{\al, \b} g]\right|\right) \ang X ^\eta d V d X \\
    & \qquad \les \iint |D ^{\al, \b} g| \cdot |V| \ang X ^{-1} \left| D ^{\le 1} \cP _M [D ^{\al, \b} g]\right| \ang X ^\eta d V d X.
\end{align*}
Recall that $|V| \ang X ^{-1} \les \Lams ^\f12$ from \eqref{eqn: interpolation-v-holder}. So 
\beq\begin{aligned}
    \label{eqn: splitting-in-PM-VDX}
    \left| \la \cP _M [V \cdot \na _X D^{\al,\b} g], D^{\al,\b} g \ra_{\cYe} \right|
    & \les \| D ^{\al, \b} g \| _{\cYe} ^\f12 \| D ^{\le 1} \cP _M [D ^{\al, \b} g] \| _{\cYe} ^\f12 \\
    & \qquad \times \| D ^{\al, \b} g \| _{\cYE} ^\f12 \| D ^{\le 1} \cP _M [D ^{\al, \b} g] \| _{\cYE} ^\f12.
\end{aligned}\eeq
By the commutator estimate \eqref{eqn: comm-P-Dab-L2} and Corollary \ref{cor: concatenate-derivative}, we can commute $\cP _M$ and $D ^{\le 1}$, $D ^{\al, \b}$ up to lower order commutator: 
\begin{align*}
    \| D ^{\le 1} \cP _M [D ^{\al, \b} g] \| _{\cYe} &\les \| \cP _M [D ^{\le 1} D ^{\al, \b} g] \| _{\cYe} + \| D ^{\al, \b} g \| _{\cYe} \\
    &\les _{\al, \b} \| \cP _M D ^{\le |\al| + |\b| + 1} g] \| _{\cYe} + \| D ^{\al, \b} g \| _{\cYe} \\
    & \les _{\al, \b} \| D ^{\le |\al| + |\b|} g \| _{\cYe} +  \| D ^{\le |\al| + |\b| + 1} \cP _M  g \| _{\cYe} \\
    & = \| D ^{\le |\al| + |\b|} g \| _{\cYe}.
\end{align*}
In the last step, we used $g = \cP _m g, \cP_M g = 0$. Similarly using \eqref{eqn: comm-P-Dab-sigma} we have 
\begin{align*}
    \| D ^{\le 1} \cP _M [D ^{\al, \b} g] \| _{\cYE} \les _{\al, \b} \| D ^{\le |\al| + |\b|} g \| _{\cYE}.
\end{align*}
Plugging into \eqref{eqn: splitting-in-PM-VDX} we obtain \eqref{eqn: PMVnaX}. As for \eqref{eqn: PMdcm-1}, note that 
\begin{align*}
    & \la \cP _M [\dcm D^{\al,\b} g], D^{\al,\b} g \ra_{\cYe} = \iint \ang X ^\eta \cP _M [D ^{\al, \b} g] \cdot \dcm D ^{\al, \b} g d V d X.
\end{align*}
Using $\dcm \les \ang X ^{-r} \ang \vc ^3 \les \Lams ^\f3{\g + 2}$ from Lemma \ref{lem: commu-derivative-projection} (4) and \eqref{eqn: interpolation-v3-holder}, we have the following bound similar to \eqref{eqn: splitting-in-PM-VDX}:
\begin{align*}
    \left| \la \cP _M [\dcm D^{\al,\b} g], D^{\al,\b} g \ra_{\cYe} \right|
    & \les \| D ^{\al, \b} g \| _{\cYe} ^\f{\g - 1}{\g + 2} \| \cP _M [D ^{\al, \b} g] \| _{\cYe} ^\f{\g - 1}{\g + 2} \| D ^{\al, \b} g \| _{\cYE} ^\f3{\g + 2} \| \cP _M [D ^{\al, \b} g] \| _{\cYE} ^\f3{\g + 2}.
\end{align*}
Thus, the conclusion follows the same proof. The case of $\td d _\cM$ is identical: thanks to Lemma \ref{lem: commu-derivative-projection} (4) and \eqref{eqn: interpolation-v3-holder} again we have $\td d _\cM \les \ang X ^{-1} \cs \ang \vc ^3 \les \Lams ^\f3{\g + 2}$.
\end{proof}

\vspace{0.5cm}

We are ready to prove the $\cYe ^k$ estimate for $\cP _M$ terms. 

\begin{cor}
    \label{cor: Yk-Pm-estimate}
    Suppose $\g \in (1, 2]$. If $g = \cP _m g$ then
    \begin{align}
        \label{eqn: Yk-free-transport}
        \la \cP _M [V \cdot \na _X g], g \ra _{\cYe ^k} &\les _k C _k \| g \| _{\cYe ^k} \| g \| _{\cYE ^k}, \\
        \label{eqn: Yk-dcm}
        \la \cP _M [(2 \dcm + \td d _\cM) g], g \ra _{\cYe ^k} &\les _k C _k \| g \| _{\cYe ^k} ^\f{2 (\g - 1)}{\g + 2} \| g \| _{\cYE ^k} ^\f6{\g + 2}.
    \end{align}
\end{cor}

\begin{proof}

    Recall the commutator $h _4$ is defined as 
    \begin{align*}
        h _4 & := 
        D ^{\al, \b} \cP _M [V \cdot \na _X g] - \cP _M [(V \cdot \na _X) D ^{\al, \b} g] \\
        & = D ^{\al, \b} \cP _M [V \cdot \na _X g] - \cP _M [D ^{\al, \b} (V \cdot \na _X ) g] \\
        & \qquad + \cP _M [D ^{\al, \b} (V \cdot \na _X ) g] - \cP _M [(V \cdot \na _X) D ^{\al, \b} g ] \\
        &=: h_{4,1} + \cP _M [h_{4,2}]. 
    \end{align*}
    We first handle $h _{4, 1}$ term. By interpolation \eqref{eqn: interpolation-v-holder} we conclude 
    \begin{align*}
        | \ang{h _{4, 1}, D^{\al, \b} g} _{\cYe}| \le & \| \Lams ^\f14 D ^{\al, \b} g \| _{\cYe} \| \Lams ^{-\f14} h _{4, 1} \| _{\cYe} \\
        \les & \| D ^{\al, \b} g \| _{\cYe} ^\f12 \|D^{\al, \b} g \|_{\cYE} ^\f12 \| \ang X ^{\f12} h _{4, 1} \| _{\cYe}.
    \end{align*}
    Using \eqref{eqn:commutator-projection} and $|V| \les \cs \ang \vc$, we know
    \begin{align}
        \notag
        \| h _{4, 1} \| _{L ^2 (V)} &\les \| \ang \vc ^{-1} D ^{\le |\al| + |\b| - 1} (V \cdot \na _X g) \| _{L ^2 (V)} \\
        \notag
        & \les \| \ang \vc ^{-1} (|V| + \cs) D ^{\le |\al| + |\b| - 1} (\na _X g) \| _{L ^2 (V)} \\
        & \les \ang X ^{-1} \| D ^{\le |\al| + |\b|} g \| _{L ^2 (V)}.
        \label{eqn: h41-L2V-est}
    \end{align}
    Thus 
    \begin{align*}
        \| \ang X ^\f12 h _{4, 1} \| _{\cYe} \les \| \ang X ^{-\f12} D ^{\le |\al| + |\b|} g \| _{\cYe} \les \| \Lams ^\f14 D ^{\le |\al| + |\b|} g \| _{\cYe} \les \| D ^{\le |\al| + |\b|} g \| _{\cYe} ^\f12 \|D ^{\le |\al| + |\b|} g \|_{\cYE} ^\f12
    \end{align*}    
    and we conclude 
    \begin{align}
        \label{eqn: main-h41}
        | \ang{h _{4, 1}, D^{\al, \b} g} _{\cYe} | \le \| D ^{\le |\al| + |\b|} g \| _{\cYe} \|D ^{\le |\al| + |\b|} g \|_{\cYE}.
    \end{align}

    Next we handle $h _{4, 2}$. By \eqref{eqn: comm-free-transp-Dab} we have 
    \begin{align*}
        |h _{4, 2}| &\les \cs \ang X ^{-1} \ang \vc |D ^{\le |\al| + |\b|} g| \les \Lams ^\f12 |D ^{\le |\al| + |\b|} g|,
    \end{align*}
    thus  
    \begin{align*}
        | \la \cP _M [h _{4, 2}], D ^{\al, \b} g \ra _{\cYe} | & = | \la h _{4, 2}, \cP _M [ D ^{\al, \b} g] \ra _{\cYe} | \\
        & \les \| \Lams ^\f14 D ^{\le |\al| + |\b|} g \| _{\cYe}  \| \Lams ^\f14 \cP _M [D ^{\al, \b} g] \| _{\cYe} \\
        & \les \| D ^{\le |\al| + |\b|} g \| _{\cYe} ^\f12 \|D ^{\le |\al| + |\b|} g \|_{\cYE} ^\f12 \| \cP _M [D ^{\al, \b} g] \| _{\cYe} ^\f12 \| \cP _M [D ^{\al, \b} g] \|_{\cYE} ^\f12.
    \end{align*}
    Use the projection bound \eqref{eqn: comm-P-Dab-L2} and \eqref{eqn: comm-P-Dab-sigma} we conclude 
    \begin{align}
        \label{eqn: main-h42}
        | \ang{\cP _M [h _{4, 2}], D^{\al, \b} g} _{\cYe} |  \les & \| D ^{\le |\al| + |\b|} g \| _{\cYe} \|D^{\le |\al| + |\b|} g \|_{\cYE}.
    \end{align}
    Summarizing, we have 
    \begin{align*}
        | \la D ^{\al, \b} \cP _M [V \cdot \na _X g], D ^{\al, \b} g \ra _{\cYe} |  &= \B|  \underbrace{\la \cP _M [(V \cdot \na _X) D ^{\al, \b} g], D ^{\al, \b} g \ra} _{\eqref{eqn: PMVnaX}} \\
        & \qquad + \underbrace{\ang{h _{4,1}, D^{\al, \b} g} _{\cYe}} _{\eqref{eqn: main-h41}} + \underbrace{\ang{\cP _M [h _{4,2}], D^{\al, \b} g} _{\cYe}} _{\eqref{eqn: main-h42}} \B| \\
        & \les \| D ^{\le |\al| + |\b|} g \| _{\cYe} \|D^{\le |\al| + |\b|} g \|_{\cYE}.
    \end{align*}
    This proves \eqref{eqn: Yk-free-transport} by H\"older inequality.

    We apply a similar splitting to $h _5$:
    \begin{align*}
        h _5 & := 
        -D ^{\al, \b} \cP _M [(2 \dcm + \td d _\cM) g] + \cP _M[ (2 \dcm + \td d _\cM) D ^{\al, \b} g] \\
        & = -D ^{\al, \b} \cP _M [(2 \dcm + \td d _\cM) g] + \cP _M [D ^{\al, \b} (2 \dcm + \td d _\cM) g] \\
        & \qquad - \cP _M [D ^{\al, \b} (2 \dcm + \td d _\cM) g] + \cP _M [(2 \dcm + \td d _\cM) D ^{\al, \b} g] \\
        &=: h_{5,1} + \cP _M [h_{5,2}].
    \end{align*}
    By interpolation \eqref{eqn: interpolation-v3-holder} we conclude 
    \begin{align*}
        | \ang{h _{5, 1}, D^{\al, \b} g} _{\cYe} | \le & \| \Lams ^\f3{2 (\g + 2)} D ^{\al, \b} g \| _{\cYe} \| \Lams ^{-\f3{2 (\g + 2)}} h _{5, 1} \| _{\cYe} \\
        \les & \| D ^{\al, \b} g \| _{\cYe} ^\f{\g - 1}{\g + 2} \|D^{\al, \b} g \|_{\cYE} ^\f3{\g + 2} \| \ang X ^{\f12} \cs ^{-\f12} h _{5, 1} \| _{\cYe}.
    \end{align*}
    By \eqref{eqn:commutator-projection} with $N = 3$, using derivative bound Lemma \ref{lem: commu-derivative-projection} (4) we have
    \begin{align*}
        \| h_{5, 1} \| _{L ^2 (V)} 
        & \les \left\| \ang \vc ^{-3} \left|D ^{< |\al| + |\b|} (2 \dcm + \td d _\cM)\right| \cdot \left|D ^{< |\al| + |\b|} g \right|  \right\| _{L ^2 (V)}
        \\
        & \les \| \ang X ^{-1} \cs D ^{< |\al| + |\b|} g \| _{L ^2 (V)}.
    \end{align*}
    Combining them, by \eqref{eqn: interpolation-v3-holder} we have
    \begin{align*}
         \| \ang X ^\f12 \cs ^{-\f12} h _{5, 1} \| _{\cYe} 
         \les \| \ang X ^{-\f12} \cs ^\f12 D ^{\le |\al| + |\b|} g \| _{L ^2 (V)} 
         \les \| D ^{\le |\al| + |\b|} g \| _{\cYe} ^\f{\g - 1}{\g + 2} \|D ^{\le |\al| + |\b|} g \|_{\cYE}^\f3{\g + 2}.
    \end{align*}
    Thus 
    \begin{align}
        \label{eqn: main-h51}
        | \ang{h _{5, 1}, D^{\al, \b} g} _{\cYe} | \le & \| D ^{\le |\al| + |\b|} g \| _{\cYe} ^\f{2 (\g - 1)}{\g + 2} \|D^{\le |\al| + |\b|} g \|_{\cYE} ^\f6{\g + 2} .
    \end{align}
    Finally, since 
    \begin{align*}
        |h _{5, 2}| \les \left|D ^{< |\al| + |\b|} (2 \dcm + \td d _\cM)\right| \cdot \left|D ^{< |\al| + |\b|} g \right| \les \Lams ^{\f3{\g + 2}} |D ^{< |\al| + |\b|} g |,
    \end{align*}
    we have
    \begin{align*}
        | \ang{\cP _M [h _{5, 2}], D^{\al, \b} g} _{\cYe}|  &= | \ang{h _{5, 2}, \cP _M D^{\al, \b} g} _{\cYe} |\\
        & \les \| D ^{< |\al| + |\b|} g \| _{\cYe} ^\f{\g - 1}{\g + 2} \|D ^{< |\al| + |\b|} g \|_{\cYE}^\f3{\g + 2} \| \cP _M D ^{\al, \b} g \| _{\cYe} ^\f{\g - 1}{\g + 2} \| \cP _M D ^{\al, \b} g \|_{\cYE}^\f3{\g + 2}.
    \end{align*}
    Use the projection bound \eqref{eqn: comm-P-Dab-L2} and \eqref{eqn: comm-P-Dab-sigma} again we conclude 
    \begin{align}
        \label{eqn: main-h52}
        |\ang{\cP _M [h _{5, 2}], D^{\al, \b} g} _{\cYe} | \les & \| D ^{\le |\al| + |\b|} g \| _{\cYe} ^\f{2 (\g - 1)}{\g + 2} \|D^{\le |\al| + |\b|} g \|_{\cYE}^\f6{\g + 2}.
    \end{align}
    Summarizing, we have 
    \begin{align*}
        \la D ^{\al, \b} \cP _M [(2 \dcm + \td d _\cM) g], D ^{\al, \b} g \ra _{\cYe} &= \underbrace{\la \cP _M [(2 \dcm + \td d _\cM) D ^{\al, \b} g], D ^{\al, \b} g \ra} _{\eqref{eqn: PMdcm-1}, \eqref{eqn: PMdcm-2}} \\
        & \qquad + \underbrace{\ang{h _{5, 1}, D^{\al, \b} g} _{\cYe}} _{\eqref{eqn: main-h51}} + \underbrace{\ang{\cP _M [h _{5, 2}], D^{\al, \b} g} _{\cYe}} _{\eqref{eqn: main-h52}} \\
        & \les \| D ^{\le |\al| + |\b|} g \| _{\cYe} ^\f{2 (\g - 1)}{\g + 2} \|D^{\le |\al| + |\b|} g \|_{\cYE}^\f6{\g + 2}.
    \end{align*}
    The proof of \eqref{eqn: Yk-dcm} is complete after applying  H\"older inequality.
\end{proof}

Using the above estimates, below, we prove \eqref{eq:micro_lin_est:lin} and \eqref{eq:micro_lin_est:lin2}. 

\begin{proof}[Proof of \eqref{eq:micro_lin_est:lin} and \eqref{eq:micro_lin_est:lin2}]
    The estimates for $\partial_s - \Lmic$ in \eqref{eq:micro_lin_est:lin} follows from Lemma \ref{lem:L_m_Dalphabeta} combined with Corollary \ref{cor: scaling-field-coersive} and Corollary \ref{cor: Yk-Pm-estimate}. Using Young's inequality and $0< \g-1 < 3$, 
    we obtain
    \[
    \bal 
        C _{k, \eta} \es \| g \|_{\cYe ^k}^2 + \f{1}{ 50 \es } \| g \|_{\cYE ^k }^2 
        & \geq  C _{k,\eta} \| g \|_{\cYe ^k}\| g \|_{\cYE ^k }, \\ 
        C _{k, \eta} \es \| g \|_{\cYe ^k}^2 + \f{1}{ 50 \es } \| g \|_{\cYE^k }^2  
        & \geq C _{k, \eta}  \left( \es  \| g \|_{\cYe ^k}^2 \right)^{ \f{\g-1}{\g+2} }   \left( \es^{-1}  \| g \|_{\cYE ^k}^2 \right)^{ \f{3}{\g+2} }
        \geq C_{k,\eta} \es^{ \f{\g-4}{\g+1}} \| g \|_{\cYe ^k}^{ \f{2(\g-1)}{\g+2}}  \| g \|_{\cYE ^k}^{ \f{6}{\g+2} }.
        \eal 
    \]
    Since $\es \leq 1$ and $2 \lam_{\eta} = \f{\bcx}{2} (\etab - \eta) $ by \eqref{eq:lam_eta},  we obtain  $ \es^{ \f{\g-4}{\g+1}} \geq 1$ and prove \eqref{eq:micro_lin_est:lin2}.
\end{proof}

\subsection{Estimates of the macro terms}
\label{sec: projection-free-transport}

Recall the decomposition \ref{eq:Pm_tran_FM}. Now we handle $\cP _m$ terms in \eqref{eq:micro_lin_est:cross_d} and \eqref{eq:micro_lin_est:cross}, which involve interaction with $\tFM$. 

\begin{proof}[Proof of \eqref{eq:micro_lin_est:cross_d}]

Thanks to commutator estimate \eqref{eqn: comm-P-Dab-L2} and derivative bound \eqref{eqn: Dab-dcm-bound} we have
\begin{align*}
    \| D ^{\al, \b} \cP _m [(2 \dcm + \td d _\cM) \tFM] \| _{L ^2 (V)} & \les _{\al, \b} \| D ^{\preceq (\al, \b)} ((2 \dcm + \td d _\cM) \tFM) \| _{L ^2 (V)} \\
    & \les _{\al, \b} \| |D ^{\preceq (\al, \b)} (2 \dcm + \td d _\cM)| \cdot |D ^{\preceq (\al, \b)} \tFM| \| _{L ^2 (V)} \\
    & \les _{\al, \b} \ang X ^{-1} \cs \| \ang \vc ^3 |D ^{\preceq (\al, \b)} \tFM| \| _{L ^2 (V)}.
\end{align*}
The weights and $V$-derivative on the macroscopic quantities are negligible, due to \eqref{eq:macro_UPB_est:b}:
\begin{align*}
    \| \ang \vc ^3 |D ^{\preceq (\al, \b)} \tFM| \| _{L ^2 (V)} \les _{\al, \b} \| D ^{\preceq (\al, 0)} \tFM \| _{L ^2 (V)}.
\end{align*}
Therefore, 
\begin{align*}
    |\la D ^{\al, \b} \cP _m [(2 \dcm + \td d _\cM) \tFM], D ^{\al, \b} \tFm \ra| _{L ^2 (V)} &\les \ang X ^{-1} \cs \| D ^{\preceq (\al, 0)} \tFM \| _{L ^2 (V)} \| D ^{\al, \b} \tFm \| _{L ^2 (V)} \\
    & \les \| D ^{\preceq (\al, 0)} \tFM \| _{L ^2 (V)} \| \Lams ^\f12 D ^{\al, \b} \tFm \| _{L ^2 (V)},
\end{align*}
where we used $\ang X ^{-1} \cs \les \Lams ^\f12$ in the last step due to \eqref{eqn: interpolation-v-holder}. We conclude \eqref{eq:micro_lin_est:cross_d} by integrating in $X$ with weight $\ang X ^\eta$ and H\"older inequality.
\end{proof}

\begin{proof}[Proof of \eqref{eq:micro_lin_est:cross}]

We first decompose the commutator $h _2$ similar to the $h _4$ term:
\begin{align*}
    h _2 &:= D ^{\al, \b} \cP _m [V \cdot \na _X \tFM] - \cP _m [(V \cdot \na _X) D ^{\al, \b} \tFM] \\
    & = D ^{\al, \b} \cP _m [(V \cdot \na _X) \tFM] - \cP _m [D ^{\al, \b} (V \cdot \na _X \tFM)] \\
    & \qquad + \cP _m [D ^{\al, \b} (V \cdot \na _X \tFM)] - \cP _m [(V \cdot \na _X) D ^{\al, \b} \tFM] \\
    & =: h _{2, 1} + \cP _m [h _{2, 2}].
\end{align*}
$h _{2, 1}$ can be handled completely analogously to $h _{4, 1}$ in \eqref{eqn: h41-L2V-est}, yielding
\begin{align*}
    \| h _{2, 1} \| _{L ^2 (V)} \les \ang X ^{-1} \| D ^{\le |\al| + |\b|} \tFM \| _{L ^2 (V)}, 
\end{align*}
For $h _{2, 2}$ we apply the commutator bound \eqref{eqn: comm-free-transp-Dab}, \eqref{eq:macro_UPB_est:b}, and 
$\cs \les 1$:
\begin{align*}
    \| h _{2, 2} \| _{L ^2 (V)} & \les \cs \ang X ^{-1} \| \ang \vc D ^{\le |\al| + |\b|} \tFM \| _{L ^2 (V)} \\
    & \les \ang X ^{-1} \| D ^{\le |\al| + |\b|} \tFM \| _{L ^2 (V)}.
\end{align*}
Therefore 
\begin{align*}
    \la h _2, D ^{\al, \b} \tFm \ra _V &\les \ang X ^{-1} \| D ^{\le |\al| + |\b|} \tFM \| _{L ^2 (V)} \| D ^{\al, \b} \tFm \| _{L ^2 (V)} \\
    & \les \| D ^{\le |\al| + |\b|} \tFM \| _{L ^2 (V)} \| \Lams ^\f12 D ^{\al, \b} \tFm \| _{L ^2 (V)},
\end{align*}
again using $\ang X ^{-1} \les \Lams ^\f12$. Integrating in $X$ with weight $\ang X ^\eta$ yields 
\begin{align*}
    \la h _2, D ^{\al, \b} \tFm \ra _{\cYe} &\les \| D ^{\le |\al| + |\b|} \tFM \| _{\cYe} \| D ^{\al, \b} \tFm \| _{\cYE},
\end{align*}
Therefore, 
\begin{align*}
    - \la \cP _m [(V \cdot \na _X) \tFM], \tFm \ra_{ \cYe ^k}  
    & = \cO _{k, \eta} (\| \tFm \|_{\cYE^{k}}  \| \tFM \|_{\cYe ^k}) \\
    & \qquad - \sum _{|\al| + |\b| \le k} \nu ^{|\al| + |\b| - k} 
        \f{ |\al|! }{ \al !}  \iint (V \cdot \na _X) D ^{\al, \b} \tFM \cdot D ^{\al, \b} \tFm \ang X ^\eta d V d X.
\end{align*}
Whenever $|\al| < k$, we can use \eqref{cor: concatenate-derivative}, \eqref{eq:macro_UPB_est:b} and get 
\begin{align*}
    \| V \cdot \na _X D ^{\al, \b} \tFM \| _{L ^2 (V)} &= \sum _{i = 1} ^3 \| V _i \na _{X _i} D ^{\al, \b} \tFM \| _{L ^2 (V)} \\
    &\les \sum _{i = 1} ^3 \ang X ^{-1} \cs \| \ang \vc D ^{\ee _i, 0} D ^{\al, \b} \tFM \| _{L ^2 (V)} \\
    &\les \sum _{i = 1} ^3 \ang X ^{-1} \cs \| \ang \vc D ^{\preceq (\al + \ee _i, \b)} \tFM \| _{L ^2 (V)} \\
    &\les \sum _{i = 1} ^3 \ang X ^{-1} \cs \| D ^{\preceq (\al + \ee _i, 0)} \tFM \| _{L ^2 (V)}.
\end{align*}
Thus we conclude 
\begin{align*}
    & \sum _{\substack{|\al| + |\b| \le k \\ |\al| < k}} \nu ^{|\al| + |\b| - k}     \f{ |\al|! }{ \al !} \iint (V \cdot \na _X) D ^{\al, \b} \tFM \cdot D ^{\al, \b} \tFm \ang X ^\eta d V d X \\
    & \qquad \les \sum _{\substack{|\al| + |\b| \le k \\ |\al| < k}} \nu ^{|\al| + |\b| - k}     \f{ |\al|! }{ \al !} \int \ang X ^{-1} \cs \|D ^{\le |\al| + 1} \tFM\| _{L ^2 (V)} \| D ^{\al, \b} \tFm \| _{L ^2 (V)} \ang X ^\eta d X \\
    & \qquad \les _k \sum _{ \substack{ |\al| + |\b| \le k  \\ |\al| < k}  } \nu ^{\f{|\b| - 1}2} \int \nu ^{\f{|\al| + 1 - k}2} \|D ^{\le |\al| + 1} \tFM\| _{L ^2 (V)} \cdot \nu ^{\f{|\al| + |\b| - k}2}\| \Lams ^\f12 D ^{\al, \b} \tFm \| _{L ^2 (V)} \ang X ^\eta d X \\
    & \qquad \les _{k} \nu ^{-\f12} \| \tFm \|_{\cYE ^{k}}  \| \tFM \|_{\cYe ^k}.
\end{align*}
Combined, we have completed the proof of \eqref{eq:micro_lin_est:cross}.
\end{proof}

\subsection{Error estimate}
\label{sec: micro_error}

We now prove the estimate \eqref{eq:micro_lin_est:err}. 
Recall $\eM = (\pa _s + \cT) \cM$ which was defined in \eqref{eq:error0}. 
By Lemma \ref{lem:ds_t}, we can write $\eM = \cM p _3 (s, X, \vc)$, where $p _3 (s, X, \vc)$ is a class $\bF ^{-r}$ polynomial of $\vc$ with degree 3 (see Definition \ref{def: polynomial} and Definition \ref{def: bF}).
Therefore,
\begin{align*}
    \cP _m [\cMM ^{-1/2} \eM] = \cs ^3 \cP _m [\cMM ^{1/2} p _3 (s, X, \vc)] = \cs ^3 \cMM ^{1/2} \td p _3 (s, X, \vc),
\end{align*}
where $\td p _3$ is another polynomial of degree 3, because $\cP _m$ is the projection orthogonal  to the space spanned by $\{1, \vc, |\vc| ^2\} \cMM ^{1/2}$.

\begin{proof}[Proof of \eqref{eq:micro_lin_est:err}]
    With the above expression, we compute its weighted derivative:
    \begin{align*}
        D ^{\al, \b} \cP _m [\cMM ^{-1/2} \eM] &= D ^{\al, \b} (\cs ^3 \cM _1 ^{1/2} \td p _3 (s, X, \vc)) \\
        &= \sum _{\substack{\al _1 + \al _2 + \al _3 = \al \\ \b _1 + \b _2 = \b}} C _{\al _i, \b _i} \cdot D ^{\al _1, \b _1} \cMM ^{1/2} \cdot D ^{\al _2, \b _2} \td p _3 (s, X, \vc)  \cdot D ^{\al _3, 0} \cs ^3.
    \end{align*}
    Apply Corollary \ref{cor: Dab-polynomial} on $\td p _3$, \eqref{eq:M11/2} on $\cMM ^{1/2}$, and Lemma \ref{lem:example_fl} (4) on $\cs ^3$, we conclude 
    \begin{align*}
        |D ^{\al, \b} \cP _m [\cMM ^{-1/2} \eM]| \les \ang \vc ^N \cMM ^{1/2} \ang X ^{-r} \cs ^{3} 
    \end{align*}
    for some $N$ depending on $\al, \b$. So 
    \begin{align*}
        & \la D ^{\al, \b} \cP _m [\cMM ^{-1/2} \eM], D^{\al, \b} \tFm \ra_{\cYe}
        \\
        & \qquad =  \iint D ^{\al, \b} \cP _m [\cMM ^{-1/2} \eM] \cdot D ^{\al, \b} \tFm \ang X ^\eta d X d V \\
        & \qquad \les _k \iint  \ang X^{-r} \ang \vc^N \cs^3 \cM_1^{1/2} |D^{\al,\b} \tFm |\ang X ^\eta d X d V\\
        & \qquad \le \left( \iint \ang X^{\eta -2r}\ang \vc^{2N} \cs^{3-\gamma} \cM_1 d X d V \right) ^\f12 \left( \iint \cs^{3+\gamma}\ang \vc^{2+\gamma} | D^{\al,\b} \tFm|^2 \ang X ^\eta d X d V \right) ^\f12 \\
        & \qquad \les _k \left( \int \ang X^{\eta - 2r} d X \right) ^\f12 \| D^{\al,\b}\tFm \| _{\cYE} \\
        & \qquad \le \left( \int \ang X^{\etab - 2r} d X \right) ^\f12 \| D^{\al,\b}\tFm \| _{\cYE}.
    \end{align*}
    We used $\g \le 3$ so $\cs ^{3 - \g} \les 1$. For $\ang X ^{\etab - 2r}$ to be integrable, we need $\etab - 2 r < -3$, $6 (r - 1) < 2 r$, $r < \f32$, which is satisfied by any $r < 3 - \sqrt 3$. We conclude the proof of \eqref{eq:micro_lin_est:err} by H\"older inequality, and the proof of Theorem \ref{thm: micro-Hk-main} is now complete.
\end{proof}

\section{Top order estimates for the cross term}\label{sec:EE_top}

In this section, we estimate the cross terms in the energy estimates, e.g. 
$\cI_i(\tFm)$ in \eqref{eq:lin_euler2} and $\cP_m [V \cdot \na_X \tFM]$ in \eqref{eqn: micro-part}. We estimate them together and exploit an integration by parts to avoid the loss of derivatives. We have the following estimates.
\begin{proposition}\label{prop:EE_cross}
Let $\kp = \f{5}{3}$, $\cI_i$ be the moments defined in \eqref{eq:moments_Fm}, $\cX _\eta ^k$, $\cYe ^k$ be the norms defined in \eqref{norm:Xk} and \eqref{norm:Y}, respectively. 
Let $\tw = (\tu, \tp, \tb)$ and $\tFM = \cF_M(\tw)$ be the macro-perturbation associated with $\tw$. For any $ \g \in [0,2]$, 
$\eta \leq \etab$, and even non-negative integer $k$, we have 
\bseq\label{eq:est_cross}
\begin{align}
    \B| 
    \kp \B\la (\tu, \tp, \tb), (- \cI_1, -\cI_2,  \cI_2) (\tFm)  \B\ra_{\cX _\eta^k}   
    - \B\la \cP _m [V \cdot \na _X \tFM], \tFm \B \ra_{ \cYe ^k} 
    \B| &\les_{k, \eta} \|\tw \|_{\cX _\eta ^k} \| \tFm \|_{\cYE ^k}, \label{eq:est_cross:a} \\
    \B| \B\la \cP _m [(2 \dcm + \td d _\cM) \tFM], \tFm \B \ra_{ \cYe ^k} 
    \B| &\les_{k, \eta} \|\tw \|_{\cX _\eta ^k} \| \tFm \|_{\cYE ^k}. 
    \label{eq:est_cross:aa}
\end{align}
Moreover, we have 
\begin{align}
   \B|  \B\la  (\tu, \tp, \tb), \, (  - \cI_1 , - \cI_2, \cI_2 )(\tFm)  \B\ra_{\cX _\eta^k} 
  \B|  &  \les_k   \| \tFm \|_{\cYE ^{ k+1 } }   \| \tw \|_{\cX_\eta ^k },  \label{eq:est_cross:b}   \\
\B| 
 \B\la \cP _m [(V \cdot \na _X + 2 \dcm + \td d _\cM) \tFM], \tFm \B\ra_{ \cYe ^k} \B| 
    & \les_k
      \| \tFm \|_{\cYE ^k}  \| \tw \|_{\cX _{\eta} ^{k+1}} . \label{eq:est_cross:c}
\end{align}
\eseq
\end{proposition}

In Section \ref{sec:main_cross_term}, we derive the main terms in 
the first inner product in \eqref{eq:est_cross}. In Section \ref{sec:IBP_top}, we prove Proposition \ref{prop:EE_cross} by applying integration by parts.

\subsection{Main terms in the macro cross terms}\label{sec:main_cross_term}

\begin{lemma}\label{lem:cross}
For any multi-indices $\al \in \Znonneg ^3$, we have 
\begin{align*}
    \kp \B( D_X^{\al} \tu \cdot D_X^{\al} \cI_1 +   D_X^{\al} 
    \tp \cdot D_X^{\al} \cI_2
    - \f{3}{2} D_X^{\al} \tb \cdot  D_X^{\al} \cI_2 \B)
    &=  \int ( V \cdot \na_X  D_X^{\al} \tFm )  \cdot  ( D_X^{\al} \tFM + \cR)  d V 
    + \EEs , %
\end{align*}
with lower order terms $\cR$ satisfying 
\[
    \| D^{\leq 1} \cR (X, \cdot) \| _{L ^2 (V)} \les _\al \| D^{\leq |\al|} \tFM \|_{L ^2 (V)},
\]
and the error term $\EEs$ satisfying 
\[
    |\EEs (X)| \les _\al \| D ^{\le |\al|} \tFM \| _{L ^2 (V)} \| D ^{\le |\al|} \tFm \| _{\s}.
\]

\end{lemma}

\begin{proof}
We fix the multi-indices $\al$ and $N \geq 0$. For any function $g$, we denote 
\bseq\label{eq:cross_pf12}
\beq
  \cJ(g) := \cMM^{-1/2} V \cdot \na_X (\cMM^{1/2} g ) = V \cdot \na _X g - \td d _\cM g, 
\eeq
where $\td d _\cM$ was introduced in \eqref{eqn: dcm} with bound \eqref{eqn: Dab-dcm-bound}. Using Leibniz's rule, we obtain 
\begin{align}
    \notag
    D _X^{\al} \cJ(g) &= V \cdot \na _X D _X ^\al g - V \cdot \na _X \log (\varphi _1 ^{|\al|}) \cdot D _X ^\al g - D _X ^\al (\td d _\cM g) \\
    \notag
    &= V \cdot \na_X  D _X^{\al} g + \cO (|\al| |V| \varphi _1 ^{-1} \na _X \varphi _1) \cdot |D _X ^\al g| + \cO (D ^{\le |\al|} \td d _\cM) \cdot |D ^{\le |\al|} g| \\
    & = V \cdot \na_X  D _X^{\al} g + \cO _\al (\ang X ^{-1} \cs \ang \vc ^3) \cdot |D ^{\le |\al|} g| \\
    & = \cO _\al (\ang X ^{-1} \cs \ang \vc ^3) \cdot |D ^{\le |\al| + 1} g|.
\end{align}
Here we used $\vp_1 \asymp \ang X$ and $|\na _X \varphi _1| \les 1$ from \eqref{eq:wg_asym}, $|V| \les \cs \ang \vc$ from \eqref{eq:V_est}, and $|D ^{\le |\al|} \td d _\cM| \les _\al \ang X ^{-1} \cs \ang \vc ^3$ from \eqref{eqn: Dab-dcm-bound}. 
\eseq

Using the identities \eqref{eq:mom_comp_FM}, \eqref{eq:mom_comp_tranFm}, we obtain 
\beq\label{eq:cross_final_pf1}
\kp D _X^{\al} \tu \cdot D_X^{\al} \cI_1 + \kp 
  D_X^{\al}  \tp \cdot  D_X^{\al} \cI_2 - \f{3}{2} \kp 
D_X^{\al} \tb \cdot D_X^{\al} \cI_2 
 = \sum _{i = 0} ^4 
  D_X^{\al} \la \tFM, \Phi_i \ra_V \cdot 
    D_X^{\al} \la \cJ(\tFm), \Phi_i \ra_V .
\eeq
Denote by $\msf{RS}_{ \eqref{eq:cross_final_pf1}}$ the right side. 
Thanks to commutator estimate \eqref{eq:proj_commu}, we know
\begin{align*}
    D _X ^\al \la \tFM, \Phi_i \ra _V &= \la D _X ^\al \tFM, \Phi_i \ra _V + \cO _\al (\| D ^{< |\al|} \tFM \| _{L ^2 (V)}) = \cO _\al (\| D ^{\le |\al|} \tFM \| _{L ^2 (V)}).
\end{align*}
Similarly, by applying \eqref{eq:proj_commu} with $N = 3$ to \eqref{eq:cross_pf12} we get   
\begin{align*}
    D _X ^\al \la \cJ (\tFm), \Phi_i \ra _V &= \la D _X ^\al \cJ (\tFm), \Phi_i \ra _V + \cO _{\al} (\| \ang \vc ^{-3} D ^{< |\al|} \cJ (\tFm) \| _{L ^2 (V)}) \\
    &= \la V \cdot \na _X D _X ^\al \tFm, \Phi_i \ra _V + \cO _{\al} (\ang X ^{-1} \cs \| D ^{\le |\al|} \tFm \| _{L ^2 (V)}) \\
    &= \la V \cdot \na _X D _X ^\al \tFm, \Phi_i \ra _V + \cO _{\al} (\| D ^{\le |\al|} \tFm \| _{\s}).
\end{align*}
In the last step, we used $\ang X ^{-1} \cs \les \Lams ^\f12$ from \eqref{eqn: interpolation-v-holder}. Combine them,
we obtain 
\[
\bal
    \msf{RS}_{ \eqref{eq:cross_final_pf1}}
    & = \sum _{i = 1} ^4 D _X ^{\al} \la \tFM, \Phi _i \ra _V \cdot \la V \cdot \na _X D _X ^\al \tFm, \Phi _i \ra_V + \EEs \\ 
    & = \sum _{i = 1} ^4
    \left( 
        \la D _X ^\al \tFM, \Phi_i \ra _V + \cR _{\al, i}(X)
    \right) \cdot 
    \la V \cdot \na_X D _X ^\al \tFm, \Phi_i \ra _V + \EEs,
\eal 
\]
where $\EEs$ satisfies the bound claimed in the lemma, and $\cR _{\al, i}(X)$ satisfies 
\beq\label{eq:commu_error_X}
    |D ^{\leq 1} \cR _{\al, i}(X) | \les _\al \| D^{\le |\al|} \tFM \| _{L ^2 (V)}.
\eeq
Use the definition of projection $\cP _M$, we rewrite the above identity as 
\[
    \msf{RS}_{ \eqref{eq:cross_final_pf1}}
    = \int (V \cdot \na_X   D_X^{\al}  \tFm) \cdot \left( \cP _M [D _X ^\al \tFM]
    + \sum _{i = 0} ^4 \cR _{\al, i}(X) \Phi _i \right) d V + \EEs,
\]
which justifies the identity in the lemma with $\cR = \sum _{i = 0} ^4 \cR _{\al, i} \Phi _i - \cP _m [D _X ^\al \tFM]$ and $\| D ^{\le 1} \Phi _i \| _{L ^2 (V)} \le C$ due to \eqref{eqn:DabPhi-L2}. Moreover, by applying Corollary \ref{cor: concatenate-derivative} and \eqref{eqn:commutator-projection} we obtain
\begin{align*}
    \| D ^{\le 1} \cP _m [D ^\al \tFM] \| _{L ^2 (V)} &\le \| \cP _m [D ^{\le 1} D ^\al \tFM] \| _{L ^2 (V)} + \cO _\al (\| D ^{\le |\al|} \tFM \| _{L ^2 (V)}) \\
    &\le \| \cP _m [D ^{\le |\al| + 1} \tFM] \| _{L ^2 (V)} + \cO _\al (\| D ^{\le |\al|} \tFM \| _{L ^2 (V)}) \\
    &= \| D ^{\le |\al| + 1} \cP _m [\tFM] \| _{L ^2 (V)} + \cO _\al (\| D ^{\le |\al|} \tFM \| _{L ^2 (V)}) \\
    &= \cO _\al (\| D ^{\le |\al|} \tFM \| _{L ^2 (V)}).
\end{align*}
Combining with \eqref{eq:commu_error_X}, we have proved the bound on $\cR$.
\end{proof}

\subsection{Proof of Proposition \ref{prop:EE_cross}}\label{sec:IBP_top}

Now, we are in a position to prove Proposition \ref{prop:EE_cross}. Since $k$ is even, we assume $k = 2n$. Denote 
\beq\label{eq:IBP_lower}
    J = \| \tFm \|_{\cYE ^k} \| \tFM \|_{\cYe ^k}.
\eeq
By $|V| \ang X ^{-1} \les \Lams ^\f12$ from \eqref{eqn: interpolation-v-holder}, we know
\begin{align}
    \label{eqn: J-bound}
    \iint |V| \la X \ra^{-1} |D_X^{\leq k} \tFm| |D_X^{\leq k} \tFM| \ang X ^\eta d V d X \les _{k, \eta} J.
\end{align}

\paragraph{Main terms in the micro cross terms}

Since the parameter $\nu$ in $\cY$-norm has been chosen in Theorem \ref{thm: micro-Hk-main}, using \eqref{eq:micro_lin_est:cross_d} and $D^{\al, 0} = D_X^{\al}$, 
we estimate the second cross term in \eqref{eq:est_cross:a} as
\beq\label{eq:prop_IBP_pf_micro}
\bal 
- \la \cP _m [V \cdot \na _X \tFM], \tFm \ra_{ \cYe ^k} 
& = - \sum _{|\al| = k} 
\f{ |\al| !}{ \al!} \iint (V \cdot \na _X) D_X^{\al} \tFM \cdot D_X ^{\al} \tFm \ang X ^\eta d V d X +  O_{k, \eta } (J).  
\eal 
\eeq 

Next, we show that the above main term can be further rewritten as 
\beq\label{eq:prop_IBP_pf1}\bal
    & - \sum_{|\al| = k}  \f{|\al|!}{\al!} \int ( V \cdot \na_X ) D_X^{\al} \tFM  
    \cdot D_X^{\al} \tFm \la X \ra ^\eta dV dX \\
    & \qquad = - \int ( V \cdot \na_X ) \D^{n} \tFM  
    \cdot \D^n \tFm \varphi _1 ^{2 k}  \ang X ^\eta dV dX + \cO _{k, \eta}(J).
\eal\eeq
To simplify notation, below, we simplify $\pa_{x_i}$ as $\pa_i$. Applying integration by parts, we obtain 
\[
\bal
    \msf{RS}_{\eqref{eq:prop_IBP_pf1}} 
    &= - \sum_{i_1,.. ,i_n, j_1,.,,j_n \in \{1,2,3\} } 
    \int ( V \cdot \na_X ) \pa_{ i_1 }^2 ... \pa_{ i_n }^2   \tFM  
    \cdot \pa_{j_1  }^2 ... \pa_{ j_n }^2 \tFm \ang X ^\eta \varphi _1 ^{2 k} dV dX + O_{k, \eta}(J) \\
    & =  \sum_{i_1,.. ,i_n, j_1,.,,j_n \in \{1,2,3\}  } 
    \int ( V \cdot \na_X ) \pa_{ i_1 } \pa_{i_2}^2 ... \pa_{ i_n }^2   \tFM  
    \cdot \pa_{i_1} \pa_{ j_1 }^2 ... \pa_{ j_n }^2 \tFm \ang X ^\eta \varphi _1 ^{2 k} dV dX  \\
    & \hspace{6em} + ( V \cdot \na_X ) \pa_{ i_1 } \pa_{i_2}^2 ... \pa_{ i_n }^2   \tFM  
    \cdot \pa_{ j_1 }^2 ... \pa_{ j_n }^2 \tFm \cdot \pa_{i_1} (\varphi _1 ^{2 k} \ang X ^\eta) dV dX   + O_{k, \eta}(J) 
\eal
\]
Since $ |V \pa_{i_1} (\varphi _1 ^{2 k} \ang X ^\eta)| \les _{k, \eta} |V| \ang X ^{-1} \varphi _1 ^{2 k} \ang X ^\eta$ (see \eqref{eq:wg_asym}), 
the integral of the second term is bounded by $\cO _{k, \eta} (J)$ thanks to \eqref{eqn: J-bound}. Similarly, applying integration by parts in $\pa_{j_1}$, we yield 
\[
    \msf{RS}_{\eqref{eq:prop_IBP_pf1}} 
    = - \sum_{i_1,.. ,i_n, j_1,.,,j_n \in \{1,2,3\}  } 
    \int ( V \cdot \na_X ) \pa_{j_1} \pa_{ i_1 } \pa_{i_2}^2 ... \pa_{ i_n }^2   \tFM  
    \cdot \pa_{i_1} \pa_{ j_1 } ... \pa_{ j_n }^2 \tFm \ang X ^\eta \varphi _1 ^{2k} d V d X + O_{k, \eta}( J).
\]
Repeating the above argument, we yield 
\[
    \msf{RS}_{\eqref{eq:prop_IBP_pf1}} 
    = - \sum_{i_1,.. ,i_n, j_1,.,,j_n \in \{1,2,3\}  } 
    \int ( V \cdot \na_X ) \pa_{i_1 }  .. \pa_{i_n} \pa_{j_1} .. \pa_{j_n}    \tFM  
    \, \cdot \,  \pa_{i_1} .. \pa_{i_n} \pa_{j_1} .. \pa_{j_n} \tFm \ang X ^\eta \varphi _1 ^{2k} d V d X + O_{k, \eta}( J).
\]
Using identity \eqref{eq:two_summation} between two summations with $(n, g_1, g_2 ) \rsa (2n, \tFM, \tFm)$, we prove \eqref{eq:prop_IBP_pf1}.

\paragraph{Proof of \eqref{eq:est_cross:a} and  \eqref{eq:est_cross:aa}}

Recall the $\cX$-norm from \eqref{norm:Xk}. Summing Lemma \ref{lem:cross} with $|\al | = k = 2n$ and integrating it over $X$ with weight $\ang X ^\eta$, we obtain 
\beq\label{eq:prop_IBP_pf2}
\bal 
& \kp \B\la  (\tu, \tp, \tb),  \,  (  - \cI_1, -  \cI_2,  \cI_2 )(\tFm)  \B\ra_{\cX _\eta^k} \\ 
& \qquad = -\int \kp \B( \D^{n} \tu \cdot \D^{ n} \cI_1 +   \D^{ n} 
\tp \cdot \D^{ n } \cI_2
- \f{3}{2} \D^{ n } \tb \cdot  \D^{ n } \cI_2 \B) \ang X ^\eta d X  + O_{k, \eta}( J) \\
& \qquad =  - \int  \B( ( V \cdot \na_X  \D^{n} \tFm )  \cdot  (  \D^{n} \tFM + \cR_{2n})  \ang X ^\eta d V \B)
+ \EEs_{2n}(X) \ang X ^\eta d X + O _{k, \eta} (J), 
\eal 
\eeq
where the $J$-term bound the $0$-th order inner product in $\cX$-norm \eqref{norm:Xk} by 
Lemma \ref{lem:cross} with $\al = 0$ and \eqref{eqn: J-bound},  and $\EEs_{2n}, \cR_{2n}$ satisfy the estimates in Lemma \ref{lem:cross} with $|\al| = 2n$. 

Combining the above estimate, \eqref{eq:prop_IBP_pf_micro}, and \eqref{eq:prop_IBP_pf1},  we obtain 
\beq
\bal
\msf{LS}_{\eqref{eq:prop_IBP_pf_micro}}+ 
\msf{LS}_{\eqref{eq:prop_IBP_pf2}} 
& = 
- \int V \cdot \na_X \D^{n} \tFM  
\cdot \D ^n \tFm \vp _1 ^{2 k} \ang X ^\eta d V d X + O _{k, \eta} (J)  \\
& \quad - \int V \cdot \na_X \D ^n \tFm \cdot  ( \vp _1 ^k \D^{n} \tFM + \cR_{2n} ) \vp _1 ^k \ang X ^\eta d V d X
+ \EEs_{2n}(X) \cdot \ang X ^\eta d X \\ 
 & = 
\underbrace{- \int (V\cdot \na_X) (\D^n \tFm \cdot \D^n \tFM) \varphi _1 ^{2 k} \ang X ^\eta d V d X} _{:=I}  + II ,
\eal 
\label{eq:prop_IBP_pf3}
\eeq
where $II$ denotes the error terms 
\[
  II = O_{k, \eta }( J) - \int ( V \cdot \na_X \D^n \tFm ) \cdot \cR_{2n} \varphi _1 ^k \ang X ^\eta d V d X
  - \int \EEs_{2n}(X)  \cdot \ang X ^\eta dX := II_1 + II_2 + II_3.
\]

For the first term $I$, applying integration by parts and using \eqref{eqn: J-bound}, we obtain 
\[
  |I| \les_{k, \eta} \int |V| \cdot |\na_X (\varphi _1 ^{2 k} \ang X ^\eta)| 
  \cdot |\D^n \tFm \D^n \tFM| d X d V
   \les_{k, \eta} J.
\]

For $II _2$, applying integration by parts, $|V| \les \cs \la \vc \ra$ from \eqref{eq:V_est}, and 
using the estimates of $\EEs_{2n},  \cR_{2n}$ in Lemma \ref{lem:cross}, we obtain 
\beq\label{eq:prop_IBP_pf_II2}
\bal
|II_2| & \les \int |V \D^n \tFm | \cdot |\na_X (\cR_{2n} \varphi _1 ^k \ang X ^\eta)| d V d X \\
& \les _{k, \eta} \int |V| \ang X ^{-1} | \varphi _1 ^k \D^n \tFm | \cdot |D ^{\le 1} \cR_{2n}| \cdot \ang X ^\eta d V d X  \\
& \les \int \| D _X^{\le k} \tFm \| _{\s}
\| D _X^{\le k} \tFM \|_{L ^2(V)} \ang X ^\eta d X \les _k J.
\eal
\eeq
Moreover, $|II _3| \les _k J$ directly follows from the bound of $\EEs_{2n}$ in Lemma \ref{lem:cross}. Combining $I$, $II _1$, $II _2$ and $II _3$, we conclude 
\begin{align*}
    \msf{LS}_{\eqref{eq:est_cross:a}} \les _{k, \eta} \| \tFm \|_{\cYE ^k} \| \tFM \|_{\cYe ^k}.
\end{align*}
By Lemma \ref{lem:macro_UPB} we know $\| \tFM \|_{\cYe ^k} \asymp \| \tw \| _{\cX _\eta ^k}$, so \eqref{eq:est_cross:a} is proven. 

Estimate \eqref{eq:est_cross:aa} follows from  \eqref{eq:micro_lin_est:cross_d} and $\| \tFM \|_{\cYe ^k} \asymp \| \tw \| _{\cX _\eta ^k}$.

\paragraph{Proof of \eqref{eq:est_cross:b}, \eqref{eq:est_cross:c}}

The proofs of \eqref{eq:est_cross:b}, \eqref{eq:est_cross:c} are similar, except that we estimate 
the main  terms in \eqref{eq:prop_IBP_pf_micro}, \eqref{eq:prop_IBP_pf2} directly, without using integration by parts. We have estimated the integral of the $ \cR_{2n}, \EEs_{2n}$ terms in the above proof of \eqref{eq:est_cross:a}, e.g.\eqref{eq:prop_IBP_pf_II2}, which are bounded by $J$ and are further bounded by the upper bounds in \eqref{eq:est_cross:b}, \eqref{eq:est_cross:c}. For $(f, g) = (\tFm, \tFM)$ or $(\tFM, \tFm)$, we have
\[
    \left| \int  V \cdot \na_X  \D^{n} f \cdot \D^{n} g \cdot \varphi _1 ^{2k} \ang X ^\eta d V d X \right|
\les \int |V| \ang X ^{-1} | D_X^{\leq k+1} f | 
\cdot | D_X^{\leq k} g | \ang X ^\eta d V d X,
\]
where $k=2 n$. Applying \eqref{eqn: interpolation-v-holder}, \eqref{eq:est_cross:aa}, the Cauchy--Schwarz inequality, and $\| \tFM \|_{\cYe ^k} \asymp \| \tw \| _{\cX _\eta ^k}$, we prove
\eqref{eq:est_cross:b}, \eqref{eq:est_cross:c}.

\section{Nonlinear estimates of collision operator in the energy space}\label{sec:non_Q}

Our main nonlinear estimates are as follows. 
\footnote{
Nonlinear estimates near the \emph{global} Maxwellian $\operatorname e ^{-|V|^2}$ on the torus $ X \in\mathbb{T}^3$, which are similar to \eqref{eq:non_Q:micro}, have been established in \cite[Theorem 3]{guo2002landau}. We refer to Section \ref{sec:compare_guo} for a discussion of the difficulties in our setting.
}

\begin{theorem} \label{theo: N nonlin_est}
    Recall $\etab = -3 + 6 (r - 1)$ from \eqref{wg:X_power} and let $\xeta < \etab$ satisfies 
    \begin{align}
        \label{eqn: etaD-constraint2}
        \etab - \xeta \le \f{(1 + \om) r}{2}.
    \end{align}
    Let $\eta, \eta _1, \eta _2 \in [\xeta, \etab]$ 
    \footnote{
        This constraint is not essential. We impose this range so that the constants related to $\eta, \eta _1, \eta _2 \in [\xeta, \etab]$ in Theorem \ref{theo: N nonlin_est} are bounded by absolute constants.
    }
    satisfy $\eta _1 + \eta _2 \ge \eta + \etab$. There exists an absolute constant $\bar C_{\cN}$ such that for $k \leq \kkl$ with $\kkl = 2 d + 16$, we get
    \bseq\label{eq:non_Q:micro}
    \begin{align}\label{eq:non_Q:micro:a} 
        |\la \cN (f, g), h \ra _{\cYe ^k}| \leq 
        \bar C_{\cN} \| f \| _{\cY _{\eta _1} ^{ \kkl  }} \| g \| _{\cYL{\eta _2} ^k} \| h \| _{\cYE ^k} .
    \end{align}
    For $k > \kkl$, we get
    \begin{align}\label{eq:non_Q:micro:b} 
        |\la \cN (f, g), h \ra _{\cYe ^k}| & \leq 
        \big( 
            \bar C _{\cN} \| f \| _{\cY _{\eta _1}^{\kkl}} \| g \|_{\cYL{\eta _2}^k}   + C_{k} \| f \| _{\cY _{\eta _1}^{k}} \|  g \| _{\cYL{\eta _2} ^{k - 1}}
        \big) \| h \| _{\cYE ^k} \notag \\
        & \leq C_{k} \| f \| _{\cY _{\eta _1} ^k} \| g \| _{\cYL{\eta _2} ^k} \| h \| _{\cYE ^k} .
    \end{align}
    \eseq
Here, the pairing $\ang{\cdot, \cdot} _{\cYe ^k}$ is associated with $\cYe ^k$ norm defined in \eqref{norm:Y}:
    \beq\label{eq:non_recall}
     \la \cN (f, g), h \ra _{\cYe ^k} = \sum _{|\al| + |\b| \le k} \nu ^{|\al| + |\b| - k} 
     \f{|\al|!}{\al!}
     \int \ang X ^\eta \la D ^{\al, \b} \cN (f, g), D ^{\al, \b} h \ra _V d X . 
    \eeq
    Furthermore, if $g = \cP _M g$, then for $k \geq \kkl$
\bseq\label{eq:non_Q:macro}
\beq\label{eq:non_Q:macro:mM}
    |\la \cN (f, g), h \ra _{\cYe ^k}| \les_{k} \| f \| _{\cYL{\eta _1} ^k} \| g \| _{\cY _{\eta _2} ^k} \| h \| _{\cYE ^k} ,
\eeq
    and \beq\label{eq:non_Q:macro:MM2}
        |\la \cN (f, g), h \ra _{\cYe ^k}| \les_{k} 
        \left( 
            \| f \| _{\cYY ^{k-4}} \| g \| _{\cYY ^k}
         + \| f \| _{\cYY ^k} \| g \| _{\cYY ^{k-4}} 
        \right)
        \| h \| _{\cYE ^k}.
    \eeq\eseq
\end{theorem}

 We will apply estimate \eqref{eq:non_Q:micro} with  $g$ being microscopic, 
estimate \eqref{eq:non_Q:macro:mM} with 
$f$ being microscopic, $g$ being macroscopic, and estimate \eqref{eq:non_Q:macro:MM2} with $f, g$ being the macroscopic.

\begin{proof}

Recall the nonlinear term from \eqref{eq:non_recall}. First, we separate the inner product into two parts:
\begin{align*}
    \la \cN (f, g), h \ra _{\cYe ^k} = \udb{ \sum _{|\al| + |\b| \le k} \nu ^{|\al| + |\b| - k} 
         \f{|\al|!}{\al!}
     \int \ang X ^\eta \la \cN(f, D ^{\al, \b} g), D ^{\al, \b} h \ra _V d X}_{:=I} + II
    \end{align*}
where $II$ denotes the lower order terms and satisfies the following estimates due to Lemma \ref{lem: N-derivative-bound}
\bseq\label{eq:non_pf_decomp}
\beq
    \label{eq:non_pf_decomp-II}
    |II| \les_k \sum _{|\al| + |\b| \le k}  \nu ^{|\al| + |\b| - k} \sum _{\substack{\al _1 + \al _2 \preceq \al \\ \b _1 + \b _2 \preceq \b \\ (\al _2, \b _2) \prec (\al, \b)}} \int \ang X ^\eta \cs ^{-3} \|  D^{\alpha_1, \beta_1} f \| _{L ^2 (V)} \| D^{\alpha_2, \beta_2}g \| _\s \| D^{\al, \b} h \| _\s d X .
\eeq
For $I$, since $\cN = \sum \cN_i$ \eqref{N(f,g)}, using Lemma \ref{lem: nonlinear}, we obtain
\beq
    \label{eq:non_pf_decomp-I}
    |I| \les \sum _{ |\al| + |\b| \le k} \nu ^{|\al| + |\b| - k}      \f{|\al|!}{\al!} \int \ang X ^\eta \cs ^{-3} \| f \|_{L ^2 (V)} \| D^{\al, \b} g \|_{\s} \| D^{\al, \b} h \|_{\s} d X ,
\eeq
\eseq
with constant independent of $k$. Following the assumption $\eta _1 + \eta _2 \ge \eta + \etab$ together with $\cs ^{-3} \les \ang X ^{3 (r - 1)}$, we have
\begin{align*}
    \ang X ^\eta \cs ^{-3} \les \ang X ^{3 (r - 1) + \f \eta2 + \f{\eta _1}2 + \f{\eta _2}2 - \f{\etab}2} = \ang X ^{\f \eta2 + \f{\eta _1}2 + \f{\eta _2}2 + \f32}.
\end{align*}

    Let us first bound $II$. If $|\al _1| + |\b _1| \le \kkl - 3$ then by weighted Sobolev embedding \eqref{eq:linf_Y}, we take supremum for $f$:
    \beq\bal
    \label{eq:non_pf_LH}
        &\int \ang X ^\eta \cs ^{-3} \|  D^{\al_1, \b_1} f \| _{L ^2 (V)} \| D^{\alpha_2, \beta_2}g \| _\s \| D^{\al, \b} h \| _\s d X \\
        & \qquad \les \sup _{X} \left\{\ang X ^{\f{\eta _1 + 3}2} \| D ^{\al_1, \b_1} f (s, X, \cdot) \| _{L ^2 (V)} \right\} \int \ang X ^{\f{\eta _2}2} \|  D^{\alpha_2, \beta_2} g \| _\sigma \ang X ^{\f \eta2} \| D ^{\al, \b} h \| _\sigma d X \\
        & \qquad \les _{\eta _1} \| f \| _{\cY _{\eta _1} ^{|\al _1| + |\b _1| + 3}} \| D^{\alpha_2, \beta_2}g \| _{\cYL{\eta _2} } \| D^{\al, \b}h \| _{\cYE} \\
        & \qquad \les \| f \| _{\cY _{\eta _1} ^{\kkl}} \| D^{\alpha_2, \beta_2}g \| _{\cYL{\eta _2} } \| D^{\al, \b}h \| _{\cYE}.
    \eal\eeq
    Otherwise, $|\al _2| + |\b _2| \le k - \kkl + 3 \le k - 5$. We take the supremum for $g$:
    \beq\bal
    \label{eq:non_pf_HL}
        &\int \ang X ^\eta \cs ^{-3} \|  D^{\al_1, \b_1} f \| _{L ^2 (V)} \| D^{\alpha_2, \beta_2}g \| _\s \| D^{\al, \b} h \| _\s d X \\
        & \qquad \les \sup _{X} \left\{\ang X ^{\f{\eta _2 + 3}2} \| D ^{\al_2, \b_2} g (s, X, \cdot) \| _{\s} \right\} \| D^{\al _1, \b _1} f \| _{\cY _{\eta _1} } \| D^{\al, \b} h \| _{\cYE}.
    \eal\eeq
    We recall that the $\sigma$ norm can be bounded from above as 
    $$
        \| f \|^2 _\s \le \int \cs^{\gamma+3} \la \vc\ra^{\gamma+2}f^2 + \cs^{\gamma+5} \la \vc\ra^{\gamma+2}|\na _V f|^2 = \int \Lam (f ^2 + |\cs \na _V f| ^2) \;dV
    $$
   where $\Lam = \cs ^{\g + 3} \ang \vc ^{\g + 2}$. Since $\cs \pa _{V _i} D^{\al, \b} f = D^{\al, \b + \ee _i}f$, by the definition of $D ^{\al, \b}$, we have that 
    \begin{align*}
        \| D^{\al _2,\b _2}g \|^2_\s &\le \int \Lam \B( |D^{\al _2,\b _2} g|^2 + \sum _i |D^{\al _2,\b _2+\ee _i}g|^2 \B) \;dV \\
        &\les \| \Lam ^{1/2} D ^{\al_2, \b_2} g \|^2 _{L ^2 (V)} + \sum _{i} \| \Lam ^{1/2} D ^{\al_2, \b_2 + \ee _i} g \|^2 _{L ^2 (V)}.
    \end{align*}
    We apply \eqref{eq:linf_Ylam} and obtain
    \begin{align*}
        \| D ^{\al_2, \b_2} g \| _\sigma \les _{\eta _2} \ang X ^{-\f{\eta _2 + 3}2} \| g \| _{\cYL{\eta _2} ^{|\al _2| + |\b _2| + 4}}  \les _k \ang X ^{-\f{\eta _2 + 3}2} \| g \| _{\cYL{\eta _2} ^{k - 1}}.
    \end{align*}
    Therefore, we can continue to bound \eqref{eq:non_pf_HL} as 
    \beq\bal
    \label{eq:non_pf_HL-2}
        &\int \ang X ^\eta \cs ^{-3} \|  D^{\al_1, \b_1} f \| _{L ^2 (V)} \| D^{\alpha_2, \beta_2}g \| _\s \| D^{\al, \b} h \| _\s d X \\
        & \qquad \les _{\eta _2} \| D^{\al _1, \b _1} f \| _{\cY _{\eta _1} } \| g \| _{\cYL{\eta _2} ^{|\al _2| + |\b _2| + 4}} \| D^{\al, \b} h \| _{\cYE}.
    \eal\eeq
    Summarizing \eqref{eq:non_pf_LH} for $|\al _1| + |\b _1| \le \kkl - 3$ with $|\al _2| + |\b _2| \le k - 1$, and \eqref{eq:non_pf_HL-2} for $|\al _2| + |\b _2| \le k - 5$, we obtain
    \begin{align*}
        |II| \les _{k, \eta _1, \eta _2} \| f \| _{\cY _{\eta _1} ^{\kkl}} \| g \| _{\cYL{\eta _2} ^{k - 1}} \| h \| _{\cYE ^k} + \| f \| _{\cY _{\eta _1} ^k} \| g \| _{\cYL{\eta _2} ^{k - 1}} \| h \| _{\cYE ^k} \les_{k} \| f \| _{\cY _{\eta _1} ^{\max\{k, \kkl\}}} \| g \| _{\cYL{\eta _2} ^{k - 1}} \| h \| _{\cYE ^k}. 
    \end{align*}

For the first term $I$ in \eqref{eq:non_pf_decomp}, applying estimates \eqref{eq:non_pf_LH} with $(\al_1, \b_1,\al_2,\b_2) = (0,0,\al, \b)$, summing over $\al, \b$, and using the definition of $\cY$-norm \eqref{norm:Y} and the Cauchy--Schwarz inequality, 
we prove
\begin{align}
    \label{eq:non_pf_I}
	|I| \les _{\eta _1} \| f \|_{\cY _{\eta _1} ^{\kkl}} \| g \| _{\cYL{\eta _2} ^k} 
	\| h\|_{\cYE ^k} , 
\end{align}
with absolute constants independent of $k$. 

For $k \leq \kkl$, the constants in $II$ can be treated as independent of $k$ constants, thus we prove \eqref{eq:non_Q:micro:a} and \eqref{eq:non_Q:micro:b} by combining $I$ and $II$.

\paragraph{Proof of \eqref{eq:non_Q:macro}} 
    If $g = \cP _M g$ is macroscopic, then 
    \begin{align*}
        \| g \| _\sigma &\les \cs ^\f{\g + 3}2 \| g \| _{L ^2 (V)}.
    \end{align*}
    Indeed, $g$ is a linear combination of $\{\Phi _i \} _{i = 1} ^5$ which are orthonormal in $L ^2 (V)$, so 
    \begin{align*}
        \| g \| _\sigma &\les \| g \| _{L ^2 (V)} \max _i \| \Phi _i \| _\s \les \cs ^{\f{\g + 3}2} \| g \| _{L ^2 (V)}.
    \end{align*}
    Here we used $\| \Phi _i \| _\sigma ^2 \les \cs ^{\g + 3}$ from \eqref{eqn:DabPhi-sigma}.
    
    Now we integrate in $X$, and recall that $|D ^{\al, 0} \cs ^\f{\g + 3}2| \les \cs ^\f{\g + 3}2$, we conclude
    \begin{align*}
        \| g \| _{\cYL{\eta _2} ^j} &\les \| \cs ^{\f{\g + 3}2} g \| _{\cY _{\eta _2} ^j},
        \quad \forall \ j \geq 0.
    \end{align*}
    Similarly, since $\| f \| _{L ^2} \le \cs ^{-\f{\g + 3}2} \| f \| _\sigma$ for every $f$ 
    and $\g \geq -2$ by Corollary \ref{cor:sig_norm}, we know
    \begin{align*}
        \| f \| _{\cY _{\eta _2} ^k} &\les \| \cs ^{-\f{\g + 3}2} f \| _{\cYL{\eta _2} ^k}.
    \end{align*}
    By definition \eqref{eq:non_nota}, $\cN(\cdot ,\cdot)$ commutes with     
 multiplication by any function $a(X)$. So we can apply \eqref{eq:non_Q:micro:b} and prove \eqref{eq:non_Q:macro:mM}:
    \begin{align*}
        |\la \cN (f, g), h \ra _{\cYe ^k}| &= |\la \cN (\cs ^{\f{\g + 3}2} f, \cs ^{-\f{\g + 3}2} g), h \ra _{\cYe ^k}| \\
        & \les \| \cs ^\f{\g + 3}2 f \| _{\cY _{\eta _1} ^k} \| \cs ^{-\f{\g + 3}2} g \| _{\cYL{\eta _2} ^k} \| h \| _{\cYE ^k} \\
        & \les \| f \| _{\cYL{\eta _1} ^k} \| g \| _{\cY _{\eta _2} ^k} \| h \| _{\cYE ^k}.
    \end{align*}
    
    Finally, note that when $|\al _1| + |\b _1| \le \f{k}{2} \leq k - 8$ or
    $|\al _2| + |\b _2| \le \f{k}{2}\leq k - 8$, we have 
    \begin{align*}
        &\int \ang X ^\eta \cs ^{-3} \|  D^{\al_1, \b_1} f \| _{L ^2 (V)} \| D^{\alpha_2, \beta_2}g \| _\s \| D^{\al, \b} h \| _\s d X \\
        & \qquad \les \int \ang X ^\eta \cs ^\f{\g-3}2 \|  D^{\al_1, \b_1} f \| _{L ^2 (V)} \| D^{\alpha_2, \beta_2}g \| _{L ^2 (V)} \| D^{\al, \b} h \| _\s d X \\
        & \qquad \les \sup _{X} \left\{\ang X ^{\eta + \f{3 - \g}2(r - 1)} \ang X ^{-\f{3 + \eta _1 + \eta _2 + \eta}2} \right\} \left( \| f \| _{\cY _{\eta _1} ^{k-4}} \| g \| _{\cY _{\eta _2} ^k} + \| f \| _{\cY _{\eta _1} ^{k}} \| g \| _{\cY _{\eta _2} ^{k-4}} \right)
        \| h \| _{\cYE ^k}.
    \end{align*}
    The supremum is bounded by 1 when 
    \begin{align*}
        \f{3 - \g}2 (r - 1) + \f{\eta - \eta _2 - \eta _1 - 3}2 \le 0 \iff \eta _1 + \eta _2 & \ge \eta + (3 - \g) (r - 1) - 3 \\
        & = \eta + \etab - (3 + \g) (r - 1) = \eta + \etab - r (\omega + 1).
    \end{align*}
    Therefore, when $\eta _1 = \eta _2 = \xeta$, the constraint is satisfied by $\eta \le \etab$ provided $2 \xeta \ge 2 \etab - r (\omega + 1)$, which reduces to \eqref{eqn: etaD-constraint2}. 
   Thus, applying the above estimate to $I, II$ in \eqref{eq:non_pf_decomp}, we prove \eqref{eq:non_Q:macro:MM2}.   
\end{proof}

\section{Construction of blowup solution}\label{sec:solu}

In this section, we prove Theorem \ref{thm:blowup} by constructing global solutions 
to \eqref{eq:LC_ss} in the vicinity of the local Maxwellian $\cM$ defined in \eqref{eq:local_max}. 
Throughout this section, we perform weighted $H^{2\kk}$ or $H^{2\kk+2}$ energy estimates with the regularity parameter $\kk$ 
chosen in \eqref{def:kk} and use the compact operator $\cK_{\kk} = \cK_{\kk, \xeta}$ \eqref{def:kk} constructed in Proposition \ref{prop:compact}. 
The implicit constants in this section may depend on $\xeta$, $\etab$, and $\kk$, and we omit these dependencies for simplicity.

\subsection{Decomposition of the solution} 
We use $F$ to denote the nonlinear solution to~\eqref{eq:LC_ss}. As in~\eqref{eq:pertb_dec} and 
 \eqref{def:proj}, we denote the perturbation $\td F $ to the profile $\cM$ and its macroscopic 
 $\tFM$ and microscopic parts $\tFm$ as
 \bseq\label{eq:non_decomp}
\beq
\cMM^{1/2} \td F := F - \cM,
\quad  \tFM := \cP_M \td F, \quad  \tFm := \cP_m \td F.
\eeq

We define the weighted hydrodynamic fields $(\trho, \tu, \tp)$ of the perturbation
and $\tb$ via \eqref{eq:pertb_macro}. 
\beq
    ( \td  \rho, \td \UU, \td P)  := \int  \cMM^{1/2} \td F \B(  1,  \f{ V - \bu}{\cs}, \f{|V -\bu|^2}{3\cs^2}  \B) d V ,
    \quad \tb = \trho - \tp, 
\eeq
and denote 
\beq
  \tw = (\tu, \tp, \tb) = \cF_E(\td F).
\eeq
\eseq

Given $\tw$, we construct the macro-perturbation via \eqref{eq:macro_UPB}: $\tFM = \cF_M(\tw)$. We recall that the perturbation $\td F$ solves \eqref{eqn: linearized-td-F} and $\tw$ solves \eqref{eq:lin_euler}. 

We further decompose the macro-perturbation as 
\[
	\tw = \tw_1 + \tw_2, \quad  \tFM = \tFMa + \tFMb
\]
with \footnote{ We emphasize that the $\td \cdot_1$ or $\td \cdot_2$ sub-index denote different parts of the perturbation, they  {\em do not represent} Cartesian coordinates.}
\[
	\tu = \tu_1 + \tu_2, \quad  \tp = \tp_1 + \tp_2, \quad \tb = \tb_1 + \tb_2 ,
\]
so that 
\beq\label{eq:init}
\bga
F   = \cM + \cMM^{1/2} ( \tFm + \tFMa + \tFMb ), 
\quad \td F = \tFm + \tFMa + \tFMb,  \\
   \tw_i = (\tu_i, \tp_i, \tb_i) , 
   \quad \td F_{M, i} = \cF_M( \tw_i ) , 
  \quad i = 1, 2 .
  \ega
\eeq

The field $\tw_i$ are defined as solutions of 
\bseq\label{eq:non_W}
\begin{align}
  \pa_s \tw_1 &= (\cL_{E, s} - \cK_{\kk} ) \tw_1  
  + (\cL_{E,s} - \cL_E) \tw_2   - (\cI_1, \cI_2,  -  \cI_2  )(\tFm)
   - ( \cs^3 \eu, \cs^3  \ep, 0), \label{eq:non_W:a} \\
  \pa_s \tw_2 &= \cL_E \tw_2 + \cK_{\kk} \tw_1  \label{eq:non_W:b}, 
\end{align}
and $\tFm$ solves \eqref{eqn: defn-Lmic}
\beq\label{eq:non_W:micro}
\bal
   \pa _s  \tFm = \Lmic  \tFm - \cP _m [(V \cdot \na _X + 2 \dcm + \td d _\cM) \tFM] 
         + \f1\es \cN (\td F, \td F) - \cP _m [\cMM ^{-1/2} \eM]. 
\eal
\eeq
where we recall $\Lmic$ from \eqref{eqn: defn-Lmic}
\beq
 \Lmic \tFm =  \f1\es \cL _\cM \tFm - \left(\cT + \dcm - \f32 \bcv \right) \tFm+ \cP _M [(V \cdot \na _X - 2 \dcm - \td d _\cM) \tFm] .
\eeq 
\eseq

Let us clarify the definitions of the operators and functions in \eqref{eq:non_W}. 
The operators $\cK_\kk, \cI_i, \cL_E$, $\cL_{E, s},  \Lmic$ are linear. 
We define $\cK_\kk = \cK_{\kk, \xeta}$ \eqref{def:kk} in Proposition~\ref{prop:compact} with parameter $\xeta$, $\cI_1, \cI_2$ in \eqref{def:lin_op}, 
$\cL_{E, s}$ in \eqref{eq:lin_euler}, $\cL_{E}$ in \eqref{eq:lin_euler_limit} with 
\[
  \cL_{E, s} = (\cL_{U, s}, \cL_{P, s}, \cL_{B, s}),
  \quad \cL_E = (\cL_{U}, \cL_{P}, \cL_{B}), 
\]
and $\cT, \Lmic$ in \eqref{def:lin_op} and \eqref{eqn: defn-Lmic}. The error terms $\eu, \erho, \eM$ in \eqref{eq:non_W} are defined in  \eqref{eq:error0} or \eqref{eq:error}.  It is clear, by definition, that a global solution $\tw_1, \tw_2, \tFm$ of \eqref{eq:non_W} provides via \eqref{eq:init} a global solution $F$ of \eqref{eq:LC_ss}.

There are a few important advantages to the decomposition
\footnote{
A similar decomposition was first developed in~\cite{ChenHou2023a} to analyze stable blowup in the 3D incompressible Euler equations,  and then generalized in  \cite{chen2024Euler,chen2024vorticity} for stability analysis 
of implosion in the compressible Euler equations. 
}
 \eqref{eq:init} and \eqref{eq:non_W}. First, the part $\tw_2$, which is used to capture unstable parts, is almost decoupled from the equations of $\tw_1$ \eqref{eq:non_W:a} 
 with a small error $(\cL_{E,s} - \cL_E) \tw_2$ (see Proposition \ref{prop:error_diff}) and $\tFm$ \eqref{eq:non_W:micro} at the linear level, and so we can obtain decay estimates for $\tw_1$ and $\tFm$ directly using energy estimates and the dissipative estimate of $\cL_{E, s} -\cK_{\kk}$ (see ~\eqref{eq:dissip}) and of the linearized operators in \eqref{eq:non_W:micro} (see Theorem \ref{thm: micro-Hk-main}), without appealing to semigroups. Second, by applying energy estimates on $\tw_1$ and $\tFm$, we can estimate the time-dependent linear operators in \eqref{eq:non_W:a} and \eqref{eq:non_W:micro}. 
Third, we can obtain a representation formula (and an estimate) for $\tw_2$ by using Duhamel's formula 
\cite{ChenHou2023a,chen2024Euler}
\footnote{
In general, the projections $\Pi_{\mathsf s}, \Pi_{\mathsf u}$ can lead to a complex-valued solution. 
We restrict to the real part of the semigroup  so that $\tw_2$ is real. 
}
: 
\bseq\label{eq:W2_form}
\begin{align}
\tw_2(s) 
& := \tw_{2,{\mathsf s}}(s) -\tw_{2, {\mathsf u}}(s) + \operatorname e^{s \cL_E}\Bigl(\tw_{2, {\mathsf u}}( 0 ) \bigl(1 - \chi\bigl(\tfrac{y}{8 R_{\xeta}}\bigr) \bigr) \Bigr), 
\label{eq:W2_form:a}\\
\tw_{2, {\mathsf s}} (s) 
& :=  \Re \int_{0}^s \operatorname e^{(s - s^\prime) \cL_E} \Pi_{\mathsf s} \cK_\kk( \tw_1 )(s^\prime) d s^\prime,  
\label{eq:W2_form:b} \\
\tw_{2, {\mathsf u}} (s) 
& :=  \Re  \int_s^{\infty} \operatorname e^{-(s^\prime-s)\cL_E} \Pi_{\mathsf u} \cK_\kk( \tw_1 )(s^\prime) d s^\prime ,
\label{eq:W2_form:c}
\end{align}
where $\chi$ is a smooth radial cutoff function with $\chi(y) = 1$ for $|y| \leq 2/3 $, and $\chi(y) = 0$ for $|y| \geq 1$, and $R_{\xeta}$ is the parameter determined in Theorem \ref{thm:coer_est}, $\Pi_{\mathsf u}$ is the orthogonal projection from $\cXC^{2\kk}$ to $\cX_{\mw{un}}^{2 \kk}$ (see~\eqref{eq:dec_X}-\eqref{eq:dec_Xu}) and $\Pi_{\mathsf s} := {\rm Id} - \Pi_{\mathsf u}$.

It is not difficult to see that \eqref{eq:W2_form} solves \eqref{eq:non_W:b} with initial data taken as 
\begin{equation}
\tw_{2, \init} 
= - \tw_{2,{\mathsf u}}( 0 ) \chi\bigl(\tfrac{y}{8 R_4}\bigr)
=-\chi\bigl(\tfrac{y}{8 R_4}\bigr)
\Re
\int_{ 0}^{\infty} \operatorname e^{- s^\prime \cL_E} \Pi_{\mathsf u} \cK_\kk( \tw_1) (s^\prime) d s^\prime.
\label{eq:W2_form:d}
\end{equation}

\eseq

The detailed representation \eqref{eq:W2_form} shows that $\tw_2$ is computed as a function of 
$\tw_1$; for later purposes it is useful to codify this relation as a map, $\cA_2$, and to denote
\begin{equation}  
\cA_2(\tw_1)  :=
\mbox{Right Side of }\eqref{eq:W2_form:a}
.
\label{eq:T2:def} 
\end{equation}

\subsection{Functional setting and parameters}

In the rest of this section, we will consider power $\eta = \etab$ defined in \eqref{wg:X_power} or $\eta = \xeta$ satisfying \eqref{eq:eta_constraint}:
\beq\bal\label{wg:X_power_recall}
    \etab &= - 3 + 6 (r-1), & 4 \om \cdot r < \etab - \xeta < \f{(1 + \om) r}2.
\eal\eeq

We introduce the spaces $\cZ^{2 \kk + 2}$, which are used for closing nonlinear estimates. Our goal is to perform both weighted $H^{2 \kk}$ and weighted $H^{2 \kk +2}$ estimates on \eqref{eq:non_W}, using the {\em same compact operator} $\cK_\kk$ and the {\em same projections} $\Pi_{\mathsf s}, \Pi_{\mathsf u}$ appearing in~\eqref{eq:non_W} and~\eqref{eq:W2_form:b}-\eqref{eq:W2_form:c}; that is, we do not wish to change $\cK_\kk$ into $\cK_{\kk+1}$ for the weighted $H^{2 \kk + 2}$ bound.

Recall the parameter $\lam_1 < \lame$ chosen in \eqref{eq:eta_constraint}. For some $\vpi^{\pr}_{\kk+1} > 0$ to be chosen sufficiently large,  using Theorem~\ref{thm:coer_est}, Proposition~\ref{prop:compact} (which in particular gives that $\cK_\kk : \cXX^{0} \to \cXX^{2 \kk+ 6}$), and the fact that by definition we have $\|\cdot\|_{\cXX^0} \les_{n, \xeta} \|\cdot\|_{\cXX^n}$, we obtain 
\begin{align*}   
& \vpi_{\kk+1}^{\pr} \la  (\cL_{E,s} - \cK_\kk) f , f \ra_{\cXX^{2\kk}}   +  \la  (\cL_{E, s} - \cK_\kk) f , f \ra_{\cXX^{ 2\kk+2  }}  \\
 &\quad \leq - \lame  \vpi_{\kk+1}^{\pr} \| f \|_{\cXX^{2\kk}}^2 
 +  \bigl( \la \cL_{E, s} f , f \ra_{\cXX^{2\kk + 2}}   
 + \| \cK_\kk f \|_{\cXX^{2\kk+2}} \| f \|_{\cXX^{2\kk+ 2}}\bigr )  \\
 &\quad \leq 
 -\lame \vpi_{\kk + 1}^{\pr} \| f \|_{\cXX^{2 \kk}}^2 
 +  \bigl( - \lame \| f \|_{\cXX^{2\kk+2} }^2 +  C_{\kk, \xeta} \|f\|_{\cXX^0}^2 + C_{\kk, \xeta} \|f\|_{\cXX^0} \|f\|_{\cXX^{2\kk+2}} \bigr)
  \\
 &\quad \leq 
 - \lam_1 \bigl(  \vpi_{\kk+1}^{\pr}  \| f \|_{\cXX^{2\kk} }^2  +  \| f \|_{\cXX^{2\kk+2}}^2 \bigr)
 +
 \bigl(- (\lame - \lam_1) \vpi_{\kk+1}^{\pr} + 
C_{\kk, \xeta, \lam_1}
 \bigr)
 \| f \|_{\cXX^{2\kk} }^2,
\end{align*}
for all $f \in \{ \WW \in \cXX^{2\kk+2} \colon \cL_{E, s} \WW \in \cXX^{2\kk+2}\}$.
Choosing $ \vpi^{\pr}_{\kk+1} >0 $ large enough in terms of $\kk, \xeta$ and $\lame$, so that 
$ - (\lame - \lam_1) \vpi_{\kk+1}^{\pr} + 
C_{\kk, \xeta, \lam_1} < 0 $, 
we obtain the coercivity estimate
\begin{align*}
\vpi_{\kk+1}^{\pr} \la  (\cL_{E, s} - \cK_\kk) f , f \ra_{\cXX^{2\kk}}   + \la  (\cL_{E, s} - \cK_\kk) f , f \ra_{\cXX^{2\kk+2}} 
\leq  
-\lam_1 \bigl(   \vpi_{\kk+1}^{\pr}   \| f \|_{\cXX^{2\kk}}^2  + \| f \|_{\cXX^{2\kk + 2}}^2 \bigr)
.
\end{align*}

In light of the above coercive bounds, with $\vpi_{\kk+1}^{\pr} >0$ chosen as above, we  define the Hilbert spaces $\cZ^{2 \kk+ 2}\subset \cXX^{2\kk + 2}$ according to the inner products
\footnote{
  We apply the $\cZ$-norm only with power $\eta = \xeta$. To simplify the notation, we do not indicate the dependence of $\cZ$  on $\xeta$.  
}
\begin{subequations}
\label{norm:Z}
\begin{align}
&\la f, g \ra_{\cZ^{2\kk+2}} 
:= \la f, g \ra_{\cXX^{2\kk+2}}   
+ \vpi_{\kk+1}^{\pr} \la f , g \ra_{\cXX^{2\kk}}, \;\;\; \| f \|_{\cZ^{2 \kk + 2}}^2 = \la f , f \ra_{\cZ^{ 2\kk + 2}} ,  
\end{align}
\end{subequations}
and obtain with $\lam_1$ determined in \eqref{eq:eta_constraint} that 
\beq\label{eq:coer_Zk}
  \la (\cL_{E, s} - \cK_\kk ) f , f \ra_{\cZ^{2 \kk+2}} \leq - \lam_1 \| f\|^2_{ \cZ^{ 2\kk + 2 }}, \;\;\; 
\eeq
for all $f \in \{ (\UU, P, B) \in \cXX^{2\kk +2} \colon \cL_{E, s}(\UU, P, B) \in \cXX^{ 2 \kk +2}\}$.
Estimate~\eqref{eq:coer_Zk} shows that we can use the same compact operator $\cK_\kk$ to simultaneously obtain coercivity estimates in weighted $H^{2\kk+2}$ and weighted $H^{2\kk}$ spaces. Moreover, we have the following equivalence.

\begin{lemma}\label{lem:norm_equiv_XZ}
For $f \in \cXX^{2\kk+2}$, we have 
\[
  \| f \|_{\cXX^{2\kk+2}} \les  \| f \|_{\cZ^{2\kk+2}} \les   \| f \|_{\cXX^{2\kk+2}} .
\]
\end{lemma}

Since $\kk$ is fixed, we treat constants depending on $\xeta, \etab, \kk$ as absolute constants.

\paragraph{Parameters}

Note that we have fixed $\xeta, \etab$ and the regularity parameter $\kk \geq k_0$. We recall from 
\eqref{eq:eta_constraint}, \eqref{norm:X_lim}, \eqref{eq:decay_stab} 
and~\eqref{eq:coer_Zk}, that the decay rates $
 \lame,  \lams, \lamu$, and $\lam_1$  are given by
\beq\label{eq:decay_para}
\qquad 
\left(\f{2}{3} -\el \right) \rE
<
 \lams < \lamu < \f{2}{3} \rE , \quad \rE <  \lam_1 < \lame .
\eeq

We will only use parameters $\lams, \lamu$ in Lemmas \ref{lem:decay_Km} and \ref{lem:A2} for the semigroup estimates. We use $\lam_1$ and $\lame$ for the energy estimates.

\subsection{Nonlinear stability and the proof of Theorem~\ref{thm:blowup}}\label{sec:non_stab}

We have the following nonlinear stability results.

\begin{theorem}[Nonlinear stability] \label{thm:non}
Let $\kk$ be the regularity index chosen in \eqref{def:kk}. There exists a sufficiently small $\d > 0$ such that
for any initial perturbation $\tw_{1,\init} = (\tu_{1 }(0), \tp_{1}(0), \tb_{1}(0) )$ 
and $\td F_{ m , \init } = \tFm( 0 )$  which are {\em smooth} enough
\footnote{
We require the $\cX^{2\kk+4}$-regularity of $\tw_{1,\init}$, a space which is stronger than $\cZ^{2\kk+2}$, in order to obtain the local-in-time existence of a $\cX^{2\kk + 4}$-solution (see 
Theorem \ref{thm:LWP}); in turn, this allows us to justify a few estimates, e.g.~\eqref{eq:coer_Zk} for $\tw_1$ which requires $\cL_{E, s}(\tw_1) \in \cX^{2\kk+2}$. Note that this regularity requirement is only \textit{qualitative}, and we only use Theorem~\ref{thm:non} with an $C^{\infty}$ initial perturbation (see~\eqref{eq:non_init}) in order to prove Theorem~\ref{thm:blowup}. The {\em quantitative} assumption on the initial data is given by~\eqref{eq:IC:small}.
} 
to ensure $\tw_{1, \init} \in  \cXb^{2 \kk + 4}, \td F_{m, \init} \in  \cYb^{2 \kk + 4} $ 
\footnote{
Since $\xeta < \etab$, using Lemma \ref{lem:Xm_chain} and the definition of $\cYe ^k$ in \eqref{norm:Y}, we also obtain $\tw_{1, \init} \in \cX^{2\kk+4}_{\xeta}, \td F_{m ,\init} \in \cYY ^{2\kk+4}$.
}
and {\em small} enough to ensure
\begin{equation}\label{eq:IC:small}
\bga
\| \tw_{1, \init} \|_{ \cXX^{2\kk+2} } 
+ \|  \td F_{m ,\init} \|_{\cYY^{2\kk+2}} 
< \d^{1/2 }, \quad 
\| \tw_{1, \init} \|_{ \cXX^{2\kk} }  < \d^{2/3 + \el},   \\
\| \tw_{1, \init} \|_{ \cXb^{2\kk+2} } 
 + \|  \td F_{m ,\init} \|_{\cYb^{2\kk+2}}  < \d^{2 \el}  , 
\ega
\end{equation} 
there exists a global solution $\tw_1$ to \eqref{eq:non_W:a} with initial data $\tw_{1, \init}$, a global solution $\tw_2$ to \eqref{eq:non_W:b}  given by~\eqref{eq:W2_form}, and a global solution to \eqref{eq:non_W:micro} with initial data $ \td F_{m ,\init}$ satisfying exponential decay bounds
\begin{subequations}
\label{eq:solu_small}
\begin{align}
\| \tw_1(s) \|_{\cXX^{2\kk+2} } + \| \td F_{m}(s) \|_{\cYY^{2\kk+2} } &\les  \es^{1/2 - \el},  \label{eq:solu_small:a}\\
  \| \tw_1(s) \|_{\cXX^{2\kk}}  & <  \es^{2/3} ,  
 \label{eq:solu_small:b}   \\
  \|  \tw_2(s)  \|_{\cXX^{2\kk+6}} & \les \es^{2/3-\ell} ,
\end{align}
\end{subequations}
and the smallness bound 
\begin{align}
 \| \tw_2( 0) \|_{ \cYb^{ n}  } &  \les_n \d^{2/3 }  , \label{eq:solu_small:W2} \\
 \| \tw_1(s) + \tw_2(s) \|_{ \cXb^{2\kk+2} } +   \| \tFm(s)  \|_{\cYb^{2\kk+2}} & \les \d^{\el}, 
 \label{eq:small_Yeta}
\end{align}
for all $s\geq 0$ and $n \geq 0$. We emphasize that we cannot {\em prescribe} the initial data $
\tw_{2, \init} = \tw_2(0)=  (\tu_{2}(0),\tp_{2}(0), \tb_{2}(0))$; rather, this data is constructed via \eqref{eq:W2_form:d} (simultaneously with $\tw_1$) to lie in a finite-dimensional subspace of $\cXX^{2\kk + 4}$.
\end{theorem}

It is important to obtain extra smallness for the lower order norm  $\| \tw_1 \|_{ \cXX^{2\kk} }$ 
 compared to estimate of the top order norm \eqref{eq:solu_small:a}. See the motivation in \hyperref[sec: step6]{\itshape \underline{Step 6}} in Section \ref{sec:idea}.

\begin{remark}[Exponential decay estimates]\label{rem:critical_no_decay}
We establish temporal exponential decay estimates of perturbation only in norms with faster decay weights, 
e.g. norms $\cXX^{2\kk}, \cYY^{2\kk}, \cZ^{2\kk}$ with parameter $\xeta$, rather than $\etab$. Note that $\xeta < \etab$. 
In the norm $\cXb^{2\kk}, \cYb^{2\kk}$ with parameter $\etab$, we prove smallness instead of
temporal decay estimates. See the motivation in \hyperref[sec: step1]{\itshape \underline{Step 1}} in Section \ref{sec:idea}.
\end{remark}

\begin{remark}[Initial data]\label{rem:data}
The initial data for $F_{\iin} = \cM + \cMM^{1/2} ( \tFm + \cF_M(\tw_1 + \tw_2) )$ is obtained from Theorem~\ref{thm:non} and the decomposition~\eqref{eq:init} at time $s= 0$. In light of Theorem~\ref{thm:non}, we identify the space $X_2$ mentioned in Remark~\ref{rem:initial:data:set} with an open ball in the weighted Sobolev space $\cYb^{2\kk +2}$ defined in~\eqref{norm:Y}.
On the other hand, the space $X_1$ mentioned in Remark~\ref{rem:initial:data:set} consists of functions which are given as the sum of an element $\tw_1$ which lies in open ball in the weighted Sobolev space $\cXb^{2\kk+2}$ (see definition~\eqref{norm:Xk}) and the element 
$\tw_2$ constructed in \eqref{eq:W2_form:d}, which lies in a finite-dimensional subspace of $\cXb^{ 2\kk+ 6}$. 
From \eqref{eq:init_sm} in the proof of Theorem \ref{thm:blowup}, one can construct a finite codimension set of positive initial data $F_{\iin}$.
\end{remark}

We defer the proof of Theorem \ref{thm:non} to Sections \ref{sec:fix_pt_setup}-Section \ref{sec:contra}. Based on Theorem \ref{thm:non},  we are in a position to prove Theorem \ref{thm:blowup}.

\begin{proof}[Proof of Theorem \ref{thm:blowup}]

The proof of Theorem \ref{thm:blowup} consists of a few steps. First, we construct initial perturbation $\tw_1, \tFm$ satisfying the assumptions in Theorem \ref{thm:non} and the initial data $F_{\iin}$ satisfying assumptions in Theorem \ref{thm:blowup}. Then we use the estimates of the perturbation from Theorem \ref{thm:non} to prove the regularity and limiting behaviors of the blowup solution in Theorem \ref{thm:blowup}.

\paragraph{Step 1: Initial data}

Recall from \eqref{eq:non_decomp} that the initial data are given by 
\[
    F_{\mw{in}} = \cM +  \cMM^{1/2} \td F = \cM + \cMM^{1/2}( \td F_{M} + \tFm ) , 
    \quad 
    \tFM = \cF_M( \tw ) = \cF_M( \tw_1 + \tw_2 ) ,
\]
with $\tw_2$ determined by \eqref{eq:W2_form} implicitly. 
Due to finite codimension stability of $\tw$, we cannot prescribe $\tw$ freely. 
To ensure that $F_{\mw{in}}>0$, we first design a specific micro-perturbation  $\FFp$ that is positive for large $|\vc|$ and has much slower decay in $\la \vc \ra$ compared to $\tFM$ so that $F_{\mw{in }}  \approx \cM + \cMM^{1/2} \FFp >0 $. 
 
\paragraph{A specific micro-perturbation}
Recall from \eqref{eq:localmax2} 
\beq\label{eq:localmax2_recall}
  \cMM = \cs^{-3} \mu(\vc),  \quad \cM = \cs^3  \cMM = \mu(\vc).
\eeq

We design a micro-perturbation as 
\footnote{
Since $\cMM$  \eqref{eq:localmax2} and $\cs$ \eqref{eq:Euler_profi_modi} depend on $s$,
We evaluate the functions  in \eqref{eq:Fm_init_pos}  at $s=0$ to construct 
time-independent function $\FFp(X,V)$.
}
\bseq\label{eq:Fm_init_pos}
\beq
\FFp(X, V) = 
\la X \ra^{- l} 
( \cs^{-3/2} \la \vc \ra^{-2} +
( c_1 + c_2 ( |\vc|^2 - \f95  ) ) \cMM^{1/2} ) \B|_{s =0} ,
\quad l >  3(r-1) > 0,
\eeq
and choose $c_1, c_2$ to ensure the orthogonal conditions
\beq\label{eq:Fm_init_pos:b}
  \la \cMM^{1/2} \FFp, h(\vc) \ra_\lvv = 0, \quad h(\vc)  =  1, \  \vc_i , \ |\vc|^2 .
\eeq
\eseq
Since $\FFp$ is radial in $\vc$, we obtain 
$\la  \cMM^{1/2} \FFp  , \vc_i \ra_\lvv = 0$. 
 Using \eqref{eq:Fm_init_pos:b}, \eqref{eq:localmax2_recall}, and a change of variable 
$V = \cs \vc + \bu$, we rewrite the equation \eqref{eq:Fm_init_pos:b} equivalently as 
\[
  0 =\int \cMM^{1/2} \FFp h(\vc) d V
  = \cs^{-3}  \ang X^{-l} \int \mu(  \vc )^{1/2} \B(  \la \vc \ra^{-2} + 
( c_1 + c_2 ( |\vc|^2 - \f95  ) ) \mu(\vc)^{1/2} \B) h(\vc)  
   d \vc
\]
for $h(\vc) = 1, |\vc|^2$. 
Dividing the factor $\cs^{-3}  \ang X^{-l}$ and changing $\vc$ to a dummy variable $z \in \R^3$, we simplify the equations of $c_1, c_2$ as 
\[
\int \mu(  z )^{1/2} \B(  \la z \ra^{-2} + 
( c_1 + c_2 ( |z|^2 - \f95  ) ) \mu(z)^{1/2} \B) h(z)  
   d  z, \quad 
h(z) = 1, |z|^2 .
\]
Since the variance of the Gaussian $\mu(z)$ defined in \eqref{eq:gauss} is $ \kp^{-1} = \f35$, 
by choosing $h(z) = 1$ and using $ \int \mu(z)  ( |z|^2 - \f95  ) d z = 0$, we obtain $c_1 $.  
By choosing $h(z) = |z|^2$,  we further obtain $c_2 $. 
Thus, we obtain constants $c_1, c_2$ \emph{independent of $X$} which satisfy \eqref{eq:Fm_init_pos:b} and 
\beq\label{eq:Fm_init_bd}
  |c_1|,  \ |c_2| \les 1.
\eeq

Using \eqref{eq:Fm_init_pos}, \eqref{eq:Fm_init_bd}, and \eqref{eq:localmax2_recall},  we obtain 
\beq\label{eq:Fm_init_pos1}
\FFp = \ang X^{-l} \cs^{-3/2} \ang \vc^{-2}  + \cE_{F_m},
\quad  |\cE_{F_m}| \les \ang X^{-l} \cs^{-3/2} \ang \vc^{-3}.
\eeq
Since the error part $ \cE_{F_m}$ has a Gaussian decay for large $\vc$, we obtain $\FFp > 0$ for large $|\vc|$. 
Moreover, since $l > 3 (r-1)$, using 
\eqref{eq:M11/2}, \eqref{eq: Dab_vcl}, and a direct computation, we obtain
\beq\label{eq:Fm_pos_Hk}
  \| \FFp \|_{\cYb^n} \les_n 1 , \quad \forall \ n \geq 0. 
\eeq

\paragraph{Initial perturbation}

Based on $\FFp$, we construct small initial perturbation that satisfies the smallness and smoothness assumptions in Theorem \ref{thm:non}.
Consider a family of initial perturbations:  
\bseq\label{eq:init_sm}
\beq
 \tFm = \d_1  \cdot  ( \FFp + \cP_m H), \quad \tw_1 = \d_1 \d^{\ell} \cdot \om ,\quad \d_1 =  b_{\kk}  \d^{2/3 }. 
\eeq
with a small constant $b_{\kk}$ to be chosen and any $H(X, V), \om(X)$ satisfying 
\beq\label{eq:init_sm:a}
\bal 
H &\in \cap_{ n \geq 0} \cYb^{n} \subset C^{\infty},  \quad \| H \|_{\cYb^{2\kk+ 2}} \leq 1, 
\quad  |H(X, V) | \leq \cs^{-3/2} \ang X^{-l} \ang \vc^{-3} , \\ 
\om &\in  \cap_{n \geq 0}\cXb^{ n} \subset  C^{\infty}, \quad \|\om\|_{ \cXb^{ 2 \kk+2} } \leq 1, \quad \supp( \om  ) \in  B(0, 8 R_{\xeta} ) .
\eal 
 \eeq
Clearly, $H=0$ or small $H$ with compact support in $X, V$ satisfies the above assumptions.  

From \eqref{eq:Fm_init_pos}, the definition of $\cY$-norm \eqref{norm:Y}, 
\eqref{eq:Fm_pos_Hk}, and Lemma \ref{lem: commu-derivative-projection}, we have 
\beq\label{eq:init_sm:b}
  \quad  \|\tFm\|_{\cYb^n} \leq \d_1  ( \| \FFp\|_{\cYb^n} 
  + \| \cP_m H\|_{\cYb^n} ) \les_n  \d_1,  \quad 
\| \tw_1(0) \|_{\cYb^n} 
= \d_1 \d^{\ell} \| \om \|_{\cYb^n} \les_n \d_1 \d^{ \ell},  \quad \forall \   n \geq 0 . 
\eeq

We take $\d_1 = b_{\kk}  \d^{2/3}$ with small constant $b_{\kk} = b( \| \FFp\|_{\cYb^{2\kk+4}}  )>0$ depending on $\kk$ so that the smallness assumptions \eqref{eq:IC:small} in Theorem \ref{thm:non} are satisfied for $\tw_1( 0 ), \tFm( 0)$. 

We construct $\tw_2$ via Theorem \ref{thm:non} and use the support of $\tw_2(0)$ \eqref{eq:W2_form} and $\tw_1$ \eqref{eq:init_sm:a} to obtain
\beq\label{eq:init_sm:c}
\bga
 \| \tw_2( 0 ) \|_{\cXb^{n}} \les_n \d^{2/3}, \quad  \| \tw_2( 0 ) \|_{\cXb^{2 \kk + 6}} \les  \d^{2/3} \les \d_1 ,  \\
  \quad \supp(\tw( 0)) , \quad \supp( \tw_1( 0)), \quad \supp(\tw_2( 0)) \subset B( 8 R_{\xeta} ) ,
 \ega 
\eeq
\eseq
for any $n \geq 0$. In particular, we have $\tw_2(0) \in C_c^{\infty}$. 
Using the above construction, we obtain a finite co-dimension set of positive initial data.

\paragraph{Gaussian upper and lower bound of $F_{\iin}$ }
Next, we show that for
\[
\d_1 = b_{\kk} \d^{2/3}
\]
with $\d$ small enough, the initial data satisfy 
\beq\label{eq:non_init}
	  F_{\iin} = \cM + \cMM^{1/2} ( \d_1 (\FFp + \cP_m H) + \tFM)( 0) \geq \tf{1}{2} \cM,
	  \quad \tFM = \cF_M(\tw_1 + \tw_2).
\eeq
for any $X, V$. In particular, $F_{\iin}$ is positive.

Since $R_{\xeta}$ is absolute constant and $l > 3(r-1)$, using the estimates of size and support in \eqref{eq:init_sm} and the embedding in Lemma \ref{lem:prod},  we obtain
\beq\label{eq:FM_init_sm3}
  |\tw(0, X) | \les    \one_{|X| \leq 8 R_{\xeta} }  \ang X^{- \f{\etab + d}{2} } \| \tw( 0) \|_{\cXb^{2\kk} }
    \les (\d_1 \d^{\ell} + \d^{2/3})  \ang X^{-3(r-1)}  \one_{|X| \leq 8 R_{\xeta} }  
     \les \d_1   \ang X^{-l}.
\eeq

Since each basis in \eqref{eq:func_Phi} satisfies $|\Phi_i| \les \ang \vc^2 \cMM^{1/2}$, using 
\eqref{eq:macro_UPB}, \eqref{eq:FM_init_sm3}, and \eqref{eq:localmax2_recall}, we obtain 
\[
\bal 
	|\tFM| \les |\tw( 0, X)| \ang \vc^2 \cMM^{1/2}
\les 
	 \d_1 \ang X^{-l} \ang \vc^{2} \cMM^{1/2}
	 =  \d_1 \ang X^{-l} \cs^{-3/2} \ang \vc^{2} \mu(\vc)^{1/2}
     \les  \d_1 \ang X^{-l} \cs^{-3/2} \ang \vc^{-3}.\\ 
\eal 
\]

Since $\cP_m$ is a projection in $V$, using \eqref{eq:init_sm:a}, \eqref{eq:macro_UPB}, and \eqref{eq:cross_pf2}, we obtain 
\[
\bal 
|\cP_m H| 
& \leq |H| + |\cP_M H| 
\les |H| +  \|\cP_M H \|_{L ^2 (V)}  \ang \vc^2 \cMM^{1/2}  \\
&\les \cs^{-3/2} ( \ang X^{-l} \ang \vc^{-3}+ \ang X^{-l}  \ang \vc^2 \mu(\vc)^{1/2} ) 
 \les \cs^{-3/2} \ang X^{-l} \ang \vc^{-3}.
\eal 
\]

Since $l > 3(r-1)$ \eqref{eq:Fm_init_pos}, using \eqref{eq:localmax2_recall}, we obtain 
\[
	\cM = \mu(\vc) = \cs^3 \cs^{-3} \mu(\vc) \gtr \ang X^{-3(r-1)} \cs^{-3} \mu(\vc) 
\gtr \ang X^{-l} \cs^{-3} \mu(\vc).
\]

Using the above three estimates, \eqref{eq:localmax2_recall}, and \eqref{eq:Fm_init_pos1}, we obtain 
\beq\label{eq:Gauss_low_initial}
\bal
	  F_{\iin} - \tf{1}{2} \cM & = \tf12 \cM + \cMM^{1/2} ( \d_1 \FFp + \d_1 \cP_m H + \tFM) \\
   &  \geq C_1 \ang X^{-l} \cs^{-3} \mu(\vc)
   + \cs^{-3/2} \mu(\vc)^{1/2}
   (  \d_1 \ang X^{-l} \cs^{-3/2}  \ang \vc^{-2}    
   - \d_1 |\cE_{F_m}| - \d_1 |\cP_m H| -  |\tFM| ) \\
  & \geq 
  C_1 \ang X^{-l} \cs^{-3} \mu(\vc)
   + \cs^{-3/2} \mu(\vc)^{1/2} \B(  \d_1 \ang X^{-l} \cs^{-3/2}  \ang \vc^{-2}  
   - C \d_1 \ang X^{-l} \cs^{-3/2}   \ang \vc^{-3}  \B) \\
& = \ang X^{-l} \cs^{-3} \mu(\vc)^{1/2}
(  C_1 \mu(\vc)^{1/2} + \d_1 \cdot \ang \vc^{-2} - C \d_1 \ang \vc^{-3} ),
\eal
\eeq
for some absolute constants $C_1, C >0$.  The above term is positive for $\d_1 = b_{\kk} \d^{2/3}$ with $\d$ small enough. Thus, we prove \eqref{eq:non_init}: $F_{\iin}( 0) > \f{1}{2} \mu(\vc) =\f{1}{2} \cM$.

Next, we prove the uniform Gaussian decay bound for $F_{\iin}$. Using $\cs|_{s=0} \gtr_{R_0} 1$ \eqref{eq:dec_S} and $\vc = \f{V- \bu}{\cs}$ and $V = v$ \eqref{eq:SS} at $t=0$, we obtain  $|\vc| \gtr c_1 \ang v - c_2$. 

Using \eqref{eq:fxv_est1} to be shown with $\al=\b=0$ and the bound on $|\vc|$, we obtain 
 Gaussian decay in $v$ that is uniformly in $x$
\beq\label{eq:Gauss_upper} 
\bal 
   |F_{\iin}( X, V) - \mu(\vc)| \les \mu(\vc)^{1/4} ,
  \quad  |F_{\iin}(X, V) | \les \mu(\vc) + \mu(\vc)^{1/4} 
\les \exp( - C \ang v^2).
\eal 
\eeq 

Recall  $c_v \equiv \bcv,   c_x \equiv \bcx, c_f \equiv 0 $ and the self-similar relation from \eqref{eq:SS} 
\beq\label{eq:ss_recall}
s = - \log( 1-t), \  X = \f{x}{(1-t)^{\bcx}}, \  V = \f{v}{(1-t)^{\bcv}}, 
 \quad  f(t, x, v) = F(s, X, V)
 = \cM +  \cMM^{1/2 } \td F. 
\eeq 
Since $f(0, x, v) = F(0, x, v)$, we prove that the initial data $f$ satisfies the Gaussian decay estimate in $v$ \eqref{eq:Gauss_upper}. 

\paragraph{Initial hydrodynamic fields}

Recall the mass, momentum, and energy density  $ (\mass, \mome, \eee)$ from \eqref{eq:density}. Using the self-similar relation \eqref{eq:SS} and \eqref{eq:pertb_macro} and $\td \rho =\tb + \tp$, we obtain 
\beq\label{eq:SS_macro}
\bal
 (\mass, \mome, \eee )(t,x) & =	\int f(t, x ,v)(1, v , |v|^2) d v  \\ 
 & = (T-t)^{3 \bcv} \int F(s, X, V) (1, (T-t)^{\bcv} V, (T-t)^{2 \bcv} |V|^2 ) d  V \\ 
  & = \B(  (T-t)^{3 \bcv } ( \brho + \td \rho),  \, (T-t)^{4 \bcv} 
  ( \brho \bu + \td \rho \bu + \cs \tu), \\
  & \quad  \qquad  (T-t)^{5 \bcv} (  3 \brho \bar \Th_s + \brho |\bu|^2
+  3 \cs^2 \tp + 2 \bu \cdot \tu \cs + \td \rho |\bu|^2  )  
    \B).
\eal
\eeq

Using $f_{\iin} = F_{\iin},x = X, v= V$ \eqref{eq:ss_recall} at $t =0$,  \eqref{eq:Gauss_low_initial} and \eqref{eq:dec_S}, we estimate the initial density 
\[
 \mass( 0, x) = \int F_{\iin} d V  \geq \f{1}{2} \int \cM dV = \f{1}{2} \brho|_{s=0}(  x) = \f{1}{2} \cs^3|_{s=0} 
 \gtr R_0^{-3 (r-1)} ,
\]
for any $x$. 
We prove  $\mass_{\iin} \geq \rm{constant}>0$ in Theorem \ref{thm:blowup}.
Using \eqref{eq:SS_macro} and \eqref{eq:solu_small_point} to be shown below, 
we obtain that the initial data have uniformly bounded hydrodynamic fields. 
We have proved all the properties of initial data in Theorem \ref{thm:blowup}.

\vspace{0.1in}

\paragraph{Step 2: Asymptotically self-similar blowup}

Since the initial perturbation $(\tw_1, \tFm)$ satisfies \eqref{eq:init_sm}, which implies \eqref{eq:IC:small} 
and \eqref{eq:LWP_small_IC}, using Theorem \ref{thm:non}, we construct a global solution 
$(\tw_1 \tw_2, \tFm)$ to \eqref{eq:non_W} with estimates \eqref{eq:solu_small} and \eqref{eq:small_Yeta}.  Since system \eqref{eq:non_W} is equivalent to the linearized Landau equation \eqref{eq:lin},  $F = \cM + \cMM^{1/2} \td F$
with $\td F = \tFm + \cF_M(\tw_1 + \tw_2) $ is a global solution to the Landau equation \eqref{eq:LC_ss} with $\td F(s) \in \cYb^{2\kk+2}$, arising from the initial perturbation $\td F(0)$. 
By requiring $\d$ small, \eqref{eq:init_sm} implies that $\td F(0)$ also satisfies 
\eqref{eq:LWP_ass_IC} for any $ k \geq 0$. Therefore, by  uniqueness of solutions, the global solution $\td F$ constructed in Theorem \ref{thm:non} and the \emph{local} solution constructed in Corollary \ref{prop:LWP_landau} from the same initial perturbation $\td F(0)$ are the same. Since estimates \eqref{eq:init_sm:b} and \eqref{eq:init_sm:c} imply $ \td F(0) \in \cap_{n\geq 0 } \cYb^n$ and since $F(0, X, V) > \f{1}{2} \cM(0,X, V)$ by \eqref{eq:Gauss_low_initial}, 
using Proposition \ref{prop:LWP_landau}, we further obtain that $\td F(s) \in \cap_{n \geq 0} \cYb^n \subset C^{\infty}$ and $F$ satisfies a Gaussian lower bound \eqref{eq:Gauss_lower} with $l =0$.

Using estimates \eqref{eq:X_linf} and \eqref{eq:linf_Y_point} in Lemma \ref{lem:prod} with $d =3$ and \eqref{eq:solu_small}, we obtain 
\bseq\label{eq:solu_small_point}
\beq
\bal
	|\tw(s, X)| & \les \ang X^{ - \f{ \xeta + 3 }{2} } \| \tw\|_{\cXX^{2\kk}}
	\les \es^{1/2- \ell} \ang X^{ - \f{ \xeta + 3 }{2} }, \\
	|\td F(s, X, V)| & \les \cs^{- \f{3}{2}} \ang X^{- \f{\xeta+ 3 }{2}} \| \td F \|_{\cYY^{2d}}
	\les \cs^{- \f{3}{2}} \ang X^{- \f{\xeta+ 3 }{2}} ( \| \tw\|_{\cXX^{2d}} + \| \tFm \|_{\cYY^{2d}}
	) \\
 & 	\les  \es^{1/2- \ell}
	\cs^{- \f{3}{2}} \ang X^{- \f{\xeta+ 3 }{2}} , \\ 
	| D^{\leq {2\kk- 2 d}  } \td F(s, X, V)| & \les \cs^{- \f{3}{2}} \| D_V^{\leq d} D^{\leq 2\kk - 2 d} \td F(X, \cdot ) \|_{L ^2 (V)}  \\
    & \les \cs^{- \f{3}{2} } \ang X^{ - \f{\etab + 3}{2}} \| \td F \|_{\cYb^{2\kk}}
    \les \d^{\ell}  \cs^{- \f{3}{2} } \ang X^{-3(r-1)} \les \d ^\ell \cs ^\f32.
\eal
\eeq
In the last step, we used $\cs \gtrsim \ang X ^{-(r - 1)}$ from Lemma \ref{lem:profile}. 

From the definition \eqref{eq:Euler_profi_modi}, we obtain 
\beq
\lim_{s \to \infty} \rs = \infty, \quad 
 \lim_{s \to \infty} \cs = \bar \sc, 
 \quad \lim_{s \to \infty} (\brho, \ths, \bp) = (\bar \rho, \bth, \bar P),
 \quad \lim_{s \to \infty} \cM = \cM_{\bar \rho, \bar \UU, \bar \Th} ,
\eeq
\eseq
where $\cM$ is the time-dependent local Maxwellian defined in \eqref{eq:localmax2},
and $\cM_{\bar \rho, \bar \UU, \bar \Th}$ is defined in \eqref{eq:local_max}.

Since $t\to 1^-$ is equivalent to $s = -\log(1-t) \to \infty$, for fixed $X,V$, using the decay estimates in \eqref{eq:solu_small_point}, $\es \to 0$ as $s \to \infty$, and the relation \eqref{eq:SS_macro}, we establish the blowup asymptotics 
\[
\bal 
 \lim_{t \to T^-}  ( (T-t)^{-3 \bcv} \mass,
 (T-t)^{-4 \bcv} \mome, (T-t)^{-5 \bcv } \eee) ( t, (T-t)^{\bcx} X, (T-t)^{\bcv} V)  & = ( \bar \rho, \bar \rho \bu,  \bar \rho (  3 \bar \Th + |\bu|^2 ) ), \\
	\lim_{t \to T^-}  f(t, (T-t)^{\bcx} X, (T-t)^{\bcv} V)
	& = \cM_{\bar \rho, \bar \UU, \bar \Th}.
\eal 
\]
Using $\bcx = \f{1}{r}, \bcv = \f{1}{r}-1$ from \eqref{eq:Euler_profi}, we prove \eqref{eq:blow_asym:macro}, 
\eqref{eq:blow_asym:micro}. Since $\bcv = \f{1}{r} - 1 <0$ \eqref{eq:Euler_profi}, 
the mass $\mass$, moments $\mome$, and the energy $\eee$ blow up at $t=1$.
We prove results  (b), (c) in Theorem \ref{thm:blowup}.

\vspace{0.1in}
\paragraph{Step 3. Estimates of blowup solution  $f$}

In this step, we study the limiting behavior of the blowup solution and its regularity away from $x = 0$.

Recall $ F = \cM + \cMM^{1/2} \td F$. Using Leibniz rule, $|D ^{\al, \b} \cMM ^{1/2}| \les _{\al, \b} \ang \vc ^{|\b| + 2 |\al|} \cMM ^{1/2} $ from \eqref{eq:M11/2}, $|D ^{\al, \b} \cM| \les _{\al, \b} \ang \vc ^{|\b| + 2 |\al|} \cM$ from \eqref{eq:DabM}, together with \eqref{eq:solu_small_point}, we obtain for any $|\al| + |\b| \le 2 \kk - 2 d$ that
\beq\label{eq:fxv_est1}
\bal 
|D^{\al, \b} (F - \cM)| &\les |D^{\preceq (\al, \b)} \cMM ^{1/2}| \cdot
|D^{\le 2 \kk - 2 d} \td F| \\
& \les \ang \vc ^{|\b| + 2 |\al|} \cMM ^{1/2} \cdot \d ^\ell \cs ^{\f32} 
= \d ^\ell \ang \vc ^{|\b| + 2 |\al|} \mu (\vc) ^{1/2} \les _{\al, \b} \d ^\ell \mu ^{1/4} (\vc), \\
|D^{\al, \b} F| &\les _{\al, \b}  |D^{\al, \b} \cM| + |D^{\al, \b} (F - \cM)| \\
& \les \ang \vc ^{|\b| + 2 |\al|} \cM  + \d ^\ell \ang \vc ^{|\b| + 2 |\al|} \mu (\vc) ^{1/2} \\
& = \ang \vc ^{|\b| + 2 |\al|} \left(\mu (\vc) + \d ^\ell \mu (\vc) ^{1/2} \right) \les _{\al, \b} \mu^{1/4}(\vc) .
\eal 
\eeq

For fixed $x$ and fixed $v$, using the self-similar relation \eqref{eq:ss_recall}, we obtain 
\beq\label{eq:vring_limit}
 \vc(t, x, v) = \f{ V- \bu(X)}{\cs(X)}
 =  \f{ v -  (1-t)^{ 1/r - 1 } \bu( \f{x}{(1- t )^{1/r } } )}{  (1-t)^{ 1/r - 1 }  \cs( \f{x}{(1-t)^{ 1/r }  } ) },
 \quad  \vc(t,0,v) =  \f{ v }{ \bc(0) \cdot (1-t)^{1/r - 1 } } ,
\eeq
For $x=0$ and $|\al| = |\b| =0$, using \eqref{eq:fxv_est1}, $X=0$ \eqref{eq:ss_recall}, 
and $\e_0 = \d$ \eqref{eq:EE_para1}, we prove 
\[
 |f(t, 0, v) -\mu(\vc) |
 = |F(s, 0, V) - \mu(\vc) | \les  \e_0^{\ell} \mu^{1/2}(\vc), \quad \vc =  \f{ v }{ \bc(0) \cdot (1-t)^{ 1/r - 1 } } ,
\]
and obtain the first estimate in Remark \ref{rem:tail_flat}.
By choosing $\e_0 = \d$ small enough, we obtain  $C \e_0^{\ell} \mu(0)^{1/2} \leq \e_0^{\ell/2} \mu(0)$.
For fixed $v $ and $x=0$, since $\vc \to 0$ as $ t \to 1^-$, we prove the second estimate in 
Remark \ref{rem:tail_flat}.

\paragraph{Smoothness away from $x=0$}
We first derive that the limit of $\vc$ is $\mr v$ as $t \to 1 ^-$.
For fixed  $x \neq 0$, using the asymptotics of $(\bu, \bc)$ \eqref{eq:UC_refine_asym} , we obtain
\[
\bal 
 \lim_{ t \to 1^-}  (1-t)^{ \f{1}{r} - 1 } \bu \left( \f{x}{(1- t )^{ 1/r } } \right)
 & = \lim_{ t \to 1^-}   (1-t)^{ \f{1}{r} - 1 } \left( C_{\bu} \ee _R  \B| \f{x}{(1- t )^{1/r } }\B|^{-(r-1)}  
 + O \left(   \B| \f{x}{(1- t )^{ 1/r } }\B|^{-2 r + 1}  \right)
 \right) \\ 
 & = C_{\bu} \ee _R |x|^{-r+1}.
 \eal 
\]
Next, we compute a similar limit for  $\cs$ \eqref{eq:Euler_profi_modi}.  Since $(1 -t)^{\bcx} = \operatorname e^{-\bcx s}$ \eqref{eq:ss_recall}, using the definition of $\rs = R_0 \operatorname e^{\bcx s}$ in \eqref{eq:S_radial}, \eqref{eq:EE_para1}, we obtain
\[
 \f{X}{\rs} = \f{ x \operatorname e^{\bcx s}}{R_0 \operatorname e^{\bcx s}} = \f{x}{R_0},
 \quad \chi_{\rs}(X) = \chi \left( \f{X}{\rs} \right) = \chi \left( \f{x}{R_0}\right) = \chi _{R _0} (x).
\]
Thus, for $x \neq 0$, using \eqref{eq:UC_refine_asym} and $\rs = R_0 \operatorname e^{\bcx s} = R_0 (1-t)^{-\f{1}{r}}$, we obtain
\beq\label{def:cRx}
\bal 
& \lim_{t \to 1^-}  (1-t)^{ \bcv } \cs \left( \f{x}{(1- t )^{1/r} } \right) \\
&\qquad = \lim_{t \to 1^-}  (1-t)^{ \f{1}{r} - 1 } \left( \bc \left( \f{x}{(1- t )^{1/r} } \right)
\chi _{R _0} (x) + \left(1 - \chi _{R _0} (x) \right) (R_0 (1-t)^{- \f{1}{r} } )^{-(r-1)} \right) \\ 
&\qquad  = C_{\bc} |x|^{-r+1} \chi _{R _0} (x)
+ \left(1 - \chi _{R _0} (x) \right)  R_0^{-(r-1)} := \msf{c}_{R_0}(x).
\eal 
\eeq
By definition $x = \f{X}{(1-t)^{\bcx}} $ and using \eqref{eq:dec_S}, for any fixed $x \neq 0$ and $t \in [0,1)$,  we obtain 
\beq\label{eq:est_cRx}
\bga 
\msf{c}_{R_0}(x) \asymp \min \{|x|, R_0\}^{-(r-1)} , \\ 
 (1-t)^{\bcv} \cs(X) \gtr 
(1-t)^{\bcv} \rs^{-(r-1)} 
\gtr R_0^{-(r-1)},  \\ 
  (1-t)^{\bcv} \cs(X)  \les (1-t)^{\bcv}  ( X^{-(r-1)} + \rs^{-(r-1)})
= |x|^{-(r-1)} + R_0^{-(r-1)} \les  \msf{c}_{R_0}(x).
\ega 
\eeq
Combining the above estimates, for fixed $x \neq 0$ and $v$,  using \eqref{eq:vring_limit}, we derive 
\beq\label{def:small_vring}
 \lim_{t \to 1^-} \vc(t,x, v) = \f{v - C_{\bu} \ee _R |x|^{-r+1}}{ \msf{c}_{R_0}(x) } := \mr{v}(x, v).
\eeq
Using \eqref{eq:fxv_est1} with $\al = \b = 0$, we have 
\[
\limsup_{t \to 1^-} |f (t, x, v) - \mu(\mr{v})| = \limsup_{t \to 1^-} |F (s, X, V) - \mu(\vc)|
\les \d^{\ell} \limsup_{t \to 1^-} \mu(\vc)^{1/2} =
\d^{\ell} \mu^{1/2}(\mr{v}) .
\]
With $\e _0 = \d$ \eqref{eq:EE_para1} we establish estimate  \eqref{eq:limit_solu} in Remark \ref{rem:limit_solution}.

Now we derive pointwise estimates for higher derivatives.
Let functions $g, G$ be related by $g(x, v) =G( \f{x}{ (1-t)^{\bcx}} , \f{v}{ (1-t)^{\bcv}}) =  G(X, V)$. Using the definition of $D^{\al, \b}$ \eqref{eq:deri_wg}, $\vp_1 \asymp \ang X$ from Lemma \ref{lem:wg}, 
and the self-similar relation \eqref{eq:ss_recall}, we obtain
\beq\label{eq:fxv_est2}
 |D^{\al, \b} G(X, V) | \asymp_{\al}  \ang X^{|\al|} \cs^{|\b|} |\pa_X^{\al} \pa_V^{\b} G(X, V)|
 = (1-t)^{\bcx |\al| + \bcv |\b|} \ang X^{|\al|} \cs^{|\b|} | \pa_x^{\al} \pa_v^{\b} g(x, v) |.
\eeq
For fixed $x \neq 0$ and $v$, using the above asymptotics, and the property that $\mu(\cdot) \in C^{\infty}$, we obtain 
\beq\label{eq:fxv_est3}
\bal  
  \lim_{t \to 1^-} (1-t)^{\bcx |\al| + \bcv |\b|} \ang X^{|\al|}  \cs^{|\b|}(X)
  &= |x|^{ |\al| }   (\msf{c}_{R_0}(x))^{|\b|} , \\ 
 \lim_{  t\to 1^-} \pa_x^{\al} \pa_v^{\b}  \mu(\vc(t, x, v)) 
 &= \pa_x^{\al} \pa_v^{\b}  \mu( \mr{v}). 
 \eal 
\eeq
Thus, for $x \neq 0, v\in \R^3$, using \eqref{eq:fxv_est1}-\eqref{eq:fxv_est3} with $(g, G) \rsa (f, F)$, we prove  %
\[
\bal 
& \limsup_{t \to 1^-} |x|^{ |\al| }   (\msf{c}_{R_0}(x))^{|\b|}  |\pa_x^{\al} \pa_v^{\b} (f - \mu(\mr{v})) |
 =  \limsup_{t \to 1^-} |x|^{ |\al| }   (\msf{c}_{R_0}(x))^{|\b|}  |\pa_x^{\al} \pa_v^{\b} (f - \mu(\vc)) |  \\
& \qquad \les_{\al} \limsup_{t \to 1^- } |  D^{\al, \b} ( F - \mu(\vc)) |
\les_{\al, \b} \limsup_{t \to 1^-}  \d^{\ell}  \mu^{1/4}( \vc)  
= \d^{\ell} \mu^{1/4}( \mr{v}).
\eal 
\]
Dividing $|x|^{ |\al| }   (\msf{c}_{R_0}(x))^{|\b|}$ in the above estimate and using  \eqref{eq:est_cRx}, we prove 
\[
\limsup_{t \to 1^-}  |\pa_x^{\al} \pa_v^{\b} (f - \mu(\mr{v})) |
\les_{\al, \b} \d^{\ell} |x|^{-|\al|} \msf{c}_{R_0}(x)^{-|\b|} \mu^{1/4}(\mr{v}) .
\]
This yields higher-order estimates for the error.

Recall $R_0 = \e_0^{-\ell_r}$ from \eqref{eq:EE_para1}.  Using \eqref{eq:fxv_est1}, \eqref{eq:fxv_est2}, and then \eqref{eq:est_cRx}
$ (1-t)^{\bcx} |X|= |x|$, for any $|x| \neq 0$, we obtain 
\[
\bal 
  \mu^{1/4}(\vc) & \gtr_{\al, \b} | D^{\al, \b} F| 
  \gtr_{\al, \b}   (1-t)^{\bcx |\al| + \bcv |\b|}  |X|^{|\al|} \cs^{|\b|} | \pa_x^{\al} \pa_v^{\b} f(t,x, v) |\\ 
&  \gtr_{\al, \b} |x|^{ |\al|} R_0^{-(r-1) |\b|} | \pa_x^{\al} \pa_v^{\b} f(t,x, v) |
\gtr_{\al, \b,\e_0}  |x|^{ |\al|}  | \pa_x^{\al} \pa_v^{\b} f(t,x, v) | .
\eal 
\]
Recall $\mu(\cdot)$ from \eqref{eq:gauss}. Since $ |\f{\bu(X) }{\cs(X)}| \les 1$, using \eqref{eq:vring_limit}, the upper bound on $\cs$ in \eqref{eq:est_cRx}, 
and \eqref{def:small_vring},  we obtain 
\[
 |\vc| \geq C_0 | \f{v}{ (1-t)^{1/r-1} \cs(X) } | - C_2
 \geq  C_1 |\f{v}{ \msf{c}_{R_0}(x) } | -  C_2
 \geq C_3 | \mr{v}| - C_4, \quad \mu( \vc)^{1/4} \leq c \exp(- C |\mr{v} |^2). 
\]
for some absolute constants $C_i > 0$. Combining the above estimates, we prove 
\beq\label{eq:smooth_away_0:pf}
 | \pa_x^{\al} \pa_v^{\b} f(t,x, v) | \les_{\al, \b, \e_0} |x|^{-|\al| } \exp( - C |\mr{v}|^2)|,
\eeq
and obtain the estimate in Remark \ref{rem:smooth_away_0}, which implies estimate \eqref{eq:smooth_away_0:1} in 
result (a)  in Theorem \ref{thm:blowup}.

\vspace{0.1in}

\paragraph{Step 4: Regularity of the blowup solution}

Using the self-similar transform \eqref{eq:SS}, we obtain 
\[
  \| f(t, \cdot , v = 0) \|_{ \dot C_x^{\al} } = (1 - t)^{- \al \bcx} \|  F(s, \cdot , V =0 ) \|_{ \dot C_X^{\al} } .
\]

From the profile equations \eqref{eq:profile_eqn} and $\bc(X) \gtr 1 $ for $|X| \leq 1$,  there exists some $ |X_0| \leq 1$ with $\bu(X_0) \neq 0$ (otherwise \eqref{eq:profile_eqn} implies $\bu(X) = 0, \bc(X) = 0 , \forall |X| \leq 1$). 
Using \eqref{eq:solu_small_point}, \eqref{eq:fxv_est1} and $F = \cM + \cMM^{1/2} \td F = \mu(\vc) + \cMM^{1/2} \td F$,  by choosing 
$\d$ small enough, we obtain
\beq\label{eq:fxv_est4}
\bal 
\|  F(s,\cdot, V=0) \|_{ \dot C^{\al}}
& \geq  | F(s, X_0,  0) -  F(s, 0, 0 ) | 
\geq  | \cM(s, X_0,  0) -  \cM(s, 0 , 0 ) | -  C \d ^\ell \\
& \geq \f{1}{2} \B| \mu\left( \f{ \bu(X_0) }{ \bc(X_0)} \right) - \mu(0) \B| \geq \bar c  ,
\eal 
\eeq
uniformly in $s$ for some $\bar c>0$.  Combining the above estimates, we prove that $ \| f(t, \cdot , v=0) \|_{C_x^{\al} }$ blows up for any $\al >0$. 

Using \eqref{eq:ss_recall} and \eqref{eq:fxv_est1}, 
we obtain $  \| f(t) \|_{L^{\infty}} = \| F(s) \|_{L^{\infty}} \les 1 $. Thus, $ \| f(t) \|_{L^{\infty}}$ is uniformly bounded for $t \in [0, 1)$. 
Estimate \eqref{eq:smooth_away_0:1}
in result (a) in Theorem \ref{thm:blowup} has been proved in \eqref{eq:smooth_away_0:pf}.

Next, we fix $v$. Using the self-similar relation 
\eqref{eq:ss_recall}, \eqref{eq:fxv_est1}, \eqref{eq:fxv_est4}, and then taking $1-t$ small enough, we obtain 
\[
\bal 
& \left| F \left( s, X_0, \f{v}{ (1-t)^{1/r-1}} \right) -
  F\left( s, 0, \f{v}{ (1-t)^{1/r-1}} \right) \right| \\ 
 & \qquad \geq  | F( s, X_0, 0 ) -
  F( s, 0, 0)| - C \B|\f{v}{ (1-t)^{1/r-1}} \B| 
\geq \bar c -  C \B| \f{v}{ (1-t)^{1/r-1}} \B| \gtr \f{1}{2} \bar c.
\eal 
\]
For $1-t>0$ small enough, using the mean-value theorem, and  \eqref{eq:ss_recall}, we prove 
\[
 \sup_{|x| \leq (1-t)^{\bcx}} |\na_x f(t,x ,v) |
 \geq \sup_{ |X| \leq 1} (1-t)^{-\bcx} \left|
 \na_X F \left(s,X, \f{v}{ (1-t)^{1/r-1}}\right)
 \right|
\gtr  (1-t)^{ - \bcx }c. 
\]
For any fixed $v$, using $\bcx = \f{1}{r}$ and taking $t \to 1^-$, we prove the gradient blowup result in result (a) in Theorem \ref{thm:blowup}. 
We complete the proof of Theorem \ref{thm:blowup}.
\end{proof}

\subsection{Setup of the fixed point problem }\label{sec:fix_pt_setup}

In subsections \ref{sec:fix_pt_setup}-\ref{sec:contra}, our goal is to prove Theorem~\ref{thm:non}.  

Since the formula for $\tw_2$ (see~\eqref{eq:W2_form:c}) involves the future of the solution $\tw_1$, and since $\tw_2$ enters the evolution \eqref{eq:non_W:a} for $\tw_1$ and 
\eqref{eq:non_W:micro} for $\tFm$ through the nonlinear term, we cannot solve for the perturbation  $ \tw_1$ directly. Instead, we reformulate \eqref{eq:non_W:a} as a fixed point problem. We fix the initial data $\tw_1 |_{s= 0} = \tw_{1,\init} \in \cX^{2\kk + 4}$ sufficiently smooth, and sufficiently small such that \eqref{eq:IC:small} holds. We define the space $Y$
\begin{align}\label{norm:fix}
    \| \tw_1 \|_{Y}
  := \sup_{ s \geq 0} \es^{- 2/3}   \| \tw_1 \|_{\cXX^{2\kk}}  ,
\end{align}
and energy 
\footnote{
  Note that we only define $\cZ$ norm in \eqref{norm:Z}
with the power $\xeta$ and do not consider similar norm with power $\etab$. Therefore, we consider $ 
\| \tw_1 \|_{\cXb^{2\kk+2}}^2$ in the definitions \eqref{norm:fix} rather than some $\cZ$ norms with power $\etab$ of $\tw_1$.
}
\bseq\label{eq:non_E}
\begin{align}
  E_{\kk+1, \xeta}(\tw_1, \tFm) & := \kp  \| \tw_1 \|_{ \cZ^{2\kk+2}}^2 +  \| \tFm \|_{\cYY^{2\kk+2}}^2 , \label{EE:onto_sub}  \\
  E_{\kk+1, \etab}(\tw, \tFm) & := \kp \| \tw \|_{ \cXb^{2\kk+2}}^2 + 
\| \tFm \|_{\cYb^{2\kk+2}}^2 . \label{EE:ontp_cric} 
\end{align}
\eseq

Note that $E_{\kk + 1, \xeta}$ controls $\tw_1$, while $E_{\kk + 1, \etab}$ controls $\tw$ rather than $\tw_1$. The parameter $\kp = \f{5}{3}$ in \eqref{norm:fix} and \eqref{eq:non_E} relates to the coupled estimates in Proposition \ref{prop:EE_cross}. Showing $\tw_1 \in Y$  implies that 
the norm $\cXX^{2\kk} $ of the perturbation $\tw_1$ 
decay with a rate $\es^{2/3}$ as $s\to \infty$.

Next, we define an operator $\cA$ (see~\eqref{eq:map_A}), whose fixed point (see~\eqref{eq:the:fixed:point}) is the desired solution of~\eqref{eq:non_W:a}, \eqref{eq:non_W:micro}. We remark that throughout the remainder of this proof, we distinguish the $W = (\UU, P, B)$-components of an input of a map  (e.g.~$\cA$, or $\cA_2$) by variables with a ``hat'' (e.g.~$\wh \WW_1 = ( \wh \UU_1, \wh P_1, \wh B_1 )$), and the output of these maps by variables with a ``tilde'' (e.g.~$\tw_1 = (\td \UU_1, \tp_1, \tb_1)$). With this notational convention in place, the two-step process is:
\begin{itemize}
\item first, for  $\hw_1 \in Y$, we define 
\begin{equation}
\tw_2 = \cA_2 (\hw_1),
\label{eq:non_fix:W2}
\end{equation} 
where the linear map $\cA_2$ is defined by \eqref{eq:T2:def}, via~\eqref{eq:W2_form};
\item second, we define $\tw_1$ as the solution of a modified version of~\eqref{eq:non_W:a} 
and $\tFm$ as the solution of \eqref{eq:non_W:micro}, namely
\begin{subequations}
\label{eq:non_fix}
\begin{align}
\pa_s \tw_1 & =  (\cL_{E, s} - \cK_\kk ) \tw_1   + (\cL_{E,s} - \cL_E) \tw_2  - (\cI_1, \cI_2,  -  \cI_2  )(\tFm)
   - ( \cs^3 \eu, \cs^3  \ep, 0),  \label{eq:non_fix:a} \\
     \pa _s  \tFm &= \Lmic  \tFm - \cP _m [(V \cdot \na _X + 2 \dcm + \td d _\cM) \tFM] 
          + \f1\es \cN( \td F, \td F ) \notag  - \cP _m [\cMM ^{-1/2} \eM]. 
  \label{eq:non_fix:b}
\end{align} 
where $\td F = \tFm + \tFM$, and we construct the macro-perturbation $\tFM$ associated with $\tw_1 + \tw_2$ using the linear operator $\cF_M$ \eqref{eq:macro_UPB}
\beq
\tFM  = \cF_M( \tw_1 + \tw_2 )
= \cF_M( \tw_1 + \cA_2( \hw_1 ) ) .
\eeq
We choose the initial data as in Theorem~\ref{thm:non}
\beq
\tw_1 |_{s = 0} = (\td \UU_1(0), \tp_1(0), \tb_1(0) ),
\quad \tFm |_{s = 0}  = \tFm(0).
\eeq
\end{subequations}

\end{itemize}

For initial data $\tw_1(0), \tFm(0)$ satisfying \eqref{eq:IC:small} with $\d$ small enough,
applying \eqref{eq:est_init} for $\tw_2$ (to be established) and $\kk \geq \kkl$, we have $\td F =\cF_M(\tw_1 + \tw_2) + \tFm \in \cYb^{2\kk+4} $ and 
 \beq\label{eq:LWP_small_IC}
 \|\td F_{\iin} \|_{\cYb^{\kkl}} \leq C ( \| \tw_{1, \iin} + \tw_2(0) \|_{\cXb^{\kkl}}
 + \| \td F_{m,\iin} \|_{\cYb^{\kkl}} ) \leq C \d^{2\ell} <  \ze_2,
 \eeq
 with $\ze_2$ chosen in Theorem \ref{thm:LWP}. Thus, $\td F_{\iin}$ satisfies assumption \eqref{eq:LWP_ass_IC}. Applying Theorem \ref{thm:LWP} with $\sss = 1$, 
 we construct local-in-time solutions $\td F(s) \in \cYb^{2\kk+4}$ and 
 $(\tw_1(s), \tFm(s) ) \in \cXb^{2\kk+4} \times \cYb^{2\kk+4}$.  We will prove estimate \eqref{eq:EE_onto:a} in Proposition \ref{prop:onto}, which ensures that $\| \td F(s)\|_{\cYb^{\kkl}}$ remains small. 
 Therefore, using the continuation criterion in Theorem \ref{thm:LWP}, we justify the global existence 
 of a solution $(\tw_1(s), \tFm(s) ) \in \cXb^{2\kk+4} \times \cYb^{2\kk+4}$ to \eqref{eq:non_fix}.

Concatenating the two steps given above defines a map with input $\hw_1$ and output the solution of \eqref{eq:non_fix}: 
\beq
 (\tw_1, \tFm ) 
\stackrel{\eqref{eq:non_fix}}{=} \cA( \hw_1)
=\bigl( \cA_W( \hw_1), \cA_{ \rm{mic} }( \hw_1) \bigr) .
\label{eq:map_A}
\eeq
Denoting by  $\cA_{W}$ the restriction of $\cA$ to the $\WW$-components, we have thus reformulated the system \eqref{eq:non_W}
as a {\em fixed point problem}: find $\tw_1$ such that 
\begin{equation}
 \tw_1 = \cA_W( \tw_1) ,
 \label{eq:the:fixed:point}
\end{equation}
with $\tw_2$ and $\tFm$ computed as $\cA_2( \tw_1 )$ and $\cA_{\rm{mic}}(\tw_1)$ , respectively.

By definition of $\tw_2$ \eqref{eq:non_fix:W2} and \eqref{eq:W2_form}, $\tw_2$ satisfies \eqref{eq:non_W:b} with the forcing $\cK_\kk \hw_1$:
\bseq\label{eq:non_fix:W}
\beq\label{eq:non_fix:Wa}
   \pa_s \tw_2 = \cL_E \tw_2 + \cK_\kk \hw_1.
\eeq
Combining the above equation and \eqref{eq:non_fix}, we derive the equation of $\tw$
\beq
  \pa_s \tw = \cL_{E, s} \tw  + \cK_\kk (\hw_1 - \tw_1)  - (\cI_1, \cI_2,  -  \cI_2 )(\tFm)
   - ( \cs^3 \eu, \cs^3  \ep, 0) .
  \eeq
\eseq

The proof of Theorem~\ref{thm:non} reduces to establishing that the operator  $\cA_W$ is a contraction with respect to the norm in~\eqref{norm:fix}, in a vicinity of the zero state as in the statement of Theorem~\ref{thm:non}. 
The proof of Theorem~\ref{thm:non} is broken down in two steps, according to Proposition~\ref{prop:onto} (which shows that the map $\cA_W$ maps the ball of radius $1$ in $Y$ into itself and into a space with higher regularity characterized by $E_{\kk+1, \xeta}, E_{\kk+1, \etab}$), 
and Proposition~\ref{prop:contra} (which shows that $\cA_W$ is a contraction for the topology $Y$).
 
\begin{proposition}\label{prop:onto}
Recall $\es = \d \operatorname e^{ -\rE s }$ from \eqref{eq:EE_para1}, the energy 
$E_{\kk, \xeta}, E_{\kk, \etab}$ from \eqref{eq:non_E}, and the space $Y$ from \eqref{norm:fix}. Let $(\tw_1, \tFm) =\cA( \hw_1 )$ and $\el$ be the parameter to be chosen in \eqref{eq:EE_para_k}.
There exists a positive $\d_0 \ll_m 1$ such that for any $\d < \d_0$ and any $ \hw_1 \in Y$ with $\| 
\hw_1\|_{Y} < 1$, we have 
\bseq\label{eq:EE_onto}
\beq\label{eq:EE_W2_est}
  \| \tw_2(0) \|_{ \cXb^n} \les_n \d^{2/3} ,
  \quad   \|\tw_2(s) \|_{\cXX^{2\kk+6}}   \les  \es^{2/3 -\ell},
\eeq 
and
\begin{gather}
 \| \tw_1(s) \|_{ Y } \leq C \d^{ \el } < 1,  \quad 
  E_{\kk+1, \xeta}(s) < \es^{ 1 - 2 \el },
 \quad E_{\kk+1, \etab}(s) < \d^{ 2 \el }, \label{eq:EE_onto:a} \\
  \int_{0}^s   \f{1}{\e _\tau} \| \tFm(\tau) \|^2_{\cYQ ^{ 2\kk+2}}  d \tau \les \d^{ 2 \el },   \label{eq:EE_onto:b} 
 \end{gather}
 \eseq
for all $s\geq 0$, $n\geq0$, and any $\th \in [0, 1 )$, with implicit constants independent of $s, \th$ and $\d$. 
\end{proposition}

Note that the norms $\| \tw_1 \|_{\cXb^{2\kk+2}}, \| \tFm \|_{\cYb^{2\kk+2}} $. may not decay in time. See Remark \ref{rem:critical_no_decay}.

\begin{proposition}\label{prop:contra}
There exists a positive $ \d_0\ll_k 1$ such that for any $\d < \d_0$ and any pairs 
$\hw_{1, a}, \hw_{1, b} \in Y$ with $\| \hw_{1, a}\|_{Y} < 1$ and  $\| \hw_{1, b}\|_{Y} < 1$
, we have 
\[
 \|\cA_{W}( \hw_{1, a} ) - \cA_W(\hw_{1, b}) \|_{Y} < \tfrac{1}{2} 
  \| \hw_{1, a} - \hw_{1, b} \|_{Y}  .
\]
\end{proposition}

From Proposition~\ref{prop:onto} and Proposition~\ref{prop:contra} we directly obtain:

\begin{proof}[Proof of Theorem~\ref{thm:non}]
Propositions~\ref{prop:onto} and~\ref{prop:contra}  allow us to apply a Banach fix-point theorem for the operator $\cA_W$, in the ball of radius $1$ around the origin in the space $Y$ \eqref{norm:fix}; this results in a unique fixed point $\tw_1$ in this ball, as claimed in~\eqref{eq:the:fixed:point}. Upon defining 
$ \tFm := \cA_{\rm{mic}}( \tw_1 )$ and $\tw_2 := \cA_2( \tw_1)$, by construction we have that 
$\tw_1$ solves \eqref{eq:non_W:a} and $\tw_2 $ solves \eqref{eq:non_W:b}. Using the definitions of the $Y$ norm in~\eqref{norm:fix} and the energies \eqref{eq:non_E} and Proposition~\ref{prop:onto}), 
we deduce that~\eqref{eq:solu_small}, \eqref{eq:solu_small:W2}, and \eqref{eq:small_Yeta} hold, thereby concluding the proof of Theorem~\ref{thm:non}. 
\end{proof}

\

The following subsections are dedicated to the proof of Propositions~\ref{prop:onto} and~\ref{prop:contra}. In subsection~\ref{sec:S2}, we obtain suitable estimates for the linear map $\cA_2$; in particular, in Lemma~\ref{lem:decay_Km} we demonstrate a smoothing effect  for $\tw_2$, which allows us to overcome the loss of a spatial derivative due to the term $ \na \tw_2$
and $( \cL_{E, s} - \cL_E ) \tw_2$ present in the first equation of \eqref{eq:non_W}. In subsection~\ref{sec:fix_pt} we prove Proposition~\ref{prop:onto}, while in subsection~\ref{sec:contra}, we prove Proposition \ref{prop:contra}.

\subsubsection{Estimates on $\cA_2$}\label{sec:S2}
Recall the decomposition~\eqref{eq:dec_X} of $\cXC^{2 \kk}$ into stable and unstable modes. In light of definitions~\eqref{eq:W2_form:b} and~\eqref{eq:W2_form:c}, we establish the following decay and smoothing estimates for the stable and unstable parts of $\cK_\kk$:
\begin{lemma}\label{lem:decay_Km}
For any real-valued $f \in \cXX^{2\kk}$, we have 
\[
\bal
 \| \Re \, ( \operatorname e^{s \cL_E} \Pi_{\mathsf s} \cK_\kk f ) \|_{\cXX^{2\kk+6}} 
 & \les \operatorname  e^{-\lams s  }  \| f \|_{\cXX^{2\kk}}, \\
 \| \Re \, ( \operatorname e^{ - s \cL_E} \Pi_{\mathsf u} \cK_\kk f ) \|_{\cXX^{2\kk+6}} 
& \les \operatorname e^{\lamu s }  \| f \|_{\cXX^{2\kk}},
\eal
\]
for all $s\geq 0$, where $\lamu$ and $\lams$ are as in~\eqref{eq:decay_para}.
\end{lemma}

The proof uses the semigroup estimates in \eqref{eq:decay_stab}, \eqref{eq:decay_unstab}. Since the proof is the same as \cite[Lemma 4.5]{chen2024Euler}, we omit it and refer the proof to 
 \cite{chen2024Euler}.

Using Lemma~\ref{lem:decay_Km} and the fact that $\cL_E$ generates a semigroup,  we obtain a direct estimate for the operator $\cA_2$, as defined in~\eqref{eq:W2_form}.

\begin{lemma}\label{lem:A2}

Recall $\lams < \f{2}{3}\rE$ from \eqref{eq:decay_para}. For $\hw_1 \in Y$ and for all $s\geq 0$ we have
\[
 \|\cA_2( \hw_1 )(s) \|_{\cXX^{2\kk+6}} \les 
\operatorname  e^{-\lams s} 
 \sup_{s \geq 0} \operatorname e^{ \f{2}{3} \rE s  } \| \hw_1(s) \|_{\cXX^{2\kk}}.
 \]

\end{lemma}

Lemma \ref{lem:A2} is an analog of \cite[Lemma 4.6]{chen2024Euler}, which was proved using decay estimates essentially the same as 
those in Lemma \ref{lem:decay_Km} and Proposition \ref{prop:far_decay}. Here, Lemma \ref{lem:decay_Km} corresponds to \cite[Lemma 4.5]{chen2024Euler}, 
 Proposition \ref{prop:far_decay} corresponds to 
 \cite[Proposition 3.8]{chen2024Euler}, and the parameters $(\lams, \lamu, \f{2}{3}\rE )$ with $\lams < \lamu < \f{2}{3} \rE$ correspond to $\eta_s < \eta < \lam_1 $ in  \cite[Sections 4.3, 4.4]{chen2024Euler}. The proof of Lemma \ref{lem:A2} is the same as that of \cite[Lemma 4.6]{chen2024Euler}. 
A minor difference is that the map $\cT_2$ used in \cite[Lemma 4.6]{chen2024Euler} depends on two variables $(\UU, \S)$, while $\cA_2$ we use here depends on three variables $( \UU, P, B )$. 
We omit the proof 
of Lemma \ref{lem:A2} and refer to \cite{chen2024Euler} for more details. Since $\lams < \lamu < \f{2}{3}\rE$ by \eqref{eq:decay_para}, we obtain a decay rate $\operatorname e^{-\lams s} $ in the above Lemma.

\subsection{Proof of Proposition \ref{prop:onto} }\label{sec:fix_pt}

In this section, we prove Proposition \ref{prop:onto} 
via a bootstrap argument. Recall the  notations from the beginning of Section \ref{sec:solu}. Per the assumption of Proposition~\ref{prop:onto}, let $  \|  \hw_1 \|_{\cXX^{2\kk}}  < \es^{2/3} $. 
Define $\tw_2$ using \eqref{eq:non_fix:W2}, and then define $\tw_1$ as the solution of~\eqref{eq:non_fix}. Denote 
\[
  \tw_2  = \cA_2( \hw_1 ),
  \quad \tw = \tw_1 + \tw_2 .
\]

\paragraph{Bootstrap assumptions}

We assume the following bootstrap bounds  
 \bseq\label{eq:EE_boot}
 \begin{align}
  E_{\kk+1, \xeta}(s)   & <  \es^{ 1 - 2 \el } =   \d^{ 1 - 2\el } \operatorname e^{- (1 - 2 \el) \rE  s  }, \label{eq:EE_boot:a} \\
   \| \tw_1 \|_{ \cXX^{2\kk} } & < \es^{ \f23 } =  \d^{ \f{2}{3} } \operatorname e^{- \f{2}{3} \rE  s  } , \label{eq:EE_boot:b} \\
  \quad E_{\kk+1, \etab}(s) & < \d^{ 2 \el } , \label{eq:EE_boot:c}
  \end{align}
  \eseq
for $s \in [0, \bar s], \bar s > 0$, where $\el = 10^{-4}$ is chosen in \eqref{eq:idea_size} and is a small parameter satisfying
\beq\label{eq:EE_para_k}
 0 < \el = 10^{-4} < \min \left\{ \f{2}{3} - \f{1}{2}, \ \f{1}{10} \right\}.
\eeq

In the following sections, our goal is to show that there exists  $\d_0 = \d_0(\kk, \xeta, \etab)$ such that these bounds can be improved  for any $\d < \d_0$ and $s \in [0, \bar s]$. Since we have fixed $\kk, \xeta, \etab, \ell_i$, the following implicit constants $C$ 
or those in the notation ``$\les$" can depend on $\kk, \xeta,\etab $ but independent of $s, \d, \es$.

\paragraph{Estimate of $\tw_2$}
Using Lemma \ref{lem:A2}, $\es = \d \operatorname e^{-\rE s }$ \eqref{eq:EE_para1}, 
$ \| \hw_1\|_{\cXX^{2\kk}} < \es^{2/3}$, and $\lams > ( \f23 - \ell ) \rE$ by \eqref{eq:decay_para}, we obtain 
\bseq\label{eq:est_W12}
\beq\label{eq:est_W12:a}
  \|\tw_2 \|_{\cXX^{2\kk+6}}   \les  
  \operatorname e^{-\lams s  } \sup_{s \geq 0} \operatorname e^{\f23 \rE s } \| \hw_1 \|_{\cXX^{2\kk}}
\les \d^{ \f23 }   \operatorname e^{-\lams s  }  \les      \es^{2/3-\ell} .  \\
\eeq
Under the bootstrap assumption \eqref{eq:EE_boot}, using \eqref{eq:non_E} and Lemmas \ref{lem:macro_UPB},  we estimate 
\beq
\bal
  \| \tFM\|_{\cYY^{2\kk+2} }
   & \les \| \tw\|_{\cXX^{2\kk+2}}
\les   E_{\kk+1, \xeta}^{1/2} +\|\tw_2 \|_{\cXX^{2\kk+2}}  \les \es^{ 1/2 - \el }, \\
  \| \tFM\|_{\cYY^{2\kk } } 
   & \les   \| \tw_1 \|_{\cXX^{2\kk}} +  \| \tw_2 \|_{\cXX^{2\kk}}
  \les \es^{ 2/3 -\ell } .
  \eal
\eeq
\eseq

Recall the initial data $\tw_2(0) , \tw_{2,{\mathsf u}}( 0 )$ from \eqref{eq:W2_form:c}, \eqref{eq:W2_form:d}, which 
 depend on $\hw_1$. Using  Lemma \ref{lem:decay_Km}, 
$ \| \hw_1 \|_{ \cXX^{2\kk}} \leq \es^{2/3}$, 
$\lamu <  \f{2}{3}  \rE$ by \eqref{eq:decay_para}, and $\es = \d \operatorname e^{-\rE s}$ \eqref{eq:EE_para1}, we obtain 
\[
\bal 
\| \tw_{2,{\mathsf u}}( 0 ) \|_{\cXX^{2 \kk+6}}
 &= \left\| \int_{ 0}^{\infty} \operatorname e^{- \cL_E s } \Pi_{\mathsf u} \cK_\kk( \hw_1) (s) d s \right\|_{\cXX^{2 \kk+6}}
 \les \int_0^{\infty} \operatorname e^{ \lamu s} \es^{ \f23} d s
\les \d^{\f23 } \int_0^{\infty} \operatorname e^{ (\lamu - \f23 \rE ) s} d s 
\les \d^{ \f23 }.
\eal 
\]
From the definition of $\tw_{2,{\mathsf u}}( 0 )$ and the projection $ \Pi_{\mathsf u}$ in \eqref{eq:W2_form:c}, 
$\tw_{2,{\mathsf u}}( 0 ) $ can be written as $\Re \, g$ for some $ g \in \cX_{\mathrm{un}}^{2\kk}$. Using \eqref{eq:X_un_smooth} and the above estimate, for any $n\geq 0$, we obtain
\[
\| \tw_{2,{\mathsf u}}( 0 ) \|_{\cXX^{n}}
= \| \Re \, g \|_{\cXX^{n}}
\les_n \| \Re \, g \|_{\cXX^{2\kk }}  = \| \tw_{2,{\mathsf u}}( 0 ) \|_{\cXX^{2\kk}} \les_n
\| \tw_{2,{\mathsf u}}( 0 ) \|_{\cXX^{2\kk+6}} \les_n \d^{2/3}.
\]
From \eqref{eq:W2_form:d}, since $\tw_2(0) = - \tw_{2,{\mathsf u}}( 0 ) \chi\bigl(\tfrac{y}{8 R_4}\bigr)$ has compact support $\supp( \tw_2(0)) \subset B(0, 8 R_{\xeta})$ and $\chi$ is a smooth cutoff function, using the definition of the $\cX^k_{\eta}$ norms \eqref{norm:Xk} and the above estimate, for any $n\geq0$, we obtain 
\bseq\label{eq:est_init}
\beq
\| \tw_2(0) \|_{ \cXb^n}
\les_n \| \tw_2(0) \|_{ \cXX^n} \les_n \d^{2/3}. 
\eeq
where the implicit constants can depend on $R_4, \kk$ (these parameters are fixed throughout this section) and $n$. 
Using \eqref{eq:est_init} and the assumption \eqref{eq:IC:small} on $\tw_1, \tFm$, we yield 
\beq
   \| \tw( 0 ) \|_{\cXb^{2\kk+2} } 
    + \| \tFm (0) \|_{\cYb^{2\kk+2}}
   \les    \| \tw_1( 0 ) \|_{\cXb^{2\kk+2} } 
   +  \| \tw_2( 0 ) \|_{\cXb^{2\kk+2} } 
   + \d^{2 \el}
   \les \d^{2\el}.
\eeq
\eseq
Combining \eqref{eq:est_W12} and \eqref{eq:est_init}, we prove estimates \eqref{eq:EE_W2_est} on $\tw_2$.

\begin{remark}[Size of perturbations]\label{rem:size}
The typical size of perturbations $\|\tw_1 \|_{\cXX^{2\kk+2}}, \| \tFm\|_{\cYY^{2\kk+2}}$ is $\es^{ \f12 -\el}$. The terms $\|\tw_2\|_{\cXX^{2\kk+6} }, \| \tw_1 \|_{\cXX^{2\kk}}$ (not $\| \tw_1 \|_{\cXX^{2\kk+2}}$) satisfy much smaller bounds $\es^{ \f23-\ell}$. From Remark \ref{lem:para_bd}, we have $\rs \ll \es^2$. 
The reader can essentially treat the terms as if
\[
\|\tw_2\|_{\cXX^{2\kk+6} } \approx 0,\quad  \| \tw_1 \|_{\cXX^{2\kk}}\approx 0, \quad  \rs \approx 0.
\]

\end{remark}

\subsubsection{Energy estimates in $E_{\kk+1, \xeta}$}\label{sec:EE_sub_norm_onto}

In light of Lemma \ref{lem:A2}, we already bound $\| \tw_2 \|_{\cXX^{2\kk+6}}$. 
In order to estimate $E_{\kk+1, \xeta} $, we perform
$\cXX^{2\kk+2}$ energy estimates on $\tw_1$ and 
$ \cYY^{2\kk+2}$ energy estimates on $\tFm$ using equations \eqref{eq:non_fix}
\bseq \label{eq:non_Y1}
\beq
\bal 
 \underbrace{ \kp \B\la  ( \pa_s - (\cL_{E,s} - \cK_k)  ) \tw_1, \tw_1    \B\ra_{\cZ^{2\kk+2} } }_{ := I_{\cL, M, 1} }
 & +  \underbrace{ \B\la (\pa_s - \Lmic) \tFm, \tFm \B\ra_{\cYY^{2\kk+2} } }_{ := I_{\cL, m, 1}} \\ 
  & = I_{\cL, M, 2} + I_{\cL, M, 3} + I_{\cL, m, 2} + I_{\cN, \xeta} + I_{\cE, \xeta} , 
  \label{eq:non_Y1:a} 
  \eal 
  \eeq
where $I_{\cL, M, \cdot}, I_{\cL,m, \cdot }$ denote macro and micro linear terms 
given by  
\beq \label{eq:non_Y1:b}
\bal
   I_{\cL, M, 2} +  I_{\cL, M, 3} &:= \kp \la  (- \cI_1, - \cI_2,  \cI_2 )(\tFm) , \  \tw_1 \ra_{\cZ^{2\kk+2}}  +  \kp \la  (\cL_{E,s} - \cL_E ) \tw_2, \tw_1 \ra_{\cZ^{2\kk+2}}  
, \\
I_{\cL, m,2} &:= %
 - \la  \cP _m [(V \cdot \na _X + 2 \dcm + \td d _\cM) \tFM] , \ \tFm  \ra_{\cYY^{2\kk+2}} , \\ 
 \eal 
\eeq
$I_{\cN, \xeta}$ is the nonlinear term 
\beq \label{eq:non_Y1:c}
 I_{\cN, \eta}  := \tf{1}{\e_s} \la \cN(\td F, \td F), \tFm \ra_{\cYe^{2\kk+2}} ,
\quad \eta = \xeta  \rm{\ or \ } \etab, 
\eeq
and $I_{\cE, \xeta}$ is the error term 
\beq \label{eq:non_Y1:d}
I_{\cE, \xeta} := - \kp \la   \tw_1,  ( \cs^3 \eu, \cs^3  \ep, 0) \ra_{\cZ^{2\kk+2}} 
- \la    \cP _m [\cMM ^{-1/2} \eM] , \tFm \ra_{\cYY^{2\kk+2}} :=  I_{\cE, M, \xeta} + I_{\cE, m, \xeta}
\eeq
\eseq

\paragraph{Estimates of linear terms}
Note that the weight in $\cX$-norm \eqref{norm:Xk} and $\cZ$-norm \eqref{norm:cZZ} are $s$-independent. Using the coercivity estimates in $\cZ$ norm \eqref{norm:Z} and  \eqref{eq:micro_lin_est:lin2} in Theorem \ref{thm: micro-Hk-main} with $ \eta = \xeta$, we estimate $I_{\cL, M, 1}, I_{\cL, m, 1}$ as 
\beq\label{eq:lin_M1}
\bal
& \f{1}{2} \f{d}{d s}  (\kp \| \tw_1 \|_{\cZZ^{2\kk+2}}^2 + \| \tFm \|_{\cYY^{2\kk+2}}^2)
+ \lam_1 \kp \| \tw_1 \|_{\cZ^{ 2\kk + 2}}^2
+  (2\lame - C \es) \| \tFm \|_{\cYY^{2\kk+2}}^2
+  \f{ \cgam } { 6 \es} \| \tFm \|_{\cYL{\xeta} ^{2\kk + 2}}^2  \\ 
  & \quad  \leq \kp \la  (\pa_s -  (\cL_{E, s} - \cK_\kk) ) \tw_1 , \tw_1 \ra_{\cZ^{2\kk+2}} 
  +  \la (\pa_s - \Lmic)  \tFm, \tFm \ra_{\cYY^{2\kk+2} }  = I_{\cL, M, 1} + I_{\cL, m, 1}. 
\eal
\eeq

Next, we estimate the interaction between the macro and micro parts in $I_{\cL, M, 2}$. 
Using the definition of $\cZ$ norm in \eqref{norm:Z} and estimate \eqref{eq:est_cross:b} with 
$ \eta = \xeta$, we estimate $I_{\cL, M, 2}$ 
\beq\label{eq:lin_M2}
\bal
  I_{\cL, M, 2} & =  \kp \big\la \tw_1, \ (- \cI_1, - \cI_2,  \cI_2 )(\tFm) \ra_{ \cX_{\xeta}^{2\kk+2}}
  +  \kp \varpi_{\kk + 1}^{\pr} \la \tw_1, \  (- \cI_1, - \cI_2, \cI_2)(\tFm) \ra_{ \cXX^{2\kk}}   \\
    & = \kp \la \tw_1,  \ (- \cI_1, - \cI_2,  \cI_2 )(\tFm)   \ra_{ \cX_{\xeta}^{2\kk+2}}
  + O(  \| \tFm\|_{\cYL{\xeta}^{2\kk+2} }  \| \tw_1 \|_{\cXX^{2\kk}} ). 
  \eal
\eeq
Recall the map \eqref{eq:macro_UPB} and $\tw = \tw_1 + \tw_2$. Using 
$\tFM = \cF_M(\tw_1) + \cF_M(\tw_2)$, Lemma \ref{lem:macro_UPB}, and estimate \eqref{eq:est_cross:b} with $ \eta = \xeta$, we estimate $I_{\cL, m, 2}$ as 
\[
\bal
  I_{ \cL, m, 2 }  & = 
  - \big\la  \cP _m [(V \cdot \na _X + 2 \dcm + \td d _\cM) \tFM(\tw_1)] , \  \tFm \big\ra_{\cYY^{2\kk+2}} \\ 
   & \qquad - \big\la  \cP _m [(V \cdot \na _X + 2 \dcm + \td d _\cM) \tFM(\tw_2)] , \ \tFm \big\ra_{\cYY^{2\kk+2}}  \\
& =    - \big\la  \cP _m [(V \cdot \na _X + 2 \dcm + \td d _\cM) \tFM(\tw_1)] , \  \tFm \big\ra_{\cYY^{2\kk+2}}
+ O(    \| \tFm\|_{\cYL{\xeta} ^{2\kk+2 } } 
  \| \tw_2 \|_{\cXX^{2\kk+4}} ). 
\eal 
\]
We estimate the main terms in $I_{\cL, M, 2}$ and $I_{\cL, m, 2}$ together using \eqref{eq:est_cross:a}
\eqref{eq:est_cross:aa} in Proposition \ref{prop:EE_cross} with $\eta = \xeta$ and combine the error terms using 
$ \|\tw_1 \|_{\cXX^{2\kk }} \les \|\tw_1 \|_{\cXX^{2\kk+2} }$:
\bseq\label{eq:lin_mM2}
\beq
  |I_{\cL, m, 2} + I_{\cL, M, 2}|
  \les   \| \tFm\|_{\cYL{\xeta} ^{2\kk + 2}  }  ( \| \tw_1 \|_{\cXX^{2\kk+2}} + \| \tw_2\|_{ \cXX^{2\kk+4}}  ) .
\eeq
Using Cauchy--Schwarz inequality and the energy $E_{\kk+1, \xeta}$ \eqref{eq:non_E}, 
we obtain 
\beq
    |I_{\cL, m, 2} + I_{\cL, M, 2}| 
    \les \es^{1/2} \left(  E_{\kk+1, \xeta} +  \| \tw_2 \|^2_{\cXX^{2\kk+4}}
+ \f{1}{\es} \| \tFm\|_{\cYL{\xeta}^{2\kk+2} }^2  \right) .
\eeq
\eseq

For $I_{\cL, M, 3}$, applying Proposition \ref{prop:error_diff}, the equivalence of norms in Lemma \ref{lem:norm_equiv_XZ}, Cauchy--Schwarz inequality, and using the energy $E_{\kk+1, \xeta}$ in \eqref{eq:non_E} and $R_s^{-r} \ll \es^2$ from Remark \ref{lem:para_bd}, we estimate
\beq\label{eq:lin_M3}
\bal
  |I_{\cL, M, 3}| & \les \| \tw_1 \|_{\cZ^{2\kk+2}}
  \| ( \cL_{E, s} - \cL_E ) \tw_2 \|_{\cXX^{2\kk+2}} \\ 
& \les R_s^{-r} E_{\kk+1, \xeta}^{1/2} \| \tw_2 \|_{\cXX^{2\kk+4}}  
\les \es^2 E_{\kk+1, \xeta}^{1/2} \| \tw_2 \|_{\cXX^{2\kk+4}}  .
\\
\eal
\eeq

\paragraph{Estimates of nonlinear terms}

Consider $\eta = \xeta$ or $\etab$. For the nonlinear terms $I_{\cN}$,  we use $\td F = \tFm + \tFM$ to decompose
\[
\bal
   \la \cN( \td F, \td F ), \tFm \ra_{\cYe^{2\kk+2}}
   & = \la  \cN(\td F, \tFm) , \tFm \ra_{\cYe^{2\kk+2}}
   + \la  \cN(\tFm, \tFM) , \tFm \ra_{\cYe^{2\kk+2}}
   + \la  \cN(\tFM, \tFM) , \tFm \ra_{\cYe^{2\kk+2}} \\
   & := I_{\cN, m}
   +I_{\cN,m M} + I_{\cN, MM}.
   \eal
\]
Applying \eqref{eq:non_Q:micro} 
in Theorem \ref{theo: N nonlin_est}
to $I_{\cN, m}$, \eqref{eq:non_Q:macro:mM} to $I_{\cN, m M}$,
and \eqref{eq:non_Q:macro:MM2} with $\eta \in \{ \xeta, \etab \}$ to $I_{\cN, MM}$, and then using \eqref{eq:macro_Y} in Lemma  \ref{lem:macro_UPB} to bound $\tFM = \cF_M(\tw)$, we obtain
\bseq\label{eq:est_non}
\beq
\bal
   |I_{\cN,m} | & \les \| \td F \|_{\cYb^{2\kk+2}} \| \tFm \|_{\cYE ^{2\kk+2 }} 
   \|   \tFm  \|_{\cYE ^{2\kk+2} } 
  \les ( \| \tFm \|_{\cYb^{2\kk+2}} + \| \tw \|_{\cXb^{2\kk+2}} )  
   \|   \tFm  \|_{\cYE^{ 2\kk+2 } }^2,  \\
  |I_{\cN, mM}| & \les 
  \| \tFM \|_{\cYb^{2\kk+2}} \| \tFm \|_{\cYE^{ 2\kk+2 }} 
   \|   \tFm  \|_{\cYE^{2\kk+2} } 
  \les \| \tw \|_{\cXb^{2\kk+2}}  \|   \tFm  \|_{\cYE^{ 2\kk+2  } }^2 ,  \\
    |I_{\cN, MM}| & 
  \les \| \tFM \|_{\cYY^{2\kk+2}  }
    \| \tFM \|_{\cYY^{2\kk - 2} } 
        \|   \tFm  \|_{\cYE^{ 2\kk+2  } } .
   \eal
\eeq
For $\eta = \etab$ or $\xeta$,
combining the above estimates and using the energy $E_{ \kk +1, \etab }$ \eqref{eq:non_E}, 
we get 
\begin{align}\label{eq:est_non:b}
     |I_{\cN, \eta}|=   \f{1}{\es}   |\la \cN( \td F, \td F ), \tFm \ra_{\cYe^{ 2\kk+2  }}|
    & \les  \f{1}{\es} \B( E_{ \kk +1, \etab }^{1/2} \| \tFm \|_{\cYE^{ 2\kk+2  }}^2 
+ \| \tFM \|_{\cYY^{2\kk+2}  }
    \| \tFM \|_{\cYY^{ 2\kk - 2 } } 
        \|   \tFm  \|_{\cYE^{ 2\kk+2  } } \B). 
\end{align}
\eseq

\paragraph{Estimate of error terms}

Using \eqref{eq:error_L2} in Lemma \ref{lem: cutoff_error} and Cauchy--Schwarz inequality, for $
\eta = \xeta$ or $\etab$, $ n = 2\kk , 2\kk+2$, 
and any function $G \in \cX^{n}_{\eta}$, we have 
\[
    |\la   G,  ( \cs^3 \eu, \cs^3  \ep, 0) \ra_{\cX_{\eta}^n} |
    \les \| G \|_{\cX_{\eta}^n} 
     \| ( \cs^3 \eu, \cs^3  \ep, 0) \|_{\cX_{\eta}^{n}}  
  \les \| G \|_{\cX_{\eta}^n}  \rs^{-r}. 
\]

Recall the definition $\cZ$ norm \eqref{norm:Z}. 
Using the energy \eqref{eq:non_E}, we estimate $I_{\cE, M , \cdot }$ \eqref{eq:non_Y1:d}
\bseq\label{eq:error_onto}
\beq\label{eq:error_onto:a}
\bal
  |\la   \tw_1,  ( \cs^3 \eu, \cs^3  \ep, 0) \ra_{\cZ^{2\kk+2}}| & \les \rs^{-r} E_{\kk+1, \xeta}^{1/2} ,
  \\
  |\la   \tw,  ( \cs^3 \eu, \cs^3  \ep, 0) \ra_{\cXb^{2\kk+2}}| & \les \rs^{-r} 
  E_{\kk+1 , \etab}^{1/2} .
  \eal 
\eeq
Note that in the second estimate, we use $\tw$ (the whole macro-perturbation)  instead of $\tw_1$.

For the other error term,  for any $\upsilon>0$, using \eqref{eq:micro_lin_est:err} in  Theorem 
\ref{thm: micro-Hk-main} with $\eta = \xeta$ or $\etab$, we have
\beq
    | \la    \cP _m [\cMM ^{-1/2} \eM] , \tFm \ra_{\cYe^{ 2\kk+2  }}| 
    \les   \| \tFm \|_{\cYE^{2\kk + 2} }  
    \les  \f{\ups}{\es}  
   \| \tFm \|_{\cYE^{ 2\kk+2  }}^2  + \ups^{-1} \es
  \eeq
  \eseq

\paragraph{Consequences of the bootstrap assumptions}

We treat all the terms except the coercive terms in \eqref{eq:lin_M1} perturbatively. Under the bootstrap assumptions \eqref{eq:EE_boot}, using the bounds \eqref{eq:est_W12}, 
\eqref{eq:EE_para_k} and $\rs^{-r} \ll \es^2$ from Remark \ref{lem:para_bd},  we simplify the estimates \eqref{eq:lin_mM2}, \eqref{eq:lin_M3}, \eqref{eq:error_onto} with $\eta = \xeta$
\beq\label{eq:non_Y1_boot}
\bal
    |I_{\cL, m, 2} + I_{\cL, M, 2}| 
   & \les \es^{1/2} (  E_{\kk+1, \xeta} +  \| \tw_2 \|_{\cXX^{2\kk+2}}^2
+ \f{1}{\es} \| \tFm\|_{\cYL{\xeta}^{2\kk+2} }^2  )  \\
& \les \es^{1/2} ( \es^{ 1 - 2 \el } + \f{1}{\es} \| \tFm\|_{\cYL{\xeta}^{2 \kk+2} }^2 ) 
\les \es + \es^{1/2} \cdot \f{1}{\es}\| \tFm\|_{\cYL{\xeta}^{2\kk+2} }^2  , \\
|I_{\cL,M,3}| & \les
\es^2 E_{\kk+1, \xeta}^{1/2} \| \tw_2 \|_{\cXX^{ 2\kk+4 }}  \les  \es^{2} \cdot \es^{ 1/2 - \el} \es^{2/3-\ell} \les \es^2 ,  \\
| I_{\cE, \xeta}| & \leq | I_{ \cE, M, \xeta}|  + |I_{\cE, m , \xeta}|
\leq  C \es^2  + \f{\ups}{\es} \| \tFm\|^2_{\cYL{\xeta}^{2\kk+2} } + \ups^{-1} \es.
\eal
\eeq

For $\eta = \etab$ or $\xeta$, since $\f{1}{6} - \el > \el$ \eqref{eq:EE_para_k} and $\es \les \d$ from \eqref{eq:EE_para1}, we simplify the estimate \eqref{eq:est_non:b} as
\beq\label{eq:non_boot}
\bal
|I_{\cN, \eta} | & \les  \d^{ \el} \cdot  \f{1}{\es}  \| \tFm \|_{\cYE^{2\kk+2}}^2 
+ \es^{ \f{1}{2}-\el + \f{2}{3} - 1}
        \|   \tFm  \|_{\cYE^{2\kk+2 } }  \\
    & \les \es +  ( \d^{\el} + \es^{2 ( \f{1}{6} -\el )   } ) \cdot  \f{1}{\es}  \| \tFm \|_{\cYE^{2\kk+2}}^2
    \les \es +   \d^{\el}  \cdot  \f{1}{\es}  \| \tFm \|_{\cYE^{2\kk+2}}^2
\eal
\eeq

\paragraph{Summary of the estimates}
Recall 
\[
E_{\kk+1, \xeta} = \kp \| \tw_1 \|_{\cZZ^{2\kk+2}}^2 + \| \tFm \|_{\cYY^{2\kk+2}}^2.
\] 
Applying the estimates \eqref{eq:lin_M1}, \eqref{eq:non_Y1_boot}, \eqref{eq:non_boot} with $\eta = \xeta$
to \eqref{eq:non_Y1}, using $\lame > \lam_1$ \eqref{eq:decay_para}, $\es \leq \d$ by \eqref{eq:EE_para1}, and choosing $\d$ small enough, we derive
\[
\bal
 \tf{1}{2} \tf{d}{d s} E_{\kk+1, \xeta}  \leq &  -\lam_1 E_{\kk+1, \xeta} 
+ 
 \f{1}{\es} ( 
 - \f{\cgam}{6} + C \es^{1/ 2} + \ups + C\d^{ \el}  ) \| \tFm \|_{\cYL{\xeta}^{2\kk+2} }^2  +   C (1 + \ups^{-1})\es  .
\eal 
\]
Choosing $\ups = \f{1}{100} \cgam$ and $\d >0$ small enough (depending on $\kk, \xeta$)  so that $\es \leq \d$ is very small (by \eqref{eq:EE_para1}), we establish 
\bseq\label{eq:EE_Y1}
\beq\label{eq:EE_Y1:a}
 \f{1}{2} \f{d}{d s} E_{\kk+1, \xeta}   \leq - \lam_1  E_{\kk+1, \xeta} 
 -  \f{ \cgam  }{ 8 \es} \| \tFm \|_{\cYL{\xeta}^{2\kk+2} }^2
 + C \es \leq  - \lam_1  E_{\kk+1, \xeta}  + C \es .
\eeq
Recall $\es = \d \operatorname e^{-\rE s  }$ from \eqref{eq:EE_para1}. Since $\lam_1 > \f{1}{2}\rE$ from \eqref{eq:decay_para}, $\el > 0$ from \eqref{eq:EE_para_k}, and $E_{\kk+1, \xeta}(0) < \d$ from \eqref{eq:IC:small}, solving the inequality and choosing $\d $ small enough, we prove 
\beq
\bal
  E_{\kk+1, \xeta} & \leq \operatorname e^{-2 \lam_1 s } E_{\kk+1, \xeta}(0) + C \d \int_{0}^s  \operatorname e^{-2\lam_1(s - \tau )} e^{-\rE \tau} d \tau 
   \leq C \d \operatorname e^{-\rE s  }= C \es \ll  \f{1}{2} \es^{ 1-2\el }. 
\eal
\eeq
\eseq
We improve the bound \eqref{eq:EE_boot:a} and proved the second estimate in \eqref{eq:EE_onto:a} for $s \in [0, \bar s]$, where the bootstrap assumptions in \eqref{eq:EE_boot} hold.

\subsubsection{Energy estimates in  $\cXX^{2\kk}$}\label{sec:sub_norm_Xk}

The energy estimates on $ \| \tw_1 \|_{\cXX^{2\kk}}$ is similar and is easier. 
We perform $\cXX^{2\kk}$ estimates on $\tw_1$ in \eqref{eq:non_fix} and use similar decompositions as in \eqref{eq:non_Y1} %
\beq 
\bal 
  \la  ( \pa_s - (\cL_{E,s} - \cK_k)  ) \tw_1, \tw_1    \ra_{\cXX^{2\kk} } %
  & =  \la  (- \cI_1, - \cI_2,  \cI_2 )(\tFm) , \  \tw_1 \ra_{\cXX^{2\kk} }  \\ 
 & \quad +   \la  (\cL_{E,s} - \cL_E ) \tw_2, \tw_1 \ra_{\cXX^{2\kk} }
 -  \la   \tw_1,  ( \cs^3 \eu, \cs^3  \ep, 0) \ra_{\cXX^{2\kk} } \\ 
  &  := I_{\cL, M, 2} + I_{\cL, M, 3} + I_{\cE, M} . 
   \eal 
  \eeq
The estimates of the left hand side, $I_{\cL, M, 3}, I_{\cE_M}$ are
 similar to 
\eqref{eq:lin_M1}, \eqref{eq:lin_M2}, \eqref{eq:lin_M3}, \eqref{eq:error_onto:a} (replacing the norm $\cZ^{2\kk+2}$ by $\cXX^{2\kk}$ and $\kp$ by $1$). Thus, we only state the estimates and use the bootstrap assumption \eqref{eq:EE_boot} to further simplify them
\bseq\label{eq:EE_Wk0}
\beq
\bal
\f{1}{2} \f{d}{ds} \| \tw_1 \|_{\cXX^{2\kk}}^2 
+ \lame \| \tw_1 \|_{\cXX^{2\kk}}^2  
& \leq \la  (\pa_s - (\cL_{E, s} - \cK_\kk)  ) \tw_1 , \tw_1 \ra_{\cXX^{2\kk}}  , \\ 
|I_{\cL, M, 3} | & =  | \la  (\cL_{E,s} - \cL_E ) \tw_2, \tw_1 \ra_{\cXX^{2\kk}} |
\les  \rs^{-r}  \| \tw_1 \|_{\cXX^{2\kk}} \| \tw_2 \|_{\cXX^{2\kk+2}}
\les \es^2 ,  \\
|I_{\cE, M}| & = 
  |\la   \tw_1,  ( \cs^3 \eu, \cs^3  \ep, 0) \ra_{\cXX^{2\kk}}| 
  \les \| \tw_1 \|_{\cXX^{2\kk}} \rs^{-r} \les \es^2.
  \eal
\eeq
For the cross term $I_{\cL, M, 2}$, we simply bound it using \eqref{eq:est_cross:b} with $\eta = \xeta$
\beq
  |I_{\cL, M, 2}| = 
   |  \la \tw_1, (- \cI_1, - \cI_2,  \cI_2 )(\tFm) \ra_{\cXX^{2\kk}} |
  \les     \| \tFm\|_{\cYL{\xeta}^{ 2\kk+2}}  
  \| \tw_1 \|_{\cXX^{2\kk}} .
\eeq
\eseq
Thus, combining the above estimates, we prove
\bseq\label{eq:EE_Wk}
\beq\label{eq:EE_Wk:a}
 \f{1}{2}  \f{d}{d s}  \| \tw_1 \|_{\cXX^{2\kk}}^2
 \leq - \lame   \| \tw_1 \|_{\cXX^{2\kk}}^2
 + \es^2  +
 C    \es^{1/2} \cdot   \f{1}{\es^{1/2}} \| \tFm\|_{\cYL{\xeta}^{2\kk+2}} 
  \| \tw_1 \|_{\cXX^{2\kk}} .
\eeq

The small factor $\es^{1/2}$ in the above estimates indicates that the estimates of $\| \tw_1 \|_{\cXX^{2\kk}}$ and $E_{\kk+1, \xeta}$ in \eqref{eq:EE_Y1} are weakly coupled. 
Recall $\es = \d \operatorname  e^{-\rE s }$ \eqref{eq:EE_para1}.  Next, we estimate  
\beq\label{eq:non_Emix}
E_{\kk,\rm{mix}} =  E_{\kk+1, \xeta} + \es^{-2( \f{1}{6} + \el )}  \| \tw_1 \|_{\cXX^{2\kk}}^2 .
\eeq
where $ \es^{-2( \f{1}{6} + \el )}$ is the difference between decay rates of $E_{\kk+1, \xeta}$ and 
$\| \tw_1 \|_{\cXX^{2\kk}}^2$ in \eqref{eq:EE_boot}.
We estimate $\es^{-2( \f{1}{6} + \el )}  \| \tw_1 \|_{\cXX^{2\kk}}^2$ 
by multiplying  \eqref{eq:EE_Wk:a} by $\es^{-2( \f{1}{6} + \el )}$ 
and using $\tf{1}{2}\f{d}{ds} \es^{-b} = -\tf{1}{2} b \rE \es^{-b}$ 
\beq
\bal
 \f{1}{2}  \f{d}{d s} \left( \es^{-2( \f{1}{6} + \el )}  \| \tw_1 \|_{\cXX^{2\kk}}^2 \right)
& \leq  - \left( \lame -  \left( \f{1}{6} + \el \right) \rE \right)  \es^{-2( \f{1}{6} + \el )}  \| \tw_1 \|_{\cXX^{2\kk}}^2
+ C \es^{2 - 2( \f{1}{6} + \el  )} \\
& \qquad + C \es^{ \f12 - ( \f{1}{6} + \el )}  
\cdot   \f{1}{\es^{1/2}} \| \tFm\|_{\cYL{\xeta}^{2\kk+2}} 
  \cdot \es^{ - (\f{1}{6} + \el)} \| \tw_1 \|_{\cXX^{2\kk}} .
  \eal
\eeq
\eseq

Combining \eqref{eq:EE_Y1:a} and \eqref{eq:EE_Wk} and using Cauchy--Schwarz inequality, we derive
\beq\label{eq:EE_Hk_W1}
\bal
\f{1}{2} \f{d}{ds} E_{\kk, \rm{mix}}
& \leq - \min \left\{ \lame - \left( \f{1}{6} + \el  \right) \rE, \lam_1 \right\} E_{\kk,  \rm{mix}}- \f{ \cgam }{ 8 \es} \| \tFm\|_{\cYL{\xeta}^{2\kk+2}}^2   \\
& \quad + C \es^{ \f{1}{2}- ( \f{1}{6} + \el )} \left( E_{\kk, \rm{mix}} + \f{1}{\es} \| \tFm\|_{\cYL{\xeta}^{ 2\kk+2 }}^2 \right) + C ( \es + \es^{2 - 2( \f{1}{6} + \el )} ).
\eal
\eeq
Since $\el = 10^{-4}$ \eqref{eq:EE_para_k}, from \eqref{eq:IC:small} and \eqref{eq:decay_para}, we obtain 
\beq\label{eq:EE_Hk_W1_init}
E_{\kk, \mw{mix}}( 0) \les \d + \d^{-2 (1/6 + \el)} 
\d^{4/3 + 2\el} \les \d 
\eeq
From \eqref{eq:decay_para}, we have 
\beq\label{eq:EE_Hk_W1_para}
  \lame - ( \f16 + \ell ) \rE > ( \f{2}{3} - \ell ) \rE > \f{7}{12} \rE, \quad \lam_1 > \f{7}{12} \rE , \quad \f{1}{2} - ( \f{1}{6} + \el ) > \f{1}{60},
  \quad 2 - 2(\f{1}{6} + \el) > 1.
\eeq
Using $\es \leq \d$ by \eqref{eq:EE_para1}, choosing $\d$ small enough,  and using \eqref{eq:EE_Hk_W1} and \eqref{eq:EE_Hk_W1_para}, we derive 
\[
\bal
\f{1}{2} \f{d}{ds} E_{\kk, \rm{mix}}
& \leq -  (\f{7}{12} \rE - C \es^{\f{1}{60}} ) E_{\kk,  \rm{mix}}- \f{ (\cgam - C \es^{ \f{1}{60} }) }{ 8 \es} \| \tFm\|^2_{\cYL{\xeta}^{2\kk+2}}  + C \es \\
& \leq -  \f{13}{24} \rE \cdot E_{\kk,  \rm{mix}}- \f{ \cgam }{ 9 \es} \| \tFm\|^2_{\cYL{\xeta}^{2\kk+2}}  + C \es .
\eal
\]
Since $\f{13}{24} \rE > \f{1}{2} \rE$, solving the above inequality similar to \eqref{eq:EE_Y1} and using \eqref{eq:EE_Hk_W1_init}, we prove 
\[
  E_{\kk, \rm{mix}} \leq C \d \operatorname  e^{-\rE s  } = C \es , 
\]
which along with \eqref{eq:non_Emix} implies 
\beq\label{eq:EE_W1k}
  \| \tw_1 \|_{\cXX^{2\kk}} 
\leq C \es^{ \f{1}{2} + \f{1}{6} + \el }  \ll \es^{ 2/3 }.
\eeq
We improve the estimate \eqref{eq:EE_boot:b} %
and prove the first bound in  \eqref{eq:EE_onto:a}
for $s \in [0, \bar s]$, where the bootstrap assumptions in \eqref{eq:EE_boot} hold.

\subsubsection{Energy estimates in $\cXb^{2\kk+2}$ and $\cYb^{2\kk+2}$}\label{sec:EE_crit_norm_onto}
 To control $E_{\kk+1, \etab}$, we estimate $\tw$ and $\tFm$. 
We recall the equation of $\tw$ from  \eqref{eq:non_fix:W}
\[
  \pa_s \tw = \cL_{E, s} \tw  + \cK_\kk (\hw_1 - \tw_1)  - (\cI_1 , \cI_2,  -  \cI_2 )(\tFm)
   - ( \cs^3 \eu, \cs^3  \ep, 0) .
\]

The energy estimates on $E_{\kk+1, \etab}$ is similar to those of $E_{\kk+1, \xeta}$ in Section \ref{sec:EE_sub_norm_onto}. We have 
\bseq
\beq
\bal 
 \underbrace{ \kp \B\la  ( \pa_s - \cL_{E,s} )  ) \tw , \tw     \B\ra_{\cXb^{2\kk+2} } }_{ := I_{\cL, M, 1} }
 & +  \underbrace{ \B\la (\pa_s - \Lmic) \tFm, \tFm \B\ra_{\cYb^{2\kk+2} } }_{ := I_{\cL, m, 1}} \\ 
  & = I_{\cL, M, 2} +  I_{\cL, M, 4} +  I_{\cL, m, 2} + I_{\cN, \etab} + I_{\cE, \etab} , 
  \label{eq:non_Y2:a} 
  \eal 
  \eeq
where $I_{\cL, M, \cdot}$ denote the macro linear terms given by
\beq
\bal
 I_{\cL, M, 2} +  I_{\cL, M, 4}  & :=      \kp \big\la  (- \cI_1 , - \cI_2 ,\cI_2 )(\tFm) , \ \tw \big\ra_{\cXb^{2\kk+2}}  
 + \kp \la \cK_\kk (\hw_1 - \tw_1) , \tw \ra_{\cXb^{2\kk+2}} , 
\eal
\eeq
and we decompose $I_{\cL, m}, I_{\cN, \etab}$ and $I_{\cE}$ in the same way as those in \eqref{eq:non_Y1} 
\begin{align}
I_{\cL, m, 2}  & = 
 -  \big\la   \cP _m [(V \cdot \na _X + 2 \dcm + \td d _\cM) \tFM] , \ \tFm  \big\ra_{\cYb^{2\kk+2}},  \label{eq:non_Y2:b} \\
 I_{\cN, \etab} & := \f{1}{\e_s} \la \cN(\td F, \td F), \tFm \ra_{\cYb^{2\kk+2}} ,  \label{eq:non_Y2:c} \\
I_{\cE, \etab} & := - \kp \la   \tw_1,  ( \cs^3 \eu, \cs^3  \ep, 0) \ra_{\cXb^{2\kk+2}} 
- \la    \cP _m [\cMM ^{-1/2} \eM], \tFm \ra_{\cYb^{2\kk+2}} :=  I_{\cE, M,\etab} + I_{\cE, m,\etab} ,
\label{eq:non_Y2:d}
\end{align}
\eseq

Using Theorem \ref{thm:coer_est} and Theorem \ref{thm: micro-Hk-main} with $ \eta = \etab$, we estimate $I_{\cL, M, 1}, I_{\cL, m, 1}$ as 
 \[
 \bal
& \f{1}{2} \f{d}{ds} ( \kp \| \tw \|_{\cXb^{2\kk+2}}^2 + \| \tFm \|_{\cYb^{2\kk+2}}^2 ) 
   - C \int  | D_X^{\leq 2\kk+2} \tw|^2 \la X \ra^{\etab- r} d X 
   - C \es  \| \tFm \|_{\cYb^{2\kk+2}}^2
 + \f{ \cgam } { 6 \es} \| \tFm \|_{\cYQ^{2\kk+2}}^2 \\ 
 &\qquad   \leq  \kp \la ( \pa_s -  \cL_{E, s}  \tw  )  , \tw  \ra_{\cXb^{2\kk+2}} 
  +  \la  (\pa_s - \Lmic) \tFm , \tFm \ra_{ \cYb^{2\kk+2} } 
  =  I_{\cL, M , 1}  + I_{\cL, m, 1} .
 \eal 
 \]
 The $\f{d}{ds}$-term  gives exactly $\f{d}{ds} E_{\kk+1, \etab}$ \eqref{eq:non_E}. 
 Note that on the left hand side, we have the $\cYb$-norm term $ - C \es \| \tFm \|_{\cYb^{2\kk+2}}^2$ rather than 
 $ - C \es \| \tFm \|_{\cYY^{2\kk+2}}^2 $. Since $\etab - r < \xeta$ by \eqref{eq:eta_constraint}, using 
 \[
 \int  | D_X^{\leq 2\kk+2} \tw|^2 \la X \ra^{\etab- r} d X \les \| \tw \|_{\cXX^{2\kk+2}}^2,
 \]
and the energy $E_{\kk+1, \xeta}, E_{\kk+1, \etab}$ \eqref{eq:non_E}, we obtain 
\[
 \f{1}{2} \f{d}{ds} E_{\kk+1, \etab}  - C E_{\kk+1, \xeta} - C \es E_{k+1, \etab} + \f{  \cgam} { 6 \es} \| \tFm \|_{\cYQ^{2\kk+2}}^2 
   \leq 
    I_{\cL, M , 1}  + I_{\cL, m, 1} .
\]

For $I_{\cL, M, 4}$, using $\cK_{\kk} = \cK_{\kk, \xeta}$ by \eqref{def:kk} and 
$\supp( \cK_{\kk} f  ) \subset B(0, 4 R_{\xeta})$ by item (a) in Proposition \ref{prop:compact}, we obtain 
\[
  |I_{\cL, M, 4} | \les \| \cK_\kk ( \tw_1 - \hw_1) \|_{\cXb^{2\kk+2}}
  \| \tw \|_{\cXb^{2\kk+2}} 
\les \| \cK_\kk ( \tw_1 - \hw_1) \|_{\cXX^{2\kk+2}}
  \| \tw \|_{\cXb^{2\kk+2}} .
\]
Using energy \eqref{eq:non_E}, item (c) in Proposition \ref{prop:compact}, 
and bounds \eqref{eq:EE_boot:b}, \eqref{eq:est_W12}, implied by the bootstrap assumptions, we obtain 
\[
    |I_{\cL, M, 4} | \les  \|  \tw_1 - \hw_1 \|_{\cXX^{2\kk}}    \| \tw \|_{\cXb^{2\kk+2}} 
\les E_{\kk+ 1, \etab}^{1/2}  
 \|  \tw_1 - \hw_1 \|_{\cXX^{2\kk}} .
\]

For $I_{\cL, M, 2 } + I_{\cL, m, 2}$, applying Proposition \ref{prop:EE_cross} with $\eta = \etab$, the energy \eqref{eq:non_E}, we obtain 
\[
  | I_{\cL, M, 2 } + I_{\cL, m, 2}|
  \les \es^{1/2} \| \tw \|_{\cXb^{2\kk+2}} 
  \cdot \f{1}{\es^{1/2}} \| \tFm\|_{\cYb^{ 2\kk+2} }  
\les \es^{1/2} ( E_{ \kk +1, \etab} +  \f{1}{\es} \| \tFm\|_{\cYb^{ 2 \kk+ 2} }^2 ) .
\]

We have estimated $I_{\cN, \etab}$ in \eqref{eq:est_non} and $I_{\cE, \etab}$ in \eqref{eq:error_onto} with $\eta = \etab$.

\paragraph{Consequences of the bootstrap assumptions}

Under the bootstrap assumptions \eqref{eq:EE_boot}, using the bounds \eqref{eq:est_W12}, 
\eqref{eq:EE_para_k} and $\rs^{-r} \ll \es^2$ from Remark \ref{lem:para_bd}, we simplify the above estimates as
\bseq\label{eq:non_Y2_boot}
\beq
\bal
 \f{1}{2} \f{d}{ds} E_{\kk+1, \etab}
 - C \es^{ 1 - 2 \el }  - C \es + \f{  \cgam} { 6 \es} \| \tFm \|_{\cYQ^{2\kk+2}}^2 
    &\leq 
    I_{\cL, M , 1}  + I_{\cL, m, 1} , \\ 
    |I_{\cL, M, 4} | & \les  \es^{2/3 } \d^{ \el }, \\
      | I_{\cL, M, 2 } + I_{\cL, m, 2}| & \les \es^{1/2} \d^{ 2 \el }
      + \es^{1/2} \cdot  \f{1}{\es} \| \tFm\|_{\cYb^{2\kk + 2} }^2 .
\eal
\eeq
Applying the estimates of $I_{\cN, \etab}$ in \eqref{eq:est_non}, \eqref{eq:non_boot} %
and estimates of $I_{\cE, \etab}$ in \eqref{eq:error_onto} with $\eta = \etab$, and estimates
$\rs^{-r} E_{\kk+1, \etab}^{1/2} \les \es^2$ from \eqref{eq:EE_boot:c},
\eqref{eq:para_bd},  we obtain
\beq
\bal
|I_{\cN, \etab}| & \les \es +   \d^{\el}  \cdot  \f{1}{\es}  \| \tFm \|_{\cYb^{2\kk+2}}^2 , \\
 I_{\cE,\etab} & \leq | I_{\cE, M,\etab}|  + |I_{\cE, m,\etab}|
 \leq  C \es^2  + \f{\ups}{\es} \| \tFm\|_{\cYQ^{2\kk+2} }^2 + \ups^{-1} \es .
\eal
\eeq
\eseq

\paragraph{Summary of the estimates}

Combining the estimates in \eqref{eq:non_Y2_boot},  we derive
\[
\bal
 \f{1}{2} \f{d}{d s} E_{\kk+1, \etab}  \leq & 
 \f{1}{\es} (- \f{\cgam}{ 6 } + C \es^{1/ 2} + C \d^{ \el } 
 + \ups
 ) \| \tFm \|_{\cYQ^{2\kk+2} }^2 \\
 & +   C (1 + \ups^{-1}) \es + C ( \es^{ 1 - 2 \el } +  \es^{2/3} \d^{ \el }
+ \es^{1/2} \d^{ 2 \el }
  ).
\eal 
\]
Recall the bounds of $\el$ from \eqref{eq:EE_para_k} and $\es = \d \operatorname e^{-\rE s  }$ from \eqref{eq:EE_para1}. We have 
\beq\label{eq:non_Y2_es}
  \max \{ \es,  \   \es^{ 1 - 2\el } +  \es^{2/3} \d^{ \el }
+ \es^{1/2} \d^{ 2 \el }
  \} \les \es^{1/2}.
\eeq
By choosing $\ups = \f{\cgam}{100 }$, then choosing $\d$ small enough (depending on $\kk, \xeta, \etab$),
and using \eqref{eq:non_Y2_es}, we obtain 
\[
   \f{1}{2} \f{d}{d s} E_{\kk+1, \etab}  
   \leq - \f{ \cgam  }{8 \es}  \| \tFm \|_{\cYQ^{2\kk+2} }^2
   + C \es^{1/2} \leq C \es^{1/2}.
\]
Integrating the above estimate in $s$, using $\es = \d \operatorname e^{-\rE  s }$, and  \eqref{eq:est_init}
\footnote{
Recall that we assume that the bootstrap assumptions \eqref{eq:EE_boot} hold for $s \in [0, \bar s] $. 
}
 , we obtain
\beq\label{eq:EE_Wk_etab}
  \f{1}{2} E_{\kk+1, \etab}(s) +  \f{ \cgam  }{8}  \int_{0}^{s} \f{1}{\es} 
   \| \tFm(\tau) \|_{\cYQ^{2\kk+2} }^2 d \tau  
  \leq   \f{1}{2} E_{\kk+1, \etab}( 0 ) + C \d^{ 1/ 2 } \les \d^{4\el} \ll  \d^{2 \el}.
\eeq
Thus, we have improved the estimate \eqref{eq:EE_boot:c} and proved the third estimate in \eqref{eq:EE_onto:a} for any $s \leq \bar s$. The above estimate also implies \eqref{eq:EE_onto:b} for $s \leq \bar s$.

Combining \eqref{eq:EE_Y1}, \eqref{eq:EE_W1k}, \eqref{eq:EE_Wk_etab}, we improve all bootstrap assumptions in \eqref{eq:EE_boot} for $s \in [0, \bar s]$.  Therefore, the bootstrap assumptions hold for $s \in [0, \bar s)$ with $\bar s = \infty$. Combining estimates \eqref{eq:est_W12:a}, \eqref{eq:est_init}, \eqref{eq:EE_Y1}, \eqref{eq:EE_W1k}, and \eqref{eq:EE_Wk_etab}, we prove Proposition \ref{prop:onto}.

\subsection{Proof of Proposition \ref{prop:contra}}\label{sec:contra}

As in the assumption of the proposition, 
let $\wh W_{1, \al} \in Y, \al \in \{a, b\} $ be such that $\hat E_{\al} = \| \wh W_{1, \al} \|_{Y} < 1$. According to \eqref{eq:non_fix:W2}, \eqref{eq:non_fix},\eqref{eq:map_A}, denote the associated solutions
\[
\bal  
\tw_{2, \al} = \cA_2(\hw_{1, \al}),
\quad  ( \tw_{1,\al}, \td{F}_{\al, m } ) = \cA( \hw_{1, \al} ),
\quad \tw_{\al} = \tw_{1,\al} + \tw_{2, \al},
\quad \td{F}_{\al, M} = \cF_M( \tw_{\al} ),
\eal
\]
for $\al \in \{a, b\}$. Throughout this proof, we use the subscript $\al \in \{a, b\}$ to denote two different solutions, and we adopt the notation introduced \eqref{eq:non_E}; e.g.~$E_{\kk+1, \eta}, \eta = \xeta, \etab$ for the ``energies'' of these two solutions. From Propositions \ref{prop:onto} 
and \eqref{eq:est_W12} and estimate \eqref{eq:macro_Y}, we obtain 
\bseq\label{eq:non_contra1}
\begin{gather}
  \| \tw_{2, \al} \|_{\cXX^{2\kk+6} }
  + \| \tw_{1, \al} \|_{\cXX^{2\kk}}  \les \es^{2/3-\ell},  
\quad \| \tw_{1,\al} \|_{\cXX^{2\kk+2} } + \| \td F_{\al, m} \|_{ \cYY^{ 2\kk+2 } } \les \es^{1/2 - \ell}, \label{eq:non_contral:a} \\
\| \td F_{\al, M} \|_{\cXX^{2\kk}} \les \| \tw_{\al} \|_{\cXX^{2\kk}}
\les \es^{2/3-\ell},
\quad \  \| \td F_{\al, M} \|_{\cXX^{2\kk+2}} \les  \| \tw_{\al} \|_{\cXX^{2\kk+2}} \les  \es^{1/2-\ell} , 
\label{eq:non_contral:FMb} \\
    \| \tw_{1, \al} \|_{\cXb^{2\kk+2}} < \d^{ \el},
  \quad  
  \| \td F_{\al, m} \|_{\cYb^{2\kk+2}}  \les \d^{ \el} ,    \\ 
  \int_{0}^{\infty} \f{1}{\es} \| \td F_{\al, m}(s) \|^2_{\cYQ^{2\kk+2}} d s \les \d^{2\el}, \quad \al \in \{ a, b \}.
   \label{eq:non_contra1:b}
\end{gather}
\eseq
Additionally, we denote the difference of two solutions by a $\Delta$-sub-index:
\begin{subequations}
\label{eq:nota_del}
\begin{gather}
 \hw_{1,\D} = \hw_{1,a} - \hw_{1, b}, \quad  
\tw_{i, \D} = \tw_{i, a} - \tw_{i, b},  \ i=1,2,  \quad  
\tw_{\D} = \tw_{a} - \tw_b ,\\
 \td F_{\D,  m } = \td{F}_{a, m} - \td{F}_{b, m}, 
\quad 
\td F_{\D, M} = \td F_{a, M} - \td F_{b, M} ,
\quad \td F_{\D} = \td F_{a} - \td F_{b}, \\
   \cN_{ \D} = \cN( \td F_a ,\td F_a) - \cN( \td F_b , \td F_b) , 
\end{gather}
and introduce the following energies for the difference
\beq \label{eq:nota_del:d}
   E_{\kk+1, \D}(s) :=  \kp \| \tw_{1,\D}  \|_{\cZ^{2\kk + 2}}^2 
 + \| \td F_{m, \D} \|_{\cYY^{2\kk + 2}}^2,
 \quad  \cE_{\D} = \| \tw_{2, \D} \|_{\cXX^{2\kk+6}}.
\eeq 
\end{subequations}
With this notation, to prove Proposition~\ref{prop:contra}, we will show  
\begin{equation}
\label{eq:contra_goal}
  \| \tw_{1,\D}(s) \|_{\cXX^{2\kk}}
< \f{1}{2 } \es^{2/3} \| \hw_{1, \D} \|_Y .
\end{equation}

Using  \eqref{eq:non_fix}, we deduce that $ \tw_{1,\D}, \td F_{m, \D}$ solves 
\beq\label{eq:F_diff}
\bal
  \pa_s \tw_{1,\D} & =  (\cL_{E, s} - \cK_\kk ) \tw_{1,\D}   + (\cL_{E,s} - \cL_E) \tw_{2,\D}  - (\cI_1, \cI_2,  -  \cI_2)(\td F_{m,\D} ) ,
  \\
    \pa_s \td F_{\D, m} & = \Lmic \td F_{\D, m}  
    -    \cP _m [(V \cdot \na _X + 2 \dcm + \td d _\cM)\td F_{\D, M} ]
  + \tf{1}{\es} \cN_{\D}  . 
  \eal
\eeq

\begin{remark}[Improved decay rates]\label{rem:no_error}
The error terms $\eu, \ep, \cE$ in \eqref{eq:non_fix} are canceled in the above equations. 
This enables us to prove that  $E_{\kk+1, \D}$ decays faster than $\es$.
\end{remark}

Similar to Sections \ref{sec:EE_sub_norm_onto}, \ref{sec:sub_norm_Xk}, we  estimate $(\tw_{1,\D}, \td F_{\D, m})$ in energy $E_{\kk+1, \xeta}$ and norm $\cXX^{2\kk}$.
Performing energy estimates on $E_{\kk+1 ,\D}$, we yield 
  \begin{subequations} \label{eq:non_del}
\beq
\bal 
 \underbrace{ \kp \B\la  ( \pa_s - (\cL_{E,s} - \cK_k)  ) \tw_{1,\D}, \tw_{1,\D}    \B\ra_{\cZ^{2\kk+2} } }_{ := I_{\cL, M, 1} }
 & +  \underbrace{ \B\la (\pa_s - \Lmic) \td F_{\D, m}, \td F_{\D, m} \B\ra_{\cYY^{2\kk+2} } }_{ := I_{\cL, m, 1}} \\ 
  & = I_{\cL, M, 2} + I_{\cL, M, 3} + I_{\cL, m, 2} + I_{\cN_{\D}} , 
  \label{eq:non_del:a}\\
\eal 
  \eeq
where $I_{\cL, M, \cdot}, I_{\cL,m,\cdot }$ are the macro and micro linear terms  given by  
\beq \label{eq:non_del:b}
\bal
 I_{\cL, M, 2} +  I_{\cL, M, 3}  & :=  
    \kp \big\la \tw_{1,\D}, \ (- \cI_1, - \cI_2,  \cI_2 )(\td F_{ \D, m} )  \big\ra_{\cZ^{2\kk + 2}}   + \kp \la  (\cL_{E,s} - \cL_E ) \tw_{2,\D}, \tw_{1,\D} \ra_{ \cZ^{2\kk+ 2} }  , \\
 I_{\cL, m, 2}  &:=
 -  \big\la  \cP _m [(V \cdot \na _X + 2 \dcm + \td d _\cM)\td F_{\D, M} ],  \  \td F_{ \D, m}  \big\ra_{\cYY^{2\kk+2}} ,  \\
 \eal 
\eeq
and $I_{\cN_{\D}}$ is the nonlinear term 
\beq \label{eq:non_del:c}
 I_{\cN_{\D}}  := \es^{-1} \la \cN_{\D}, \td F_{ \D, m}\ra_{\cYY^{2\kk+2}} .
\eeq
\end{subequations}

\paragraph{Estimates of linear terms}

The estimates of the linear terms are the same as those in Section \ref{sec:EE_sub_norm_onto}. We apply the linear estimates \eqref{eq:lin_M1}, \eqref{eq:lin_mM2} and \eqref{eq:lin_M3} 
with $(\tw_1, \tw_2, \tFm, E_{\kk+1, \xeta})$ replaced by  $( \tw_{1,\D}, \tw_{2,\D}, \td F_{\D, m}, E_{\kk+1, \D} )$ and use the energy $E_{\kk+1, \D}$ \eqref{eq:nota_del}   to obtain 
\bseq\label{eq:contra_coer}
\begin{align}
\f{1}{2}\f{d}{ds}  \big( \kp \| \tw_{1,\D} \|_{\cZ^{2\kk+2}}^2 + \| \td F_{\D, m} \|_{\cYY^{2\kk+2}}^2  \big) & + \lam_1 \kp \| \tw_{1,\D} \|_{\cZ^{2\kk+2}}^2 + (2 \lame - C \es) \| \td F_{\D, m} \|_{\cYY^{2\kk+2}}^2
+ \f{ \cgam } {6 \es} \| \td F_{\D, m} \|_{\cYL{\xeta}^{2\kk+2}}^2  \notag \\ 
& \leq 	I_{\cL, M,1} + I_{\cL, m,1}  ,  
\end{align}
and 
\beq 
\bal
|  I_{\cL, M, 2} + I_{\cL, m, 2}|
& \les \es^{1/2} ( E_{\kk + 1, \D} + \| \tw_{2,\D} \|_{\cXX^{2\kk+4}}^2 + \f{1}{\es} \| \td F_{m, \D}\|_{\cYL{\xeta}^{2\kk+2} }^2 ) , \\
|I_{\cL, M, 3}| & \les \es^2 E_{\kk+1, \D}^{1/2} \| \tw_{2,\D}\|_{\cXX^{2\kk+4}}
\les \es^2 ( E_{\kk+1, \D} + \| \tw_{2,\D}\|_{\cXX^{2\kk+4}}^2).
\eal
\eeq
Recall the energy $E_{\kk+1, \D}$ from \eqref{eq:nota_del}. Using $\lame > \lam_1$ \eqref{eq:decay_para}, 
$\es \leq \d$ \eqref{eq:EE_para1}, and choosing $\d$ small enough, we simplify the first estimate as 
\beq 
\f{1}{2}\f{d}{ds} E_{\kk+1, \D}  + \lam_1 E_{\kk+1, \D} 
+ \f{ \cgam } { 6 \es} \| \td F_{\D, m} \|_{\cYL{\xeta}^{2\kk+2}}^2  \leq 	I_{\cL, M,1} + I_{\cL, m,1}   . 
\eeq
\eseq

\paragraph{Estimates of nonlinear terms}

The estimate of nonlinear terms are more difficult. 
Since $\cN(\cdot ,\cdot)$ \eqref{eq:non_nota} is bilinear, using the definition of $\cN_{\D}$ in 
\eqref{eq:nota_del}, we obtain 
\[
     \cN_{\D} = \cN( \td F_a - \td F_b, \td F_a ) + \cN( \td F_b, \td F_a - \td F_b )
   = \cN(\td F_{\D}, \td F_a) + \cN( \td F_b , \td F_{\D} ),
   \]
  and further decompose $\td F_{\D}, \td F_a, \td F_b$ into the macro and micro perturbation
\bseq\label{eq:non_contra2}
\beq
\bal 
\cN(\td F_{\D}, \td F_a)
&= \cN(\td F_{\D, m}, \td F_{a,M})
+ \cN(\td F_{\D, M}, \td F_{a, M})
+ \cN(\td F_{\D}, \td F_{a, m}) \\
&:= I_{a, m M} + I_{a, MM} + I_{a, m} , \\
\cN( \td F_b , \td F_{\D} ) 
& = \cN( \td F_b , \td F_{\D, m} ) 
+ \cN( \td F_{b, m} , \td F_{\D, M} ) 
+ \cN( \td F_{b, M} , \td F_{\D, M} ) \\
& := I_{b, m} + I_{b, m M} + I_{b, MM} .
   \eal 
\eeq
We estimate these terms using Theorem \ref{theo: N nonlin_est} with $\eta = \xeta$. Applying \eqref{eq:non_Q:macro} (micro-macro) to $I_{a, m M}$ with $(\eta_1, \eta_2) = (\xeta, \etab)$, 
\eqref{eq:non_Q:micro} ($*$-micro) to  $I_{b, m}$ with $(\eta_1, \eta_2) = (\etab, \xeta)$, and using 
the estimates of $(\tw_{1,\al}, \tw_{2,\al}, \td F_{\al, m})$ with $\al \in \{ a, b\}$ in \eqref{eq:non_contra1}, we obtain
\beq 
\bal 
  |\la  \cN(\td F_{\D, m}, \td F_{a,M}), \td F_{\D, m} \ra_{\cYY^{2\kk + 2}} |  %
 & \les 
\| \td F_{\D, m} \|_{\cYL{\xeta}^{2\kk + 2}}^2 
\| \td F_{a, M} \|_{\cYb^{2\kk + 2}}
\les  \d^{\el}\| \td F_{\D, m} \|_{\cYL{\xeta}^{2\kk + 2} }^2 , \\
 | \la  \cN( \td F_b , \td F_{\D, m} )  ,  \td F_{\D, m} \ra_{\cYY^{2\kk + 2}} |   %
& \les \| \td F_b \|_{\cYb^{2\kk + 2} } \| \td F_{\D, m} \|_{\cYL{\xeta}^{2\kk + 2} }^2
\les \d^{\el} \| \td F_{\D, m} \|_{\cYL{\xeta}^{2\kk + 2} }^2 . %
\eal
\eeq
Applying \eqref{eq:non_Q:macro:MM2} (macro-macro) to $I_{a, MM}$ and $I_{b, MM}$ with $\eta = \etab$ 
and then using $\| q \|_{\cYe^{2\kk-2}} \les \| q \|_{\cYe^{2\kk}}$ and the bound \eqref{eq:non_contral:FMb}, we obtain 
\begin{align}
& 
| \la  \cN(\td F_{\D, M}, \td F_{a, M}) +   \cN( \td F_{b, M} , \td F_{\D, M} ),  \td F_{\D, m}  \ra_{\cYY^{2\kk + 2}} | 
\notag \\
& \les   \B( (\| \td F_{a, M} \|_{\cYY^{2\kk }} + \| \td F_{b, M} \|_{\cYY^{2\kk }} )
\| \td F_{\D, M} \|_{\cYY^{2\kk + 2}}  
  + (\| \td F_{a, M} \|_{\cYY^{2\kk+2 }} + \| \td F_{b, M} \|_{\cYY^{2\kk +2}} )
\| \td F_{\D, M} \|_{\cYY^{2\kk }}  \B)
 \| \td F_{\D, m} \|_{\cYL{\xeta}^{2\kk + 2}}  \notag \\
 & \les  ( \es^{2/3-\ell} \| \td F_{\D, M} \|_{\cYY^{2\kk + 2}} 
 + \es^{1/2-\ell} \| \td F_{\D, M} \|_{\cYY^{2\kk }}  )   \| \td F_{\D, m} \|_{\cYL{\xeta}^{2\kk + 2}}.
\end{align}
\eseq
To estimate $I_{a,m}, I_{b, m, M}$, we need the extra smallness on the dissipation in \eqref{eq:EE_onto:b}. Applying \eqref{eq:non_Q:micro} ($*$-micro)  to $I_{a, m}$ 
and \eqref{eq:non_Q:macro} (micro-macro) to $I_{b, m M}$ with $(l_1, l_2) = (\xeta, \etab)$, we obtain 
\[
\bal
  |\la I_{a, m}, \td F_{\D, m} \ra_{\cYY^{2\kk + 2}} |  
 & = |\la  \cN(\td F_{\D}, \td F_{a, m}), \td F_{\D, m} \ra_{\cYY^{2\kk + 2}} |  
\les \| \td F_{\D} \|_{\cYY^{2\kk + 2}} 
\| \td F_{a, m} \|_{\cYQ^{2\kk + 2}  }  \| \td F_{\D, m} \|_{\cYY^{2\kk + 2}}, \\
|\la I_{b, m M} , \td F_{\D, m} \ra_{\cYY^{2\kk + 2}} |  
& = | \cN( \td F_{b, m} , \td F_{\D, M} ) ,  \td F_{\D, m} \ra_{\cYY^{2\kk + 2}} |  
\les \| \td F_{b, m} \|_{\cYQ^{2\kk + 2} } 
\| \td F_{\D, M} \|_{\cYY^{2\kk + 2}} \| \td F_{\D, m} \|_{\cYY^{2\kk + 2}},
\eal
\]
Using the energy $E_{\kk+1, \D}$ and $\cE_{\D}$ in \eqref{eq:non_del}, we obtain 
\[
\bal 
  \| \td F_{\D} \|_{\cYY^{2\kk + 2}} 
& \les \| \td F_{\D, M}  \|_{\cYY^{2\kk + 2}}  
+ \| \td F_{\D, m}  \|_{\cYY^{2\kk + 2}}   \\ 
& \les \| \td W_{2, \D}  \|_{\cXX^{2\kk + 2}}  
+ \| \td W_{1, \D}  \|_{\cXX^{2\kk + 2}}  
+ \| \td F_{\D, m}  \|_{\cYY^{2\kk + 2}}  \les E_{\kk + 1, \D}^{1/2} + \cE_{\D}, \\
\| \td F_{\D, M} \|_{\cYY^{2\kk}}
& \les \| \td W_{2, \D}  \|_{\cXX^{2\kk }}  
+ \| \td W_{1, \D}  \|_{\cXX^{2\kk }}  
\les \cE_{\D} + \| \td W_{1, \D}  \|_{\cXX^{2\kk }}  .
\eal 
\]
Using \eqref{eq:non_contra2}, the above estimates, and $\es^{2/3-\ell} \leq \es^{2/3-2\ell}$, we derive
\[
\bal
 |\la \cN_{\D}, \td F_{\D, m} \ra_{\cYY^{2\kk}} |
& \les \d^{\el}
\| \td F_{\D, m} \|_{\cYY^{2\kk + 2}}^2
+ \bigl( \es^{2/3- 2\ell } 
( E_{ \kk + 1, \D}^{1/2} + \cE_{\D} )  + \es^{1/2-\ell} ( \| \tw_{1,\D}\|_{\cXX^{2\kk}} + \cE_{\D})  \bigr)\| \td F_{\D, m} \|_{\cYL{\xeta}^{2\kk + 2} } \\
& \qquad + ( \| \td F_{a, m} \|_{\cYb^{2\kk + 2} } 
+ \| \td F_{b, m} \|_{\cYb^{2\kk + 2} } 
) ( E_{\kk + 1, \D}^{1/2} + \cE_{\D} ) \| \td F_{\D, m} \|_{\cYL{\xeta}^{2\kk+2} } .
\eal
\]
Using $\e$-Young's inequality and definition of $I_{\cN_{\D}}$ \eqref{eq:non_del:c}, for any $\upsilon >0$, we establish 
\beq\label{eq:non_contra3}
\bal 
 |I_{\cN, \D} | & = \es^{-1} |\la \cN_{\D}, \td F_{\D, m} \ra_{\cYY^{2\kk + 2}} | \\
& \leq
  \f{ (\ups + C\d^{\el} ) }{\es} \| \td F_{\D, m}\|_{ \cYL{\xeta}^{2\kk + 2} }^2 
+  \f{C}{\es \ups} \B( g(s) ( E_{\kk + 1, \D}  + \cE_{\D}^2) + \es^{1 - 2\ell} ( \| \tw_{1,\D}\|^2_{\cXX^{2\kk}} + \cE_{\D}^2) \B),
\eal 
\eeq
where we denote 
\beq\label{eq:non_contra_g}
 g(s) = \es^{4/3 - 4 \ell } +  \| \td F_{a, m} \|_{\cYb^{2\kk+2} }^2  + \| \td F_{b, m} \|_{\cYb^{2\kk + 2} }^2.
\eeq

\paragraph{Energy estimates in $E_{\kk+1, \D}$}

Plugging  linear estimates \eqref{eq:contra_coer} and \eqref{eq:non_contra3} in \eqref{eq:non_del:a} and using $\es \leq 1$,  we prove 
\[
\bal 
\f{1}{2}\f{d}{ds} E_{\kk+1, \D}  + \lam_1 E_{\kk+1, \D} 
& + \f{ \cgam } { 6 \es} \| \td F_{\D, m} \|_{\cYL{\xeta}^{2\kk+2}}^2  \leq 	I_{\cL, M,1} + I_{\cL, m,1}
= I_{\cL, M, 2} + I_{\cL, M, 3} + I_{\cL, m, 2} + I_{\cN, \D}  \\
& \leq 
C \es^{1/2} ( E_{\kk + 1, \D} + \| \tw_{2,\D} \|_{\cXX^{2\kk+4}}^2  ) 
+ \f{ (\ups + C\d^{\el} + C \es^{1/2} ) }{\es} \| \td F_{\D, m}\|_{ \cYL{\xeta}^{2\kk + 2} }^2  \\
& \quad +  \f{C}{\es \ups} \B( g(s) ( E_{\kk + 1, \D}  + \cE_{\D}^2) + \es^{1 - 2\ell} ( \| \tw_{1,\D}\|^2_{\cXX^{2\kk}} + \cE_{\D}^2) \B) .
\eal
\]
Since $\tw_{2,\D} = \cA_2( \hw_{1,\D} )$, applying Lemma \ref{lem:A2} and the definition of $Y$ norm 
\eqref{norm:fix}, we obtain 
\beq\label{eq:est_W2_del}
\bal 
  \cE_{\D} = \| \tw_{2,\D} \|_{\cXX^{2\kk+6}} 
   & \les \operatorname e^{-  \lams s}  \sup_{s\geq 0} \operatorname e^{ \f23  \rE s  } 
 \|  \hw_{1, \D} \|_{\cXX^{2\kk}}  \\
 & \les \operatorname e^{-  \lams s  }  \d^{ \f23 } \|  \hw_{1,\D}\|_{Y}
 \les \es^{ \f23-\ell}  \|  \hw_{1,\D}\|_{Y} . 
\eal 
\eeq

Combining the above two estimates, bounding $\es^{1/2} \les g(s) / \es$ due to \eqref{eq:non_contra_g}, 
choosing $\ups = \f{1}{100}\cgam $, 
and $\d$ small enough so that $\es$ is small by \eqref{eq:EE_para1}, we obtain 
\begin{align}\label{eq:non_contra_EE}
  \f{1}{2} \f{d}{d s} E_{\kk + 1, \D} 
    & \leq \left( - \lam_1 + \f{C}{\es }  g(s)   \right)  E_{\kk + 1, \D}  
- \f{\cgam}{8\es} 
\| \td F_{\D, m}\|^2_{\cYL{\xeta}^{2\kk+2}}   \\
    & \qquad   + C (g(s) + \es^{1-2\ell}) \es^{ \f13 -2\ell }  \|  \hw_{1,\D}\|_{Y}^2
    + C \es^{-2\ell} \| \tw_{1,\D} \|_{\cXX^{2\kk}}^2, \notag 
\end{align}
    where $C$ is some absolute constant (depending on $\kk, \etab, \xeta$). 

\paragraph{Energy estimates in $ \cXX^{2\kk }$ }

Note that the equation of $\tw_{1,\D}$ \eqref{eq:F_diff} is linear and has the same form as that of $\tw_1$ in \eqref{eq:non_W:a} except that we do not have the error term in \eqref{eq:F_diff}. 
Applying the estimate of $I_{\cL, M, 1}, I_{\cL, M, 2}, I_{\cL, M, 3}$ in \eqref{eq:EE_Wk0} and 
\eqref{eq:EE_Wk} with $(\tw_1, \tw_2, \tFm)$ replaced by $(\tw_{1,\D}, \tw_{2,\D}, \td F_{\D, m})$, 
and using $R_s^{-r} \les \es^2$ from Remark \ref{lem:para_bd}, we obtain  $ \cXX^{2\kk }$ estimates of $\tw_{1,\D}$
\[
\bal
 \tf{1}{2}	\tf{d}{ds} \| \tw_{1, \D} \|_{ \cXX^{2\kk} }^2 
 \leq & - \lame   \| \tw_{1,\D} \|_{\cXX^{2\kk}}^2 +
 C \es^2 \| \tw_{1,\D}\|_{ \cXX^{2\kk} } \| \tw_{2,\D} \|_{ \cXX^{2\kk} } %
 + C    \es^{ \f12 } \cdot \f{1}{\es^{\f12}} \| \td F_{\D, m}\|_{\cYL{\xeta}^{2\kk+2}}   \| \tw_{1,\D} \|_{\cXX^{2\kk}} .
 \eal
 \]

We do not have an error term similar to $I_{\cE, M}$ in \eqref{eq:EE_Wk0} since there is no error term in \eqref{eq:F_diff}. Using \eqref{eq:est_W2_del} and Cauchy--Schwarz inequality, we obtain 
\beq\label{eq:non_contra_EE2}
\bal 
	 \tf{1}{2}	\tf{d}{ds} \| \tw_{1, \D} \|_{ \cXX^{2\kk} }^2 
	& \leq (-\lame + \es^2 + \es^{\ell }) \| \tw_{1, \D} \|_{ \cXX^{2\kk} }^2 
	 + C \es^2 \| \hw_{1,\D}\|_Y^2  \\
  & \qquad  + C \es^{1- \ell} 
\cdot \tf{1}{\es} \| \td F_{\D, m}\|_{\cYL{\xeta}^{2\kk+2}}^2 .
\eal 
\eeq

\paragraph{Summary of the estimates}

Recall $\es = \d \operatorname e^{-\rE s  }$. We estimate the mix energy 
\beq\label{eq:EE_contra_norm}
	E_{\mw{mix},\D} := \es^{- 1 } E_{\kk+1, \D} + \es^{-4/3}  \| \tw_{1,\D}\|_{\cXX^{2\kk}}^2.
\eeq

From \eqref{eq:EE_para1}, for any $b$, we have 
\[
\tf12 \tf{d}{ds} \es^{-b} =\tf12 b \rE \es^{-b} .
\]
Combining \eqref{eq:non_contra_EE} $\times \es^{-1}$ and \eqref{eq:non_contra_EE2} $\times \es^{-4/3}$, we obtain 
\[
\bal
	\f{1}{2} \f{d}{ds} E_{ \mw{mix}, \D}
	\leq & \left(-\lam_1 + \f{1}{2}\rE + \f{C}{\es} g(s)\right) \f{1}{\es}  E_{\kk+1, \D} 
	+ \left( - \lame + \f{2}{3} \rE + C \es^{\ell} \right) \es^{-\f43}  \| \tw_{1,\D}\|_{\cXX^{2\kk}}^2 
	- \f{\cgam}{ 8 \es^{2} } \| \td F_{\D, m}\|^2_{\cYL{\xeta}^{2\kk+2}} \\
&  	+ C \es^{ 1 -\ell - \f43} \cdot  \f{1}{\es} \| \td F_{\D, m}\|_{\cYL{\xeta}^{2\kk+2}}^2 
	 + C\left( \f{ g(s)}{\es} \cdot \es^{ \f13- 2\ell}  + \es^{\f13-4\ell} +  \es^{\f23} \right) \| \hw_{1,\D}\|_{Y}^2 .
\eal
\]

 Recall $\ell = 10^{-4}$ from \eqref{eq:EE_para_k}. Using $\lam_1 > \rE$ \eqref{eq:decay_para}, 
\[
 \f{1}{3}- 2 \ell > 0, \quad  \es^{ 1 -\ell - \f43} < \es^{\ell -1}, \quad \es^{2/3} +  \es^{1/3-4\ell} \les g(s) \es^{ - 1}
\]
by \eqref{eq:non_contra_g}, $\es \leq \d$ by \eqref{eq:EE_para1}, and taking $\d$ small enough, we obtain 
\[
\bal
 	\f{1}{2} \f{d}{ds} E_{\mw{mix}, \D} 
 & \leq  (-\lam_1 + \f{2}{3} \rE + \f{C}{\es} g(s) + C \es^{\ell} ) E_{\mw{mix}, \D}
  - \f{ \cgam - C \es^{ \ell} }{ 8 \es^{2}} \| \td F_{\D, m}\|^2_{\cYL{\xeta}^{2\kk+2}} + \f{C g(s)}{\es} \|  \hw_{1,\D}\|_{Y}^2 \\  
  & 
 \leq   \f{C g(s)}{\es}    E_{\mw{mix}, \D} + \f{C g(s)}{\es} \|  \hw_{1,\D}\|_{Y}^2 .
  \eal
\]

Denote $G(\tau) = C \int_{0}^{\tau} \f{1}{\es} g(s) d s$. 
For any $ 0 \leq s_1 < s_2$, using \eqref{eq:non_contra1:b} and $\el < \f{1}{3}$ \eqref{eq:EE_para_k}, we obtain
\[
0 \leq  G(s_2) - G(s_1)
= \int_{s_1}^{s_2} \f{1}{\es} g( s ) d s 
  \les \int_{0}^{\infty} \es^{1/3-2\ell} 
  + \d^{\el} 
\les \d^{1/3-2\ell} d s + \d^{\el} \les \d^{\el}.
\]
Since $(\tw_{1,\D}, \td F_{\D, m} )|_{s = 0} = 0 $, we have $E_{\kk+1, \D}(0) = 0, E_{\mw{mix}, \D}(0)=0$. Using Gr\"onwall's inequality,  we establish  
\[ 
\bal
E_{\mw{mix}, \D}(s)
  &\leq C \int_{0}^s \operatorname e^{G(s) - G(\tau)  } 
   \cdot \f{g(\tau)}{\e_{\tau}} d \tau  \cdot  \|  \tw_{1,\D}\|_{Y}^2 
    \leq C \int_{0}^s 
  \cdot \f{g(\tau)}{\e_{\tau}} d \tau 
  \cdot 
  \|  \tw_{1,\D}\|^2_{Y} \leq C \d^{\el}   \|  \tw_{1,\D}\|_{Y}^2 .%
  \eal
\]
Using the definition \eqref{eq:EE_contra_norm} and taking $\d$ small enough, we prove 
\[
	\| \tw_{1,\D}\|_{\cXX^{2\kk}} \leq  \es^{2/3} E_{\mw{mix}, \D}^{1/2} \leq C \es^{2/3} \d^{\ell/2}  \|  \tw_{1,\D}\|_{Y} 
	\ll \f{1}{2}  \es^{2/3} \|  \tw_{1,\D}\|_{Y} .
\]

The decay estimates of $E_{ \mw{mix}, \D}$ can be improved, but we do not need such an improvement.
We have proved \eqref{eq:contra_goal} and concluded the proof of Proposition \ref{prop:contra}.

\section{Local well-posedness of the fixed point equations}\label{app:LWP_fixed}

In this section, we show that the fixed point equations \eqref{eq:non_fix} and the Landau equation 
\eqref{eq:Landau} admit a local-in-time solution, by constructing a solution to the following system with an appropriate initial value:
\beq\bal
    \label{eq:fix_F_1}
    \left(\pa_s + \cT + \dcm - \f32 \bcv\right) \td F  
    &= \f1\es \left[ 
        \cL _\cM (\td F)
        + \cN (\td F, \td F)  
    \right] - \cMM ^{-1/2} \eM \\
    & \qquad +  \sss \cdot \cF_M \circ \cK_k (\hw _1 + \cA _2 (\hw _1) - \cF_E(\td F) ).
\eal\eeq
Here the data $\hw _1$ is given, $\cA _2$ is defined in \eqref{eq:T2:def}, and parameter $\sss \in \{0, 1\}$. Note that the linearized self-similar Landau equation \eqref{eqn: linearized-td-F} corresponds to $\sss = 0$. 

\subsection{Reformulation of the fixed-point equations}

Firstly, we show that given  $\hw_1$, the fixed point equations \eqref{eq:non_fix} are equivalent to \eqref{eq:fix_F_1} with $\sss = 1$. We consider \eqref{eq:fix_F_1} since it is easier to establish the local well-posedness.

Given $\hw_1$, we recall the fixed-point equations of $(\tw_1, \tw_2, \tFm)$ \eqref{eq:non_fix} as follows 
  \bseq\label{eq:fix_recall}
  \begin{align}
    \tw_2 &= \cA_2(\hw_1),  \label{eq:fix_recall:a}  \\  
\pa_s \tw_1 & =  (\cL_{E, s} - \cK_k ) \tw_1   + (\cL_{E,s} - \cL_E) \tw_2  
- (\cI_1, \cI_2,  -  \cI_2  )(\tFm)
   - ( \cs^3 \eu, \cs^3  \ep, 0) ,  \label{eq:fix_recall:b}  \\
  \pa_s \tFm & = \Lmic \tFm - \cP _m [(V \cdot \na _X + 2 \dcm + \td d _\cM) \tFM]
  + \f{1}{\es} \cN( \td F, \td F ) - \cP_m(\cMM^{-1/2} \cE_M) ,
   \label{eq:fix_recall:c}
 \end{align}
where $\Lmic$ is defined in \eqref{eqn: defn-Lmic}, 
and $\cA_2$ is 
defined in 
\eqref{eq:T2:def}. Using $\tw = \tw_1 + \tw_2$,
\[
  \cK_k (\hw_1 - \tw_1) = \cK_k( \hw_1 + \tw_2 - \tw) , %
\]
we rewrite the equation of $\tw$ \eqref{eq:non_fix:W} as
\beq\label{eq:fix_recall:d}
    \tag{\ref*{eq:fix_recall}b'}
      \pa_s \tw = \cL_{E, s} \tw  + \sss \cdot \cK_k( \hw_1 + \tw_2 - \tw )
       - (\cI_1, \cI_2,  -  \cI_2 )(\tFm) 
   - ( \cs^3 \eu, \cs^3  \ep, 0) ,
\eeq
\eseq
with  $\sss = 1$. %
Given $\hw_1$, we construct $\tw_2$ using \eqref{eq:fix_recall:a}. Then the system of 
$(\tw_1, \tFm)$ in \eqref{eq:fix_recall:b}, \eqref{eq:fix_recall:c} is equivalent to that of $(\tw, \tFm)$ in  \eqref{eq:fix_recall:d}, \eqref{eq:fix_recall:c}. 

We argue that the above system \eqref{eq:fix_recall:d}, \eqref{eq:fix_recall:c} is equivalent to 
\bseq\label{eq:fix_F}
\begin{align}
   (\pa_s + \cT) (\cMM ^{\f12} \td F)
    = & \f{1}{\es} \B[ 
    Q(\cM, \cMM^{\f12} \td F) 
    + Q(\cMM^{\f12} \td F, \cM)
    + Q(\cMM^{\f12} \td F, \cMM^{\f12} \td F)  
    \B] - \eM + I
\end{align}
via $\td F = \cF _M (\tw) + \tFm$, where $I$ is defined as 
\beq
	 I :=  \sss \cdot \cMM^{\f12} \cF_M \circ \cK_k(\hw_1 + \tw_2 - \cF_E(\td F) ).
\eeq
\eseq
Note that \eqref{eq:fix_F} differs from \eqref{eq:lin} by the $I$ terms. 

First, we show \eqref{eq:fix_F} implies the system \eqref{eq:fix_recall:d}-\eqref{eq:fix_recall:c}.
In fact, since $I$ is purely macroscopic, following the derivations in Lemma 
\ref{lem:Deri-Fm} for the equations of $\tFm$ by first dividing $\cMM^{1/2}$ and then applying projection $\cP_m$, we obtain \eqref{eq:fix_recall:c} from \eqref{eq:fix_F}. 
Using the map $\cF_E$ \eqref{eq:pertb_macro} and the relations \eqref{eq:macro_UPB}, \eqref{eq:macro_UPB_inv}, we obtain 
\beq\label{eq:fix_Kk}
\bal
	 \B\la \sss \cdot \cMM^{1/2} \cF_M \circ \cK_k(\hw_1 + \tw_2 - \cF_E(\td F) ),& \  
	   \B( \tf{V -\bu}{\cs},  \tf{| V - \bu|^2}{3 \cs^2} ,    1 - \tf{| V - \bu|^2}{3 \cs^2}  \B)  \B\ra_V \\
	 &  = \sss \cdot \cF_E \circ \cF_M \circ \cK_k(\hw_1 + \tw_2 - \cF_E(\td F) ) \\ 
	 & = \sss \cdot  \cK_k(\hw_1 + \tw_2 - \cF_E(\td F) ).
\eal
\eeq

Note that equations \eqref{eq:lin_euler} are derived by integrating \eqref{eqn: lin-with-M11/2} 
against $1, \f{V-\bu}{\cs}, \f{|V-\bu|^2}{3 \cs^2}, 1 -  \tf{|V-\bu|^2}{3 \cs^2}$ over $V$ (see \eqref{eq:pertb_macro}). Integrating \eqref{eq:fix_F} against $1, \f{V-\bu}{\cs}, \f{|V-\bu|^2}{3 \cs^2}, 1 -  \tf{|V-\bu|^2}{3 \cs^2}$ over $V$, applying the same derivations (see Appendix \ref{app:euler_deri}), and using the integrals of $I$ over $V$ in \eqref{eq:fix_Kk}, 
we derive \eqref{eq:fix_recall:d}. 
Thus, \eqref{eq:fix_F} implies the system \eqref{eq:fix_recall:d}-\eqref{eq:fix_recall:c}. 

Using the $\tw$-equation \eqref{eq:fix_recall:d} and the relation \eqref{eq:macro_UPB}, we derive the equations of $\cF_M(\tw)$. Along with the equation of $\tFm$ \eqref{eq:fix_recall:c}, we can derive the equation \eqref{eq:fix_F}. The derivations are similar and are thus omitted.

Dividing  \eqref{eq:fix_F} by $\cMM ^{1/2}$ and using the notation $\cL_{\cM}, \cN$ \eqref{def:T-LM}, \eqref{eq:non_nota}, we obtain
\eqref{eq:fix_F_1}. With $\td F$ being a solution to the {\em{nonlinear}} problem \eqref{eq:fix_F_1}, we construct the solution $\tw_1 = \cF _E \cP _M \td F - \tw_2$ and $\tFm = \cP_m \td F$ to the system \eqref{eq:fix_recall}.

The main result in this section is the following local existence theorem.

\begin{theorem}\label{thm:LWP}
There exists absolute constants $ 0< \ze_2 < \ze_1 < 1$  such that the following statement holds. Consider equation
\eqref{eq:fix_F_1} with $\sss \in \{0, 1 \}$, $k\geq \kkl$ with initial data
\beq\label{eq:LWP_ass_IC}
	\td F(0) \in \cYb^{k} , \quad \| \td F(0)\|_{\cYb^{\kkl}} < \ze_2. 
\eeq
When $\sss=1$, we further assume that $\hw_1$ satisfies 
$\| \hw_1 (s) \|_{\cXX^{k}} <\es^{2/3} \d _0 ^{2\ell}$ with $\d _0$ given in Proposition \ref{prop:contra}. 
There exists a unique local solution $ \td F \in L^{\infty}( [0,  T], \cYb^{ k } ) 
\cap L^2( [0,  T], \cYQ^{k})$ to \eqref{eq:fix_F_1} with $T \asymp \min \{\e_0, 1 \}$ and 
\beq\label{eq:LWP_thm_EE}
 	\| \td F(s) \|_{\cYb^{\kkl}} \leq \min\left\{ \ze_1 , \, C ( \| \td F(0) \|_{\cYb^{\kkl}} + \e_0^{-1} s )  \right\} ,
    \quad s \in [ 0 ,   T].
\eeq
Moreover, the solution can be continued beyond $s \in [ 0,   T_*)$ in the same regularity class if 
  \beq\label{eq:LWP_blowup}
   \sup_{s \in [ 0,  T_*)} \|\td F(s)\|_{\cYb^{\kkl}} < \ze_2 .
\eeq

\end{theorem}

The solution satisfies the energy estimates \eqref{eq:LWP_EE_non}. Since we develop much sharper estimates on $\td F$ when $\sss=1$ in Section \ref{sec:solu}, we do not derive the explicit bounds 
in \eqref{eq:LWP_EE_non} when $j > \kkl$.
Note that we only require smallness in $\cYb^{\kkl}$ norm, but \emph{not} the higher order $\cYb^k$ norm. 

From the above assumption on $\hw_1$ and \eqref{eq:est_W12},  $\hw_1, \tw_2$ 
satisfy the following estimate  for any $s$
\beq\label{eq:LWP_W1W2}
	\| \hw_1(s) \|_{\cXb^{ k }} < \es^{2/3} \les 1, \quad  \| \tw_2 \|_{\cXb^{ k +6} } \les  \es^{2/3- \ell} \les \es^{1/3} \les 1.
\eeq

Based on Theorem \ref{thm:LWP}, we establish the following local existence results for the Landau equation with 
a solution satisfying a Gaussian lower bound.

\begin{proposition}\label{prop:LWP_landau}
Consider $F = \cM + \cMM^{1/2} \td F$. Suppose that $\td F(0)$ satisfies \eqref{eq:LWP_ass_IC} with some $k \geq \kkl$ and $F(0) > 0$. The Landau equation \eqref{eq:LC_ss} admits a unique local solution $ \td F \in L^{\infty}( [ 0,  T], \cYb^{ k } ) 
\cap L^2( [ 0,   T], \cYQ^{k})$ to \eqref{eq:fix_F_1} with $T \asymp \min \{ \e_0, 1 \}$
and $F(s) \geq 0$. The nonnegative solution extends beyond $s \in [ 0,  T_*)$ in the same regularity class 
as long as condition \eqref{eq:LWP_blowup} is satisfied.

Let $\cM$ be the time-dependent local Maxwellian constructed in \eqref{eq:local_max}. There exists $a _0 \geq 1$ such that the following holds. If the initial data satisfy
\beq\label{eq:Gauss_lower_IC}
F(0, X, V)  \geq c \cdot \ang X^{-l} \cM (0, X ,V) ^{\aaa} ,\quad \forall \ (X, V) \in \R^6, 
\eeq
for some $l \in [0, 100],  \aaa \geq \aaa_0$, $c > 0$, then there exists $\bbb \gtr \aaa^2$, such that
\beq\label{eq:Gauss_lower}
F(s, X, V) \geq c \cdot \exp(  \bbb ( \e_0^{-1} - \es^{-1}  ) ) \ang X^{-l} \cM(s,X,V)^\aaa ,
\eeq
 for any $(X, V)\in \R^6$ and $s \in [ 0 ,  T_*)$, where $[ 0,  T_*)$ is the maximal interval on which \eqref{eq:LWP_blowup} holds. 
\end{proposition}

\begin{remark}[\bf Local $C^{\infty}$ solutions in the physical variables]
Since assumption \eqref{eq:LWP_ass_IC} imposes smallness \emph{only} on $\cYb^{\kkl}$ norm 
for a fixed $\kkl$, we can choose $\td F(0) \in \cap_{k \geq \kkl} \cYb^k$ in Proposition \ref{prop:LWP_landau} in the case $\sss =0$.  Since $\cMM^{1/2}, \cM \in C^{\infty}$, using the embedding estimates in Lemma \ref{lem:prod}, we obtain 
\beq\label{eq:LWP_upper}
 |D^{ \leq {k-2d}} F(s,X,V)| \les_{k} \mu(\vc)^{1/4} ( 1 + \| \td F(s) \|_{\cYb^{k }} ) ,
\eeq
where $\cM(s, X, V) = \mu(\vc)$ and $\mu(\cdot)$ is the Gaussian defined in \eqref{eq:gauss}.
See the estimates in \eqref{eq:solu_small_point} and \eqref{eq:fxv_est1}. 
As a result, the local solution $F$ corresponding to the initial perturbation $\td F(0)$ is $C^{\infty}$.
For any $ s < \infty$, since the physical time $t$ \eqref{eq:SS} satisfies $t< 1$,  we obtain 
\[ 
|X| \asymp_s |x|, \quad |V| \asymp_s |v|, \quad  \operatorname e^{-C_{1,s} |v|^2}\les_s \mu(\vc) \les_s \operatorname e^{- C_{2,s} |v|^2} .
\]
Using the self-similar transform \eqref{eq:SS} and Proposition \ref{prop:LWP_landau}, we construct a local smooth solution $f$ to the Landau equation \eqref{eq:Landau} with the uniform Gaussian decay \eqref{eq:LWP_upper} and the Gaussian lower bound \eqref{eq:Gauss_lower}.
\end{remark}

In Section \ref{sec:LWP_iterative}, we reformulate solving the nonlinear equation \eqref{eq:fix_F_1} as a fixed point problem 
and perform uniform energy estimates. 
We prove the local existence of solution in Section \ref{sec:LWP_exist}, the continuation criterion \eqref{eq:LWP_blowup} in Section \ref{sec:LWP_continue}, and Proposition \ref{prop:LWP_landau} in Section \ref{sec:LWP_landau}.

\subsection{Iterative scheme and uniform energy estimates}\label{sec:LWP_iterative}
Let us rewrite the linear operator $\cL _\cM (\td F)$ defined in \eqref{def:T-LM} as: 
\begin{align*}
\cL _\cM (\td F) &= \cMM ^{-1/2} Q (\cM, \cMM ^{1/2} \td F) +  \cMM ^{-1/2} Q (\cMM ^{1/2} \td F, \cM) \\
&= \cN (\rhos \cMM ^{1/2}, \td F) + \cN ( \td F, \rhos \cMM ^{1/2}) \\
&= \cN _{1} (\rhos \cMM ^{1/2}, \td F) + \dots + \cN _{6} (\rhos \cMM ^{1/2}, \td F) + \cN ( \td F, \rhos \cMM ^{1/2}). 
\end{align*}
where $\cN _i$, $1 \le i \le 6$, are defined in \eqref{N(f,g)}.
Solution to \eqref{eq:fix_F_1} can be regarded as a solution to the linear equation below with $\td G = \td F$:
\begin{align}
    \label{eqn: linear-transport_1}
    \left(\pa _s + \cT + \dcm - \f32 \bcv\right) \td F = \f1{\es} \left[(\cN _1 + \cN _5)(\rhos \cMM ^{1/2}, \td F) \right] + \f1{\es} \cN (\td G, \td F) + \td H,
\end{align}
with 
\begin{align}\label{tildaH}
    \td H  = \td H (\td G, \cMM) =  %
    \f1\es \left[ (\cN _2 + \cN _3 + \cN _4 + \cN _6)(\rhos \cMM ^{1/2}, \td G) \right] + \f1\es \cN (\td G, \rhos \cMM ^{1/2}) \\
    - \cMM ^{-1/2} \eM + \sss \cdot \cF_M \circ \cK_k(\hw_1 + \tw_2 - \cF_E (\td G)) . \nonumber
\end{align}
Note that 
\begin{align*}
    (\cN _1 + \cN _5)(\rhos \cMM ^{1/2}, \td F) = \div (A [\cM] \na _V \td F) - \kappa _2 ^2 \cs ^{-2} A [\cM \vc \otimes \vc] \td F,
\end{align*}
so
\begin{align}
\label{eqn: N1+N5=sigma}
\left \la (\cN _1 + \cN _5)(\rhos \cMM ^{1/2}, \td F), \td F \right\ra _V &= - \| \td F \|^2_\s,
\end{align}
and by the same computation as  \cite[Page 396]{guo2002landau}, 
\begin{align}
    \label{eqn: N2346}
    (\cN _2 + \cN _3 + \cN _4 + \cN _6)(\rhos \cMM ^{1/2}, \td G) = \kp _2 \cs ^{-1} (\div A) [\cM \vc] \td G.
\end{align}

\subsubsection{Functional spaces}

We find the solution to (\ref{eq:fix_F_1}) as the fixed point of the map 
\bseq\label{eq:LWP_map}
\beq
\mathscr T: \td G \in \boldsymbol J _\ze ^k \mapsto \td F \in \boldsymbol J _\ze ^k 
\eeq
where for some $k \ge \kkl$, $\ze > 0$, $T > 0$ small to be chosen, we denote
\beq
    \boldsymbol J_\ze^k :=  \B\{u \in \boldsymbol{\cY} _{\etab} ^k: \| u \| _{\boldsymbol{\cY} _{\etab} ^{\kkl}} \le \ze \B\}, \qquad 
    \boldsymbol{\cY} _{\etab} ^k = L^\infty (0, T; \cYb ^k) \cap L ^2 (0, T; \cYQ ^k),
   \eeq
where we define the $T$-dependent $\boldsymbol{\cY} _{\etab} ^k$ norm as 
\beq\label{eq:LWP_map:c}
    \| u \|^2 _{\boldsymbol{\cY} _{\etab} ^k} := \sup _{s \in [0, T]} \| u (s) \| _{\cYb ^k} ^2 + \int _0 ^T \f1\es \| u (s) \| _{\cYQ ^k} ^2 d s.
 \eeq
 \eseq
Since we only use the weight with exponent $\etab$ throughout this section, we do not indicate its dependence in the new norms and functional spaces, e.g. $\bbJ_{\ze}^k$. We will choose the life span $T$ depending on the size of $\e _0$, as in Theorem \ref{thm:LWP}.

For $\mathscr T$ to have a fix point, we need to show (1) $\eqref{eqn: linear-transport_1}$ has a unique solution in $\bbY^k$; (2) $\|\td F \| _{ \bbY^\kkl} \le \ze$ whenever $\| \td G \| _{\bbY ^\kkl} \le \ze $; and (3) $\mathscr T$ is contractive in $\bbJ_{\ze} ^k$.

\subsubsection{Localization and regularization}

Recall that $\chi _R: \R ^3 \to \R$ is a smooth, radial cutoff function supported in $B _{2 R}$ with $\chi _R = 1$ in $B _R$ defined in Section \ref{sec: modified-profile}. 
Specifically, we define the cutoff function by
\begin{align}
    \label{eqn: defn-varphi}
    \varphi _R (s, X, V) := \chi _R (X) \chi _R (\vc).
\end{align}
Then, inside the support of $\varphi _R$, it holds
\begin{align*}
    |X| ^2 + |V| ^2 \le (2 R) ^2 + C \cs ^2 \ang \vc ^2 \le 4 R ^2 + C (4 R ^2 + 1) \le C \ang R ^2 =: R ^{*2}
\end{align*}
so $\varphi _R$ is supported in $[0, \infty) \times B _{R ^*}$, where $B _{R ^*}$ denotes a ball in $\R ^6$ with radius $R ^*$. 

We now compute the derivatives of the cut-off. The gradient in $V$ of $\varphi _R$ is
\begin{equation}
\label{eqn: nabla_V_varphi}
\begin{aligned}
    \na _V \varphi _R (s, X, V) &= 
    \chi _R (X) \na _V \chi (\vc / R) 
    \\
    &= \chi _R (X) \cdot \f1{\cs R} \cdot \na \chi (\vc / R) \\
    &= \chi _R (X) \cdot \f1{\cs R} \cdot \chi _\xi (\vc / R) \cdot \f{\vc}{|\vc|} \\
    &= \td \varphi _R (s, X, V) \cdot \kp _2 \cs ^{-1} \vc,
\end{aligned}
\end{equation}
where $\chi _\xi$ means the radial derivative of $\chi$, and we introduce
\begin{align}
    \label{eqn: defn-tdvarphi}
    \td \varphi _R (s, X, V) := \f1{\kp _2 R ^2} \cdot \chi _R (X) \cdot \f{\chi _\xi (\vc / R)}{|\vc| / R}.
\end{align}

Next, we show that the smooth cut off function enjoys the derivative bound:
\begin{align}
    \label{eqn: uniform-bound-cutoff}
    |D ^{\al, \b} \varphi _R| + R ^2 |D ^{\al, \b} \td \varphi _R| \les _{\al, \b} 1.
\end{align}
In particular, the upper bound is uniform in $R$ for any $ R> 0$. To prove this, note that 
\begin{align*}
    D ^\al _X \chi _R (X) = R ^{-|\al|} \ang X ^{|\al|} \pa _X ^\al \chi (X / R).
\end{align*}
For any $\al \succ 0$, $\pa _X ^\al \chi$ is supported in $B _2 \setminus B _1$, so $D ^\al _X \chi _R$ is bounded. Similarly, when $|\al| + |\b| = 1$,  
\begin{align*}
    D ^{\al, \b} \chi _{R} (\vc) = \na \chi (\vc / R) \cdot D ^{\al, \b} \vc / R.
\end{align*}
Note that $|D ^{\al, \b} \vc| \les _{\al, \b} \ang \vc$ (see Lemma \ref{lem:est_Dxvring} for $\b = 0$ and Remark \ref{rem: Dab-Vc} for $|\b| > 0$), so it is bounded by $C R$ in the support of $\na \chi (\vc / R)$, thus $D ^{\al, \b} \chi _{R} (\vc)$ is bounded. By induction in view of Corollary \ref{cor: concatenate-derivative}, for general multi-index $\al, \b$ with $|\al| + |\b| = k$, we have $D ^{\al, \b} \chi _R (\vc)$ is bounded, and our claim \eqref{eqn: uniform-bound-cutoff} is proven for $\varphi _R$. The proof for $\td \varphi _R$ is identical so we do not repeat here. 

Regarding the material derivative of $\phi _R$, i.e. $(\pa _s + \cT) \phi _R$, it equals to  
\begin{align*}
    (\pa _s + \cT) \chi _R (\vc) = (\pa _s + \cT) \vc \cdot \na \chi _R (\vc) &= (\pa _s + \cT) \vc \cdot \f1R \chi _\xi (\vc / R) \cdot \f{\vc}{|\vc|}.%
\end{align*}
Lemma \ref{lem:ds_t} gives
\begin{align*}
    |(\pa _s + \cT) \chi _R (\vc)| \les R ^{-1} \ang X ^{-r} \ang \vc ^2 \les R ^{-1} \ang \vc ^2.
\end{align*}
Moreover, by direct computation 
\begin{align*}
    (\pa _s + \cT) \chi _R (X) = \f{\bcx X + \bu + \cs \vc}R \cdot \na \chi (X / R) \les R ^{-1} \ang X \ang \vc 
\end{align*}
using $|\vc| \le 2 R$ in the support of $\chi _R (\vc)$. We can iterate the estimate for higher derivatives $\al, \b$, and summarize
\begin{align}
    \label{eqn: ds_T_cutoff}
    |D ^{\al, \b} (\pa _s + \cT) \vp _R| \les R ^{-1} \ang X \ang \vc ^2.
\end{align}

We consider localized initial data 
\beq\label{eq:LWP_init}
    \td F_{\init, R}(X, V) = \td F_{\init}(X, V) \cdot \varphi_R (0, X, V). 
\eeq
Clearly, $\td F_{\init, R} \in C ^\infty _0 (B _{R ^*})$ and $\td F_{\init, R} \to \td F_{\init}$ in $\cYb ^k$ as $R \to \infty$.  

Next, we show that (\ref{eqn: linear-transport_1}) equipped with regularizing terms $\th \D _{X, V}$ with weight has an unique smooth solution $\td F _{\th, R}$ for $(X, V) \in B _R$ such that $\td F _{\th, R} =0$ on $\partial B _R$. 
Specifically, we consider 
\beq\begin{aligned}\label{eq: linear_viscous}
    \left(\pa _s + \cT + \dcm - \f32 \bcv \right) \td F _{\th, R} &= \f1{\es} \left[(\cN _1 + \cN _5)(\rhos \cMM ^{1/2}, \td F_{\th, R}) \right] + \f1{\es} \cN (\td G, \td F _{\th, R}) + \td H + \th \D _W \td F _{\th, R}
\end{aligned}\eeq
with $\th > 0$, where the weighted Laplacian is defined as 
\begin{align}\label{eq: weighted-laplacian}
    \D _W F & = -   \nu ^{-1} \ang X ^{2} \ang \vc ^{4} F   +   \sum _{|\al_1|  +  |\b_1| = 1}  \ang X ^{1-\etab} \ang \vc ^2  
    \pa_X^{\al_1} \pa_V^{\b_1}  \B( \vp_1^{2 |\al_1| } \cs^{2 |\b_1|} \ang X^{\etab}  \pa_X^{\al_1} \pa_V^{\b_1}  (\ang X \ang \vc ^2 F)  \B)  ,
\end{align}
and $\nu$ is the parameter for the $\cYb$-norm \eqref{norm:Y} chosen in 
Theorem \ref{thm: micro-Hk-main}. While $\D_W F$ may seem complicated, it is in divergence form relative to the $\cYb$ norm.

For this weighted diffusion term, we have the following estimate. The proof of  Lemma \ref{lem: diffusion-innerprod} 
follows by applying integration by parts and tracking the main terms. We defer it to Appendix \ref{sec: proof-diffusion-innerprod}.

\begin{lemma}[Weighted diffusion]
    \label{lem: diffusion-innerprod}
    There exist constants $C_k \ge 0$ such that for compactly supported $h \in \cYb ^{k + 1}$, we have
    \footnote{
The constant $C_k$ depends on the parameter $\nu$ chosen in the $\cYb$-norm \eqref{norm:Y}. Since we have 
determined $\nu$ as some small constant in Theorem \ref{thm: micro-Hk-main}, $C_k$ can be treated as a constant that depends only on $k$.
    }
    \bseq        \label{eqn: DW-Yk}
    \begin{align}
        \la \D_W h , h \ra_{\cYb} & = - \| \ang X \ang \vc ^2 h \|_{\cYb^1} ^2 ,  \\ 
        \left\la 
            \D _W h, h
        \right\ra _{\cYb ^k} & \le - \f12 \| \ang X \ang \vc ^2 h \| _{\cYb ^{k + 1}} ^2
        +  C _k \one_{k>0} \| \ang X \ang \vc ^2 h \| _{\cYb ^k} ^2. 
    \end{align}
    \eseq
\end{lemma}

Next, we show that \eqref{eqn: linear-transport_1} is parabolic in $V$.

\begin{lemma}[Parabolicity]
\label{lem:parab}
There exists $\ze_0 \in (0, 1)$ such that for any $\| \td G \| _{\cYb^{\kkl}} \leq \ze_0$, we have 
\[
  \f{1}{2} A[\cM] \preceq	A[\cM + \cMM^{1/2} \td G ] \preceq  \f{3}{2} A[\cM].
\]
\end{lemma}

\begin{proof}
Recall $\etab = -3 + 6 (r - 1)$. Since $\kkl \geq d$, using \eqref{eq:linf_Y} in Lemma \ref{lem:prod}, $ \ang X^{-r+1} \les \cs$ from \eqref{eq:dec_S}, $\brho = \cs^3$ from \eqref{eq:Euler_profi}, and $\| \td G \| _{\cYb^{\kkl}} \leq \ze_0$, we obtain  
\begin{align}
    \label{eqn: smallness-G-L2}
	\| \td G(X, \cdot ) \|_{L ^2 (V)} \les \ang X^{- \f{\etab + d}2} \| G \| _{\cYb^{\kkl}} 
\les  \ang X^{-3(r-1)} 	\| \td G \| _{\cYb^{\kkl}} \les \cs^3 \| \td G \| _{\cYb^{\kkl}}
\leq C_1 \cs^3 \ze_0 =  C_1\ze_0 \brho.
\end{align}
Using Lemma \ref{lem: A-pointwise-bound} and taking $\ze_0$ small enough, we prove
\begin{align*}
    -A_1  \preceq  A [\cMM ^{1/2} \td G]  \preceq  A_1, \quad A_1 =   C \| \td G(X, \cdot ) \|_{L ^2 (V)} \cs^{-3} \bS \preceq C C_1 \ze_0 \bS \preceq \f12 A [\cM].
\end{align*}
which proves the desired estimates.
\end{proof}
\vspace{0.5cm}
Classical parabolic theory then implies that (\ref{eq: linear_viscous}) has a unique smooth solution in $B _{R ^*}$. A more precise statement is summarized in the next theorem:

\begin{lemma}
    \label{lem: existence-FthR}
    Let $\td F_{\init,R}$ be as in \eqref{eq:LWP_init}, and $\td G \in \boldsymbol J _\ze ^k$ with $\ze \leq \ze_0$, where $\ze_0$ is chosen in Lemma \ref{lem:parab}. There exists a unique solution $\td F _{\th, R}$ to (\ref{eq: linear_viscous}) in the space
    \begin{align*}
        \td F _{\th, R} \in C ([0, T]; H _0 ^1 (B _{R ^*}) \cap H ^k (B _{R ^*})) \cap L^2 (0, T; H ^{k+1} (B _{R ^*})) , 
    \end{align*}
    with $\pa _s \td F _{\th, R} \in L ^2 (0, T; H ^{k - 1} (B _{R ^*}))$.

\end{lemma}

\begin{proof}
First note that $\th \Delta_W F$ in bounded domain is equivalent to $ \th \cs^2 \vp_1^2 \ang X \ang \vc ^2  \Delta_{X,V}  F$ (plus first order or $0$ order terms), where $\D_{X, V}F$ is the standard Laplacian in $\R^6$. Since $\ze \leq \ze_0$, from Lemma \ref{lem:parab}, we obtain that $ A[\cM+ \cMM^{1/2}\td G] \succeq \f{1}{2} A [\cM] \succeq 0$ pointwise in $X$ and $V$.
It follows that the equation (\ref{eq: linear_viscous}) is uniformly parabolic in both $X,V$ variables.

Next, we analyze the regularity of the coefficients, starting with the coefficients of the diffusion term. 
By \eqref{eq:M11/2} we know for any $|\al| + |\b| \le k$, 
\begin{align*}
    |D ^{\al, \b} (\cMM ^{1/2} \td G)| \les \cMM ^{1/2} \ang \vc ^{|\b| + 2 |\al|} |D ^{\le k} \td G|.
\end{align*}
By Lemma \ref{lem: A-pointwise-bound}, we know
\begin{align*}
    |D ^{\al, \b} A [\cMM ^{1/2} \td G]| = |A [D ^{\al, \b} (\cMM ^{1/2} \td G)]| \les \| D ^{\le k} \td G \| _{L ^2 (V)} \cs ^{\g + 2} \ang \vc ^{\g + 2},
\end{align*}
which is bounded in $B _{R ^*}$ since $\td G \in \boldsymbol J _\ze ^k \subset L ^\infty (0, T; \cYb ^k)$. As the weight $\ang X ^{|\al|} \cs ^{|\b|}$ in $D ^{\al, \b}$ is bounded from above and below in $B _{R ^*}$, $\pa _{X, V} ^{\al, \b} A [\cMM ^{1/2} \td G]$ is also bounded in $B _{R ^*}$. Therefore, the coefficients of second order derivatives in \eqref{eq: linear_viscous} are in $L ^\infty (0, T; W ^{k, \infty} (B _{R ^*}))$. 

The coefficients in transport terms are the $V$ and $X$ in $\cT$, $\div A [\cM]$ in $\cN _1$, and $\cs$, $\vc$, $A [\cMM ^{1/2} \td G]$, $\div A [\cMM ^{1/2} \td G]$ in $\cN _1$, $\cN _2$, $\cN _3$, $\cN _4$, which are all in $L ^\infty (0, T; W ^{k, \infty} (B _{R ^*}))$. Similarly, coefficients of the reaction terms in $\dcm$, $\bcv$, $\cN _{i}$, are in $L ^\infty (0, T; W ^{k, \infty} (B _{R ^*}))$.

Next, we analyze the regularity of the forcing term $\td H$ defined in \eqref{tildaH}. In particular, by \eqref{eqn: N2346} we have
\begin{align*}
    \| (\cN _2 + \cN _3 + \cN _4 + \cN _6)(\rhos \cMM ^{1/2}, \td G) \| _{H ^k (B _{R ^*})} \le \| D ^{\le k} (\kp _2 \cs ^{-1} (\div A) [\cM \vc]) \| _{L ^\infty (B _{R ^*})} \| D ^{\le k} \td G \| _{L ^2 (B _{R ^*})}.
\end{align*}
Therefore $(\cN _2 + \cN _3 + \cN _4 + \cN _6)(\rhos \cMM ^{1/2}, \td G) \in L ^\infty (0, T; H ^k (B _{R ^*}))$. Note that for functions supported in $B _{R ^*}$, $\sigma$ norm and $H ^1 (V)$ norm are equivalent. By Lemma \ref{lem: N-derivative-bound}, we see for any $|\al| + |\b| \le k$, 
\begin{align*}
    \| D ^{\al, \b} \cN (\td G, \rhos \cMM ^{1/2}) \| _{H ^{-1} (V)} \les _R \| D ^{\le k} \td G \| _{L ^2 (V)} \| D ^{\le k} (\rhos \cMM ^{1/2}) \| _\s,
\end{align*}
which is bounded. Therefore $D ^{\al, \b} \cN (\td G, \rhos \cMM ^{1/2}) \in L ^\infty (0, T; H ^{k - 1} (B _{R ^*}))$.
Finally, from Lemma \ref{lem:ds_t} we know the term $-\cMM ^{-1/2} \eM$ is smooth and bounded in $B _{R ^*}$, and $\sss \cdot \cF_M \circ \cK_k(\hw_1 + \tw_2 - \cF_E (\td G))$ is in $L ^\infty (0, T; H ^{2 k + 6} (B _{R ^*}))$ in view of Proposition \ref{prop:compact}
and Lemma \ref{lem:macro_UPB}.

Summarizing, \eqref{eq: linear_viscous} is a  linear, uniformly parabolic equation with $W ^{k, \infty}$ coefficients and $H ^{k - 1}$ forcing, and therefore has a unique regular solution $F _{\th, R}$ in $(0, T) \times B _{R ^*}$. This concludes the proof of the lemma.   
\end{proof}

\subsubsection{Uniform weighted $L ^2$ estimate}
Before passing to the limit $R \to \infty$ and $\th \to 0$, we need some energy estimates of $\td F _{\th, R}$ uniform in $\th$ and $R$. 

\begin{lemma} \label{lem: Unif_inR_lin}
    Let $\td F _{\th, R}$ be the solution to \eqref{eq: linear_viscous} obtained in Lemma \ref{lem: existence-FthR} and $\ze_0$ be the parameter chosen in Lemma \ref{lem:parab}. Suppose that $\td G \in \boldsymbol J _\ze ^k$ with $\ze < \ze_0$. Then 
    \begin{align} \label{est: Unif_inR_lin}
    &  \f12  \frac{d}{ds} \; \| \td F _{\th, R} \| _{\cYb} ^2 + \f\th2  \| \ang X \ang \vc ^2 \td F _{\th, R} \|^2_ {\cYb ^1} + \f1\es\left( \f12  - \bar C_{\cN,0} \| \td G \| _{ \cYb ^{\kkl}} \right)  \| \td F _{\th, R} \| _{\cYQ} ^2  \\ 
    & \qquad \le 
    C \| \td F _{\th, R} \|^2_ {\cYb} + \f{C}{\es} \| \td G \| _{\cYb} ^2 + \f{C}{\es} \| \td  G\|_{\cYb}^{ \f{4 }{\g+2}} \| \td  G\|_{\cYQ}^{ \f{2 \g }{\g+2}} + C .  \nonumber
\end{align}

\end{lemma}

\begin{proof}    
 We define $\td F _{\th, R} = 0$ outside $B _{R ^*}$. Since $\tFR|_{\pa B_{R^*}} =0$ and $\tFR \in H_0^1(B_{R^*} )$,  this zero extension defines a function on $\R^6$ satisfying $\tFR \in \cYQ(\R^6)$.  
 
 Next, we perform $\cYb$ estimates on \eqref{eq: linear_viscous} by estimating 
$ \iint_{B_{R^*}}  \eqref{eq: linear_viscous} \cdot \tFR \ang X^{\etab} d X d V$
\begin{align}
\label{eq: linear_viscous_weak_0}
& \hspace{-4em}  \iint_{B_{R^*}} \B( \partial_s + \cT - \dcm + \f32 \bcv  \B) \td {F}_{\theta,R} \cdot \tFR  \ang X^{\etab} d X d V  \\ 
\label{eq: linear_viscous_weak_2}
 = &  \f1{\es} \iint _{B _{R^*}} 
 \left[(\cN _1 + \cN _5)(\rhos \cMM ^{1/2}, \td F_{\th, R}) \right] \cdot \tFR \ang X^{\etab} \;dXdV \\
\label{eq: linear_viscous_weak_3}
& +   \f1{\es} \iint_{B_{R^*}}   \cN (\td G, \td F _{\th, R}) \cdot \tFR \ang X^{\etab} \;d X d V \\
\label{eq: linear_viscous_weak_4}
& +  \iint  _{B _{R^*}} \td H \cdot \tFR \ang X^{\etab} \;dXdV \\
\label{eq: linear_viscous_weak_5}
& +  \iint_{B_{R^*}} \th \Delta_W \td F _{\th, R} \cdot \tFR \ang X^{\etab} \; d X d V .
\end{align}
Recall that $\td H$ is defined in \eqref{tildaH}.

For \eqref{eq: linear_viscous_weak_0}, %
we apply Corollary \ref{cor: scaling-field-coersive} with $k = 0$, 
$l = \etab$, and $\d$ large enough and get 
\begin{align*}
    \f12 \f d{ds} \| \td F _{\th, R} \| _{\cYb} ^2 \le %
     \iint_{B_{R ^*}} \B( \partial_s + \cT - \dcm + \f32 \bcv  \B) \td {F}_{\theta,R} \cdot \tFR  \ang X^{\etab} d X d V
    + C \| \td F _{\th, R} \| _{\cYb} ^2  + \f{1}{8 }  \| \td F _{\th, R} \| _{\cYQ} ^2 . 
\end{align*}

For \eqref{eq: linear_viscous_weak_2}, using the regularity of $\tFR$ from Lemma \ref{lem: existence-FthR} and $\tFR \in \cYQ(\R^6)$, we apply \eqref{eqn: N1+N5=sigma}:
\begin{align*}
    \eqref{eq: linear_viscous_weak_2} %
    &= - \f1\es \int \la X \ra^{\etab} \| \td F _{\th, R}\| _\sigma ^2 d X  = - \f1\es \| \td F _{\th, R} \| _{\cYQ} ^2 .
\end{align*}

For \eqref{eq: linear_viscous_weak_3} we use Lemma \ref{lem: nonlinear} and Sobolev embedding \eqref{eqn: smallness-G-L2} to get 
\begin{align*}
\eqref{eq: linear_viscous_weak_3} &=   \f1{\es} \iint_{B_{R^*}}   \cN (\td G, \td F _{\th, R}) \cdot \tFR \ang X^{\etab} \;d X d V \\ 
& \le C  \sup _{ X} \left\{\cs ^{-3} \| \td G (s, X, \cdot) \| _{L ^2 (V)}\right\}  \f1\es \int \ang X ^{\etab} \| \td F _{\th, R}\| _\sigma ^2 d X  \\
& \leq  \f{\bar C_{\cN,0}}{\es} \| \td G \| _{ \cYb ^{\kkl}}   \| \td F _{\th, R} \|^2 _{\cYQ } .
\end{align*}
In the proof of Lemma \ref{lem: nonlinear}, we only use integration by parts once. Since $\tFR |_{\pa B_R^*}=0$, we obtain the same proof and estimates.

For the diffusion term \eqref{eq: linear_viscous_weak_5}, applying Lemma \ref{lem: diffusion-innerprod}, we have 
\begin{align*}
    \eqref{eq: linear_viscous_weak_5} %
    = \th \la \D _W \td F _{\th, R}, \td F _{\th, R} \ra _{\cYb}  = -\th \| \ang X \ang \vc ^2 \td F _{\th, R} \| _{\cYb ^1} ^2.
\end{align*}

We now estimate \eqref{eq: linear_viscous_weak_4}: 
\bseq
\begin{align}
    \notag
    \eqref{eq: linear_viscous_weak_4} &= \iint \la X \ra^{\etab} \td H \td F _{\th, R} \;dXdV \\
    \label{eq: linear_viscous_weak_4a}
    \tag{\ref*{eq: linear_viscous_weak_4}a}
    & = \f1\es \iint \kp _2 \cs ^{-1} \la X \ra^{\etab} (\div A) [\cM \vc] \td G \td F _{\th, R} \;dXdV \\
    \label{eq: linear_viscous_weak_4b}
    \tag{\ref*{eq: linear_viscous_weak_4}b}
    & \qquad + \f1\es \iint \la X \ra^{\etab} \cN (\td G, \rhos \cMM ^{1/2})\td F _{\th, R} \;dXdV \\
    \label{eq: linear_viscous_weak_4c}
    \tag{\ref*{eq: linear_viscous_weak_4}c}
    & \qquad - \iint \la X \ra^{\etab} \cMM ^{-1/2} \eM \td F _{\th, R} \;dXdV \\
    \label{eq: linear_viscous_weak_4d}
    \tag{\ref*{eq: linear_viscous_weak_4}d}
    & \qquad + \iint  \la X \ra^{\etab} \sss \cdot \cF_M \circ \cK_k(\hw_1 + \tw_2 - \cF_E (\td G)) \td F _{\th, R} \;dXdV.
\end{align}
\eseq

We start with the first integral. Recall $\cMM = \cs^{-3} \mu(\vc)$ and $\cM = \mu(\vc) $ from \eqref{eq:localmax2}.
Using \eqref{eqn: divA-pointwise} with $i=1, j=0$, $ \cM \vc = \cMM^{1/2} \cdot \cs^{3/2} \mu(\vc)^{1/2}\vc$, and \eqref{eq:cross_pf2:a}, we obtain
\beq\label{eq:LWP_divA_point}
 |(\div A) [\cM \vc] | 
 \les \cs^{\g+1} \ang \vc^{\g+1} \| \cs^{3/2} \mu(\vc)^{1/2} \vc \|_{L^2(V)}
 \les \cs^{\g+4} \ang \vc^{\g+1} .
\eeq
Using \eqref{eq:LWP_divA_point}, we estimate the first integral as 
\begin{align*}
     \eqref{eq: linear_viscous_weak_4a} &= \f1\es \iint \kp _2 \cs ^{-1} \la X \ra^{\etab} (\div A) [\cM \vc] \td G \td F _{\th, R} \;dXdV \\
     & \leq \f{C}{\es} \iint  \ang X ^{\etab} \cs ^{\g + 3} \ang \vc ^{\g + 1} |\td G \td F _{\th, R}| \;dXdV
     \\
     & \le  \udb{\f{C}{\es} \iint \ang X ^{\etab} \cs ^{\g + 3} \ang \vc ^\g |\td G| ^2 dV dX }_{:=I} + \f1{8 \es} \iint \ang X ^{\etab} \cs ^{\g + 3} \ang \vc ^{\g + 2} |\td F _{\th, R} | ^2 dV dX .
     \end{align*}

Since $\cs \les 1$ and $\g \geq 0$, using H\"older's inequality, we bound 
\beq\label{eq:linear_L2_G_hol}
I \les  \f{1}{\es} \B( \iint \ang X^{\etab} \cs^{\g+3} |\td G|^2  \B)^{ \f{ 4 }{\g+2} } 
\B( \iint \ang X^{\etab}  \cs^{\g+3} \ang \vc^{\g+2} |\td G|^2  \B)^{ \f{  2 \g}{\g+2} } 
\les \f{1}{\es} \| \td  G\|_{\cYb}^{ \f{4 }{\g+2}} \| \td  G\|_{\cYQ}^{ \f{2 \g }{\g+2}} . 
\eeq

Summarizing we obtain: 
\begin{align*}
    \eqref{eq: linear_viscous_weak_4a} \le  \f{C}{\es} \| \td  G\|_{\cYb}^{ \f{4 }{\g+2}} \| \td  G\|_{\cYQ}^{ \f{2 \g }{\g+2}} 
    + \f1{ 8 \es} \| \td F _{\th, R} \| _{\cYQ} ^2.
\end{align*}

For the second integral: using $\| \cMM ^{1/2} \| _\sigma \les \cs ^{\f{3 + \g}2} \les 1$ and $\brho = \cs^3$ \eqref{eq:Euler_profi_modi}, we have 
\begin{align*}
\eqref{eq: linear_viscous_weak_4b} &= \f1\es  \int \la X \ra^{\etab} \la \cN (\td G, \rhos \cMM ^{1/2}), \td F _{\th, R} \ra _V \;dX \\
& \le \f1\es  \int  \cs ^{-3} \ang X ^{\etab} \| \td G (s, X, \cdot) \| _{L ^2 (V)} \| \td F _{\th, R}\| _\sigma \| \rhos \cMM ^{1/2} \| _\sigma  d X\\
& \les \f1\es  \int \ang X ^{\etab} \| \td G (s, X, \cdot) \| _{L ^2 (V)} \| \td F _{\th, R}\| _\sigma d X\\
& \le \f1{ 8 \es}  \int \ang X ^{\etab} \| \td F _{\th, R}\|^2 _\sigma  + \f C\es \int \ang X ^{\etab} \| \td G (s, X, \cdot) \|^2 _{L ^2 (V)} d X \\
& \le \f1{ 8 \es} \| F _{\th, R}\|^2 _{\cYQ} + 
\f{C}{\es} \| \td G \| _{\cYb} ^2 . 
\end{align*}

For the third integral, recall from Lemma \ref{lem:ds_t} we know
   $ |\eM| \les \cM \ang X ^{-r} \ang \vc ^3.$
Since $\cM = \cs^3 \cMM$ \eqref{eq:localmax2}, we obtain 
\begin{align*}
    \| \cMM ^{-1/2} \eM \| _{\cYb} ^2 \les
    \iint \cs ^6 \cMM \ang X ^{-2 r} \ang \vc ^6 \ang X ^{\etab} d V d X \les \int _{\R ^3} \ang X ^{-2 r + \etab} d X \le C.
\end{align*}
Here we used $- 2 r + \etab < -3$. Hence, using Cauchy--Schwarz inequality, we have 
\begin{align*}
    \eqref{eq: linear_viscous_weak_4c} &=  \iint \cMM ^{-1/2} \eM \td F _{\th, R} \ang X ^{\etab} \;dX dV  \\ 
    & \les \| \cMM^{-1/2}  \eM \|_{\cYb}^2  + \| \tFR \|_{\cYb}^2 
\leq   \| \td F  _{\th, R} \| _{\cYb} ^2  + C . 
\end{align*}

For the last term, using Cauchy--Schwarz inequality, we have: 
\begin{align*}
     \eqref{eq: linear_viscous_weak_4d} &=  \iint  \la X \ra^{\etab} \sss \cdot \cF_M \circ \cK_k(\hw_1 + \tw_2 - \cF_E (\td G)) \td F _{\th, R} \;dXdV \\
     & \le  \| \td F  _{\th, R} \| _{\cYb} ^2  + \sss  \| \cF_M \circ \cK_k(\hw_1 + \tw_2 - \cF_E (\td G)) \| _{\cYb} ^2 .
\end{align*}

Recall that $\cK _k: \cX _\xeta \to \cX _\xeta$ is defined in Proposition \ref{prop:compact}, with parameter $\xeta$. Moreover, its image is compactly supported in $B _{4 R_{\xeta}}$ which depends only on $\xeta$. 
Applying Lemma \ref{lem:macro_UPB}, \eqref{eq:LWP_W1W2}, %
we obtain
\begin{align*}
    \| \cF_M \circ \cK_k(\hw_1 + \tw_2 - \cF_E (\td G)) \| _{\cYb} &\les \| \cK_k(\hw_1 + \tw_2 - \cF_E (\td G)) \| _{\cX _\etab} \\
    &\les _\xeta \| \cK_k(\hw_1 + \tw_2 - \cF_E (\td G)) \| _{\cX _\xeta} \\ 
    &\les _\xeta \| \hw_1 + \tw_2 - \cF_E (\td G) \| _{\cX _\xeta} \\
    &\les _\xeta \| \hw_1 + \tw_2 \| _{\cXX} + \| \td G \| _{\cYY} \\
    & \les_{\xeta} 1 +  \| \td G \| _{\cYb} . 
\end{align*}

For $\sss \in \{0, 1\}$, summarizing we get 
\begin{align*}
    \underbrace{ \f12 \frac{d}{ds}  \;    \| \td F _{\th, R} \| _{\cYb} ^2 }_{ %
    \eqref{eq: linear_viscous_weak_0}} &\le \underbrace{ C \| \td F _{\th, R} \| _{\cYb} ^2  + \f{1}{8}  \| \td F _{\th, R} \| _{\cYQ} ^2 } _{%
    \eqref{eq: linear_viscous_weak_0}} - \underbrace{ \f1\es \| \td F _{\th, R} \| _{\cYQ} ^2 }_{\eqref{eq: linear_viscous_weak_2}}  
    + \underbrace{\f C\es \| \td G \| _{ \cYb ^{\kkl}} \| \td F _{\th, R} \|^2 _{\cYQ}\;}_{\eqref{eq: linear_viscous_weak_3}} \\  %
    & \qquad -\underbrace{\f\th2  \| \ang X \ang \vc ^2 \td F _{\th, R} \|^2_ {\cYb ^1}\;}_{\eqref{eq: linear_viscous_weak_5}}   %
    + \underbrace{\f{C}{\es} \| \td  G\|_{\cYb}^{ \f{4 }{\g+2}} \| \td  G\|_{\cYQ}^{ \f{2 \g }{\g+2}} + \f1{8 \es} \|\td F _{\th, R} \|_{\cYQ} ^2 }_{\eqref{eq: linear_viscous_weak_4a}} \\
    & \qquad 
    + \underbrace{\f1{8 \es} \| \td F _{\th, R}\|^2 _{\cYQ}  + \f{C}{\es} \| \td G \| _{\cYb} ^2} _{\eqref{eq: linear_viscous_weak_4b}} 
    + \underbrace{ \| \td F  _{\th, R} \| _{\cYb} ^2  + C_0}_{\eqref{eq: linear_viscous_weak_4c}} 
    + \underbrace{\| \td F _{\th, R} \| _{\cYb} ^2 + 
    C ( 1
    + \| \td G \| _{\cYb} ^2)
    }_{\eqref{eq: linear_viscous_weak_4d}}.
\end{align*}

Since $\es \leq 1$ \eqref{eq:EE_para1}, after reorganization we get 
\begin{align*}
    & \f12  \frac{d}{ds} \;  \| \td F _{\th, R} \| _{\cYb} ^2 + \f\th2  \| \ang X \ang \vc ^2 \td F _{\th, R} \|^2_ {\cYb ^1} + \f1\es\left( \f12  - C \| \td G \| _{ \cYb ^{\kkl}} \right)  \| \td F _{\th, R} \| _{\cYQ} ^2  \\ 
    & \qquad \le 
       C \| \td F _{\th, R} \|^2_ {\cYb} + \f{C}{\es} \| \td G \| _{\cYb} ^2 + \f{C}{\es} \| \td  G\|_{\cYb}^{ \f{4 }{\g+2}} \| \td  G\|_{\cYQ}^{ \f{2 \g }{\g+2}} + C,
\end{align*}
where $C$ is some absolute constant. 

This completes the $L^2$ energy estimate.
\end{proof}

\subsubsection{Uniform weighted $H ^j$ estimate}

Next, we derive interior $\cYb ^j$ estimates for $\td F_{\theta, R}$ uniformly in $R$, with $1 \le j \le k$. 
We introduce a sequence of cutoff functions as follows. We define 
\beq\label{eq:LWP_phi_j}
\phi _0 := 1, \ \td \phi _0 = 0,   \quad 
\phi _j := \varphi _{2 ^{-j} R}, \ \td \phi _j = \td \varphi _{2 ^{-j} R} , \  \forall j > 0, 
\eeq 
where $\varphi _R$ is defined in \eqref{eqn: defn-varphi} and $\td \varphi _R$ is defined in \eqref{eqn: defn-tdvarphi}. Consequently, $\phi _j$ for $j \ge 1$ are all supported in $B _{R ^*}$. To simplify notation, we omit the dependence of $\phi_j$ on $R$.
We define 
\begin{align}\label{eq:LWP_hj}
    h _j := \td F _{\th, R} \cdot \phi _j.    
\end{align}
Note that $h _0 = \td F _{\th, R}$. Since $\phi _j = 1$ over the support of $\phi _i$ for all $i <  j$, we can write $h _j = h _{j - 1} \cdot \phi _j$. 

Applying   \eqref{eqn: uniform-bound-cutoff} with $( \vp_R, \td \vp_R, R) \rsa ( \phi_j = \vp_{2^{-j}R}, \td \phi_j = \td \vp_{2^{-j} R}, 2^{-j} R)$ we obtain 
\beq\label{eq:uniform-cutoff-phij}
    |D ^{\al, \b} \phi _j| \les_{\al, \b} 1, \quad   |D ^{\al, \b} \td \phi _j| \les _{\al, \b}  ( 2^{-j} R)^{-2} 
    \les_{\al, \b, j} R^{-2}.
\eeq

In the next lemma we calculate the equation satisfied by $h _j$. 

\begin{lemma}
    For $j \ge 1$, $h _j$ solves the following equation in $\R ^6$: 
    \begin{align}\label{eq: h_linear_viscous}
        \left(\pa _s + \cT + \dcm - \f32 \bcv\right) h _j & = \f1{\es} \left[(\cN _1 + \cN _5) (\rhos \cMM ^{1/2}, h _j) \right] + \f1{\es} \cN (\td G, h _j) + \td H _{j} + \th \D _W h _j,
    \end{align}
    where 
    \beq\bal
    \label{eq: h_linear_viscous_2}
    \td H _j &:= \phi _j \td H + \left(\pa _s + \cT \right) \phi _j \cdot h _{j - 1} - \th [\D _W, \phi _j ] h _{j - 1}
    \\
    & \qquad - \f1\es \left( \cN _1 (\rhos \cMM ^{1/2}, \phi _j) h _{j - 1} + \cN _1 (\td G, \phi _j) h _{j - 1} \right) \\
    & \qquad + \f1\es \td \phi _j \left( 2 \cN _4 (\rhos \cMM ^{1/2}, h _{j - 1}) - (\cN _6 - 2 \cN _4 - 2 \cN _5) (\td G, h _{j - 1}) \right).
\eal\eeq
\end{lemma}

\begin{proof}
We omit the subscript $j$. First, we apply Leibniz rule and $\tFR = h_{j-1}$ on the support of $\phi_j$:
\begin{align*}
    \left(\pa _s + \cT \right) h &= \phi \left(\pa _s + \cT \right) \td F _{\th, R} + \left(\pa _s + \cT \right) \phi \cdot \td F _{\th, R}, \\
    \D _W h &= \D _W (\td F _{\th, R} \phi) = \phi \D _W \td F _{\th, R} + [ \D _W, \phi] \td F _{\th, R}
    =  \phi \D _W \td F _{\th, R} + [ \D _W, \phi] h_{j-1} .
\end{align*}
For the collision terms, recall the definition of $\cN _i$, $1 \le i \le 6$, in \eqref{N(f,g)}. By Leibniz rule and \eqref{eqn: nabla_V_varphi}, we have 
\begin{align*}
    \cN _1 (\td G, \td F _{\th, R} \phi) &= \phi \cN _1 (\td G, \td F _{\th, R}) + \cN _1 (\td G, \phi) \td F _{\th, R} + 2 A [\cMM ^{1/2} \td G] \na _V \phi \cdot \na \td F _{\th, R}, \\
    &= \phi \cN _1 (\td G, \td F _{\th, R}) + \cN _1 (\td G, \phi) \td F _{\th, R} - 2 \td \phi \cN _4 (\td G, \td F _{\th, R}), \\
    \cN _2 (\td G, \td F _{\th, R} \phi) &= \phi \cN _2 (\td G, \td F _{\th, R}) - \na _V \phi \cdot \div A [\cMM ^{1/2} \td G] \cdot \td F _{\th, R}, \\
    &= \phi \cN _2 (\td G, \td F _{\th, R}) + \td \phi \cN _6 (\td G, \td F _{\th, R}), \\
    \cN _{3, 4} (\td G, \td F _{\th, R} \phi) &= \phi \cN _{3, 4} (\td G, \td F _{\th, R}) - A [\cMM ^{1/2} \td G] \na _V \phi \cdot \kp _2 \cs ^{-1} \vc \cdot \td F _{\th, R}, \\
    &= \phi \cN _{3, 4} (\td G, \td F _{\th, R}) - \td \phi \cN _5 (\td G, \td F _{\th, R}), \\
    \cN _{5, 6} (\td G, \td F _{\th, R} \phi) &= \phi \cN _{5, 6} (\td G, \td F _{\th, R}).
\end{align*}
Therefore, 
\begin{align*}
    \cN (\td G, h) = \phi \cN (\td G, \td F _{\th, R}) + \cN _1 (\td G, \phi) \td F _{\th, R} + \td \phi (\cN _6 - 2 \cN _4 - 2 \cN _5) (\td G, \td F _{\th, R}).
\end{align*}
Similarly, 
\begin{align*}
    (\cN _1 + \cN _5) (\rhos \cMM ^{1/2}, h) = \phi (\cN _1 + \cN _5) (\rhos \cMM ^{1/2}, \td F _{\th, R}) + \cN _1 (\rhos \cMM ^{1/2}, \phi) \td F _{\th, R} - 2 \td \phi \cN _4 (\rhos \cMM ^{1/2}, \td F _{\th, R}). 
\end{align*}
We thus obtain the equation for $h$ in $\R ^6$, from multiplying \eqref{eq: linear_viscous} by $\phi$:
\begin{align*}
    \left(\pa _s + \cT + \dcm - \f32 \bcv\right) h &= \f1{\es} \left[(\cN _1 + \cN _5) (\rhos \cMM ^{1/2}, h) \right] + \f1{\es} \cN (\td G, h) + \td H _j + \Delta_W h
\end{align*}
with a new forcing term $\td H _j$ in \eqref{eq: h_linear_viscous_2}.
\end{proof}

Recall $h_j$ from \eqref{eq:LWP_hj}. Since $  h_j =  \phi_j \phi_{i} \tFR = \phi_j h_i$ for any $ i < j$, a straightforward consequence of 
    \eqref{eqn: uniform-bound-cutoff} is that by Leibniz rule, 
\begin{align}\label{eq:Hj_chain}
    \| h _j \| _{\cY _l ^k} \les_{k}      \| h _i \| _{\cY _l ^k} \les_{  k} \|  \td F _{\th, R} \| _{\cY _l ^k},  \quad 
    \| h _j \| _{\cYE ^k} \les_{ k}   \| h _i \| _{\cYE ^k} \les_{ k} \| \td F _{\th, R} \| _{\cYE ^k} ,
\end{align}
with constants \emph{only} depending on $k$. This can be directly verified by Leibniz rule so we do not go into details. 

For any $i, j \in  \R$, since the commutator between the weight  $ \ang X^i \ang \vc^j $ and 
derivatives $D^{\al, \b}$ consists of terms with $\leq |\al| + |\b|-1$ derivatives,  
using induction on $k$, 
and $|D^{\leq l} \ang X^i \ang \vc^j | \les_{i, j, l} \ang X^i \ang \vc^j $, we have 
\bseq
\label{eq:norm_diffu_equiv}
\beq
| D^{\leq k} (\ang X^i \ang \vc^j f) | \asymp_{i, j, k}  \ang X^i \ang \vc^j | D^{\leq k} f  |,
\quad \forall \,  k \geq 0, i, j \in \R. 
\eeq
Recall $\Lams$ from \eqref{eq:coer_coe}. Similarly, we have
\beq 
|D ^{\le k} (\Lams ^\f12 f)| \asymp _k \Lams ^\f12 |D ^{\le k} f|,
\quad \Lams = \cs^{\g+3} \ang \vc^{\g+2}.
\eeq
\eseq
We omit the proof.

To bound the collision commutators in $\td H _j$ \eqref{eq: h_linear_viscous_2}, we need the following estimates. 
\begin{lemma}\label{lem:Q_commu}
For any $ i \leq j$, we have 
\[
\bal 
    I _{i, j} &:= \int \ang X ^{\etab} \cs ^\g \ang \vc ^{\g + 2} \cdot \|D ^{\le i} \td G \| _{L ^2 (V)} 
    \cdot |D ^{\le j - i} f| \cdot |D ^{\le j} g| d V d X \\
    &\les _j  \left( 
    \one_{i>0} \| \td G \|_{\cYb^{\max \{\kkl, i \}} }  \| f\|_{\cYQ^{j-1}}
        + \| \td G\|_{\cYb^{\kkl}}  \| \Lam ^\f12 f \|_{\cYb^j} 
    \right) \| g\|_{\cYQ ^j}.
\eal 
\]
\end{lemma}

\begin{proof}
Recall the $\s$-norm from \eqref{eqn: sigma-norm-approx}. We have 
\begin{align*}
    \int \cs ^\g \ang \vc ^{\g + 2} |D^{\le j - i} f| \cdot |D ^{\le j} g| d V \les \cs ^{-3} \| D ^{\le j - i} f \| _\s \| D ^{\le j} g \| _\s.
\end{align*}
For $i \geq 1$, $I _{i, j}$ has the same structure as the $II$ term in \eqref{eq:non_pf_decomp-II}. Using \eqref{eq:non_pf_LH}-\eqref{eq:non_pf_HL-2}, we estimate:
\begin{align*}
    I _{i, j} \les_j \| \td G \|_{\cYb^{\max{\{\kkl, i\}}}}  \| f\|_{\cYQ^{j-1}} \| g\|_{\cYQ^j} .
\end{align*}

Note that $\cs \les 1$, so instead of $\s$ norm we can also bounded $I _{i,j}$ by weighted diffusion
\begin{align*}
    \int \cs ^\g \ang \vc ^{\g + 2} |D^{\le j - i} f| \cdot |D ^{\le j} g| d V \les \cs ^{-3} \| \Lam ^\f12 D ^{\le j - i} f \| _{L ^2 (V)} \| D ^{\le j} g \| _\s.
\end{align*}
For $i = 0$, $I _{i, j}$ has the same structure as the $I$ term in \eqref{eq:non_pf_decomp-I}, which is bounded using \eqref{eq:non_pf_I} 
and then \eqref{eq:norm_diffu_equiv} by 
\begin{align*}
    I _{0, j} \les \| \td G \|_{\cYb ^\kkl} \| \Lam ^\f12 f\|_{\cYb ^{j}} \| g\|_{\cYQ^j}.
\end{align*}
We complete the proof.
\end{proof}

\begin{lemma}\label{lem:LWP_Hj}
Let $ \td G \in \JJ_{\ze}^{k}$ with $k > \kkl$ and $\ze < 1$. For any $ j \geq 1$, we have the following energy estimates for $h_j = \phi_j \tFR$
\beq\label{eq:LWP_EE_Hj}
\bal 
  \f12 \f{d}{ds} \| h _j \| _{\cYb ^j} ^2 
& \le  C_j  \| h _j \| _{\cYb ^j} ^2 
-  \f{1}{\es} \big( \f{ 1}{2} - \cnon \| \td G \|_{\cYb^{\kkl}} \big) \| h _j \| _{\cYQ ^j} ^2 + \one_{j>0} \f{  \bar C_{j,1}}{\es}   \| h _{j - 1} \| _{\cYQ ^{j - 1}} ^2 \\
& \quad +   \one_{j > \kkl} \f{C_j}{\es}   \| \td G \|_{ \cYb^{ j }}  
               \| h_{j-1} \| _{\cYQ ^{j-1}}   \| h _j \| _{\cYQ ^j}  
            +  \f{C_j}{\es} \big( 
      \| \td G \| _{\cYb ^j}^{ \f{4}{\g+2}} \| \td G \| _{\cYQ ^j}^{ \f{2\g}{\g+2}} 
    +\| \td G \| _{\cYb ^j} ^2 \big) +  C _j     \\
& \quad  - \f\theta 4 \|\ang X \ang \vc ^2 h_j \| _{\cYb ^{j + 1}} ^2   + \one_{j>0} \bar C_{j, 2} \left(  \theta + \f{1}{R^2 \es} \right)  \|\ang X \ang \vc ^2 h_{j-1} \| _{\cYb ^{j }} ^2  
\eal 
\eeq
for some absolute constant $\bar C_{j, i}$ with $\bar C_{0, i} = 0, i=1,2$, where $\cnon$ is an absolute constant determined in Theorem \ref{theo: N nonlin_est}. 

Suppose that $ R^{-2} \leq \th$. We define the $R$-dependent norms $\cZZ_R^j,  \cZZQR^j$: \footnote{
Note that these parameters $\varpi_{Z,j}$ are different from those in \eqref{norm:Xk}.
}
\beq\label{norm:cZZ}
\bga
 \| f\|_{{\cZZ_R^j}}^2 :=   \sum_{ 0\leq i\leq j} \varpi_{Z,i}  \| \phi_i f\|^2_{\cYb^i} , \quad 
  \| f \|_{ \cZZQR^j}^2 :=  \sum_{0\leq  i\leq j} \varpi_{Z,i} \| \phi_i f\|_{ \cYQ^i}^2, \\ 
  \varpi_{Z, 0} = 1, \quad  \varpi_{Z,j} = \prod_{i \leq j} \f{1}{16 (1+ \bar C_{j,1} + \bar C_{j, 2} )} ,
 \ega 
\eeq
which depend on $R$ via the cutoff functions $\phi_i$ defined in \eqref{eq:LWP_phi_j}.

Then, we have the following estimates uniformly in $R, \theta$ that satisfy $ R \cdot \es \geq 1,  R^{-1} \leq \th$: 
\footnote{
To construct a local solution in the time-interval $s \in [ 0,  T]$ for some finite $T$, 
we choose $R$ large enough and then $\th$ small enough so that  the  assumptions on $R, \th$ are satisfied.
}
\begin{align}\label{eq:LWP_EE_Hj2}
	\f{1}{2} \f{d}{d s} \| \tFR \|_{\cZZ_R^j}^2 
	& \leq C_j \|  \tFR \|_{ \cZZ_R^j }^2
	-  \f{ 1}{\es}  (  \f{3 }{8} - \cnon \| \td G\|_{ \cYb^\kkl }  ) \| \tFR \|_{ \cZZQR^j}^2  \\
    & \quad +\one_{j>\kkl} \f{  C_j}{\es} \| \td G\|_{ \cYb^{ j}}     
	\| \tFR \|_{ \cZZQR^{  j-1 }}  \| \tFR \|_{ \cZZQR^j} 
            +   \f{C_j}{\es} \| \td  G\|_{\cYb^j}^{ \f{4 }{\g+2}} \| \td  G\|_{\cYQ^j}^{ \f{2 \g }{\g+2}}
     +  \f{C_j}{\es} \| \td G\|_{\cYb^j}^2 + C_j , \notag   
\end{align}
 for any $j \geq 0$, where the constants $C_j$ may change from line to line.

\end{lemma}

\begin{proof}
We recall that 
$$
\la f_1 , f_2 \ra _{\cY _{\etab} ^k} = \sum _{|\al| + |\b| \le k} \nu ^{|\al| + |\b| - k} 
\f{|\al|!}{\al!}
\int \ang X ^\etab \la D ^{\al, \b} f_1, D ^{\al, \b} f_2 \ra _{L ^2 (V)} d X.
$$
We perform $\cYb^j$ estimates on \eqref{eq: h_linear_viscous} by estimating 
$\la \textrm{\ref{eq: h_linear_viscous}} , h_j \ra _{\cYb ^j}$: %
\begin{align}
    \label{eq: h_linear_viscous_III}
    \left\la \left(\pa _s + \cT + \dcm - \f32 \bcv\right) h _j, h_j \right\ra _{\cYb ^j} 
    & = \f1{\es} \left\la (\cN _1 + \cN _5)(\rhos \cMM ^{1/2}, h _j), h_j \right\ra _{\cYb ^j} \\
    & \qquad + \f1{\es} \left \la \cN (\td G, h _j), h_j \right\ra _{\cYb ^j} + \la \td H _{j},  h_j\ra _{\cYb ^j} + \th \la \D _W h_j, h_j \ra _{\cYb ^j}. \notag
\end{align}

\paragraph{Proof of \eqref{eq:LWP_EE_Hj}} 

We analyze inner products in \eqref{eq: h_linear_viscous_III} term by term. This will be analogous to the $L ^2$ estimate in Lemma \ref{lem: Unif_inR_lin}.

\begin{itemize}
\setlength{\leftskip}{-2em}
\item \textit{Viscosity}. Apply Lemma \ref{lem: diffusion-innerprod} to $h \rsa h_j$, we have
\beq%
    - \th \la  \D _W h_j, h_j \ra _{\cYb ^j} \ge \f \theta 2 \| \ang X \ang \vc ^2 h_j \| _{\cYb ^{j + 1}} ^2
    - C_j \theta  \one_{j>0} \| \ang X \ang \vc ^2 h_j \| _{\cYb ^{j }} ^2 , \quad \forall j \geq 0 . \nonumber
\eeq 
Note that for $j=0$ we have $h_0 = F_{\theta, R}$.
Using $ h_j = h_{j-1} \phi_j$ for $j \ge 1$, and \eqref{eq:Hj_chain}, we obtain
\begin{align}\label{eq:LWP_Lap1}
	     -\th \la \D _W h_j, h_j \ra _{\cYb ^j} 
	    \geq \f\theta 2 \|\ang X \ang \vc ^2 h_j \| _{\cYb ^{j + 1}} ^2
	    - C_j \th \one_{j>0} \|\ang X \ang \vc ^2 h_{j-1} \| _{\cYb ^{j }} ^2 .
\end{align}

\item \textit{Transport}.
Using Corollary \ref{cor: scaling-field-coersive} for the left hand side of \eqref{eq: h_linear_viscous_III} with $\d$ large enough, we obtain
\begin{align} \label{trans_lin}
    \f12 \f{d}{ds} \| h _j \| _{\cYb ^j} ^2 \le \left\langle \left(\pa _s + \cT + \dcm - \f32 \bcv\right) h _j, h _j \right\rangle _{\cYb ^j} + C_j \| h _j \| _{\cYb ^j} ^2 + \f{1}{16} \| h _j \| _{\cYQ ^j} ^2.
\end{align}

\item \textit{Main collision}. Now we handle the first term on the right hand side of \eqref{eq: h_linear_viscous_III}.
For $\cN _1, \cN _5$, we first apply Lemma \ref{lem: N-derivative-bound} to get the lower order terms: for $|\al| + |\b| \le j$, 
\begin{align*}
    & \left\la D ^{\al, \b} \cN _i (\rhos \cMM ^{1/2}, h _j) - \cN _i (\rhos \cMM ^{1/2}, D ^{\al, \b} h _j), D ^{\al, \b} h _j \right\ra _{\cYb}\nonumber \\
    & \le \cs ^{-3} \|D ^{\le j} (\rhos \cMM) ^{1/2}\| _{L ^2 (V)}
    \| D ^{\prec (\al, \b)} h _j \| _{\cYQ} \| D ^{\al, \b} h _j \| _{\cYQ} \nonumber\\
    & \le \f{1}{8} \| D ^{\al, \b} h _j \| _{\cYQ} ^2 + C _{j} \| D ^{\prec (\al, \b)} h _j \| _{\cYQ} ^2. %
\end{align*}
The leading order term, using \eqref{eqn: N1+N5=sigma} with $\td F = D ^{\al, \b} h _j$, reads
\begin{align*}
    \left\la (\cN _1 + \cN _5)(\rhos \cMM ^{1/2}, D ^{\al, \b} h _j), D ^{\al, \b} h _j \right\ra _{\cYb} = -\| D ^{\al, \b} h _j \| _{\cYQ} ^2.
\end{align*}
Combined and taking summation over $\al, \b$, and then using \eqref{eq:Hj_chain}, for $j\geq1$, we conclude 
\beq  \label{main_coll_main}
\bal 
    \left\la (\cN _1 + \cN _5)(\rhos \cMM ^{1/2}, h _j) , h_j \right\ra _{\cYb ^j} & \le -\f78 \| h _j \| _{\cYQ ^j} ^2 + C _j \one_{j>0} \| h _j \| _{\cYQ ^{j - 1}} ^2 \\
  &  \leq  -\f78 \| h _j \| _{\cYQ ^j} ^2 + C _j \one_{j>0} \| h _{j-1} \| _{\cYQ ^{j - 1}} ^2
\eal 
\eeq 
When $j =0$, we obtain the above estimate from \eqref{eqn: N1+N5=sigma}, and we do not have the lower order term.

\item \textit{Secondary collision}. We handle the second term on the right hand side of \eqref{eq: h_linear_viscous_III}.
Using estimate \eqref{eq:non_Q:micro:a}, \eqref{eq:non_Q:micro:b} in Theorem \ref{theo: N nonlin_est} with $l _1 = l _2 = l = \etab$, we get 
\begin{align*}
        \B|\B\la \cN (\td G, h _j), h_j \B\ra _{\cYb ^j} \B| \leq 
        \B( \cnon \| \td G \| _{ \cYb^{\kkl} } 
        \| h _j \| _{\cYQ ^j} 
        +  C_j \one_{j > \kkl} \| \td G \|_{ \cYb^{ j }}  
               \| h _j \| _{\cYQ ^{j-1}}  \B) \| h _j \| _{\cYQ ^j}  .
\end{align*}
For $j \leq \kkl$, we only need the first term on the right hand side to bound the nonlinear terms. %
Again, using \eqref{eq:Hj_chain} and $j> j-1$, we further bound the nonlinear terms as 
\begin{align}\label{sec_coll_lin}
	  \B|  \B\la \cN (\td G, h _j), h_j \B\ra _{\cYb ^j} \B| \leq 
        \B( \cnon \| \td G \|_{ \cYb^{\kkl} } 
        \| h _j \| _{\cYQ ^j} 
              +  C_j \one_{j > \kkl} \| \td G \|_{ \cYb^{ j }}  
               \| h_{j-1} \| _{\cYQ ^{j-1}}  \B) \| h _j \| _{\cYQ ^j}  .
\end{align}

\item \textit{Forcing}. Recall the forcing term $\td H_j$ from \eqref{eq: h_linear_viscous_2}
\beq\label{eq:Hj_recall}
\bal
    \td H _j &:= \phi _j \td H + \left(\pa _s + \cT \right) \phi _j \cdot h _{j - 1} - \th [\D _W, \phi _j] h _{j - 1} \\
    & \qquad - \es^{-1}\cN _1 (\rhos \cMM ^{1/2}, \phi _j) h _{j - 1} + 2  \es^{-1} \td \phi _j \cN _4 (\rhos \cMM ^{1/2}, h _{j - 1}) \\
    &\qquad - \es^{-1} \cN _1 (\td G, \phi _j) h _{j - 1} - \es^{-1} \td \phi _j (\cN _6 - 2 \cN _4 - 2 \cN _5) (\td G, h _{j - 1}),
\eal
\eeq 
Let us analyze them term by term, first with the main forcing term, then the commutators.

\item[$\circ$] \textit{Main forcing}. For the term $\phi _j \td H$, recall the definition of $\td H$: 
\begin{align*}
    \td H (\td G,\cMM) =  %
    &  \f1{\es} \left[ (\cN _2 + \cN _3 + \cN _4 + \cN _6)(\rhos \cMM ^{1/2}, \td G) + \cN (\td G, \rhos \cMM ^{1/2})\right] \\
     &- \cMM ^{-1/2} \eM + \sss \cdot \cF_M \circ \cK_k(\hw_1 + \tw_2 - \cF_E (\td G)) . \nonumber
\end{align*}
Recall $|D ^{\le j} \phi _j| \les_j 1$.  For the term $(\cN _2 + \cN _3 + \cN _4 + \cN _6)(\rhos \cMM ^{1/2}, \td G) = \kp _2 \cs ^{-1} (\div A) [\cM \vc] \td G$ (see \eqref{eqn: N2346}), given $|\al| + |\b| \le j$, 
using $|(\div A) [\cM \vc]| \les \cs ^{\g + 4} \ang \vc ^{\g + 1}$ by \eqref{eq:LWP_divA_point}, we have 
\begin{align*}
    & \left\la D ^{\al, \b} (\phi _j  \kp _2 \cs ^{-1} (\div A) [\cM \vc] \td G) , D ^{\al, \b} h _j \right\ra _{\cYb} \\
    & \qquad \le C_j \iint \kp _2 \cs ^{-1} \la X \ra^{\etab} |(\div A) [D ^{\le j} (\cM \vc)]| \cdot |D ^{\preceq (\al, \b)} \td G| \cdot |D ^{\al, \b} h _j| dXdV \\
    & \qquad \le C_j \iint \ang X ^{\etab} \cs ^{\g + 3} \ang \vc ^\g |D ^{\preceq (\al, \b)} \td G| ^2 dV dX + \f1{16} \iint \ang X ^{\etab} \cs ^{\g + 3} \ang \vc ^{\g + 2} |D ^{\al, \b} h _j| ^2 dV dX \\
    & \qquad \le C_j  \B( \iint \ang X ^{\etab} |D ^{\preceq (\al, \b)} \td G| ^2 dV dX \B)^{ \f{4 }{ \g+2}} 
        \B( \iint \ang X ^{\etab} \cs^{\g+3} \ang \vc^{\g+2} |D ^{\preceq (\al, \b)} \td G| ^2 dV dX \B)^{ \f{ 2 \g }{ \g+2}}  \\ 
     & \qquad \qquad + \f1{16} \iint \ang X ^{\etab} \cs ^{\g + 3} \ang \vc ^{\g + 2} |D ^{\le j} h _j| ^2 dV dX.
\end{align*}
where we have applied the H\"older's inequality and $\cs \les 1$ similar to \eqref{eq:linear_L2_G_hol} in the last inequality. 
The integrals of $\td G$ are further bounded by the $\cY$-norm.
Summing up $\al$ and $\b$, we obtain 
\begin{align*}
    \left\la \phi _j  \kp _2 \cs ^{-1} (\div A) [\cM \vc] \td G , h _j \right\ra _{\cYb ^j}\le 
     C_j \| \td G \| _{\cYb ^j}^{ \f{4}{\g+2}} \| \td G \| _{\cYQ ^j}^{ \f{2\g}{\g+2}}  + \f1{16} \| h _j \| _{\cYQ ^j} ^2.
\end{align*}

For the second term in $\tH$, note that $\| D ^{\le j} \cMM ^{1/2} \| _\sigma 
\les _j \cs^{\f{\g+3}{2}} \les 1$ using \eqref{eqn:DabPhi-sigma} with $i = 0$, $\Phi_0 = \cMM^{1/2}$ 
\eqref{eq:func_Phi}, $|D ^{\le j} \rhos| = | D^{\le j} \cs^3 | \les_j \cs ^3$ by \eqref{eq:dec_S} and Leibniz rule.
Following estimates in Theorem \ref{theo: N nonlin_est} and using \eqref{eqn: uniform-bound-cutoff}, we obtain for any $|\al| + |\b| \le j$ that
\begin{align*}
    & 
    \la D ^{\al, \b} (\phi _j \cN (\td G, \rhos \cMM ^{1/2} )), D ^{\al, \b} h _j \ra _{\cYb} 
    \\
    & \qquad \le C_j \int  \cs ^{-3} \ang X ^{\etab} \| D ^{\preceq (\al, \b)} \td G (s, X, \cdot) \| _{L ^2 (V)} \| D ^{\al, \b} h _j \| _\sigma \| D ^{\le j} (\rhos \cMM ^{1/2}) \| _\sigma  d X \\
    & \qquad \le \f1{16} \int \ang X ^{\etab} \| D ^{\al, \b} h _j \|^2 _\sigma d X + C_j \int \ang X ^{\etab} \| D ^{\preceq (\al, \b)} \td G \|^2 _{L ^2 (V)} d X \\
    & \qquad \le \f1{16} \| D ^{\al, \b} h _j \|^2 _{\cYQ} + C_j \| D ^{\preceq (\al, \b)} \td G \| _{\cYb} ^2.
\end{align*}
Taking summation in $\al, \b$ we conclude 
\begin{align*}
    \la \phi _j \cN (\td G, \rhos \cMM ^{1/2}), h _j \ra _{\cYb ^j} \le \f1{16} \| h _j \| _{\cYQ ^j} ^2 + C_j \| \td G \| _{\cYb ^j} ^2.
\end{align*}
For the remaining two terms, by Cauchy--Schwarz and \eqref{eq:LWP_W1W2} we have
\begin{align*}
     \la \phi _j \cMM ^{-1/2} \eM, h_j \ra_{\cYb ^j} &\le \| \cMM ^{-1/2} \eM \| _{\cYb ^j} \| h _j \| _{\cYb ^j} \le C _j + \| h _j \| _{\cYb ^j} ^2, \\
     \sss \la \phi _j \cF_M \circ \cK_k(\hw_1 + \tw_2 - \cF_E (\td G)) , h_j \ra_{\cYb ^j} &\le \sss \| \cF_M \circ \cK_k(\hw_1 + \tw_2 - \cF_E (\td G))  \| _{\cYb ^j} \| h _j \| _{\cYb ^j} \\
    & %
    \le \| h _j \| _{\cYb ^j} ^2 + C _j \| \hw_1 + \tw_2 \| _{\cXX ^j} ^2 + C _j \| \td G \| _{\cYY ^j} ^2 \\
    & %
    \le \| h _j \| _{\cYb ^j} ^2 + C _j %
    + C _j \| \td G \| _{\cYb ^j} ^2 .
\end{align*}
Combined we get 
\begin{align}
    \la \phi _j \tH, h _j \ra _{\cYb ^j} &\le
    \f1\es \left( %
     C_j  \| \td G \| _{\cYb ^j}^{ \f{4}{\g+2}} \| \td G \| _{\cYQ ^j}^{ \f{2\g}{\g+2}} 
    + 
    \f18 \| h _j \|^2 _{\cYQ ^j} + C_j \| \td G \| _{\cYb ^j} ^2 \right) + C_j \| h _j \| _{\cYb ^j} ^2  + C _j. \label{htilde_high}
\end{align}
This concludes the estimate for $\phi_j \tH$.

\item[$\circ$] \textit{Transport commutators}. 
We now deal with the second term of $\tH_j$. Take $|\al| + |\b| \le j$. 
Using \eqref{eqn: ds_T_cutoff} with $( \vp_R, R) \rsa (\phi_j, 2^{-j} R), \phi_j = \vp_{2^{-j}R}$
and \eqref{eq:norm_diffu_equiv},  we get
\begin{align*}
    & \la D ^{\al, \b} \left((\pa _s + \cT) \phi _j \cdot h _{j - 1}\right), D ^{\al, \b} h _j \ra _{\cYb} \\
    & \qquad \le \iint \ang X ^{\etab} |D ^{\le j} (\pa _s + \cT) \phi _j| \cdot |D ^{\preceq (\al, \b)} h _{j - 1}| \cdot |D ^{\al, \b} h _j| dX dV \\
    & \qquad \leq C_j \iint \ang X ^{\etab} \cdot 2^j  R ^{-1} \ang X \ang \vc ^2 |D ^{\preceq (\al, \b)} h _{j - 1}| \cdot |D ^{\al, \b} h _j| dX dV \\
    & \qquad \le C _j \| D ^{\al, \b} h _j \| _{\cYb} ^2 + C_j R ^{-2} \| D ^{\preceq (\al, \b)} (\ang X \ang \vc ^2 h _{j - 1}) \| _{\cYb} ^2 .
\end{align*}

Take summation in $\al, \b$ we conclude 
\begin{align}
    \la (\pa _s + \cT) \phi _j \cdot h _{j - 1}, h _j \ra _{\cYb ^j} & \le C _j  \| h _j \| _{\cYb ^j} ^2 + C_j R ^{-2} \| \ang X \ang \vc ^2 h _{j - 1} \| _{\cYb ^j} ^2. \label{trans_comm_lin}
\end{align}

\item[$\circ$] \textit{Diffusion commutator}.
Recall the weighted Laplacian from \eqref{eq: weighted-laplacian}
\[
    \D _W F  = -   \nu_{\etab}^{-1} \ang X ^{2M} \ang \vc ^{2N} F   +   \sum _{|\al_1|  +  |\b_1| = 1}  \ang X ^{M-\etab} \ang \vc ^N  
    \pa_X^{\al_1} \pa_V^{\b_1}  \B( \vp_1^{2 |\al_1| } \cs^{2 |\b_1|} \ang X^{\etab}  \pa_X^{\al_1} \pa_V^{\b_1}  (\ang X ^{M} \ang \vc ^N F)  \B)  \notag , \\
\]
and $h_{j-1} = \tFR \cdot \phi_{j-1}$.  Next, we compute the commutator $  [\D_W, \phi_j] h_{j-1}$. 
In the support of $h_j$, we have $\tFR = h_{j-1}$. Denote $g_{j-1} = \ang X^M \ang \vc^N h_{j-1}$. For each $|\al_1| + |\b_1| = 1$, we have 
\begin{align*}
  \pa_X^{\al_1} \pa_V^{\b_1}  \B( \vp_1^{2 |\al_1| } \cs^{2 |\b_1|} \ang X^{\etab}  \pa_X^{\al_1} \pa_V^{\b_1}  
    (  g_{j-1} \phi_j )\B)  
    = &   \pa_X^{\al_1} \pa_V^{\b_1}  \B( \vp_1^{2 |\al_1| } \cs^{2 |\b_1|} \ang X^{\etab}    g_{j-1}\pa_X^{\al_1} \pa_V^{\b_1}  \phi_j   \B)  \\ 
  & + \pa_X^{\al_1} \pa_V^{\b_1} \B( 
   \vp_1^{2 |\al_1| } \cs^{2 |\b_1|} \ang X^{\etab}   \pa_X^{\al_1} \pa_V^{\b_1}  
     g_{j-1} \B)  \cdot  \phi_{j} \\ 
    &   +    \vp_1^{2 |\al_1| } \cs^{2 |\b_1|} \ang X^{\etab}   \pa_X^{\al_1} \pa_V^{\b_1}  
     g_{j-1}   \cdot   \pa_X^{\al_1} \pa_V^{\b_1}  \phi_{j} 
   := \sum_{1\leq i\leq 3} I_{i,\al_1, \b_1 } . %
\end{align*}
 Since the commutator associated with the first term in $\D_W F$ is $0$ and the term $I_{2,\al_1, \b_1}$ 
 is canceled in the commutator, using Leibniz rule, we obtain 
\[
    \th [\D_W, \phi_j] h_{j-1} = \th  \sum _{|\al_1|  +  |\b_1| = 1}  \ang X ^{M-\etab} \ang \vc ^N (I_{1,\al_1, \b_1 } + I_{3,\al_1, \b_1 }  ) .
\]

For $|\al_1| + |\b_1| \leq 1, i=1,3$, since $I_{1i,\al_1,\b_1}$ involves at most one derivative acting on $g_{j-1}$, using integration by parts, Leibniz rule, and \eqref{eq:uniform-cutoff-phij}, we conclude 
\begin{align}
    \th  \la [\D _W , \phi _j] h _{j - 1}, h _j \ra _{\cYb ^j}  \le C _j \th \| \ang X \ang \vc ^2 h _{j - 1} \| _{\cYb ^j} ^2 + \f\th4 \| \ang X \ang \vc ^2 h _j \| _{\cYb ^{j + 1}} ^2 \label{comm_lap+lin}
\end{align}

\item[$\circ$] \textit{Collision commutator: $\cN _1$}. Finally, we bound the last four terms in $\td H _j$ \eqref{eq:Hj_recall}. Start with the two terms involving $\cN _1$:
\[
    - \cN _1 (\rhos \cMM ^{1/2}, \phi _j) h _{j - 1} - \cN _1 (\td G, \phi _j) h _{j - 1}.
\]

By \eqref{eqn: nabla_V_varphi}, we get
\begin{align*}
    \cN _1 (\td G, \phi _j) &= 
     \div_V (A [\cMM^{1/2} \td G] \na _V \phi_j) = 
    \div _V (A [\cMM ^{1/2} \td G \vc] \cdot \kp _2 \cs ^{-1} \td \phi _j) \\
    &= \div_V A [\cMM ^{1/2} \td G \vc] \cdot \kp _2 \cs ^{-1} \td \phi _j + A [\cMM ^{1/2} \td G \vc] \kp _2 \cs ^{-2} \sum _{i = 1} ^3 D ^{0, \ee _i} \td \phi _j \ee _i.
\end{align*}
Using \eqref{eqn: divA-pointwise} with $i=0$ and \eqref{eq:uniform-cutoff-phij}, we get %
\begin{align*}
    |D ^{\le i} \cN _1 (\td G, \phi _j)| &\les _i \cs ^{-2} |D ^{\le i} A [\cMM ^{1/2} \td G \vc]| \cdot |D ^{\le i + 1} \td \phi _j| + \cs ^{-1} |D ^{\le i} \div A [\cMM ^{1/2} \td G \vc]| \cdot |D ^{\le i} \td \phi _j| \\
    & \les _i \cs ^{\g + 2 - 2} \| D ^{\le i} \td G \| _{L ^2 (V)} R ^{-2} \ang \vc ^{\g + 2} .
\end{align*}
Using Leibniz rule and Lemma \ref{lem:Q_commu} with $(f, g) \rsa (h_{j-1}, h_j)$, we obtain 
\bseq\label{eq:LWP_commu_N1}
\begin{align}
    & \left|\la \cN _1 (\td G, \phi _j) \cdot h _{j - 1}, h _j \ra _{\cYb ^j} \right| 
    \notag \\
    \label{eq:LWP_commu_N1:a}
    & \qquad \les_j \sum_{i\leq j} \int \ang X^{\etab}    |D ^{\le i} \cN _1 (\td G, \phi _j) | \cdot |D^{\leq j-i} h_{j-1}| \cdot | D^{\leq j} h _j |  \\
    & \qquad \les_j R^{-2} \left( 
        \| \td G \|_{\cYb^{\max \{\kkl, j\}} }  \| h_{j-1}\|_{\cYQ^{j-1}} 
        + \| \td G\|_{\cYb^{\kkl}}  \|\Lam ^\f12 h _{j - 1} \| _{\cYb^j} 
    \right) \|  h_j \|_{\cYQ^j}.\notag 
\end{align}
When applying the above estimate with $\td G$ replaced by $\brho \cMM^{1/2}$, we use 
$\rhos = \cs^3$ \eqref{eq:Euler_profi_modi}, $\cMM = \cs^{-3} \mu(\vc)$ \eqref{eq:localmax2},
and  estimates \eqref{eq:dec_S}, \eqref{eq:M11/2}, \eqref{eq:cross_pf2:a} to obtain 
\beq\label{eq:LWP_L2_mass1}
\| D ^{\le j} (\rhos \cMM ^{1/2}) \| _{L ^2 (V)} 
\les_j \| \cs^3 \ang \vc^{2j} \cMM^{1/2} \|_{L^2(V)}
=\cs^{3/2} \| \ang \vc^{2 j} \mu(\vc)^{1/2}  \|_{L^2(V)} \les \cs^3. 
\eeq
 So
\begin{align*}
    |D ^{\le j} \cN _1 (\rhos \cMM ^{1/2}, \phi _j)| &\les _j \cs ^{\g + 5 - 2} R ^{-2} \ang \vc ^{\g + 2} = R ^{-2} \Lam.
\end{align*}
Plugging the above estimate into \eqref{eq:LWP_commu_N1:a} and using the definition of $\cYQ$-norm \eqref{norm:Y}, we get 
\beq\bal
    \label{eq:LWP_commu_N1:b}
    \left|\la \cN _1 (\rhos \cMM ^{1/2}, \phi _j) \cdot h _{j - 1}, h _j \ra _{\cYb ^j} \right| &\le 
    C _j R^{-2}  ( \| \Lam ^\f12 h _{j - 1} \|_{\cYb ^j} + \| h_{j-1} \|_{\cYQ^{j-1}}  ) \| h_j \|_{\cYQ ^j}  
    \\
    &\le \f{1}{16} \| h _j \| _{\cYQ ^j} ^2 + \f{C_j}{R ^4} \| \Lam ^\f12 h _{j - 1} \| _{\cYb ^j} ^2.
\eal\eeq
\eseq

\item[$\circ$] \textit{Collision commutator: $\cN _{4,5,6}$}.
We now estimate the commutator terms related to $\cN_4, \cN_5, \cN_6$, which are
\begin{align*}
    \td \phi _j \left( 2 \cN _4 (\rhos \cMM ^{1/2}, h _{j - 1}) - (\cN _6 - 2 \cN _4 - 2 \cN _5) (\td G, h _{j - 1}) \right).
\end{align*}
Recall from \eqref{N(f,g)}
\beq\label{eq:LWP_N_recall}
\bal  
\cN _4 (f, g) &=  - \kp _2 \cs ^{-1} \vc ^\top A [\cMM ^{1/2} f] \nabla_V g \\ 
\cN _5 (f, g) &= \kp _2 ^2 \cs ^{-2} g \vc ^\top A [\cMM ^{1/2} f] \vc \\
\cN _6 (f, g) &= \kp _2 \cs ^{-1} g \div_V A [\cMM ^{1/2} f] \cdot \vc.
\eal 
\eeq 
Below, we let $l=4,5,6$.  For any $|\al| +| \b | \leq j$,  we use Leibniz rule, 
 $|D ^{\le j} \td \phi _j| \les _j R ^{-2}$ by \eqref{eqn: uniform-bound-cutoff}, Lemma \ref{lem: nonlinear}, and Lemma \ref{lem: N-derivative-bound} to bound the lower order derivatives on $h _{j - 1}$:
\begin{align*}
    & \left\la D ^{\al,\b} ( \td \phi _j \cN _l (\td G, h _{j - 1}) ) - \td \phi_j \cN_l( \td G, D^{\al, \b} h_{j-1})  , D ^{\al, \b} h _j  \right\ra _V \\
    & \qquad \les _j \sum _{1 \le i \le j} R ^{-2} \cs ^{-3} \| D ^{\le i} G \| _{L ^2 (V)} \| D ^{\le j - i} h _{j - 1} \| _\s \| D^{\leq j} h _j \| _\s.
\end{align*}
This is the same situation as $I _{i, j}$ for $i \ge 1$ in Lemma \ref{lem:Q_commu}, so they are bounded in the same way as $\cN _1$ terms. 

For the main term, applying estimate of $\cN_{4,5,6}$ in  Lemma \ref{lem: nonlinear} and using \eqref{eq:uniform-cutoff-phij}, we obtain 
\begin{align*}
    \B| \B\la \td \phi _j \cN _l (\td G, D ^{\al, \b} h _{j - 1}), D ^{\al, \b} h _j \B\ra _V \B|
    \les R^{-2} \cs^{-3} \| \td G\|_{L ^2 (V)}  \| \Lams^{1/2} D ^{\al, \b} h _{j - 1} \|_{L ^2 (V)}  \|  D ^{\al, \b} h _j  \|_{\s}. 
\end{align*}
With the same idea as in the estimate of $I_{0, j}$ in Lemma \ref{lem:Q_commu} we obtain 
\begin{align*}
    & \B| \left\la \td \phi_j \cN _{l} (\td G, D ^{\al, \b} h _{j - 1}), D ^{\al, \b} h _j \right \ra _{\cYb ^j} \B|  \les R^{-2} \| \td G \|_{\cYb^\kkl} \|\Lam ^\f12 h _{j - 1} \|_{\cYb^j} \| h _j \|_{\cYQ ^j}.
\end{align*}

Summarize in $\al,\b$, we conclude the same upper bound for $\cN_{l},l=4,5,6$ terms as $\cN _1$ terms: %
\bseq\label{eq:LWP_commu_N4}
\beq\label{eq:LWP_commu_N4:a}
    \B| \big\la \td \phi _j \cN _l (\td G, h _{j - 1}), h _j \big\ra _{\cYb ^j} \B|  
     \les_j R^{-2} \left( 
        \| \td G \|_{\cYb^{\max \{\kkl, j\}} }  \| h_{j-1}\|_{\cYQ^{j-1}} 
        + \| \td G\|_{\cYb^{\kkl}}  \|  \Lam ^\f12 h _{j - 1} \| _{\cYb^j} 
    \right) \| h _j \|_{\cYQ^j}.
\eeq 
Replacing the above estimate of $\td G$ by $\rhos \cMM ^{1/2}$ and replacing the bound of $\| D^{\leq j} \td G\|_{L ^2 (V)}$ 
(based on  \eqref{eq:linf_Y} or \eqref{eqn: smallness-G-L2}) by $\| D^{\leq j} (\brho \cMM^{1/2} ) \|_{L ^2 (V)} \les_j \cs^3$ (see \eqref{eq:LWP_L2_mass1}),  we have the same bound as \eqref{eq:LWP_commu_N1:b}
\beq\label{eq:LWP_commu_N4:b}
\bal 
    \big \la \td \phi _j \cN _4 (\rhos \cMM ^{1/2}, h _{j - 1}), h _j \big \ra _{\cYb ^j} & \le \f{1}{16} \| h _j \| _{\cYQ ^j} ^2 + \f{C_j}{R ^4} \| \Lam ^\f12 h _{j - 1} \| _{\cYb ^j} ^2.
\eal 
\eeq
\eseq 
\end{itemize}

\paragraph{Summary} 
Summarizing, %
we obtain 
\begin{align*}
& \hspace{-3em} \underbrace{
    \f12 \f{d}{ds} \| h _j \| _{\cYb ^j} ^2 
    - C_j \| h _j \| _{\cYb ^j} ^2 
    - \f{1}{16} \| h _j \| _{\cYQ ^j} ^2
}_{\eqref{trans_lin}} 
+ \underbrace{ 
    \f\theta 2 \|\ang X \ang \vc ^2 h_j \| _{\cYb ^{j + 1}} ^2
	- C_j \theta \|\ang X \ang \vc ^2 h_{j-1} \| _{\cYb ^{j }} ^2  
}_{\eqref{eq:LWP_Lap1}}\\
\le & - \underbrace{
    \f7{8\es} \| h _j \| _{\cYQ ^j} ^2 
    + \f{C_j}{\es}  \| h _{j - 1} \| _{\cYQ ^{j - 1}} ^2
}_{\eqref{main_coll_main}} 
+ \underbrace{ 
    \f1\es \B( 
        \cnon \| \td G \| _{\cYb^{\kkl}} \| h _j \| _{\cYQ ^j} 
        + C_j \one_{j > \kkl} \| \td G \|_{ \cYb^{ j }} \| h_{j-1} \| _{\cYQ ^{j-1}}  
    \B) \| h _j \| _{\cYQ ^j} 
}_{\eqref{sec_coll_lin}} \\
& + \underbrace{
    \f1\es \left( %
        C_j \| \td G \| _{\cYb ^j}^{ \f{4}{\g+2}} \| \td G \| _{\cYQ ^j}^{ \f{2\g}{\g+2}} 
        + \f18 \| h _j \|^2 _{\cYQ ^j} + C _j \| \td G \| _{\cYb ^j} ^2 \right) + C_j \| h _j \| _{\cYb ^j} ^2  
        + C _j 
}_{\eqref{htilde_high}} \\
& + \underbrace{
    C _j  \| h _j \| _{\cYb ^j} ^2 + \f{C _j}{R ^2} \| \ang X \ang \vc ^2 h _{j - 1} \| _{\cYb ^j} ^2
}_{\eqref{trans_comm_lin}}
+ \underbrace{ 
    C _j \th \| \ang X \ang \vc ^2 h _{j - 1} \| _{\cYb ^j} ^2 + \f\th4 \| \ang X \ang \vc ^2 h _j \| _{\cYb ^{j + 1}} ^2 
}_{\eqref{comm_lap+lin}} \\
& + \underbrace{ 
    \f{C_j}{\es R^2} \left( 
        \| \td G \|_{\cYb^{\max \{\kkl, j\}} }  \| h_{j-1}\|_{\cYQ^{j-1}} 
        + \| \td G\|_{\cYb^{\kkl}}  \|\Lam ^\f12 h _{j - 1} \| _{\cYb^j} 
    \right) \| h _j \|_{\cYQ^j}
}_{\eqref{eq:LWP_commu_N1:a}, \eqref{eq:LWP_commu_N4:a}} \\
& 
+ \underbrace{ 
    \f{1}{8\es} \| h_j \|_{\cYQ^j}^2 
    + \f{C_j}{\es R^4} \| \Lam ^\f12 h_{j - 1} \|_{ \cYb^j}^2 
}_{\eqref{eq:LWP_commu_N1:b}, \eqref{eq:LWP_commu_N4:b}}.
\end{align*}

Note that we have multiplied the estimate in  \eqref{main_coll_main}, \eqref{sec_coll_lin}, \eqref{eq:LWP_commu_N1}, \eqref{eq:LWP_commu_N4} by $\es^{-1}$, which is associated with the $\cN_i$-term in \eqref{eq: h_linear_viscous}. For the upper bound in \eqref{eq:LWP_commu_N1:a}, \eqref{eq:LWP_commu_N4:a}, since $\td G \in \JJ_{\ze}^{k}$ with $ \ze < 1 $ and $R>1$,  we obtain $\| \td G\|_{\cYb^{\kkl}} \leq 1$ and 
\[
\bal 
  \f{C_j}{\es R^2}  
        \| \td G \|_{\cYb^{\max \{\kkl, j\}} }  \| h_{j-1}\|_{\cYQ^{j-1}} \| _{\cYb^j}  \| h _j \|_{\cYQ^j}
& \leq \one_{j > \kkl} 
  \f{C_j}{\es R^4}  
        \| \td G \|_{\cYb^{ j } }  \| h_{j-1}\|_{\cYQ^{j-1}} \| _{\cYb^j}  \| h _j \|_{\cYQ^j} \\ 
& \quad + \f{C_j}{\es}  \| h _{j-1} \|_{\cYQ^{j-1} }^2 
+ \f{1}{32 \es} \| h _{j} \|_{\cYQ^{j} }^2 , \\
 \f{C_j }{\es R^2} \| \td G\|_{\cYb^{\kkl}}  \|\Lam ^\f12 h _{j - 1} \| _{\cYb^j}  \| h _j \|_{\cYQ^j}
 & \leq  \f{C_j }{\es R^4}   \|\Lam ^\f12 h _{j - 1} \| _{\cYb^j}^2 +  \f{1}{32 \es} \| h _j \|_{\cYQ^j}^2 .
\eal 
\]

Since $2 \geq \f{\g + 2}2, \g+3\geq 0$, we obtain $\Lam ^\f12 \les \ang \vc^{ (\g+2)/2} \les \ang X \ang \vc ^2$. Using \eqref{eq:norm_diffu_equiv}, we bound 
\begin{align*}
    \| \Lam ^{1/2} h _{j - 1} \| _{\cYb^j} \les \| \ang X \ang \vc ^2 h _{j - 1} \| _{\cYb ^j}.
\end{align*} 
Combining similar terms in the above two estimates, using the diffusion and $\es \leq 1, R > 1$, we simplify the energy estimates as
\beq\label{eq:LWP_EE_Hj3}
\bal 
    \f12 \f{d}{ds} \| h _j \| _{\cYb ^j} ^2 
    & \le C _j \| h _j \| _{\cYb ^j} ^2 
    - \f{1}{\es} \left( 
        \f12 - \cnon \| \td G \| _{\cYb^{\kkl}} 
    \right) \| h _j \| _{\cYQ ^j} ^2 + \f{C_j}{\es} \| h _{j - 1} \| _{\cYQ ^{j - 1}} ^2 \\
    & \quad + \one_{j > \kkl} \f{C_j}{\es} \| \td G \|_{\cYb ^j}  
    \| h_{j-1} \| _{\cYQ ^{j - 1}} \| h _j \| _{\cYQ ^j} + \f{C_j}{\es} \left( 
        \| \td G \| _{\cYb ^j} ^\f{4}{\g+2} \| \td G \| _{\cYQ ^j} ^\f{2\g}{\g+2}
        + \| \td G \| _{\cYb ^j} ^2 
    \right) + C _j \\
    & \quad - \f\th4 \|\ang X \ang \vc ^2 h _j \| _{\cYb ^{j + 1}} ^2 
    + C _j \left(\th + \f1{R ^2} + \f1{\es R^4} \right) \|\ang X \ang \vc ^2 h_{j-1} \| _{\cYb ^j} ^2.
\eal 
\eeq 
Changing the constants for the diffusion terms to other absolute constants, we prove \eqref{eq:LWP_EE_Hj}.

\paragraph{Proof of \eqref{eq:LWP_EE_Hj2}}
Summing the  $L^2$ estimates \eqref{est: Unif_inR_lin} and the weighted $H^j$ estimates \eqref{eq:LWP_EE_Hj} 
(or see above \eqref{eq:LWP_EE_Hj3}) with weight $\varpi_{Z,j}$, we obtain
\[
\bal
	& \f{1}{2} \f{d}{d s} \sum_{j \leq k} \varpi_{Z,j} \| h_j \|_{\cYb^j}^2   \\ 
	 \leq &  \sum_{j\leq k}  \varpi_{Z,j} \B( C_j \| h_j \|_{\cYb^j}^2
- \f{1}{\es} \big(  \f{1}{2} -  (\cnon + \bar C_{\cN,0}) \| \td G\|_{ \cYb^\kkl }   \big)	\| h_j\|_{\cYQ^j}^2    + \udb{ \one_{j>0} \f{  \bar C_{j,1}}{\es}  \| h _{j - 1} \| _{\cYQ ^{j - 1}} ^2 }_{:=I_{j,2}}   \\ 
&     + \one_{j > \kkl} \f{ C_j}{\es} \| \td G\|_{\cYb^{  j}}  \| h_{j-1} \|_{\cYQ^{  j-1 }}  \| h_j \|_{\cYQ^j} 
  +   \f{C_j}{\es} ( \| \td  G\|_{\cYb^j}^{ \f{4 }{\g+2}} \| \td  G\|_{\cYQ^j}^{ \f{2 \g }{\g+2}} 
  +  \| \td G \| _{\cYb ^j} ^2 ) + C_j \\
	&  - \underbrace{ \f\th4 \|\ang X \ang \vc ^2 h_{j} \| _{\cYb ^{j + 1}} ^2 }_{:= I_{j, 3}}  + \underbrace{ \one_{j>0}  \bar C_{j,2} (  \th + \f{1}{ \es R^2})  \|\ang X \ang \vc ^2 h_{j-1} \| _{\cYb ^{j }} ^2 }_{:= I_{j,4} } \B).
   \eal
\]

Recall the norms $\cZZ^j, \cZZQ^j$ and weight $\varpi_{Z,j}$ from \eqref{norm:cZZ}. 
By definition, we obtain $\varpi_{Z,0} = 1, \bar C_{0, i} = 0$.
Since $ R \es \geq 1, R^{-1} \leq \th$, we obtain 
\[ 
\bar C_{j, 2} (\th + \f{1}{R^2 \es} ) \varpi_{Z, j}
\leq  \bar C_{j, 2} \cdot 2 \th \cdot  \varpi_{Z, j} <  \f{\th }{8} \varpi_{Z, j-1},  \quad 
\bar C_{j, 1}  \varpi_{Z, j} 
<  \f{ 1}{8} \varpi_{Z, j-1} ,  \quad 
 \forall  j\geq 1.
\]

For the weighted sum of the diffusion terms $I_{j,i},i=2,3,4$, using the above estimates, we obtain
\[
\bal 
\sum_{ j \leq k} \varpi_{Z, j}  I_{j,2} 
& \leq \sum_{ 0< j \leq k }  \f{  \varpi_{Z,j-1} }{ 8\es} 
\| h _{j - 1} \| _{\cYQ ^{j - 1}}^2
=  \sum_{0\leq j\leq k-1}  \f{ \varpi_{Z,j} }{ 8\es} 
\| h _{j } \| _{\cYQ ^{j }}^2 , \\ 
\eal 
\] 
and 
\[
\bal 
\sum_{ j\leq k}  \varpi_{Z, j}  I_{j, 4}
 \leq 
 \sum_{ 0<  j\leq k }  \f{\th}{8} \varpi_{Z,j-1} 
  \|\ang X \ang \vc ^2 h_{j-1} \| _{\cYb ^{j }} ^2 
 =  \sum_{ 0\leq j\leq k -1 }   \f{\th}{8} \varpi_{Z,j} 
  \|\ang X \ang \vc ^2 h_{j} \| _{\cYb ^{j +1 }} ^2
\leq \f{1}{2} \sum_{ j\leq k}  \varpi_{Z, j}  I_{j, 3}. 
\eal 
\]

Thus, in the above weighted sum, the diffusion terms have the negative sign up to the term $  (\cnon + \bar C_{\cN,0}) \| \td G\|_{ \cYb^\kkl } \| h_j \|_{\cYQ^j}^2$.
By definition \eqref{norm:cZZ}, for any $i\geq 0$, we have
\beq\label{eq:LWP_sum:hi}
	\| h_i \|_{\cYb^i} \les_i \| \tFR \|_{\cZZ_R^i} ,
	\quad 	\| h_i \|_{\cYQ^i} \les_i \| \tFR \|_{\cZZQR^i} .
\eeq

Using \eqref{eq:LWP_sum:hi}, the $Z$-norm \eqref{norm:cZZ}, 
and dropping the weighted diffusion in $I_{j,3},I_{j,4}$, we prove
\[
\bal
	\f{1}{2} \f{d}{d s} \| \tFR\|_{\cZZ_R^k}^2
&	\leq 
C_k \| \tFR \|_{\cZZ_R^k }^2 + \f{1}{\es} \big(  -  \f{3 }{8 } + (\cnon + \bar C_{\cN,0})  \| \td G\|_{ \cYb^\kkl }   \big) \| \tFR \|_{\cZZQR^k}^2 \\
& \quad +    \one_{ k > \kkl} \f{ C_k}{\es} \| \td G\|_{\cYb^{  k}} 
	\| \tFR \|_{\cZZQR^{  k-1 }}  \| \tFR \|_{\cZZQR^k} 
  +   \f{C_k}{\es} ( \| \td  G\|_{\cYb^k}^{ \f{4 }{\g+2}} \| \td  G\|_{\cYQ^k}^{ \f{2 \g }{\g+2}} + 
     \| \td G \| _{\cYb ^k} ^2  )  +  C_k . 
	\eal
\]
Changing the dummy variable $k$ to $j$, we prove \eqref{eq:LWP_EE_Hj2}.
\end{proof}

\subsubsection{Uniform energy estimates}

Recall the constant $\cnon$ from \eqref{eq:non_Q:micro} Theorem \ref{theo: N nonlin_est}, $\bar C_{\cN,0}$ from Lemma \ref{lem: Unif_inR_lin}, and $\ze_0$ from Lemma \ref{lem:parab}.  We choose 
\bseq\label{def:LWP_ze}
\beq\label{def:LWP_ze1}
	\ze_1 =  \min \left\{ 
        \f{1}{ 8 ( \cnon + \bar C_{\cN,0}) } , \ze_0, 1  
    \right\}.
\eeq
Next, we assume $ \td G \in \bbJ_{\ze_1}^k$. From \eqref{def:LWP_ze1}, we obtain 
\beq
	 - \f{3  }{8} +   (\cnon + \bar C_{\cN, 0}) \| \td G \|_{\cYb^{\kkl}} 
	 \leq 	 - \f{3 }{8} +    (\cnon + \bar C_{\cN, 0})  \ze_1 < -\f{ 1}{4} .
\eeq
\eseq

Using \eqref{def:LWP_ze} and %
Lemma \ref{lem:LWP_Hj} with $j \geq \kkl$, we obtain
\[
\bal
	\f{1}{2} \f{d}{d s} \| \tFR \|_{\cZZ^j}^2 
	& \leq C_j \|  \tFR \|_{ \cZZ^j }^2
	 - \f{  1 }{ 4 \es}  \| \tFR \|_{ \cZZQ^j}^2  
	 	+ \one_{j > \kkl} \f{  C_j}{\es} \| \td G\|_{ \cYb^{  j}}  \| \tFR \|_{ \cZZQ^{  j-1 }}  \| \tFR \|_{ \cZZQ^j}\\
 & \qquad + C_j  
 + C_j \es^{-1} ( \| \td  G\|_{\cYb^j}^{ \f{4 }{\g+2}} \| \td  G\|_{\cYQ^j}^{ \f{2 \g }{\g+2}} + 
     \| \td G \| _{\cYb ^j} ^2  ) \\
\eal
\]
Applying Young's inequality 
\beq\label{eq:LWP_Young}
	  \f{C_j }{\es} \| \td G\|_{ \cYb^{  j}}  \| \tFR \|_{ \cZZQ^{  j-1 }}  \| \tFR \|_{ \cZZQ^j}
\leq \f{ 1 }{ 8\es} \| \tFR \|_{ \cZZQ^j}^2
+ 	\f{  C_j}{  \es} \| \td G\|_{ \cYb^{  j}}^2  \| \tFR \|_{ \cZZQ^{  j-1 }}^2,
\eeq
we bound 
\bseq\label{eq:LWP_gron1}
\beq 
\bal
	\f{1}{2} \f{d}{d s} \| \tFR \|_{\cZZ^j}^2 
	 \leq C_j \|  \tFR \|_{ \cZZ^j }^2
	 - \f{   1 }{ 8 \es}  \| \tFR \|_{ \cZZQ^j}^2   + C_j \cR_j  ,
\eal 
\eeq 
where $\cR_j$ denotes the forcing terms 
\beq 
  \cR_j :=     \one_{j > \kkl} \es^{-1} \| \td G\|_{ \cYb^{  j}}^2 
	\| \tFR \|_{ \cZZQ^{  j-1 }}^2    + 1   +  \es^{-1}  ( \| \td  G\|_{\cYb^j}^{ \f{4 }{\g+2}} \| \td  G\|_{\cYQ^j}^{ \f{2 \g }{\g+2}} + \| \td G \| _{\cYb ^j} ^2  ) .
\eeq 
\eseq

Integrating \eqref{eq:LWP_gron1} over $s$, we obtain
\beq\label{eq:LWP_gron2}
\bal
&	\f{1}{2} \| \tFR(s) \|_{\cZZ^j}^2 +\int_0^{s} \f{ 1 }{ 8 \e_{\tau}}
	 \| \tFR \|_{ \cZZQ^j}^2   \leq  	\f{1}{2} \| \tFR( 0 ) \|_{\cZZ^j}^2  
 +  C_j \int_{ 0 }^s  ( \|  \tFR \|_{ \cZZ^j }^2 + \cR_j)(\tau) d \tau .
\eal
\eeq
Applying Gr\"onwall's inequality to $	\| \tFR(s) \|_{\cZZ^{ j }}$, 
and using $ \e_{\tau}^{-1} \leq \es^{-1} $ for $\tau \leq s$, we obtain 
\beq\label{eq:LWP_gron3}
\bal
	\| \tFR(s) \|_{\cZZ^{ j }}^2 & \leq \operatorname e^{ C_j s } ( \| \tFR(0) \|_{\cZZ^{j}}^2
	+    C_j \int_0^s \cR_j(\tau) d \tau  ). 
\eal
\eeq
Using \eqref{eq:LWP_init}, \eqref{eqn: uniform-bound-cutoff} and \eqref{norm:cZZ}, we obtain  
\beq\label{eq:FR_IC}
\| \tFR( 0 ) \|_{\cZZ^j}^2  \les_j \| \td F(0) \|_{\cYb^j }^2. 
\eeq 
Applying \eqref{eq:LWP_gron3} in the upper bound in \eqref{eq:LWP_gron2},
and using  \eqref{eq:FR_IC},  we prove 
\bseq\label{eq:LWP_gron_Hj1}
\beq 
\bal
	 \| \tFR(s) \|_{\cZZ^j}^2 +\int_0^{s} \f{1}{\e_{\tau}}
	 \| \tFR \|_{ \cZZQ^j}^2
	& \leq C_j \operatorname e^{ C_j s } \left(  \| \td F( 0 ) \|_{\cYb^j}^2 
	 +   \int_0^s \cR_j(\tau) d \tau \right) , \\
\eal
\eeq
 where $C_j$ may change from line to line.

 For $\cR_j$ \eqref{eq:LWP_gron1}, using $\es \geq \e_0 \operatorname e^{-C s}$ \eqref{eq:EE_para1} and H\"older's inequality, 
 for $s \leq T$, we obtain 
  \beq
\bal 
  \int_0^{s} \f{1}{\e_{\tau}}   \| \td  G\|_{\cYb^j}^{ \f{4 }{\g+2}} \| \td  G\|_{\cYQ^j}^{ \f{2 \g }{\g+2}} d \tau 
& \leq \f{C \operatorname e^{Cs}}{\e_0}  \| \td G\|_{L^{2}(0, s; \cYb^{j} ) }^{ \f{4}{\g+2}}
\| \td G\|_{L^{2}( 0, s; \cYQ^{j} ) }^{ \f{2\g}{\g+2}}  . 
\eal 
 \eeq

Combining the above two estimates and using \eqref{eq:LWP_gron1}, for $s \leq T$, we obtain 
\begin{align} 
	 \| \tFR(s) \|_{\cZZ^j}^2 +\int_0^{s} \f{1}{\e_{\tau}}
	 \| \tFR \|_{ \cZZQ^j}^2
& \leq 
C_j \operatorname e^{ C_j s } \B( \| \td F( 0 ) \|_{\cYb^j}^2 
    + \int_0^{s}  
   \f{  1 }{\e_{\tau}} (  \one_{j > \kkl}  \| \tFR \|_{ \cZZQ^{  j-1 }}^2  + 1) \| \td G\|_{ \cYb^{  j}}^2(\tau) d \tau \notag \\
  & \hspace{6em}  + %
  \f{1}{\e_0}  \| \td G\|_{L^{2}(0, s; \cYb^{j} ) }^{ \f{4}{\g+2}}
\| \td G\|_{L^{2}( 0, s; \cYQ^{j} ) }^{ \f{2\g}{\g+2}}  
  +  s   \B) . 
  \end{align}
\eseq 
 Note that when $j \leq \kkl$, the $\tFR$-term on the right hand side vanishes.

Thus, given $\td G \in \bbJ_{\ze_1}^k$  with  $ k \geq \kkl$, %
using \eqref{eq:LWP_gron_Hj1} with $\nu = 1$ inductively for $j =\kkl, \kkl+1, \kkl+2,..,k$, we obtain $\tFR \in L^{\infty}( (0, T), \cZZ^{j} ) \cap  L^2( (0, T); \cZZQ^{j} )$ with
\beq\label{eq:LWP_uniform_Hj}
	 \| \tFR(s) \|_{\cZZ^j}^2 +\int_0^{s} \f{1}{\e_{\tau}}
	 \| \tFR \|_{ \cZZQ^j}^2
	 \leq %
	 C_j\B(s, \| \td F( 0 ) \|_{\cYb^j},  \; \sup_{ \tau \leq s } \| \td G(\tau)\|_{\cYb^j}^2 
     +  \int_0^{\tau} \f{1}{\e_{\tau}} \| \td G \|_{\cYQ^j}^2 
     \B) , 
\eeq
for any $j\leq k$ and any $s \leq T$. Here, we do not require $T$ to be small.

\subsubsection{Convergence}

Suppose $\td G \in \bbJ_{\ze_1}^k$. We consider $ s \in [ 0 ,  1]$. We take $R_n = n, \th_n = n^{-1}$. 
For large $n \geq N_s$, assumptions $R \cdot \es \geq 1, R^{-1} \leq \th$ in Lemma \ref{lem:LWP_Hj} are satisfied.  Recall the cutoff function $\phi_{L}$ from \eqref{eqn: defn-varphi} and the uniform estimates \eqref{eqn: uniform-bound-cutoff}. Since for fixed $L = m$, for $n$ large enough, 
we have $ \phi_L = \phi_n \phi_{L} $. Using estimate similar to \eqref{eq:Hj_chain}, \eqref{eq:uniform-cutoff-phij}, and \eqref{norm:cZZ}, for any $ j \geq 0$, we obtain
\beq\label{eq:Hj_chain2}
	 \|\phi_L \td F_{\th_n, R_n}\|_{\cYb^j} 
	\les_j \|\td F_{\th_n, R_n} \|_{ \cZZ_{R_n}^j } 
    \les_j \| \td F_{\th_n, R_n}\|_{\cYb^j}  , 
	\quad 
		 \|\phi_L \td F_{\th_n, R_n}\|_{\cYQ^j}  \les_j \|\td F_{\th_n, R_n} \|_{ \cZZ_{\Lams, R_n}^j } .
	\eeq

Note that \eqref{eq:LWP_uniform_Hj} is uniform in $\th, R$ for $\th \leq R^{-1}, R \es \geq 1$. 
A subsequence of $\phi_L \td F_{\th_n, R_n} $ converges weakly to some limit in $\td F_L \in \cYb^k$. 
Note that $\td F_L = \td F_{L^{\prime}}$ on the ball of radius $\min\{L, L^{\prime}\}$ and $ s \leq T$. Using a diagonalization argument, we can take $L \to \infty$ and extract a subsequence such that 
$\td F_{\th_{n_i}, R_{n_i}}  \rightharpoonup  \td F$ weakly in $\cYb^k$ for any $s \leq T$
and compact sets, and $\td F_{\th_{n_i}, R_{n_i}}  \rightharpoonup  \td F$ weakly in $L^2( (0, T), \cYQ^k )$ on compact sets. From \eqref{eq:LWP_gron_Hj1} and \eqref{eq:LWP_uniform_Hj}, 
we obtain $\td F \in \bbY^k $ (see \eqref{eq:LWP_map}) and it satisfies the  energy estimates %
\bseq\label{eq:LWP_gron_Hj}
\beq\label{eq:LWP_gron_Hj:a}
\bal 
	 \| \td F(s) \|_{\cYb^j}^2 +\int_0^{s} \f{1}{\e_{\tau}}
	 \| \td F \|_{ \cYQ^j}^2
    & \leq C_j \operatorname e^{ C_j s } \B( \| \td F( 0 ) \|_{\cYb^j}^2 
    +  \int_0^{s}  
   \f{  1 }{\e_{\tau}} (  \one_{j > \kkl}  \| \td F\|_{ \cYQ^{  j-1 }}^2  + 1)  \| \td G\|_{ \cYb^{  j}}^2(\tau) \\
  & \qquad  + %
  \f{1}{\e_0}  \| \td G\|_{L^{2}(0, s; \cYb^{j} ) }^{ \f{4}{\g+2}}
\| \td G\|_{L^{2}( 0, s; \cYQ^{j} ) }^{ \f{2\g}{\g+2}}  
  +  s   \B) , 
     \eal 
\eeq
for any $s \leq T $ and $j \leq k$. %
The $\td F$ term on the right hand side vanishes when $j = \kkl$.

Recall the $\bbY^{\kkl}$ norm from \eqref{eq:LWP_map}. %
For $ j = \kkl$ and $ \| \td G\|_{\bbY^{\kkl}} < \ze_1 < 1$, since $\es \leq \e_0 \leq 1, \f{2}{2+\g} \in [0, 1]$, we bound the norm of $\td G$ using $\| \td G\|_{\bbY^{\kkl}}$:
\[
\int_0^{s} \f{1}{\e_{\tau}} \| \td G\|_{\cYb^{\kkl}}^2
\les \f{ C \operatorname e^{C s} }{\e_0} s, 
\quad \| \td G\|_{L^{2}(0, s; \cYb^{ \kkl} ) }^{ \f{4}{\g+2}}
\| \td G\|_{L^{2}( 0, s; \cYQ^{ \kkl} ) }^{ \f{2\g}{\g+2}}  
\leq  s^{ \f{2}{\g+2}} 
\| \td G\|_{\bbY^{\kkl}}^2  \les  s^{ \f{2}{\g+2}} , 
\]
and then take supremum over $s \leq T \leq 1$ to yield 
\begin{align}\label{eq:LWP_gron_Hj:b}
	 \| \td F\|_{\bbY^{\kkl}}^2   = \sup_{s \leq T} \| \td F(s) \|_{\cYb^{\kkl}}^2 +\int_0^{s} \f{1}{\e_{\tau}}
	 \| \td F(\tau) \|_{ \cYQ^{\kkl}}^2 
	  & \leq    C \operatorname e^{ C T}  
     \big( \|\td F( 0 )  \|_{\cYb^{\kkl} }^2 + 
     T + ( T^{ \f{2}{2 + \g}} + T) \f{1}{\e_0} 
     \big)  \notag  \\ 
     & \leq 
       \bar C_1 
     \big( \|\td F( 0 )  \|_{\cYb^{\kkl} }^2 + 
      \f{1}{\e_0} T^{ \f{2}{2 + \g}}    
     \big)  ,  
\end{align}
\eseq
for some absolute constant $\bar C_1$. 

Thus, to ensure that the map $\TTs$ satisfies the property $\TTs : \bbJ_{\ze_1}^k \to \bbJ_{\ze_1}^k$, 
for some $\ze_2, T$ with $ \ze_2 < \ze_1, T\leq 1$ determined in Section \ref{sec:LWP_para},  we first impose 
\beq\label{eq:LWP_para_ineq1}
  \| \td F( 0)  \|_{\cYb^{\kkl}} < \ze_2 ,\quad 
   \bar C_1   \big( \ze_2^2 +  
       T^{ \f{2}{2 + \g}}  \f{1}{\e_0}  ) < \f{\ze_1^2}{4}. 
\eeq 
From the above estimates, we obtain 
\beq\label{eq:LWP_onto_pf}
	 \| \td F\|_{\bbY^{\kkl}} =  \| \TTs(\td G) \|_{\bbY^{\kkl}} < \tf{1}{2} \ze_1 , \quad \forall  \ 
     \| \td G \|_{\bbY^{\kkl}}  < \ze_1. 
\eeq
Estimates \eqref{eq:LWP_gron_Hj} and \eqref{eq:LWP_onto_pf} imply the property $\TTs : \bbJ_{\ze_1}^k \to \bbJ_{\ze_1}^k$.

\subsection{Contraction estimates and local existence}\label{sec:LWP_exist}

In this section, we first establish the contraction estimates and then choose small $\ze, T$ so that  
 the map $\TTs$ is contraction in $\bbY^{\kkl}$.

Suppose $ G_1,  G_2 \in \bbJ^{k}_{\ze_1}$. Let $\td F_i$ be the solution to 
 \eqref{eqn: linear-transport_1} associated with $G_i$: $\td F_i = \TTs G_i$. Denote 
 \[
 	 \td F_{\D} = \td F_1 - \td F_2,  \quad G_{\D} = G_1 - G_2, 
 	 \quad \cN_{\D} = \cN( G_1, \td F_1 ) - \cN(G_2, \td F_2). 
 \]
 Since the error term $\cE_{\cM}$ and $\cK_{\kk}$ term in $\td H$ in \eqref{eqn: linear-transport_1} does not depend on $G_i$, the $\cN_i$-operators are bilinear, we obtain the following equation for $\td F_{\D} = \td F_1 - \td F_2$:
 \[
\bal 
    \left(\pa _s + \cT + \dcm - \f32 \bcv \right) \td F_{\D} = \f1{\es} \left[(\cN _1 + \cN _5)(\rhos \cMM ^{1/2}, \td F_{\D}) \right] + \f1{\es} \cN_{\D}  + \td H_{\D},
\eal 
\]
where $ \td H_{\D}$ and $ \cN_{\D}$ are defined as 
\[
\bal 
   \td H_{\D} & := \f1\es  (\cN _2 + \cN _3 + \cN _4 + \cN _6)(\rhos \cMM ^{1/2},  G_{\D})  + \f1\es \cN ( G_{\D}, \rhos \cMM ^{1/2})  , \\ 
   	 \cN_{\D} & := \cN(G_1, \td F_1) - \cN(G_2, \td F_2) =  \cN( G_1 - G_2, \td F_1 ) + \cN(G_2, \td F_1 - \td F_2) \\
     & = \cN(G_{\D}, \td F_1) + \cN(G_2, \td F_{\D}).
    \eal 
\]

Below, we perform $\cYb^{\kkl}, \cZZ_{\infty}^{\kkl}$ energy estimates on $\td F_{\D}$. 
We bound $\td F_{\D}$ using $\cZZ, \cZZ_{\Lam}$-norms \eqref{norm:cZZ} and bound $ \td F_i, G_i, G_{\D}$ using $\cYb, \cYQ$ norms.

Applying Theorem \ref{theo: N nonlin_est} with $k= \kkl$,  we obtain
\[
	\B| \la \cN(G_{\D}, \td F_1) + \cN(G_2, \td F_{\D}), \td F_{\D} \ra_{\cYb^{\kkl}} \B|
	\leq  C \| G_{\D} \|_{\cYb^{\kkl}} \| \td F_1 \|_{\cYQ^{\kkl}}
	\| \td F_{\D} \|_{\cYQ^{\kkl}}
	+ \cnon \| G_2 \|_{\cYb^{\kkl}}  	\| \td F_{\D} \|_{\cYQ^{\kkl}}^2 .
\]
Since the $\cZZ_{\infty}^j, \cZZQinf^j$ norms \eqref{norm:cZZ} are the linear combinations of $\cYb^j$ norms and are equivalent to 
$\cYb^j, \cYQ^j$ norms, 
we further obtain 
\[
	| \la \cN(G_{\D}, \td F_1) + \cN(G_2, \td F_{\D}), \td F_{\D} \ra_{\cZZ_{\infty}^{\kkl}} |
\leq C \| G_{\D} \|_{\cYb^{\kkl}} \| \td F_1 \|_{\cYQ^{\kkl}}
	\| \td F_{\D} \|_{\cZZQinf^{\kkl}}
	+ \cnon \| G_2 \|_{\cYb^{\kkl}}  	\| \td F_{\D} \|_{\cZZQinf^{\kkl}}^2 .
\]

The linear terms satisfy the same estimates as those in \eqref{eq:LWP_EE_Hj2} with $R = \infty$ (the estimates in the whole space) 
without the error terms, nonlinear terms and the $\sss$-forcing terms.
Combining the linear estimates and the above nonlinear estimates, for $j = \kkl$, we obtain
\[
\bal
	 \f{1}{2} \f{d}{d s} \| \td F_{\D}\|_{\cZZ_{\infty}^{\kkl}}^2 
	& \leq C \| \td F_{\D} \|_{\cZZ_{\infty}^{\kkl} }^2
	 -  \f{ 1 }{ 2 \es}   \| \td F_{\D} \|_{\cZZQinf^{\kkl} }^2   
     + 
     \f{C}{\es}  \big( \|  G_{\D} \|_{\cYb^{\kkl} }^2 
     +       \|  G_{\D} \| _{\cYb ^{\kkl}}^{ \f{4}{\g+2}} \| G_{\D} \| _{\cYQ ^{\kkl}}^{ \f{2\g}{\g+2}}   \big)    
     \\
 	& \qquad  + C \es^{-1} \| G_{\D} \|_{\cYb^{\kkl}} \| \td F_1 \|_{\cYQ^{\kkl}}
	\| \td F_{\D} \|_{\cZZQinf^{\kkl}}
	+ \cnon \es^{-1} \| G_2 \|_{\cYb^{\kkl}}  	\| \td F_{\D} \|_{\cZZQinf^{\kkl}}^2.
  \eal
\]
Since $\cnon \| G_2 \|_{\cYb^{\kkl}} \leq \cnon \ze_1 < \f{ 1}{8}$ by \eqref{def:LWP_ze1}, using Young's 
 inequality similar to \eqref{eq:LWP_Young}, we bound 
\[
\bal
		 \f{1}{2} \f{d}{d s} \| \td F_{\D}\|_{\cZZ_{\infty}^{\kkl}}^2 
	 \leq C \| \td F_{\D} \|_{\cZZ_{\infty}^{\kkl} }^2
	 -  \f{ 1 }{ 4 \es}   \| \td F_{\D} \|_{\cZZQinf^{\kkl} }^2   +    \f{C}{\es}  \big( \|  G_{\D} \|_{\cYb^{\kkl} }^2 
     +       \|  G_{\D} \| _{\cYb ^{\kkl}}^{ \f{4}{\g+2}} \| G_{\D} \| _{\cYQ ^{\kkl}}^{ \f{2\g}{\g+2}}   \big)    
     +  \f{C}{\es} \| G_{\D} \|_{\cYb^{\kkl}}^2 \| \td F_1 \|_{\cYQ^{\kkl}}^2. 
	\eal
\]
Note that $\td F_1(0) = \td F_2(0) = \td F(0)$, we obtain $\td F_{\D}(0) = 0$. 
Applying Gr\"onwall's inequality similar to \eqref{eq:LWP_gron2}-\eqref{eq:LWP_gron_Hj1}, we obtain 
\[
\bal 
\| \td F_{\D}\|_{\cZZ_{\infty}^{\kkl}}^2 
+ \int_0^s  \f{1}{\e_{\tau}} \| \td F_{\D}\|_{\cZZQinf^{\kkl}} d \tau 
& \leq  
C \operatorname e^{ C s } \B(   %
  \f{1}{\e_0} s^{ \f{2}{\g+2}} \|  G_{\D}\|_{L^{\infty}(0, s; \cYb^{ \kkl} ) }^{ \f{4}{\g+2}}
\| G_{\D}\|_{L^{2}( 0, s; \cYb^{ \kkl} ) }^{ \f{2\g}{\g+2}} \\
& \qquad +    \int_0^s  \f{1}{\e_{\tau}} (  \| \td F_1\|_{ \cYQ^{\kkl}}^2 + 1 )  \|  G_{\D}(\tau) \|_{\cYb^{\kkl} }^2  \B) .  
\eal 
\]
Since $\| f \|_{\cZZ_{\infty}^{\kkl}} \asymp \| f \|_{ \cYb^{\kkl}}$,  
$\| f \|_{\cZZQinf^{\kkl}} \asymp \| f \|_{ \cYQ^{\kkl}}$ for any $f$, bounding the above upper bounds by 
$ \| \td G\|_{\bbY^{\kkl}}$ (see \eqref{eq:LWP_map}), 
and using $\f{2}{2+\g} \leq 1$, for $s \leq T \leq 1$, we establish 
 \[
\| \td F_{\D}(s) \|_{\cYb^{\kkl}}^2 
+ \int_0^s  \f{1}{\e_{\tau}} \| \td F_{\D}\|_{\cYQ^{\kkl}} d \tau  \leq C \B( \f{1}{\e_0} s^{\f{2}{\g+2}} + \int_0^s  \f{1}{\e_{\tau}} 
 \| \td F_1\|_{ \cYQ^{\kkl}}^2  \B) \| G_{\D} \|_{\bbY^{\kkl}}^2 .
 \]
Taking supremum over $s \in [0, T]$ and using estimate 
\eqref{eq:LWP_gron_Hj:b}, \eqref{eq:LWP_onto_pf}  on $\td F_1$, we further prove 
\beq\label{eq:LWP_contra1}
\bal
\| \td F_{\D}(s) \|_{\bbY^{\kkl}}^2
 &%
\leq \bar C_2  \B( \| \td F_1(0) \|_{\cYb^{\kkl} }^2 + 
 \f{1}{\e_0} T^{ \f{2}{2 + \g}}  \B) \cdot \| G_{\D} \|^2_{\bbY^{\kkl}}   , 
\eal
\eeq
for some absolute constant $\bar C_2$ independent of $T, \e_0$.

\subsubsection{ Choosing $\ze_2, T$}\label{sec:LWP_para}

Recall $\ze_1$ from \eqref{def:LWP_ze1}, and constants $\bar C_1, \bar C_2$ from  \eqref{eq:LWP_gron_Hj:b} and \eqref{eq:LWP_contra1}. We choose $\ze_2$ and $T>0$ as:
\bseq\label{def:ze2}
\beq
\ze_2 = \bar C_3 \ze_1, \quad T(\e_0)^{ \f{2}{2+\g}} = \min \{\bar C_3^2 \e_0 \ze_1^2, \, 1\} >0,
\eeq
with some absolute constant $\bar C_3$ small enough such that
\beq
\bar C_3 < 1, \quad 
   ( \bar C_1  + \bar C_2   ) ( \ze_2^2 + \f{1}{\e_0} T^{ \f{2}{2 + \g} } ) 
\leq 2 \bar C_3^2 \ze_1^2 
   ( \bar C_1  + \bar C_2 ) 
   <\min \left\{ \f{1}{4} \ze_1^2, \ \f{1}{4} \right\}. 
\eeq
\eseq

For initial data $\td F(0) \in \cYb^k$ with $\|\td F(0) \|_{\cYb^{\kkl}} < \ze_2$, the above parameters imply \eqref{eq:LWP_para_ineq1}
and the estimates \eqref{eq:LWP_onto_pf}. 
For any $G_1, G_2 \in \bbJ^{k}_{\ze_1}$, estimates \eqref{eq:LWP_contra1} and \eqref{def:ze2} imply
\[
	\| \TTs(G_1) - \TTs(G_2) \|_{\bbY^{\kkl}}^2
	\leq 
 \bar C_2   ( \| \td F(0) \|_{ \cYb^{\kkl} }^2 +  \f{1}{\e_0} T^{ \f{2}{2 + \g}}  )  \| G_1 - G_2 \|_{\bbY^{\kkl}}^2 < \f{1}{4}
	\| G_1 - G_2 \|_{\bbY^{\kkl}}^2 .
\]
Thus $\TTs$ is a contraction mapping in $\bbJ^{\kkl}_{\ze_1}$. 

Using the Banach fixed point theorem with the map $\TTs$ in the space $\bbJ^{\kkl}_{\ze_1}$, we construct a unique fixed point  $\td F = \TTs(\td F)$ , which solves the nonlinear equations 
\eqref{eq:fix_F_1}. We prove the existence and uniqueness of local-in-time solution in Theorem \ref{thm:LWP}.
Moreover, since $\td F$ satisfies the estimates \eqref{eq:LWP_gron_Hj:a} with $j = \kkl$ for $s \leq T<1$ and 
$ \td F\in \bbJ^{\kkl}_{\ze_1}$, we prove \eqref{eq:LWP_thm_EE}.

\paragraph{Energy estimates}

Since $\td F = \TTs(\td F)$, it satisfies the energy estimates in 
\eqref{eq:LWP_gron_Hj} with $\td G= \td F$ provided $ \td F \in \bbJ_{\ze_1}^k$. In particular, for $ j \geq \kkl, j \leq k $, using 
\eqref{eq:LWP_gron_Hj:a}, we obtain 
\beq\label{eq:LWP_EE:LHS}
\bal 
	 \| \td F(s) \|_{\cYb^j}^2 +\int_0^{s} \f{1}{\e_{\tau}}
	 \| \td F \|_{ \cYQ^j}^2
    & \leq C_j \operatorname e^{ C_j s } \B( \| \td F( 0 ) \|_{\cYb^j}^2 
    +  \int_0^{s}  
   \f{  1 }{\e_{\tau}} (  \one_{j > \kkl}  \| \td F\|_{ \cYQ^{  j-1 }}^2  + 1)  \| \td F\|_{ \cYb^{  j}}^2(\tau) \\
  & \hspace{6em}  + %
  \f{1}{\e_0}  \| \td F\|_{L^{2}(0, s; \cYb^{j} ) }^{ \f{4}{\g+2}}
\| \td F\|_{L^{2}( 0, s; \cYQ^{j} ) }^{ \f{2\g}{\g+2}}  
  +  s   \B) , 
\eal 
\eeq
Applying $\e$-Young's inequality to $ \| \td F\|_{L^{2}(0, s; \cYb^{j} ) }^{ \f{4}{\g+2}}
\| \td F\|_{L^{2}( 0, s; \cYQ^{j} ) }^{ \f{2\g}{\g+2}} $ and using $  \e_0 \operatorname e^{- Cs}  \les \es  $, we obtain 
\[
 C_j \operatorname e^{C_j s} \f{1}{\e_0} 
  \| \td F\|_{L^{2}(0, s; \cYb^{j} ) }^{ \f{4}{\g+2}}
\| \td F\|_{L^{2}( 0, s; \cYQ^{j} ) }^{ \f{2\g}{\g+2}}  
\leq \f12 \int_0^{s} \f{1}{\e_{\tau}}
	 \| \td F(\tau) \|_{ \cYQ^j}^2
     +  C_j \operatorname e^{C_j s}  \int_0^s \f{1}{ \e_{\tau}} \| \td F(\tau) \|_{\cYb^{j}}^2 ,     
\]
where $C_j$ may change from line to line. We absorb the first term on the right hand side using 
the dissipation on the left hand side of \eqref{eq:LWP_EE:LHS}. Combining the above two estimates, we prove 
\beq\label{eq:LWP_EE_non0}
\bal 
	 \| \td F(s) \|_{\cYb^j}^2 +\int_0^{s} \f{1}{\e_{\tau}}
	 \| \td F \|_{ \cYQ^j}^2
    & \leq C_j \operatorname e^{ C_j s } \B( \| \td F( 0 ) \|_{\cYb^j}^2 
    +  \int_0^{s}  
   \f{  1 }{\e_{\tau}} (  \one_{j > \kkl}  \| \td F\|_{ \cYQ^{  j-1 }}^2  + 1)  \| \td F\|_{ \cYb^{  j}}^2(\tau)  
   +   s \B)  ,
  \eal 
\eeq 
whenever $\td F \in L^{\infty}( 0, s ;  \cYb^k)  \cap L^2(0, s; \cYQ^k)$.

\subsection{Continuation criterion}\label{sec:LWP_continue}

Below, we show that if
\beq\label{eq:LWP_blowup2}
   \sup_{s \in [ 0,  T_*)} \| \td F(s)\|_{\cYb^{\kkl}} < \ze_2 ,
\eeq
we can extend the solution to \eqref{eq:fix_F_1} in $L^{\infty}( ( 0,   T_2) , \cYb^{k}  )
\cap L^{2}( ( 0 ,  T_2) , \cYQ^{k}  )$ for some $T_2> T_*$.

Suppose \eqref{eq:LWP_blowup2} holds true. Firstly, we show that $\td F \in \bbY^k$. 
Since $\ze_2 < \ze_1$,  $\td F$ satisfies the energy estimates in \eqref{eq:LWP_EE_non0}.
For $j = \kkl$, since $\one_{j > \kkl } = 0$, using \eqref{eq:LWP_blowup2} and \eqref{eq:LWP_EE_non0}, we bound $\| \td F \|_{\bbY^{\kkl}}$ in terms of  $\td F(0)$ and $s$ 
 \bseq\label{eq:LWP_EE_non}
 \beq\label{eq:LWP_EE_non:a}
	 \| \td F(s) \|_{\cYb^j}^2 +\int_0^{s} \f{1}{\e_{\tau}}
	 \| \td F \|_{ \cYQ^j}^2
\leq C \operatorname e^{C s}  ( \| \td F( 0 ) \|_{\cYb^j}^2  + \e_0^{-1} s   ). 
 \eeq  
For $j > \kkl$, applying Gr\"onwall's inequality to  $\| \td F(s) \|_{\cYb^j}^2$, we obtain 
\beq
			 \| \td F(s) \|_{\cYb^j}^2  +\int_0^{s} \f{1}{\e_{\tau}}
	 \| \td F(\tau) \|_{ \cYQ^j}^2
			 \leq 
			 C_j \exp\B(  C_j \operatorname e^{C_j s} \int_0^s  \f{1}{\e_{\tau}} \| (\td F\|_{\cYQ^{j-1}}^2 + 1 ) d \tau  \B)
			 \cdot  \B(  \| \td F( 0 ) \|_{\cYb^j}^2  +  s  \B) .
\eeq
For $j=\kkl + 1$, the right hand side is uniformly bounded for $s <  T_*$ due to 
the estimate \eqref{eq:LWP_EE_non:a} with $j = \kkl$. 
Applying the above estimates inductively on $j$, we establish the uniform boundedness 
for any $ j\leq k$ :
\beq
	 \sup_{s \in ( 0,  T_*) } \| \td F(s) \|_{\cYb^j}^2  +\int_0^{s} \f{1}{\e_{\tau}}
	 \| \td F(\tau) \|_{ \cYQ^j}^2 
\leq C(j , T_*, \| \td F(0) \|_{\cYb^j}  ) .
\eeq
\eseq

Since $\td F(s_1) \in \cYb^k$ and $\| \td F(s_1) \|_{\cYb^{\kkl}} < \ze_2$ for any $s_1 <  T_*$, we  apply the fixed-point construction in previous sections with initial data $\td F(s_1)$ and extend the solution in $ \bbJ_{\ze_1}^k$  to $s \in [s_1, s_1 + T(s_1))$ with $T(s_1)
\geq \bar C_4 \min \{ \e_{s_1}, 1 \} ^\f{\g + 1}2$ chosen in \eqref{def:ze2}, where $\bar C_4$ is some absolute constant. Note that we need to change $\e_0$ in \eqref{def:ze2} to $\e_{s_1}$ due to the change of initial time.
Since 
\[
	\lim_{ s_1 \to (  T_*)^- }  s_1 + \bar C_4 \min \{ \e_{s_1}, 1 \}^{  \f{\g+1}{2}  } >  T_* , 
\]
choosing $s_1$ close to $ T_*$ so that $s_1 + T(s_1) >  T_*$, we extend the solution beyond $  T_*$. We complete the proof of the continuation criterion. We complete the proof of Theorem \ref{thm:LWP}.

\subsection{Local existence of solution to the Landau equation}\label{sec:LWP_landau}

Consider $F =  \cM + \cMM^{1/2} \td F$ with initial data  $\td F(0)$ satisfying \eqref{eq:LWP_ass_IC} and $F(0, X, V)> 0$. Using Theorem \ref{thm:LWP} with $\sss = 0$, we construct a local solution 
$ \td F$ to \eqref{eq:fix_F_1} with $\sss = 0$. 
Equation \eqref{eq:fix_F_1} with $\sss=0$ for $\td F$ is equivalent to \eqref{eq:LC_ss} 
for $F =\cM + \cMM^{1/2} \td F $. We obtain a local solution to \eqref{eq:LC_ss} with 
$\td F \in L^{\infty}( [ 0,  T], \cYb^{k} ) \cap L^2( [ 0,  T], \cYQ^{k} )$. Moreover, the solution  
satisfies \eqref{eq:LWP_thm_EE} and the continuation criterion \eqref{eq:LWP_blowup}.

\paragraph{Gaussian lower bound}
Below, we prove the positivity of $F$. Since $\ze_1 < \ze_0$, using $F = \cM + \cMM^{1/2} \td F$, Lemma \ref{lem:parab}, 
and Lemma \ref{lem: A-pointwise-bound}, we obtain
\beq\label{eq:LWP_A_lower}
C_1 \bS
\preceq \f{1}{2} A[\cM] \preceq  A[F] 
\preceq  C_2 \bS
\eeq
for some $C_1, C_2 > 0$, where $\bS$ is defined in \eqref{eqn: defn-matrix}. Using \eqref{eqn: divA-pointwise} with $N=1,i=2, f = \brho \cMM^{1/2} + \td F$, the estimate \eqref{eqn: smallness-G-L2}, and $\ze_0 < 1$, we obtain 
\beq\label{eq:LWP_A_up}
  | \na_V^2  A[\cM + \cMM^{1/2} \td F | 
  \les \cs^{\g} \ang \vc^{\g}  \|  \brho \cMM^{1/2} + \td F  \|_{L ^2 (V)}
  \les \cs^{\g} \ang \vc^{\g} ( \cs^{3} + \cs^3  )
  \les \cs^{\g+3} \ang \vc^{\g}.
\eeq

We can treat \eqref{eq:LC_ss}, \eqref{eq:Q} as a linear parabolic equation with the following operator
\beq\label{eq:LWP_max}
\bal
	L_F g &:= (\pa_s + \cT) g 
	- \es^{-1} ( A(F) : \na_V^2 g - \div_V ( \div_V A(F) )   \cdot  g) , \\
	 \cT g & = (\bcx X \cdot \na_X + \bcv V \cdot \na_X + c_v V \cdot \na_V) g.
\eal 
\eeq
Below, we construct a barrier function $H > 0$  and show that $F - C H \geq 0$ for some $C \geq 0$.

\paragraph{Barrier function}

Recall $\cM = \mu(\vc) = \exp(-\kp_2 |\vc|^2)$ and $\cMM$ form \eqref{eq:localmax2} and \eqref{eq:gauss}. Let $ l \in [0, 100], \bbb >0, \aaa \geq 1$ be parameters to be chosen.
\footnote{We impose an upper bound for $l$ so that we do not need to track constant depending on $l$.} 
We construct the following barrier function 
\beq\label{eq:barrier1}
	 H_{l, \aaa, \bbb} = \operatorname e^{ - \bbb \cdot \es^{-1}  }  \ang X^{-l} \cM^{\aaa}
	 = \operatorname e^{- \bbb \cdot \es^{-1} }  \ang X^{-l} \exp( - \aaa \kp_2 |\vc|^2 ),
\eeq
where we recall $\es = \e_0 \operatorname e^{-\rE s}$ from \eqref{eq:EE_para1}. 

Next, we show that $L_F H \leq  -\f{\bbb}{2} \rE \es^{-1} H_{l, \aaa, \bbb} $ for some $\aaa, \bbb$ large enough.  A direct calculation yields 
\[
\bal 
	L_F H_{l,\aaa, \bbb}  & : = \udb{\operatorname e^{-\bbb \es^{-1}} \ang X^{-l} (\pa_s + \cT) \cM^{\aaa} }_{:=I_1} + \udb{ \ang X^{-l} \cM^{\aaa} \pa_s \operatorname e^{-\bbb \es^{-1}} }_{:=I_2} +  \udb{ \operatorname e^{-\bbb \es^{-1}} \cM^{\aaa}  \cT \ang X^{-l}  }_{:=I_3} \\ 
	  & \quad  \udb{   - \es^{-1} \operatorname e^{-\bbb \es^{-1}} \ang X^{-l} ( A[F] : \na^2_{V} \cM^{\aaa}
	- \div_V ( \div_V A(F) )  \cdot \cM^{\aaa} ) }_{:= II}.
\eal 
\]
Recall the function class $\FF^{-r}$ from Definition \ref{def: bF}. Using \eqref{eq:error} and Lemma \ref{lem:ds_t}, we obtain
\[
\bal
|I_1| & =\aaa | \operatorname e^{-\bbb \es^{-1}} \ang X^{-l}  \cM^{\aaa-1} (\pa_s + \cT) \cM|  
= \aaa | \operatorname e^{-\bbb \es^{-1}} \ang X^{-l}  \cM^{\aaa-1} \eM|  \\
& 	\les \aaa \operatorname e^{-\bbb \es^{-1}} \ang X^{-l}  \ang X^{-r} \ang \vc^{3}  \cM^{\aaa}
\les \aaa \ang X^{-r} \ang \vc^{3}   H_{l, a, b} .
\eal
\]
For $I_2$, from \eqref{eq:EE_para1}, we obtain $\pa_s \es^{-1} =  \rE \es^{-1}$ and 
\[
	I_2 = \ang X^{-l} \cM^{\aaa} \cdot (-\bbb \rE \es^{-1} ) \operatorname e^{ - \bbb \es^{-1} } 
	=  -\bbb \rE \es^{-1}  H_{l, a, b}.
\]
For $I_3$, using $\vc = \f{V -\bu}{\cs}$, $|V| \les \cs \ang \vc$, 
and $l \leq 100$, we obtain
\[
\bal
		|I_3| & =| \operatorname e^{-\bbb \es^{-1}} \cM^{\aaa} (  \bcx X \cdot \na_X + V \cdot \na_X ) \ang X^{-l} | 
		\les  | \operatorname e^{-\bbb \es^{-1}} \cM^{\aaa}  
			  ( \ang X^{-l} + |V| \ang X^{-l-1} ) | \\
& \les   | \operatorname e^{-\bbb \es^{-1}} \cM^{\aaa}  \ang X^{-l}
			  ( 1 +  \cs \ang X^{-1}   \ang \vc  ) |
			  \les   ( 1 +  \cs \ang X^{-1}   \ang \vc )   H_{l, a, b}. 
	\eal
\]
For the collision part, since $\cM^{\aaa} = \exp(-\kmu \aaa |\vc|^2) = \exp( -\kmu \aaa \f{|V-\bu|^2}{\cs^2} )$, we yield 
\[
	\pa_{V_i V_j } \cM^{\aaa}
	= 4 \aaa^2 \kmu^2 \cs^{-2} \vc_i \vc_j  \cM^{\aaa}
	- 2 \aaa \kmu  \cs^{-2} \d_{ij} \cM^{\aaa}.
\]
Using \eqref{eq:LWP_A_lower}, \eqref{eq:LWP_A_up}, $\aaa \geq 1$, and the above calculation 
on $\na_V^2 \cM^a$, we yield 
\[
\bal
 		A[F]: \na^2_{V} \cM^{\aaa} 
		 & - \div_V ( \div_V A[F]) \cdot \cM^{\aaa} \\ 
	 &  	\geq \B( \cs^{\g+5} \cdot C_1 \aaa^2 \cs^{-2} \ang \vc^{\g} |\vc|^2 - \cs^{\g+5}\cdot  C_2 \aaa \cs^{-2} \ang \vc^{\g+2}
		-C_3 \cs^{\g+3} \ang \vc^{\g} \B) \cM^{\aaa} \\ 
		&   \geq \cs^{\g+3} (C_1 \aaa^2 |\vc|^2 \ang \vc^{\g} - C_4 \aaa  \ang \vc^{\g+2}   ) \cM^{\aaa} ,
\eal
\]
for some absolute constant $C_i$. Thus, we estimate $II$ as 
\[
\bal
	II & \leq - \es^{-1} \operatorname e^{-\bbb \es^{-1}} \ang X^{-l}  
	\cs^{\g+3} (C_1 \aaa^2 |\vc|^2 \ang \vc^{\g} - C_4 \aaa  \ang \vc^{\g+2}   ) \cM^{\aaa} \\
	& = \cs^{\g+3} \es^{-1} (- C_1 \aaa^2 |\vc|^2 \ang \vc^{\g} + C_4 \aaa  \ang \vc^{\g+2}   ) H_{l, \aaa, \bbb}.
\eal
\]

Combining the above estimates, we prove
\[
\bal 
	L_F H_{l,\aaa,\bbb}  \leq &  \B( C  + C   \cs \ang X^{-1} \ang \vc + C \aaa \ang X^{-r} \ang \vc^3
	- \bbb \rE \es^{-1} \\ 
    & \quad +\es^{-1} \cs^{\g+3} (-  C_1 \aaa^2 |\vc|^2 \ang \vc^{\g} + C_4 \aaa  \ang \vc^{\g+2}   ) \B)  H_{l,\aaa, \bbb } .
\eal 
\]
Using Lemma \ref{lem:V/X}, we obtain
\[
		L_F H_{l,\aaa,\bbb} 
		\leq  \B( C_5 + C_5 \cs^{\g+3} \ang \vc^{\g+2} 
- \bbb \rE \es^{-1} +\es^{-1} \cs^{\g+3} (-  C_1 \aaa^2 |\vc|^2 \ang \vc^{\g} + C_4 \aaa  \ang \vc^{\g+2}   ) \B)  H_{l,\aaa, \bbb } .
\]

There exist absolute constants $ \aaa_*, \bar C$ large enough, such that for any 
\beq\label{eq:LWP_ab_range}
	\aaa \geq \aaa_* > \max \left\{ \f{ 2 C_4}{C_1}, 1 \right\},  \quad \bbb \geq \bar C \aaa^2, 
\eeq
using $\es^{-1} \gtr 1$ and $|\vc|^2 \geq \f{1}{2} \ang \vc^2$ for $|\vc| \gtr 1$, we obtain
\[
\bal
C_5  + \es^{-1} ( -C_1 \aaa^2 |\vc|^2 \ang \vc^{\g} +  C_4 \aaa  \ang \vc^{\g+2} )
 & \leq \es^{-1} (- \f{C_1}{2} \aaa^2 + C_4 \aaa )\ang \vc^{\g+2} 
 +\f{1}{3} \bbb \rE \es^{-1}
\leq \f{1}{3} \bbb \rE \es^{-1}.
\eal
\]
By further requiring $\aaa_*$ large in \eqref{eq:LWP_ab_range}, and using $\cs \les 1$, 
we prove 
\beq\label{eq:barrier2}
	L_F H_{l,\aaa,\bbb} \leq( C_5  + \f13 \bbb \rE \es^{-1} -  \bbb \rE \es^{-1} )   H_{l,\aaa,\bbb} 
\leq -\f{1}{2} \bbb \rE \es^{-1} H_{l,\aaa,\bbb}.
\eeq
uniformly for any $\aaa, \bbb$ satisfying \eqref{eq:LWP_ab_range}, and $l \in [0, 100]$. 
Using $L_F F = 0$ and \eqref{eq:barrier2}, for any $ C \geq 0$, we obtain 
\beq\label{eq:barrier3}
	L_F( F - C H_{l,\aaa,\bbb} ) \geq \tf{1}{2} \bbb  C \rE \es^{-1}.
\eeq

\paragraph{Decay at infinity}
Using \eqref{eq:linf_Y_point} with $\eta \rsa \etab = -3 + 6(r-1)$, \eqref{eq:LWP_thm_EE} and $\kkl \geq 2 d$,  we obtain
\[
\bal
	F(s, X, V) & \geq \cM - |\td F| \cMM^{1/2}
	\geq \mu(\vc) - C |\td F| \cs^{-3/2} \mu(\vc)^{1/2} \\
	& \geq - C \cs^{-3} \ang X^{- \f{\etab+3}{2} } \mu(\vc)^{1/2} \| \td F\|_{\cYb^{\kkl}}
	\geq - C \cs^{-3} \ang X^{-3(r-1)} \mu(\vc)^{1/2} . 
\eal
\]
For $H_{l,\aaa,\bbb}$ \eqref{eq:barrier1}, since $\aaa \geq 1$, we have 
\[
 H_{l, \aaa ,\bbb}(s, X, V)  \leq \ang X^{-l} \mu( \vc)^{\aaa}
 \leq \ang X^{-l} \mu( \vc).
\]

When $l > 0$, since $\cs \gtr_s 1$, $ \mu(\vc) \to 0$ as $|V-\bu| \to \infty$,
and $|\bu(X)| \to 0, \ang X^{-3(r-1)} \to 0$ as $|X| \to \infty$ (see \eqref{eq:dec_U}), we obtain 
\beq\label{eq:LWP_max_BC}
	F(s, X, V) \geq - c_1(R), \quad 	|H(s, X, V) |\leq c_l(R),  \quad \mw{for \ } |(X, V)| \leq R , 
\eeq
with $c_1(R), c_l(R) >0$ and $c_1(R), c_l(R) \to 0$ as $R \to \infty$, uniformly in $s\in[ 0 ,  T]$. 
For initial data satisfying \eqref{eq:Gauss_lower_IC}, we obtain $\psi := F - c \operatorname e^{\bbb \e_0^{-1}} H_{l,\aaa, \bbb} |_{s =0} >0$. Applying the maximum principle to $\psi$ in the domain $\Om_R = \{(s,X,V): s \in [ 0,  T], |(X, V)| \leq R \}$, we prove 
\[
\psi(s, X, V)  = (F - c \operatorname e^{\bbb \e_0^{-1}} H_{l,\aaa, \bbb} )(s, X, V)  \geq -c_1(R) - c \operatorname e^{\bbb \e_0^{-1}}  \cdot c_l(R) , \quad \forall (s, X, V) \in \Om_R .
\]
Taking $R \to \infty$, we prove $\psi(s, X, V) \geq 0
$ and obtain \eqref{eq:Gauss_lower}.

When $l =0$, under the assumption \eqref{eq:Gauss_lower}, 
since $H_{l, a, b}$ is decreasing in $l$ \eqref{eq:barrier1}, we have 
\[
 F - c \operatorname e^{ \bbb \e_0^{-1}} H_{q,\aaa, \bbb} |_{s =0} 
\geq F - c \operatorname  e^{ \bbb \e_0^{-1}} H_{0,\aaa, \bbb} |_{s =0}  > 0, \quad \forall \ q > 0.
\]
Since \eqref{eq:barrier3} holds for  any $l\in [0, 100]$, applying the maximum principle to $  F - c \operatorname e^{ \bbb \e_0^{-1}} H_{q,\aaa, \bbb}$ and then taking $q\to 0^+$, we prove 
$ F - c \operatorname e^{ \bbb \e_0^{-1}} H_{0,\aaa, \bbb} \geq 0$.  We complete the proof of Proposition \ref{prop:LWP_landau}.

\appendix

\section{Derivation of the linearized Euler equations and error estimate}\label{app:euler_deri}

In this appendix, we estimate the macro-error of the profile $\erho, \eu, \ep, \ec$ defined in \eqref{eq:error0}, and derive the linearized Euler equations \eqref{eq:lin_euler2} from \eqref{eq:lin}.

\subsection{Estimate of the macro-error}\label{app:error_est}

In this section, we estimate the macro-error $\erho, \eu, \ep, \ec$ defined in \eqref{eq:error0}. First, we recall the definitions from \eqref{eq:error0}
\bseq\label{eq:error}
\beq
\bga
    \eM = (\pa_s + \bcx X \cdot \na_X + \bcx V \cdot \na_V + V \cdot \na_X) \cM, \\
    \erho = \cs^{-3} \la \eM, 1 \ra_\lvv, \quad 
    \eu   = \cs^{-4} \la \eM, V - \bu \ra_\lvv, \quad  
    \ep   = \cs^{-5} \left\la \eM, \f{1}{3} |V - \bu|^2 \right\ra _\lvv, \\
\ega
\eeq
and 
\beq\label{eq:error_a}
	   \ec =\cs^{-1} \B( [\pa_s + ( \bcx X + \bu) \cdot \na] \cs + \f13 \cs (\na \cdot \bu) -  \bcv \cs \B) .
\eeq

Using the derivations \eqref{eq:Euler1} and $\pa_s \bu = 0$, we obtain the following formulas 
\beq\label{eq:error_b}
\bal
    \erho &= \cs^{-3} \B( 
    [\pa_s + (\bcx X + \bu) \cdot \na] \rhos + \rhos (\na \cdot \UU) -3 \bcv \rhos \B) \\
    &= [\pa_s + ( \bcx X + \bu) \cdot \na] \log \rhos + (\na \cdot \bu) - 3 \bcv = 3 \ec, \\
    \eu & = \cs^{-4} ( \la \eM, V \ra_V - \bu \la \eM, 1 \ra_V ) \\ 
    & = \cs^{-1} \B((\bcx X + \bu) \cdot \na \bu - \bcv \bu  + \brho^{-1}  \na \bp \B) \\ 
    & = \cs^{-1} \B((\bcx X + \bu) \cdot \na \bu - \bcv \bu  + 3 \cs \na \cs \B), \\
    \ep & = \cs^{-5} \B( \f{1}{3} ( \la \eM, |V|^2 \ra - 2 \bu \cdot \la \eM , V \ra + |\bu|^2 \la \eM, 1\ra ) \B) \\
    & = \cs^{-5} \B( [\pa_s + ( \bcx X + \bu) \cdot \na] \bp + \kp \bp (\na \cdot \bu) - 5 \bcv \bp \B) \\
    & = \f1\kp [\pa_s + ( \bcx X + \bu) \cdot \na] \log \bp + (\na \cdot \bu) - \f5\kp \bcv = \erho = 3 \ec.  
\eal
\eeq
\eseq
Here we used $\rhos = \cs ^3$ and $P _s = \f1\kp \cs ^5$ in deducing $\erho = \ep = 3 \ec$.

Using the equation of $\bu$ from \eqref{eqn: rhoUP}, $\bar \rho = \bar {\sc}^3$ and $\bar P = \f{1}{\kp}\bc^5$ from \eqref{eq:Euler_profi_rhoPth}, we obtain
\bseq\label{eq:error_U}
\beq
\bal
    (\bcx X + \bu) \cdot \na \bu - \bcv \bu &= - \bar \rho^{-1} \na \bar P = - \bar \sc ^{-3} \cdot \na \left(\f{3}{5} \bar \sc^5\right) = - 3 \bar \sc \na \bar \sc .  \\
  \eal
\eeq
Similarly, we obtain $\brho^{-1} \na \bp = 3 \cs \na \cs$. Combining these two estimates, we obtain 
\beq\label{eq:error_Ub}
\eu  = \cs^{-1}( - 3 \bar \sc \na \bar \sc + 3 \cs \na \cs ) 
= \f{3}{2} \cs^{-1} \na( \cs^2 - \bc^2 ).
\eeq

\eseq

\begin{lemma}[Cut-off error]
    \label{lem: cutoff_error}
Let $\ec, \eu, \ep, \erho$ be defined in \eqref{eq:error0} (or \eqref{eq:error}). We have 
\beq\label{eq:error_id}
	\erho = \ep =  3 \ec. 
\eeq
    For $ k \geq 0$, we have the following estimates of $\ec$ and $\eu$
	\begin{align}
        \label{eq:error_ec}
		|\nabla ^k \ec| & \lesssim _k \ang X ^{-r - k} \one _{\{|X| \ge \rs\}}, \\
        \label{eq:error_eu}
		|\nabla ^k \eu| & \lesssim _k \ang X ^{-r - k} \one _{\{|X| \ge \rs\}}.
	\end{align}
    Recall the $\cX$-norm from \eqref{norm:Xk}. For any $k\geq0$, $\eta \leq \etab =- 3 + 6(r-1) $, and $\cE = \eu , \ec, \ep$, or $\erho$, we have 
    \beq\label{eq:error_L2}
 \|  \cs^3 \cE \|_{\cX^k_{\eta}} \les_{k, \eta}  \rs^{  \f{\eta - \etab}{2} - r  } \les_{k, \eta} \rs^{-r }.
\eeq
    
\end{lemma}

\begin{proof}[Proof of Lemma \ref{lem: cutoff_error} ]

The identity \eqref{eq:error_id} follows from \eqref{eq:error_b}.

	Note that $(\bu, \bc)$ solves \eqref{eq:Euler3} precisely. Since $\bc = \cs$ in $\{|X| \le \rs\}$, the errors are zero inside the ball $\{|X| \le \rs\}$. Subtract \eqref{eq:Euler3} from  
    \eqref{eq:error_a}, we have 
	\begin{align*}
		\cs \ec & = \left[\partial _s + (\bcx X + \bu) \cdot \nabla - \bcv + \f13 \na \cdot \bu\right] (\cs - \bc) \\
		\notag
		& = \left[\partial _s + ( \bcx X + \bu) \cdot \nabla - \bcv + \f13 \na \cdot \bu\right] [(\rs ^{-r + 1} - \bc) (1 - \chi _{\rs})] \\
		\notag
		& = (1 - \chi _{\rs}) (\partial _s - \bcv + \f13 \na \cdot \bu) \rs ^{-r + 1}  - (\rs ^{-r + 1} - \cs) \left[\partial _s + ( \bcx X + \bu) \cdot \nabla \right] \chi _{\rs}.
	\end{align*}
	Within the first term, due to our choice of cut-off function, it holds that
	\begin{align*}
		(\partial _s - \bcv) \rs ^{-r + 1} = [(-r + 1) \bcx - \bcv] \rs ^{-r + 1} = 0.
	\end{align*}
	There is also cancellation on the second term:
	\begin{align*}
		\left(\partial _s + \bcx X \cdot \nabla \right) \chi _{\rs} = \left[-\frac{X}{\rs ^2} \partial _s \rs + \frac1{\rs} \bcx X \right] \cdot \nabla \chi = 0.
	\end{align*}
	With these cancellations, we have
	\begin{align*}
		\cs \cE _{\sc} = \f13 (1 - \chi _{\rs}) (  \na \cdot \bu) \rs ^{-r + 1} - (\rs ^{-r + 1} - \cs) \bu \cdot \nabla \chi _\rs.
	\end{align*}
    To prove \eqref{eq:error_ec}, We take a multi-index $\al$, and compute 
    \begin{align}
        \label{eq:dx_cs_ec}
        \pa _X ^\al (\cs \ec) = \rs ^{-r + 1} \pa _X ^\al \left[
            \f13 (1 - \chi _{\rs}) (  \na \cdot \bu)
        \right] - \pa _X ^\al \left[(\rs ^{-r + 1} - \cs) \bu \cdot \nabla \chi _\rs
        \right].
    \end{align}
    Note that for any multi-index $\al$, it holds that
    \begin{align*}
        |\pa _X ^\al \chi _\rs| &\les \rs ^{-|\al|} |\na ^{|\al|} \chi| \les \ang X ^{-|\al|}, &        
        |\pa _X ^\al \na \cdot \bu| &\les \ang X ^{-r - |\al|}.
    \end{align*}
    We used that $X \approx \rs$ in the support of $\na \chi$.
    By Leibniz rule we know 
    \begin{align}
        \label{eq:dx_cs_ec_1}
        \pa _X ^\al \left[
            (1 - \chi _{\rs}) (  \na \cdot \bu)
        \right] \les \ang X ^{-r - |\al|}.
    \end{align}
    Similarly, because 
    \begin{align*}
        |\pa _X ^\al (\rs ^{-r + 1} - \cs)| &\les \cs \ang X ^{-|\al|}, &
        |\pa _X ^\al \bu| &\les \ang X ^{-r + 1 - |\al|}, &
        |\pa _X ^\al \na \chi _\rs| & \les \ang X ^{-|\al| - 1},
    \end{align*}
    by Leibniz rule we conclude 
    \begin{align}
        \label{eq:dx_cs_ec_2}
        |\pa _X ^\al \left[(\rs ^{-r + 1} - \cs) \bu \cdot \nabla \chi _\rs
        \right]| \les \cs \ang X ^{-r - |\al|}.
    \end{align}
    Combining \eqref{eq:dx_cs_ec_1}, \eqref{eq:dx_cs_ec_2} with \eqref{eq:dx_cs_ec} and $\rs ^{-r + 1} \les \cs$, we have shown 
    \begin{align*}
        \pa _X ^\al (\cs \ec) \les \cs \ang X ^{-r - |\al|}.
    \end{align*}
    This proves \eqref{eq:error_ec} with $k = 0$. Using Leibniz rule again and \eqref{eq:dec_S}, \eqref{eq:error_ec} follows by induction:
    \begin{align*}
        \cs |\na ^k \ec| &\les |\na ^k (\cs \ec)| + \sum _{k' < k} |\na ^{k - k'} \cs| |\na ^{k'} \ec| \\
        &\les \cs \ang X ^{-r - k} + \ang X ^{-r + 1 - (k - k')} \ang X ^{-r - k'} \les \cs \ang X ^{-r - k}.
    \end{align*}

\paragraph{Proof of \eqref{eq:error_eu} and \eqref{eq:error_L2}}
    As for \eqref{eq:error_eu}, using \eqref{eq:error_b} and \eqref{eq:error_Ub}, we directly compute 
    \begin{align*}
        |\na ^k (\cs \eu)| \les  |\na ^k (\cs \na \cs - \bc \na \bc)| \les \cs \ang X ^{-r - k} + \bc \ang X ^{-r - k} \les \cs \ang X ^{-r - k}.
    \end{align*}
    \eqref{eq:error_eu} follows by the same Leibniz rule and induction.

    Finally, for $\cE = \eu$ or $\ec$ and any $k \geq 0$, 
    by Leibniz rule and \eqref{eq:dec_S}, \eqref{eq:error_ec}, \eqref{eq:error_eu}, we obtain %
    \begin{align*}
        | \ang X^k \na ^{k} (\cs ^3 \cE)| \les _k \ang X^k \cs ^3 \ang X ^{-r-k} \one _{\{|X| \ge \rs\}} \les_k \rs ^{-3 (r - 1)} \ang X ^{-r } \one _{\{|X| \ge \rs\}}.
    \end{align*}
    Therefore, using the definition of $\cX$-norm \eqref{norm:Xk} and the above estimate, we obtain
    \begin{align*}
     \| \cs^3 \cE \|_{\cX_{\eta}^k}
    &  \les_{k, \eta}
         \| \ang X ^{\f\eta2 + k} \na ^{k} (\cs ^3 \cE) \| _{L ^2} 
         +          \| \ang X ^{\f\eta2}  (\cs ^3 \cE) \| _{L ^2}  \\ 
         & \les _{k, \eta} \rs ^{-3 (r - 1)} \B( \int _{\{|X| \ge \rs\}} \ang X ^{-2 r  + \eta} d X \B) ^\f12 \les_{k,\eta} \rs ^{-3 (r - 1) + \f{-2 r  + \eta + 3}2} 
        = \rs ^{\f{\eta - \etab}2  - r}.
    \end{align*} 
 For $\eta \leq \etab$, since $ - 2r + \eta \leq -2r + \etab 
 =  4 r -9 < -3$ by Remark \ref{lem:para_bd}, the above integral is integrable.
    Given $\eta \le \etab$, \eqref{eq:error_L2} follows directly.
\end{proof}

\subsection{Derivation of the linearized Euler equations}
\label{sec: derivation-euler}

We need the following basic results for the orthogonality of certain polynomials in Gaussian weighted $L ^2 (V)$ space.

\begin{lemma}[Orthogonality]\label{lem: orthogonality}
Recall the basis $\Phi_i$ from \eqref{eq:func_Phi} and $\cMM = \cs ^{-3} \mu (\vc)$ 
from \eqref{eq:localmax2}, where $\mu (x) = \left(\frac\kappa{2 \pi}\right) ^{3/2} \exp \left(
- \f{\kp |x| ^2}2 \right)$ is a Gaussian with variance $\f1\kp = \f35$ given in \eqref{eq:gauss}.  Define 
\begin{align*}
    \mathbf A (\vc) &= \left(\vc \otimes \vc - \f13 |\vc| ^2 \Id\right) \cMM ^{1/2}, &
    \mathbf b (\vc) &= \left(|\vc| ^2 - 3\right) \vc \cMM ^{1/2}.
\end{align*}
Then $\mathbf A _{ij}, \mathbf b _j \perp \cMM^{1/2} p(\vc)$ for any $p(\vc) \in \operatorname{Span} \{1, \vc _i, |\vc| ^2\}$ in $L ^2 (V)$ for all $1 \le i, j \le 3$.
In particular, $\mathbf A _{ij}, \mathbf b _j \perp \Phi_k$ in $L ^2 (V)$ 
for all $1 \le i, j \le 3$ and $0 \le k \le 4$.
\end{lemma}
The proof follows standard computations of the normal distribution and is therefore omitted. We refer to \cite[Eq. (3.64)]{golse_boltzmann_2005}, where a similar result is stated for the standard Gaussian with variance $1$.
\footnote{
In our case, the variance is $\f{1}{\kp} =  \f{3}{5}$, which leads to the term
$|\vc|^2 - 3$ in $  \mathbf b (\vc)$ instead of $|\vc|^2 - 5$ in the unit-variance case \cite[Eq. (3.64)]{golse_boltzmann_2005}. Since the matrix $ \vc \otimes \vc - \f13 |\vc| ^2 \Id$ in $  \mathbf A (\vc)$ is homogeneous 
in $\vc$,  the change of variance does not affect the orthogonality $\mathbf A _{ij} \perp \Phi_k$ in $L^2(V)$.
}


Recall the linearized equation \eqref{eq:lin}:
\begin{align*}
    (\pa _s + \cT) (\cMM ^{1/2} \td F) = \f1\es Q (\cM + \cMM ^{1/2} \td F, \cM + \cMM ^{1/2} \td F) - \eM.
\end{align*}
Let $p = p (\vc) \in \operatorname{Span} \{1, \vc _i, |\vc| ^2\}$ be a polynomial of $\vc$. Then it is orthogonal to $Q$. Taking inner product with $p$ on both sides, we obtain
\begin{align}\label{eqn: inner-prod-with-p}
    \int (\pa _s + \cT) (\cMM ^{1/2} \td F) \cdot p \, d V = -\la \eM, p \ra _\lvv.
\end{align}
We separate $\td F = \tFM + \tFm$. By product rule, we have
\begin{align*}
    \int (\pa _s + \cT) (\cMM ^{1/2} \tFM) \cdot p \,d V = \udb{\int (\pa _s + \cT) (p \cMM ^{1/2} \tFM) d V} _{I} \udb{-\int (\na p)(\vc) \cdot (\pa _s + \cT) \vc \cdot \cMM ^{1/2} \tFM d V} _{II}.
\end{align*}
Note that for any function $g$ we have 
\begin{align*}
    \bcv V \cdot \na _V g &= \bcv \div _V (g V) - 3 \bcv g, \\
    V \cdot \na _X g &= \bu \cdot \na _X g + (V - \bu) \cdot \na _X g = \bu \cdot \na _X g + (\div \bu) g + \div _X [(V - \bu) g].
\end{align*}
Recall that $V - \bu = \cs \vc$. Therefore 
\begin{align*}
    (\pa _s + \cT) g = [\pa _s + (\bcx X + \bu) \cdot \na _X + \div \bu - 3 \bcv] g + \div _V [\bcv V g] + \div _X (\cs \vc g).
\end{align*}
Apply this to $g \rightsquigarrow p \cMM ^{1/2} \tFM$ and integrate over $V$, we obtain
\begin{align*}
    I &= \int (\pa _s + \cT) (p \cMM ^{1/2} \tFM) d V \\
    &= [\pa _s + (\bcx X + \bu) \cdot \na _X + \div \bu - 3 \bcv] \int p \cMM ^{1/2} \tFM d V + \div _X \left(
        \cs \int p \vc \cMM ^{1/2} \tFM d V
    \right).
\end{align*}
For $II$, we use \eqref{eq:ds_t_vc}
\begin{align*}
    (\pa _s + \cT) \vc = -\left( \eu - 3 \na \cs + \vc \cdot \na \bu \right) - \left(\ec - \f13 \na \cdot \bu + \vc \cdot \na \cs\right) \vc.
\end{align*}
Therefore 
\begin{align*}
    II &= (\eu - 3 \na \cs) \cdot \int \na p \cdot \cMM ^{1/2} \tFM d V 
    + \left(\ec - \f13 \na \cdot \bu\right) \int \na p \cdot \vc \cMM ^{1/2} \tFM d V \\
    & \qquad + \int (\vc \cdot \na \bu) \cdot \na p \cdot \cMM ^{1/2} \tFM d V + \int (\vc \cdot \na \cs) \na p \cdot \vc \cMM ^{1/2} \tFM d V.
\end{align*}

Recall $\td F = \tFm + \tFM$. The terms $I, II$ account for the contribution from $\tFM$. 
\eqref{eqn: inner-prod-with-p} becomes
\begin{align}\label{eq:moments_I_II}
    I + II + \udb{\left\la (\pa _s + \cT) (\cMM ^{1/2} \tFm), p \right\ra _\lvv} _{III} = -\la \eM, p \ra _\lvv.
\end{align}

Since $\cMM^{1/2} \tFm$ is orthogonal to $1, \vc, |\vc|^2$, and since the scaling fields 
$X \cdot \na_X, V \cdot \na_V$ and $\pa_s$ preserve the orthogonality, which follows from \eqref{eq:proj_scaling}, using the notations $\cI_i$ from \eqref{eq:moments_Fm}, we obtain 
\beq\label{eq:moments_III}
\bal 
\left\la (\pa _s + \cT) (\cMM ^{1/2} \tFm), \, \left(1, \vc, \f{1}{3} |\vc|^2 \right) \right\ra _\lvv 
& = \left\la V \cdot \na_X (\cMM ^{1/2} \tFm), \, \left(1, \vc, \f{1}{3} |\vc|^2 \right) \right\ra _\lvv \\ 
& = \left(0, \cI_1(\tFm), \cI_2(\tFm) \right).
\eal 
\eeq

\paragraph{Equation of $\trho$}

Set $p (\vc) = 1$, then $\na p = 0$, $II = 0$. 
Recall $\la \eM, 1 \ra _\lvv = \cs ^3 \erho$ from \eqref{eq:error}.
Using \eqref{eq:moments_I_II} and \eqref{eq:moments_III} (the first component), we obtain the equation of $\trho$ 
in \eqref{eq:lin_euler}
\begin{align*}
    [\pa _s + (\bcx X + \bu) \cdot \na _X + \div \bu - 3 \bcv] \trho + \div _X (\cs \tilde \UU) = -\cs ^3 \erho.
\end{align*}

\paragraph{Equation of $\tu$}

Now let $p (\vc) = \vc$, then $\na p = \Id$. 
We first compute $I$:
\begin{align*}
    I = [\pa _s + (\bcx X + \bu) \cdot \na _X + \div \bu - 3 \bcv] \tilde \UU + \div _X \left(
        \cs \int \vc \otimes \vc \cMM ^{1/2} \tFM d V
    \right).
\end{align*}
Note that by orthogonality (see Lemma \ref{lem: orthogonality})
\begin{align*}
    \int \left(\vc \otimes \vc - \f13 |\vc| ^2 \Id \right) \cMM ^{1/2} \tFM d V = 0 .
\end{align*}
Therefore we have 
\begin{align*}
    \int \vc \otimes \vc \cMM ^{1/2} \tFM d V = \left(\int \f{|\vc| ^2}3 \cMM ^{1/2} \tFM d V\right) \Id = \tp \cdot \Id,
\end{align*}
and 
\begin{align*}
    I = [\pa _s + (\bcx X + \bu) \cdot \na _X + \div \bu - 3 \bcv] \tilde \UU + \na _X (\cs \tp).
\end{align*}
Next, we compute $II$. Because $\na p = \Id$, we have
\begin{align*}
    II = \int (\pa _s + \cT) \vc \cMM ^{1/2} \tFM d V &= \left(\eu - 3 \na \cs\right) \trho + \left(\ec - \f13 \na \cdot \bu\right) \tu \\
    & \qquad + \tu \cdot \na \bu + \na \cs : \int \vc \otimes \vc \cMM ^{1/2} \tFM d V \\
    &= \left(\eu - 3 \na \cs\right) \trho + \left(\ec - \f13 \na \cdot \bu\right) \tu + \tu \cdot \na \bu + \tp \na \cs.
\end{align*}

Recall that $\la \eM, \vc \ra = \cs ^3 \eu$ from \eqref{eq:error}.
Combining $I$, $II$, and the derivation of $III$ in 
\eqref{eq:moments_III} (the second component), we derive: 
\begin{align*}
    & [\pa _s + (\bcx X + \bu) \cdot \na _X + \div \bu - 3 \bcv] \tu + \na _X (\cs \tp) \\
    & \qquad + \left(\eu - 3 \na \cs\right) \trho + \left(\ec - \f13 \na \cdot \bu\right) \tu + \tu \cdot \na \bu + \tp \na \cs + \cI _1 (\tFm) = -\cs ^3 \eu.
\end{align*}
Recall $\tb= \trho -\tp$. Using $\trho = \tb + \tp$ and \eqref{eq:error_Ub}, we obtain 
\[
(\eu - 3 \na \cs)\trho = ( 3  \na \cs - 3 \cs^{-1} \bc \na \bc - 3 \na \cs ) \trho
=  - 3 \cs^{-1} \bc \na \bc (\tb + \tp). 
\]
Collecting similar terms, we derive the $\tu$-equation in \eqref{eq:lin_euler}.

\paragraph{Equation of $\tp$} Now set $p (\vc) = \f13 |\vc| ^2$, then $(\na p)(\vc) = \f23 \vc$. 
We compute $I$: 
\begin{align*}
    I = & [\pa _s + (\bcx X + \bu) \cdot \na _X + \div \bu - 3 \bcv] \tp + \div _X \left(
        \cs \int \f13 |\vc| ^2 \vc \cMM ^{1/2} \tFM d V
    \right).
\end{align*}
Note that $(|\vc| ^2 - 3) \vc \cMM ^{1/2} \perp \Phi _i$ by Lemma \ref{lem: orthogonality}, so 
\begin{align*}
    \int \f13 |\vc| ^2 \vc \cMM ^{1/2} \tFM d V = \int \vc \cMM ^{1/2} \tFM d V = \tu.
\end{align*}
Therefore 
\begin{align*}
    I & = [\pa _s + (\bcx X + \bu) \cdot \na _X + \div \bu - 3 \bcv] \tp + \div _X (\cs \tu).
\end{align*}
Next, we compute $II$. Using $\na p = \f23 \vc$ and 
$\int (|\vc|^2 - 3) \vc \cMM^{1/2} \tFM d V = 0$ by Lemma \ref{lem: orthogonality}, we have
\begin{align*}
    II & = \f23 \left(\eu - 3 \na \cs\right) \cdot \int \vc \cMM ^{1/2} \tFM d V + \f23 \left(\ec - \f13 \na \cdot \bu\right) \int |\vc| ^2 \cMM ^{1/2} \tFM d V \\
    & \qquad + \f23 \int (\vc \cdot \na \bu) \cdot \vc \cMM ^{1/2} \tFM d V + \f23 \int (\vc \cdot \na \cs) |\vc| ^2 \cMM ^{1/2} \tFM d V \\
    & = \f23 \left(\eu - 3 \na \cs\right) \cdot \tu + 2 \left(\ec - \f13 \na \cdot \bu\right) \tp + \f23 \tp \div \bu + 2 \na \cs \cdot \tu \\
    & = \f23 \eu \cdot \tu + 2 \ec \tp.
\end{align*}

 Recall that $\la \eM, \f13 |\vc| ^2 \ra = \cs ^3 \ep$ from \eqref{eq:error}. Combining $I$, $II$, and the derivation of $III$ in \eqref{eq:moments_III} (the third component), we derive the equation of $\tp$
in \eqref{eq:lin_euler}.
\begin{align*}
    & [\pa _s + (\bcx X + \bu) \cdot \na _X + \div \bu - 3 \bcv] \tp + \div _X (\cs \tu) + \f23 \eu \cdot \tu + 2 \ec \tp + \cI _2 (\tFm) = -\cs ^3 \ep.
\end{align*}

\paragraph{Equation of $\tb$} 

Recall $\tb= \trho -\tp$ and $\erho = \ep$ from Lemma \ref{lem: cutoff_error}. Taking the difference between the equation of $\trho$ and that of $\tp$,
we derive the equation of $\tb$ in \eqref{eq:lin_euler}:
\begin{align*}
    [\pa _s + (\bcx X + \bu) \cdot \na _X + \div \bu - 3 \bcv] \tb - \f23 \eu \cdot \tu - 2 \ec \tp - \cI _2 (\tFm) = -\cs ^3 (\erho - \ep) = 0.
\end{align*}

\section{Functional Inequalities}\label{app:ineq}

The goal of this appendix is to gather a few functional analytic bounds that are used throughout the paper. Lemmas \ref{lem:leib}-\ref{lem:norm_equiv} were established in \cite[Appendix C]{chen2024Euler}, and we refer there for details.

First, we record a Leibniz rule for radially symmetric vectors/scalars.

\begin{lemma}[Lemma A.4 \cite{buckmaster2022smooth}]\label{lem:leib}
	Let $f, g$ be radially symmetric scalar functions over $\R^d$ and let $\FF = F \ee _R = (F_1, .., F_d)$ and $\GG = G \ee _R = (G_1, .., G_d)$ be radially symmetric vector fields over $\R^d$. For integers $m\geq 1$ we have
	\[
		\bal
		| \D^m(\FF \cdot \na G_i)- \FF \cdot \na \D^m G_i - 2 m \pa_{\xi} F \; \D^m G_i |
		& \les_m \sum_{1 \leq j \leq 2 m} | \na^{2m +1-j} \FF| \cdot |\na^{j} G_i | , \\
		|\D^m (f \na g) - f \na \D^m g - 2 m \na f \D^m g |
		& \les_m \sum_{1 \leq j \leq 2 m} | \na^{2m +1-j} f| \cdot |\na^{j} g| , \\
		|\D^m(\FF \cdot \na g) - \FF \cdot \na \D^m g - 2 m \pa_{\xi} F \; \D^m g |
		& \les_m \sum_{1 \leq j \leq 2 m} | \na^{2m +1-j} \FF| \cdot |\na^{j} g | ,
		\\
		| \D^m(f \div(\GG) - f \div(\D^m \GG) - 2 m \na f \cdot \D^m \GG|
		& \les_m \sum_{1 \leq j \leq 2 m} | \na^{2m +1- j} f| \cdot |\na^{j} \GG|,
		\eal
	\]
	whenever $f, g, \{F_i\}_{i=1}^d, \{G_i\}_{i=1}^d$ are sufficiently smooth.
\end{lemma}

Next, we focus on Gagliardo-Nirenberg-type interpolation bounds with {\em weights}. In all of the following lemmas, we do not assume that the functions are radially symmetric.
\begin{lemma}[Lemma C.2 \cite{chen2024Euler}]\label{lem:interp_wg}
	Let $\d_1 \in (0, 1]$ and $\d_2 \in \R$. For integers $n\geq 0$ and sufficiently smooth functions $f$ on $\R^d$, we denote
	\[
		\b_n := 2 n \d_1 + \d_2,
		\qquad
		I_n := \int |\na^n f(y) |^2 \la y \ra^{\b_n} d y ,
	\]
	where as usual we let $\la y \ra = (1+ |y|^2)^{1/2}$.
	Then, for $n<m$ and for $\nu>0$, there exists a constant $C_{\nu, n,m} = C(\nu, n,m,\d_1,\d_2,d)>0$ such that
	\begin{equation}\label{eq:interp_convex}
		I_n \leq \nu I_m + C_{\nu,n,m} I_0.
	\end{equation}

\end{lemma}

\begin{lemma}\label{lem:norm_equiv}
	Let $\d_1 \in (0, 1], \d_2 \in \R$, and define $\b_n = 2 n \d_1 + \d_2$. Let $\psi_{n}$ be a weight satisfying the pointwise properties $\psi_{n}(y) \asymp_n \la y \ra^{\b_{n}}$ and $|\na  \psi_{n}(y)| \les_n \la y \ra^{\b_{n}-1}$. Then, for any $\nu > 0$ and $n\geq 0$, there exists a constant $C_{\nu,n} = C(\nu,n,\d_1,\d_2,d)>0$ such that\footnote{Throughout the paper we denote by $|\nabla^k f|$ the Euclidean norm of the $k$-tensor $\nabla^k f$, namely, $|\nabla^k f| = (\sum_{|\alpha|=k} |\partial^\alpha f|^2)^{1/2}$.}
	\bseq\label{eq:grad_est}
	\begin{align}	
		\int |\na^{2n} f|^2 \psi_{2n}
		& \leq (1 + \nu)\int | \D^n f |^2 \psi_{2n} +
		C_{\nu, n} \int |f|^2 \la y \ra^{\b_0}, \label{eq:grad_est:even} \\
	\int |\na^{2n+1} f|^2  \psi_{2n+1}
		& \leq (1 + \nu)\int | \na \D^n f |^2 \psi_{2n+1} +
		C_{\nu, n} \int |f|^2 \la y \ra^{\b_0}, \label{eq:grad_est:odd}
	\end{align}	
	\eseq
	for any function $f$ on $\R^d$ which is sufficiently smooth and has suitable decay at infinity.
\end{lemma}

\begin{proof}

The inequality \eqref{eq:grad_est:even} has been established in \cite[Lemma C.3]{chen2024Euler}. Below, we prove \eqref{eq:grad_est:odd}. We adopt the notation $I_n$ from Lemma \ref{lem:interp_wg}. 
Denote 
\beq\label{eq:interp_odd:nota}
	\th_n = \b_{n+1} = 2 n \d_1 + (2 \d_1 + \d_2)  , \quad g_{n} =  \psi_{ n+ 1 } ,
	\quad I_n = \int |\na^n f|^2 \psi_n .
\eeq

From the assumption of $(\b_n, \psi_n)$,  $(\th_n, g_n)$ satisfies the same assumptions 
as those of $(\b_n, \psi_n)$ in Lemma \ref{lem:norm_equiv}. Thus, 
for any $\nu > 0$ applying \eqref{eq:grad_est:even} with $(f, \psi_n, \b_n)
\rightsquigarrow  ( \pa_i f, g_n ,\th_n)$, we obtain 
\beq\label{eq:interp_odd1}
	\int |\na^{2n} \pa_i f |^2 g_{2n} 
\leq (1 + \nu ) \int |\D^n \pa_i f|^2 g_{2n} 
+ C_{\nu, n} \int |\pa_i f|^2 \la y \ra^{ \th_0 } 
:= J_{1,i} + J_{2,i}.
\eeq

For the second term, since $\th_0 = \b_1$ and $ \psi_n(y) \asymp_n \ang y^{\b_n} $, for any $ \nu_1 > 0$, applying Lemma \ref{lem:interp_wg}, we obtain 
\[
	J_{2,i } \leq \nu_1 \int |\na^{2n+1} f|^2 \ang y^{\b_{2n+1}}
	+ C(\nu_1, \nu, n) \int f^2 \ang y^{\b_0} \leq C_n \nu_1  I_{2n+1} + C(\nu, \nu_1, n) I_{0},
	\]
where $I_j$ is defined in \eqref{eq:interp_odd:nota}, and $C_n$ is some constant depending on $n$. 
Recall $g_{2n} = \psi_{2n+1}$. Combining the above two estimates and summing these estimates over $i$, we prove 
\[
\bal
I_{2n+1} &= 	\int |\na^{2n+1} f|^2 g_{2n}
\leq (1 + \nu ) \int |\D^n \na f|^2 g_{2n}  + C_n \nu_1 I_{2n+1} + C(\nu, \nu_1, n) I_0.
\eal
\]
Since $\nu, \nu_1 >0$ are arbitrary parameters, taking $\nu_1$ small enough so that $C_n \nu_1 < 1$, and then rewriting the above inequality, we prove \eqref{eq:grad_est:odd}. 
\end{proof}

We record the following estimates for the functional spaces $\cX _\eta ^m$ defined in~\eqref{norm:Xk} 
and $\cYe ^m$ defined in \eqref{norm:Y}. It is convenient to state estimates for a general dimension $d$, not just for $d=3$. We recall from Lemma~\ref{lem:wg} that the weights $\vp_m$ satisfy $\vp_m(y) \asymp_m \la y \ra^{m}$, and $|\na \vp_m(y)|\les_m \vp_{m-1}(y)$.

\begin{lemma}\label{lem:prod}
Suppose that $\eta \in [-100, 100]$.  

\begin{enumerate}[\upshape(1)]
    \item For any $f \in \cX^k_{\eta}$, $0\leq i\leq k-d$, and $X \in \R^d$, we have the pointwise estimate 
    \beq\label{eq:X_linf}
      \la X \ra^i | \na_X^i f(X)| \les_k \la X \ra^{- \f{\eta+d}{2} } \| f \|_{\cX_\eta^k}.
    \eeq

    \item Recall $D^{\al, \b} = \vp_1^{|\al|} \cs^{|\b|} \pa_X^{\al} \pa_V^{\b}$ from \eqref{eq:deri_wg}. 
    Suppose that the weight $\psi(X, V) > 0$ satisfies 
    \beq\label{eq:L2V_linf_assum}
      | \pa_X^{\al} \psi | \les_{\al}  \la X \ra^{-|\al|} \psi ,
      \quad   | \pa_V \psi| \les \cs^{-1} \psi .
    \eeq
    for any multi-indices $\al$ with $|\al| \leq d$. 
    
    For $f : \R^d \times \R^d \to R$  and multi-indices $\al, \b$ with 
    $|\al| + |\b| \leq k - d $, we have
    \bseq\label{eq:L2V_linf}
    \beq\label{eq:L2V_linfW}
          \| \psi(X, \cdot)^{1/2} D^{\al, \b} f(X, \cdot) \|_{L ^2 (V)}
    \les_{\al, \b} \la X \ra^{- \f{ \eta +d}{2} } 
    \sum_{ |p| + |q| \leq k }
    \| \psi(X, V)^{1/2} \la X \ra^{\eta/2} D^{p, q} f \|_{L^2(X, V)}
    \eeq
    pointwise for $X \in \R^d$. In particular, for $\al, \b$ with $|\al| + |\b| \leq k - d$, we have
    \begin{align}
       \|  D^{\al, \b}  f(X, \cdot) \|_{L ^2 (V)}
       &\les_k \la X \ra^{- \f{\eta +d}{2} }  \| f \|_{\cYe^k},  \label{eq:linf_Y} \\
          \| \Lams(X, \cdot)^{1/2}  D^{\al, \b} f(X, \cdot) \|_{L ^2 (V)}
        &  \les_k \la X \ra^{- \f{\eta + d}{2} } \| f \|_{\cYE^k}. 
        \label{eq:linf_Ylam}
    \end{align}
    \eseq

    Moreover, we have the pointwise estimate 
    \beq\label{eq:linf_Y_point}
    	|f(X, V)| \les \cs^{- \f{d}{2}} \| D_V^{\leq d} f(X, \cdot) \|_{L ^2 (V)}
    	\les \cs^{- \f{d}{2}} \ang X^{- \f{\eta +d}{2}} \| f\|_{\cYe^{2 d}}.
    \eeq
\end{enumerate}
\end{lemma}

Result (1) is essentially the same as \cite[Lemma C.4]{chen2024Euler}. 
For completeness, we present the proof. %

\begin{proof}

We first obtain a pointwise estimate of weighted derivatives of $f$. Consider the cone with vertex at $X$ extending towards infinity: $\Om(X) := \{ z \in \R^d \colon z_j \sgn(X_j) \geq |X_j|, \forall 1\leq j\leq d \}$. 
For any fixed $V$, by integrating on rays extending to infinity, we have
\begin{align}\label{eq:linf_IXV}
I(X, V) & := 
  \psi(X, V) \la X \ra^{ \eta+d  } | D^{\al , \b} f |^2 =  \psi(X, V) \la X \ra^{ \eta +d + 2 |\al| }  \cs^{2 \b} (\pa_X^{\al} \pa_V^{\b} f(X, V))^2  \notag \\
& \les \int_{ Y  \in \Om(X) }
\B| \pa_{X_1} \pa_{X_2} .. \pa_{X_d}  
\B(  \psi(X, V) \la X \ra^{\eta +d + 2 |\al| }  \cs^{2 \b} \pa_X^{\al} \pa_V^{\b} f(X, V) \B)^2 \B|  d X \\
 & \les \sum_{ |\th_1| + |\th_2|  + |\th_3|= d }  \int_{ Y  \in \Om(X) }
 \B| \pa_X^{ \al + \th_1 } \pa_V^{\b} f 
 \cdot  \pa_X^{ \al + \th_2} \pa_V^{\b} f 
\cdot \pa_X^{\th_3 } \B( \psi(X, V) \la X \ra^{ \eta +d + 2 |\al| }  \cs^{2 \b} \B) \B|
  d X . \notag 
 \end{align}

Using the estimates of $\cs$ in \eqref{eq:dec_S} and assumption \eqref{eq:L2V_linf_assum} on $\psi$, we obtain 
\[
\bal
  & | \pa_X^{\th_3 } ( \psi(X, V) \la X \ra^{\eta +d + 2 |\al| }  \cs^{2 \b} )| 
  \les_{ \al}  \psi(X, V) \la X \ra^{\eta +d + 2 |\al| - |\th_3| }  \cs^{2 \b} \\
   =& \psi(X, V) \la X \ra^{ \eta + 2 |\al| + |\th_1| + |\th_2|} \cs^{2 \b} 
  \les_{ \al}  \psi(X, V)  \ang X ^\eta \vp_1^{ 2 |\al| + |\th_1| + |\th_2| } \cs^{2 \b} .
  \eal 
\]

Recall $D^{\al, \b} = \vp_1^{|\al|} \cs^{|\b|} \pa_X^{\al} \pa_V^{\b}$ from \eqref{eq:deri_wg}.
Combining the above estimates and applying $\vp_1 \asymp \la X \ra$ from Lemma \ref{lem:wg}, 
we establish 
\beq\label{eq:prod_pf1}
  I(X, V) \les  \sum_{ |\th_1| + |\th_2| \leq  d } 
  \int_{ Y \in \Om(X) } | D^{\al + \th_1 ,\b} f 
  \cdot D^{\al + \th_2, \b} f | \psi(X, V)  \ang X ^\eta d X  .
\eeq

\paragraph{Proof of \eqref{eq:X_linf}}
For function $f$ independent of $V$, applying the above estimate with $\b = 0, \psi \equiv 1$, any $\al$ with $|\al| = k$, and Cauchy--Schwarz inequality, we establish 
\[
   \la X \ra^{\eta +d + 2 |\al| }   ( \pa_X^{\al}  f(X) )^2
   = I(X) \les \sum_{ p \leq k } \| \la X \ra^{\eta/2} D_X^{p} f   \|_{L^2}^2 .
\]

Recall the $\cX _\eta ^n$ norm from \eqref{norm:Xk}. Applying 
the interpolation in Lemma \ref{lem:interp_wg} and Lemma \ref{lem:norm_equiv} with $\psi_n = \vp_1^{2 n}$ and $\b_n = 2 n + \eta $, we further obtain 
\[
	   \la X \ra^{\eta+d + 2 |\al| }   ( \pa_X^{\al}  f(X) )^2   \les_k \| f \|_{\cX _\eta ^k}^2.
\]

Multiplying $\la X \ra^{- (\eta+d)}$ on both sides of the above estimate,  we prove \eqref{eq:X_linf}. 

\paragraph{Proof of \eqref{eq:L2V_linfW}}
Integrating \eqref{eq:prod_pf1} over $V \in \R^d$, 
using Cauchy--Schwarz inequality, and $|\al| + |\b| + |\th| \leq k$, we establish 
\[
\bal
  \int I(X, V) d V
  & \les 
  \sum_{ |\th_1| + |\th_2|  \leq d } 
  \int_{ Y \in \Om(X) } | D^{\al + \th_1 ,\b} f 
  \cdot D^{\al + \th_2, \b} f |  \psi(X, V) \ang X ^\eta d X d V  \\
& \les \sum_{ |p| + |q| \leq k }   
 \| \psi(X, V)^{1/2} \la X \ra^{\eta/2} D^{p, q} f \|_{L^2}^2.
\eal 
\]

Multiplying $\la X \ra^{- (\eta+d)}$ on both sides of the above estimate,  we prove \eqref{eq:L2V_linfW}.

\paragraph{Proof of \eqref{eq:linf_Y} and \eqref{eq:linf_Ylam}}

 Recall the norm $\cYe^k$ from \eqref{norm:Y}.  Since $\psi(X, V) \equiv 1$ satisfies assumptions in  \eqref{eq:L2V_linf_assum}, using \eqref{eq:L2V_linfW} with $\psi \equiv 1$, we prove \eqref{eq:linf_Y}.

 Recall the norm $\cYE^k$ from \eqref{norm:Y}. From estimate \eqref{eq:Lams_est} in Lemma \ref{lem:basis},  
$\psi(X, V) = \Lams(X, V) = \cs^{\g+3} \la \vc \ra^{\g + 2}$ satisfies the assumptions in \eqref{eq:L2V_linf_assum}.  Using \eqref{eq:L2V_linfW} with $\psi = \Lams$, we prove \eqref{eq:linf_Ylam}.

\paragraph{Proof of \eqref{eq:linf_Y_point}}

The second inequality in \eqref{eq:linf_Y_point} follows from \eqref{eq:linf_Y} with $\al = 0,  |\b| \leq d$. 
To prove the first inequality, we fix $X, V$ and introduce $\Om(V) := \{ z \in \R^d \colon z_j \sgn(V_j) \geq |V_j|, \forall 1\leq j\leq d \}$. Following the argument in \eqref{eq:linf_IXV} and using Cauchy--Schwarz inequality, we estimate 
\[
\bal
	\cs^{d} f^2(X, V)
	& \les  \int_{\Om(V)} \B| \pa_{V_1} \pa_{V_2} .. \pa_{V_d} ( 	\cs^{d} f^2(X, V) )  \B| d V
\les \sum_{|\th_1| + |\th_2| = d}  \int  | \cs^{|\th_1|} \pa_V^{\th_1} f
\cdot \cs^{|\th_2|} \pa_V^{\th_2} f | d V \\
& \les \| D_V^{ \leq d } f(X,\cdot) \|_{L ^2 (V)}^2,
\eal 
\]
and prove the first inequality in \eqref{eq:linf_Y_point}.

We conclude the proof of Lemma \ref{lem:prod}. 
\end{proof}

\section{Estimates of projections and related functions}

In this appendix, we estimate the projections 
defined in \eqref{def:proj} and their related functions.

\subsection{Estimate functions of \texorpdfstring{$\vc$}{V°}}\label{sec:non_func_V}

To facilitate our proof, we introduce the following class of functions with algebraic bound.

\begin{definition}
    \label{def: bF}
    We say a function $f \in C ^\infty ((s _0, \infty) \times \R ^3$) has good decay property if for any multi-index $\al$ it satisfies
    \begin{align*}
        |\pa _X ^\al f (s, X)| \les _\al \ang X ^{-|\al|}.
    \end{align*}
    Denote $\bF$ the class of functions with good decay property. It is straightforward to verify that $\bF$ forms an algebra. For $\eta \in \R$ we define $\bF ^\eta = \ang X ^\eta \bF$, then $f \in \bF ^\eta$ iff $\pa _X ^\al f (X) \les \ang X ^{\eta - |\al|}$ for all multi-index $\al$ using simple induction. Note that $\bF ^{\eta _1} \cdot \bF ^{\eta _2} \subset \bF ^{\eta _1 + \eta _2}$. Vector-valued function is said to be of class $\bF ^\eta$ if each component is in $\bF ^\eta$. By definition, we have 
    $\bF^0 = \bF$. 
\end{definition}

\begin{lemma}
    \label{lem:example_fl}
    The following examples are in class $\bF ^\eta$:
    \begin{enumerate}[\upshape (1)]
        \item $\bu \in \bF ^{-r + 1}$, $\na \cs \in \bF ^{-r}$.
        \item $\cs ^{-1} \in \bF ^{r - 1}$.
        \item $\na \log \cs \in \bF ^{-1}$.
        \item $\cs ^l \pa ^\al _X \cs ^{-l} \in \bF ^{-|\al|}$ for any $l \in \R$ and any multi-index $\al$.
        \item $\cs ^{-l} \in \bF ^{l (r - 1)}$ for any $l > 0$.
        \item $\varphi _1 \in \bF ^1$.
    \end{enumerate}
\end{lemma}

\begin{proof}
    (1) is a direct consequence of \eqref{eq:dec_U} and \eqref{eq:dec_S}. We prove the rest. 
    
    \begin{enumerate}[\upshape (1)]
        \setcounter{enumi}{1}
        \item To see $\cs ^{-1} \in \bF ^{r - 1}$, we first prove
    \begin{align}
        \label{eq:decay_cs_-1}
        |\cs \pa ^\al _X \cs ^{-1}| \les _\al \ang X ^{- |\al|}
    \end{align}
    inductively. \eqref{eq:decay_cs_-1} clearly holds for $\alpha = 0$. Moreover, if \eqref{eq:decay_cs_-1} holds for any $\al' \prec \al$ then
    \begin{align*}
        0 = \pa ^\al _X (\cs ^{-1} \cs) &= \cs \pa ^\al _X \cs ^{-1} + \sum _{\al' \prec \al} C ^\al _{\al'} \cdot \pa ^{\al - \al'} _X \cs \cdot \pa ^{\al'} _X \cs ^{-1} \\ 
        &= \cs \pa ^\al _X \cs ^{-1} + \sum _{\al' \prec \al} C ^\al _{\al'} \cdot \frac{\pa ^{\al - \al'} _X \cs}{\cs} \cdot \cs \pa ^{\al'} _X \cs ^{-1}.
    \end{align*}
    Note that $|\pa ^{\al - \al'} _X \cs| \les \ang X ^{-r + 1 - |\al - \al'|} \les \ang X ^{-|\al| + |\al'|} \cs$, so \eqref{eq:decay_cs_-1} holds by the inductive assumption. From \eqref{eq:decay_cs_-1} and the fact $\cs \gtrsim \ang X ^{-r + 1}$, we know
    \begin{align*}
        |\pa ^\al _X \cs ^{-1}| \les \cs ^{-1} \ang X ^{-|\al|} \les \ang X ^{r - 1 - |\al|},
    \end{align*}
    so $\cs ^{-1} \in \bF ^{r - 1}$. 
    
    \item $\nabla \log \cs = \cs ^{-1} \na \cs$ which is in $\bF ^{-1}$ by the previous two conclusions.

    \item We first prove 
    \begin{align}
        \label{eq:decay_cs_-l}
        |\cs ^l \pa ^\al _X \cs ^{-l}| \les _{\al, l} \ang X ^{- |\al|} \qquad \forall l \in \R
    \end{align}
    inductively. Again, it is true for $\al = 0$. Suppose $\al = \al' + \ee _i$, then 
    \begin{align*}
        \cs ^l \pa ^\al _X \cs ^{-l} &= \cs ^l \pa ^{\al'} _X \left[
            \cs ^{-l} \cs ^l \pa _{X _i} \cs ^{-l}
        \right] \\
        &= \cs ^l \pa ^{\al'} _X \left[
            \cs ^{-l} \cdot (-l) \pa _{X _i} \log \cs
        \right] \\
        &= -l \sum _{\al'' \preceq \al'} C ^{\al'} _{\al''} \cdot \cs ^l \pa ^{\al''} _X \cs ^{-l} \cdot \pa ^{\al' - \al''} _X \pa _{X _i} \log \cs.
    \end{align*}
    If \eqref{eq:decay_cs_-l} holds for any $\al'' \preceq \al' \prec \al$, then together with $\na \log \cs \in \bF ^{-1}$ we conclude 
    \begin{align*}
        |\cs ^l \pa ^\al _X \cs ^{-l}| \les \ang X ^{-|\al''|} \ang X ^{-|\al' - \al''| - 1} = \ang X ^{-|\al|},
    \end{align*}
    so \eqref{eq:decay_cs_-l} is proved. 

    To prove (4), we need to show 
    \begin{align}
        \label{eq:decay_cs_-l_b}
        \pa ^\b _X (\cs ^l \pa ^\al _X \cs ^{-l}) \les _{\al, l, \b} \ang X ^{-|\al| - |\b|}.
    \end{align}
    \eqref{eq:decay_cs_-l} proves the $\b = 0$ case, and we now show \eqref{eq:decay_cs_-l_b} for $\b > 0$ inductively. Suppose $\b = \b' + \ee _i$, then 
    \begin{align*}
        \pa ^\b _X (\cs ^l \pa ^\al _X \cs ^{-l}) &= \pa ^{\b'} \left[
            \pa _{X _i} \cs ^l \pa ^\al _X \cs ^{-l} + \cs ^l \pa ^\al _X \pa _{X _i} \cs ^{-l}
        \right] \\
        &= \pa ^{\b'} \left[
            l \cdot \pa _{X _i} \log \cs \cdot \cs ^l \pa ^\al _X \cs ^{-l}
        \right] + \pa ^{\b'} _X (\cs ^l \pa ^\al _X \pa _{X _i} \cs ^{-l}) \\
        &= l \sum _{\b'' \preceq \b'} C ^{\b'} _{\b''} \cdot \pa ^{\b''} _X \pa _{X _i} \log \cs \cdot \pa ^{\b' - \b''} _X \cs ^l \pa ^\al _X \cs ^{-l} + \pa ^{\b'} (\cs ^l \pa ^\al _X \pa _{X _i} \cs ^{-l}).
    \end{align*}
    Therefore, if \eqref{eq:decay_cs_-l_b} holds for $\b'$ and $\b''$, then it holds for $\b$ as well.
    
    \item When $l > 0$, by \eqref{eq:decay_cs_-l} we have 
    \begin{align*}
        |\pa ^\al _X \cs ^{-l}| \les \cs ^{-l} \ang X ^{-|\al|} \les \ang X ^{l (r - 1) - |\al|},
    \end{align*}
    so $\cs ^{-l} \in \bF ^{l (r - 1)}$.

    \item From Lemma \ref{lem:wg} we see $\varphi _1 \in C ^\infty (\R ^3)$ is smooth and $\varphi _1 (X) = C (1 + c _3 \ang X)$ for $|X| \ge R _2 + 1$, where $C$, $c _3$, $R _2$ are constants. Verification of the good decay property is straightforward.
    \end{enumerate} 
\end{proof}

This lemma enables us to derive the following corollary, which shows we can concatenate derivatives up to lower order corrections.
\begin{cor}
    \label{cor: concatenate-derivative}
    Let $f$ be a function of $X, V$. For multi-index $\al, \b, \al _1, \b _1$, we have
    \begin{align*}
        D ^{\al _1, \b _1} D ^{\al, \b} f -D ^{\al _1 + \al, \b _1 + \b} f = \sum _{\substack{\al' \prec \al _1 + \al \\ \b' \preceq \b _1 + \b}} c _{\al', \b', \al _1, \b _1, \al, \b} D ^{\al', \b'} f
    \end{align*}
    where $c _{\al', \b', \al _1, \b _1, \al, \b} \in \bF$. With a slight abuse of notation, we can write 
    \begin{align*}
        D ^{\al _1, \b _1} D ^{\al, \b} f -D ^{\al _1 + \al, \b _1 + \b} f = c D ^{\prec (\al _1 + \al, \b _1 + \b)} f,
    \end{align*}
    where $c$ is a class $\bF$ tensor.
\end{cor}

\begin{proof}
    Since the weight does not depend on $V$, we have $D ^{\al, \b + \b _1} = D ^{0, \b _1} D ^{\al, \b}$ and
    \begin{align*}
        D ^{\al _1, \b _1} D ^{\al, \b} f &= D ^{0, \b _1} D ^{\al _1, 0} D ^{\al, \b} f, &
        D ^{\al _1 + \al, \b _1 + \b} f &= D ^{0, \b _1} D ^{\al _1 + \al, \b} f.
    \end{align*}
    Thus we can assume $\b _1 = 0$ without loss of generality. 
    By induction in $|\al _1|$, we only need to prove the case $|\al _1| = 1$. Note that 
    \begin{align*}
        D ^{\ee _i, 0} D ^{\al, \b} f = D ^{\ee _i, 0} \varphi _1 ^{|\al|} \cs ^{|\b|} \pa _X ^\al \pa _V ^\b f = D ^{\al + \ee _i, \b} f + D ^{\ee _i, 0} \log (\varphi _1 ^{|\al|} \cs ^{|\b|}) \cdot D ^{\al, \b} f.
    \end{align*}
    We conclude the proof since $D ^{\ee _i, 0} \log (\varphi _1 ^{|\al|} \cs ^{|\b|}) \in \bF$ (see Lemma \ref{lem:example_fl} (3) and (6)). 
\end{proof}

\begin{definition}
    \label{def: polynomial}
    We say $p \in C ^\infty ((s _0, \infty) \times \R ^3 \times \R ^3)$ is a class $\bF$ polynomial of $\vc$  if 
    \begin{align*}
        p (s, X, V) = \sum _{|\alpha| \le N} c _\alpha (t, X) \vc ^\alpha 
    \end{align*}
    with $c _\alpha \in \bF$ and $N \ge 0$. The smallest $N$ is called the degree of $p$. We define class $\bF ^\eta$ polynomial of $\vc$ similarly when the coefficients are of class $\bF ^\eta$.
\end{definition}

\begin{lemma}[Estimate of $\vc$]\label{lem:est_Dxvring}
For any multi-index $\al$, we have $D ^{\al, 0} \vc _i$ is a class $\bF$ polynomial of $\vc$ with degree at most $1$. In particular,
\bseq
\begin{align}
    \label{eq:dx_vring}
    \pa _X ^\al \vc & = c _{\al, 1} \vc + \boldsymbol c _{\al, 2}, \qquad c _{\al, 1}, \boldsymbol c _{\al, 2} \in \bF ^{-|\al|} .
\end{align}
\eseq
\end{lemma}

\begin{proof}
    We prove \eqref{eq:dx_vring} by induction. The case $\al = 0$ is trivial. For $|\al| = 1$, we use 
    \begin{align*}
        \pa _{X _i} \vc = - \frac{V - \bu}{\cs ^2} \pa _{X _i} \cs - \frac{\pa _{X _i} \bu}{\cs}= - \pa _{X _i} \log \cs \vc - \cs ^{-1} \pa _{X _i} \bu =: c _{\ee _i, 1} \vc + \boldsymbol c _{\ee _i, 2}.
    \end{align*}
    By Lemma \ref{lem:example_fl}, $c _{\ee _i, 1}, \boldsymbol c _{\ee _i, 2} \in \bF ^{-1}$. For $|\al| \ge 2$, we can write $\al = \al' + \ee _i$ and 
    \begin{align*}
        \pa _X ^\al \vc 
        = \pa _{X _i} \pa _X ^{\al'} \vc 
        &= \pa _{X _i} (c _{\al', 1} \vc + \boldsymbol c _{\al', 2}) 
        \\
        &= \pa _{X _i} c _{\al', 1} \vc - c _{\al', 1} (c _{\ee _i, 1} \vc + \boldsymbol c _{\ee _i, 2}) + \pa _{X _i} \boldsymbol c _{\al', 2}.
    \end{align*}
    so $c _{\al, 1} = \pa _{X _i} c _{\al', 1} - c _{\al', 1} c _{\ee _i, 1} \in \bF ^{-|\al'| - 1} = \bF ^{-|\al|}$, $\boldsymbol c _{\al, 2} = \pa _{X _i} \boldsymbol c _{\al', 2} - c _{\al', 1} \boldsymbol c _{\ee _i, 2} \in \bF ^{-|\al'| - 1} = \bF ^{-|\al|}$. Since $\varphi _1 \in \bF ^1$, we know $\varphi _1 ^{|\al|} \pa _X ^\al \vc = \varphi _1 ^{|\al|} c _{\al, 1} \vc + \varphi _1 ^{|\al|} \boldsymbol c _{\al, 2}$, where $\varphi _1 ^{|\al|} c _{\al, 1}, \varphi _1 ^{|\al|} \boldsymbol c _{\al, 2}$ are both of class $\bF$. Therefore, $D ^{\al, 0} \vc _i$ is a class $\bF$ polynomial of degree at most 1. 
\end{proof}

\begin{remark}
    \label{rem: Dab-Vc}
    Note that $D ^{0, \ee _i} \vc = \ee _i$, and $D ^{\al, \b} \vc = 0$ when $|\b| \ge 2$.
\end{remark}

\begin{cor}
    \label{cor: Dab-polynomial}
    If $p (s, X, V)$ is a class $\bF$ polynomial of $\vc$ of degree $N$, then $D ^{\al, \b} p (s, X, V)$ is also a class $\bF$ polynomial of $\vc$, with degree at most $N - |\b|$. Recall that degree of $0$ is $-\infty$. 
\end{cor}

\begin{proof}
    Without loss of generality, assume $p (s, X, V) = c (s, X) \vc ^{\b'}$ for some coefficient $c \in \bF$ and multi-index $\b'$. Then 
    \begin{align}
        \label{eqn: split-al-b}
        D ^{\al, \b} p = \cs ^{|\b|} D ^{\al, 0} (\cs ^{-|\b|} D ^{0, \b} p) = \sum _{\al' \preceq \al} C ^\al _{\al'} \cdot \cs ^{|\b|} D ^{\al - \al', 0} \cs ^{-|\b|} \cdot D ^{\al', 0} D ^{0, \b} p.
    \end{align}
    By Lemma \ref{lem:example_fl} (4) we know $\cs ^{|\b|} D ^{\al - \al', 0} \cs ^{-|\b|} \in \bF$, so it remains to verify $D ^{\al', 0} D ^{0, \b} p$ is a class $\bF$ polynomial of $\vc$.
    
    Note that $D ^{0, \b} (\vc ^{\b'}) = 0$ if $\b \not\preceq \b'$, so $D ^{\al', \b} p = 0$ whenever $|\b| > |\b'|$. When $\b \preceq \b'$ we have
    \begin{align*}
        D ^{0, \b} (\vc ^{\b'}) = C ^{\b'} _\b \vc ^{\b' - \b}.
    \end{align*}
    Therefore, $D ^{0, \b} p = c (s, X) \cdot C ^{\b'} _\b \vc ^{\b - \b'}$, and by Lemma \ref{lem:est_Dxvring} with product rule we know $D ^{\al, 0} D ^{0, \b} p$ is a class $\bF$ polynomial with degree $|\b - \b'|$. 
\end{proof}

Next, we estimate functions involving $\mu(\vc)$.

\begin{lemma}\label{lem:H_deri}
Let $H (s, X, V) = \mu(\vc)^{1/2} p (s, X, V)$, where $\mu$ is the Maxwellian defined in \eqref{eq:gauss}, and $p (s, X, V)$ is a class $\bF$ polynomial of $\vc$ with degree $d _p$. Then 
\begin{align*}
    D ^{\al, \b} H (s, X, V) = \mu (\vc) ^{1/2} \td p (s, X, \vc)
\end{align*}
where $\td p (s, X, \vc)$ is another class $\bF$ polynomial of $\vc$, with degree $d _p + |\b| + 2 |\al|$. 
\end{lemma}

\begin{proof}
    By the product rule and Corollary \ref{cor: Dab-polynomial}, it suffices to verify the case $p \equiv 1$. That is, 
    \begin{align}
        \label{eqn: Dab-mu1/2} 
        D ^{\al, \b} \mu (\vc) ^{1/2} = \mu (\vc) ^{1/2} p _{\al, \b} (s, X, \vc)
    \end{align}
    where $p _{\al, \b} (s, X, \vc)$ is some class $\bF$ polynomial of $\vc$ with degree $|\b| + 2 |\al|$.

    We use induction. Assume \eqref{eqn: Dab-mu1/2} holds for all multi-index $(\al, \b)$ with $|\al| + |\b| \le k$. Now we want to show it also holds for $(\al + \al _1, \b + \b _1)$ where $|\al _1| + |\b _1| = 1$. By Corollary \ref{cor: concatenate-derivative}, we have 
    \begin{align*}
        D ^{\al + \al _1, \b + \b _1} \mu (\vc) ^{1/2} = D ^{\al, \b} D ^{\al _1, \b _1} \mu (\vc) ^{1/2} + c D ^{\prec (\al + \al _1, \b + \b _1)} \mu (\vc) ^{1/2}.
    \end{align*}
    Using inductive assumption, the lower order term $c D ^{\prec (\al + \al _1, \b + \b _1)} \mu (\vc) ^{1/2} = p \mu (\vc) ^{1/2}$ with some class $\bF$ polynomial $p$ with degree at most $|\al + \al _1| + 2 |\b + \b _1| - 1$. For the leading term, note we have 
    \begin{align*}
        D ^{\al _1, \b _1} \mu (\vc) ^{1/2} = \f12 \mu (\vc) ^{1/2} D ^{\al _1, \b _1} \log \mu (\vc) = -\kp _2 \mu (\vc) ^{1/2} D ^{\al _1, \b _1} |\vc| ^2.
    \end{align*}
    By Corollary \ref{cor: Dab-polynomial}, $D ^{\al _1, \b _1} |\vc| ^2$ is a class $\bF$ polynomial of $\vc$ with degree 2 if $\al _1 = 1$, with degree 1 if $\b _1 = 1$. Therefore, $D ^{\al, \b} D ^{\al _1, \b _1} \mu (\vc) ^{1/2} = p \mu (\vc) ^{1/2}$ where $p$ has degree at most $|\b| + 2 |\al| + 2$ if $\al _1 = 1$, $|\b| + 2 |\al| + 1$ if $\b _1 = 1$. In either case, $p$ has degree at most $2 |\al + \al _1| + |\b + \b _1|$. The induction is completed.
\end{proof}

We estimate the transport operator applied to $\vc$. Recall the transport operator defined in \eqref{def:lin_op}: $\partial _s + \cT = \partial _s + \bcx X \cdot \nabla _X + \bcv V \cdot \nabla _V + V \cdot \nabla _X$.

\begin{lemma}
    \label{lem:ds_t}
    $(\pa _s + \cT) \vc$ is a class $\bF ^{-r}$ polynomial of $\vc$ with degree 2, which equals to 
    \begin{align}
        (\pa _s + \cT) \vc = -\eu + 3 \na \cs - \vc \cdot \na \bu - \left(\ec - \f13 \na \cdot \bu + \vc \cdot \na \cs\right) \vc = \cO (\ang X ^{-r} \ang {\vc} ^2) .
        \label{eq:ds_t_vc}
    \end{align}
    The error term $\eM$ defined in 
    \eqref{eq:error} equals to 
    \begin{align*}
        \eM
        &= -\kappa \cM \vc \cdot \left(-\eu + 3 \nabla \cs -\vc\cdot \nabla \bu - \left(\ec - \f13 \na \cdot \bu + \vc \cdot \na \cs\right) \vc \right) = \cM p _3 (s, X, V),
    \end{align*}
    where $p _3 (s, X, V)$ is a class $\bF ^{-r}$ polynomial of $\vc$ with degree 3. 
\end{lemma}

\begin{proof}
	Recall $\ec$ and $\eu$ were computed in \eqref{eq:error_a} and \eqref{eq:error_b}. We first apply $\partial _s + \cT$ to $\cs$:
	\begin{align*}
		(\partial _s + \cT) \cs & = [\partial _s + (\bcx X + V) \cdot \nabla] \cs \\
		& = [\partial _s + (\bcx X + \bu) \cdot \nabla] \cs + (V - \bu) \cdot \nabla \cs \\
		& = \left[\partial _s + (\bcx X + \bu) \cdot \nabla - \bcv + \f13 \na \cdot \bu \right] \cs + \left(\bcv - \f13 \na \cdot \bu\right) \cs \\
		& \qquad + \cs \vc \cdot \nabla \cs \\
		& = \cs \ec + \left(\bcv - \f13 \na \cdot \bu\right) \cs + \cs \vc \cdot \nabla \cs.
	\end{align*}
    Dividing $\cs$, by bound \eqref{eq:dec_S} and \eqref{eq:error_ec} we know
	\begin{align}
        \label{eqn: ds_t_logcs}
		(\pa _s + \cT) \log \cs &= \ec + \bcv - \f13 \na \cdot \bu + \vc \cdot \na \cs = \bcv + \cO (\ang X ^{-r} \ang{\vc}).
	\end{align}
    In fact, $(\pa _s + \cT) \log \cs - \bcv$ is a class $\bF ^{-r}$ polynomial of $\vc$ with degree 1. 
	Next, we apply $\partial _s + \cT$ to $\bu$:
	\begin{align*}
		(\pa _s + \cT) \bu = \cT \bu & = [\pa _s + (\bcx X + V) \cdot \nabla] \bu \\
		& = [\pa _s + (\bcx X + \bu) \cdot \nabla] \bu + (V - \bu) \cdot \nabla \bu \\
		& = \cs \eu + \bcv \bu - 3 \cs \nabla \cs + \cs \vc \cdot \nabla \bu.
	\end{align*}
	Finally, we apply $\partial _s + \cT$ to $\vc$ to get 
	\begin{align*}
		(\partial _s + \cT) \vc & = (\partial _s + \cT) \left(\frac{V - \bu}{\cs}\right) \\
		& = \frac1{\cs} \left(
		\cT V - \cT \bu - \frac{V - \bu}{\cs} (\partial _s + \cT) \cs
		\right) \\
		& = \frac1{\cs} \bigg[
		\bcv V - \left(\cs \eu + \bcv \bu - 3 \cs \nabla \cs + \cs \vc \cdot \nabla \bu\right) \\
		& \qquad \qquad
			- \left(\ec + \bcv - \f13 \na \cdot \bu + \vc \cdot \na \cs\right) \cs \vc
		\bigg] \\
		& = -\left( \eu - 3 \na \cs + \vc \cdot \na \bu\right) - \left(\ec - \f13 \na \cdot \bu + \vc \cdot \na \cs \right) \vc  \\
		& = \cO (\ang X ^{-r} \ang{\vc} ^2),
	\end{align*}
    using the decay estimates \eqref{eq:dec_U}, \eqref{eq:dec_S}, \eqref{eq:error_ec} and  \eqref{eq:error_eu}. Now we compute $\eM = (\pa _s + \cT) \cM$:
	\begin{align*}
		(\partial _s + \cT) \cM = \cM (\partial _s + \cT) \log \cM
		& = \cM (\partial _s + \cT) \left(-\frac\kappa2 |\vc| ^2 \right) = -\kappa \cM \vc \cdot (\partial _s + \cT) \vc = \cM p _3 (s, X, V).
	\end{align*}
	The proof is completed.
\end{proof}

\subsection{Commutators between \texorpdfstring{$\cP$, $D ^{\al, \b}$, and $\cT$}{P, D α,β and T}}
\label{sec: commutator}

In this subsection, we justify the following commutator estimates.

\begin{lemma}[Commutator estimate]
\label{lem: commu-derivative-projection}
Let $f \in C ^\infty ((s _0, \infty) \times \R ^3 \times \R ^3)$.
\begin{enumerate}[\upshape (1)]
    \item \emph{Commuting $\pa _s + \cT$ and $\cP _m$:} 
    \begin{align*}
        \cP _m[(\pa _s + \cT) f] - (\pa _s + \cT) \cP _m f &= \cP _m [(V \cdot \na _X + \dcm + \td d _\cM) \cP _M f] \\
        & - \cP _M [(V \cdot \na _X - \dcm - \td d _\cM) \cP _m f].
    \end{align*}
    
    \item \emph{Commuting $\pa _s + \cT$ and $D ^{\al, \b}$:}
    \begin{align}
        \label{eqn: comm-free-transp-Dab}
        [V \cdot \na _X, D ^{\al, \b}] f &= \cO _{\al, \b} (\cs \ang X ^{-1} \ang \vc) \sum _{|\al'| + |\b'| = |\al| + |\b|} |D ^{\al', \b'} f|, \\
        \label{eqn: comm-transp-Dab}
        [\pa _s + \cT - V \cdot \na _X, D ^{\al, \b}] f 
        &= \cO _{\al, \b} (\cs \ang X ^{-1} \ang \vc + \ang X^{-1})  |D ^{\al, \b} f|.
    \end{align}

    \item \emph{Commuting $D ^{\al, \b}$ and $\cP _M$:}
    recall that $\Phi _i$ is defined in \eqref{eq:func_Phi}. We define $\cR_{\al, \b}$ and $\cR_{\al, i}$ by 
    \bseq \label{eq:proj_commu}
    \begin{align}
         D ^{\al,0} \ang{f, \Phi _i} _V &= \ang{D ^{\al,0} f, \Phi _i} _V + \cR _{\al, i} (s, X) , 
     \label{eq:proj_commu:a} \\
        D ^{\al, \b} \cP _M f (s, X, V) &= \cP _M D ^{\al, \b} f (s, X, V) + \cR_{\al, \b} (s, X, V) .
         \label{eq:proj_commu:b} 
    \end{align}
    We have the following pointwise estimate on $\cR_{\al, i}$
    for any $N \geq 0$ 
    \begin{align}
        |D^{\al', 0} \cR _{\al, i} (s, X)| &\les _{\al,\al', N} \| 
        \la \vc \ra^{-N} D ^{< |\al| + |\al'|} f (s, X, \cdot) \| _{L ^2 (V)}. \label{eq:proj_commu:c} 
    \end{align}
\eseq 
Moreover, for each $N \ge 0$ and any $\al, \b$, we have 
    \begin{align}
        \label{eqn:commutator-projection}
        \|\ang \vc ^N \cR_{\al, \b} (s, X, \cdot) \| _{L ^2 (V)} &\les _{N, \al, \b} \| \ang \vc ^{-N} D ^{< |\al| + |\b|} f (s, X, \cdot) \| _{L ^2 (V)}, \\
        \label{eqn:commutator-projection-sigma}
        \|\ang \vc ^N \cR_{\al, \b} (s, X, \cdot) \| _\s &\les _{N, \al, \b} \| \ang \vc ^{-N} D ^{< |\al| + |\b|} f (s, X, \cdot) \| _\s.
    \end{align}
    As a consequence, we have the following bound:
    \begin{align}
        \label{eqn: comm-P-Dab-L2}
        \|\cR_{\al,\b} \| _{\cY _l} &\les _{\al, \b} \| D ^{ < |\al| + |\b|} f \| _{\cY _l}, \\
        \label{eqn: comm-P-Dab-sigma}
        \|\cR_{\al,\b} \| _{\cYE} &\les _{\al, \b} \| D ^{< |\al| + |\b|} f \| _{\cYE}.
    \end{align}

    \item \emph{Commuting $D ^{\al, \b}$ and $\dcm$, $\td d _\cM$:} 
    recall that $\dcm$ and $\td d _\cM$ are defined in \eqref{eqn: dcm}. The derivative of $\dcm$ and $\td d _\cM$ can be bounded by 
    \begin{align}
        \left|D ^{\al, \b} \dcm \right| &\les \ang X ^{-r} \ang \vc ^{3}, & 
        \left|D ^{\al, \b} \td d _\cM \right| &\les \ang X ^{-1} \cs \ang \vc ^3.
        \label{eqn: Dab-dcm-bound}
    \end{align}

\end{enumerate}
\end{lemma}

\begin{remark}
    Because $\cP _M + \cP _m = \id$, the commutator with $\cP _M$ is just the negative of the commutator with $\cP _m$.
\end{remark}

Before we prove Lemma \ref{lem: commu-derivative-projection}, we establish the following basic derivative bounds. %

\begin{lemma}[Estimates of the basis and weight]
    \label{lem:basis}
Recall $\cMM = \cs^{-3} \mu(\vc)$ from \eqref{eq:localmax2}. For any multi-indices $\al$, $\b$ and $l \in \R$, we have the following pointwise estimates 
    \begin{align}
        \label{eq: Dab_cs}
        |D ^{\al, \b} \cs ^l| &\les _{\al, \b} \cs ^l, \\
        \label{eq: Dab_vcl}
        |D ^{\al, \b} \ang \vc ^l| &\les _{\al, \b} \ang \vc ^l, \\
        \label{eq:Lams_est}
        |D ^{\al, \b} \Lam| &\les _{\al, \b} \Lam, \\
        \label{eq:cross_pf1}
        |D ^{\al, \b} \Phi _i| &\les _{\al, \b} \cs ^{-3/2} \ang \vc ^{2 + |\b| + 2 |\al|} \mu (\vc) ^{1/2}, \\
        \label{eq:M11/2}
        |D ^{\al, \b} \cMM ^{\pm 1/2}| &\les_{\al, \b} \ang \vc ^{|\b| + 2 |\al|} \cMM ^{\pm 1/2}, \\
        \label{eq:DabM}
        |D ^{\al, \b} \cM| &\les_{\al, \b} \ang \vc ^{|\b| + 2 |\al|} \cM, \\
        \label{eq:DablogM1}
        |D ^{\al, \b} \log \cMM| &\les_{\al, \b} \ang \vc ^2 \text{ if } |\al| + |\b| > 0.
    \end{align}
    For any function $f$, integer $N \ge 0$, and multi-indices 
    $\al, \b$, we have
    \bseq\label{eq:cross_pf2}
    \begin{align}
        & \| \ang \vc ^N \mu (\vc) ^{1/2} \| _{L ^2 (V)} \les _N \cs ^{3/2}, 
        \label{eq:cross_pf2:a} \\
        & |\la f , D ^{\al, \b} \Phi_i \ra _V| \les _{N, \al, \b} \| \ang \vc ^{-N} f \| _{L ^2 (V)}.
        \label{eq:cross_pf2:b}
    \end{align}
    \eseq
\end{lemma}

\begin{proof}
    We start with the proof of \eqref{eq: Dab_cs}-\eqref{eq:Lams_est}. $\cs$ is $V$-independent, so \eqref{eq: Dab_cs} follows Lemma \ref{lem:example_fl} (4). For \eqref{eq: Dab_vcl}, we use induction on $|\al| + |\b|$. Assume \eqref{eq: Dab_vcl} holds for any $l \in \R$ and $|\al| + |\b| \le k$. We will show it holds for $(\al + \al _1, \b + \b _1)$ with $|\al _1| + |\b _1| = 1$. By Corollary \ref{cor: concatenate-derivative}, 
    \begin{align*}
        D ^{\al + \al _1, \b + \b _1} \ang \vc ^l = D ^{\al, \b} D ^{\al _1, \b _1} \ang \vc ^l + c D ^{\prec (\al + \al _1, \b + \b _1)} \ang \vc ^l.
    \end{align*}
    By inductive assumption, the lower order term is bounded as 
    \begin{align*}
        |c D ^{\prec (\al + \al _1, \b + \b _1)} \ang \vc ^l| \les _{\al, \b} \ang \vc ^l.    
    \end{align*}
    For the top order term, note that 
    \begin{align*}
        D ^{\al _1, \b _1} \ang \vc ^l = l \ang \vc ^{l - 2} \vc \cdot D ^{\al _1, \b _1} \vc.
    \end{align*}
    By induction, for any $|\al'| + |\b'| \le k$ we have
    \begin{align*}
        |D ^{\al', \b'} \ang \vc ^{l - 2}| \les _{k} \ang \vc ^{l - 2}.
    \end{align*}
    Together with $|D ^{\al', \b'} \vc| \les _k \ang \vc$ and $|D ^{\al', \b'} D ^{\al _1, \b _1} \vc| \les _k \ang \vc$ using Corollary \ref{cor: Dab-polynomial}, we conclude by Leibniz rule that
    \begin{align*}
        D ^{\al, \b} D ^{\al _1, \b _1} \ang \vc ^l \les \ang \vc ^{l - 2 + 1 + 1} = \ang \vc ^l.
    \end{align*}
    Combined with lower order term, we proved \eqref{eq: Dab_vcl} for $\al + \al _1, \b + \b _1$ and the induction is completed. Because $\Lams = \cs ^{\g + 3} \ang \vc ^{\g + 2}$, \eqref{eq:Lams_est} follows by Leibniz rule.

    Next, we prove \eqref{eq:cross_pf1}-\eqref{eq:DablogM1}. Recall that
    \begin{align*}
        \Phi _i = \cs ^{-3/2} p _i (\vc) \mu (\vc) ^{1/2}
    \end{align*}
    where $p _i$ is a polynomial of degree $\operatorname{deg} p _i \le 2$. By \eqref{eq: Dab_cs} and Lemma \ref{lem:H_deri}, we have 
    \begin{align*}
        D ^{\al, \b} \Phi _i &= \sum _{\substack{\al' \preceq \al\\ \b' \preceq \b}} C ^\al _{\al'} C ^\b _{\b'} \cdot D ^{\al', \b'} \cs ^{-3/2} \cdot D ^{\al - \al', \b - \b'} (p _i (\vc) \mu (\vc) ^{1/2}) \\
        &\les _{\al, \b} \cs ^{-3/2} \ang \vc ^{\operatorname{deg} p _i + |\b| + 2 |\al|} \mu (\vc) ^{1/2}.
    \end{align*}
    \eqref{eq:cross_pf1} is proven. As for \eqref{eq:M11/2}, we have shown $D ^{\al, \b} \cMM ^{1/2} \lesssim \ang \vc ^{|\b| + 2 |\al|} \cMM ^{1/2}$ because $\cMM ^{1/2} = \Phi _0$. By the product rule we have for $|\al| + |\b| > 0$ that
    \begin{align*}
        0 = D ^{\al, \b} (\cMM ^{1/2} \cMM ^{-1/2}) = \sum _{\al' \preceq \al, \b' \preceq \b} C ^\al _{\al'} C ^\b _{\b'} \cdot D ^{\al', \b'} \cMM ^{1/2} \cdot D ^{\al - \al', \b - \b'} \cMM ^{-1/2}.
    \end{align*}
    We can conclude \eqref{eq:M11/2} by induction. By writing $\cM = \cs ^3 \cdot \cMM ^{1/2} \cdot \cMM ^{1/2}$, \eqref{eq:DabM} follows by Leibniz rule and \eqref{eq: Dab_cs}, \eqref{eq:M11/2}. Finally, $\log \cMM = -3 \log \cs + \log \mu (\vc) = -3 \log \cs - \kp _2 |\vc| ^2$, so \eqref{eq:DablogM1} follows Lemma \ref{lem:example_fl} (3) and Corollary \ref{cor: Dab-polynomial}.

    To prove \eqref{eq:cross_pf2:a}, we verify it using a change of variable:
    \begin{align*}
        \int \ang \vc ^{2 N} \mu (\vc) d V = \cs ^3 \int \ang \vc ^{2 N} \mu (\vc) d \vc \les _{N} \cs ^3. 
    \end{align*}
    Estimate \eqref{eq:cross_pf2:b} follows \eqref{eq:cross_pf1}, \eqref{eq:cross_pf2:a} and the Cauchy--Schwarz inequality:
    \begin{align*}
        \int f \cdot D ^{\al, \b} \Phi_i d V \les \cs ^{-3/2} \int |f (V)| \ang \vc ^{-N} \ang \vc ^{N + 2 + |\b| + 2 |\al|} \mu (\vc) ^{1/2} d V \les _{N, \al, \b} \|\ang \vc ^{-N} f \| _{L ^2 (V)}.
    \end{align*}
    We have completed the proof.
\end{proof}

\begin{proof}[Proof of Lemma \ref{lem: commu-derivative-projection}]
\ 
\begin{enumerate}[(1)]
    \item We separate $V \cdot \na _X$ from other terms in $\pa _s + \cT$: 
    \begin{align*}
		[\cP _m, \pa _s + \cT] f &= [\cP _m, \pa _s + \bcx X \cdot \na _X + \bcv V \cdot \na _V] f + [\cP _m, V \cdot \na _X] f.
	\end{align*}
	The second commutator can be computed directly as 
	\begin{align*}
		[\cP _m, V \cdot \na _X] f & = \cP _m [(V \cdot \na _X) f] - V \cdot \na _X \cP _m f \\
		& = \cP _m [(V \cdot \na _X) (\cP _m + \cP _M) f] - (\cP _M + \cP _m) [(V \cdot \na _X) \cP _m f] \\
		& = \cP _m [(V \cdot \na _X) \cP _M f] - \cP _M [(V \cdot \na _X) \cP _m f].
	\end{align*}
	For the first commutator, we observe that the projection operator commutes with the scaling field and time derivative:
	\beq\label{eq:proj_scaling}
	\bal
	\cMM ^{1/2} \pa _s (\cMM ^{-1/2} \cP _M f) &= \cP _M (\cMM ^{1/2} \pa _s (\cMM ^{-1/2} f)), \\
	\cMM ^{1/2} (X \cdot \na_X) (\cMM ^{-1/2} \cP _M f) &= \cP _M (\cMM ^{1/2} (X \cdot \na_X) (\cMM ^{-1/2} f)), \\
	\cMM ^{1/2} (V \cdot \na_V) (\cMM ^{-1/2} \cP _M f) &= \cP _M (\cMM ^{1/2} (V \cdot \na_V) (\cMM ^{-1/2} f)).
	\eal
	\eeq
	This is because $\cMM ^{1/2} \Phi _i$ is in $\operatorname{Span} \{1, V _i, |V| ^2\}$ (with $X$-dependence), and $\pa _s$, $X \cdot \na _X$, $V \cdot \na _V$ maps this span to itself. Recall $\dcm$ and $\td d _\cM$ defined in \eqref{eqn: dcm}. 
	\begin{align*}
		& \cMM ^{1/2} (\pa _s + \bcx X \cdot \na _X + \bcv V \cdot \na _V) \cMM ^{-1/2} f \\
		& \qquad = (\pa _s + \bcx X \cdot \na _X + \bcv V \cdot \na _V) f -\frac12 (\pa _s + \bcx X \cdot \na _X + \bcv V \cdot \na _V) \log \cMM \cdot f \\
		& \qquad = (\pa _s + \bcx X \cdot \na _X + \bcv V \cdot \na _V) f - \left( \dcm + \td d _\cM - \frac32 \bcv \right) f.
	\end{align*}
	As shown in \eqref{eq:proj_scaling}, the operator $\cMM ^{1/2} (\pa _s + \bcx X \cdot \na _X + \bcv V \cdot \na _V) \cMM ^{-1/2}$ commutes with the projection $\cP _M$, so it also commutes with $\cP _m$. We deduce
	\begin{align*}
		[\cP _m, \pa _s + \bcx X \cdot \na _X + \bcv V \cdot \na _V] f &= \left[\cP _m, \dcm + \td d _\cM - \frac32 \bcv\right] f \\ 
		&= \left[\cP _m, \dcm + \td d _\cM \right] f,
	\end{align*}
	because scalar multiplication commutes with $\cP _m$. The remaining computation is the same as the $V \cdot \na _X$ part.

    \item Note that $D ^{\al, \b} V = \cs \ee _i$ when $\al = 0$, $\b = \ee _i$, and $D ^{\al, \b} V = 0$ when $\al \ne 0$ or when $|\b| \ge 2$. Therefore, 
    \begin{align*}
        D ^{\al, \b} (V \cdot \na _X f) &= \sum _{\substack{\al' \preceq \al\\\b' \preceq \b}} C ^{\al} _{\al'} C ^{\b} _{\b'} \cdot D ^{\al', \b'} V \cdot D ^{\al - \al', \b - \b'} \na _X f \\
        &= V \cdot D ^{\al, \b} \na _X f + \sum _i C ^\b _{\ee _i} D ^{0, \ee _i} V \cdot D ^{\al, \b - \ee _i} \na _X f \\
        &= V \cdot \na _X D ^{\al, \b} f - V \cdot \na _X \log (\varphi _1 ^{|\al|} \cs ^{|\b|}) D ^{\al, \b} f + \sum _i \b _i \cs D ^{\al, \b - \ee _i} \pa _{X _i} f \\
        &= V \cdot \na _X D ^{\al, \b} f + \cO (|V| \ang X ^{-1}) D ^{\al, \b} f + \cs \varphi _1 ^{-1}  \sum _i \b _i D ^{\al + \ee _i, \b - \ee _i} f.
    \end{align*}
    The bound for the first commutator comes from $\varphi _1 ^{-1} \les \ang X ^{-1}$ and $|V| \les \cs \ang \vc$. Note that we only sum $D ^{\al + \ee _i, \b - \ee _i}$ for $\b _i > 0$, therefore $\b -\ee _i \succeq 0$.
    
    Similarly, for the second commutator we have 
    \begin{align*}
         [\pa _s + \cT - & V \cdot \na _X, D ^{\al, \b}] f  \\ 
            &= \left\{(\pa _s + \cT - V \cdot \na _X) \log (\varphi _1 ^{|\al|} \cs ^{|\b|}) - \bcx |\al| - \bcv |\b| \right\} D ^{\al, \b} f \\
            &= \left\{|\b|(\pa_s+ \cT-V\cdot \nabla_X) \log \cs+ \bcx |\al| X \cdot \na _X \log \varphi _1 - \bcx |\al| - \bcv |\b| \right\} D^{\al,\b} f\\
            & = \left\{|\b|(\bcv +\cO( \ang{X}^{-r}\ang{\vc}))-|\b|V\cdot \nabla_X \log \cs + \bcx |\al| (1 + \cO (\ang X ^{-1}) - \bcx |\al| - \bcv |\b| \right\} D^{\al,\b} f\\
            & = \left\{
                \cO( \ang{X}^{-r}\ang{\vc}) + \cO (|V| \ang X ^{-1}) +\cO (\ang X ^{-1}) %
            \right\} D^{\al,\b} f\\
            &= \cO (\cs \ang X ^{-1} \ang \vc + \ang X^{-1})  D ^{\al, \b} f.
        \end{align*}
        We used the equation \eqref{eq:ds_t_vc} and the estimate of $\vp_1$ from \eqref{eq:wg_asym}:
        \begin{align*}
                    X \cdot \na _X \log \varphi _1 = 1 + \cO (\ang X ^{-1}).
        \end{align*}

    \item 
    First, we consider pure $X$ derivative, i.e. the case $|\b| = 0$. 
    Recall $\cR_{\al, i}$ from \eqref{eq:proj_commu:a}. By Leibniz rule, 
    \begin{align*}
        \cR _{\al, i} = \sum _{\al' \prec \al} C ^\al _{\al'} \la D ^{\al', 0} f, D ^{\al - \al', 0} \Phi _i \ra _V .
    \end{align*}
    Therefore 
    \begin{align*}
        D ^{\al _1, 0} \cR _{\al, i} &= \sum _{\al _1' \preceq \al _1} \sum _{\al' \prec \al} C ^{\al _1} _{\al _1'} \cdot C ^\al _{\al'} \la D ^{\al _1', 0} D ^{\al', 0} f, D ^{\al _1 - \al _1'} D ^{\al - \al', 0} \Phi _i \ra _V \\
        &= \sum _{\substack{\al _2 \prec \al _1 + \al \\ \al _3 \preceq \al _1 + \al}} c _{\al _2, \al _3} \la D ^{\al _2, 0} f, D ^{\al _3, 0} \Phi _i \ra _V,
    \end{align*}
    with $c _{\al _2, \al _3} \in \bF$, thanks to Corollary \ref{cor: concatenate-derivative}. We bound the pairing using the $L ^2$ norm and the $\sigma$ norm:
    \begin{align*}
        \la D ^{\al _2, 0} f, D ^{\al _3, 0} \Phi _i \ra _V &\le \| \ang \vc ^{-N} D ^{\al _2, 0} f \| _{L ^2 (V)} \| \ang \vc ^N D ^{\al _3, 0} \Phi _i \| _{L ^2 (V)}, \\
        \la D ^{\al _2, 0} f, D ^{\al _3, 0} \Phi _i \ra _V &\le \| \ang \vc ^{-N} D ^{\al _2, 0} f \| _\s \| \Lams ^{-\f12} \ang \vc ^ND ^{\al _3, 0} \Phi _i \| _{L ^2 (V)}.
    \end{align*}
    By \eqref{eq:cross_pf1} and \eqref{eq:cross_pf2:a} we know 
    \begin{align}
        \label{eqn:DabPhi-L2}
        \| \ang \vc ^N D ^{\al, \b} \Phi _i \| _{L ^2 (V)} \les _{\al, \b, N} 1.
    \end{align}
    Let us also compute the $\sigma$ norm and weighted norm for the derivative of the basis:
    \beq\bal
        \label{eqn:DabPhi-sigma}
        \| \ang \vc ^N D ^{\al, \b} \Phi _i \| _\sigma ^2 &= \cs ^{-3} \| \ang \vc ^N p (\vc) \mu (\vc) ^{1/2} \| _\sigma ^2 \\
        &\le \cs ^{-3} \B(
            \cs ^{\gamma + 5} \int \ang \vc ^{\gamma + 2} |\na _V (\ang \vc ^N p (\vc) \mu (\vc) ^{1/2}| ^2 d V \\
            & \hspace{4em} + \cs ^{\gamma + 3} \int \ang \vc ^{\gamma + 2} |(\ang \vc ^N p (\vc) \mu (\vc) ^{1/2}| ^2 d V 
        \B) \\
        & \les _{\al, \b, N} \cs ^{\gamma + 3}.
    \eal\eeq
    Similarly, 
    \beq\bal
        \label{eqn:DabPhi-sigma*}
        \| \Lams ^{-\f12} \ang \vc ^N D ^{\al, \b} \Phi _i \| _{L ^2} ^2
        &\les \cs ^{-3} \left\{
            \cs ^{-\g - 3} \int \ang \vc ^{-\g - 2} |\ang \vc ^N p (\vc) \mu (\vc) ^{1/2}| ^2 d V 
        \right\} \les \cs ^{-\g - 3}.
    \eal\eeq
    In summary, we conclude 
    \begin{align}
        \label{eqn:commutator-bracket}
        | D ^{\al _1, 0} \cR _{\al, i} (s, X) | 
        & \les _{N, \al, \b} \| \ang \vc ^{-N} D ^{< |\al| + |\al _1|} f (s, X, \cdot) \| _{L ^2 (V)} \\
        \label{eqn:commutator-bracket-sigma}
        | D ^{\al _1, 0} \cR _{\al, i} (s, X) | 
        & \les_{N, \al, \b} \cs ^{-\f{\g + 3}2} \| \ang \vc ^{-N} D ^{< |\al| + |\al _1|} f (s, X, \cdot) \| _\s.
    \end{align}
    which proves \eqref{eq:proj_commu:c}. 
    
    Now we compute
    \begin{align*}
        \cR_{\al, 0} (s, X, V) &= D ^{\al,0} \cP _M f - \cP _M D ^{\al,0} f \\
        &= \sum _i D ^{\al,0} (\ang{f, \Phi _i} \Phi _i) - \ang{D ^{\al,0} f, \Phi _i} \Phi _i \\
        &= \sum _i \sum _{\al' \preceq \al} C ^\al _{\al'} \cdot D ^{\al',0} \ang{f, \Phi _i} D ^{\al - \al',0} \Phi _i - \ang{D ^{\al,0} f, \Phi _i} \Phi _i \\
        &= \sum _i \left(\sum _{\al' \prec \al} C ^\al _{\al'} \cdot \ang{D ^{\al',0} f, \Phi _i} \cdot D ^{\al - \al',0} \Phi _i + \sum _{\al' \preceq \al} C ^\al _{\al'} \cdot \cR _{\al', i} \cdot D ^{\al - \al',0} \Phi _i\right).
    \end{align*}

For any $\al_1, \b_1$, applying Leibniz rule, \eqref{eq:cross_pf1} to $\Phi_i$, \eqref{eq:proj_commu:c} to $\cR_{\al^{\pr}, i}$, and \eqref{eq:cross_pf2:b}, 
we obtain 
\[
\bal 
|D^{\al_1, \b_1}  ( \ang{D ^{\al',0} f, \Phi _i} \cdot D ^{\al - \al',0} \Phi _i ) | 
& \les_{\al, \al_1, \b_1}  \| D^{ \leq |\al^{\pr}| + |\al_1| + |\b_1|} f \|_{L ^2 (V)} 
\cs^{-3/2} \ang \vc^{2 + 2 |\al| + 2 |\al_1| + |\b_1|}  \mu(\vc)^{1/2} ,  \\ 
 |D^{\al_1, \b_1} ( \cR_{\al^{\pr}, i}  \cdot D ^{\al - \al',0} \Phi _i  )  | &
 \les_{\al, \al_1, \b_1}   \| D^{ < |\al^{\pr}| + |\al_1| + |\b_1|} f \|_{L ^2 (V)}  \cs^{-3/2} \ang \vc^{2 + 2 |\al| + 2 |\al_1| + |\b_1|}  \mu(\vc)^{1/2} .
\eal 
\]
    Estimate \eqref{eqn:commutator-projection} follows \eqref{eq:cross_pf2:b}, \eqref{eqn:commutator-bracket}, and \eqref{eqn:DabPhi-L2}, whereas \eqref{eqn:commutator-projection-sigma} follow \eqref{eqn:DabPhi-sigma} and \eqref{eqn:commutator-bracket-sigma}.
    
    Now we consider derivatives in both $X$ and $V$. Note that 
    \begin{align*}
        D ^{\al, \b} \cP _M f &= \sum _i \sum _{\al' \preceq \al} C ^\al _{\al'} \cdot D ^{\al',0} \ang{f, \Phi _i} \cdot D ^{\al - \al', \b} \Phi _i \\
        &= \sum _i \sum _{\al' \preceq \al} C ^\al _{\al'} \cdot (\ang{D ^{\al',0} f, \Phi _i} + \cR _{\al', i}) \cdot D ^{\al - \al', \b} \Phi _i,
    \end{align*}
    so we can obtain the weighted $L ^2$ bound for $D ^{\al, \b} \cP _M f$. Using integration by parts, 
    \begin{align*}
        \cP _M D ^{\al, \b} f = \sum _i \ang{D ^{\al, \b} f, \Phi _i} \Phi _i = (-1) ^{|\b|} \sum _i \ang{D ^{\al,0} f, D ^{0, \b} \Phi _i} \Phi _i,
    \end{align*}
    so we can obtain the weighted $L ^2$ bound for $\cP _M D ^{\al, \b} f$ as 
    \begin{align}
        \label{eqn:commutator-projection-v}
        \|\ang \vc ^N D ^{\al, \b} \cP _M f\| _{L ^2 (V)} + \|\ang \vc ^N \cP _M D ^{\al, \b} f \| _{L ^2 (V)} \les _N \| \ang \vc ^{-N} D ^{\le |\al|} f (s, X, \cdot) \| _{L ^2 (V)}.
    \end{align}
    So \eqref{eqn:commutator-projection} holds provided $\beta \succ 0$. \eqref{eqn:commutator-projection-sigma} holds similarly. \eqref{eqn: comm-P-Dab-L2}, \eqref{eqn: comm-P-Dab-sigma} follow \eqref{eqn:commutator-projection}, \eqref{eqn:commutator-projection-sigma} after integration in $X$ with weight $\ang X ^\eta$.

    \item By a direct computation, 
    \begin{align*}
        \dcm &= \f12 (\pa _s + \cT) \log \cMM + \f32 \bcv \\
        &= \f12 (\pa _s + \cT) \log \cM - \f32 (\pa _s + \cT) \log \cs + \f32 \bcv = \f12 \cM ^{-1} \eM - \f32 [(\pa _s + \cT) \log \cs - \bcv].
    \end{align*}
    Lemma \ref{lem:ds_t} shows $\eM / \cM$ is a is a class $\bF ^{-r}$ polynomial of $\vc$ of degree 3. By \eqref{eqn: ds_t_logcs}, we know $(\pa _s + \cT) \log \cs - \bcv$ is a class $\bF ^{-r}$ polynomial of $\vc$ of degree 1. Therefore, $\dcm$ is a class $\bF ^{-r}$ polynomial of $\vc$ of degree 3. Using Corollary \ref{cor: Dab-polynomial}, we know $D ^{\al, \b} \dcm$ is a class $\bF$ polynomial of $\vc$ with degree at most $3 - |\b|$, thus the conclusion follows.

Since the weight $\vp_1$ in \eqref{eq:deri_wg} satisfies $\vp_1 \asymp \ang X$ by \eqref{eq:wg_asym}, the computation for $\td d _\cM$ is straightforward:
    \begin{align*}
        D ^{\al, \b} \td d _\cM &= \f12 V \cdot D ^{\al, \b} \na _X \log \cMM + \f12 \sum _{i} C ^\b _{\ee _i} \cdot D ^{0, \ee _i} V \cdot D ^{\al, \b - \ee _i} \na _X \log \cMM
        \\
        &= \f12 \sum _i V _i \vp_1^{-1} \cdot D ^{\al + \ee _i, \b} \log \cMM + \f12 \sum _{i} \b _i \cs \vp_1^{-1} \cdot D ^{\al + \ee _i, \b - \ee _i} \log \cMM \\
        &= \cO (\ang \vc \cs \ang X ^{-1} |D ^{\le |\al| + |\b| + 1} \log \cMM|).
    \end{align*}
    By using \eqref{eq:DablogM1} we conclude the proof.
    \end{enumerate}
\end{proof}

\subsection{Estimates between \texorpdfstring{$\cX$}{X} and \texorpdfstring{$\cY$}{Y} norms}

\begin{lemma}\label{lem:macro_UPB}
Let $\tw = (\tu, \tp, \tb)$, $\cF_{M}$ be the operator defined in \eqref{eq:macro_UPB} and 
$\tFM = \cF_M(\tw)$.  For any multi-indices $\al, \b $ and $N \in \R$, we have the following relation 
\bseq\label{eq:macro_UPB_est}
\beq\bal
\int | D_X^{\al} \cF_M (\tu, \tp, \tb) |^2  d V
& = \kp \left(  
|D_X^{\al} \tu|^2 + 
|D_X^{\al} \tp|^2 + \f{3}{2} | D_X^{\al} \tb|^2 \right) 
\\
& \qquad + \one_{|\al| \geq 1} O( | D_X^{\leq |\al|} \tw| \cdot  | D_X^{\leq |\al|-1} \tw| )  ,
\label{eq:macro_UPB_est:a}
\eal\eeq
where $\kp = \f{5}{3}$,  and the following equivalence
\begin{align}
\int \ang \vc^{2N}  | D^{\al, \b}  \cF_M(\tu, \tp, \tb) |^2  d V
  \les_{\al, \b, N}   | D_X^{ \leq |\al| } (\tu, \tp, \tb) |^2 
  \les_{\al }  \int  | D_X^{ \leq |\al|} \cF_M (\tu, \tp, \tb) |^2  d V. 
  \label{eq:macro_UPB_est:b}
 \end{align}
 \eseq
 In particular, for any $\eta \in \R, k \in \BZ_+$, we have the following estimates between the $\cX$-norm and $\cY$-norm for the macro-perturbation 
\bseq\label{eq:macro_Y}
\begin{align} 
\| \cF_M(\tu, \tp, \tb) \|_{\cYe ^k} & \les  \| (\tu, \tp, \tb) \|_{\cX_{\eta}^k}, 
\label{eq:macro_Y:a} \\
\| \cF_M(\tu, \tp, \tb) \|_{\cYE ^k} 
  & \les  \| (\tu, \tp, \tb) \|_{\cX_{\eta}^k } .
\label{eq:macro_Y:b}
\end{align}
\eseq
\end{lemma}

\begin{proof}

Denote $\tw = (\tu, \tp, \tb ) $ and  $\tFM = \cF_M( \tu, \tp, \tb )$. 
Since $\tw$ only depends on $X$, using the relation \eqref{eq:macro_UPB}, the Leibniz rule, and \eqref{eq:cross_pf1}, for any multi-indices $\al, \b$, we obtain 
\bseq\label{eq:macro_UPB_pf1}
\beq\label{eq:macro_UPB_pf1:a}
	D^{\al,\b} \tFM = 
	D_X^{\al} ( \tp + \tB ) \cdot D_V^{\b} \Phi_0 + \kp^{1/2} 	D_X^{\al} \tu_i  \cdot D_V^{\b} \Phi_i 
  + \sqrt{\f{1}{6}} 	D_X^{\al} (2 \td P - 3 \tb)  \cdot  D^{\b} \Phi_4
  + I ,
\eeq
where the error term $I$ satisfies 
\beq\label{eq:macro_UPB_pf1:b}
	|I| \les_{\al, \b}  \sum_{ 0 \leq i \leq |\al | -1}  
	 | D_X^{ i} \tw | \cdot ( \sum_{ 1\leq j \leq 5} |D_{X, V}^{\leq |\al| + |\b| - i} \Phi_j | ) 
\les_{\al, \b} | D_X^{\leq |\al| - 1 } \tw| \cdot \cs^{-3/2} \ang V^{2 + |\b| + 2|\al|} \mu(\vc)^{1/2}.
\eeq
For any $ N \in \R$, using \eqref{eq:cross_pf2:a}, we obtain 
\beq\label{eq:macro_UPB_pf1:c}
	\| \ang V^{N} I\|_{L ^2 (V)} \les_{\al, \b, N} \one_{|\al| \geq 1} | D_X^{\leq |\al| - 1 } \tw| \les   \one_{|\al| \geq 1} | D_X^{\leq |\al|  } \tw|.
\eeq
\eseq

\paragraph{Proof of \eqref{eq:macro_UPB_est:a}}
 Since $\Phi_i$ are orthonormal \eqref{eq:func_Phi} and $\kp = \f53$, applying  \eqref{eq:macro_UPB_pf1} with $\b = 0$, we prove 
\[
\bal
 \|	D_X^{\al} \tFM\|_{L ^2 (V)}^2 & = |	D_X^{\al} ( \tp + \tB )|^2 
 +  \sum_i \kp  	|D_X^{\al} \tu_i|^2
 +  \f16	| D_X^{\al} (2 \td P - 3 \tb)  |^2 + 
  O( \| I\|_{L ^2 (V)}  | D_X^{\leq |\al|} \tw | + \| I\|_{L ^2 (V)}^2  )  \\
  & =  \kp (  |D_X^{\al} \tu|^2 + 
| D_X^{\al} \tp|^2 + \f{3}{2} | D_X^{\al} \tb|^2 )
+  \one_{|\al| \geq 1}  O(  | D_X^{\leq |\al|} \tw | \cdot  | D_X^{\leq |\al| - 1}  \tw |   ) .
 \eal
\]
Thus, we prove \eqref{eq:macro_UPB_est:a}. 
Using induction on $ k \geq 0$  and \eqref{eq:macro_UPB_est:a}, we obtain
\beq\label{eq:macro_UPB_equiv1}
   \| D_X^{ \leq k} \tFM \|_{L ^2 (V)}^2 \asymp_k \| D_X^{\leq k} \tw| .
\eeq

\paragraph{Proof of \eqref{eq:macro_UPB_est:b}}
For the main term on the right hand side of \eqref{eq:macro_UPB_pf1:a}, applying estimates similar to $I$ in  \eqref{eq:macro_UPB_pf1:b}, \eqref{eq:macro_UPB_pf1:c}, we prove the first bound in  \eqref{eq:macro_UPB_est:b}: 
\[
	\| \ang V^N D^{\al, \b} \tFM\|_{L ^2 (V)} 
	\les_{\al, \b}
	| D_X^{\leq |\al|  } \tw| \cdot \cs^{-3/2}  \| \ang V^{2 + |\b| + 2|\al| + N}  \mu(\vc)^{1/2} \|_{L ^2 (V)}
	  \les_{\al, \b, N}	| D_X^{\leq |\al|  } \tw| .
\]
The second bound in \eqref{eq:macro_UPB_est:b} follows from \eqref{eq:macro_UPB_equiv1}.

\paragraph{Proof of \eqref{eq:macro_Y}}

Recall the $\cX$-norm from \eqref{norm:Xk}, $\cY, \cYE ^k$-norms from \eqref{norm:Y}. 
For the coefficient in $\cYE^k$ norm, we note that $\cs \les 1$ \eqref{eq:dec_S}. 
Using \eqref{eq:macro_UPB_est:b} with $N= 0$ and $N =\f{\g+3}{2}$, we prove 
\[
\bal
	\| \cF_M(\tw) \|_{ \cYe^k }^2
	+ 	\| \cF_M(\tw) \|_{ \cYE^k }^2
	& \les_k \sum_{ |\al| \leq k, |\b| \leq k+1 } \| \ang X^{\eta/2} \ang \vc^{  \f{\g+2}{2} } 
D^{\al ,\b} \cF_M(\tw) \|_{L^2}^2  \\
	& \les_k \| \la X \ra^{ \eta/2 }  D_X^{\leq k} \tw \|_{L^2(X)}^2
	\les \| \tw\|_{\cX^k_{\eta}}^2,
	\eal
\]
where in the last inequality we have used \eqref{eq:grad_est} with $\psi_n = \ang X^{\eta}  \vp_1^{n}$
with $\vp_1$ defined in Lemma \ref{lem:wg}, $\b_n = 2 n + \eta$ and $\nu = 1$. This completes the proof.
\end{proof}

\subsection{Proof of Lemma \ref{lem: diffusion-innerprod} on weighted diffusion term}
\label{sec: proof-diffusion-innerprod}

Recall the weighted operator from \eqref{eq:deri_wg} and the weighted diffusion from \eqref{eq: weighted-laplacian}:
\begin{align} \label{eq: weighted-laplacian-recall}
    \D _W F & = -  \nu ^{-1} \ang X ^{2} \ang \vc ^{4} F   +   \sum _{|\al_1|  +  |\b_1| = 1}  \ang X ^{1-\etab} \ang \vc ^2  
    \pa_X^{\al_1} \pa_V^{\b_1}  \B( \vp_1^{2 |\al_1| } \cs^{2 |\b_1|} \ang X^{\etab}  \pa_X^{\al_1} \pa_V^{\b_1}  (\ang X \ang \vc ^2 F)  \B)  \notag \\ 
    & := \D_{L^2} F + \D_{H^1} F .
\end{align}

\begin{proof}[Proof of Lemma \ref{lem: diffusion-innerprod}]
Denote 
\beq\label{eq:diff_g}
\quad g = \ang X \ang \vc ^2 h . 
\eeq

\paragraph{Estimate of $k=0$}

First, consider $k = 0$. Recall the $\cYb$ norms from \eqref{norm:Y}. By definition, we yield 
\[
  \la  \D_W h, h \ra_{\cYb} %
  =  -\nu^{-1} \iint \ang X^{2 + \etab} \ang \vc^{4} h^2 d X d V
  + \iint \D_{H^1} h \cdot h \ang X^{ \etab} dX dV := I_1 + I_2.
\]

For $I_2$, using the notation $g$ \eqref{eq:diff_g} and integration by parts, we obtain 
\[
\bal 
I_2 & = \sum_{|\al_1| + |\b_1| = 1}  \iint   \pa_X^{\al_1} \pa_V^{\b_1}  \B( \vp_1^{2 |\al_1| } \cs^{2 |\b_1|} \ang X^{\etab}  \pa_X^{\al_1} \pa_V^{\b_1}  (\ang X \ang \vc ^2 h)  \B)    \cdot  \ang X  \ang \vc^2 h   \\ 
& = - \sum_{|\al_1| + |\b_1| = 1} \iint \vp_1^{2 |\al_1| } \cs^{2 |\b_1|} \ang X^{\etab} 
\pa_X^{\al_1} \pa_V^{\b_1} g \cdot \pa_X^{\al_1} \pa_V^{\b_1} g d X d V \\
& =  - \sum_{|\al_1| + |\b_1| = 1}  \iint |D^{\al_1, \b_1} g |^2 \ang X^{\etab}  d X d V.
\eal  
\]
Combining the above estimates, using the definition of $\cYb^1$ \eqref{norm:Y} 
and $1! = 0! = 1$,  we prove Lemma \ref{lem: diffusion-innerprod} for $k=0$:
\[
  \la  \D_W h, h \ra_{\cYb}  = 
  - \nu^{-1}  \int |g|^2  \ang X^{\etab}  d X d V 
   - \sum_{|\al_1| + |\b_1| = 1}  \iint |D^{\al_1, \b_1} g |^2 \ang X^{\etab}  d X d V 
  = -\| g \|_{\cYb^1}^2.
\]

\paragraph{Estimate for $k \geq 1$}

For higher order estimates, by Leibniz rule, we  rewrite $\D_W$ \eqref{eq: weighted-laplacian-recall} as 
\beq\label{eq:diffusion_main}
\bal 
\D_W h & = 
\sum _{|\al_1|  +  |\b_1| = 1} \ang X \ang \vc ^2  
 \vp_1^{2|\al_1|}\cs^{ 2 |\b_1|} \pa_X^{2 \al_1} \pa_V^{2 \b_1} g 
 + c_1  \ang X \ang \vc ^2 D^{\preceq (\al_1, \b_1 ) } g \\ 
 & = \sum _{|\al_1|  +  |\b_1|  = 1} \ang X  \ang \vc ^2 D^{2\al_1, 2 \b_1} g   + c_{1}  \ang X  \ang \vc ^2 D^{\preceq (\al_1, \b_1 ) } g ,
\eal 
\eeq
where $c_{1}$ denotes generic bounded functions containing functions like %
$\f{D ^{\al, \b} f}{  f}, f = \ang X, \ang \vc, \cs, \vp_1$  with bounds only depending on $k$. 
We use similar notations $c$ below, which may change from line to line. Note that by applying the Leibniz rule iteratively, we obtain
\begin{align}\label{commh1}
    \ang X \ang \vc ^2 D ^{\al, \b} h &= D ^{\al, \b} (\ang X \ang \vc ^2 h) + c D ^{\prec (\al, \b)} (\ang X \ang \vc ^2 h) = D ^{\al, \b} g + c D ^{\prec (\al, \b)} g.
\end{align}

For the main term in \eqref{eq:diffusion_main}, we take a single term $|\al _1| + |\b _1| = 1$ and obtain 
\begin{align}
    &\left\la 
        \ang X \ang \vc ^2 D ^{2 \al_1, 2 \b_1} g,  h
    \right\ra _{\cYb ^k} 
    \label{eqn: one-term-in-yk-0} \\
    & \qquad = \sum _{|\al| + |\b| \le k} \nu ^{|\al| + |\b| - k}      \f{|\al|!}{\al!}  \int \ang X ^{\etab} D ^{\al, \b} \left[\ang X  \ang \vc ^2 D ^{2 \al _1, 2 \b _1} g\right] \cdot D ^{\al, \b}\; h d V d X.
    \label{eqn: one-term-in-yk}
\end{align}
We start by commuting $D ^{\al, \b}$ with the weights $\ang X  \ang \vc ^2$: 
\begin{align}
    D ^{\al, \b} &\left[ 
        \ang X  \ang \vc ^2 D ^{2 \al _1,  2 \b _1} g
    \right] \notag \\ 
    &= \ang X  \ang \vc ^2 D ^{\al, \b} D ^{2 \al _1,  2 \b _1} g  %
    + c \ang X  \ang \vc ^2 D ^{\prec (\al, \b)} D ^{2 \al _1,  2 \b _1} g \nonumber\\
    & = \ang X  \ang \vc ^2 D ^{\al + 2 \al _1, \b + 2 \b _1} g %
    + c \ang X  \ang \vc ^2 D ^{\prec (\al + 2 \al _1, \b + 2 \b _1)}  g \nonumber\\
    & = \ang X  \ang \vc ^2 D ^{\al _1, \b _1} D ^{\al + \al _1, \b + \b _1} g %
    + c \ang X  \ang \vc ^2 D ^{\prec (\al + 2 \al _1, \b + 2 \b _1)}  g.\label{first_comm}
\end{align}

Therefore, one term in \eqref{eqn: one-term-in-yk} can be computed as 
\begin{align}\label{visc_Dalphabeta_1}
    &  \int \ang X ^{\etab} D ^{\al, \b} \left[ 
        \ang X  \ang \vc ^2 D ^{2 \al _1,  2 \b _1} g
    \right] D ^{\al, \b} h \; d V d X  \\
    \notag
    &= \int \ang X ^{\etab} \left(  D ^{\al _1, \b _1} D ^{\al + \al _1, \b + \b _1} g  +c  D ^{\prec (\al + 2 \al _1, \b + 2 \b _1)}  g \right)   {\underbrace{\ang X  \ang \vc ^2 D ^{\al, \b} h}_{\eqref{commh1}}} d V d X \\ 
    \notag
    &= \int \ang X ^{\etab} (  D ^{\al , \b } g +c D ^{\prec (\al, \b)} g) \left(  D ^{\al _1, \b _1} D ^{\al + \al _1, \b + \b _1} g  +c  D ^{\prec (\al + 2 \al _1, \b + 2 \b _1)}  g \right) d V d X \\ 
    &=
    \int \ang X ^{\etab} D ^{\al _1, \b _1} D ^{\al + \al _1, \b + \b _1} g \cdot D ^{\al, \b} g \; d V d X  \label{visc_Dalphabeta_2}\\
    & \qquad
    + \int \ang X ^{\etab} D ^{\al _1, \b _1} D ^{\al + \al _1, \b + \b _1} g \cdot c D ^{\prec (\al, \b)} g \; d V d X  \label{visc_Dalphabeta_3}\\
    & \qquad 
    + \int \ang X ^{\etab} D ^{\prec (\al + 2 \al _1, \b + 2 \b _1)} g \cdot c D ^{\preceq (\al, \b)} g \; d V d X.\label{visc_Dalphabeta_4}
\end{align}
We now deal with \eqref{visc_Dalphabeta_2}
using integration by parts: 
\begin{align*}
   \eqref{visc_Dalphabeta_2} = &\int \ang X ^{\etab} D ^{\al _1, \b _1} D ^{\al + \al _1, \b + \b _1} g D ^{\al, \b} g \; d V d X \\
    & = - \int \ang X ^{\etab} D ^{\al + \al _1, \b + \b _1} g D ^{\al _1, \b _1} D ^{\al, \b} g \; d V d X %
    + \int c \ang X ^{\etab} D ^{\al + \al _1, \b + \b _1} g D ^{\al, \b} g \; d V d X \\
    & = - \int \ang X ^{\etab} | D ^{\al + \al _1, \b + \b _1} g | ^2 \; d V d X  %
    + \int c \ang X ^{\etab} D ^{\al + \al _1, \b + \b _1} g D ^{\prec (\al + \al _1, \b + \b _1)} g \; d V d X.
\end{align*}
We apply the Cauchy--Schwarz inequality in the last integral, and we have 
\begin{align}\nonumber
    & \int c \ang X ^{\etab} D ^{\al + \al _1, \b + \b _1} g D ^{\prec (\al + \al _1, \b + \b _1)} g \; d V d X \\
    & \qquad \le \f18 \int \ang X ^{\etab} | D ^{\al + \al _1, \b + \b _1} g | ^2 \; d V d X + C _{\al, \b} \int \ang X ^{\etab} | D ^{\prec (\al + \al _1, \b + \b _1)} g | ^2 \; d V d X. \label{remainder} 
\end{align}
Summarizing we get %
\begin{align*}
    \eqref{visc_Dalphabeta_2} \le  & -\f78 \int \ang X ^{\etab} | D ^{\al + \al _1, \b + \b _1} g | ^2 \; d V d X + C _{\al, \b} \int \ang X ^{\etab} |D ^{\le |\al| + |\b|} g | ^2 \; d V d X. \nonumber
\end{align*}
Cauchy--Schwarz for \eqref{visc_Dalphabeta_3} and \eqref{visc_Dalphabeta_4} yield a similar bound as in (\ref{remainder}). We get 
\begin{align}
    \eqref{visc_Dalphabeta_1} \le  & -\f34  \int \ang X ^{\etab} | D ^{\al + \al _1, \b + \b _1} g | ^2 \; d V d X +3  C _{\al, \b} \int \ang X ^{\etab} | D ^{\le |\al| + |\b|} g | ^2 \; d V d X.   \label{visc_1} 
\end{align}
Using this, we take the summation in \eqref{eqn: one-term-in-yk}, we have 
\begin{align*}
    \eqref{eqn: one-term-in-yk} &\le \sum _{|\al| + |\b| \le k} \nu ^{|\al| + |\b| - k} \left(
    -\f12  \int \ang X ^{\etab} | D ^{\al + \al _1, \b + \b _1} g | ^2 \; d V d X +3  C _{\al, \b} \int \ang X ^{\etab} | D ^{\le |\al| + |\b|} g | ^2 \; d V d X\right) \\
    &\le -\f34 \sum _{|\al| + |\b| \le k} \nu ^{|\al| + |\b| - k}\int \ang X ^{\etab} | D ^{\al + \al _1, \b + \b _1} g | ^2 \; d V d X + C _k \| g \| _{\cYb ^k} ^2.
\end{align*}

The lower order term $\ang X \ang \vc^2 D^{\preceq (\al_1, \b_1)} g$ in \eqref{eq:diffusion_main} can be estimated similarly by 
interpolation 
\beq\label{eq:diffusion_low_est}
 \B\la  c_1  \ang X  \ang \vc ^2 D^{\preceq (\al_1, \b_1 ) } g, h \B\ra_{\cYb^k}
 \leq  \f{1}{32} \| g \|_{\cYb^{k+1}}^2 + C_k  \| g \|_{\cYb^{k}}^2.
\eeq
Since $\nu$ is a fixed constant, we treat it as an absolute constant independent of $k$. 

We now summarize \eqref{eqn: one-term-in-yk-0} in $\al _1, \b _1$ and combine \eqref{eq:diffusion_low_est} to conclude 
\begin{align}
    &\B\la 
      c_1  \ang X  \ang \vc ^2 D^{\preceq (\al_1, \b_1)} g  +   \sum _{|\al_1|  +  |\b_1| = 1} \ang X \ang \vc ^2 D ^{2 \al_1, 2 \b_1} g, h
    \B\ra _{\cYb ^k}  \label{eq:diffusion_pf1}\\
    & \le -\f58 \sum _{|\al _1| + |\b _1| = 1} \sum _{|\al| + |\b| \le k} \underbrace{\nu ^{|\al| + |\b| - k}}_{ =\nu ^{|\al+ \al_1| + |\b + \b_1| -( k+1)}} 
\f{|\al|!}{\al!}    
    \iint \ang X ^{\etab} | D ^{\al + \al _1, \b + \b _1} g | ^2 \; d V d X +
     \f{1}{8} \| g \|_{\cYb^{k+1}}^2 + 
    C _k \| g \| _{\cYb ^k} ^2  \notag . 
\end{align}

For any multi-indices $\al^{\pr}, \b^{\pr}$ with  $1\leq |\al^{\pr}|+ |\b^{\pr}| \leq k+1$ and $|\al^{\pr}| \geq 1$, we have 
\[
\bal 
 I &= \sum_{ |\al_1| + |\b_1| = 1} \sum_{ |\al| + |\b| \leq k} \one_{  (\al + \al_1, \b+ \b_1) = (\al^{\pr}, \b^{\pr} ) } 
 \f{|\al|!}{\al!}   
 \geq  \sum_{ |\al_1| = 1} \sum_{ |\al| = |\al^{\pr}| - 1} \one_{  \al + \al_1 = \al^{\pr} } 
 \f{|\al|!}{\al!}    \\
 & = \sum_{1\leq i \leq 3} \one_{  \al = \al^{\pr} - \ee _i} \one_{ \al^{\pr}_i \geq 1} \f{ (|\al^{\pr}| - 1) !   }{ \al^{\pr} ! } 
 \cdot \al_i^{\pr} 
 = \f{ (|\al^{\pr}| - 1) !   }{ \al^{\pr} ! } 
 \cdot |\al^{\pr}| = \f{ |\al^{\pr}|  !   }{ \al^{\pr} !  }  .
 \eal 
\]

If $\al^{\pr} =0$, we obtain $I \geq 1 = \f{ |\al^{\pr}|  !   }{ \al^{\pr} !  }$. Thus, 
denoting $\al^{\pr} = \al +\al_1, \b^{\pr} = \b + \b_1$,  we further bound 
\eqref{eq:diffusion_pf1} as 
\begin{align*}  
  \eqref{eq:diffusion_pf1}  &  \le -\f58 \sum _{1 \le |\al^{\pr}| + |\b^{\pr}| \le k + 1} \nu ^{|\al^{\pr}| + |\b^{\pr}| - (k + 1)} 
  \f{ |\al^{\pr}|  !   }{ \al^{\pr} !  } \iint \ang X ^{\etab} | D ^{\al, \b} g | ^2 \; d V d X + \f{1}{8} \| g \|_{\cYb^{k+1}}^2 +  C _k \| g \| _{\cYb ^k} ^2 \\
    &  =  -\f58 \| g \| _{\cYb ^{k + 1}} ^2 + \f58 \nu ^{-(k + 1)} \| g \| _{\cYb} ^2 +  \f{1}{8} \| g \|_{\cYb^{k+1}}^2 + C _k \| g \| _{\cYb ^k} ^2 \\
    &  \le -\f12 \| g \| _{\cYb ^{k + 1}} ^2 + C _k  \| g \| _{\cYb ^k} ^2.
\end{align*}

 In the last equality, we used that the lowest order term in $\cYb ^k$ norm is $\nu ^{-k} \| g \| _{\cYb} ^2$. %

Using the decomposition in \eqref{eq:diffusion_main}, we complete the proof of Lemma \ref{lem: diffusion-innerprod} for $k>0$.
\end{proof}

\bibliographystyle{plain}
\bibliography{selfsimilar}

\end{document}